\documentclass[a4paper,12pt,amsart,frenchb]{article}

\usepackage{amsmath,amsbsy,amsfonts,amssymb}
\usepackage[french]{babel}
\newcommand{\ass}[2]{\vskip0.3cm\noindent
{\bf {#1}}. { \sl {#2}}\vskip0.3cm\noindent
}
  
\begin{document}

    \title{Stabilisation de la formule des traces tordue IX: propri\'et\'es des int\'egrales orbitales pond\'er\'ees $\omega$-\'equivariantes sur le corps r\'eel}
\author{J.-L. Waldspurger}
\date{8 d\'ecembre 2014}
\maketitle

{\bf Introduction}

Dans l'article pr\'ec\'edent [VIII], on a \'etabli le r\'esultat suivant. Soient $F$ un corps local non-archim\'edien de caract\'eristique nulle, $G$ un groupe r\'eductif connexe d\'efini sur $F$, $\tilde{G}$ un espace tordu sous $G$ et ${\bf a}$ un \'el\'ement de $H^1(W_{F};Z(\hat{G}))$ auquel est associ\'e un caract\`ere $\omega$ de $G(F)$. Soit de plus $\tilde{M}$ un espace de Levi de $\tilde{G}$. Pour toute fonction $f\in C_{c}^{\infty}(\tilde{G}(F))$, nulle au voisinage de certaines classes de conjugaison dites exceptionnelles, il existe une fonction cuspidale $\epsilon_{\tilde{M}}(f)$ sur $\tilde{M}(F)$ telle que
$$I^{\tilde{M}}(\gamma,\omega,\epsilon_{\tilde{M}}(f))=I_{\tilde{M}}^{\tilde{G},{\cal E}}(\gamma,\omega,f)-I_{\tilde{M}}^{\tilde{G}}(\gamma,\omega,f)$$
pour tout $\gamma\in \tilde{M}(F)$ qui est fortement r\'egulier dans $\tilde{G}$.  On renvoie aux articles pr\'ec\'edents pour les d\'efinitions de ces termes.  La fonction $\epsilon_{\tilde{M}}(f)$ est localement constante mais on n'a pas prouv\'e  qu'elle \'etait  \`a support compact. Mais elle est "de Schwartz" en un sens que l'on a d\'efini en [VIII] (en fait, on prouvera ult\'erieurement que $\epsilon_{\tilde{M}}(f)=0$ mais ce n'est pas encore d'actualit\'e). 

Dans le pr\'esent article,  on se propose de d\'emontrer  essentiellement le m\^eme r\'esultat, le corps de base \'etant maintenant ${\mathbb R}$. Ce changement de corps de base induit plusieurs diff\'erences. D'abord, on doit souvent travailler non pas avec un triplet $(G,\tilde{G},{\bf a})$, mais avec $(KG,K\tilde{G},{\bf a})$, o\`u $K\tilde{G}$ est un $K$-espace tordu, cf. [I] 1.11. On consid\`ere un $K$-espace de Levi $K\tilde{M}$ de $K\tilde{G}$ et on va construire une application $f\mapsto \epsilon_{K\tilde{M}}(f)$. L'introduction des $K$-espaces n'est qu'une complication mineure. 
Plus s\'erieusement, pour que l'application $f\mapsto \epsilon_{K\tilde{M}}(f)$ soit utilisable, on doit montrer qu'elle v\'erifie des propri\'et\'es suppl\'ementaires. D'une part, elle doit \^etre \'equivariante en un sens facile \`a pr\'eciser pour l'action du centre $\mathfrak{Z}(G)$ de l'alg\`ebre enveloppante de l'alg\`ebre de Lie de $G$. D'autre part,   on d\'efinit facilement la notion
de fonction $K$-finie sur $K\tilde{G}({\mathbb R})$ ou de fonction $K^M$-finie sur  $K\tilde{M}({\mathbb R})$. Si $f$ est $K$-finie, la fonction $\epsilon_{K\tilde{M}}(f)$ doit \^etre $K^M$-finie. En fait, l'action de $\mathfrak{Z}(G)$ comme les questions de $K$-finitude vont intervenir non seulement dans le r\'esultat mais dans la construction elle-m\^eme de l'application $\epsilon_{K\tilde{M}}$. En particulier, $\epsilon_{K\tilde{M}}(f)$ n'est d\'efinie que pour une fonction $f$ qui est $K$-finie. Il y a enfin  une derni\`ere diff\'erence avec le cas non-archim\'edien, qui porte sur les propri\'et\'es  des caract\`eres pond\'er\'es. Dans une s\'erie d'articles, Arthur avait d\'efini ceux-ci en utilisant des op\'erateurs d'entrelacement normalis\'es. En cons\'equence,  ces op\'erateurs, que l'on peut faire d\'ependre d'une variable  parcourant un espace vectoriel complexe, \'etaient m\'eromorphes avec un nombre fini d'hyperplans polaires. Dans l'article [A5], Arthur a introduit une  nouvelle d\'efinition des caract\`eres pond\'er\'es. Pour stabiliser ceux-ci,  il est n\'ecessaire d'utiliser cette nouvelle d\'efinition. F\^acheusement, sur le corps de base r\'eel, ces nouveaux op\'erateurs  peuvent avoir un nombre infini d'hyperplans polaires, parce qu'il intervient des fonctions $\Gamma$ qui ont une infinit\'e de p\^oles. Cela conduit \`a modifier substantiellement la construction de l'application $\epsilon_{K\tilde{M}}$.

Dans les deux premi\`eres sections de l'article, on \'etudie comment se transforment sous l'action de $\mathfrak{Z}(G)$ les int\'egrales orbitales pond\'er\'ees $\omega$-\'equivariantes et leurs avatars endoscopiques. Ces actions se r\'ealisent par des op\'erateurs diff\'erentiels et on montre que ces op\'erateurs sont les m\^emes pour les deux types d'int\'egrales. Il s'agit de la version tordue du r\'esultat qu'Arthur d\'emontre dans [A2]. Notre d\'emonstration est diff\'erente de celle d'Arthur. Elle consiste \`a prouver d'abord l'\'egalit\'e des op\'erateurs diff\'erentiels associ\'es \`a l'\'el\'ement de Casimir de $\mathfrak{Z}(G)$, ce que l'on fait par un calcul explicite. On montre ensuite que, pour un \'el\'ement g\'en\'eral de $\mathfrak{Z}(G)$, les op\'erateurs diff\'erentiels associ\'es \`a cet \'el\'ement v\'erifient des propri\'et\'es formelles de commutation aux m\^emes op\'erateurs associ\'es au Casimir. Ces propri\'et\'es sont suffisamment fortes pour que l'\'egalit\'e prouv\'ee par ces derniers se propage en l'\'egalit\'e des op\'erateurs diff\'erentiels associ\'es \`a un \'el\'ement quelconque de $\mathfrak{Z}(G)$. 

La section 3 traite des majorations locales v\'erifi\'ees par les  int\'egrales orbitales pond\'er\'ees $\omega$-\'equivariantes et leurs avatars endoscopiques.  Pour les premi\`eres, on dispose des majorations \'etablies par Arthur dans [A1]. On montre que celles-ci sont aussi v\'erifi\'ees par les variantes endoscopiques. 

La section 4 \'etudie les sauts des  int\'egrales orbitales pond\'er\'ees $\omega$-\'equivariantes et de leurs avatars endoscopiques. De nouveau, ces sauts sont connus gr\^ace \`a Arthur pour les premiers types d'int\'egrales. On montre que les avatars endoscopiques satisfont les m\^emes relations de sauts. La d\'emonstration ressemble beaucoup \`a celle de Shelstad concernant les int\'egrales orbitales non pond\'er\'ees, cf. [S], ainsi qu'\`a celle de l'article [A8] d'Arthur.

Dans la section 5, on d\'efinit les variantes "compactes" des int\'egrales orbitales pond\'er\'ees $\omega$-\'equivariantes. Ici, on peut travailler de nouveau avec un triplet $(G,\tilde{G},{\bf a})$. Pour $\gamma\in \tilde{M}({\mathbb R})$ qui est fortement r\'egulier dans $\tilde{G}$ et pour $f\in C_{c}^{\infty}(\tilde{G}({\mathbb R}))$,  on note $   {^cI}_{\tilde{M}}^{\tilde{G}}(\gamma,\omega,f)$ cette int\'egrale. Elle a beaucoup de propri\'et\'es communes avec l'int\'egrale plus usuelle $I_{\tilde{M}}^{\tilde{G}}(\gamma,\omega,f)$. Mais, pour $f$ fix\'ee, elle est \`a support compact en $\gamma\in \tilde{M}({\mathbb R})$, \`a conjugaison pr\`es par $M({\mathbb R})$. Pour d\'efinir $   {^cI}_{\tilde{M}}^{\tilde{G}}(\gamma,\omega,f)$, on ne peut pas utiliser la m\^eme construction que dans le cas non-archim\'edien (cette derni\`ere \'etait directement tir\'ee d'Arthur). La possible infinit\'e des hyperplans polaires des caract\`eres pond\'er\'es invalide cette construction. On commence par d\'efinir de nouveaux caract\`eres pond\'er\'es en supprimant dans la d\'efinition de [A5] les fonctions $\Gamma$ qui peuvent cr\'eer une infinit\'e de p\^oles. Ces nouveaux termes n'ont plus qu'un nombre fini d'hyperplans polaires.  La d\'efinition est adapt\'ee pour que ces nouveaux termes se stabilisent aussi ais\'ement que ceux d\'efinis par Arthur. On doit avouer que cette d\'efinition est un peu artificielle. On pourrait sans doute la rendre plus conceptuelle en utilisant les r\'ecents r\'esultats de Mezo ([Me]). Notre unique excuse pour ne pas utiliser ceux-ci est que l'on a commenc\'e ce travail alors que l'article de Mezo n'\'etait pas encore disponible.
Dans une seconde \'etape, on reprend la construction du cas non-archim\'edien, en l'appliquant \`a nos nouveaux caract\`eres pond\'er\'es. 

Les sections 6 et 7 sont consacr\'es \`a la stabilisation des int\'egrales  $   {^cI}_{\tilde{M}}^{\tilde{G}}(\gamma,\omega,f)$. La m\'ethode est similaire \`a celle du cas non-archim\'edien.

Dans la section 8, on construit l'application $\epsilon_{K\tilde{M}}$ qui est le but principal de l'article. Sa d\'efinition est similaire \`a celle du cas non-archim\'edien. Montrer que l'application ainsi construite est \'equivariante pour les actions de $\mathfrak{Z}(G)$ est facile en utilisant les r\'esultats des deux premi\`eres sections. Par contre, prouver qu'elle conserve la $K$-finitude n\'ecessite du travail suppl\'ementaire. On utilise pour cela la version sym\'etrique de la formule des traces locale, d\'evelopp\'ee dans [A9] dans le cas non tordu et \'etendue au cas g\'en\'eral dans [Moe]. Celle-ci permet d'exprimer $  I_{K\tilde{M}}^{K\tilde{G}}(\gamma,\omega,f)$  en termes de l'image de $f$ dans l'espace de Paley-Wiener. C'est-\`a-dire que l'on obtient grosso-modo une expression
$$  I_{K\tilde{M}}^{K\tilde{G}}(\gamma,\omega,f)=\int \xi(\gamma,\tilde{\pi})I^{K\tilde{G}}(\tilde{\pi},f) d\tilde{\pi},$$
o\`u $\tilde{\pi}$ parcourt les $\omega$-repr\'esentations temp\'er\'ees de $K\tilde{G}({\mathbb R})$. 
  La fonction $(\gamma,\tilde{\pi})\mapsto \xi(\gamma,\tilde{\pi})$ poss\`ede des singularit\'es qui sont assez raisonnables. Un r\'esultat analogue vaut pour  l'int\'egrale endoscopique $I_{K\tilde{M}}^{K\tilde{G},{\cal E}}(\gamma,\omega,f)$ donc aussi pour $I^{K\tilde{M}}(\gamma,\omega,\epsilon_{K\tilde{M}}(f))$. De ce r\'esultat et de l'\'equivariance pour les actions de $\mathfrak{Z}(G)$ se d\'eduit la conservation de la $K$-finitude par des m\'ethodes standard. 
  
  Notons que $\epsilon_{K\tilde{M}}(f)$ est d\'efinie pour tout $f$: contrairement au cas non-archim\'edien, on n'impose pas que $f$ s'annule au voisinage des classes de conjugaison exceptionnelles. 
  
  Dans tout l'article, le corps de base est ${\mathbb R}$ mais nos r\'esultats s'\'etendent au cas du corps de base ${\mathbb C}$. On renvoie \`a [V] section 7 pour la justification \'el\'ementaire de cette affirmation.
  
  Enfin, signalons que l'on impose dans cet article  la condition que {\bf  $\omega$ est un caract\`ere unitaire}.  Cette hypoth\`ese \'etait \'egalement n\'ecessaire dans plusieurs articles pr\'ec\'edents bien que l'on craigne de l'avoir parfois oubli\'ee. 
  
  \bigskip

 \section{Stabilisation d'une famille d'\'equations diff\'erentielles}
 \bigskip

 \subsection{Op\'erateurs diff\'erentiels}
 Dans cette section, le corps de base est ${\mathbb R}$.  Soit $(G,\tilde{G},{\bf a})$ un triplet comme en [IV] 1.1. Fixons un tore tordu maximal $\tilde{T}$ de $\tilde{G}$.  On note $\theta$ l'automorphisme $ad_{\gamma}$ de $T$ pour n'importe quel \'el\'ement $\gamma\in \tilde{T}$. On suppose que $\omega$ est trivial sur $T({\mathbb R})^{\theta}$. On pose $\tilde{T}_{\tilde{G}-reg}=\tilde{T}\cap \tilde{G}_{reg}$. 
 
 Tout \'el\'ement $H\in \mathfrak{t}({\mathbb R}) $ d\'efinit un op\'erateur diff\'erentiel $\partial_{H}$ sur $T({\mathbb R})$. Pr\'ecis\'ement, nous privil\'egions les actions \`a gauche. C'est-\`a-dire que, pour une fonction $f$ sur $T({\mathbb R})$ on a $(\partial_{H}f)(t)=\frac{d}{dx}f(exp(-xH)t)_{\vert x=0}$. 
  L'application $H\mapsto \partial_{H}$ s'\'etend lin\'eairement \`a $\mathfrak{t}=\mathfrak{t}({\mathbb C})$, puis s'\'etend  en un isomorphisme not\'e $\partial$ de l'alg\`ebre $Sym(\mathfrak{t})$ sur l'alg\`ebre des op\'erateurs diff\'erentiels sur $T({\mathbb R})$  invariants par translations \`a  gauche par $T({\mathbb R})$. 
  En fixant un \'el\'ement $\gamma\in \tilde{T}({\mathbb R})$, on peut identifier $T({\mathbb R})$ \`a $\tilde{T}({\mathbb R})$ par $t\mapsto t\gamma $. Ainsi $ Sym(\mathfrak{t})$   est aussi isomorphe  \`a l'alg\`ebre des op\'erateurs diff\'erentiels   sur $\tilde{T}({\mathbb R})$ invariants par translations \`a  gauche par $T({\mathbb R})$. 
  
  Notons $C^{\infty}(\tilde{T} ({\mathbb R}))^{\omega-inv}$ l'espace des fonctions  $\varphi:\tilde{T}({\mathbb R})\to {\mathbb C}$ qui sont $C^{\infty}$ et v\'erifient la relation  $\varphi(t^{-1}\gamma t)=\omega(t^{-1})\varphi(\gamma)$ pour tous $\gamma\in \tilde{T} ({\mathbb R})$ et $t\in T({\mathbb R})$.  Un op\'erateur diff\'erentiel sur $\tilde{T}({\mathbb R})$ invariant par translations \`a  gauche par $T({\mathbb R})$ conserve cet espace. On note $Diff^{cst}(\tilde{T}({\mathbb R}))^{\omega-inv}$ l'espace des restrictions \`a $C^{\infty}(\tilde{T}({\mathbb R}))^{\omega-inv}$ d'op\'erateurs diff\'erentiels invariants par translations \`a gauche. 
  On a introduit en [IV] 1.2 un \'el\'ement $\tilde{\mu}(\omega)\in \mathfrak{h}^*$. L'espace $\mathfrak{h}^*$ s'identifie naturellement \`a $\mathfrak{t}^*$, on peut donc consid\'erer que $\tilde{\mu}(\omega)$ appartient \`a $\mathfrak{t}^*$. C'est un \'el\'ement central, invariant par le groupe de Weyl $W$ de $T({\mathbb C})$ dans $G({\mathbb C})$. Pour $H\in (1-\theta)( \mathfrak{t})$ et $\varphi\in C^{\infty}(\tilde{T}_{\tilde{G}-reg}({\mathbb R}))^{\omega-inv}$, on a l'\'egalit\'e $\partial_{H}\varphi=<H,\tilde{\mu}(\omega)>\varphi$. Notons $I_{\theta,\omega}$ l'id\'eal de $Sym(\mathfrak{t})$ engendr\'e par les $H-<H,\tilde{\mu}(\omega)>$ pour $H\in (1-\theta)(\mathfrak{t})$ et notons $Sym(\mathfrak{t})_{\theta,\omega}$ le quotient $Sym(\mathfrak{t})/I_{\theta,\omega}$. L'application $H\mapsto \partial_{H}$ se quotiente en un isomorphisme de   $Sym(\mathfrak{t})_{\theta,\omega}$ sur $Diff^{cst}(\tilde{T}({\mathbb R}))^{\omega-inv}$. Afin de ne pas surcharger les notations,  pour $H\in Sym(\mathfrak{t})$, on notera encore  $\partial_{H}$ l'image de  cet op\'erateur dans  $Diff^{cst}(\tilde{T}({\mathbb R}))^{\omega-inv}$. Remarquons que $Sym(\mathfrak{t})_{\theta,\omega}$ s'identifie \`a l'alg\`ebre des polyn\^omes sur l'espace affine $\tilde{\mu}(\omega)+\mathfrak{t}^{*,\theta}$.  Remarquons aussi que l'homomorphisme naturel $Sym(\mathfrak{t}^{\theta})\to Sym(\mathfrak{t})_{\theta,\omega}$ est un isomorphisme. 
  
  En rempla\c{c}ant $\tilde{T}({\mathbb R})$ par $\tilde{T}_{\tilde{G}-reg}({\mathbb R})$, on d\'efinit comme ci-dessus l'espace $C^{\infty}(\tilde{T}_{\tilde{G}-reg} ({\mathbb R}))^{\omega-inv}$.
 Notons $C^{\infty}(\tilde{T}_{\tilde{G}-reg}({\mathbb R}))^{inv}$ l'espace des fonctions  $\varphi:\tilde{T}_{\tilde{G}-reg}({\mathbb R})\to {\mathbb C}$ qui sont $C^{\infty}$ et v\'erifient la relation  $\varphi(t^{-1}\gamma t)= \varphi(\gamma)$ pour tous $\gamma\in \tilde{T}_{\tilde{G}-reg}({\mathbb R})$ et $t\in T({\mathbb R})$. On pose
 $$Diff^{\infty}(\tilde{T}_{\tilde{G}-reg}({\mathbb R}))^{\omega-inv}=C^{\infty}(\tilde{T}_{\tilde{G}-reg}({\mathbb R}))^{inv}\otimes_{{\mathbb C}}Diff^{cst}(\tilde{T}({\mathbb R}))^{\omega-inv} .$$
 Cet espace agit naturellement sur $C^{\infty}(\tilde{T}_{\tilde{G}-reg}({\mathbb R}))^{\omega-inv}$: pour $\varphi'\in C^{\infty}(\tilde{T}_{\tilde{G}-reg}({\mathbb R}))^{inv}$, $D\in  Diff^{cst}(\tilde{T}({\mathbb R}))^{\omega-inv}$ et $\varphi\in C^{\infty}(\tilde{T}_{\tilde{G}-reg}({\mathbb R}))^{\omega-inv}$, $(\varphi'\otimes D)(\varphi)$ est la fonction d\'efinie par $((\varphi'\otimes D)(\varphi))(\gamma)=\varphi'(\gamma)(D\varphi)(\gamma)$. Il est clair que l'on peut munir $Diff^{\infty}(\tilde{T}_{\tilde{G}-reg}({\mathbb R}))^{\omega-inv}$ d'une unique structure d'alg\`ebre telle que cette action devienne une action d'alg\`ebres. On note $(D,D')\mapsto D\circ D'$ le produit pour cette structure.
 
 {\bf Notation.} Un \'el\'ement $\delta\in Diff^{\infty}(\tilde{T}_{\tilde{G}-reg}({\mathbb R}))^{\omega-inv}$ sera plut\^ot consid\'er\'e comme une fonction $C^{\infty}$ de $\tilde{T}_{\tilde{G}-reg}({\mathbb R})$ dans $ Diff^{cst}(\tilde{T}({\mathbb R}))^{\omega-inv}$. On notera $\delta(\gamma)$ sa valeur en un point $\gamma$. Pour $\varphi\in C^{\infty}(\tilde{T}_{\tilde{G}-reg}({\mathbb R}))^{\omega-inv}$, on notera $\delta(\gamma)\varphi(\gamma)$ la valeur en $\gamma$ de la fonction $\delta \varphi$. 
 
 \bigskip

 Les ensembles $T({\mathbb C})$ et $\tilde{T}({\mathbb C})$  sont des vari\'et\'es alg\'ebriques complexes. On peut donc parler de fonctions polynomiales, rationnelles, holomorphes ou m\'eromorphes sur ces ensembles. Par exemple, l'espace des polyn\^omes  sur $T({\mathbb C})$ est engendr\'e lin\'eairement par les caract\`eres alg\'ebriques, c'est-\`a-dire les \'el\'ements de $X^*(T)$.   Introduisons le tore $T'=T/(1-\theta)(T)$ et  l'espace $\tilde{T}'=\tilde{T}/(1-\theta)(T)$. On peut de m\^eme d\'efinir les espaces de fonctions polynomiales, rationnelles etc... sur $T'({\mathbb C})$ ou $\tilde{T}'({\mathbb C})$. L'ensemble $\tilde{T}_{\tilde{G}-reg}$ est invariant par conjugaison par $T$, donc aussi par produit avec $(1-\theta)(T)$. Cela permet d'introduire le sous-ensemble $\tilde{T}'_{\tilde{G}-reg}=\tilde{T}_{\tilde{G}-reg}/(1-\theta)(T)$ de $\tilde{T}'$. On appelle fonction rationnelle r\'eguli\`ere sur $\tilde{T}'_{\tilde{G}-reg}({\mathbb C})$ une fonction rationnelle sur $\tilde{T}'({\mathbb C})$ qui n'a pas de p\^ole dans $\tilde{T}'_{\tilde{G}-reg}({\mathbb C})$. Notons $Pol(\tilde{T}'_{\tilde{G}-reg}({\mathbb C}))$ l'espace de ces fonctions.  La restriction de $\tilde{T}_{\tilde{G}-reg}({\mathbb C})$ \`a $\tilde{T}_{\tilde{G}-reg}({\mathbb R})$ d\'efinit  une application
 $$Pol (\tilde{T}'_{\tilde{G}-reg}({\mathbb C}))\to C^{\infty}(\tilde{T}_{\tilde{G}-reg}({\mathbb R}))^{inv}.$$
 Elle est injective. Son image est conserv\'ee par l'action de tout op\'erateur diff\'erentiel sur $\tilde{T}({\mathbb R})$ invariant par translations \`a gauche. Posons
 $$Diff^{reg}(\tilde{T}_{\tilde{G}-reg}({\mathbb R}))^{\omega-inv}=Pol(\tilde{T}'_{\tilde{G}-reg}({\mathbb C}))\otimes_{{\mathbb C}}Diff^{cst}(\tilde{T}({\mathbb R}))^{\omega-inv} .$$
 L'application ci-dessus permet d'identifier cet espace \`a une sous-alg\`ebre de $Diff^{\infty}(\tilde{T}_{\tilde{G}-reg}({\mathbb R}))^{\omega-inv}$.

   \bigskip

  \subsection{Les \'equations diff\'erentielles}
 Pour la suite de la section, on fixe  un espace de Levi $\tilde{M}$ de $\tilde{G}$ et   un sous-tore tordu maximal $\tilde{T}$ de $\tilde{M}$.  On suppose que $\omega$ est trivial sur $T({\mathbb R})^{\theta}$.
 
 Rappelons que l'on note $\mathfrak{U}(G)$ l'alg\`ebre enveloppante de la complexifi\'ee de l'alg\`ebre de Lie $\mathfrak{g}$ de $G$ et que l'on note $\mathfrak{Z}(G)$ le centre de $\mathfrak{U}(G)$. On dispose d'homomorphismes
 $$\begin{array}{ccccccc}\mathfrak{Z}(G)&\to& \mathfrak{Z}(M),&\,\,& \mathfrak{Z}(G)&\to&\mathfrak{Z}(T)\\ z&\mapsto&z_{M},&\,\,&z&\mapsto&z_{T}.\\ \end{array}$$
 L'alg\`ebre $\mathfrak{Z}(T)$ s'identifie \`a  $Sym(\mathfrak{t})$.  L'homomorphisme $z\mapsto z_{T}$ identifie $\mathfrak{Z}(G)$ \`a la sous-alg\`ebre des invariants $Sym(\mathfrak{t})^W$, o\`u $W$ est le groupe de Weyl de $G$ relatif \`a $T$.
 L'alg\`ebre $\mathfrak{U}(G)$ agit sur $C_{c}^{\infty}(\tilde{G}({\mathbb R}))$ via les translations \`a gauche. Il s'en d\'eduit une action de $\mathfrak{Z}(G)$ sur $I(\tilde{G}({\mathbb R}),\omega)$. Cette action se quotiente en l'action d'une alg\`ebre quotient $\mathfrak{Z}(G)_{\theta,\omega}$. Celle-ci est isomorphe \`a $Sym(\mathfrak{t})_{\theta,\omega}^{W^{\theta}}$.

  Pour simplifier, on fixe des mesures de Haar sur tous les groupes intervenant. Pour $\gamma\in \tilde{M}({\mathbb R})\cap \tilde{G}_{reg}({\mathbb R})$ et pour $f\in C_{c}^{\infty}(\tilde{G}({\mathbb R}))$, on sait d\'efinir l'int\'egrale orbitale pond\'er\'ee $\omega$-\'equivariante $I_{\tilde{M}}^{\tilde{G}}(\gamma,\omega,f)$, cf. [W2] 6.5. La fonction $\gamma\mapsto I_{\tilde{M}}^{\tilde{G}}(\gamma,\omega,f)$ appartient \`a $C^{\infty}(\tilde{T}_{\tilde{G}-reg}({\mathbb R}))^{\omega-inv}$. 

 Arthur d\'emontre  en [A1] proposition 11.1 et [A2] paragraphe 1 qu'il existe une  unique application lin\'eaire 
  $$\begin{array}{ccc}  \mathfrak{Z}(G)&\to& Diff^{\infty}(\tilde{T}_{\tilde{G}-reg}({\mathbb R}))^{\omega-inv}\\ z&\mapsto& \delta_{\tilde{M}}^{\tilde{G}}(z)\\ \end{array}$$
    qui v\'erifie la propri\'et\'e suivante:

  - pour tout $f\in C_{c}^{\infty}(\tilde{G}({\mathbb R}))$, tout $\gamma\in \tilde{T}_{\tilde{G}-reg}({\mathbb R})$ et tout $z\in \mathfrak{Z}({\mathbb R})$, on a l'\'egalit\'e
  $$(1) \qquad I_{\tilde{M}}^{\tilde{G}}(\gamma,\omega,zf)=\sum_{\tilde{L}\in {\cal L}(\tilde{M})}\delta_{\tilde{M}}^{\tilde{L}}(\gamma,z_{\tilde{L}})I_{\tilde{L}}^{\tilde{G}}(\gamma,\omega,f).$$
  On a not\'e $\delta_{\tilde{M}}^{\tilde{L}}(\gamma,z_{\tilde{L}})$ la valeur au point $\gamma$ de $\delta_{\tilde{M}}^{\tilde{L}}(z_{\tilde{L}})$ et on a utilis\'e la notation introduite dans le paragraphe pr\'ec\'edent.

 On v\'erifie formellement  que    $\delta_{\tilde{M}}^{\tilde{G}}(z)$ ne d\'epend que de l'image de $z$ dans $\mathfrak{Z}(G)_{\theta,\omega}$. On  pourra donc consid\'erer que $\delta_{\tilde{M}}^{\tilde{G}}(z)$ est d\'efini pour $z\in \mathfrak{Z}(G)_{\theta,\omega}$.

 D'autre part, d'apr\`es [A1] lemme 12.4,
 
 (2) si  $\tilde{M}=\tilde{G}$, on a simplement $\delta_{\tilde{G}}^{\tilde{G}}(\gamma,z)=\partial_{z_{T}}$ pour tout $z\in \mathfrak{Z}(G)$ et tout $\gamma\in \tilde{T}_{\tilde{G}-reg}({\mathbb R})$.
  
   Montrons que
 
 (3) pour $z,z'\in \mathfrak{Z}(G)$, on a l'\'egalit\'e
 $$\delta_{\tilde{M}}^{\tilde{G}}(zz')=\sum_{\tilde{L}\in {\cal L}(\tilde{M})}\delta_{\tilde{M}}^{\tilde{L}}(z_{\tilde{L}})\circ\delta_{\tilde{L}}^{\tilde{G}}(z').$$
 
 Preuve. Notons $\underline{\delta}_{\tilde{M}}^{\tilde{G}}(zz')$ le membre de droite de cette \'egalit\'e.  Pour $f\in  C_{c}^{\infty}(\tilde{G}({\mathbb R}))$, notons $\psi_{\tilde{M}}(f)$ la fonction $\gamma\mapsto I_{\tilde{M}}^{\tilde{G}}(\gamma,\omega,f)$. L'\'egalit\'e (1) prend la forme
 $$(4) \qquad \psi_{\tilde{M}}(zf)=\sum_{\tilde{L}\in {\cal L}(\tilde{M})}\delta_{\tilde{M}}^{\tilde{L}}(z_{\tilde{L}})\psi_{\tilde{L}}(f).$$
 Rempla\c{c}ons dans cette \'egalit\'e $f$ par $z'f$. D\'eveloppons ensuite chaque terme $\psi_{\tilde{L}}(z'f)$ par la m\^eme \'egalit\'e o\`u l'on remplace $z$ par $z'$ et $\tilde{M}$ par $\tilde{L}$. On obtient facilement l'\'egalit\'e 
 $$\psi_{\tilde{M}}(zz'f)=\sum_{\tilde{L}\in {\cal L}(\tilde{M})}\underline{\delta}_{\tilde{M}}^{\tilde{L}}((zz')_{\tilde{L}})\psi_{\tilde{L}}(f).$$
 En raisonnant par r\'ecurrence, on peut supposer que 
$$\underline{\delta}_{\tilde{M}}^{\tilde{L}}((zz')_{\tilde{L}}) =\delta_{\tilde{M}}^{\tilde{L}}((zz')_{\tilde{L}})$$
pour tout $\tilde{L}\not=\tilde{G}$. En comparant l'\'egalit\'e pr\'ec\'edente avec l'\'egalit\'e (1) o\`u $z$ est remmplac\'e par $zz'$, on en d\'eduit l'\'egalit\'e
$$(5) \qquad (\delta_{\tilde{M}}^{\tilde{G}}(zz')-\underline{\delta}_{\tilde{M}}^{\tilde{G}}(zz'))\psi_{\tilde{G}}(f)=0.$$
Soit $\gamma\in \tilde{T}_{\tilde{G}-reg}({\mathbb R})$. L'espaces des germes au point $\gamma$ des fonctions $\psi_{\tilde{G}}(f)$ quand $f$ d\'ecrit $C_{c}^{\infty}(\tilde{G}({\mathbb R}))$ est \'egal \`a celui des germes des fonctions $\varphi$ pour $\varphi\in C^{\infty}(\tilde{T}_{\tilde{G}-reg}({\mathbb R}))^{\omega-inv}$. On en d\'eduit que $0$ est le seul \'el\'ement de $ Diff^{\infty}(\tilde{T}_{\tilde{G}-reg}({\mathbb R}))^{\omega-inv}$ qui annule toute fonction $\psi_{\tilde{G}}(f)$. L'\'egalit\'e (5) entra\^{\i}ne alors la conclusion de (3). $\square$ 

\bigskip 
 
 \subsection{Propri\'et\'es des op\'erateurs $\delta_{\tilde{M}}^{\tilde{G}}(z)$}

 On note $\Sigma(T)$ l'ensemble des racines de $T$ dans $G$. A tout \'el\'ement $\alpha\in \Sigma(T)$, on associe le plus petit entier $n_{\alpha}\geq1$ tel que $\theta^{n_{\alpha}}(\alpha)=\alpha$. On d\'efinit  un \'el\'ement $N\alpha\in X^*(T)$ par $N\alpha=\sum_{k=0,...n_{\alpha}-1}\theta^k(\alpha)$. Il se descend en un \'el\'ement de $X^*(T')$. On pose $\alpha_{T'}=N\alpha$ si $\alpha$ est de type 1 ou 3, $\alpha_{T'}=2N\alpha$ si $\alpha$ est de type 2 (cf. [I] 1.6). On pose $\Sigma(T')=\{\alpha_{T'}; \alpha\in \Sigma(T)\}$. C'est un syst\`eme de racines en g\'en\'eral non r\'eduit.  Fixons un sous-groupe de Borel $B$ de $G$ de tore maximal $T$ et invariant par $\theta$. Notons  $\Delta$ la base de $\Sigma(T)$ associ\'ee \`a $B$.  On dira qu'une suite $(t'_{k})_{k\in {\mathbb N}}$ d'\'el\'ements de $T'({\mathbb C})$ tend vers l'infini selon $\Delta$ si et seulement si, $lim_{k\to \infty}\vert \alpha_{T'}(t'_{k})\vert =\infty$ pour tout $\alpha\in \Delta$.  Consid\'erons une fonction $\varphi:\tilde{T}'_{\tilde{G}-reg}({\mathbb C})\to {\mathbb C}$ et  un nombre complexe $c$.  Introduisons la condition
 
 (1) pour tout $\gamma'\in \tilde{T}'({\mathbb C})$ et pour toute suite $(t'_{k})_{k\in {\mathbb N}}$ d'\'el\'ements de $T'({\mathbb C})$ tendant vers l'infini selon $\Delta$ et telle que $t'_{k}\gamma'\in \tilde{T}'_{\tilde{G}-reg}({\mathbb C})$ pour tout $k$, on a  $lim_{k\to\infty}\varphi(t'_{k}\gamma')=c$.
 
 On peut remplacer "pour tout $\gamma'\in \tilde{T}'({\mathbb C})$" par "il existe  $\gamma'\in \tilde{T}'({\mathbb C})$ tel que...", on obtient une condition \'equivalente. 
 On notera $lim_{\gamma'\to_{\Delta}\infty}\varphi(\gamma')=c$ si cette condition est v\'erifi\'ee. 
 
 Fixons une base ${\cal B}$ de $Sym(\mathfrak{t})_{\theta,\omega}$. Pour tout $z\in \mathfrak{Z}(G)$, l'op\'erateur $\delta_{\tilde{M}}^{\tilde{G}}(z)$ s'\'ecrit de fa\c{c}on unique comme une somme finie
 $$\delta_{\tilde{M}}^{\tilde{G}}(z)=\sum_{U\in {\cal B}}q_{\tilde{M}}^{\tilde{G}}(U;z)\partial_{U},$$
 o\`u les $q_{\tilde{M}}^{\tilde{G}}(U;z)$ sont des  \'el\'ements de $C^{\infty}(\tilde{T}_{\tilde{G}-reg}({\mathbb R}))^{\omega-inv}$.  On note $q_{\tilde{M}}^{\tilde{G}}(U;\gamma,z)$ la valeur de cette fonction en un point $\gamma$. 
 
 \ass{Proposition}{ Supposons $\tilde{M}\not=\tilde{G}$. Soient $z\in \mathfrak{Z}(G)$ et $U\in {\cal B}$. 
 
 (i) La fonction $q_{\tilde{M}}^{\tilde{G}}(U;z)$ est la restriction \`a $\tilde{T}_{\tilde{G}-reg}({\mathbb R})$ d'une fonction rationnelle r\'eguli\`ere sur $\tilde{T}'_{\tilde{G}-reg}({\mathbb C})$, que l'on note encore 
 $q_{\tilde{M}}^{\tilde{G}}(U;z)$.
 
 (ii) Pour toute base $\Delta$ de $\Sigma(T)$, on a $lim_{\gamma'\to_{\Delta}\infty}q_{\tilde{M}}^{\tilde{G}}(U;\gamma',z)=0$.}
 
 Cette proposition sera prouv\'ee en 1.5. Remarquons que le (i) est vrai aussi si $\tilde{M}=\tilde{G}$ en vertu de 1.2(2).

 \bigskip

\subsection{Rappels sur l'action adjointe}
Pour deux racines $\alpha,\beta\in \Sigma(T)$, disons qu'elles sont dans la m\^eme orbite si et seulement s'il existe $m\in {\mathbb N}$ tel que $\beta=\theta^m\alpha$. On note $(\alpha) $ l'orbite de $\alpha$.   Fixons un ensemble de repr\'esentants $\Sigma(T)_{\theta}\subset \Sigma(T)$ des orbites. On compl\`ete $(B,T)$ en une paire de Borel \'epingl\'ee ${\cal E}$ en fixant des \'el\'ements non nuls $E_{\alpha}\in \mathfrak{g}_{\alpha}$ pour $\alpha\in \Delta$, o\`u $\mathfrak{g}_{\alpha}\subset \mathfrak{g}$ est la droite radicielle associ\'ee \`a $\alpha$. On fixe $e\in Z(\tilde{G},{\cal E})$ et on d\'efinit l'automorphisme $\theta=ad_{e}$ de $G$.
 Pour $\alpha\in \Sigma(T)_{\theta}$, fixons un \'el\'ement non nul $E_{\alpha}$ de  $\mathfrak{g}_{\alpha}$. On suppose que c'est l'\'el\'ement d\'ej\`a fix\'e si $\alpha\in \Delta$. Pour $m=1,...,n_{\alpha}-1$, posons $E_{\theta^m\alpha}=\theta^m(E_{\alpha})$. C'est  un \'el\'ement non nul de $\mathfrak{g}_{\theta^m\alpha}$. D'apr\`es [KS] 1.3, on a $ \theta^{n_{\alpha}}(E_{\alpha})=\epsilon_{\alpha}E_{\alpha}$, o\`u $\epsilon_{\alpha}=1$ si $\alpha$ est de type 1 ou 2 et $\epsilon_{\alpha}=-1$ si $\alpha$ est de type 3. Soit $\gamma\in \tilde{T}_{\tilde{G}-reg}({\mathbb C})$. On \'ecrit $\gamma=t_{\gamma}e$, avec $t_{\gamma}\in T$. On pose $(\alpha)(\gamma)=\epsilon_{\alpha}\prod_{m=0,...,n_{\alpha}-1}\theta^{m}(\alpha)(t_{\gamma})$.  On  en fixe une racine $n_{\alpha}$-i\`eme $\nu_{(\alpha)}(\gamma)$ et on note $\boldsymbol{\zeta}_{n_{\alpha}}$ le groupe des racines $n_{\alpha}$-i\`emes de $1$. Pour $\zeta\in \boldsymbol{\zeta}_{n_{\alpha}}$, posons
 $$E((\alpha),\gamma,\zeta)=\sum_{m=0,...,n_{\alpha}-1}\nu_{(\alpha)}(\gamma)^{-m}\zeta^{-m}(\prod_{j=0,...,m}\theta^{j}(\alpha)(t_{\gamma}))E_{\theta^{m}(\alpha)}.$$
 On v\'erifie que $ad_{\gamma}(E((\alpha),\gamma,\zeta))=\nu_{(\alpha)}(\gamma)\zeta E((\alpha),\gamma,\zeta)$. La famille $(E(\alpha),\gamma,\zeta))_{\zeta\in \boldsymbol{\zeta}_{n_{\alpha}}}$ est une base de l'espace $\mathfrak{g}_{(\alpha)}=\sum_{m=0,...,n_{\alpha}-1}\mathfrak{g}_{\theta^m(\alpha)}$. Notons $\mathfrak{q}$ le sous-espace de $\mathfrak{g}$ engendr\'e par les alg\`ebres de Lie des radicaux unipotents de $B$ et du sous-groupe de Borel oppos\'e $\bar{B}$. La famille  $(E(\alpha),\gamma,\zeta))_{\alpha\in \Sigma(T)_{\theta},\zeta\in \boldsymbol{\zeta}_{n_{\alpha}}}$ est une base de $\mathfrak{q}$.

 Si l'on remplace $\gamma$ par $t\gamma $, avec $t\in T^{\theta,0}$, on peut supposer $\nu_{(\alpha)}(t\gamma )=\alpha_{res}(t)\nu_{(\alpha)}(\gamma) $, o\`u $\alpha_{res}$ est la restriction de $\alpha$ \`a $T^{\theta,0}$. On a alors $E((\alpha),t\gamma ,\zeta)=E((\alpha),\gamma,\zeta)$. Les familles ci-dessus sont donc ind\'ependantes de $t$. Remarquons que $E((\alpha),\gamma,\zeta)$ est vecteur propre pour l'action $ad_{t\gamma }$, de valeur propre $\zeta\alpha_{res}(t)\nu_{(\alpha)}(\gamma) $ et est aussi vecteur propre pour l'action $ad_{t}$, de valeur propre $\alpha_{res}(t)$.

  Posons
 $$D_{\star}^{\tilde{G}}(\gamma)=\prod_{\alpha\in \Sigma(T)_{\theta}}\prod_{\zeta\in \boldsymbol{\zeta}_{n_{\alpha}}}(1-\zeta\nu_{(\alpha)}(\gamma))=\prod_{\alpha\in \Sigma(T)_{\theta}}(1-(\alpha)(\gamma)).$$
 
 {\bf Remarque.} Ce terme est d\'efini pour $\gamma\in \tilde{T}_{\tilde{G}-reg}({\mathbb C})$ et appartient \`a ${\mathbb C}^{\times}$.  Il appartient \`a ${\mathbb R}^{\times}$ si $\gamma\in \tilde{T}({\mathbb R})$. En [I] 2.4, on a d\'efini un terme $D^{\tilde{G}}(\gamma)$ pour $\gamma\in \tilde{T}_{\tilde{G}-reg}({\mathbb R})$. On a l'\'egalit\'e
 $$D^{\tilde{G}}(\gamma)=\vert D_{\star}^{\tilde{G}}(\gamma) det((1-\theta)_{\mathfrak{t}/\mathfrak{t}^{\theta}})\vert_{{\mathbb R}} .$$

 \bigskip 
 
 Puisque la restriction de  $D_{\star}^{\tilde{G}}$ \`a $\tilde{T}_{\tilde{G}-reg}({\mathbb R})$ est \`a valeurs r\'eelles, on peut choisir une racine carr\'ee $(D_{\star}^{\tilde{G}})^{1/2}$ qui soit une fonction $C^{\infty}$ sur  $\tilde{T}_{\tilde{G}-reg}({\mathbb R})$. Pour $H\in Sym(\mathfrak{t})$ ou $H\in Sym(\mathfrak{t})_{\theta,\omega}$, on d\'efinit l'op\'erateur $\partial^{\star}_{H}\in Diff^{\infty}(\tilde{T}_{\tilde{G}-reg})^{\omega-inv}$ par l'\'egalit\'e
 $$\partial^{\star}_{H}=(D_{\star}^{\tilde{G}})^{1/2}\circ \partial_{H}\circ (D_{\star}^{\tilde{G}})^{-1/2}.$$
 Cela ne d\'epend pas du choix de la racine carr\'ee $(D_{\star}^{\tilde{G}})^{1/2}$. Ecrivons
 $$\partial^{\star}_{H}=\sum_{U\in {\cal B}}b(U;H)\partial_{U}.$$
 Les termes $b(U;H)$ sont des \'el\'ements de $C^{\infty}(\tilde{T}_{\tilde{G}-reg}({\mathbb R}))^{\omega-inv}$. On note $b(U;\gamma,H)$ leurs valeurs en un point $\gamma$.
 
 \ass{Lemme}{Soient $H\in Sym(\mathfrak{t})$ et $U\in {\cal B}$. Alors la fonction $b(U;H)$ est combinaison lin\'eaire de  produits de fonctions $\gamma\mapsto \frac{1}{(1-(\alpha)(\gamma))^m}$ pour $m\in {\mathbb N}$ et $\alpha\in \Sigma(T)_{\theta}$, $\alpha>0$.}
 
 Preuve. Supposons d'abord $H\in \mathfrak{t}$.  Pour $\alpha\in \Sigma(T)_{\theta}$ et $\gamma\in \tilde{T}_{\tilde{G}-reg}({\mathbb R})$, on calcule 
 $$(1) \qquad \partial_{H}(\alpha)(\gamma)=-<N\alpha,H>(\alpha)(\gamma).$$
  Pour tout $m\in {\mathbb N}$, on en d\'eduit l'\'egalit\'e
 $$(2) \qquad \partial_{H} \frac{1}{(1-(\alpha)(\gamma))^m}= \frac{-m<N\alpha,H>(\alpha)(\gamma)}{(1-(\alpha)(\gamma))^{m+1}}= \frac{m<N\alpha,H>}{(1-(\alpha)(\gamma))^m}-\frac{m<N\alpha,H>}{(1-(\alpha)(\gamma))^{m+1}}.$$
  En introduisant l'ordre sur les racines relatif \`a la base $\Delta$ fix\'ee, on peut r\'ecrire
 $$D_{\star}^{\tilde{G}}(\gamma)=\prod_{\alpha\in \Sigma(T)_{\theta}, \alpha>0}(1-(\alpha)(\gamma))(1-(\alpha)(\gamma)^{-1})=\prod_{\alpha\in \Sigma(T)_{\theta}, \alpha>0}-(\alpha)(\gamma)^{-1}(1-(\alpha)(\gamma))^2.$$
 Fixons $\gamma$. Pour $\gamma'$ voisin de $\gamma$, on peut choisir des racines carr\'ees $(\alpha)(\gamma')^{1/2}$ qui sont $ C^{\infty}$ en $\gamma'$.   On peut supposer que
 $$  D_{\star}^{\tilde{G}}(\gamma')^{1/2}=\prod_{\alpha\in \Sigma(T)_{\theta}, \alpha>0}i(\alpha)(\gamma')^{-1/2}(1-(\alpha)(\gamma'))$$
 au voisinage de $\gamma$, 
 o\`u $i$ est une racine carr\'ee de $-1$.   En utilisant cette formule  et la relation (1), on obtient
 $$(3) \qquad \partial^{\star}_{H}(\gamma)=\partial_{H}+\sum_{\alpha\in \Sigma(T)_{\theta},\alpha>0}<N\alpha,H>(\frac{1}{2}-\frac{1}{1-(\alpha)(\gamma)}).$$
Cela v\'erifie le lemme pour $H\in \mathfrak{t}$.  
 
 En raisonnant par r\'ecurrence sur le d\'egr\'e de $H$, il reste \`a v\'erifier l'assertion suivante. Soient $ H\in \mathfrak{t}$ et $H'\in Sym(\mathfrak{t})$. Supposons que le lemme soit v\'erifi\'e pour $H'$.  Alors il l'est pour $HH'$. La formule (2) montre que $\partial_{H}$ conserve l'ensemble des fonctions qui sont combinaisons lin\'eaires de  produits de fonctions $\gamma\mapsto \frac{1}{(1-(\alpha)(\gamma))^m}$ pour $m\in {\mathbb N}$ et $\alpha\in \Sigma(T)_{\theta}$, $\alpha>0$. L'assertion r\'esulte alors facilement de (3). $\square$

 \bigskip
 
 \subsection{Une application d'Harish-Chandra}
 Notons $Tens(\mathfrak{q})$ l'alg\`ebre tensorielle de $\mathfrak{q}$. Il contient le sous-espace $Sym(\mathfrak{q})$. Soit $\gamma\in \tilde{T}_{\tilde{G}-reg}({\mathbb C})$. On d\'efinit une application lin\'eaire
$$\Gamma_{\gamma}: Tens(\mathfrak{q})\otimes_{{\mathbb C}}Sym(\mathfrak{t})\to \mathfrak{U}(G)$$
par
$$(1) \qquad \Gamma_{\gamma}(X_{1}...X_{a}\otimes H)=(R_{ad_{\gamma}(X_{1})}-L_{X_{1}})...(R_{ad_{\gamma}(X_{a})}-L_{X_{a}})H$$
pour $X_{1},...,X_{a}\in \mathfrak{q}$ et $U\in Sym(\mathfrak{t})$. On a not\'e $L_{X}$ et $R_{X}$ les applications $Y\mapsto XY$ et $Y\mapsto YX$ de $\mathfrak{U}(G)$ dans lui-m\^eme. D'apr\`es [HC] lemme 22, cette application se restreint en un isomorphisme d'espaces vectoriels
$$Sym(\mathfrak{q})\otimes_{{\mathbb C}}Sym(\mathfrak{t})\to \mathfrak{U}(G).$$
Fixons des bases $(Y_{k})_{k\in {\mathbb N}}$ de $Sym(\mathfrak{ q})$ et $(H_{l})_{l\in {\mathbb N}}$ de $Sym(\mathfrak{t})$.   Pour tout $V\in \mathfrak{U}(G)$ et tout $\gamma\in \tilde{T}_{\tilde{G}-reg}({\mathbb C})$, on  peut \'ecrire de fa\c{c}on unique
$$(2) \qquad V=\sum_{k,l}a_{k,l}(\gamma,V)\Gamma_{\gamma}(Y_{k}\otimes H_{l})$$
 o\`u $a_{k,l}(\gamma,V)\in {\mathbb C}$ et $a_{k,l}(\gamma,V)=0$ pour presque tout $(k,l)$. 
 
 {\bf Remarques.} (3)  Les preuves de Harish-Chandra concernent le cas non tordu, mais leur extension au cas tordu est imm\'ediate. On les reprendra d'ailleurs partiellement ci-dessous. 

(4) Ici et dans la suite, il y a de l\'eg\`eres diff\'erences avec les r\'ef\'erences cit\'ees dues au fait que nous provil\'egions les actions par translations \`a gauche alors que les auteurs cit\'es utilisent les translations \`a droite. 
\bigskip

Supposons que $Y_{0}=1$ tandis  que le terme constant de $Y_{k}$ est nul si $k\geq1$.  

\ass{Lemme}{Soient $V\in \mathfrak{U}(G)$ et $k,l\in {\mathbb N}$. 

(i) La fonction $\gamma\mapsto a_{k,l}(\gamma,V)$ est rationnelle sur $\tilde{T}_{\tilde{G}-reg}({\mathbb C})$. Plus pr\'ecis\'ement, il existe un entier $n_{k,l}\in {\mathbb N}$ tel que la fonction $\gamma\mapsto D_{\star}^{\tilde{G}}(\gamma)^{n_{k,l}}a_{k,l}(\gamma,V)$ se prolonge en un polyn\^ome sur $\tilde{T}({\mathbb C})$. 

(ii) Supposons que $k\geq1$ et que $V$ est invariant par l'action adjointe de $T^{\theta,0}({\mathbb C})$. Soient $\gamma\in \tilde{T}_{\tilde{G}-reg}({\mathbb C})$, $\Delta$ une base de $\Sigma(T)$ et $(t_{j})_{j\in {\mathbb N}}$ une suite d'\'el\'ements de $T^{\theta,0}({\mathbb C})$. On suppose que $t_{j}\gamma\in \tilde{T}_{\tilde{G}-reg}({\mathbb C})$ pour tout $j$ et que l'image de la suite dans $T'({\mathbb C})$ tend vers l'infini selon $\Delta$. Alors $lim_{j\to \infty}a_{k,l}(t_{j}\gamma,V)=0$.}

Preuve. Si on change les bases $(Y_{k})_{k\in {\mathbb N}}$  et $(H_{l})_{l\in {\mathbb N}}$ en d'autres bases $(Y'_{k})_{k\in {\mathbb N}}$ et $(H'_{l})_{l\in {\mathbb N}}$, on obtient de nouveaux coefficients $a'_{k,l}(\gamma,V)$ qui se d\'eduisent des pr\'ec\'edents par un syst\`eme d'\'equations lin\'eaires. Il est clair que les assertions de l'\'enonc\'e   pour ces nouveaux coefficients sont \'equivalentes aux m\^emes assertions pour les anciens. On a donc le choix des  bases  $(Y_{k})_{k\in {\mathbb N}}$ et $(H_{l})_{l\in {\mathbb N}}$. Fixons $\gamma\in \tilde{T}_{\tilde{G}-reg}({\mathbb C})$. La fonction $t\mapsto a_{k,l}(t\gamma,V)$ est d\'efinie pour presque tout $t\in T ({\mathbb C})$ (pr\'ecis\'ement pour les $t$ tels que $t\gamma\in \tilde{T}_{\tilde{G}-reg}({\mathbb C})$). On la restreint \`a $T^{\theta,0}({\mathbb C})$. On va commencer par prouver

(5) la fonction $t\mapsto a_{k,l}(t\gamma,V)$ est rationnelle sur $T^{\theta,0}({\mathbb C})$; il existe un entier $n_{k,l}\in {\mathbb N}$ tel que la fonction $t\mapsto D_{\star}^{\tilde{G}}(t\gamma)^{n_{k,l}}a_{k,l}(t\gamma,V)$ soit r\'eguli\`ere sur $T^{\theta,0}({\mathbb C})$. 

 Comme ci-dessus, cette assertion ne d\'epend pas du choix des bases. 
En cons\'equence, on suppose que $(H_{l})_{l\in {\mathbb N}}$ est form\'e d'\'el\'ements homog\`enes et que la base $(Y_{k})_{k\in {\mathbb N}}$   est form\'ee des sym\'etris\'es des \'el\'ements $E((\alpha),\gamma,\zeta)$ introduits dans le paragraphe pr\'ec\'edents. Rappelons que l'homomorphisme de sym\'etrisation identifie $Sym(\mathfrak{q})$ \`a un sous-espace de $\mathfrak{U}(G)$ et que, modulo cette identification, on a l'isomorphisme $\mathfrak{U}(G)=Sym(\mathfrak{t})\otimes_{{\mathbb C}}Sym(\mathfrak{q})$. Ainsi, la famille $(H_{l}Y_{k})_{k,l\in {\mathbb N}}$ est une base de $\mathfrak{U}(G)$. Elle est form\'ee de vecteurs propres pour l'action de $T^{\theta,0}({\mathbb C})$.   Par lin\'earit\'e, on peut supposer que $V$ est l'un de ces \'el\'ements de base $H_{l_{0}}Y_{k_{0}}$,  et on peut raisonner par r\'ecurrence sur le degr\'e de cet \'el\'ement. Ecrivons simplement $Y=Y_{k_{0}}$, $U=U_{l_{0}}$. Si $k_{0}=0$, on a simplement $V=\Gamma_{t\gamma }(Y_{0}\otimes U)$ pour tout $t\in T^{\theta,0}({\mathbb C})$ et l'assertion (5) est  claire. Remarquons que $V$ est   invariant par l'action adjointe de $T^{\theta,0}({\mathbb C})$ et que l'assertion (ii) du lemme est tout aussi claire. Supposons $k_{0}\geq1$. L'\'el\'ement $Y$ est le sym\'etris\'e de $X_{1}...X_{a}$, o\`u chaque $X_{i}$ est un \'el\'ement $E((\alpha),\gamma,\zeta)$. Ainsi, pour chaque $i=1,...,a$, il existe une racine $\alpha_{i}\in \Sigma(T)$ et un \'el\'ement $\nu_{i}\in {\mathbb C}^{\times}$ de sorte que $ad_{t}(X_{i})=\alpha_{i,res}(t)X_{i}$ pour tout $t\in T^{\theta,0}({\mathbb C})$  et $ad_{\gamma }(X_{i})=\nu_{i}X_{i}$. On a
$$\Gamma_{t\gamma}(X_{1}...X_{a}\otimes U)=(R_{ad_{t\gamma}(X_{1})}-L_{X_{1}})...(R_{ad_{t\gamma}(X_{a})}-L_{X_{a}})U$$
$$=(R_{ad_{t\gamma}(X_{1})}-R_{X_{1}}+R_{X_{1}}-L_{X_{1}})...(R_{ad_{t\gamma}(X_{a})}-R_{X_{a}}+R_{X_{a}}-L_{X_{a}})U$$
$$=(c_{1}(t)R_{X_{1}}+R_{X_{1}}-L_{X_{1}})...(c_{a}(t)R_{X_{a}}+R_{X_{a}}-L_{X_{a}})U,$$
o\`u $c_{i}(t)=\nu_{i}\alpha_{i,res}(t)-1$. On d\'eveloppe cette expression en s\'eparant les termes $c_{i}(t)R_{X_{i}}$ des $R_{X_{i}}-L_{X_{i}}$. Ces deux op\'erateurs envoient un vecteur propre pour l'action adjointe de $T^{\theta,0}({\mathbb C})$ sur un tel vecteur propre. De plus, le second op\'erateur fait baisser strictement le degr\'e. On obtient une expression
$$\Gamma_{t\gamma }(X_{1}...X_{a}\otimes U)=(\prod_{i=1,...,a}c_{i}(t))UX_{a}...X_{1}+\sum_{b=0,...,a-1}\sum_{1\leq i_{1}<...i_{b}\leq a}(\prod_{j=1,...,b}c_{i_{j}}(t))V_{i_{1},...,i_{b}},$$
o\`u les $V_{i_{1},...,i_{b}}$ sont des \'el\'ements de $\mathfrak{U}(G)$ de degr\'e strictement inf\'erieur \`a celui de $V$ et qui sont propres pour l'action adjointe de $T^{\theta,0}$. En sym\'etrisant, on obtient une formule analogue
$$\Gamma_{t\gamma }(Y\otimes U)=(\prod_{i=1,...,a}c_{i}(t))V+\sum_{b=0,...,a-1}\sum_{1\leq i_{1}<...i_{b}\leq a}(\prod_{j=1,...,b}c_{i_{j}}(t))V'_{i_{1},...,i_{b}},$$
ou encore
$$(6) \qquad V=(\prod_{i=1,...,a}c_{i}(t))^{-1}\Gamma_{t\gamma }(Y\otimes U)-\sum_{b=0,...,a-1}\sum_{1\leq i_{1}<...i_{b}\leq a}(\prod_{i=1,...,a; i\not=i_{1},...,i_{b}}c_{i}(t))^{-1}V'_{i_{1},...,i_{b}}.$$
Pour tous $k,l\in {\mathbb N}$, posons
$$a'_{k,l}(t\gamma ;V)=\sum_{b=0,...,a-1}\sum_{1\leq i_{1}<...i_{b}\leq a}(\prod_{i=1,...,a; i\not=i_{1},...,i_{b}}c_{i}(t))^{-1}a_{k,l}(t\gamma ;V'_{i_{1},...,i_{b}}).$$
Alors
$$(7) \qquad a_{k,l}(t\gamma ;V)=\left\lbrace\begin{array}{cc}-a'_{k,l}(t\gamma ;V);&\text{ si }(k,l)\not=(k_{0},l_{0}),\\ (\prod_{i=1,...,a}c_{i}(t))^{-1}-a'_{k_{0},l_{0}}(t\gamma ;V),&\text{ si }(k,l)=(k_{0},l_{0}).\\ \end{array}\right.$$
Pour tout $i=1,...,a$, la fonction $c_{i}(t)$  est une fonction   rationnelle sur $T^{\theta,0}({\mathbb C})$.  La fonction $D_{\star}^{\tilde{G}}(t\gamma)c_{i}(t)$ est r\'eguli\`ere sur cet ensemble. Jointe aux propri\'et\'es des $a_{k,l}(t\gamma ;V'_{i_{1},...,i_{b}})$ connues par r\'ecurrence, cette propri\'et\'e entra\^{\i}ne que $a_{k,l}(t\gamma ;V)$ est aussi une fonction rationnelle sur $T^{\theta,0}({\mathbb C})$ et qu'il existe un entier $n_{k,l}\in {\mathbb N}$ tel que $D_{\tilde{G}}(t\gamma)^{n_{k,l}}a_{k,l}(t\gamma;V)$ est r\'eguli\`ere sur cet ensemble. Cela prouve (5). 

Prouvons maintenant le (ii) de l'\'enonc\'e. On reprend le calcul ci-dessus en supposant $V=U_{l_{0}}Y_{k_{0}}=UY$. On a d\'ej\`a remarqu\'e que (ii) \'etait v\'erifi\'ee si $k_{0}=0$. On suppose $k_{0}\geq1$. 
Puisque $V$ est invariant par l'action adjointe de $T^{\theta,0}({\mathbb C})$, on a 

(8) $\sum_{i=1,...,a}\alpha_{i,res}=0$. 

\noindent Il est clair sur la d\'efinition (1) qu'alors $\Gamma_{t\gamma }(Y\otimes U)$ est aussi invariant par l'action adjointe de $T^{\theta,0}({\mathbb C})$. Puisque les $V'_{i_{1},...,i_{b}}$ sont propres pour cette action, on peut aussi bien supprimer de la formule (6) ceux qui ne sont pas invariants. Les \'el\'ements restants v\'erifient alors le  (ii) de l'\'enonc\'e par r\'ecurrence. 
On a
$$(9) \qquad lim_{ j\to \infty}c_{i}(t_{j})^{-1}=\left\lbrace\begin{array}{cc}0,&\text{ si }\alpha_{i}>0,\\ -1,&\text{ si }\alpha_{i}<0,\\ \end{array}\right.$$
o\`u la positivit\'e des racines est relative \`a $\Delta$. Des propri\'et\'es des $a_{k,l}(t\gamma ;V'_{i_{1},...,i_{b}})$ connues par r\'ecurrence  r\'esulte alors que, si $k\geq1$, on a $lim_{ j\to \infty}a'_{k,l}(t_{j}\gamma ;V)=0$. Pour conclure, il reste \`a prouver que $lim_{ j\to \infty}(\prod_{i=1,...,a}c_{i}(t_{j}))^{-1}=0$. Mais cela r\'esulte de (8) et (9), car (8) assure qu'il y a au moins un $i$ tel que $\alpha_{i}>0$. Cela  prouve le (ii) de l'\'enonc\'e.

Venons-en \`a la preuve du (i). De nouveau, on a le choix de la base. On suppose maintenant que  $H_{l}$ est homog\`ene pour tout $l$ et que $Y_{k}$ est propre pour l'action de $T({\mathbb C})$, de caract\`ere propre $y_{k}$. Par lin\'earit\'e, on peut aussi supposer que $V$ est propre pour cette action, de caract\`ere propre $v$. Pour $\gamma\in \tilde{T}_{\tilde{G}-reg}({\mathbb C})$ et $t\in T({\mathbb C})$, on d\'eduit de (2) l'\'egalit\'e
$$ad_{t}(V)=\sum_{k,l}a_{k,l}(\gamma,V)\Gamma_{ad_{t}(\gamma)}(ad_{t}(Y_{k})\otimes ad_{t}(H_{l})),$$
ou encore
$$V=\sum_{k,l}a_{k,l}(\gamma,V)v(t)^{-1}y_{k}(t)\Gamma_{(1-\theta)(t)\gamma}(Y_{k}\otimes H_{l}).$$
En comparant avec (2) appliqu\'ee \`a $(1-\theta)(t)\gamma$, on obtient
$$(10) \qquad  a_{k,l}((1-\theta)(t)\gamma,V)=v(t)^{-1}y_{k}(t)a_{k,l}(\gamma,V).$$
Ceci implique que $a_{k,l}(\gamma,V)=0$ si les restrictions de $v$ et $y_{k}$ \`a $T^{\theta}({\mathbb C})$ sont diff\'erentes. Supposons qu'elles soient \'egales. Alors le caract\`ere $v^{-1}y_{k}$ se descend en un caract\`ere de $(1-\theta)(T({\mathbb C}))$ que l'on note $y'_{k}$.  Fixons $\gamma_{0}\in \tilde{T}_{\tilde{G}-reg}({\mathbb C})$. L'homomorphisme
$$\begin{array}{ccc}(1-\theta)(T({\mathbb C}))\times T^{\theta,0}({\mathbb C})&\to& \tilde{T}({\mathbb C})\\ (t_{1},t_{2})&\mapsto& t_{1}t_{2}\gamma_{0}\\ \end{array}$$
est surjectif.  Soit $n_{k,l}$ un entier v\'erifiant (5). La fonction $\gamma\mapsto D_{\star}^{\tilde{G}}(\gamma)^{n_{k,l}}a_{k,l}(\gamma,V)$ se rel\`eve en une fonction sur $(1-\theta)(T({\mathbb C}))\times T^{\theta,0}({\mathbb C})$. Au point $(t_{1},t_{2})$, celle-ci vaut $y'_{k}(t_{1})D_{\star}^{\tilde{G}}(t_{2}\gamma_{0})^{n_{k,l}}a_{k,l}(t_{2}\gamma_{0},V)$. L'assertion (5) entra\^{\i}ne que cette fonction se prolonge en une fonction polynomiale sur $(1-\theta)(T({\mathbb C}))\times T^{\theta,0}({\mathbb C})$.     Donc la fonction $\gamma\mapsto D_{\star}^{\tilde{G}}(\gamma)^{n_{k,l}}a_{k,l}(\gamma,V)$  se prolonge elle aussi en une fonction polynomiale sur $\tilde{T}({\mathbb C})$. 
 $\square$

Les fonctions $\gamma\mapsto a_{k,l}(\gamma,V)$ sont rationnelles. On sait qu'en tout $\gamma$, il n'y a qu'un nombre fini de couples $(k,l)$ pour lesquels $a_{k,l}(\gamma,V)\not=0$. Il en r\'esulte qu'il n'y a qu'un nombre fini de couples $(k,l)$ pour lesquels la fonction $\gamma\mapsto a_{k,l}(\gamma,V)$ n'est pas nulle.

 \bigskip
 
 \subsection{Preuve de la proposition 1.3}
 On va rappeler la construction de nos op\'erateurs diff\'erentiels, d'apr\`es [A1] paragraphe 12.  
 Soit $\tilde{P}\in {\cal P}(\tilde{M})$. Notons $\mathfrak{p}^1$ l'alg\`ebre de Lie du sous-groupe de $P$ dont les points r\'eels sont les $x\in P({\mathbb R})$ tels que $H_{\tilde{P}}(x)=0$. Notons $\mathfrak{u}_{\bar{P}}$ le radical nilpotent de l'alg\`ebre de Lie du sous-groupe parabolique oppos\'e $\bar{P}$ et notons $\mathfrak{U}_{\bar{P}}$ la sous-alg\`ebre qu'il engendre dans $\mathfrak{U}(G)$. On a la d\'ecompositon en somme directe
$$\mathfrak{U}(G)=\mathfrak{p}^1\mathfrak{U}(G)\oplus Sym(\mathfrak{a}_{\tilde{M}})\mathfrak{U}_{\bar{P}}\mathfrak{u}_{\bar{P}}\oplus Sym(\mathfrak{a}_{\tilde{M}}).$$
Elle est compatible \`a l'action adjointe de $T$. Pour tout $Y\in \mathfrak{U}(G)$, \'ecrivons $Y=Y'+Y''+\mu_{\tilde{P}}(Y)$ conform\'ement \`a cette d\'ecomposition ci-dessus.    Notons $d$ la diff\'erence des dimensions de ${\cal A}_{\tilde{M}}$ et ${\cal A}_{\tilde{G}}$. Notons $\mu_{\tilde{P},d}(Y)$ la composante homog\`ene de degr\'e $d$ de $\mu_{\tilde{P}}(Y)$. On peut consid\'erer $Sym(\mathfrak{a}_{\tilde{M}})$ comme l'alg\`ebre des polyn\^omes sur ${\cal A}_{\tilde{M},{\mathbb C}}^*$. Pour $\nu\in {\cal A}_{\tilde{M},{\mathbb C}}^*$, notons $<\mu_{\tilde{P},d}(Y),\nu>$ la valeur en $\nu$ du polyn\^ome associ\'e \`a $\mu_{\tilde{P},d}(Y)$. Posons
$$c_{\tilde{M}}^{\tilde{G}}(Y)=\sum_{\tilde{P}\in {\cal P}(\tilde{M})}<\mu_{\tilde{P},d}(Y),\nu>\epsilon_{\tilde{P}}^{\tilde{G}}(\nu),$$
o\`u $\nu$ est un \'el\'ement assez r\'egulier de ${\cal A}_{\tilde{M},{\mathbb C}}^*$ et $\epsilon_{\tilde{P}}^{\tilde{G}}$ est la fonction usuelle de la th\'eorie des $(\tilde{G},\tilde{M})$-familles (la notation d'Arthur est $(\theta_{\tilde{P}}^{\tilde{G}})^{-1}$; notre notation  est reprise de [LW] p. 28). L'expression ci-dessus ne d\'epend pas de $\nu$.

Comme dans le paragraphe pr\'ec\'edent, fixons des bases   $(H_{l})_{l\in {\mathbb N}}$ de $Sym(\mathfrak{t})$ et $(Y_{k})_{k\in {\mathbb N}}$ de $Sym(\mathfrak{q})$. On suppose que $H_{l}$ est homog\`ene pour tout $l$, que $Y_{0}=1$ et que, pour $k\geq1$, $Y_{k}$ est de terme constant nul et est propre pour l'action adjointe de $T({\mathbb C})$, de caract\`ere propre $y_{k}$.

Soient $z\in \mathfrak{Z}(G)$ et $\gamma\in \tilde{T}_{\tilde{G}-reg}({\mathbb R})$. D'apr\`es [A1] lemme 12.1, on a l'\'egalit\'e
$$(1) \qquad  \delta_{\tilde{M}}^{\tilde{G}}(\gamma,z)=\sum_{k\geq1, l\geq0}c_{\tilde{M}}^{\tilde{G}}(Y_{k}) a_{k,l}(\gamma,z)\partial^{\star}_{H_{l}}(\gamma).$$
Puisque la projection $Y\mapsto \mu_{\tilde{P}}(Y)$ est \'equivariante pour les actions de $T$, on a  $\mu_{\tilde{P}}(Y_{k})=0$ si $y_{k}\not=1$. A fortiori, le terme index\'e par $(k,l)$ dans la somme ci-dessus est nul si $y_{k}\not=1$.  En utilisant les notations de 1.4, on obtient
$$\delta_{\tilde{M}}^{\tilde{G}}(\gamma,z)=\sum_{U\in {\cal B}}\sum_{k\geq1; y_{k}=1}\sum_{l\geq0}c_{\tilde{M}}^{\tilde{G}}(Y_{k}) a_{k,l}(\gamma,z)b(U;\gamma,H_{l})\partial_{U}.$$
Donc, pour tout $U\in {\cal B}$, 
$$(2) \qquad q_{\tilde{M}}^{\tilde{G}}(U;\gamma,z)=\sum_{k\geq1; y_{k}=1}\sum_{l\geq0}c_{\tilde{M}}^{\tilde{G}}(Y_{k}) a_{k,l}(\gamma,z)b(U;\gamma,H_{l}).$$

Fixons $U\in {\cal B}$.  Consid\'erons un couple $(k,l)$ intervenant ci-dessus. Le terme $a_{k,l}(\gamma,z)$ est d\'efini pour $\gamma\in \tilde{T}_{\tilde{G}-reg}({\mathbb C})$. Puisque $y_{k}=1$ et que $z$ est invariante par l'action adjointe de $T({\mathbb C})$, le lemme 1.5(i) et la relation 1.5(10) entra\^{\i}nent que $\gamma\mapsto a_{k,l}(\gamma,z)$ se quotiente en une fonction rationnelle r\'eguli\`ere sur $\tilde{T}'_{\tilde{G}-reg}({\mathbb C})$. Le lemme 1.4 montre qu'il en est de m\^eme de la fonction $\gamma\mapsto b(U;\gamma,H_{l})$. Donc la fonction $\gamma\mapsto q_{\tilde{M}}^{\tilde{G}}(U;\gamma,z)$ se prolonge en une fonction rationnelle r\'eguli\`ere sur $\tilde{T}'_{\tilde{G}-reg}({\mathbb C})$. 
 
  Fixons $\gamma$. L'homomorphisme 
$$\begin{array}{ccc}T^{\theta,0}({\mathbb C})&\to& \tilde{T}'({\mathbb C})\\ t&\mapsto& t\gamma\\ \end{array}$$
est surjectif.  Pour d\'emontrer le (ii) de la proposition 1.3, il suffit de prouver la propri\'et\'e (3) suivante. Fixons une base $\Delta$ de $\Sigma(T)$ et une suite $(t_{j})_{j\in {\mathbb N}}$ d'\'el\'ements de $T^{\theta,0}({\mathbb C})$. Supposons que $t_{j}\gamma\in \tilde{T}_{\tilde{G}-reg}({\mathbb C})$ pour tout $j$ et que l'image de la suite $(t_{j})_{j\in {\mathbb N}}$ dans $T'({\mathbb C})$ tend vers l'infini selon $\Delta$. Alors  on a 

(3) $lim_{j\to \infty} q_{\tilde{M}}^{\tilde{G}}(U;t_{j}\gamma,z)=0$.  

Consid\'erons un couple $(k,l)$ intervenant dans (2). Puisque $k\geq1$, le lemme 1.5(ii) implique que $lim_{j\to \infty}a_{k,l}(t_{j}\gamma,z)=0$. D'apr\`es le lemme 1.4, le terme $b(U;t_{j}\gamma,H_{l})$ est combinaison lin\'eaire de produits de termes 
  $(1- (\alpha)(\gamma)\alpha_{res}(t_{j})^{n_{\alpha}})^{-m}$ pour des $\alpha\in \Sigma(T)_{\theta}$, $\alpha>0$ et des $m\in {\mathbb N}$. Un tel produit  a une limite quand $j$ tend vers l'infini: sa limite est nulle sauf si la fonction est constante. Ces r\'esultats et la formule (2) impliquent (3), ce qui ach\`eve la d\'emonstration de la proposition. 
   $\square$
   
   Une  cons\'equence formelle de la formule (1) et de l'alg\'ebricit\'e des constructions est le comportement de nos op\'erateurs diff\'erentiels par conjugaison. Enon\c{c}ons cette propri\'et\'e sous la forme g\'en\'erale qui nous servira plus loin. Supposons que $\tilde{G}$ soit une composante connexe d'un $K$-espace $K\tilde{G}$, cf. [I] 1.11, et que $\tilde{M}$ soit une composante connexe d'un $K$-espace de Levi $K\tilde{M}$, cf. [I] 3.5. Disons que $\tilde{G}$ et $\tilde{M}$ sont les composantes $\tilde{G}_{p}$ et $\tilde{M}_{p}$. Soit $q\in \Pi^{\tilde{M}}$. On a alors une composante $\tilde{G}_{q}$ de $K\tilde{G}$, une composante $\tilde{M}_{q}$ de $K\tilde{M}$ et un isomorphisme   $\tilde{\phi}_{p,q}:\tilde{G}_{q}\to \tilde{G}_{p}$ d\'efini sur ${\mathbb C}$ qui envoie $\tilde{M}_{q}$ sur $\tilde{M}_{p}$.
   Supposons donn\'e  un sous-tore maximal $\tilde{T}_{q}$ de $\tilde{M}_{q}$ et un \'el\'ement $x\in M_{p}$ tel que    $ad_{x}\circ \tilde{\phi}_{p,q}(\tilde{T}_{q})=\tilde{T}$. On note $\tilde{\iota}:\tilde{T}_{q}\to \tilde{T}$ la restriction de $ad_{x}\circ \tilde{\phi}_{p,q}$. On suppose que $\tilde{\iota}$ est d\'efini sur ${\mathbb R}$. 
    Les centres $\mathfrak{Z}(G)$ et $\mathfrak{Z}(G_{q})$ s'identifient. On peut donc d\'efinir $\delta_{\tilde{M}_{q}}^{\tilde{G}_{q}}(\gamma_{q},z)$ pour $\gamma_{q}\in \tilde{T}_{q,\tilde{G}_{q}-reg}({\mathbb R})$. On a alors l'\'egalit\'e
$$(4) \qquad \tilde{\iota}\circ \delta_{\tilde{M}_{q}}^{\tilde{G}_{q}}(\gamma_{q},z)\circ \tilde{\iota}^{-1}=\delta_{\tilde{M}}^{\tilde{G}}(\tilde{\iota}(\gamma_{q}),z)$$
pour tout $\gamma_{q}\in \tilde{T}_{q,\tilde{G}_{q}-reg}({\mathbb R})$.

\bigskip

\subsection{L'op\'erateur de Casimir}
On a fix\'e en [IV] 1.1 une forme bilin\'eaire sur $X_{*}(T)\otimes_{{\mathbb Z}}{\mathbb R}$ d\'efinie positive, invariante par le groupe de Weyl $W$ et par $\theta$.   On la prolonge en une forme ${\mathbb C}$-bilin\'eaire sur $\mathfrak{t}$. On peut alors identifier  le dual $\mathfrak{t}^*$ \`a $\mathfrak{t}$. Pour $X\in X_{*}(T)\otimes_{{\mathbb Z}}{\mathbb R}\subset \mathfrak{t}$ ou $X\in X^{*}(T)\otimes_{{\mathbb Z}}{\mathbb R}\subset \mathfrak{t}^*$, on a $(X,X)\geq0$ et on note $\vert X\vert $ la racine carr\'ee positive ou nulle de $(X,X)$.

D'autre part, la forme sur $\mathfrak{t}$ se prolonge en une unique forme ${\mathbb C}$-bilin\'eaire sur $\mathfrak{g}$ invariante par l'action adjointe de $G$. Elle est aussi invariante par $ad_{\gamma}$ pour tout $\gamma\in \tilde{G}$. On la note $(.,.)$. De cette forme se d\'eduit un \'el\'ement de Casimir $\Omega\in \mathfrak{Z}(G)$, cf. [Va] I.2.3.  Fixons  une base $(H_{j})_{i=1,...,n}$ de $X_{*}(T)\otimes_{{\mathbb Z}}{\mathbb R}$. Supposons-la orthonorm\'ee. On a fix\'e
en 1.4 des \'el\'ements $E_{\alpha}$ pour $\alpha\in \Sigma(T)$. On a forc\'ement $(E_{\alpha},E_{\beta})=0$ si $\beta\not=-\alpha$ et on peut supposer $(E_{\alpha},E_{-\alpha})=1$. Alors
$$\Omega=(\sum_{j=1,...,n}H_{j}^2)+(\sum_{\alpha\in \Sigma(T), \alpha>0}E_{\alpha}E_{-\alpha}+E_{-\alpha}E_{\alpha}).$$
On calcule classiquement
$$(1) \qquad \Omega_{T}=\sum_{j=1,...,n}H_{j}^2-\rho_{B}(H_{j})^2,$$
o\`u $\rho_{B}$ est la demi-somme des racines positives.

Notons plus pr\'ecis\'ement $\Sigma^{G}(T)$ l'ensemble not\'e jusqu'alors $\Sigma(T)$. L'ensemble $\Sigma^{M}(T)$ est invariant par $\theta$ et on peut supposer $\Sigma^{M}(T)_{\theta}\subset
\Sigma^{G}(T)_{\theta}$. Ainsi $\Sigma^{G}(T)_{\theta}-\Sigma^{M}(T)_{\theta}$ est un ensemble de repr\'esentants des orbites dans $\Sigma^{G}(T)-\Sigma^{M}(T)$. On peut aussi supposer que $\Sigma^{G}(T)_{\theta}$ est invariant par $\alpha\mapsto -\alpha$ et que $\nu_{(-\alpha)}(\gamma)=\nu_{(\alpha)}(\gamma)^{-1}$ pour tout $\gamma\in \tilde{T}_{\tilde{G}-reg}({\mathbb C})$. Posons $d=a_{\tilde{M}}-a_{\tilde{G}}$. Supposons $d\geq1$. D\'efinissons la fonction $C_{\tilde{M}}^{\tilde{G}}$ sur $\tilde{T}_{\tilde{G}-reg}({\mathbb C})$ par 

- si $d\geq2$, $C_{\tilde{M}}^{\tilde{G}}=0$;

- si $d=1$,
$$C_{\tilde{M}}^{\tilde{G}}(\gamma)=-\sum_{\alpha\in \Sigma^{G}(T)_{\theta}-\Sigma^{M}(T)_{\theta}}\vert (N\alpha)_{\vert {\cal A}_{\tilde{M}}}\vert n_{\alpha}(1- (\alpha)(\gamma))^{-1}(1- (\alpha)(\gamma)^{-1})^{-1}.$$
Le terme $(\alpha)(\gamma)$ a \'et\'e d\'efini en 1.4. On a not\'e $(N\alpha)_{\vert {\cal A}_{\tilde{M}}}$ la projection orthogonale de $N\alpha$ sur ${\cal A}_{\tilde{M}}$, ce dernier espace s'identifiant \`a un sous-espace de $\mathfrak{t}$. 

\ass{Proposition}{On suppose $d\geq1$. Alors, pour tout $\gamma\in \tilde{T}_{\tilde{G}-reg}({\mathbb R})$, l'op\'erateur $\delta_{\tilde{M}}^{\tilde{G}}(\gamma,\Omega)$ est la multiplication par $C_{\tilde{M}}^{\tilde{G}}(\gamma)$.}

Preuve. Soit $\gamma\in \tilde{T}_{\tilde{G}-reg}$. Utilisons les \'el\'ements $E((\alpha),\gamma,\zeta)$ de 1.4. On v\'erifie que, pour tout $\alpha\in \Sigma(T)_{\theta}$, on a l'\'egalit\'e
$$\sum_{m=0,...,n_{\alpha}-1}E_{\theta^m(\alpha)}E_{-\theta^m(\alpha)}=n_{\alpha}^{-1}\sum_{\zeta\in \boldsymbol{\zeta}_{n_{\alpha}}} E((\alpha),\gamma,\zeta)E((-\alpha),\gamma,\zeta^{-1}).$$
Ainsi
 $$\Omega=(\sum_{j=1,...,n}H_{j}^2)+\sum_{\alpha\in \Sigma(T)_{\theta}, \alpha>0}n_{\alpha}^{-1}
\sum_{\zeta\in \boldsymbol{\zeta}_{n_{\alpha}}} Y((\alpha),\gamma,\zeta),$$
o\`u
$$Y((\alpha),\gamma,\zeta)= E((\alpha),\gamma,\zeta)E((-\alpha),\gamma,\zeta^{-1})+ E((-\alpha),\gamma,\zeta^{-1})E((\alpha),\gamma,\zeta).$$
Pour tous $\alpha$, $\zeta$, on calcule
$$\Gamma_{\gamma}( Y((\alpha),\gamma,\zeta)\otimes 1)=(1-\nu_{(\alpha)}(\gamma)\zeta)(1-\nu_{(\alpha)}(\gamma)^{-1}\zeta^{-1}) Y((\alpha),\gamma,\zeta)$$
$$+(\nu_{(\alpha)}(\gamma)\zeta-\nu_{(\alpha)}(\gamma)^{-1}\zeta^{-1})[E((\alpha),\gamma,\zeta),E((-\alpha),\gamma,\zeta^{-1})].$$
Le dernier terme appartient \`a $\mathfrak{t}$ puisqu'il est fixe par $ad_{\gamma}$. On en d\'eduit
$$\Omega=U+\sum_{\alpha\in \Sigma(T)_{\theta},\alpha>0}n_{\alpha}^{-1}$$
$$\sum_{\zeta\in \boldsymbol{\zeta}_{n_{\alpha}}} (1-\nu_{(\alpha)}(\gamma)\zeta)^{-1}(1-\nu_{(\alpha)}(\gamma)^{-1}\zeta^{-1})^{-1}\Gamma_{\gamma}( Y((\alpha),\gamma,\zeta)\otimes 1),$$
avec $U\in Sym(\mathfrak{t})$. En appliquant 1.6(1), on obtient que $\delta_{\tilde{M}}^{\tilde{G}}(\gamma,\Omega)$ est la multiplication par
$$(2) \qquad \sum_{\alpha\in \Sigma(T)_{\theta},\alpha>0}n_{\alpha}^{-1}\sum_{\zeta\in \boldsymbol{\zeta}_{n_{\alpha}}} (1-\nu_{(\alpha)}(\gamma)\zeta)^{-1}(1-\nu_{(\alpha)}(\gamma)^{-1}\zeta^{-1})^{-1}c_{\tilde{M}}^{\tilde{G}}(Y((\alpha),\gamma,\zeta)).$$
Soit $\tilde{P}\in {\cal P}(\tilde{M})$, $\alpha\in \Sigma(T)_{\theta}$ et $\zeta\in \boldsymbol{\zeta}_{n_{\alpha}}$. Si $\alpha\in\Sigma^{M}(T)$, on a $E(\pm(\alpha),\gamma,\zeta)\in \mathfrak{p}^1$ et $\mu_{\tilde{P}}(Y((\alpha),\gamma,\zeta))=0$. Supposons $\alpha\not\in \Sigma^{M}(T)$. Soit $\xi\in\{\pm 1\}$ tel que $\xi\alpha$ soit positif pour $P$. Alors $E(\xi(\alpha),\gamma,\zeta)\in \mathfrak{p}^1$ et on calcule 
$$\mu_{\tilde{P}}(Y((\alpha),\gamma,\zeta))=\xi[E(-(\alpha),\gamma,\zeta^{-1}),E((\alpha),\gamma,\zeta)].$$
Puisqu'on a d\'ej\`a remarqu\'e que ce terme appartenait \`a $\mathfrak{t}$, on le calcule facilement:
$$\mu_{\tilde{P}}(Y((\alpha),\gamma,\zeta))=\xi \sum_{m=0,...,n_{\alpha}}[E_{-\theta^m(\alpha)},E_{\theta^m(\alpha)}].$$
Soit $\beta\in \Sigma(T)$. Parce que notre forme bilin\'eaire est invariante par l'action adjointe, on a
$$([E_{-\beta},H],E_{\beta})+ (H,[E_{-\beta},E_{\beta}])=0$$
pour tout $H\in \mathfrak{t}$. On a aussi $[E_{-\beta},H]=\beta(H)E_{-\beta}$. Puisqu'on a choisi nos \'el\'ements de sorte que $(E_{-\beta},E_{\beta})=1$, on obtient $(H,[E_{-\beta},E_{\beta}])=-\beta(H)$. Modulo notre identification de $\mathfrak{t}^*$ \`a $\mathfrak{t}$, on obtient $[E_{-\beta},E_{\beta}]=-\beta$. D'o\`u
$$\mu_{\tilde{P}}(Y((\alpha),\gamma,\zeta))=-\xi\sum_{m=0,...,n_{\alpha}}\theta^m(\alpha)=-\xi N\alpha.$$
Ce terme est homog\`ene de degr\'e $1$. Si $d\geq2$, on en d\'eduit $c_{\tilde{M}}^{\tilde{G}}(Y((\alpha),\gamma,\zeta))=0$ ce qui d\'emontre l'\'enonc\'e dans ce cas. Supposons $d=1$. Fixons $\tilde{P}\in {\cal P}(\tilde{M})$ et un \'el\'ement $X\in {\cal A}_{\tilde{M}}^{\tilde{G}}$ positif pour $P$ et tel que $\vert X\vert =1$. Il r\'esulte des d\'efinitions que
$$c_{\tilde{M}}^{\tilde{G}}(Y((\alpha),\gamma,\zeta))=2(\mu_{\tilde{P}}(Y((\alpha),\gamma,\zeta),X))=-2\xi(N\alpha,X).$$
Il r\'esulte de la d\'efinition de $X$ que ceci est \'egal \`a $-2\vert (N\alpha)_{\vert {\cal A}_{\tilde{M}}}\vert $. Cela transforme (2) en
$$-2\sum_{\alpha\in \Sigma^{G}(T)_{\theta}-\Sigma^{M}(T)_{\theta},\alpha>0}\vert (N\alpha)_{\vert {\cal A}_{\tilde{M}}}\vert n_{\alpha}^{-1}\sum_{\zeta\in \boldsymbol{\zeta}_{n_{\alpha}}} (1-\nu_{(\alpha)}(\gamma)\zeta)^{-1}(1-\nu_{(\alpha)}(\gamma)^{-1}\zeta^{-1})^{-1}.$$
On montre ais\'ement que, pour presque tout $x\in {\mathbb C}$, on a l'\'egalit\'e
$$n_{\alpha}^{-1}\sum_{\zeta\in \boldsymbol{\zeta}_{n_{\alpha}}}(1-x\zeta)^{-1}(1-x^{-1}\zeta^{-1})^{-1}=n_{\alpha}(1-x^{n_{\alpha}})^{-1}(1-x^{-n_{\alpha}})^{-1}.$$
On se rappelle que $\nu_{(\alpha)}(\gamma)^{n_{\alpha}}= (\alpha)(\gamma)$. Alors l'expression ci-dessus devient
$$-2\sum_{\alpha\in \Sigma^{G}(T)_{\theta}-\Sigma^{M}(T)_{\theta},\alpha>0}\vert (N\alpha)_{\vert {\cal A}_{\tilde{M}}}\vert n_{\alpha}(1-(\alpha)(\gamma))(1-(\alpha)(\gamma)^{-1})^{-1}.$$
Les expressions \'etant sym\'etriques en $\alpha$ et $-\alpha$, on peut supprimer le facteur $2$ en supprimant la restriction $\alpha>0$ dans la sommation. Le terme ci-dessus devient $C_{\tilde{M}}^{\tilde{G}}(\gamma)$.
Cela prouve la proposition. $\square$

{\bf Remarque.} L'alg\`ebre $\mathfrak{Z}(G)$ contient $Sym(\mathfrak{z}(G))$, o\`u $\mathfrak{z}(G)$ est l'alg\`ebre de Lie du centre de $G$.  Il r\'esulte facilement de  1.6(1) que, pour $d\geq1$, on a $\delta_{\tilde{M}}^{\tilde{G}}(\gamma,z)=0$ pour $z\in Sym(\mathfrak{z}(G))$.  La proposition pr\'ec\'edente reste vraie si l'on remplace $\Omega$ par n'importe quel \'el\'ement de $\Omega+Sym(\mathfrak{z}(G))$. 

\bigskip

\subsection{Variante avec caract\`ere central}
On suppose $(G,\tilde{G},{\bf a})$ quasi-d\'eploy\'e et \`a torsion int\'erieure. On consid\`ere des extensions compatibles
$$1\to C_{\natural}\to G_{\natural}\to G\to 1\text{ et }\tilde{G}_{\natural}\to \tilde{G}$$
o\`u $C_{\natural}$ est un tore central induit et o\`u $\tilde{G}_{\natural}$ est encore \`a torsion int\'erieure. Soit $\lambda_{\natural}$ un caract\`ere de $C_{\natural}({\mathbb R})$. On note $\tilde{M}_{\natural}$ et $\tilde{T}_{\natural}$ les images r\'eciproques dans $\tilde{G}_{\natural}$ de $\tilde{M}$ et $\tilde{T}$. 

Puisqu'il n'y a pas de torsion, on note simplement $Diff^{cst}(\tilde{T}_{\natural}({\mathbb R}))$ l'alg\`ebre not\'ee $Diff^{cst}(\tilde{T}_{\natural}({\mathbb R}))^{\omega-inv}$ en 1.1. Elle agit sur l'espace des fonctions $C^{\infty}$ sur $\tilde{T}_{\natural}({\mathbb R})$. Restreignons-nous aux fonctions $f$ qui v\'erifient pour tout $\gamma_{\natural}\in \tilde{T}_{\natural}({\mathbb R})$ la propri\'et\'e 

(1) $f(\gamma_{\natural}c)=\lambda_{\natural}(c)^{-1}f(\gamma_{\natural})$ pour tout $c\in C_{\natural}({\mathbb R})$.

\noindent Alors l'action de $Diff^{cst}(\tilde{T}_{\natural}({\mathbb R}))$ se quotiente en une action d'une alg\`ebre quotient $Diff^{cst}(\tilde{T}_{\natural}({\mathbb R}))_{\lambda_{\natural}}$.

Notons $\mu(\lambda_{\natural})\in \mathfrak{c}_{\natural}^*$ le param\`etre de $\lambda_{\natural}$ (c'est-\`a-dire que $\lambda_{\natural}(exp(H))=e^{<H,\mu(\lambda_{\natural})>}$ pour $H\in \mathfrak{c}_{\natural}({\mathbb R})$ assez proche de $0$). Notons $I_{\lambda_{\natural}}$ l'id\'eal de $Sym(\mathfrak{t}_{\natural})$ engendr\'e par les $H+<H,\mu(\lambda_{\natural})>$ pour $H\in \mathfrak{c}_{\natural}$. Posons $Sym(\mathfrak{t}_{\natural})_{\lambda_{\natural}}=Sym(\mathfrak{t}_{\natural})/I_{\lambda_{\natural}}$.  Introduisons l'ensemble $\mathfrak{t}_{\natural,\lambda_{\natural}}^*$ des \'el\'ements de $\mathfrak{t}_{\natural}^*$ qui se projettent sur $\mu(\lambda_{\natural})\in \mathfrak{c}_{\natural}^*$. C'est un espace affine sous $\mathfrak{t}$ et $Sym(\mathfrak{t}_{\natural})_{\lambda_{\natural}}$ s'identifie \`a l'alg\`ebre des polyn\^omes sur $\mathfrak{t}_{\natural,\lambda_{\natural}}^*$. L'application $H\mapsto \partial_{H}$ se quotiente en un isomorphisme not\'e de m\^eme de $Sym(\mathfrak{t}_{\natural})_{\lambda_{\natural}}$ sur $Diff^{cst}(\tilde{T}_{\natural}({\mathbb R}))_{\lambda_{\natural}}$.

Rappelons que l'on note $C_{c,\lambda_{\natural}}^{\infty}(\tilde{G}_{\natural}({\mathbb R}))$ l'espace des fonctions $f$ sur $\tilde{G}_{\natural}({\mathbb R})$ qui sont $C^{\infty}$,  \`a support compact modulo $C_{\natural}({\mathbb R})$ et qui v\'erifient la condition  (1) pour tous $\gamma_{\natural}\in \tilde{G}_{\natural}({\mathbb R})$. On note $I_{\lambda_{\natural}}(\tilde{G}_{\natural}({\mathbb R}))$ le quotient de $C_{c,\lambda_{\natural}}^{\infty}(\tilde{G}_{\natural}({\mathbb R}))$ par le sous-espace des fonctions dont les orbitales orbitales sont nulles en tout point fortement r\'egulier. L'action de $\mathfrak{Z}(G_{\natural})$ sur $C_{c,\lambda_{\natural}}^{\infty}(\tilde{G}_{\natural}({\mathbb R}))$ ou $I_{\lambda_{\natural}}(\tilde{G}_{\natural}({\mathbb R}))$ se quotiente en l'action d'une alg\`ebre quotient $\mathfrak{Z}(G_{\natural})_{\lambda_{\natural}}$. Celle-ci est isomorphe \`a  $Sym(\mathfrak{t}_{\natural})_{\lambda_{\natural}}^W$.

Soient $\gamma_{\natural}\in \tilde{T}_{\natural,\tilde{G}-reg}({\mathbb R})$ et $f\in C_{c,\lambda_{\natural}}^{\infty}(\tilde{G}_{\natural}({\mathbb R}))$. On sait d\'efinir $I_{\tilde{M}_{\natural},\lambda_{\natural}}^{\tilde{G}_{\natural}}(\gamma_{\natural},f)$ (cf. [II] 1.10). Par exemple,  fixons une fonction $\dot{f}\in C_{c}^{\infty}(\tilde{G}_{\natural}({\mathbb R}))$ telle que
$$f(\gamma_{\natural})=\int_{C_{\natural}({\mathbb R})}\dot{f}(c\gamma_{\natural})\lambda_{\natural}(c)\,dc.$$
Alors 
$$I_{\tilde{M}_{\natural},\lambda_{\natural}}^{\tilde{G}_{\natural}}(\gamma_{\natural},f)=\int_{C_{\natural}({\mathbb R})}I_{\tilde{M}_{\natural}}^{\tilde{G}_{\natural}}(\gamma_{\natural},\dot{f}^c)\lambda_{\natural}(c)\,dc,$$
o\`u $\dot{f}^c$ est d\'efinie par $\dot{f}^c(\gamma'_{\natural})=\dot{f}(c\gamma'_{\natural})$. En appliquant 1.2, on a pour $z\in \mathfrak{Z}(G_{\natural})$ une \'egalit\'e
$$(2) \qquad I_{\tilde{M}_{\natural},\lambda_{\natural}}^{\tilde{G}_{\natural}}(\gamma_{\natural},zf)=\sum_{\tilde{L}\in {\cal L}(\tilde{M})}\delta_{\tilde{M}_{\natural}}^{\tilde{L}_{\natural}}(\gamma_{\natural},z_{L_{\natural}})I_{\tilde{L}_{\natural},\lambda_{\natural}}^{\tilde{G}_{\natural}}(\gamma_{\natural},f).$$
Puisque les op\'erateurs diff\'erentiels s'appliquent ici \`a des fonctions v\'erifiant (1), on peut les remplacer par leurs images dans l'espace  $Diff^{cst}(\tilde{T}_{\natural}({\mathbb R}))_{\lambda_{\natural}}$. On note $\delta_{\tilde{M}_{\natural},\lambda_{\natural}}^{\tilde{G}_{\natural}}(\gamma_{\natural},z)$ l'image de $\delta_{\tilde{M}_{\natural}}^{\tilde{G}_{\natural}}(\gamma_{\natural},z)$ dans  cet espace. Il r\'esulte formellement de l'\'egalit\'e (2) que

-  $\delta_{\tilde{M}_{\natural},\lambda_{\natural}}^{\tilde{G}_{\natural}}(\gamma_{\natural},z)$ ne d\'epend que de l'image de $z$ dans $\mathfrak{Z}(G_{\natural})_{\lambda_{\natural}}$;  

- $\delta_{\tilde{M}_{\natural},\lambda_{\natural}}^{\tilde{G}_{\natural}}(\gamma_{\natural},z)$ ne d\'epend que de l'image $\gamma\in \tilde{T}({\mathbb R})$ de $\gamma_{\natural}$.

On peut donc noter $\delta_{\tilde{M}_{\natural},\lambda_{\natural}}^{\tilde{G}_{\natural}}(\gamma,z)$ notre op\'erateur. Il est d\'efini pour $\gamma\in \tilde{T}_{\tilde{G}-reg}({\mathbb R})$ et $z\in \mathfrak{Z}(\tilde{G})_{\lambda_{\natural}}$ et prend ses valeurs dans  $Diff^{cst}(\tilde{T}_{\natural}({\mathbb R}))_{\lambda_{\natural}}$. La proposition 1.3 reste vraie pour ces op\'erateurs. En particulier,  la fonction $\gamma\mapsto \delta_{\tilde{M}_{\natural},\lambda_{\natural}}^{\tilde{G}_{\natural}}(\gamma,z)$ est rationnelle et r\'eguli\`ere. Posons
$$Diff^{reg}(\tilde{T}_{\tilde{G}-reg}({\mathbb R})_{\lambda_{\natural}}=Pol(\tilde{T}_{\tilde{G}-reg}({\mathbb C}))\otimes_{{\mathbb C}}Diff^{cst}(\tilde{T}_{\natural}({\mathbb R}))_{\lambda_{\natural}}.$$
On a alors un \'el\'ement $\delta_{\tilde{M}_{\natural},\lambda_{\natural}}^{\tilde{G}_{\natural}}(z)\in 
Diff^{reg}(\tilde{T}_{\tilde{G}-reg}({\mathbb R})_{\lambda_{\natural}}$ dont la valeur en un point $\gamma$ est $\delta_{\tilde{M}_{\natural},\lambda_{\natural}}^{\tilde{G}_{\natural}}(\gamma,z)$.

On a muni $X_{*}(T)\otimes_{{\mathbb Z}}{\mathbb R}$ d'une forme quadratique d\'efinie positive. Fixons un scindage
$$X_{*}(T_{\natural})\otimes_{{\mathbb Z}}{\mathbb R}=X_{*}(C_{\natural})\otimes_{{\mathbb Z}}{\mathbb R}\oplus X_{*}(T)\otimes_{{\mathbb Z}}{\mathbb R}.$$
Munissons le premier facteur d'une forme quadratique d\'efinie positive et munissons l'espace total de la forme somme directe des formes sur les deux facteurs.  On dira que cettte derni\`ere forme est compatible avec celle sur $X_{*}(T)\otimes_{{\mathbb Z}}{\mathbb R}$. On d\'efinit alors l'op\'erateur de Casimir $\Omega^{\tilde{G}_{\natural}}$. Il d\'epend des constructions ci-dessus, mais la classe $\Omega^{\tilde{G}_{\natural}}+Sym(\mathfrak{z}(G_{\natural}))$ n'en d\'epend pas.  Pour tout \'el\'ement $z$ de cette classe, l'op\'erateur $\delta_{\tilde{M}_{\natural},\lambda_{\natural}}^{\tilde{G}_{\natural}}(\gamma,z)$ est donn\'e par la formule de la proposition 1.7 (pourvu que $d\geq1$). 

Consid\'erons d'autres donn\'ees
$$1\to C_{\flat}\to G_{\flat}\to G\to 1\text{ et }\tilde{G}_{\flat}\to \tilde{G}$$
ainsi qu'un caract\`ere $\lambda_{\flat}$ de $C_{\flat}({\mathbb R})$ v\'erifiant les m\^emes conditions que pr\'ec\'edemment. On introduit les produits fibr\'es $G_{\natural,\flat}$ et $\tilde{G}_{\natural,\flat}$ de $G_{\natural}$ et $G_{\flat}$ au-dessus de $G$, resp. de $\tilde{G}_{\natural}$ et $\tilde{G}_{\flat}$ au-dessus de $\tilde{G}$. On suppose donn\'es un caract\`ere $\lambda_{\natural,\flat}$ de $G_{\natural,\flat}({\mathbb R})$ et une fonction $\tilde{\lambda}_{\natural,\flat}$ sur $\tilde{G}_{\natural,\flat}({\mathbb R})$ tels que

- $\lambda_{\natural,\flat}( c_{\natural}x_{\natural},c_{\flat}x_{\flat})=\lambda_{\natural}(c_{\natural})\lambda_{\flat}(c_{\flat})^{-1}\lambda_{\natural,\flat}(x_{\natural},x_{\flat})$ pour tous $(x_{\natural},x_{\flat})\in G_{\natural,\flat}({\mathbb R})$, $c_{\natural}\in C_{\natural}({\mathbb R})$, $c_{\flat}\in C_{\flat}({\mathbb R})$;

- $\tilde{\lambda}_{\natural,\flat}(x_{\natural}\gamma_{\natural} ,x_{\flat}\gamma_{\flat})=\tilde{\lambda}_{\natural,\flat}(\gamma_{\natural},\gamma_{\flat})\lambda_{\natural,\flat}(x_{\natural},x_{\flat})$ pour tous $(\gamma_{\natural},\gamma_{\flat})\in \tilde{G}_{\natural,\flat}({\mathbb R})$ et $(x_{\natural},x_{\flat})\in G_{\natural,\flat}({\mathbb R})$.

On d\'efinit un isomorphisme
$$\begin{array}{ccc}C_{c,\lambda_{\natural}}^{\infty}(\tilde{G}_{\natural}({\mathbb R}))&\to&C_{c,\flat}^{\infty}(\tilde{G}_{\flat}({\mathbb R}))\\ f_{\natural}&\mapsto&f_{\flat}\\ \end{array}$$
par $f_{\flat}(\gamma_{\flat})=\tilde{\lambda}_{\natural,\flat}(\gamma_{\natural},\gamma_{\flat})f_{\natural}(\gamma_{\natural})$ o\`u $\gamma_{\natural}$ est n'importe quel \'el\'ement de $\tilde{G}_{\natural}({\mathbb R})$ tel que $(\gamma_{\natural},\gamma_{\flat})\in \tilde{G}_{\natural,\flat}({\mathbb R})$. On a de m\^eme un isomorphisme entre l'espace des fonctions sur $\tilde{T}_{\natural}({\mathbb R})$ v\'erifiant (1) et l'espace de fonctions analogue sur $\tilde{T}_{\flat}({\mathbb R})$. De cet isomorphisme se d\'eduit un isomorphisme  
$$(3) \qquad Diff^{cst}(\tilde{T}_{\natural}({\mathbb R}))_{\lambda_{\natural}} \simeq Diff^{cst}(\tilde{T}_{\flat}({\mathbb R}))_{\flat}.$$
Au caract\`ere $\lambda_{\natural,\flat}$ est associ\'e un param\`etre $\mu(\lambda_{\natural,\flat})\in \mathfrak{t}_{\natural,\flat}^*$, o\`u $T_{\natural,\flat}$ est l'image r\'eciproque de $T$ dans $G_{\natural,\flat}$. On a
$$\mathfrak{t}_{\natural,\flat}^*=(\mathfrak{t}_{\natural}^*\oplus \mathfrak{t}_{\flat}^*)/diag_{-}(\mathfrak{t}^*).$$
La projection de $\mu(\lambda_{\natural,\flat})$ dans $\mathfrak{c}_{\natural}^*\oplus \mathfrak{c}_{\flat}^*$ est $(\mu(\lambda_{\natural}),-\mu(\lambda_{\flat}))$. Pour $\nu_{\natural}\in \mathfrak{t}_{\natural,\lambda_{\natural}}^*$ il existe un unique $\nu_{\flat}\in \mathfrak{t}_{\flat,\lambda_{\flat}}$ tel que $(-\nu_{\natural},\nu_{\flat})$ ait pour projection $\mu(\lambda_{\natural,\flat})$ dans $\mathfrak{t}_{\natural,\flat}^*$. L'application $\nu_{\natural}\mapsto \nu_{\flat}$
est un isomorphisme de $\mathfrak{t}_{\natural,\lambda_{\natural}}^*$ sur $\mathfrak{t}_{\flat,\lambda_{\flat}}^*$. Il est compatible \`a (3) et au fait que les alg\`ebres d'op\'erateurs diff\'erentiels s'identifient \`a des alg\`ebres de polyn\^omes sur ces deux espaces affines. Il s'en d\'eduit un isomorphisme 
$$(4)\qquad \mathfrak{Z}(G_{\natural})_{\lambda_{\natural}}\simeq \mathfrak{Z}(G_{\flat})_{\lambda_{\flat}}.$$
 Il est formel de v\'erifier que, si $z_{\natural}$ et $z_{\flat}$ se correspondent par (4), alors $\delta_{\tilde{M}_{\natural},\lambda_{\natural}}^{\tilde{G}_{\natural}}(\gamma,z_{\natural})$ et $\delta_{\tilde{M}_{\flat},\lambda_{\flat}}^{\tilde{G}_{\flat}}(\gamma,z_{\flat})$ se correspondent par (3). Remarquons que l'image dans $\mathfrak{Z}(G_{\natural})_{\lambda_{\natural}}$ de la classe $\Omega^{\tilde{G}_{\natural}}+Sym(\mathfrak{z}(G_{\natural}))$ correspond par (4) \`a l'image dans 
 $\mathfrak{Z}(G_{\flat})_{\lambda_{\flat}}$ de la classe $\Omega^{\tilde{G}_{\flat}}+Sym(\mathfrak{z}(G_{\flat}))$.
 
 Revenons au cas o\`u $(G,\tilde{G},{\bf a})$ est quelconque. Consid\'erons des donn\'ees endoscopiques ${\bf G}'=(G',{\cal G}',\tilde{s})$ de $(G,\tilde{G},{\bf a})$ et ${\bf M}'=(M',{\cal M}',\tilde{\zeta})$ de $(M,\tilde{M},{\bf a}_{M})$. Supposons que $M'$ s'identifie \`a un Levi de $G'$ et que ${\bf M}'$ soit la donn\'ee d\'eduite de ${\bf G}'$ via cette identification (cf. [I] 3.4). Supposons aussi que ${\bf M}'$ soit relevante (donc aussi ${\bf G}'$). Fixons des donn\'ees auxiliaires $G'_{1},\tilde{G}'_{1},C_{1},\hat{\xi}_{1}, \Delta_{1}$ pour ${\bf G}'$, qui se restreignent en des donn\'ees auxiliaires $M'_{1},...,\Delta_{1}$ pour ${\bf M}'$. Soient $\gamma_{0}\in \tilde{T}_{\tilde{G}-reg}({\mathbb R})$ et $\delta_{0}\in \tilde{M}'({\mathbb R})$ correspondant \`a $\gamma_{0}$. 
 Notons $\tilde{T}'$ le sous-tore tordu maximal de $\tilde{M}'$ tel que $\delta\in \tilde{T}'({\mathbb R})$. Alors on dispose d'un isomorphisme $\xi:T/(1-\theta)(T)\simeq T'$ dont on d\'eduit un isomorphisme $\tilde{\xi}:\tilde{T}/(1-\theta)(T)\simeq \tilde{T}'$ tel que $\tilde{\xi}(t\gamma_{0})=\xi(t)\delta_{0}$. Modulo ces isomorphismes, les notations $T'$ et $\tilde{T}'$ sont coh\'erentes avec celles du paragraphe 1.1. 
 
 Pour $z'_{1}\in \mathfrak{Z}(G'_{1})_{\lambda_{1}}$ et $\delta\in \tilde{T}'_{\tilde{G}'-reg}({\mathbb R})$, on construit comme plus haut l'op\'erateur $\delta_{\tilde{M}'_{1},\lambda_{1}}^{\tilde{G}'_{1}}(\delta,z'_{1})\in Diff^{cst}(\tilde{T}'_{1}({\mathbb R}))_{\lambda_{1}}$.  Faisons varier les donn\'ees auxiliaires. Les alg\`ebres  $\mathfrak{Z}(G'_{1})_{\lambda_{1}}$ se recollent en une alg\`ebre que l'on a not\'ee $\mathfrak{Z}({\bf G}')$ en [IV] 2.1.  Les espaces affines $\mathfrak{t}_{1,\lambda_{1}}^{'*}$ se recollent en un espace affine isomorphe \`a $\mathfrak{t}_{\theta,\omega}^*$, cf. [IV] 2.1. Donc les  alg\`ebres $ Diff^{cst}(\tilde{T}'_{1}({\mathbb R}))_{\lambda_{1}}$ se recollent en une alg\`ebre isomorphe \`a $Diff(\tilde{T}({\mathbb R}))_{\theta,\omega}$. D'apr\`es les consid\'erations ci-dessus, pour $z'\in \mathfrak{Z}({\bf G}')$, les op\'erateurs $\delta_{\tilde{M}'_{1},\lambda_{1}}^{\tilde{G}'_{1}}(\delta,z'_{1})$  (o\`u $z'_{1}$ correspond \`a $z'$) se recollent en un op\'erateur que l'on peut noter $\delta_{{\bf M}'}^{{\bf G}'}(\delta,z')\in Diff^{cst}(\tilde{T}({\mathbb R}))_{\theta,\omega}$. Les alg\`ebres de fonctions rationnelles et r\'eguli\`eres $Pol(\tilde{T}'_{1,\tilde{G}'-reg}({\mathbb C}))$ se recollent aussi en une sous-alg\`ebre de l'alg\`ebre not\'ee $Pol(\tilde{T}'_{\tilde{G}-reg}({\mathbb C}))$ en 1.1. Ce n'est pas forc\'ement cette alg\`ebre tout enti\`ere car la condition de r\'egularit\'e relative \`a $\tilde{G}'$ est plus faible que celle relative \`a $\tilde{G}$. En tout cas, les termes $\delta_{\tilde{M}'_{1},\lambda_{1}}^{\tilde{G}'_{1}}(z'_{1})$ se recollent en un \'el\'ement $\delta_{{\bf M}'}^{{\bf G}'}(z')\in Diff^{reg}(\tilde{T}_{\tilde{G}-reg}({\mathbb R}))^{\omega-inv}$. 
 
\bigskip

\section{Versions stables et endoscopiques des op\'erateurs diff\'erentiels}

\bigskip

\subsection{Version stable des op\'erateurs diff\'erentiels}
  On  suppose dans ce paragraphe $(G,\tilde{G},{\bf a})$ quasi-d\'eploy\'e et \`a torsion int\'erieure. On fixe comme en 1.1 un espace de Levi $\tilde{M}$ et un sous-tore tordu maximal $\tilde{T}$ de $\tilde{M}$, ainsi que des mesures de Haar sur tous les groupes intervenant.   On sait d\'efinir $S_{\tilde{M}}^{\tilde{G}}(\boldsymbol{\delta},f)$ pour une distribution stable $\boldsymbol{\delta}$ sur $\tilde{M}({\mathbb R})$ \`a support $\tilde{G}$-r\'egulier et pour $f\in C_{c}^{\infty}(\tilde{G}({\mathbb R}))$, cf. [V] 1.4.  A un \'el\'ement $\delta\in \tilde{T}_{\tilde{G}-reg}({\mathbb R})$ est associ\'ee une telle distribution $\boldsymbol{\delta}$, \`a savoir l'int\'egrale orbitale stable sur $\tilde{M}({\mathbb R})$ associ\'ee \`a la classe de conjugaison stable de $\delta$. On pose simplement $S_{\tilde{M}}^{\tilde{G}}(\delta,f)=S_{\tilde{M}}^{\tilde{G}}(\boldsymbol{\delta},f)$. Ce terme devient ainsi une fonction de $\delta$, qui est clairement $C^{\infty}$ sur $\tilde{T}_{\tilde{G}-reg}({\mathbb R})$.

\ass{Proposition}{  Il existe une unique application lin\'eaire
$$\begin{array}{ccc}  \mathfrak{Z}(G)&\to& Diff^{reg}(\tilde{T}_{\tilde{G}-reg}({\mathbb R}))\\
z&\mapsto& S\delta_{\tilde{M}}^{\tilde{G}}(z)\\ \end{array}$$  
qui v\'erifie la propri\'et\'e suivante

 -  pour tout $f\in C_{c}^{\infty}(\tilde{G}({\mathbb R}))$, tout $\delta\in \tilde{T}_{\tilde{G}-reg}({\mathbb R})$ et tout $z\in \mathfrak{Z}(G)$, on a l'\'egalit\'e
$$S_{\tilde{M}}^{\tilde{G}}(\delta,zf)=\sum_{\tilde{L}\in {\cal L}(\tilde{M})} S\delta_{\tilde{M}}^{\tilde{L}}(\delta,z_{L})S_{\tilde{L}}^{\tilde{G}}(\delta,f).$$}

 Preuve. Soit $\delta\in \tilde{T}_{\tilde{G}-reg}({\mathbb R})$. On va commencer par d\'efinir les op\'erateurs $S\delta_{\tilde{M}}^{\tilde{G}}(\delta,z)$ et on montrera ensuite qu'ils v\'erifient les propri\'et\'es requises. Soit $s\in Z(\hat{M})^{\Gamma_{{\mathbb R}}}/Z(\hat{G})^{\Gamma_{{\mathbb R}}}$, $s\not=1$. On fixe des donn\'ees auxiliaires $G'_{1}(s),...,\Delta_{1}(s)$ pour ${\bf G}'(s)$. L'\'el\'ement $z$ d\'etermine un \'el\'ement $z^{{\bf G}'(s)}\in \mathfrak{Z}({\bf G}'(s))$, cf. [III] 2.1.  Fixons un \'el\'ement  $z_{1}(s)\in \mathfrak{Z}(G'_{1}(s))$ d'image $z^{{\bf G}'(s)}$ dans $ \mathfrak{Z}({\bf G}'(s))$. Notons $\tilde{M}_{1}(s)$ et  $\tilde{T}_{1}(s)$ les images r\'eciproques de  $\tilde{M}$ et $\tilde{T}$ dans $\tilde{G}'_{1}(s)$. Soit $\delta\in \tilde{T}_{\tilde{G}-reg}({\mathbb R})$ et fixons un \'el\'ement $\delta_{1}(s)\in \tilde{T}_{1}(s)$ se projetant sur $\delta$.  Puisque $dim(G_{1,SC}'(s))<dim(G_{SC})$, on peut supposer la proposition connue pour $\tilde{G}'_{1}(s)$.  On dispose donc d'un op\'erateur $S\delta_{\tilde{M}_{1}(s)}^{\tilde{G}'_{1}(s)}(\delta_{1}(s),z_{1}(s))$.  Les m\^emes formalit\'es que l'on a d\'evelopp\'ees en 1.8 pour les op\'erateurs $\delta_{\tilde{M}}^{\tilde{G}}(\gamma,z)$ s'appliquent. On est dans le cas particulier simple o\`u la donn\'ee endoscopique de $(M,\tilde{M},{\bf a}_{M})$ est la donn\'ee endoscopique maximale ${\bf M}=(M,{^LM},1)$.  On voit que, quand on fait varier les donn\'ees auxiliaires, les op\'erateurs $S\delta_{\tilde{M}_{1}(s)}^{\tilde{G}'_{1}(s)}(\delta_{1}(s),z_{1}(s))$ se recollent en des op\'erateurs $S\delta_{{\bf M}}^{{\bf G}'(s)}(\delta,z^{{\bf G}'(s)})\in Diff^{cst}(\tilde{T}({\mathbb R}))$. On dispose aussi de l'op\'erateur $\delta_{\tilde{M}}^{\tilde{G}}(\delta,z)$ de 1.2. On d\'efinit alors
 $$(1) \qquad S\delta_{\tilde{M}}^{\tilde{G}}(\delta,z)=\delta_{\tilde{M}}^{\tilde{G}}(\delta,z)-\sum_{s\in Z(\hat{M})^{\Gamma_{{\mathbb R}}}/Z(\hat{G})^{\Gamma_{{\mathbb R}}}, s\not=1}i_{\tilde{M}}(\tilde{G},\tilde{G}'(s))S\delta_{{\bf M}}^{{\bf G}'(s)}(\delta,z^{{\bf G}'(s)}).$$
 
 Par r\'ecurrence, pour $s\not=1$, l'application $\delta\mapsto S\delta_{{\bf M}}^{{\bf G}'(s)}(\delta,z^{{\bf G}'(s)})$ est  la restriction d'une fonction rationnelle et r\'eguli\`ere  sur $\tilde{T}_{\tilde{G}'(s)-reg}({\mathbb C})$, a fortiori sur $\tilde{T}_{\tilde{G}-reg}({\mathbb C})$. Il en est de m\^eme de $\delta\mapsto \delta_{\tilde{M}}^{\tilde{G}}(\delta,z)$ d'apr\`es la proposition 1.3. Donc $\delta\mapsto S\delta_{\tilde{M}}^{\tilde{G}}(\delta,z) $ est la restriction d'une fonction rationnelle et r\'eguli\`ere sur $\tilde{T}_{\tilde{G}-reg}({\mathbb C})$.

 Posons
 $$\Lambda_{\tilde{M}}^{\tilde{G}}(\delta,f)=\sum_{\gamma\in \dot{{\cal X}}(\delta)}I_{\tilde{M}}^{\tilde{G}}(\gamma,f),$$
 o\`u $\dot{{\cal X}}(\delta)$ est un ensemble de repr\'esentants des classes de conjugaison par $M({\mathbb R})$ dans la classe de conjugaison stable de $\delta$. 
  Par d\'efinition
 $$(2) \qquad S_{\tilde{M}}^{\tilde{G}}(\delta,zf)= \Lambda_{\tilde{M}}^{\tilde{G}}(\delta,zf)-\sum_{s\in Z(\hat{M})^{\Gamma_{{\mathbb R}}}/Z(\hat{G})^{\Gamma_{{\mathbb R}}}, s\not=1} i_{\tilde{M}}(\tilde{G},\tilde{G}'(s))S_{{\bf M}}^{{\bf G}'(s)}(\delta,(zf)^{{\bf G}'(s)}).$$
 Cette formule requiert une explication.  On a d\'efini l'espace de distributions stables $D_{g\acute{e}om}^{st}({\bf M})$ par recollement des espaces $D_{g\acute{e}om,\lambda_{1}}^{st}(\tilde{M}_{1}({\mathbb R}))$ pour les diff\'erentes donn\'ees auxiliaires $M_{1}$, $\tilde{M}_{1}$ etc... Mais pour notre donn\'ee ${\bf M}$, on peut prendre pour donn\'ees auxiliaires les donn\'ees triviales $M_{1}=M$, $\tilde{M}_{1}=\tilde{M}$ etc... L'\'el\'ement $\delta$ ayant d\'ej\`a \'et\'e identifi\'e \`a une distribution stable sur $\tilde{M}({\mathbb R})$, il s'identifie aussi \`a un \'el\'ement de $D_{g\acute{e}om}^{st}({\bf M})$. Cela donne un sens \`a la formule (2) ci-dessus.
  
 Etudions $\Lambda_{\tilde{M}}^{\tilde{G}}(\delta,zf)$. Fixons $\delta$. Pour tout $\gamma\in \dot{{\cal X}}(\delta)$, il existe un tore tordu $\tilde{T}_{\gamma}$ et un \'el\'ement $x_{\gamma}\in M$ tel que $ad_{x_{\gamma}}(\tilde{T})=\tilde{T}_{\gamma}$, $ad_{x_{\gamma}}(\delta)=\gamma$ et l'isomorphisme $ad_{x_{\gamma}}:\tilde{T}\to \tilde{T}_{\gamma}$ soit d\'efini sur ${\mathbb R}$. Pour tout $\delta' \in\tilde{T}_{\tilde{G}-reg}({\mathbb R})$, on peut supposer 
 $$\dot{{\cal X}}(\delta')=\{ad_{x_{\gamma}}(\delta'); \gamma\in \dot{{\cal X}}(\delta)\}.$$
 Alors
  $$\Lambda_{\tilde{M}}^{\tilde{G}}(\delta',f)=\sum_{\gamma\in \dot{{\cal X}}(\delta)}I_{\tilde{M}}^{\tilde{G}}(ad_{x_{\gamma}}(\delta'),f).$$
  Rempla\c{c}ons $z$ par $zf$. En appliquant 1.2, on obtient
$$\Lambda_{\tilde{M}}^{\tilde{G}}(\delta,zf)=\sum_{\gamma\in \dot{{\cal X}}(\delta)} \sum_{\tilde{L}\in {\cal L}(\tilde{M})}\delta_{\tilde{M}}^{\tilde{L}}(ad_{x_{\gamma}}(\delta),z_{L})I_{\tilde{L}}^{\tilde{G}}(ad_{x_{\gamma}}(\delta),f),$$
o\`u le terme  $\delta_{\tilde{M}}^{\tilde{L}}(ad_{x_{\gamma}}(\delta),z_{L})I_{\tilde{L}}^{\tilde{G}}(ad_{x_{\gamma}}(\delta),f)$ est  ici la valeur en $\gamma'=ad_{x_{\gamma}}(\delta)$ de la fonction  $\gamma'\mapsto \delta_{\tilde{M}}^{\tilde{L}}(\gamma',z_{L})I_{\tilde{L}}^{\tilde{G}}(\gamma',f)$ d\'efinie sur $\tilde{T}_{\gamma}({\mathbb R})$. La propri\'et\'e 1.6(4) implique qu'il revient au m\^eme d'\'evaluer d'appliquer l'op\'erateur $\delta_{\tilde{M}}^{\tilde{L}}(\delta',z_{L})$ \`a la fonction $\delta'\mapsto I_{\tilde{M}}^{\tilde{G}}(ad_{x_{\gamma}}(\delta'),f)$, puis d\'evaluer la fonction obtenue en $\delta$. On obtient alors l'\'egalit\'e
$$(3) \qquad \Lambda_{\tilde{M}}^{\tilde{G}}(\delta,zf)= \sum_{\tilde{L}\in {\cal L}(\tilde{M})}\delta_{\tilde{M}}^{\tilde{L}}(\delta,z_{L})\Lambda_{\tilde{L}}^{\tilde{G}}(\delta,f).$$
Fixons $s\in Z(\hat{M})^{\Gamma_{{\mathbb R}}}/Z(\hat{G})^{\Gamma_{{\mathbb R}}}$, $s\not=1$. Comme plus haut, on  fixe des donn\'ees auxiliaires $G'_{1}(s),...,\Delta_{1}(s)$ pour ${\bf G}'(s)$ ainsi que des \'el\'ements   $z_{1}(s)\in \mathfrak{Z}(G'_{1}(s))$ d'image $z^{{\bf G}'(s)}$ dans $ \mathfrak{Z}({\bf G}'(s))$ et $\delta_{1}(s)\in \tilde{T}_{1}(s)$ se projetant sur $\delta$. En tant qu'\'el\'ement de $D_{g\acute{e}om}^{st}({\bf M})$, $\delta$ s'identifie \`a un   certain multiple de la distribution stable associ\'ee \`a $\delta_{1}(s)$, disons que c'est $c_{1}(s)$ fois cette distribution.  Identifions aussi $f^{{\bf G}'}$ \`a un \'el\'ement $f^{\tilde{G}'_{1}(s)}\in C_{c,\lambda_{1}(s)}^{\infty}(\tilde{G}'_{1}(s;{\mathbb R}))$. Alors
$$S_{{\bf M}}^{{\bf G}'(s)}(\delta,(zf)^{{\bf G}'(s)})=c_{1}(s)S_{\tilde{M}_{1}(s)}^{\tilde{G}'_{1}(s)}(\delta_{1}(s),z_{1}(s)f^{\tilde{G}'_{1}(s)}).$$
  Puisque $dim(G_{1,SC}'(s))<dim(G_{SC})$, on peut  appliquer la proposition par r\'ecurrence. On obtient l'\'egalit\'e
 $$ S_{{\bf M}}^{{\bf G}'(s)}(\delta,(zf)^{{\bf G}'(s)})=c_{1}(s)\sum_{\tilde{L}'_{1,s}\in {\cal L}^{\tilde{G}'(s)}(\tilde{M})}S\delta_{\tilde{M}_{1}(s)}^{\tilde{L}'_{1,s}}(\delta_{1}(s),z_{1}(s)_{L'_{1,s}})S_{\tilde{L}'_{1,s}}^{\tilde{G}'_{1}(s)}(\delta_{1}(s),f^{\tilde{G}'_{1}(s)}).$$
 Comme on l'a dit plus haut, quand on fait varier les donn\'ees auxiliaires, les op\'erateurs $S\delta_{\tilde{M}_{1}(s)}^{\tilde{G}'_{1}(s)}(\delta_{1}(s),z_{1}(s))$ se recollent en des op\'erateurs $S\delta_{{\bf M}}^{{\bf G}'(s)}(\delta,z^{{\bf G}'(s)})\in Diff^{cst}(\tilde{T}({\mathbb R}))$. L'\'egalit\'e ci-dessus se r\'ecrit
$$(4) \qquad S_{{\bf M}}^{{\bf G}'(s)}(\delta,(zf)^{{\bf G}'(s)})=\sum_{\tilde{ L}'_{s}\in {\cal L}^{\tilde{G}'(s)}(\tilde{M})}S\delta_{{\bf M}}^{{\bf L}'(s)}(\delta,(z^{{\bf G}'(s)})_{{\bf L}'(s})S_{{\bf L}'(s)}^{{\bf G}'(s)}(\delta,f^{{\bf G}'(s)}).$$
On sait que $\tilde{L}'_{s}$ d\'etermine un espace de Levi $\tilde{L}\in {\cal L}(\tilde{M})$ et une donn\'ee endoscopique de $(L,\tilde{L})$. On a not\'e ${\bf L}'(s)$ cette donn\'ee endoscopique.

Maintenant, la preuve reprend celle de la proposition 2.5 de [II]. On ins\`ere la formule (4) dans la somme en $s$ de (2). Chaque couple $(s,\tilde{L}'_{s})$ d\'etermine un espace de Levi $\tilde{L}$ de $\tilde{G}$ et, comme on vient de le dire, une donn\'ee endoscopique   ${\bf L}'(s)$ de $\tilde{L}$, dont la classe ne d\'epend que de l'image de $s$ dans $Z(\hat{M})^{\Gamma_{{\mathbb R}}}/Z(\hat{L})^{\Gamma_{{\mathbb R}}}$. Le terme $(z^{{\bf G}'(s)})_{{\bf L}'_{s}})$ est \'egal \`a $(z_{L})^{{\bf L}'(s)}$. La somme en $s$ de (2) se transforme en
$$\sum_{\tilde{L}\in {\cal L}(\tilde{M})}\sum_{s\in Z(\hat{M})^{\Gamma_{{\mathbb R}}}/Z(\hat{L})^{\Gamma_{{\mathbb R}}},{\bf L}'(s)\text{ elliptique}}\sum_{t\in sZ(\hat{L})^{\Gamma_{{\mathbb R}}}/Z(\hat{G})^{\Gamma_{{\mathbb R}}}, t\not=1}i_{\tilde{M}}(\tilde{G},\tilde{G}'(t))$$
$$S\delta_{{\bf M}}^{{\bf L}'(s)}(\delta,(z_{L})^{{\bf L}'(s)})S_{{\bf L}'(s)}^{{\bf G}'(t)}(\delta,f^{{\bf G}'(t)}).$$
Remarquons que la somme en $\tilde{L}$ et $s$  est en fait limit\'ee \`a $(\tilde{L},s)\not=(\tilde{G},1)$ puisque pour ce terme exceptionnel, la somme en $t$ est vide. L'op\'erateur $S\delta_{{\bf M}}^{{\bf G}}(\delta,z^{{\bf G}})$ correspondant au couple  $(\tilde{G},1)$ n'intervient donc pas (ce qui est heureux puisqu'on ne l'a pas encore d\'efini).  Mais, pour la commodit\'e du calcul qui suit, il convient de poser 
  formellement $S\delta_{{\bf M}}^{{\bf G}}(\delta,z^{{\bf G}})=S\delta_{\tilde{M}}^{\tilde{G}}(\delta,z)$. 
Pour $s,t$ intervenant ci-dessus, on v\'erifie l'\'egalit\'e
$$i_{\tilde{M}}(\tilde{G},\tilde{G}'(t))=i_{\tilde{M}}(\tilde{L},\tilde{L}'(s))i_{\tilde{L}'(s)}(\tilde{G},\tilde{G}'(t)).$$
De plus, la non-nullit\'e du membre de droite ci-dessus implique l'ellipticit\'e de ${\bf L}'(s)$. 
La somme ci-dessus devient
$$(5) \qquad \sum_{\tilde{L}\in {\cal L}(\tilde{M})}\sum_{s\in Z(\hat{M})^{\Gamma_{{\mathbb R}}}/Z(\hat{L})^{\Gamma_{{\mathbb R}}}}i_{\tilde{M}}(\tilde{L},\tilde{L}'(s))X(\tilde{L},s),$$
o\`u $X(\tilde{L},s)$ s'obtient en appliquant l'op\'erateur diff\'erentiel $S\delta_{{\bf M}}^{{\bf L}'(s)}(\delta,(z_{L})^{{\bf L}'(s)})$ \`a la fonction de $\delta$
$$(6) \qquad \sum_{t\in sZ(\hat{L})^{\Gamma_{{\mathbb R}}}/Z(\hat{G}^{\Gamma_{{\mathbb R}}}, t\not=1} i_{\tilde{L}'(s)}(\tilde{G},\tilde{G}'(t))S_{{\bf L}'(s)}^{{\bf G}'(t)}(\delta,f^{{\bf G}'(t)}).$$
 Supposons d'abord $\tilde{L}\not=\tilde{M}$ et $s\not=1$. 
Dans la formule (6), la restriction $t\not=1$ est superflue et  (6) n'est autre que $I_{\tilde{L}}^{\tilde{G},{\cal E}}({\bf L}'(s),\delta,f)$. Le transfert de $\delta$ vu comme une distribution stable n'est autre que la somme des int\'egrales orbitales sur les \'el\'ements de $\dot{{\cal X}}(\delta)$.  On a prouv\'e en [V] proposition 1.13 l'\'egalit\'e $I_{\tilde{L}}^{\tilde{G},{\cal E}}({\bf L}'(s),\delta,f)=\Lambda_{\tilde{L}}^{\tilde{G}}(\delta,f)$. Supposons maintenant  $\tilde{L}\not=\tilde{M}$ et $s=1$. Alors (6) est \'egal \`a $I_{\tilde{L}}^{\tilde{G},{\cal E}}({\bf L},\delta,f)-S_{\tilde{L}}(\delta,f)$. Ou encore, par le m\^eme argument, \`a $\Lambda_{\tilde{L}}^{\tilde{G}}(\delta,f)-S_{\tilde{L}}^{\tilde{G}}(\delta,f)$. Si enfin $\tilde{L}=\tilde{M}$, la somme (6) est encore \'egale \`a $\Lambda_{\tilde{M}}^{\tilde{G}}(\delta,f)-S_{\tilde{M}}^{\tilde{G}}(\delta,f)$: c'est la formule (1) pour $z=1$. Remarquons que la contribution des termes pour lesquels $s=1$ se simplifie: on a simplement $i_{\tilde{M}}(\tilde{L},\tilde{L}'(1))=i_{\tilde{M}}(\tilde{L},\tilde{L})=1$ et $S\delta_{{\bf M}}^{{\bf L}'(1)}(\delta,(z_{L})^{{\bf L}'(1)})=S\delta_{\tilde{M}}^{\tilde{L}}(\delta,z_{L})$. On transforme (5) conform\'ement \`a ces calculs.  On obtient
$$\sum_{\tilde{L}\in {\cal L}(\tilde{M})}\sum_{s\in Z(\hat{M})^{\Gamma_{{\mathbb R}}}/Z(\hat{L})^{\Gamma_{{\mathbb R}}}}i_{\tilde{M}}(\tilde{L},\tilde{L}'(s))S\delta_{{\bf M}}^{{\bf L}'(s)}(\delta,(z_{L})^{{\bf L}'(s)})\Lambda_{\tilde{L}}^{\tilde{G}}(\delta,f)$$
 $$-\sum_{\tilde{L}\in {\cal L}(\tilde{M})}S\delta_{\tilde{M}}^{\tilde{L}}(\delta,z_{L})S_{\tilde{L}}^{\tilde{G}}(\delta,f).$$
 En utilisant la formule (1) en y rempla\c{c}ant $\tilde{G}$ par $\tilde{L}$,  la premi\`ere somme ci-dessus devient simplement
 $$\sum_{\tilde{L}\in {\cal L}(\tilde{M})}\delta_{\tilde{M}}^{\tilde{L}}(\delta,z_{L})\Lambda_{\tilde{L}}^{\tilde{G}}(\delta,f),$$
 ce qui est le membre de droite de (3).  En appliquant (3), on obtient que (5) est \'egal \`a
 $$\Lambda_{\tilde{M}}^{\tilde{G}}(\delta,zf)-\sum_{\tilde{L}\in {\cal L}(\tilde{M})}S\delta_{\tilde{M}}^{\tilde{L}}(\delta,z_{L})S_{\tilde{L}}^{\tilde{G}}(\delta,f).$$
 On se rappelle que (5) est la somme en $s$ intervenant dans (2) et que, dans cette formule, elle est affect\'ee du signe $-$.  En appliquant la formule ci-dessus, la formule (1) devient celle de l'\'enonc\'e. $\square$ 
 
 \bigskip
 
 \subsection{Propri\'et\'es des versions stables des op\'erateurs diff\'erentiels}
 On conserve la situation du paragraphe pr\'ec\'edent.  
 La formule de d\'efinition (1) de ce paragraphe  et une r\'ecurrence imm\'ediate montrent que les op\'erateurs $S\delta_{\tilde{M}}^{\tilde{G}}(\delta,z)$ v\'erifient une proposition analogue \`a 1.3.

    On a 
  
  (1) si $\tilde{M}=\tilde{G}$, alors $S\delta_{\tilde{M}}^{\tilde{G}}(\delta,z)=\partial_{z_{T}}$.
  
  En effet, la d\'efinition 2.1(1) donne simplement $S\delta_{\tilde{M}}^{\tilde{G}}(\delta,z)=\delta_{\tilde{M}}^{\tilde{G}}(\delta,z)$, d'o\`u le r\'esultat d'apr\`es 1.2(2).

 Comme en 1.6, introduisons l'op\'erateur de Casimir $\Omega$.  Posons $d=a_{\tilde{M}}-a_{\tilde{G}}=a_{M}-a_{G}$ et supposons $d\geq1$.  Supposons d'abord $d=1$.  Remarquons  que, dans notre situation \`a torsion int\'erieure, la d\'efinition de la fonction $C_{\tilde{M}}^{\tilde{G}}$ de 1.7 se simplifie. On a simplement
$$(2) \qquad C_{\tilde{M}}^{\tilde{G}}(\delta)=-\sum_{\alpha\in \Sigma^{G}(T)-\Sigma^{M}(T)}\vert\alpha_{\vert {\cal A}_{M}}\vert (1- \alpha(\delta))^{-1}(1- \alpha(\delta)^{-1})^{-1}.$$
On a not\'e simplement $\alpha(\delta)$ le terme $\alpha(t_{\delta})$, o\`u $t_{\delta}$ est d\'efini en   1.4. Une d\'efinition \'equivalente dans notre situation \`a torsion int\'erieure est d'envoyer $\delta$ en un \'el\'ement $\delta_{ad}\in T_{ad}\subset G_{AD}$ et de poser $\alpha(\delta)=\alpha(\delta_{ad})$.
 Introduisons comme toujours une paire de Borel \'epingl\'ee de $\hat{G}$ pour laquelle on r\'ealise $\hat{M}$ comme un Levi standard. On note $\hat{T}$ le tore de cette paire. Il y a une bijection $\alpha\mapsto \hat{\alpha}$ de $\Sigma^{G}(T)$ sur $\Sigma^{\hat{G}}(\hat{T})$, qui se restreint en une bijection de $\Sigma^{G}(T)-\Sigma^{M}(T)$ sur $\Sigma^{\hat{G}}(\hat{T})-\Sigma^{\hat{M}}(\hat{T})$. Puisque $d=1$, les restrictions \`a $Z(\hat{M})^{\Gamma_{{\mathbb R}}}$ d'\'el\'ements de $\Sigma^{\hat{G}}(\hat{T})-\Sigma^{\hat{M}}(\hat{T})$ sont toutes proportionnelles. Pr\'ecis\'ement, pour $\beta\in \Sigma^{\hat{G}}(\hat{T})-\Sigma^{\hat{M}}(\hat{T})$, notons $\beta_{*}$ sa restriction \`a  $Z(\hat{M})^{\Gamma_{{\mathbb R}}}$. On peut fixer $\beta_{1}\in \Sigma^{\hat{G}}(\hat{T})-\Sigma^{\hat{M}}(\hat{T})$ et un entier $N\geq1$ tels que l'ensemble des $\beta_{*}$ pour 
$\beta\in \Sigma^{\hat{G}}(\hat{T})-\Sigma^{\hat{M}}(\hat{T})$ soit exactement l'ensemble $\{\pm k\beta_{1,*}; k=1,...,N\}$. Pour tout $\beta\in \Sigma^{\hat{G}}(\hat{T})-\Sigma^{\hat{M}}(\hat{T})$, on note $k^{\hat{G}}(\beta)\in {\mathbb N}$ l'entier tel que  $\beta_{*}=\pm k^{\hat{G}}(\beta)\beta_{1,*}$. On d\'efinit une fonction $SC_{\tilde{M}}^{\tilde{G}}$ sur $\tilde{T}_{\tilde{G}-reg}$ par
$$(3) \qquad SC_{\tilde{M}}^{\tilde{G}}(\delta)=-\sum_{\alpha\in \Sigma^{G}(T)-\Sigma^{M}(T)}k^{\hat{G}}(\hat{\alpha})^{-1}\vert\alpha_{\vert {\cal A}_{M}}\vert (1- \alpha(\delta))^{-1}(1- \alpha(\delta)^{-1})^{-1}.$$

Si maintenant $d\geq2$, on pose $SC_{\tilde{M}}^{\tilde{G}}=0$.

\ass{Proposition}{On suppose $d\geq1$. Alors, pour tout $\delta\in \tilde{T}_{\tilde{G}-reg}({\mathbb R})$, l'op\'erateur $S\delta_{\tilde{M}}^{\tilde{G}}(\delta,\Omega)$ est la multiplication par $SC_{\tilde{M}}^{\tilde{G}}(\delta)$.}

Preuve.  On va plut\^ot prouver que cela est vrai quand on remplace $\Omega$ par n'importe quel \'el\'ement $\underline{\Omega}\in \Omega+Sym(\mathfrak{z}(G))$.  
Utilisons la d\'efinition 2.1(1) pour $z=\Omega$. Le premier terme $\delta_{\tilde{M}}^{\tilde{G}}(\delta,\underline{\Omega})$  est la multiplication par $C_{\tilde{M}}^{\tilde{G}}(\delta)$, cf. la remarque de 1.7. Soit $s\in Z(\hat{M})^{\Gamma_{{\mathbb R}}}/Z(\hat{G})^{\Gamma_{{\mathbb R}}}$ avec $s\not=1$. On introduit des donn\'ees auxiliaires $G'_{1}(s),...,\Delta_{1}(s)$. En notant $\tilde{T}_{1}(s)$ l'image r\'eciproque de $\tilde{T}\subset \tilde{M}\subset \tilde{G}'(s)$ dans $\tilde{G}'_{1}(s)$, on fixe une forme quadratique d\'efinie positive sur $X_{*}(T_{1}(s))\otimes_{{\mathbb Z}}{\mathbb R}$ compatible avec celle d\'ej\`a fix\'ee sur $X_{*}(T)\otimes_{{\mathbb Z}}{\mathbb R}$, au sens expliqu\'e en 1.8. On introduit l'op\'erateur de Casimir $\Omega^{G'_{1}(s)}$ relatif \`a cette forme. Il n'est pas vrai en g\'en\'eral que l'image de $\Omega^{G'_{1}(s)}$ dans $\mathfrak{Z}({\bf G}'(s))$ soit \'egale \`a celle de $\Omega=\Omega^{G}$. Mais on v\'erifie facilement que l'image de l'espace affine $\Omega^{G'_{1}(s)}+Sym(\mathfrak{z}(G'_{1}(s)))$ contient celle de l'espace affine $\Omega^G+Sym(\mathfrak{z}(G))$. On peut donc choisir $\underline{\Omega}^{G'_{1}(s)}\in \Omega^{G'_{1}(s)}+Sym(\mathfrak{z}(G'_{1}(s)))$ dont l'image dans $\mathfrak{Z}({\bf G}'(s))$ soit \'egale \`a celle de $\underline{\Omega}$. Pour $\delta_{1}\in \tilde{T}_{1}(s;{\mathbb R})$ au-dessus de $\delta$, l'op\'erateur $\delta_{\tilde{M}_{1}(s)}^{\tilde{G}'_{1}(s)}(\delta_{1},\underline{\Omega}^{G'_{1}(s)})$ est la multiplication par $SC_{\tilde{M}_{1}(s)}^{\tilde{G}'_{1}(s)}(\delta_{1})$. Ce terme co\"{\i}ncide avec $SC_{\tilde{M}}^{\tilde{G}'(s)}(\delta)$. On obtient que
$S\delta_{{\bf M}}^{{\bf G}'(s)}(\delta,\underline{\Omega}^{{\bf G}'(s)})$ est la multiplication par $SC_{\tilde{M}}^{\tilde{G}'(s)}(\delta)$. En revenant \`a la d\'efinition 2.1(1), on voit que, pour prouver la proposition, il suffit de prouver l'\'egalit\'e
$$SC_{\tilde{M}}^{\tilde{G}}(\delta)=C_{\tilde{M}}^{\tilde{G}}(\delta)-\sum_{s\in Z(\hat{M})^{\Gamma_{{\mathbb R}}}/Z(\hat{G})^{\Gamma_{{\mathbb R}}}, s\not=1}i_{\tilde{M}}(\tilde{G},\tilde{G}'(s))SC_{\tilde{M}}^{\tilde{G}'(s)}(\delta),$$
ou encore
$$(4) \qquad C_{\tilde{M}}^{\tilde{G}}(\delta)=\sum_{s\in Z(\hat{M})^{\Gamma_{{\mathbb R}}}/Z(\hat{G})^{\Gamma_{{\mathbb R}}}}i_{\tilde{M}}(\tilde{G},\tilde{G}'(s))SC_{\tilde{M}}^{\tilde{G}'(s)}(\delta).$$
Cela est clair si $d\geq2$: tout est nul. Supposons $d=1$. Pour tout $s$, l'ensemble $\Sigma^{G'(s)}$ est un sous-ensemble de $\Sigma^G$. Tous les termes sont donc des sommes sur l'ensemble $\Sigma^G(T)-\Sigma^M(T)$. On peut fixer un \'el\'ement $\alpha$ de cet ensemble et montrer que le terme index\'e par $\alpha$ dans le membre de gauche de (4) est \'egal \`a la somme des termes index\'es par $\alpha$ dans le membre de droite. Dans chacun de ces termes appara\^{\i}t un facteur commun $-\vert\alpha_{\vert {\cal A}_{\tilde{M}}}\vert (1- \alpha(\delta))^{-1}(1- \alpha(\delta)^{-1})^{-1}$. On peut se limiter aux autres termes. D'apr\`es les formules (2) et (3), l'\'egalit\'e restant \`a prouver se r\'eduit \`a
 $$(5)\qquad 1=\sum_{s\in Z(\hat{M})^{\Gamma_{{\mathbb R}}}/Z(\hat{G})^{\Gamma_{{\mathbb R}}},\alpha\in \Sigma^{G'(s)}(T)}k^{\hat{G}'(s)}(\hat{\alpha})^{-1}i_{\tilde{M}}(\tilde{G},\tilde{G}'(s)).$$
Fixons $\beta_{1}\in \Sigma^{\hat{G}}-\Sigma^{\hat{M}}$ dont la restriction  $\beta_{1,*}$ \`a $Z(\hat{M})^{\Gamma_{{\mathbb {R}}}}$ soit indivisible. Quitte \`a changer $\beta_{1}$ en $-\beta_{1}$, on peut supposer que  $\hat{\alpha}_{*}=k^{\hat{G}}(\hat{\alpha})\beta_{1,*}$. L'application $s\mapsto  \beta_{1}(s)$ identifie  $Z(\hat{M})^{\Gamma_{{\mathbb R}}}/Z(\hat{G})^{\Gamma_{{\mathbb R}}}$ \`a ${\mathbb C}^{\times}$.  Pour $\alpha'\in \Sigma^{G}(T)-\Sigma^M(T)$, on a $\alpha'\in \Sigma^{G'(s)}(T)$ si et seulement si $\beta_{1}(s)^{k^{\hat{G}}(\alpha')}=1$. 
Il en r\'esulte d'abord que la somme sur $s$ de l'expression (5) s'identifie \`a une somme sur $\zeta\in \boldsymbol{\zeta}_{k}$, o\`u $k=k^{\hat{G}}(\alpha)$ et $\boldsymbol{\zeta}_{k}$ est le groupe des racines $k$-i\`emes de l'unit\'e dans ${\mathbb C}^{\times}$. Ensuite, soit $\zeta\in \boldsymbol{\zeta}_{k}$ et soit $s$  tel que $\beta_{1}(s)=\zeta$. Notons $l(\zeta)$ l'ordre de $\zeta$. Il divise $k$. Le calcul ci-dessus montre que $\Sigma^{G(s)}(T)-\Sigma^M(T)$ est form\'e des $\alpha'\in \Sigma^{G}(T)-\Sigma^M(T)$ tels que $\hat{\alpha}'_{*}\in l(\zeta){\mathbb Z}\beta_{1,*}$. Cela entra\^{\i}ne que $k^{\hat{G}'(s)}(\alpha)=k/l(\zeta)$ et que $(Z(\hat{M})^{\Gamma_{{\mathbb R}}}\cap Z(\hat{G}'(s)))/Z(\hat{G})^{\Gamma_{{\mathbb R}}}$ s'identifie au groupe des $\zeta'\in {\mathbb C}^{\times}$ tels que $(\zeta')^{l(\zeta)}=1$. Ce groupe a $l(\zeta)$ \'elements. Par construction de $G'(s)$, les actions galoisiennes sur $Z(\hat{G}'(s))$  et sur $Z(\hat{M})$ s'identifient sur $Z(\hat{G}'(s))\cap Z(\hat{M})$. Le groupe pr\'ec\'edent est donc \'egal \`a $Z(\hat{G}'(s))^{\Gamma_{{\mathbb R}}}/Z(\hat{G})^{\Gamma_{{\mathbb R}}}$. En se rappelant la d\'efinition de $i_{\tilde{M}}(\tilde{G},\tilde{G}'(s))$ ([II] 1.10), on obtient $i_{\tilde{M}}(\tilde{G},\tilde{G}'(s))=l(\zeta)^{-1}$.  L'\'egalit\'e (5) \`a prouver se r\'ecrit
$$1=\sum_{\zeta\in \boldsymbol{\zeta}_{k}}\frac{l(\zeta)}{k}l(\zeta)^{-1}.$$
Cette \'egalit\'e est \'evidente. $\square$

\bigskip

\subsection{Variante endoscopique des op\'erateurs diff\'erentiels }
 Au lieu d'un triplet $(G,\tilde{G},{\bf a})$, on consid\`ere maintenant un triplet $(KG,K\tilde{G},{\bf a})$ comme en [I] 1.11. On fixe un $K$-espace de Levi minimal $K\tilde{M}_{0}$.  Soit $K\tilde{M}\in {\cal L}(K\tilde{M}_{0})$, cf. [I] 3.5. Au lieu d'un tore tordu $\tilde{T}$, on fixe des familles finies $(\tilde{T}_{i})_{i\in I}$ et $(\tilde{\iota}_{i,j})_{i,j\in I}$ v\'erifiant les conditions suivantes (on renvoie \`a [I] 1.11 et [I] 3.5 pour les notations).  On peut supposer que $\tilde{\phi}_{p,q}(\tilde{M}_{q})=\tilde{M}_{p}$ pour tous $p,q\in \Pi^M$. Pour tout $i\in I$, il existe $p(i)\in \Pi^M$ tel que $\tilde{T}_{i}$ soit un sous-tore tordu maximal de $\tilde{M}_{p(i)}$.   Pour deux \'el\'ements $i,j\in I$, $\tilde{\iota}_{i,j}$ est un isomorphisme d\'efini sur ${\mathbb R}$ de $\tilde{T}_{j}$ sur $\tilde{T}_{i}$. On suppose qu'il existe $x_{i,j}\in M_{p(i)}$ tel que $\tilde{\iota}_{i,j}$ soit la restriction de $ad_{x_{i,j}}\circ \tilde{\phi}_{p(i),p(j)}$. On note $\iota_{i,j}:T_{j}\to T_{i}$ la restriction de $ad_{x_{i,j}}\circ \phi_{p(i),p(j)}$, qui est un isomorphisme d\'efini sur ${\mathbb R}$. On suppose que $\tilde{\iota}_{i,j} \circ\tilde{\iota}_{j,k}=\tilde{\iota}_{i,k}$ pour tous $i,j,k\in I$. On suppose enfin que, pour tout $i\in I$ et tout $\gamma\in \tilde{T}_{i}({\mathbb R})$ qui est fortement r\'egulier dans $\tilde{M}_{p(i)}$, l'ensemble $\{\tilde{\iota}_{j,i}(\gamma);j\in I\}$ est un ensemble de repr\'esentants des classes de conjugaison dans la classe de conjugaison stable de $\gamma$ dans $K\tilde{M}({\mathbb R})$.
 
 {\bf Remarque.} Pour tout $p\in \Pi^M$ et tout sous-tore tordu maximal $\tilde{T}$ de $\tilde{M}_{p}$, on peut compl\'eter $\tilde{T}$ en une famille $(\tilde{T}_{i})_{i\in I}$ et d\'efinir une famille d'isomorphismes $(\tilde{\iota}_{i,j})_{i,j\in I}$ de sorte que les conditions ci-dessus soient v\'erifi\'ees.
 
 \bigskip
 
 On note $\tilde{T}$ l'ensemble des familles $\gamma_{I}=(\gamma_{i})_{i\in I}$ telles que, pour tout $i\in I$, $\gamma_{i}\in \tilde{T}_{i}$ et, pour tous $i,j\in I$, on ait $\gamma_{i}=\tilde{\iota}_{i,j}(\gamma_{j})$. On d\'efinit le sous-ensemble \'evident $\tilde{T}_{\tilde{G}-reg}$. Pour $i\in I$, on a introduit en 1.1 divers espaces de fonctions ou op\'erateurs diff\'erentiels $C_{c}^{\infty}(\tilde{T}_{i,\tilde{G}_{i}-reg}({\mathbb R}))^{\omega-inv}$, $Diff^{cst}(\tilde{T}_{i}({\mathbb R}))^{\omega-inv}$ etc... Pour $i,j\in I$, les isomorphismes $\tilde{\iota}_{i,j}$ identifient les espaces index\'es par $i$ \`a ceux index\'es par $j$. On note simplement ces espaces $C_{c}^{\infty}(\tilde{T}_{\tilde{G}-reg}({\mathbb R}))^{\omega-inv}$, $Diff^{cst}(\tilde{T}({\mathbb R}))^{\omega-inv}$ etc... Pour chacun de ces espaces, disons pour l'espace $E$, on note $Mat_{I}(E)$ l'espace des matrices carr\'ees $I\times I$ \`a coefficients dans $E$.

 De nouveau, on fixe des mesures sur tous les groupes intervenant. On suppose que les mesures sur les $G_{p}({\mathbb R})$ sont coh\'erentes en ce sens que, pour $p,q\in \Pi$, les mesures sur $G_{p}({\mathbb R})$ et $G_{q}({\mathbb R})$ se correspondent via le torseur int\'erieur $ \phi_{p,q}$. De m\^eme pour les mesures sur les $M_{p}({\mathbb R})$. On suppose aussi que, pour $i,j\in I$, l'isomorphisme $\iota_{i,j}$ transporte la mesure sur $T_{j}({\mathbb R})$ en celle sur $T_{i}({\mathbb R})$.
 
 Soient $\gamma_{I}=(\gamma_{i})_{i\in I}\in \tilde{T}_{\tilde{G}-reg}({\mathbb R})$ et $f\in C_{c}^{\infty}(K\tilde{G}({\mathbb R}))$. Pour $i\in I$, on d\'efinit l'int\'egrale orbitale pond\'er\'ee $\omega$-\'equivariante endoscopique $I_{K\tilde{M}}^{K\tilde{G},{\cal E}}(\gamma_{i},\omega,f)$ comme en [V] 1.8.

 A partir de maintenant, on va utiliser les hypoth\`eses de r\'ecurrence telles qu'on les a pos\'ees en [V] 1.1.

 \ass{Proposition}{  Il existe une unique application lin\'eaire 
 $$\begin{array}{ccc}  \mathfrak{Z}(G)&\to&Mat_{I}(Diff^{reg}(\tilde{T}_{\tilde{G}-reg}({\mathbb R}))^{\omega-inv})\\ z&\mapsto& \delta_{K\tilde{M}}^{K\tilde{G},{\cal E}}(z)=(\delta_{K\tilde{M};i,j}^{K\tilde{G},{\cal E}}(z))_{i,j\in I}\\ \end{array}$$
 qui v\'erifie la propri\'et\'es suivante:

 - pour tout $f\in C_{c}^{\infty}(\tilde{G}({\mathbb R}))$, tout $ \gamma_{I}\in \tilde{T}_{\tilde{G}-reg}({\mathbb R})$, tout $z\in \mathfrak{Z}(G)$ et tout $i\in I$, on a l'\'egalit\'e
 $$I_{K\tilde{M}}^{K\tilde{G},{\cal E}}(\gamma_{i},\omega,zf)=\sum_{K\tilde{L}\in {\cal L}(K\tilde{M})}\sum_{j\in I}\delta_{K\tilde{M};i,j}^{K\tilde{L},{\cal E}}(\gamma_{j},z_{L})I_{K\tilde{L}}^{K\tilde{G},{\cal E}}(\gamma_{j},\omega,f).$$}
 
 Preuve.   On peut s'autoriser \`a conjuguer chaque $\tilde{T}_{i}$ par un \'el\'ement de $M_{p(i)}({\mathbb R})$. On peut donc supposer qu'il existe un $K$-espace de Levi $K\tilde{R}\in {\cal L}(K\tilde{M}_{0})$ contenu dans  $K\tilde{M}$ tel que, pour tout $i$, $\tilde{T}_{i}$ soit un sous-tore tordu maximal elliptique de $\tilde{R}_{p(i)}$. Pour simplifier, on note simplement  $\tilde{R}_{i}=\tilde{R}_{p(i)}$, $\tilde{M}_{i}=\tilde{M}_{p(i)}$, $\tilde{G}_{i}=\tilde{G}_{p(i)}$. En appliquant la description de [I] 4.9, on voit que l'on peut fixer des familles $({\bf R}'_{i})_{i\in I}$ et $(\tilde{T}'_{i})_{i\in I}$ v\'erifiant les conditions qui suivent. Pour tout $i\in I$, ${\bf R}'_{i}$ est une donn\'ee endoscopique elliptique et relevante de $(KR,K\tilde{R},{\bf a}_{R})$ et $\tilde{T}'_{i}$ est un sous-tore tordu maximal de $\tilde{R}'_{i}$. Pour tous $i,j\in I$, il y a un homomorphisme $\xi_{i,j}:T_{j}\to T'_{i}$ et une application compatible $\tilde{\xi}_{i,j}:\tilde{T}_{j}\to \tilde{T}'_{i}$, tous deux d\'efinis sur ${\mathbb R}$, qui se quotientent en des isomorphismes $T_{j}/(1-\theta)(T_{j})\simeq T'_{i}$ et $\tilde{T}_{j}/(1-\theta)(T_{j})\simeq \tilde{T}'_{i}$ (on note uniform\'ement $\theta$ les automorphismes bien qu'ils d\'ependent \'evidemment du tore en question). On a $\tilde{\xi}_{i,j}\circ \tilde{\iota}_{j,k}=\tilde{\xi}_{i,k}$ pour tous $i,j,k\in I$. Pour un \'el\'ement $\gamma_{j}\in \tilde{T}_{j}({\mathbb R})$ fortement r\'egulier dans $\tilde{R}_{i}$, la classe de conjugaison stable de $\tilde{\xi}_{i,j}(\gamma_{j})$ dans $\tilde{R}'_{i}({\mathbb R})$ correspond \`a celle $\gamma_{j}$ dans $\tilde{R}_{j}({\mathbb R})$. Notons $\tilde{T}'$ l'ensemble des $\delta_{I}=(\delta_{i})_{i\in I}$ tels que $\delta_{i}\in \tilde{T}'_{i}$ pour tout $i$ et il existe $\gamma_{I}\in \tilde{T}$ de sorte que $\delta_{i}=\tilde{\xi}_{i,j}(\gamma_{j})$ pour tous $i,j\in I$. Des applications $\tilde{\xi}_{i,j}$ se d\'eduit une application surjective $\tilde{\xi}:\tilde{T}\to \tilde{T}'$. Pour tout $i\in I$, fixons des donn\'ees auxiliaires $R'_{i,1},...,\Delta_{i,1}$ pour ${\bf R}'_{i}$. Notons $\tilde{T}'_{1}$ l'ensemble des $\delta_{I,1}=(\delta_{i,1})_{i\in I}$ o\`u $\delta_{i,1}\in \tilde{T}'_{i,1}$ pour tout $i$, tels que, en notant $\delta_{i}$ l'image de $\delta_{i,1}$ dans $\tilde{T}'_{i}$, la famille $\delta_{I}=(\delta_{i})_{i\in I}$ appartienne \`a $\tilde{T}'$. Soient $\varphi=(\varphi_{i})_{i\in I}\in C_{c}^{\infty}(K\tilde{R}({\mathbb R}))$, $\gamma_{I}=(\gamma_{i})_{i\in I}\in \tilde{T}_{\tilde{G}-reg}({\mathbb R})$ et $\delta_{I,1}=(\delta_{i,1})_{i\in I}\in \tilde{T}'_{1}({\mathbb R})$ un \'el\'ement au-dessus de $\tilde{\xi}(\gamma_{I})$. On a pour tout $i\in I$ les formules d'inversion
 $$(1) \qquad S^{\tilde{R}'_{i,1}}(\delta_{i,1},\varphi^{\tilde{R}'_{i,1}})=c\sum_{j\in I}\Delta_{i,1}(\delta_{i,1},\gamma_{j})I^{\tilde{R}_{j}}(\gamma_{j},\omega,\varphi_{j}),$$
 $$(2) \qquad I^{\tilde{R}_{i}}(\gamma_{i},\omega,\varphi_{i})=c^{-1}\vert I\vert ^{-1}\sum_{j\in I}\Delta_{j,1}(\delta_{j,1},\gamma_{i})^{-1}S^{\tilde{R}'_{j,1}}(\delta_{j,1},\varphi^{\tilde{R}'_{j,1}}),$$
 o\`u $c$ est une constante non nulle. Pour tout $i$, r\'ealisons la donn\'ee endoscopique ${\bf R}'_{i}$ de $K\tilde{R}$ comme une "donn\'ee de Levi"  d'une donn\'ee endoscopique elliptique ${\bf M}'_{i}=(M'_{i},{\cal M}'_{i},\tilde{\zeta}_{i})$ de $K\tilde{M}$ et  supposons que les donn\'ees auxiliaires ci-dessus sont les restrictions de donn\'ees auxiliaires pour les ${\bf M}'_{i}$. Alors, par simple induction,  on peut remplacer les $R$ par des $M$ dans les formules d'inversion ci-dessus. D\'esormais, on oublie le $K$-espace $K\tilde{R}$ qui ne nous a servi qu'\`a appliquer les consid\'erations de [I]  4.9. 
 
 Soient $\gamma_{I}=(\gamma_{i})_{i\in I}\in \tilde{T}_{\tilde{G}-reg}({\mathbb R})$ et $\delta_{I,1}=(\delta_{i,1})_{i\in I}\in \tilde{T}'_{1}({\mathbb R})$ un \'el\'ement au-dessus de $\tilde{\xi}(\gamma_{I})$.  Pour tout $i\in I$, on peut identifier comme en 2.1 l'\'element $\delta_{i,1}$ \`a  un \'el\'ement de $D_{g\acute{e}om,\lambda_{i,1}}^{st}(\tilde{M}'_{i,1}({\mathbb R}))$, \`a savoir l'int\'egrale orbitale stable associ\'ee \`a $\delta_{i,1}$. On l'identifie ensuite \`a un \'el\'ement de $D_{g\acute{e}om}^{st}({\bf M}'_{i})$. La formule (2) ci-dessus (avec $R$ remplac\'e par $M$) signifie que l'int\'egrale orbitale dans $\tilde{M}_{i}$ associ\'ee \`a $\gamma_{i}$ est  \'egale \`a $c^{-1}\vert I\vert ^{-1}\sum_{j\in I}\Delta_{j,1}(\delta_{j,1},\gamma_{i})^{-1}transfert(\delta_{j,1})$. Par d\'efinition, on a alors
 $$I_{K\tilde{M}}^{K\tilde{G},{\cal E}}(\gamma_{i},\omega,f)=c^{-1}\vert I\vert ^{-1}\sum_{j\in I}\Delta_{j,1}(\delta_{j,1},\gamma_{i})^{-1}I_{K\tilde{M}}^{K\tilde{G},{\cal E}}({\bf M}'_{j},\delta_{j,1},f^{{\bf M}'_{j}}).$$
  $$=c^{-1}\vert I\vert ^{-1}\sum_{j\in I}\Delta_{j,1}(\delta_{j,1},\gamma_{i})^{-1}\sum_{\tilde{s}_{j}\in \tilde{\zeta}_{j}Z(\hat{M})^{\Gamma_{{\mathbb R}},\hat{\theta}}/Z(\hat{G})^{\Gamma_{{\mathbb R}},\hat{\theta}}}i_{\tilde{M}_{j}'}(\tilde{G},\tilde{G}'_{j}(\tilde{s}_{j}))S_{{\bf M}'_{j}}^{{\bf G}'_{j}(\tilde{s}_{j})}(\delta_{j,1},f^{{\bf G}'_{j}(\tilde{s}_{j})}).$$
  Rempla\c{c}ons $f$ par $zf$ et appliquons la proposition 2.1 \`a chaque terme  $S_{{\bf M}'_{j}}^{{\bf G}'_{j}(\tilde{s}_{j})}(\delta_{j,1},(zf)^{{\bf G}'_{j}(\tilde{s}_{j})})$.  Ce terme se d\'eveloppe en une somme sur des espaces de Levi $\tilde{L}'_{j,\tilde{s}_{j}}\in {\cal L}^{\tilde{G}'_{j}}(\tilde{s}_{j})(\tilde{M}'_{j})$. Comme toujours, aux donn\'ees $j$, $\tilde{s}_{j}$ et $\tilde{L}'_{j,\tilde{s}_{j}}$ est associ\'e un $K$-espace de Levi $K\tilde{L}\in {\cal L}(K\tilde{M})$. Le Levi  $\tilde{L}'_{j,\tilde{s}_{j}}$ s'identifie \`a l'espace $\tilde{L}'_{j}(\tilde{s}_{j})$ de la donn\'ee endoscopique ${\bf L}'_{j}(\tilde{s}_{j})$ de $(KL,K\tilde{L},{\bf a}_{L})$ d\'eduite de ${\bf M}'_{j}$ et $\tilde{s}_{j}$. En regroupant les termes selon le $K$espace $K\tilde{L}$, on obtient
   $$I_{K\tilde{M}}^{K\tilde{G},{\cal E}}(\gamma_{i},\omega,zf)=\sum_{K\tilde{L}\in {\cal L}(K\tilde{M})}c^{-1}\vert I\vert ^{-1}\sum_{j\in I}\Delta_{j,1}(\delta_{j,1},\gamma_{i})^{-1}\sum_{\tilde{s}_{j}\in \tilde{\zeta}_{j}Z(\hat{M})^{\Gamma_{{\mathbb R}},\hat{\theta}}/Z(\hat{G})^{\Gamma_{{\mathbb R}},\hat{\theta}}, {\bf L}'_{j}(\tilde{s}_{j})\text{ elliptique}}$$
   $$i_{\tilde{M}_{j}'}(\tilde{G},\tilde{G}'_{j}(\tilde{s}_{j}))S\delta_{{\bf M}'_{j}}^{{\bf L}'_{j}(\tilde{s}_{j})}(\delta_{j},(z_{L})^{{\bf L}'_{j}(\tilde{s}_{j})})S_{{\bf L}'_{j}(\tilde{s}_{j}}^{{\bf G}'_{j}(\tilde{s}_{j})}(\delta_{j,1},f^{{\bf G}'_{j}(\tilde{s}_{j})}).$$
 On a pass\'e 
  quelques consid\'erations formelles similaires \`a celles d\'evelopp\'ees en 1.8 et qu'on laisse au lecteur. On d\'ecompose la somme en $\tilde{s}_{j}$ en une double somme sur $\tilde{s}_{j}\in \tilde{\zeta}_{j}Z(\hat{M})^{\Gamma_{{\mathbb R}},\hat{\theta}}/Z(\hat{L})^{\Gamma_{{\mathbb R}},\hat{\theta}}$ tel que ${\bf L}'_{j}(\tilde{s}_{j})$ soit elliptique et $\tilde{t}_{j}\in \tilde{s}_{j}Z(\hat{L})^{\Gamma_{{\mathbb R}},\hat{\theta}}/Z(\hat{G})^{\Gamma_{{\mathbb R}},\hat{\theta}}$. Comme en [II] 2.6, on a l'\'egalit\'e
  $$i_{\tilde{M}_{j}'}(\tilde{G},\tilde{G}'_{j}(\tilde{t}_{j}))=i_{\tilde{M}_{j}'}(\tilde{L},\tilde{L}'_{j}(\tilde{s}_{j}))i_{\tilde{L}_{j}'(\tilde{s}_{j})}(\tilde{G},\tilde{G}'_{j}(\tilde{t}_{j}))$$
  et la non-nullit\'e de ce dernier terme implique l'ellipticit\'e de ${\bf L}'_{j}(\tilde{s}_{j})$. D'o\`u
   $$(3) \qquad I_{K\tilde{M}}^{K\tilde{G},{\cal E}}(\gamma_{i},\omega,zf)=\sum_{K\tilde{L}\in {\cal L}(K\tilde{M})}c^{-1}\vert I\vert ^{-1}\sum_{j\in I}\Delta_{j,1}(\delta_{j,1},\gamma_{i})^{-1}\sum_{\tilde{s}_{j}\in \tilde{\zeta}_{j}Z(\hat{M})^{\Gamma_{{\mathbb R}},\hat{\theta}}/Z(\hat{L})^{\Gamma_{{\mathbb R}},\hat{\theta}}}i_{\tilde{M}_{j}'}(\tilde{L},\tilde{L}'_{j}(\tilde{s}_{j}))$$
   $$\sum_{\tilde{t}_{j}\in \tilde{s}_{j}Z(\hat{L})^{\Gamma_{{\mathbb R}},\hat{\theta}}/Z(\hat{G})^{\Gamma_{{\mathbb R}},\hat{\theta}}}i_{\tilde{L}_{j}'(\tilde{s}_{j})}(\tilde{G},\tilde{G}'_{j}(\tilde{t}_{j}))S\delta_{{\bf M}'_{j}}^{{\bf L}'_{j}(\tilde{s}_{j})}(\delta_{j},(z_{L})^{{\bf L}'_{j}(\tilde{s}_{j})})S_{{\bf L}'_{j}(\tilde{s}_{j})}^{{\bf G}'_{j}(\tilde{t}_{j})}(\delta_{j,1},f^{{\bf G}'_{j}(\tilde{t}_{j})}).$$
   Pour exploiter cette formule, on a besoin de quelques pr\'eliminaires. On a
   
   (4) pour $i,j,k\in I$, le terme $\Delta_{j,1}(\delta_{j,1},\gamma_{i})^{-1}\Delta_{j,1}(\delta_{j,1},\gamma_{k})$ ne d\'epend ni de $\gamma_{I}$, ni de $\delta_{j,1}$ mais seulement de $i,j,k$.
   
 Remarquons que  
 $$\Delta_{j,1}(\delta_{j,1},\gamma_{i})^{-1}\Delta_{j,1}(\delta_{j,1},\gamma_{k})=\boldsymbol{\Delta}_{j,1}(\delta_{j,1},\gamma_{k};\delta_{j,1},\gamma_{i}).$$
  Fixons $x_{k}\in M_{k}$ tel que $\tilde{\iota}_{i,k}=\phi_{i,k}\circ ad(x_{k})$. D\'ecomposons $x_{k}$ en $x_{k}=\pi(x_{k,sc})z_{k}$, o\`u $z_{k}\in Z(M_{k})$ et $x_{k,sc}\in M_{k,SC}$. On note $T_{k,sc}$ l'image r\'eciproque de $T_{k}$ dans $M_{k,SC}$ et $\hat{T}_{k,ad}$ l'image de $\hat{T}_{k}$ dans $\hat{M}_{AD}$. Pour $\sigma\in \Gamma_{{\mathbb R}}$, posons $V_{k,i}(\sigma)=x_{k,sc}^{-1}\nabla_{k,i}(\sigma)\sigma(x_{k,sc})$, cf. [I] 1.11 pour la d\'efinition de $\nabla_{k,i}$. On v\'erifie que $V_{k,i}$ est un cocycle \`a valeurs dans $T_{k,sc}$. Le couple $(V_{k,i},(1-\theta)(z_{k})^{-1})$  d\'efinit un \'el\'ement de $H^{1,0}(\Gamma_{{\mathbb R}}; T_{k,sc}\stackrel{1-\theta}{\to }(1-\theta)(T_{k}))$. Le cocycle $t_{T_{k}}$ d\'efini en [I] 2.2 se pousse en un cocycle de $W_{{\mathbb R}}$ dans $\hat{T}_{k}/\hat{T}_{k}^{\hat{\theta},0}$.  Rappelons que ${\bf M}'_{k}=(M'_{k},{\cal M}'_{k},\tilde{\zeta}_{k})$ et que l'on suppose implicitement que $\tilde{\zeta}_{k}=\zeta_{k}\hat{\theta}$, avec $\zeta_{k}\in  \hat{T}_{k}$.  Le couple $(t_{T_{k}},\zeta_{k})$ d\'efinit un \'el\'ement de $H^{1,0}(W_{{\mathbb R}};\hat{T}_{k}/\hat{T}_{k}^{\hat{\theta},0}\stackrel{1-\hat{\theta}}{\to}\hat{T}_{k,ad})$. En d\'evissant les d\'efinitions de [I] 2.2, on v\'erifie que 
  $$\boldsymbol{\Delta}_{j,1}(\delta_{j,1},\gamma_{k};\delta_{j,1},\gamma_{i})=<(V_{k,i},(1-\theta)(z_{k})^{-1}), (t_{T_{k}},\zeta_{k})>,$$
  le produit \'etant celui sur
  $$H^{1,0}(\Gamma_{{\mathbb R}}; T_{k,sc}\stackrel{1-\theta}{\to }(1-\theta)(T_{k}))\times H^{1,0}(W_{{\mathbb R}};\hat{T}_{k}/\hat{T}_{k}^{\hat{\theta},0}\stackrel{1-\hat{\theta}}{\to}\hat{T}_{k,ad}).$$
  Cela prouve (4).

  Pour $j,k\in I$, l'alg\`ebre $Diff^{cst}(\tilde{T}'_{j,1}({\mathbb R}))_{\lambda_{j,1}}$ s'identifie \`a $Diff(\tilde{T}_{k}({\mathbb R}))^{\omega-inv}$, cf. 1.8.  L'identification a \'et\'e faite de sorte que la propri\'et\'e suivante soit v\'erifi\'ee. Fixons $\gamma_{k}\in \tilde{T}_{k,\tilde{G}_{k,reg}}({\mathbb R})$ et $\delta_{j,1}\in \tilde{T}'_{j,1}({\mathbb R})$ dont la projection dans $\tilde{T}'_{j}({\mathbb R})$ soit $\tilde{\xi}_{j,k}(\gamma_{k})$. 
  Soit $\varphi$ une fonction sur $\tilde{T}_{k}({\mathbb R})$ telle que $\varphi(x^{-1}\gamma_{k}' x)=\omega(x)^{-1}\varphi(\gamma_{k}')$ pour tous $x\in T_{k}({\mathbb R})$ et $\gamma_{k}'\in \tilde{T}_{k}({\mathbb R})$. On en d\'eduit une fonction $\varphi'$ sur un voisinage de $\delta_{j,1}$ dans  $\tilde{T}'_{j,1}({\mathbb R})$ par la formule $\varphi'(\delta'_{1})=\Delta_{j,1}(\delta'_{1},\gamma')\varphi(\gamma'_{k})$ o\`u $\gamma'_{k}$ est un \'el\'ement de $\tilde{T}_{k}({\mathbb R})$ proche de $\gamma_{k}$  tel que $\tilde{\xi}_{j,k}(\gamma'_{k})$ est \'egal \`a la projection de $\delta'_{1}$ dans $\tilde{T}'_{j}({\mathbb R})$. Soit $D' \in Diff^{cst}(\tilde{T}'_{j,1}({\mathbb R}))_{\lambda_{j,1}}$ et notons $D$ l'\'el\'ement de $Diff^{cst}(\tilde{T}_{k}({\mathbb R}))^{\omega-inv}$ auquel il s'identifie. Alors on a l'\'egalit\'e  $(D'\varphi')(\delta_{j,1})=(D\varphi)'(\delta_{j,1})$. En composant les isomorphismes
  $$Diff^{cst}(\tilde{T}'_{j,1}({\mathbb R}))_{\lambda_{j,1}}\simeq Diff^{cst}(\tilde{T}_{k}({\mathbb R}))^{\omega-inv}\simeq Diff^{cst}(\tilde{T}({\mathbb R}))^{\omega-inv},$$
  on obtient un isomorphisme qui ne d\'epend pas de $k$: cela r\'esulte de la propri\'et\'e ci-dessus et de (4). D\'esormais, on identifie $Diff^{cst}(\tilde{T}'_{j,1}({\mathbb R}))_{\lambda_{j,1}}$ \`a $Diff^{cst}(\tilde{T}({\mathbb R}))^{\omega-inv}$.

  Revenons \`a l'expression (3).  On reconna\^{\i}t la derni\`ere somme en $\tilde{t}_{j}$: c'est
  $$S\delta_{{\bf M}'_{j}}^{{\bf L}'_{j}(\tilde{s}_{j})}(\delta_{j},(z_{L})^{{\bf L}'_{j}(\tilde{s}_{j})})I_{K\tilde{L}}^{K\tilde{G},{\cal E}}({\bf L}'_{j}(\tilde{s}_{j}),\delta_{j,1},f),$$
  ou encore
  $$S\delta_{{\bf M}'_{j}}^{{\bf L}'_{j}(\tilde{s}_{j})}(\delta_{j},(z_{L})^{{\bf L}'_{j}(\tilde{s}_{j})})I_{K\tilde{L}}^{K\tilde{G},{\cal E}}(transfert(\delta_{j,1}),\omega,f).$$
  La formule (1) d\'ecrit le transfert de $\delta_{j,1}$. On obtient
  $$ S\delta_{{\bf M}'_{j}}^{{\bf L}'_{j}(\tilde{s}_{j})}(\delta_{j},(z_{L})^{{\bf L}'_{j}(\tilde{s}_{j})})c \sum_{k\in I}  \Delta_{j,1}(\delta_{j,1},\gamma_{k})I_{K\tilde{L}}^{K\tilde{G},{\cal E}}(\gamma_{k},\omega,f).$$
  En utilisant la description ci-dessus de l'isomorphisme $Diff^{cst}(\tilde{T}'_{j,1}({\mathbb R}))_{\lambda_{j,1}} \simeq Diff^{cst}(\tilde{T}({\mathbb R}))^{\omega-inv}$, on transforme cette expression en
 $$c\sum_{k\in I }\Delta_{j,1}(\delta_{j,1},\gamma_{k}) S\delta_{{\bf M}'_{j}}^{{\bf L}'_{j}(\tilde{s}_{j})}(\delta_{j},(z_{L})^{{\bf L}'_{j}(\tilde{s}_{j})})I_{K\tilde{L}}^{K\tilde{G},{\cal E}}(\gamma_{k},\omega,f).$$
 En reportant cela dans la formule (3), on obtient
$$ I_{K\tilde{M}}^{K\tilde{G},{\cal E}}(\gamma_{i},\omega,zf)=\sum_{K\tilde{L}\in {\cal L}(K\tilde{M})}\sum_{k\in I}\vert I\vert ^{-1}\sum_{j\in I}\Delta_{j,1}(\delta_{j,1},\gamma_{i})^{-1}\Delta_{j,1}(\delta_{j,1},\gamma_{k})$$
$$\sum_{\tilde{s}_{j}\in \tilde{\zeta}_{j}Z(\hat{M})^{\Gamma_{{\mathbb R}},\hat{\theta}}/Z(\hat{L})^{\Gamma_{{\mathbb R}},\hat{\theta}}}i_{\tilde{M}_{j}'}(\tilde{L},\tilde{L}'_{j}(\tilde{s}_{j})) S\delta_{{\bf M}'_{j}}^{{\bf L}'_{j}(\tilde{s}_{j})}(\delta_{j},(z_{L})^{{\bf L}'_{j}(\tilde{s}_{j})})I_{K\tilde{L}}^{K\tilde{G},{\cal E}}(\gamma_{k},\omega,f).$$
En posant gr\^ace \`a (4)
   $$\kappa_{j}(i,k)=\Delta_{j,1}(\delta_{j,1},\gamma_{i})^{-1}\Delta_{j,1}(\delta_{j,1},\gamma_{k}),$$
   on obtient
$$(5) \qquad  I_{K\tilde{M}}^{K\tilde{G},{\cal E}}(\gamma_{i},\omega,zf)=\sum_{K\tilde{L}\in {\cal L}(K\tilde{M})}\sum_{k\in I}\vert I\vert ^{-1}\sum_{j\in I}\kappa_{j}(i,k)$$
$$\sum_{\tilde{s}_{j}\in \tilde{\zeta}_{j}Z(\hat{M})^{\Gamma_{{\mathbb R}},\hat{\theta}}/Z(\hat{L})^{\Gamma_{{\mathbb R}},\hat{\theta}}}i_{\tilde{M}_{j}'}(\tilde{L},\tilde{L}'_{j}(\tilde{s}_{j})) S\delta_{{\bf M}'_{j}}^{{\bf L}'_{j}(\tilde{s}_{j})}(\delta_{j},(z_{L})^{{\bf L}'_{j}(\tilde{s}_{j})})I_{K\tilde{L}}^{K\tilde{G},{\cal E}}(\gamma_{k},\omega,f).$$
  D\'efinissons maintenant les op\'erateurs de l'\'enonc\'e par
$$(6) \qquad \delta_{K\tilde{M},i,k}^{K\tilde{G},{\cal E}}(\gamma_{k},z)=\vert I\vert ^{-1}\sum_{j\in I}\kappa_{j}(i,k)\sum_{\tilde{s}_{j}\in \tilde{\zeta}_{j}Z(\hat{M})^{\Gamma_{{\mathbb R}},\hat{\theta}}/Z(\hat{G})^{\Gamma_{{\mathbb R}},\hat{\theta}}}$$
$$i_{\tilde{M}_{j}'}(\tilde{G},\tilde{G}'_{j}(\tilde{s}_{j})) S\delta_{{\bf M}'_{j}}^{{\bf G}'_{j}(\tilde{s}_{j})}(\delta_{j},z^{{\bf G}'_{j}(\tilde{s}_{j})}).$$
Alors la formule (5) devient celle de l'\'enonc\'e. 

Le tore tordu  $\tilde{T}'({\mathbb C})$ s'identifie \`a $\tilde{T}({\mathbb C})/(1-\theta)(T({\mathbb C}))$ (c'est-\`a-dire \`a $\tilde{T}_{k}({\mathbb C})/(1-\theta)(T_{k}({\mathbb C}))$ pour tout $k$). Les fonctions $\delta_{j}\mapsto  S\delta_{{\bf M}'_{j}}^{{\bf G}'_{j}(\tilde{s}_{j})}(\delta_{j},z^{{\bf G}'_{j}(\tilde{s}_{j})})$ s'\'etendent en des fonctions rationnelles sur ce tore. Chacune est r\'eguli\`ere sur un sous-ensemble qui d\'epend de $j$: l'ensemble des  points $\tilde{G}'_{j}$-r\'eguliers. Mais cet ensemble contient l'ensemble commun des points qui sont $\tilde{G}$-r\'eguliers. Il en r\'esulte que la fonction $\gamma_{k}\mapsto \delta_{K\tilde{M},i,k}^{K\tilde{G},{\cal E}}(\gamma_{k},z)$ s'\'etend en une fonction rationnelle r\'eguli\`eres sur $\tilde{T}_{k,\tilde{G}_{k}-reg}({\mathbb C})/(1-\theta)(T_{k}({\mathbb C}))$. Cela ach\`eve la preuve. 
$\square$
\bigskip

\subsection{Propri\'et\'es des op\'erateurs diff\'erentiels endoscopiques}
Il r\'esulte de la d\'efinition 2.3(6)  et de ce que l'on a dit en 2.2 que les op\'erateurs $\delta_{K\tilde{M},i,k}^{K\tilde{G},{\cal E}}(\gamma_{I},z)$ v\'erifient des propri\'et\'es analogues \`a celles de la proposition 1.3. 

Les diff\'erents homomorphismes $\mathfrak{Z}(G)\to \mathfrak{Z}(T_{i}) \to Sym(\mathfrak{t}_{i})_{\theta,\omega}$, pour $i\in I$, co\"{\i}ncident modulo les identifications de $Sym(\mathfrak{t}_{i})_{\theta,\omega}$ \`a une alg\`ebre commune $Sym(\mathfrak{t})_{\theta,\omega}$. Pour $z\in \mathfrak{Z}(G)$, on note $z_{T}$ l'image de $z$ dans cette alg\`ebre. On a

(1) si $K\tilde{M}=K\tilde{G}$, la matrice $\delta_{K\tilde{G}}^{K\tilde{G}}(\gamma_{I},z)$ est diagonale, de termes diagonaux tous \'egaux \`a $\partial_{z_{T}}$.

En effet, la d\'efinition 2.3(6) se simplifie en
$$\delta_{K\tilde{G},i,k}^{K\tilde{G},{\cal E}}(\gamma_{I},z)=\vert I\vert ^{-1}\sum_{j\in I}\kappa_{j}(i,k)S^{{\bf G}'_{j}}(\delta_{j},z^{{\bf G}'_{j}}).$$
D'apr\`es 2.2(1), l'op\'erateur $S^{{\bf G}'_{j}}(\delta_{j},z^{{\bf G}'_{j}})$ est simplement  $\partial_{z^{{\bf G}'_{j}}}$, o\`u on identifie $z^{{\bf G}'_{j}}$ \`a son image naturelle  dans $Sym(\mathfrak{t}'_{j,1})_{\lambda_{j,1}}$ (avec les notations de la preuve pr\'ec\'edente). Modulo l'identification  $Sym(\mathfrak{t}'_{j,1})_{\lambda_{j,1}}\simeq Sym(\mathfrak{t})_{\theta,\omega}$, ce n'est autre que $\partial_{z_{T}}$. La formule devient
$$\delta_{K\tilde{G},i,k}^{K\tilde{G},{\cal E}}(\gamma_{I},z)=\partial_{z_{T}}\vert I\vert ^{-1}\sum_{j\in I}\kappa_{j}(i,k).$$
Il est sous-jacent aux formules d'inversion de [I] 4.9 que la correspondance entre les $\gamma_{i}$ et les $\delta_{i}$ est essentiellement une transformation de Fourier entre un groupe ab\'elien fini et son dual, de sorte que la somme en $j$ ci-dessus vaut $\vert I\vert$  si $i=k$ et $0$ sinon. Cette propri\'et\'e est largement d\'evelopp\'ee dans [L] et [KS] et on ne la d\'etaillera pas. Elle entra\^{\i}ne la conclusion (1).

Posons $d=a_{\tilde{M}}-a_{\tilde{G}}$, supposons $d\geq1$. Pour tout $i\in I$, on a d\'efini en 1.7 une fonction $C_{\tilde{M}_{i}}^{\tilde{G}_{i}}$ sur $\tilde{T}_{i,\tilde{G}-reg}({\mathbb C})$. Ces fonctions s'identifient par les isomorphismes $\tilde{\iota}_{i,j}$. On en d\'eduit une fonction $C_{K\tilde{M}}^{K\tilde{G}}$ sur $\tilde{T}_{\tilde{G}-reg}({\mathbb C})$. Soit $\Omega$ l'op\'erateur de Casimir.

\ass{Proposition}{Pour $d\geq1$, la matrice $\delta_{K\tilde{M}}^{K\tilde{G}}(\gamma_{I},\Omega)$ est diagonale. Ses coefficients diagonaux sont tous \'egaux \`a la multiplication par $C_{K\tilde{M}}^{K\tilde{G}}(\gamma_{I})$.}

Preuve. On va plus g\'en\'eralement prouver que ces propri\'et\'es sont v\'erifi\'ees si l'on remplace $\Omega$ par un \'el\'ement $\underline{\Omega}\in \Omega+Sym(\mathfrak{z}(G))$. On consid\`ere la formule 2.3(6). Les m\^emes consid\'erations formelles que dans la preuve de la proposition 2.2 montre que, pour $j$ et $\tilde{s}_{j}$ intervenant dans cette formule, l'op\'erateur $S\delta_{{\bf M}'_{j}}^{{\bf G}'_{j}(\tilde{s}_{j})}(\delta_{j},\underline{\Omega}^{{\bf G}'_{j}(\tilde{s}_{j})})$ est la multiplication par la fonction $SC_{\tilde{M}'_{j}}^{\tilde{G}'_{j}(\tilde{s}_{j})}(\delta_{j})$. Ainsi, pour $i,k\in I$, $\delta_{K\tilde{M},i,k}^{K\tilde{G}}(\gamma_{I},\underline{\Omega})$ est la multiplication par la fonction
$$\vert I\vert ^{-1}\sum_{j\in I}\kappa_{j}(i,k)X_{j}(\delta_{j}),$$
o\`u
$$X_{j}(\delta_{j})=\sum_{\tilde{s}_{j}\in \tilde{\zeta}_{j}Z(\hat{M})^{\Gamma_{{\mathbb R}},\hat{\theta}}/Z(\hat{G})^{\Gamma_{{\mathbb R}},\hat{\theta}}}i_{\tilde{M}_{j}'}(\tilde{G},\tilde{G}'_{j}(\tilde{s}_{j})) SC_{\tilde{ M}'_{j}}^{\tilde{ G}'_{j}(\tilde{s}_{j})}(\delta_{j}).$$
On va prouver

(2) on a $X_{j}(\delta_{j})=C_{K\tilde{M}}^{K\tilde{G}}(\gamma_{I})$ pour tout $j\in I$. 

En admettant cela, on obtient que $\delta_{K\tilde{M},i,k}^{K\tilde{G}}(\gamma_{I},\underline{\Omega})$ est la multiplication par la fonction
$$C_{K\tilde{M}}^{K\tilde{G}}(\gamma_{I})\vert I\vert ^{-1}\sum_{j\in I}\kappa_{j}(i,k).$$
La proposition s'en d\'eduit par le m\^eme argument d'inversion utilis\'e dans la preuve de (1).

Prouvons (2). Cette \'egalit\'e est \'evidente si $d\geq2$ car alors tout est nul. On suppose $d=1$. Pour simplifier, on fixe $j$ et on abandonne les indices $j$. De m\^eme, on fixe un certain $i\in I$, on  identifie 
$C_{K\tilde{M}}^{K\tilde{G}}(\gamma_{I})$ \`a $C_{\tilde{M}_{i}}^{\tilde{G}_{i}}(\gamma_{i})$ et on abandonne l'indice $i$. Soit $\tilde{s}\in \tilde{\zeta}Z(\hat{M})^{\Gamma_{{\mathbb R}},\hat{\theta}}/Z(\hat{G})^{\Gamma_{{\mathbb R}},\hat{\theta}}$, que l'on \'ecrit $\tilde{s}=s\hat{\theta}$. D'apr\`es [KS] 1.3, repris en [I] 1.6, les \'el\'ements de $\Sigma^{G'(\tilde{s})}(T')-\Sigma^{M'}(T')$ sont les

 $N\alpha$ pour $\alpha\in \Sigma^G(T)_{\theta}-\Sigma^M(T)_{\theta}$ de type 1 tel que $N\hat{\alpha}(s)=1$;

 $2N\alpha$ pour $\alpha\in \Sigma^G(T)_{\theta}-\Sigma^M(T)_{\theta}$ de type 2 tel que $N\hat{\alpha}(s)=1$;

 $N\alpha$ pour $\alpha\in \Sigma^G(T)_{\theta}-\Sigma^M(T)_{\theta}$ de type 3 tel que $N\hat{\alpha}(s)=-1$.

Fixons comme toujours une paire de Borel \'epingl\'ee de $\hat{G}$ pour laquelle $\hat{M}$ est standard et dont on note le tore $\hat{T}$. Dualement, la racine de $\hat{G}(\tilde{s})$ correspondant \`a $N\alpha$, pour $\alpha$ de type 1 ou 3, ou \`a $2N\alpha$ pour $\alpha$ de type 2,  est la restriction de $\hat{\alpha}$ \`a $\hat{T}^{\hat{\theta},0}$. 

On se rappelle que, d'une racine $\alpha$ de type 2, se d\'eduit une racine $\beta=\alpha+\theta^{n_{\alpha}/2}(\alpha)$ de type 3. On peut choisir l'ensemble $\Sigma^G(T)_{\theta}$ de sorte que, de cette correspondance,  se d\'eduise une bijection entre les \'el\'ements de $\Sigma^G(T)_{\theta}$ de type 2 et les \'el\'ements de $\Sigma^G(T)_{\theta}$ de type 3. On a $\beta\in \Sigma^M(T)$ si et seulement si $\alpha\in \Sigma^M(T)$. On peut donc supposer que  la bijection ci-dessus se restreint en une bijection entre les \'el\'ements de $\Sigma^G(T)_{\theta}-\Sigma^M(T)_{\theta}$ de type 2 et ceux de type 3. 

Ainsi, $X(\delta)$ comme $C_{\tilde{M}}^{\tilde{G}}(\gamma)$ apparaissent comme des sommes sur les $\alpha\in \Sigma^G(T)_{\theta}-\Sigma^M(T)_{\theta}$.  On va prouver

(3) pour $\alpha \in \Sigma^G(T)_{\theta}-\Sigma^M(T)_{\theta}$ de type 1, les termes index\'es par $\alpha$ dans $X(\delta)$ et dans $C_{\tilde{M}}^{\tilde{G}}(\gamma)$ sont \'egaux;

(4) soit $\alpha \in \Sigma^G(T)_{\theta}-\Sigma^M(T)_{\theta}$ de type 2; soit $\beta\in \Sigma^G(T)_{\theta}-\Sigma^M(T)_{\theta}$ l'\'el\'ement de type 3 correspondant; alors la somme des termes index\'es par $\alpha$ et $\beta$ dans $X(\delta)$ est \'egale \`a la somme des termes index\'es par $\alpha$ et $\beta$ dans $C_{\tilde{M}}^{\tilde{G}}(\gamma)$.

Soit $\alpha \in \Sigma^G(T)_{\theta}-\Sigma^M(T)_{\theta}$ de type 1. En utilisant les d\'efinitions de 1.7 et 2.2, prouver (3) se ram\`ene \`a prouver l'\'egalit\'e
$$-\vert N\alpha_{\vert {\cal A}_{\tilde{M}}}\vert n_{\alpha}(1-(\alpha)(\gamma))^{-1}(1-(\alpha)(\gamma)^{-1})^{-1}=-\sum_{\tilde{s}\in \tilde{\zeta}Z(\hat{M})^{\Gamma_{{\mathbb R}},\hat{\theta}}/Z(\hat{G})^{\Gamma_{{\mathbb R}},\hat{\theta}}, N\alpha\in \Sigma^{\tilde{G}'(\tilde{s})}}$$
$$i_{\tilde{M}'}(\tilde{G},\tilde{G}'(\tilde{s}))\vert N\alpha_{\vert {\cal A}_{\tilde{M}'}}\vert k^{\hat{G}'(\tilde{s})}(\hat{\alpha})^{-1}(1-N\alpha(\delta))^{-1}(1-N\alpha(\delta)^{-1})^{-1}.$$
Puisque ${\cal A}_{\tilde{M}}={\cal A}_{\tilde{M}'}$, les termes $\vert N\alpha _{\vert {\cal A}_{\tilde{M}}}\vert $ et $\vert N\alpha_{\vert {\cal A}_{\tilde{M}'}}\vert$ sont \'egaux. De m\^eme, il r\'esulte des d\'efinitions que $1- (\alpha)(\gamma)^{\pm 1}=1- N\alpha(\delta)^{\pm 1}$. L'\'egalit\'e \`a prouver se simplifie en
$$(5) \qquad n_{\alpha}=\sum_{\tilde{s}\in \tilde{\zeta}Z(\hat{M})^{\Gamma_{{\mathbb R}},\hat{\theta}}/Z(\hat{G})^{\Gamma_{{\mathbb R}},\hat{\theta}}, N\alpha\in \Sigma^{\tilde{G}'(\tilde{s})}}i_{\tilde{M}'}(\tilde{G},\tilde{G}'(\tilde{s}))  k^{\hat{G}'(\tilde{s})}(\hat{\alpha})^{-1}.$$
Posons $\hat{{\cal T}}=(Z(\hat{M})^{\Gamma_{{\mathbb R}}}\cap \hat{T}^{\hat{\theta},0})/(Z(\hat{G})^{\Gamma_{{\mathbb R}}}\cap \hat{T}^{\hat{\theta},0})$. 
Remarquons que les homomorphismes de la suite suivante
$$Z(\hat{M})^{\Gamma_{{\mathbb R}},\hat{\theta},0} /( Z(\hat{M})^{\Gamma_{{\mathbb R}},\hat{\theta},0}\cap Z(\hat{G}))\to 
 \hat{{\cal T}}\to Z(\hat{M})^{\Gamma_{{\mathbb R}},\hat{\theta}}/Z(\hat{G})^{\Gamma_{{\mathbb R}},\hat{\theta}}$$
 sont des isomorphismes. Ce groupe commun est un tore complexe de dimension $1$. La restriction de $\hat{\alpha}$ \`a ce tore n'est pas triviale. Il en r\'esulte que l'on peut trouver $\tilde{s}\in \tilde{\zeta}Z(\hat{M})^{\Gamma_{{\mathbb R}},\hat{\theta}}/Z(\hat{G})^{\Gamma_{{\mathbb R}},\hat{\theta}}$ tel que $N\hat{\alpha}(s)=1$, autrement dit la somme en $\tilde{s}$ ci-dessus est non vide. On ne perd alors rien \`a supposer que $\tilde{\zeta}$ lui-m\^eme v\'erifie cette condition. En posant $\tilde{s}=\tilde{\zeta}t$ avec $t\in  \hat{{\cal T}}$, la somme sur les $\tilde{s}$ tels que $N\alpha\in \Sigma^{\tilde{G}'(\tilde{s})}$ devient une somme sur les $t$ tels que $N\hat{\alpha}(t)=1$. Puisque $t$ est invariant par $\hat{\theta}$, cette derni\`ere condition \'equivaut \`a $\hat{\alpha}(t)^{n_{\alpha}}=1$. Consid\'erons un tel $t$ et $\tilde{s}=\tilde{\zeta}t$ l'\'el\'ement correspondant. Pour $\beta\in \Sigma^G(T)_{\theta}-\Sigma^M(T)_{\theta}$, notons $\hat{\beta}_{*}$ la restriction de $\hat{\beta}$ \`a $Z(\hat{M}')^{\Gamma_{{\mathbb R}}}$.  Parmi les $\beta$ tels que  $\hat{\beta}\in \Sigma^{\hat{G}'(\tilde{s})} $,  on peut en fixer un tel que que $\hat{\beta}_{*}$   soit indivisible et que $\hat{\alpha}_{*}=k^{\hat{G}'(\tilde{s})}(\hat{\alpha})\hat{\beta}_{*}$. Le groupe $Z(\hat{G}'(\tilde{s}))^{\Gamma_{{\mathbb R}}}$ est \'egal \`a celui des $z\in Z(\hat{M}')^{\Gamma_{{\mathbb R}}}$ tel que $\hat{\beta}(z)=1$. Il contient le groupe des $z\in Z(\hat{M})^{\Gamma_{{\mathbb R}}}\cap \hat{T}^{\hat{\theta},0}$ tels que $\hat{\beta}(z)=1$. Donc 
 $$[Z(\hat{G}'(\tilde{s}))^{\Gamma_{{\mathbb R}}}:Z(\hat{G})^{\Gamma_{{\mathbb R}}}\cap \hat{T}^{\hat{\theta},0}]=[\{z\in Z(\hat{M}')^{\Gamma_{{\mathbb R}}}; \hat{\beta}(z)=1\}:\{z\in Z(\hat{M})^{\Gamma_{{\mathbb R}}}\cap \hat{T}^{\hat{\theta},0}; \hat{\beta}(z)=1\}]$$
 $$[\{z\in Z(\hat{M})^{\Gamma_{{\mathbb R}}}\cap \hat{T}^{\hat{\theta},0}; \hat{\beta}(z)=1\}:Z(\hat{G})^{\Gamma_{{\mathbb R}}}\cap \hat{T}^{\hat{\theta},0}].$$
  Puisque $\hat{\beta}$ n'est pas trivial sur $Z(\hat{M})^{\Gamma_{{\mathbb R}}}\cap \hat{T}^{\hat{\theta},0}$, l'homomorphisme
 $$\{z\in Z(\hat{M}')^{\Gamma_{{\mathbb R}}}; \hat{\beta}(z)=1\}/\{z\in Z(\hat{M})^{\Gamma_{{\mathbb R}}}\cap \hat{T}^{\hat{\theta},0}; \hat{\beta}(z)=1\}\to Z(\hat{M}')^{\Gamma_{{\mathbb R}}}/(Z(\hat{M})^{\Gamma_{{\mathbb R}}}\cap \hat{T}^{\hat{\theta},0})$$
 est un isomorphisme. D'autre part, $\hat{\beta}$ se quotiente en un caract\`ere du tore $\hat{{\cal T}}$.  Il en r\'esulte que 
$$[Z(\hat{G}'(\tilde{s}))^{\Gamma_{{\mathbb R}}}:Z(\hat{G})^{\Gamma_{{\mathbb R}}}\cap \hat{T}^{\hat{\theta},0}]=[Z(\hat{M}')^{\Gamma_{{\mathbb R}}}:Z(\hat{M})^{\Gamma_{{\mathbb R}}}\cap \hat{T}^{\hat{\theta},0}] $$
$$\vert \{z\in  \hat{{\cal T}}; \hat{\beta}(z)=1\}\vert .$$
En se rappelant la d\'efinition de $i_{\tilde{M}'}(\tilde{G},\tilde{G}'(\tilde{s}))$,  cf. [II] 1.7, on obtient
$$i_{\tilde{M}'}(\tilde{G},\tilde{G}'(\tilde{s}))=\vert \{z\in \hat{{\cal T}}; \hat{\beta}(z)=1\}\vert ^{-1}.$$
On a un homomorphisme surjectif
$$\hat{{\cal T}}\to Z(\hat{M}')^{\Gamma_{{\mathbb R}}}/Z(\hat{G}'(\tilde{s}))^{\Gamma_{{\mathbb R}}}.$$
Puisque la  restriction de $\hat{\alpha}$ au tore de droite est    \'egale \`a $k^{\hat{G}'(\tilde{s})}(\hat{\alpha})$ fois celle de $\hat{\beta}$, on a la m\^eme relation entre les restrictions de ces racines \`a $\hat{{\cal T}}$. Puisque ce tore est de dimension $1$, on obtient
$$ k^{\hat{G}'(\tilde{s})}(\hat{\alpha})^{-1}i_{\tilde{M}'}(\tilde{G},\tilde{G}'(\tilde{s}))=\vert \{z\in \hat{{\cal T}}; \hat{\alpha}(z)=1\}\vert ^{-1},$$
ou encore, par le m\^eme argument
$$ k^{\hat{G}'(\tilde{s})}(\hat{\alpha})^{-1}i_{\tilde{M}'}(\tilde{G},\tilde{G}'(\tilde{s}))=n_{\alpha}\vert \{z\in \hat{{\cal T}}; N\hat{\alpha}(z)=1\}\vert ^{-1}.$$
La relation (5) \`a prouver se r\'eduit \`a
$$n_{\alpha}=\sum_{t\in \hat{{\cal T}},N\hat{\alpha}(t)=1 }n_{\alpha}\vert \{z\in \hat{{\cal T}}; N\hat{\alpha}(z)=1\}\vert ^{-1}.$$
C'est \'evident et cela prouve (3).

Soit $\alpha\in \Sigma^G(T)_{\theta}-\Sigma^M(T)_{\theta}$ de type 2 et soit $\beta=\alpha+\theta^{n_{\alpha}/2}(\alpha)$ l'\'el\'ement de type 3 associ\'e.   L'\'egalit\'e (4) \`a prouver se ram\`ene \`a $$(6) \qquad A_{\alpha}+A_{\beta}=B_{\alpha}+B_{\beta},$$
 o\`u
$$A_{\alpha}=-\vert N\alpha _{\vert {\cal A}_{\tilde{M}}}\vert n_{\alpha}(1-(\alpha)(\gamma))^{-1}(1-(\alpha)(\gamma)^{-1})^{-1} ,$$
$$A_{\beta}=-\vert N\beta _{\vert {\cal A}_{\tilde{M}}}\vert n_{\beta}(1-(\beta)(\gamma))^{-1}(1-(\beta)(\gamma)^{-1})^{-1},$$
$$B_{\alpha}=-\sum_{\tilde{s}\in \tilde{\zeta}Z(\hat{M})^{\Gamma_{{\mathbb R}},\hat{\theta}}/Z(\hat{G})^{\Gamma_{{\mathbb R}},\hat{\theta}}, 2N\alpha\in \Sigma^{\tilde{G}'(\tilde{s})}}i_{\tilde{M}'}(\tilde{G},\tilde{G}'(\tilde{s}))\vert 2N\alpha_{\vert {\cal A}_{\tilde{M}'}}\vert k^{\hat{G}'(\tilde{s})}(\hat{\alpha})^{-1}$$
$$(1-(2N\alpha)(\delta))^{-1}(1-(2N\alpha)(\delta)^{-1})^{-1},$$
$$B_{\beta}=-\sum_{\tilde{s}\in \tilde{\zeta}Z(\hat{M})^{\Gamma_{{\mathbb R}},\hat{\theta}}/Z(\hat{G})^{\Gamma_{{\mathbb R}},\hat{\theta}}, N\beta\in \Sigma^{\tilde{G}'(\tilde{s})}}i_{\tilde{M}'}(\tilde{G},\tilde{G}'(\tilde{s}))\vert N\beta_{\vert {\cal A}_{\tilde{M}'}}\vert k^{\hat{G}'(\tilde{s})}(\hat{\beta})^{-1}$$
$$(1-N\beta(\delta))^{-1}(1-N\beta(\delta)^{-1})^{-1}.$$
On a $N\alpha=N\beta$ et $\vert N\alpha _{\vert {\cal A}_{\tilde{M}}}\vert =\vert N\alpha_{\vert {\cal A}_{\tilde{M}'}}\vert $. D'apr\`es les d\'efinitions, $(\alpha)(\gamma)=N\alpha(\delta)$ et $(\beta)(\gamma)=-N\alpha(\delta)$. En posant $x=N\alpha(\delta)$ et $c=-\vert N\alpha _{\vert {\cal A}_{\tilde{M}}}\vert$, on obtient
$$A_{\alpha}=cn_{\alpha}(1-x)^{-1}(1-x^{-1})^{-1},$$
$$A_{\beta}=cn_{\beta}(1+x)^{-1}(1+x^{-1})^{-1},$$
$$B_{\alpha}=2c (1-x^2)^{-1}(1-x^{-2})^{-1}b_{\alpha},$$
$$B_{\beta}=c(1-x)^{-1}(1-x^{-1})^{-1}b_{\beta},$$
o\`u
$$b_{\alpha}=\sum_{\tilde{s}\in \tilde{\zeta}Z(\hat{M})^{\Gamma_{{\mathbb R}},\hat{\theta}}/Z(\hat{G})^{\Gamma_{{\mathbb R}},\hat{\theta}}, 2N\alpha\in \Sigma^{\tilde{G}'(\tilde{s})}}i_{\tilde{M}'}(\tilde{G},\tilde{G}'(\tilde{s}))k^{\hat{G}'(\tilde{s})}(\hat{\alpha})^{-1},$$
$$b_{\beta}=\sum_{\tilde{s}\in \tilde{\zeta}Z(\hat{M})^{\Gamma_{{\mathbb R}},\hat{\theta}}/Z(\hat{G})^{\Gamma_{{\mathbb R}},\hat{\theta}}, N\beta\in \Sigma^{\tilde{G}'(\tilde{s})}}i_{\tilde{M}'}(\tilde{G},\tilde{G}'(\tilde{s}))k^{\hat{G}'(\tilde{s})}(\hat{\beta})^{-1}.$$
Les conditions $2N\alpha\in \Sigma^{\tilde{G}'(\tilde{s})}$ ou  $N\beta\in \Sigma^{\tilde{G}'(\tilde{s})}$ \'equivalent \`a $N\hat{\alpha}(s)=1$, resp. $N\hat{\beta}(s)=-1$. Le m\^eme calcul que dans la preuve de (3) conduit aux \'egalit\'es $b_{\alpha}=n_{\alpha}$, $b_{\beta}=n_{\beta}$. On a l'\'egalit\'e $n_{\alpha}=2n_{\beta}$. L'\'egalit\'e (6) \'equivaut alors \`a
$$2(1-x)^{-1}(1-x^{-1})^{-1}+(1+x)^{-1}(1+x^{-1})^{-1}=4(1-x^2)^{-1}(1-x^{-2})^{-1}+(1-x)^{-1}(1-x^{-1})^{-1}.$$
Cela r\'esulte d'un calcul imm\'ediat. Cela prouve (4) et la proposition. $\square$

\bigskip

\subsection{Le r\'esultat de stabilisation}
Les hypoth\`eses et notations sont comme en 2.3. Soit  $z\in \mathfrak{Z}(G)$. D'apr\`es 1.6(4), les diff\'erents op\'erateurs $\delta_{\tilde{M}_{i}}^{\tilde{G}_{i}}(z)$ s'identifient tous \`a un m\^eme \'el\'ement de $Diff^{reg}(\tilde{T}_{\tilde{G}-reg}({\mathbb R}))^{\omega-inv}$. On note  $\delta_{K\tilde{M}}^{K\tilde{G}}(z)$ l'\'el\'ement  diagonal de $Mat_{I}(Diff^{reg}(\tilde{T}_{\tilde{G}-reg}({\mathbb R}))^{\omega-inv})$ dont les coefficients diagonaux sont tous \'egaux \`a cet \'el\'ement.
\ass{Proposition}{Pour tout  $z\in \mathfrak{Z}(G)$,  on a l'\'egalit\'e
$$\delta_{K\tilde{M}}^{K\tilde{G},{\cal E}}(z)=\delta_{K\tilde{M}}^{K\tilde{G}}(z).$$}

Preuve. C'est clair si $K\tilde{M}=K\tilde{G}$. On suppose $K\tilde{M}\not=K\tilde{G}$. Soient $z,z'\in \mathfrak{Z}(G)$.  La relation 1.2(3) donne
$$(1) \qquad \delta_{K\tilde{M}}^{K\tilde{G}}(zz')=\sum_{K\tilde{L}\in {\cal L}(K\tilde{M})}\delta_{K\tilde{M}}^{K\tilde{L}}(z_{L})\delta_{K\tilde{L}}^{K\tilde{G}}(z').$$
Il s'agit ici de produits de matrices d'op\'erateurs diff\'erentiels. Une m\^eme preuve s'applique aux op\'erateurs endoscopiques et donne la relation similaire:
$$(2) \qquad \delta_{K\tilde{M}}^{K\tilde{G},{\cal E}}(zz')=\sum_{K\tilde{L}\in {\cal L}(K\tilde{M})}\delta_{K\tilde{M}}^{K\tilde{L},{\cal E}}(z_{L})\delta_{K\tilde{L}}^{K\tilde{G},{\cal E}}(z').$$
On raisonne par r\'ecurrence sur $dim(G_{SC})$, mais aussi par r\'ecurrence sur $a_{\tilde{M}}-a_{\tilde{G}}$. On peut donc supposer que l'\'egalit\'e de l'\'enonc\'e est v\'erifi\'ee si l'on remplace $\tilde{G}$ par $\tilde{L}$ avec $\tilde{L}\not=\tilde{G}$ et aussi si l'on remplace $\tilde{M}$ par $\tilde{L}$ pour $\tilde{L}\not=\tilde{M}$. Alors les termes des deux membres de droite ci-dessus index\'es par des $\tilde{L}$ diff\'erents de $\tilde{M}$ et $\tilde{G}$ sont \'egaux. D'autre part, d'apr\`es 1.2(2) et 2.4(1), on a
$$\delta_{K\tilde{G}}^{K\tilde{G}}(z')=\delta_{K\tilde{G}}^{K\tilde{G},{\cal E}}(z')=diag(\partial_{z'_{T}}),$$
$$\delta_{K\tilde{M}}^{K\tilde{M}}(z_{M})=\delta_{K\tilde{M}}^{K\tilde{M},{\cal E}}(z_{M})=diag(\partial_{z_{T}}).$$
Pour $D\in Diff^{cst}(\tilde{T}({\mathbb R}))^{\omega-inv}$, on a ici identifi\'e $ D$ \`a la fonction constante sur $\tilde{T}_{\tilde{G}-reg}({\mathbb R})$ de valeur $ D$, et on a not\'e $diag(D)$  la matrice diagonale de coefficients diagonaux tous \'egaux \`a $ D$.
Posons 
$$A(z)=\delta_{K\tilde{M}}^{K\tilde{G},{\cal E}}(z)-\delta_{K\tilde{M}}^{K\tilde{G}}(z).$$
Par diff\'erence de (2) et (1), on obtient
$$ A(zz')=diag(\partial_{z_{T}})A(z')+A(z)diag(\partial_{z'_{T}}).$$
Puisque le membre de gauche est sym\'etrique en $z$ et $z'$, on obtient
$$(3) \qquad diag(\partial_{z_{T}})A(z')+A(z)diag(\partial_{z_{T}})=diag(\partial_{z'_{T}})A(z)+A(z')diag(\partial_{z_{T}}).$$

De m\^eme que l'on a d\'efini $\tilde{T}$ comme l'ensemble des familles $\gamma_{I}=(\gamma_{i})_{i\in I}$ telles que $\gamma_{i}\in \tilde{T}_{i}$ pour tout $i$ et $\gamma_{i}=\tilde{\iota}_{i,j}(\gamma_{j})$ pour tous $i,j$, on peut d\'efinir les ensembles $T$, $\mathfrak{t}$ et $X_{*}(T)$. On a une forme quadratique d\'efinie positive sur chaque $X_{*}(T_{i})\otimes_{{\mathbb Z}}{\mathbb R}$. Elles se correspondent par les isomorphismes $\iota_{i,j}$. Donc $X_{*}(T)\otimes_{{\mathbb Z}}{\mathbb R}$ est muni d'une telle forme. Cet espace est aussi muni d'un automorphisme $\theta$. Fixons une  base orthonorm\'ee $\{H_{1},...,H_{m}\}$  de $X_{*}(T)^{\theta}\otimes_{{\mathbb Z}}{\mathbb R}$. C'est aussi une base sur ${\mathbb C}$ de $\mathfrak{t}^{\theta}({\mathbb C})$. On dispose d'un homomorphisme $H\mapsto \partial_{H}$ de $\mathfrak{t}^{\theta}({\mathbb R})$ dans $Diff^{cst}(\tilde{T}({\mathbb R}))^{\omega-inv}$. Par lin\'earit\'e, il s'\'etend en un isomorphisme de $Sym(\mathfrak{t}^{\theta}({\mathbb C}))$ sur $Diff^{cst}(\tilde{T}({\mathbb R}))^{\omega-inv}$. Pour $\underline{l}=(l_{1},...,l_{m})\in {\mathbb N}^m$, on pose 
$$\partial^{\underline{l}}=\partial_{H_{1}}^{l_{1}}...\partial_{H_{m}}^{l_{m}}.$$
Ainsi, pour $i,j\in I$ et  $\gamma_{I}\in \tilde{T}_{\tilde{G}-reg}({\mathbb R})$, on peut \'ecrire de fa\c{c}on unique
$$(4) \qquad A_{i,j}(\gamma_{I},z)=\sum_{\underline{l}\in {\mathbb N}^m}a_{i,j;\underline{l}}(\gamma_{I}) \partial^{\underline{l}},$$
o\`u les $a_{i,j;\underline{l}}$ sont des fonctions rationnelles en $\gamma_{I}$. 

D\'efinissons l'op\'erateur
$$\nabla=\partial_{H_{1}}^2+...+\partial_{H_{m}}^2.$$
 Introduisons le Casimir $\Omega$.  Il r\'esulte de 1.7(1) que l'on peut choisir $\underline{\Omega}\in \Omega+Sym(\mathfrak{z}(G))$ tel que   $\partial_{\underline{\Omega}_{T}}=\nabla$. Les propositions 1.7 et 2.4 impliquent que $A(\underline{\Omega})=0$. Pour  $z'=\underline{\Omega}$, la relation (3) se simplifie en
$$(5) \qquad A(z)diag(\nabla)=diag(\nabla)A(z).$$
Fixons $\gamma_{I}\in \tilde{T}_{\tilde{G}-reg}({\mathbb R})$. Pour $H\in \mathfrak{t}^{\theta}({\mathbb R})$, on d\'efinit l'\'el\'ement $exp(H)\in T({\mathbb R})$. Si $H$ est assez voisin de $0$, le produit $exp(H)\gamma_{I}$ appartient encore \`a $\tilde{T}_{\tilde{G}-reg}({\mathbb R})$. Pour $i,j\in I$ et $\underline{l}\in {\mathbb N}^m$, posons simplement $b_{i,j;\underline{l}}(H)=a_{i,j;\underline{l}}(exp(H)\gamma_{I})$. Alors
$$A_{i,j}(exp(H)\gamma_{I},z)=\sum_{\underline{l}\in {\mathbb N}^m}b_{i,j;\underline{l}}(H) \partial^{\underline{l}}.$$
Montrons que

(6) pour tous $i,j\in I$ et $\underline{l}\in {\mathbb N}^m$, $b_{i,j;\underline{l}}(H)$ est un polyn\^ome en $H$.

Tout \'el\'ement $H\in \mathfrak{t}^{\theta}({\mathbb C})$ agit sur les fonctions d\'efinies sur $\mathfrak{t}^{\theta}({\mathbb R})$ par l'op\'erateur de d\'erivation usuel $\partial_{H}$. Pour $k=1,...,m$, notons $1_{k}$  l'\'el\'ement de ${\mathbb N}^m$ dont toutes les coordonn\'ees sont nulles, sauf la $k$-i\`eme qui vaut $1$. Les r\`egles usuelles de d\'erivation montrent que l'\'egalit\'e (5) \'evalu\'ee au point $exp(H)\gamma_{I}$ \'equivaut aux \'egalit\'es 
$$\sum_{\underline{l}\in {\mathbb N}^m}\sum_{k=1,...,m}(\partial_{H_{k}}^2b_{i,j;\underline{l}})(H)\partial^{\underline{l}}+2\sum_{\underline{l}\in {\mathbb N}^m}\sum_{k=1,...,m}(\partial_{H_{k}}b_{i,j;\underline{l}})(H)\partial^{\underline{l}+1_{k}}=0,$$
pour tous $i,j\in I$. 
Fixons $i$ et $j$ et abandonnons ces indices pour simplifier la notation. L'\'egalit\'e pr\'ec\'edente se r\'ecrit
$$\sum_{\underline{l}\in {\mathbb N}^m}\sum_{k=1,...,m}\left((\partial_{H_{k}}^2b_{\underline{l}})(H)+2(\partial_{H_{k}}b_{\underline{l}-1_{k}})(H)\right)\partial^{\underline{l}}=0,$$
avec la convention $b_{\underline{l}}=0$ si une coordonn\'ee de $\underline{l}$ est stritement n\'egative. Ou encore, pour tout $\underline{l}\in {\mathbb N}^m$,
$$(7) \qquad \sum_{k=1,...,m}\partial_{H_{k}}^2b_{\underline{l}}+2\partial_{H_{k}}b_{\underline{l}-1_{k}}=0.$$
Pour tout $\underline{l}$, posons $S(\underline{l})=l_{1}+...+l_{m}$. On peut supposer que les $b_{\underline{l}}$ ne sont pas tous nuls, sinon (6) est clair. On note alors $N$ le plus grand entier pour lequel il existe $\underline{l}$ avec $S(\underline{l})=N$ et $b_{\underline{l}}\not=0$. On va prouver que, pour tout $n=1,...,m$ et tout $\underline{l}$ tel que $S(\underline{l})\leq N$, on a

(8) $\partial_{H_{n}}^{N-l_{n}+1}b_{\underline{l}}=0$. 

On raisonne par r\'ecurrence descendante sur $S=S(\underline{l})$ et, cet entier \'etant fix\'e, sur $l_{n}\in \{0,...,S\}$. Soit $\underline{l}\in {\mathbb N}^m$ tel que $S(\underline{l})=S\leq N$. On applique (7) \`a $\underline{l}+1_{n}$. On applique l'op\'erateur $\partial_{H_{n}}^{N-l_{n}}$ \`a l'\'egalit\'e obtenue. L'\'el\'ement $b_{\underline{l}+1_{n}}$ est annul\'e par $\partial_{H_{n}}^{N-l_{n}}$: si $S=N$, on a simplement $b_{\underline{l}+1_{n}}=0$ et si $S<N$, c'est l'hypoth\`ese de r\'ecurrence. Donc la premi\`ere somme dispara\^{\i}t. Pour $k\not=n$, il appara\^{\i}t dans le deuxi\`eme somme le terme $\partial_{H_{n}}^{N-l_{n}}\partial_{H_{k}}b_{\underline{l}+1_{n}-1_{k}}$. C'est nul si $l_{k}=0$. Supposons $l_{k}\not=0$. Alors l'hypoth\`ese de r\'ecurrence pour $l_{n}+1$  assure que le terme est nul. Il reste seulement le terme de la deuxi\`eme somme index\'e par $k=n$, qui est donc nul lui-aussi. Or ce terme est $\partial_{H_{n}}^{N-l_{n}+1}b_{\underline{l}}$. Cela prouve (8). Evidemment, (8) implique (6). 

Pour $i,j\in I$ et $\underline{l}\in {\mathbb N}^m$, la fonction $\gamma_{I}\mapsto a_{i,j;\underline{l}}(\gamma_{I},z)$ est rationnelle. Pour $\gamma_{I}$ fix\'e, la fonction $H\mapsto a_{i,j,\underline{l}}(exp(H)\gamma_{I};z)$ est donc rationnelle en $exp(H)$. Or (6) montre que c'est un polyn\^ome en $H$. 
  Un r\'esultat \'el\'ementaire sur les polyn\^omes exponentiels montre que ces deux propri\'et\'es ne sont v\'erifi\'ees que pour les constantes. On obtient que $\gamma_{I}\mapsto a_{i,j;\underline{l}}(\gamma_{I})$ est localement constant. Puisque c'est une fraction rationnelle, elle est vraiment constante. On utilise maintenant le (ii) de la proposition 1.3, qui se propage lui-aussi aux variantes de nos op\'erateurs diff\'erentiels: quand $\gamma_{I}$ tend vers l'infini au sens expliqu\'e en 1.3, $a_{i,j;\underline{l}}(\gamma_{I})$ tend vers $0$. Puisque c'est une constante, cette constante est nulle. On conclut que $A_{i,j}(\gamma_{I},z)=0$. C'est ce qu'il fallait d\'emontrer. $\square$

\bigskip

\section{Majorations}

\bigskip

\subsection{Quelques consid\'erations formelles}
Soient $(G,\tilde{G},{\bf a})$ un triplet comme en 1.1,  $\tilde{T}$ un sous-tore tordu maximal de $\tilde{G}$ et $\eta$ un \'el\'ement de $\tilde{T}({\mathbb R})$. On suppose $\omega$ trivial sur $T({\mathbb R})^{\theta}$. 

Ainsi qu'on l'a dit en 1.1, on a des isomorphismes 
$$Sym(\mathfrak{t}^{\theta})\simeq Sym(\mathfrak{t})/I_{\theta,\omega}=Sym(\mathfrak{t})_{\theta,\omega}\simeq Diff(\tilde{T}({\mathbb R}))^{\omega-inv}.$$ 

 Soit ${\bf G}'=(G',{\cal G}',\tilde{s})$ une donn\'ee endoscopique relevante de $(G,\tilde{G},{\bf a})$. Supposons donn\'e un diagramme $(\epsilon,B',T',B,T,\eta)$, o\`u $\epsilon$ est un \'el\'ement semi-simple de $\tilde{G}'({\mathbb R})$, cf. [I] 1.10. Fixons des donn\'ees auxiliaires $G'_{1}$,...,$\Delta_{1}$ pour ${\bf G}'$ et un point $\epsilon_{1}\in G'_{1}({\mathbb R})$ au-dessus de $\epsilon$. Posons $\tilde{T}'=T'\epsilon$  et notons $\tilde{T}'_{1}$ l'image r\'eciproque de $\tilde{T}' $ dans $\tilde{G}'_{1}$. Fixons une d\'ecomposition d'alg\`ebres de Lie
 $$(1) \qquad \mathfrak{g}'_{1}=\mathfrak{c}_{1}\oplus \mathfrak{g}'.$$
 Remarquons que $G'_{1,SC}=G'_{SC}$ donc la section $\mathfrak{g}'_{SC}\to \mathfrak{g}'_{1}$ est uniquement d\'etermin\'ee. Choisir une d\'ecomposition comme ci-dessus \'equivaut \`a choisir une d\'ecomposition 
 $$\mathfrak{z}(G'_{1})=\mathfrak{c}_{1}\oplus \mathfrak{z}(G').$$
 Il se d\'eduit de (1) une d\'ecomposition
$$(2) \qquad \mathfrak{t}'_{1}=\mathfrak{c}_{1}\oplus \mathfrak{t}'.$$
  Le diagramme fournit un isomorphisme
$$\xi:\mathfrak{t}^{\theta}\to \mathfrak{t}'.$$
Consid\'erons la fonction
$$X\mapsto \Delta_{1}(exp(\xi(X))\epsilon,exp(X)\eta)$$
d\'efinie presque partout dans un voisinage de $0$ dans $\mathfrak{t}^{\theta}({\mathbb R})$ (dans le premier terme, $\xi(X)$ est identifi\'e \`a un \'el\'ement de $\mathfrak{t}'_{1}({\mathbb R})$ gr\^ace \`a la d\'ecomposition (2)). Pr\'ecis\'ement, cette fonction est d\'efinie en tout $X$ proche de $0$ tel que $exp(X)\eta$ soit fortement r\'egulier. Le lemme [I] 2.8 montre qu'elle est localement proportionnelle \`a la fonction $X\mapsto e^{<b,\xi(X)\oplus X>}$, avec les notations de ce lemme. D'apr\`es [I] 2.8(1), $b$ est un \'el\'ement de $\mathfrak{z}(G'_{1})^*\oplus (1-\theta)(\mathfrak{t}^*)$. Puisque $X$ est fixe par $\theta$, son produit avec la deuxi\`eme composante de $b$ est nul. En notant $b_{1}$ la premi\`ere composante, on obtient

(3) $ \Delta_{1}(exp(\xi(X))\epsilon,exp(X)\eta)$ est le produit de $e^{<b_{1},\xi(X)>}$ et d'une fonction qui est localement constante sur les \'el\'ements $X\in \mathfrak{t}^{\theta}({\mathbb R})$ proches de $0$ et tels que $exp(X)\eta$ soit fortement r\'egulier. 

Remarquons que, puisque $b_{1}\in \mathfrak{z}(G'_{1})^*$, on a $<b_{1},\xi(X)>=0$ si $\xi(X)\in \mathfrak{t}'_{sc}$, en notant $T'_{sc}$ l'image r\'eciproque de $T'$ dans $G'_{SC}$.

De (2) se d\'eduit un isomorphisme
 $Sym(\mathfrak{t}')\to Sym(\mathfrak{t}'_{1})_{\lambda_{1}}$. D'apr\`es ce que l'on a dit en 1.8, on a donc des isomorphismes
$$ Sym(\mathfrak{t}^{\theta})\simeq  Sym(\mathfrak{t})_{\theta,\omega}\simeq Diff(\tilde{T}({\mathbb R}))^{\omega-inv}\simeq Diff(\tilde{T}'_{1}({\mathbb R}))_{\lambda_{1}}\simeq Sym(\mathfrak{t}'_{1})_{\lambda_{1}}\simeq  Sym(\mathfrak{t}').$$
On note $\Xi_{1}:Sym(\mathfrak{t}^{\theta})\to Sym(\mathfrak{t}')$ l'isomorphisme compos\'e. 
  On prendra garde que l'isomorphisme $\Xi_{1}$ n'est pas en g\'en\'eral celui d\'eduit de $\xi$.  Avec les notations ci-dessus, l'isomorphisme $\Xi_{1}$ envoie un \'el\'ement $X\in \mathfrak{t}^{\theta}$ sur l'\'el\'ement $\xi(X) -<b_{1},\xi(X)>$. En particulier, $\Xi_{1}$ d\'epend de la d\'ecomposition (1). 

Consid\'erons d'autres donn\'ees auxiliaires $G'_{2}$, etc... et une d\'ecomposition 
$$ \mathfrak{t}'_{2}=\mathfrak{c}_{2}\oplus \mathfrak{t}'$$
d\'efinie de la m\^eme fa\c{c}on que ci-dessus. 
On se souvient de l'application de transition $\tilde{\lambda}_{1,2}$ d\'efinie en [I] 2.5. On a rappel\'e en 1.8 l'isomorphisme $C_{c,\lambda_{1}}^{\infty}(\tilde{G}'_{1}({\mathbb R}))\simeq C_{c,\lambda_{2}}^{\infty}(\tilde{G}'_{2}({\mathbb R}))$.   On d\'eduit de $\tilde{\lambda}_{12}$ une fonction $\tau_{12}$ sur $\mathfrak{t}'({\mathbb R})$ par $\tau_{12}(X)=\tilde{\lambda}_{12}(exp(X)\epsilon_{1},exp(X)\epsilon_{2})^{-1}$. La fonction $\tilde{\lambda}_{12}$ est $C^{\infty}$, il en est donc de m\^eme de la fonction $\tau_{12}$. Par d\'efinition de $\tilde{\lambda}_{12}$, on a l'\'egalit\'e
$$\tau_{12}(X)=\Delta_{1}(exp(X)\epsilon_{1},exp(\xi^{-1}(X))\eta)\Delta_{2}(exp(X)\epsilon_{2},exp(\xi^{-1}(X))\eta)^{-1}$$
si $exp(\xi^{-1}(X))\eta$ est fortement r\'egulier. Il r\'esulte de ce qui pr\'ec\`ede qu'il existe un \'el\'ement $b_{12}\in \mathfrak{z}(G')^*$ tel que $\tau_{12}(X)$ soit le produit d'une constante non nulle avec la fonction $X\mapsto e^{<b_{12},X>}$ ($b_{12}$ est la diff\'erence entre les projections de $b_{1}$ et $b_{2}$ dans $\mathfrak{z}(G')^*$).  En particulier, $\tau_{12}$ prend la m\^eme valeur sur des points stablement conjugu\'es. Fixons une mesure de Haar sur $T'({\mathbb R})$. Pour $X$ r\'egulier, l'application duale de la pr\'ec\'edente envoie l'int\'egrale orbitale associ\'ee \`a $exp(X)\epsilon_{1}$ et \`a cette mesure sur l'int\'egrale orbitale associ\'ee \`a $exp(X)\epsilon_{2}$ et \`a cette m\^eme mesure, multipli\'ee par $\tau_{12}(X)$. Pour $H\in Sym(\mathfrak{t}^{\theta})$, notons $H_{1}=\Xi_{1}(H)$ et $H_{2}=\Xi_{2}(H)$. Il se d\'eduit de $H_{1}$ et $H_{2}$ des op\'erateurs diff\'erentiels $\partial_{H_{1}}$ et $\partial_{H_{2}}$ sur $\mathfrak{t}'({\mathbb R})$. On a l'\'egalit\'e
$$\partial_{H_{2}}=\tau_{12}^{-1}\circ \partial_{H_{1}}\circ \tau_{12}.$$

\bigskip

\subsection{Majoration des int\'egrales orbitales pond\'er\'ees $\omega$-\'equivariantes}
Soient   $\tilde{M}$ un espace de Levi de $\tilde{G}$,  $\tilde{T}$ un sous-tore tordu maximal de $\tilde{M}$ et $\eta$ un \'el\'ement de $\tilde{T}({\mathbb R})$. On suppose $\omega$ trivial sur $T({\mathbb R})^{\theta}$. Fixons des voisinages  $\mathfrak{u}$ et $\mathfrak{u}'$ de $0$ dans $\mathfrak{t}^{\theta,0}({\mathbb R})$ v\'erifiant les conditions suivantes

(1) l'application $X\mapsto exp(X)\eta(1-\theta)(T({\mathbb R})$ est un isomorphisme de $\mathfrak{u}'$ sur un voisinage de l'image de $\eta$ dans $\tilde{T}({\mathbb R})/(1-\theta)(T({\mathbb R})$;

(2) $\mathfrak{u}$ et $\mathfrak{u}'$ sont ouverts; la cl\^oture de $\mathfrak{u}$ est compacte et est contenue dans $\mathfrak{u}'$;

(3)   pour $X\in \mathfrak{u}'$, $X$ est r\'egulier dans $\mathfrak{g}_{\eta}({\mathbb R})$ si et seulement si $exp(X)\eta$ est  r\'egulier dans $\tilde{G}({\mathbb R})$. 

{\bf Remarque.} Le groupe $Z_{G}(\eta)$ agit naturellement sur $\mathfrak{g}_{\eta}$. A la place de (3), on pourrait \'enoncer

(3)'  pour $X\in \mathfrak{u}'$, le stabilisateur de $X$ dans $Z_{G}(\eta)$ est \'egal \`a $T^{\theta}$ si et seulement si $exp(X)\eta$ est fortement  r\'egulier dans $\tilde{G}({\mathbb R})$. 
\bigskip

 On fixe des mesures de Haar sur $G({\mathbb R})$, $M({\mathbb R})$ et  $T^{\theta,0}({\mathbb R})$ qui permettent de d\'efinir les int\'egrales orbitales qui suivent. On a

(4) pour tout $U\in Sym(\mathfrak{t})_{\theta,\omega}$, il existe un entier $N$ et, pour tout $f\in C_{c}^{\infty}(\tilde{G}({\mathbb R}))$, il existe $c>0$ de sorte que l'on ait la majoration
$$\vert \partial_{U}I_{\tilde{M}}^{\tilde{G}}(exp(X)\eta,\omega,f)\vert \leq c D^{G_{\eta}}(X)^{-N}$$
pour tout $X\in \mathfrak{u} $ tel que $exp(X)\eta$ soit fortement r\'egulier dans $\tilde{G}({\mathbb R})$.

Cf. [A1] lemme 13.2.

\bigskip

\subsection{Majoration des int\'egrales orbitales pond\'er\'ees stables}
On conserve la m\^eme situation mais on suppose $(G,\tilde{G},{\bf a})$ quasi-d\'eploy\'e et \`a torsion int\'erieure.  On a

\ass{Lemme}{ Pour tout $U\in Sym(\mathfrak{t}) $, il existe un entier $N$ et, pour tout $f\in C_{c}^{\infty}(\tilde{G}({\mathbb R}))$, il existe $c>0$ de sorte que l'on ait la majoration
$$\vert \partial_{U}S_{\tilde{M}}^{\tilde{G}}(exp(X)\eta,f)\vert \leq cD^{G_{\eta}}(X)^{-N}$$
pour tout $X\in \mathfrak{u} $ tel que $exp(X)\eta$ soit fortement r\'egulier dans $\tilde{G}({\mathbb R})$.}

Preuve. Pour tout $s\in Z(\hat{M})^{\Gamma_{{\mathbb R}}}/Z(\hat{G})^{\Gamma_{{\mathbb R}}}$ avec $s\not=1$, on fixe des donn\'ees auxiliaires $G_{1}'(s)$,...,$\Delta_{1}(s)$. On note $\tilde{M}_{1}(s)$ l'image r\'eciproque de $\tilde{M}$ dans $\tilde{G}'_{1}(s)$ et $\tilde{T}_{1}(s)$ celle de $\tilde{T}$. On fixe  un point $\eta_{1}(s)\in \tilde{T}_{1}(s;{\mathbb R})$ au-dessus de $\eta$ et une d\'ecomposition
$$\mathfrak{t}_{1}(s)=\mathfrak{c}_{1}(s)\oplus \mathfrak{t}$$
selon la recette de 3.1.
Comme on l'a dit dans ce paragraphe (et en adaptant les notations de fa\c{c}on compr\'ehensible), l'int\'egrale orbitale sur $\tilde{M}({\mathbb R})$ associ\'ee \`a $exp(X)\eta$ s'identifie \`a celle sur $\tilde{M}_{1}(s;{\mathbb R})$ associ\'ee \`a $exp(X)\eta_{1}(s)$, multipli\'ee par une fonction  que l'on note $\tau(s;X)$. L'\'el\'ement $U$ se transf\`ere en un element  que l'on note $\Xi_{1}(s;U)\in Sym(\mathfrak{t})$. On identifie le transfert $f^{{\bf G}'(s)}$ \`a  une fonction $f^{\tilde{G}'_{1}(s)}$ sur $\tilde{G}'_{1}(s;{\mathbb R})$. On introduit le terme $\Lambda_{\tilde{M}}^{\tilde{G}}(exp(X)\eta,f)$ de 2.1. C'est la somme des $I_{\tilde{M}}^{\tilde{G}}(exp(X')\eta',f)$ quand $exp(X')\eta'$ d\'ecrit un ensemble de repr\'esentants des classes de conjugaison par $M({\mathbb R})$ dans la classe de conjugaison stable de $exp(X)\eta$ dans $\tilde{M}({\mathbb R})$. 
On  a alors l'\'egalit\'e
$$(1) \qquad \partial_{U}S_{\tilde{M}}^{\tilde{G}}(exp(X)\eta,f)=\partial_{U}\Lambda_{\tilde{M}}^{\tilde{G}}(exp(X)\eta,f)$$
$$-\sum_{s\in Z(\hat{M})^{\Gamma_{{\mathbb R}}}/Z(\hat{G})^{\Gamma_{{\mathbb R}}},s\not=1} i_{\tilde{M}}(\tilde{G},\tilde{G}'(s))\tau(s;X)\partial_{\Xi_{1}(s;U)}S_{\tilde{M}_{1}(s),\lambda_{1}(s)}^{\tilde{G}_{1}'(s)}(exp(X)\eta_{1}(s),f^{\tilde{G}'_{1}(s)}).$$
Il est clair que le fait que $f^{\tilde{G}'_{1}(s)}$ ne soit pas \`a support compact mais se transforme selon le caract\`ere $\lambda_{1}(s)$ de $C_{1}(s;{\mathbb R})$ ne change rien aux propri\'et\'es locales de ses int\'egrales orbitales stables. Pour $s\not=1$, on peut donc appliquer par r\'ecurrence la majoration du lemme au terme $\partial_{\Xi_{1}(s;U)}S_{\tilde{M}_{1}(s),\lambda_{1}(s)}^{\tilde{G}_{1}'(s)}(exp(X)\eta_{1}(s),f^{\tilde{G}'_{1}(s)})$. On obtient que celui-ci est major\'e par une puissance n\'egative de $D^{G_{1}'(s)_{\eta_{1}(s)}}(X)$. Ce terme est \'egal \`a $D^{G'(s)_{\eta}}(X)$. Or le syst\`eme de racines de $G'(s)_{\eta}$ est inclus dans celui de $G_{\eta}$. Dans un voisinage compact de $0$, $D^{G'(s)_{\eta}}(X)^{-1}$ est donc major\'e par $D^{G_{\eta}}(X)^{-1}$.
Les m\^emes arguments qu'en 2.1 montrent que  la majoration 3.2(4) pour $\partial_{U}I_{\tilde{M}}^{\tilde{G}}(exp(X)\eta,f)$ s'\'etend  au terme $\partial_{U}\Lambda_{\tilde{M}}^{\tilde{G}}(exp(X)\eta,f)$. Tous les termes du membre de droite de (1) sont donc major\'es par une puissance n\'egative de $D^{G_{\eta}}(X)$. L'assertion de l'\'enonc\'e s'ensuit. $\square$

\bigskip

\subsection{Majoration des int\'egrales orbitales endoscopiques}
On revient \`a la situation de 3.2 o\`u $(G,\tilde{G},{\bf a})$ est quelconque. On suppose que $\tilde{G}$ est une composante connexe d'un $K$-espace $K\tilde{G}$. Alors $\tilde{M}$ est une composante connexe d'un $K$-espace de Levi $K\tilde{M}$ et on suppose $K\tilde{M}\in {\cal L}(K\tilde{M}_{0})$. Pour $f\in C_{c}^{\infty}(K\tilde{G}({\mathbb R}))$ et $X\in \mathfrak{u}$ tel que $exp(X)\eta$ soit fortement r\'egulier, on sait d\'efinir l'int\'egrale orbitale endoscopique $I^{K\tilde{G},{\cal E}}_{K\tilde{M}}(exp(X)\eta,\omega,f)$. 

\ass{Lemme}{Pour tout $U\in Sym(\mathfrak{t}) $, il existe un entier $N$ et, pour tout $f\in C_{c}^{\infty}(K\tilde{G}({\mathbb R}))$, il existe $c>0$ de sorte que l'on ait la majoration
$$\vert \partial_{U}I_{K\tilde{M}}^{K\tilde{G},{\cal E}}(exp(X)\eta,f)\vert \leq c D^{G_{\eta}}(X)^{-N}$$
pour tout $X\in \mathfrak{u} $ tel que $exp(X)\eta$ soit fortement r\'egulier dans $\tilde{G}({\mathbb R})$.}

Preuve. On peut fixer les objets suivants:

- une famille finie $({\bf M}'_{i})_{i=1,...,n}$ de donn\'ees endoscopiques elliptiques et relevantes de $(M,\tilde{M},{\bf a})$, munies de donn\'ees auxiliaires $M'_{i,1}$,...,$\Delta_{i,1}$;

- pour tout $i=1,...,n$, un diagramme $(\epsilon_{i},B^{M'_{i}},T^{M'_{i}},B_{i}^M,T,\eta)$ joignant $\eta$ \`a un \'el\'ement semi-simple $\epsilon_{i}\in \tilde{M}'_{i}({\mathbb R})$;  on note $\xi_{i}:\mathfrak{t}^{\theta}\to \mathfrak{t}^{M'_{i}}$ l'isomorphisme d\'eduit du diagramme;

- pour tout $i=1,...,n$, un rel\`evement $\epsilon_{i,1}\in \tilde{M}'_{i;1}({\mathbb R})$ de $\epsilon_{i}$ et, en notant $T^{M'_{i}}_{1}$ l'image r\'eciproque de $T^{M'_{i}}$ dans $M'_{i,1}$, une d\'ecomposition
$$\mathfrak{t}^{M'_{i}}_{1}=\mathfrak{c}_{i,1}\oplus \mathfrak{t}^{M'_{i}}$$
selon la recette de 3.1;

- un nombre $c_{\tilde{T}}>0$;

\noindent ces donn\'ees v\'erifiant la condition suivante. Soit $\varphi\in C_{c}^{\infty}(K\tilde{M}({\mathbb R}))$. Pour tout $i$, notons $\varphi^{\tilde{M}'_{i,1}}$ son transfert \`a $\tilde{M}'_{i,1}$. Alors on a l'\'egalit\'e
$$I^{K\tilde{M}}(exp(X)\eta,\omega,\varphi)=c_{\tilde{T}}\sum_{i=1,...,n}\Delta_{i,1}(exp(\xi_{i}(X)\epsilon_{i,1},exp(X)\eta)^{-1}S^{\tilde{M}'_{i,1}}_{\lambda_{i,1}}(exp(\xi_{i}(X))\epsilon_{i,1},\varphi^{\tilde{M}'_{i,1}}).$$
Pour tout $i=1,...,n$, on a d\'efini l'int\'egrale $I_{K\tilde{M}}^{K\tilde{G},{\cal E}}({\bf M}'_{i},\boldsymbol{\delta},f)$ pour tout $\boldsymbol{\delta}\in D^{st}_{g\acute{e}om,\tilde{G}-\acute{e}qui}({\bf M}')$. On note simplement  $I_{K\tilde{M}}^{K\tilde{G},{\cal E}}({\bf M}'_{i},exp(\xi_{i}(X))\epsilon_{i,1},f)$ ce terme quand $\boldsymbol{\delta}$ est l'\'el\'ement de $D^{st}_{g\acute{e}om,\tilde{G}-\acute{e}qui}({\bf M}')$  auquel s'identifie l'int\'egrale orbitale stable sur $\tilde{M}'_{i,1}({\mathbb R})$ associ\'ee au point $exp(\xi_{i}(X))\epsilon_{i,1}$. 
 Il r\'esulte des d\'efinitions que l'on a aussi
$$I_{K\tilde{M}}^{K\tilde{G},{\cal E}}(exp(X)\eta,f) =c_{\tilde{T}}\sum_{i=1,...,n}\Delta_{i,1}(exp(\xi_{i}(X)\epsilon_{i,1},exp(X)\eta)^{-1}I_{K\tilde{M}}^{K\tilde{G},{\cal E}}({\bf M}'_{i},exp(\xi_{i}(X))\epsilon_{i,1},f).$$
Pour tout $i$, on a un homomorphisme $\Xi_{i}:Sym(\mathfrak{t}^{\theta})\to Sym(\mathfrak{t}^{M'_{i}})$ comme en 3.1. On a alors 
$$\partial_{U}I_{K\tilde{M}}^{K\tilde{G},{\cal E}}(exp(X)\eta,f) =c_{\tilde{T}}\sum_{i=1,...,n}\Delta_{i,1}(exp(\xi_{i}(X)\epsilon_{i,1},exp(X)\eta)^{-1}$$
$$\partial_{\Xi_{i}(U)}I_{K\tilde{M}}^{K\tilde{G},{\cal E}}({\bf M}'_{i},exp(\xi_{i}(X))\epsilon_{i,1},f).$$
Pr\'ecisons que $\partial_{\Xi_{i}(U)}I_{K\tilde{M}}^{K\tilde{G},{\cal E}}({\bf M}'_{i},exp(\xi_{i}(X))\epsilon_{i,1},f)$ est la valeur en $\xi_{i}(X)$ de la fonction $Y\mapsto \partial_{\Xi_{i}(U)}I_{K\tilde{M}}^{K\tilde{G},{\cal E}}({\bf M}'_{i},exp(Y)\epsilon_{i,1},f)$ sur $\mathfrak{t}^{M'_{i}}({\mathbb R})$. 
La formule ci-dessus nous permet de fixer $i$ et de majorer $\partial_{U'}I_{K\tilde{M}}^{K\tilde{G},{\cal E}}({\bf M}'_{i},exp(\xi_{i}(X)))\epsilon_{i,1},f)$ pour $U'\in Sym(\mathfrak{t}^{M'_{i}})$. Puisque $i$ est fix\'e, on le fait dispara\^{\i}tre de la notation et on pose  ${\bf M}' =(M',{\cal M}',\tilde{\zeta})$. Pour tout $\tilde{s}\in Z(\hat{M})^{\Gamma_{{\mathbb R}},\hat{\theta}}/Z(\hat{G})^{\Gamma_{{\mathbb R}},\hat{\theta}}$, on introduit des donn\'ees auxiliaires $G'_{1}(\tilde{s})$,...,$\Delta_{1}(\tilde{s})$. On introduit les images r\'eciproques $\tilde{M}'_{1}(\tilde{s})$ et $\tilde{T}^{M'}_{1}(\tilde{s})$ de $\tilde{M}'$ et $\tilde{T}^{M'}$ dans $\tilde{G}'_{1}(\tilde{s})$. On  fixe un point $\epsilon_{1}(\tilde{s})\in \tilde{T}^{M'}_{1}(\tilde{s};{\mathbb R})$ au-dessus de $\epsilon$ et une d\'ecomposition
$$\mathfrak{t}^{M'}_{1}(\tilde{s})=\mathfrak{c}_{1}(\tilde{s})\oplus \mathfrak{t}^{M'}$$
selon la recette de 3.1.
D'apr\`es 3.1, pour $Y\in \mathfrak{t}^{M'}({\mathbb R})$,  l'int\'egrale orbitale dans $\tilde{M}'_{1}({\mathbb R})$ associ\'ee \`a $exp(Y)\epsilon_{1}$ s'identifie \`a l'int\'egrale orbitale dans $\tilde{M}'_{1}(\tilde{s};{\mathbb R})$ associ\'ee \`a $exp(Y)\epsilon_{1}(\tilde{s})$, multipli\'ee par une fonction $\tau(\tilde{s},Y)$. Un \'el\'ement $U'\in Sym(\mathfrak{t}^{M'})$ se transf\`ere en un \'el\'ement $\Xi(\tilde{s},U')$ de la m\^eme alg\`ebre. On a alors
$$\partial_{U'}I_{K\tilde{M}}^{K\tilde{G},{\cal E}}({\bf M}',exp(Y)\epsilon_{1},f)=\sum_{\tilde{s}\in Z(\hat{M})^{\Gamma_{{\mathbb R}},\hat{\theta}}/Z(\hat{G})^{\Gamma_{{\mathbb R}},\hat{\theta}} }i_{\tilde{M}'}(\tilde{G},\tilde{G}'(\tilde{s}))$$
$$\tau(\tilde{s},Y)\partial_{\Xi(\tilde{s},U')}S_{\tilde{M}'_{1}(\tilde{s}),\lambda_{1}(\tilde{s})}^{\tilde{G}'_{1}(\tilde{s})}(exp(Y)\epsilon_{1}(\tilde{s}),f^{\tilde{G}'_{1}(\tilde{s})}).$$
D'apr\`es le lemme 3.3, le terme index\'e par $\tilde{s}$ est major\'e par une puissance n\'egative de $D^{G'(\tilde{s})_{\epsilon}}(Y)$. On sait que, par l'homomorphisme $\xi:\mathfrak{t}^{\theta}\to \mathfrak{t}^{M'}$, les racines du groupe $G'(\tilde{s})_{\epsilon}$ s'identifie \`a des multiples rationnels de racines du groupe $G_{\eta}$. Il en r\'esulte que $D^{G'(\tilde{s})_{\epsilon}}(\xi(X))^{-1}$ est major\'e par $D^{G_{\eta}}(X)^{-1}$ dans un voisinage de $0$. On en d\'eduit la majoration cherch\'ee pour $\partial_{U'}I_{K\tilde{M}}^{K\tilde{G},{\cal E}}({\bf M}'_{i},exp(\xi_{i}(X))\epsilon_{i,1},f)$. $\square$

\bigskip

\section{Propri\'et\'es locales}

\bigskip

\subsection{Sauts des int\'egrales orbitales pond\'er\'ees $\omega$-\'equivariantes}
Soient   $\tilde{M}$ un espace de Levi de $\tilde{G}$,  $\tilde{T}$ un sous-tore tordu maximal de $\tilde{M}$ et $\eta$ un \'el\'ement de $\tilde{T}({\mathbb R})$. On suppose $\omega$ trivial sur $T({\mathbb R})^{\theta}$. On impose dans toute la section les hypoth\`eses (1), (2) et (3) ci-dessous:

(1) $\tilde{T}$ est un sous-tore tordu elliptique de $\tilde{M}$;

(2) $G_{\eta,SC}$ est isomorphe \`a $SL(2)$;

(3) l'image r\'eciproque $T_{d}$  de $T^{\theta,0}$ dans $G_{\eta,SC}$ est un tore d\'eploy\'e.

Fixons un sous-tore maximal $T_{c}$   non d\'eploy\'e de $G_{\eta,SC}$. Il existe un unique sous-tore maximal   de $G_{\eta}$ dont l'image r\'eciproque dans $G_{\eta,SC}$ soit $T_{c}$. On note $\underline{T}$ son commutant dans $G$ et $\underline{\tilde{T}}=\underline{T}\eta$. L'ensemble $\underline{\tilde{T}}$ est un sous-tore tordu maximal de $\tilde{G}$.    Notons $\underline{\tilde{M}}$ le plus petit espace de Levi dans lequel $\underline{\tilde{T}}$ est contenu, c'est-\`a-dire le commutant de $A_{\underline{\tilde{T}}}$ dans $\tilde{G}$. On a

(4) $M_{\eta}=T^{\theta,0}$, $\underline{M}_{\eta}=G_{\eta}$, $A_{\underline{\tilde{M}}}=A_{G_{\eta}}$, $\tilde{M}$ est un sous-espace de Levi propre et maximal de $\underline{\tilde{M}}$ et $\underline{\tilde{T}}$ est un sous-tore tordu elliptique de $\underline{\tilde{M}}$.

Preuve.  Pour un sous-tore $T'$ de $G_{\eta}$ qui est d\'eploy\'e, on a deux possibilit\'es.  Ou bien son image r\'eciproque dans $G_{\eta,SC}$ est triviale. Dans ce cas $T'$ est contenu dans $Z(G_{\eta})$ et donc (puisque c'est un tore d\'eploy\'e) dans $A_{G_{\eta}}$. Ou bien l'image r\'eciproque de $T'$ dans $G_{\eta,SC}$ est de rang $1$ et poss\`ede au signe pr\`es une seule racine $\alpha$. La composante neutre du noyau de cette racine est alors un sous-tore $T''$ de $T'$ qui est du premier cas, donc contenu dans $A_{G_{\eta}}$. Et on a $dim(T'')=dim(T')-1$.  Puisque $T^{\theta,0}$ et $\underline{T}^{\theta,0}$ sont des sous-tores maximaux de $G_{\eta}$, on a les inclusions $A_{G_{\eta}}\subset A_{T^{\theta,0}}$, $A_{G_{\eta}}\subset A_{\underline{T}^{\theta,0}}$. Puisque les images r\'eciproques de $T^{\theta,0}$ et $\underline{T}^{\theta,0}$ dans $G_{\eta,SC}$ sont respectivement d\'eploy\'ees et non d\'eploy\'ees, la discussion pr\'ec\'edente montre que $A_{G_{\eta}}=A_{\underline{T}^{\theta,0}}$, tandis que $dim(A_{T^{\theta,0}})=dim(A_{G_{\eta}})+1$. La premi\`ere \'egalit\'e entra\^{\i}ne que le commutant  $\underline{M}$ de $A_{\underline{T}^{\theta,0}}$ dans $G$ contient $G_{\eta}$. A fortiori $\underline{M}_{\eta}=G_{\eta}$. On a aussi $A_{\underline{T}^{\theta,0}}=A_{G_{\eta}}\subset A_{T^{\theta,0}}$. En prenant les commutants dans $\tilde{G}$ et en se rappelant que le commutant de $A_{T^{\theta,0}}$ est $\tilde{M}$ d'apr\`es (1), on obtient que $\tilde{M}\subset \underline{\tilde{M}}$.  Puisque l'image r\'eciproque de $A_{T^{\theta,0}}$ dans $G_{\eta,SC}$ est $T_{d}$, $A_{T^{\theta,0}}$ ne commute pas \`a $G_{\eta}$. Donc ce groupe  $G_{\eta}$ n'est pas contenu dans $M$. Donc l'inclusion $M_{\eta}\subset G_{\eta}$ est stricte. Puisqu'il s'agit d'une inclusion d'un Levi dans un groupe de rang semi-simple 1, $M_{\eta}$ est un tore, donc \'egal \`a $T^{\theta,0}$. Puisque $\underline{M}$ contient $G_{\eta}$ tandis que $M$ ne le contient pas, l'inclusion $\tilde{M}\subset \underline{\tilde{M}}$ est stricte. Les inclusions
$$A_{G_{\eta}}\subset A_{\underline{\tilde{M}}}\subset A_{\tilde{M}}=A_{T^{\theta,0}}$$
jointes aux faits que la deuxi\`eme inclusion est stricte et que la diff\'erence des dimensions entre le tore de droite et celui de gauche est 1 entra\^{\i}ne que la premi\`ere inclusion est une \'egalit\'e et que $\tilde{M}$ est maximal parmi les espaces de Levi propres de $\underline{\tilde{M}}$. Enfin, que $\underline{\tilde{T}}$ soit elliptique dans $\underline{\tilde{M}}$ r\'esulte simplement de la d\'efinition de cet espace de Levi. $\square$

 On fixe l'une des racines de $T^{\theta,0}$ dans $G_{\eta}$, que l'on note $\alpha$. On introduit la coracine associ\'ee, que l'on peut voir comme un \'el\'ement  $\check{\alpha}\in X_{*}(T^{\theta,0})$.  La preuve de (4) montre que l'on a m\^eme $\check{\alpha}\in X_{*}(A_{\tilde{M}})$. On se rappelle que l'on a fix\'e une forme quadratique d\'efinie positive sur ${\cal A}_{\tilde{M}}=X_{*}(A_{\tilde{M}})\otimes_{{\mathbb Z}}{\mathbb R}$, qui nous sert \`a d\'efinir la mesure sur cet espace. On dispose donc de la norme $\vert \check{\alpha}\vert $ pour cette forme. Cela \'etant, on d\'efinit au voisinage $0$ dans $\mathfrak{t}^{\theta}({\mathbb R})$ une modification de l'int\'egrale orbitale pond\'er\'ee $\omega$-\'equivariante de la fa\c{c}on suivante. Pour $f\in C_{c}^{\infty}(\tilde{G}({\mathbb R}))$ et $X\in \mathfrak{t}^{\theta}({\mathbb R}) $ assez proche de $0$, on pose
 $$I_{\tilde{M}}^{\tilde{G},mod}(exp(X)\eta,\omega,f)=I_{\tilde{M}}^{\tilde{G}}(exp(X)\eta,\omega,f)+\vert \check{\alpha}\vert log(\vert \alpha(X)\vert )I_{\underline{\tilde{M}}}^{\tilde{G}}(exp(X)\eta,\omega,f).$$

 On fixe un \'el\'ement $u\in G_{\eta,SC}({\mathbb C})$ tel que $ad_{u}(T_{d})=T_{c}$. Fixons un \'el\'ement non nul $H_{d}\in \mathfrak{t}_{d}({\mathbb R})$. On v\'erifie que $ad_{u}(H_{d})$ est un \'el\'ement de $i\mathfrak{t}_{c}({\mathbb R})$ et on introduit l'\'el\'ement $H_{c}$ de $\mathfrak{t}_{c}({\mathbb R})$ tel que $ad_{u}(H_{d})=iH_{c}$. De l'application $ad_{u}$ se d\'eduit un isomorphisme $C:Sym(\mathfrak{t}^{\theta,0})\to Sym(\underline{\mathfrak{t}}^{\theta,0})$. On note $w_{d}$, resp. $w_{c}$, l'unique \'el\'ement non trivial du groupe de Weyl de $G_{\eta}$ relativement \`a $T^{\theta,0}$, resp. $\underline{T}^{\theta,0}$. L'isomorphisme $C$  entrelace l'action de $w_{d}$ sur le premier espace avec celle de $w_{c}$ sur le second. On a fix\'e une mesure sur $T^{\theta,0}({\mathbb R})$ pour d\'efinir nos int\'egrales orbitales. Elle se d\'eduit d'une forme diff\'erentielle de rang maximal sur $T^{\theta,0}$, autrement dit d'un \'el\'ement  $d\in \wedge^n\mathfrak{t}^{\theta,*}({\mathbb R})$ o\`u $n=dim(T^{\theta,0})$. De $ad_{u}$ se d\'eduit un isomorphisme $\wedge^n\mathfrak{t}^{\theta,*}\simeq \wedge^n\underline{\mathfrak{t}}^{\theta,*}$. On v\'erifie qu'il se restreint en un isomorphisme de $\wedge^n\mathfrak{t}^{\theta,*}({\mathbb R})$ sur $i\wedge^n\underline{\mathfrak{t}}^{\theta,*}({\mathbb R})$. On introduit l'\'el\'ement $\underline{d}\in \wedge^n\underline{\mathfrak{t}}^{\theta,*}({\mathbb R})$ tel que 
 $ad_{u}(d)=i\underline{d}$.  De $\underline{d}$ se d\'eduit alors une mesure sur $\underline{T}^{\theta,0}({\mathbb R})$ qui nous sert \`a d\'efinir les int\'egrales orbitales dans ce qui suit. 
 
 Il y a deux possibilit\'es pour la classe de conjugaison stable de $H_{c}$ dans $\mathfrak{g}_{\eta}({\mathbb R})$. Soit elle est r\'eduite \`a la classe de conjugaison de $H_{c}$ par $G_{\eta}({\mathbb R})$, auquel cas on pose $c(\eta)=0$. Soit elle se d\'ecompose en deux classes de conjugaison par $G_{\eta}({\mathbb R})$, celles de $H_{c}$ et de $-H_{c}$, auquel cas on pose $c(\eta)=1$.
 
 Selon l'usage, pour un voisinage $V$ de $0$ dans ${\mathbb R}$ et pour une fonction $\varphi$ d\'efinie  sur $V-\{0\}$, on note $lim_{r\to 0+}\varphi(r)$ et $lim_{r\to 0-}\varphi(r)$ les limites de $\varphi(r)$ quand $r$ tend vers $0$ en restant strictement positif, resp. n\'egatif (ces limites existent ou pas).
 
 \ass{Proposition}{Soient $f\in C_{c}^{\infty}(\tilde{G}({\mathbb R}))$.  
 
 (i) Pour tout  $U\in Sym(\mathfrak{t}^{\theta})$, les limites 
 
 $lim_{r\to 0+}\partial_{U}I_{\tilde{M}}^{\tilde{G},mod}(exp(rH_{d})\eta,\omega,f)$ et $lim_{r\to 0-}\partial_{U}I_{\tilde{M}}^{\tilde{G},mod}(exp(rH_{d})\eta,\omega,f)$
 
\noindent existent.
 
 (ii) Si $w_{d}(U)=U$, ces limites sont \'egales.
 
 (iii) Pour tout $\underline{U}\in Sym(\underline{\mathfrak{t}}^{\theta})$, les limites 
 
 $lim_{r\to 0+}\partial_{\underline{U}}I_{\underline{\tilde{M}}}^{\tilde{G}}(exp(rH_{c})\eta,\omega,f)$ et $lim_{r\to 0-}\partial_{\underline{U}}I_{\underline{\tilde{M}}}^{\tilde{G}}(exp(rH_{c})\eta,\omega,f)$ 
 
 \noindent existent.
 
 (iv) Soit $U\in Sym(\mathfrak{t}^{\theta})$, supposons $w_{d}(U)=-U$ et  posons $\underline{U}=C(U)$. Alors on a les \'egalit\'es
 $$lim_{r\to 0+}\partial_{U}I_{\tilde{M}}^{\tilde{G},mod}(exp(rH_{d})\eta,\omega,f)-lim_{r\to 0-}\partial_{U}I_{\tilde{M}}^{\tilde{G},mod}(exp(rH_{d})\eta,\omega,f)$$
 $$=2^{1+c(\eta)}\pi i\vert \check{\alpha}\vert \,lim_{r\to 0+}\partial_{\underline{U}}I_{\underline{\tilde{M}}}^{\tilde{G}}(exp(rH_{c})\eta,\omega,f)$$
 $$=-2^{1+c(\eta)}\pi i\vert \check{\alpha}\vert \,lim_{r\to 0-}\partial_{\underline{U}}I_{\underline{\tilde{M}}}^{\tilde{G}}(exp(rH_{c})\eta,\omega,f).$$}

 Cf. [A1] lemme 13.1.  La diff\'erence entre notre constante $2^{1+c(\eta)}\pi i\vert \check{\alpha}\vert $ et celle d'Arthur provient de nos normalisations diff\'erentes. 
 \bigskip
 
 \subsection{Sauts des int\'egrales orbitales pond\'er\'ees stables}
 
 On conserve les m\^emes hypoth\`eses que dans le paragraphe pr\'ec\'edent. On impose de plus que $(G,\tilde{G},{\bf a})$ est quasi-d\'eploy\'e et \`a torsion int\'erieure. Fixons  un sous-groupe de Borel de $G$ contenant $T$ et pour lequel $M$ et $\underline{M}$ sont standard. Compl\'etons $(B,T)$ en une paire de Borel \'epingl\'ee    ${\cal E}$ de $G$. Fixons une paire de Borel \'epingl\'ee  $\hat{{\cal E}}$ de $\hat{G}$  conserv\'ee par l'action galoisienne, dont on note le tore  $\hat{T}$ (l'action galoisienne sur ce tore n'en fait pas le tore dual de $T$).   On en d\'eduit des r\'ealisations de  $\hat{M}$ et $\underline{\hat{M}}$ comme des groupes de Levi standard de $\hat{G}$. De la coracine $\check{\alpha}$ se d\'eduit une racine $\hat{\alpha}$ de $\hat{T}$ qui se restreint en un caract\`ere non trivial de $Z(\hat{M})^{\Gamma_{{\mathbb R}}}$. Notons $Z(\hat{M})^{\hat{\alpha}}$ le groupe des $s\in Z(\hat{M})^{\Gamma_{{\mathbb R}}}$ tels que $\hat{\alpha}(s)=1$. Il r\'esulte des d\'efinitions que $Z(\underline{\hat{M}})^{\Gamma_{{\mathbb R}}}$ est inclus dans  ce groupe. Plus pr\'ecis\'ement, $Z(\underline{\hat{M}})^{\Gamma_{{\mathbb R}}}/Z(\hat{G})^{\Gamma_{{\mathbb R}}}$ est la composante neutre de $Z(\hat{M})^{\hat{\alpha}}/Z(\hat{G})^{\Gamma_{{\mathbb R}}}$. On pose
 $$i^{\tilde{G}}_{\tilde{M}}(\eta)=[Z(\hat{M})^{\hat{\alpha}}:Z(\underline{\hat{M}})^{\Gamma_{{\mathbb R}}}]^{-1}.$$
 
 Pour tout  \'el\'ement semi-simple $\epsilon\in \tilde{G}({\mathbb R})$, posons
 $${\cal Y} (\epsilon)=\{y\in G({\mathbb C}); \forall \sigma\in \Gamma_{{\mathbb R}}, y\sigma(y)^{-1}\in G_{\epsilon}\}.$$
Pour un element $y$ de cet ensemble, on pose $\epsilon[y]=y^{-1}\epsilon y$. On fixe un ensemble de repr\'esentants $\dot{{\cal Y}}(\epsilon)$ de l'ensemble de doubles classes
 $$G_{\epsilon}\backslash {\cal Y}(\epsilon)/G({\mathbb R}).$$
 Remarquons que l'application qui, \`a $y\in \dot{{\cal Y}}(\epsilon)$, associe la classe du cocycle 
 $\sigma\mapsto y\sigma(y)^{-1}$, est une bijection de $\dot{{\cal Y}}(\epsilon)$ sur le noyau de l'application
 $$H^1(\Gamma_{{\mathbb R}};G_{\epsilon})\to H^1(\Gamma_{{\mathbb R}};G).$$
 Soulignons que les ensembles ci-dessus ne sont pas des groupes. Le noyau est l'ensemble des \'el\'ements de $H^1(\Gamma_{{\mathbb R}};G_{\epsilon})$ qui deviennent des cobords dans $G$. 
 Si $\epsilon$ est fortement r\'egulier dans $\tilde{G}({\mathbb R})$, la famille $(\eta[y])_{y\in \dot{{\cal Y}}(\epsilon)}$ est un ensemble de repr\'esentants des classes de conjugaison par $G({\mathbb R})$ dans la classe de conjugaison stable de $\epsilon$. Dans ce cas, modulo le choix de mesures de Haar sur les groupes intervenant, on a d\'efini l'int\'egrale orbitale stable associ\'ee \`a $\epsilon$ comme la forme lin\'eaire qui, \`a $f\in C_{c}^{\infty}(\tilde{G}({\mathbb R}))$, associe
 $$S^{\tilde{G}}(\epsilon,f)=\sum_{y\in \dot{{\cal Y}}(\epsilon)}I^{\tilde{G}}(\epsilon[y],f).$$
 Si $\epsilon$ est r\'egulier mais pas fortement r\'egulier,  la famille $(\eta[y])_{y\in \dot{{\cal Y}}(\epsilon)}$ n' est plus, en g\'en\'eral, un ensemble de repr\'esentants des classes de conjugaison par $G({\mathbb R})$ dans la classe de conjugaison stable de $\epsilon$. On d\'efinit n\'eanmoins $S^{\tilde{G}}(\epsilon,f)$ par la formule ci-dessus. Pour $Y\in \mathfrak{g}_{\epsilon}({\mathbb R})$ voisin de $0$ et en position g\'en\'erale, l'\'el\'ement $exp(Y)\epsilon$ est fortement r\'egulier et on  voit que l'on peut choisir $\dot{{\cal Y}}(exp(Y)\epsilon)=\dot{{\cal Y}}(\epsilon)$. Parce que les int\'egrales orbitales ordinaires sont continues en un point r\'egulier, on a la propri\'et\'e de continuit\'e
 $$S^{\tilde{G}}(\epsilon,f)=lim_{Y\to 0}S^{\tilde{G}}(exp(Y)\epsilon,f).$$
 Cela entra\^{\i}ne que la distribution de gauche est stable. 
 Soient $\tilde{L}$ un espace de Levi de $\tilde{G}$ et $\epsilon$ un \'el\'ement semi-simple de $\tilde{L}({\mathbb R})$, r\'egulier dans cet ensemble. On d\'efinit les ensembles ${\cal Y}^{L}(\epsilon)$ et $\dot{{\cal Y}}^{L}(\epsilon)$ en rempla\c{c}ant $G$ par $L$ dans les d\'efinitions ci-dessus. On vient de d\'efinir une distribution stable $\boldsymbol{\delta}_{\epsilon}=S^{\tilde{L}}(\epsilon,.)$ sur $\tilde{L}({\mathbb R})$ associ\'ee \`a $\epsilon$, qui est une combinaison lin\'eaire d'int\'egrales orbitales. On sait donc d\'efinir la distribution $f\mapsto S_{\tilde{L}}^{\tilde{G}}(\boldsymbol{\delta}_{\epsilon},f)$ sur $\tilde{G}({\mathbb R})$. On la note simplement $f\mapsto S_{\tilde{L}}^{\tilde{G}}(\epsilon,f)$.

 Pour $f\in C_{c}^{\infty}(\tilde{G}({\mathbb R}))$ et $X\in \mathfrak{t} ({\mathbb R})\cap \mathfrak{g}_{\eta,reg}({\mathbb R})$ assez proche de $0$, on pose
 $$S_{\tilde{M}}^{\tilde{G},mod}(exp(X)\eta,f)=S_{\tilde{M}}^{\tilde{G}}(exp(X)\eta,f)+i^{\tilde{G}}_{\tilde{M}}(\eta)\vert \check{\alpha}\vert log(\vert \alpha(X)\vert )S_{\underline{\tilde{M}}}^{\tilde{G}}(exp(X)\eta,f).$$
 Remarquons que l'\'el\'ement $exp(X)\eta$ est r\'egulier dans $\tilde{M}({\mathbb R})$ et dans $\underline{\tilde{M}}({\mathbb R})$. 
 
  \ass{Proposition}{Soient $f\in C_{c}^{\infty}(\tilde{G}({\mathbb R}))$.  
 
 (i) Pour tout  $U\in Sym(\mathfrak{t})$, les limites 
 
 $lim_{r\to 0+}\partial_{U}S_{\tilde{M}}^{\tilde{G},mod}(exp(rH_{d})\eta,f)$ et $lim_{r\to 0-}\partial_{U}S_{\tilde{M}}^{\tilde{G},mod}(exp(rH_{d})\eta,f)$
 
\noindent existent.
 
 (ii) Si $w_{d}(U)=U$, ces limites sont \'egales.
 
 (iii) Pour tout $\underline{U}\in Sym(\underline{\mathfrak{t}})$, les limites 
 
 $lim_{r\to 0+}\partial_{\underline{U}}S_{\underline{\tilde{M}}}^{\tilde{G}}(exp(rH_{c})\eta,f)$ et $lim_{r\to 0-}\partial_{\underline{U}}S_{\underline{\tilde{M}}}^{\tilde{G}}(exp(rH_{c})\eta,f)$ 
 
 \noindent existent.
 
 (iv) Soit $U\in Sym(\mathfrak{t})$, supposons $w_{d}(U)=-U$ et  posons $\underline{U}=C(U)$. Alors on a les \'egalit\'es
 $$lim_{r\to 0+}\partial_{U}S_{\tilde{M}}^{\tilde{G},mod}(exp(rH_{d})\eta,f)-lim_{r\to 0-}\partial_{U}S_{\tilde{M}}^{\tilde{G},mod}(exp(rH_{d})\eta,f)$$
 $$=2\pi i\vert \check{\alpha}\vert i_{\tilde{M}}^{\tilde{G}}(\eta) \,lim_{r\to 0+}\partial_{\underline{U}}S_{\underline{\tilde{M}}}^{\tilde{G}}(exp(rH_{c})\eta,f)$$
 $$=-2\pi i\vert \check{\alpha}\vert i_{\tilde{M}}^{\tilde{G}}(\eta) \,lim_{r\to 0-}\partial_{\underline{U}}S_{\underline{\tilde{M}}}^{\tilde{G}}(exp(rH_{c})\eta,f).$$}
 
 Preuve.   Consid\'erons un \'el\'ement $X\in \mathfrak{t}({\mathbb R})\cap \mathfrak{g}_{\eta,reg}({\mathbb R})$ proche de $0$.  Remarquons que, puisque $M_{\eta}=T=M_{exp(X)\eta}$, on peut supposer $\dot{{\cal Y}}^M(exp(X)\eta)=\dot{{\cal Y}}^M(\eta)$.   
 Montrons que
 
 (1) on peut supposer
 $$\dot{{\cal Y}}^{\underline{\tilde{M}}}(exp( X)\eta)=\dot{{\cal Y}}^M(\eta).$$
 
 Puisque $\underline{M}_{exp( X)\eta}=T$, on a simplement ${\cal Y}^M(exp(X)\eta)={\cal Y}^{\underline{M}}(exp(X)\eta)\cap M\subset  {\cal Y}^{\underline{M}}(exp(X)\eta)$. On veut montrer que de  cette injection se d\'eduit une bijection
 $$T\backslash {\cal Y}^M(exp(X)\eta)/ M({\mathbb R})\simeq T\backslash {\cal Y}^{\underline{M}}(exp(X)\eta)/ \underline{M}({\mathbb R}).$$
 D'apr\`es ce que l'on a dit plus haut, il suffit de voir que les noyaux des applications
 $$H^1(\Gamma_{{\mathbb R}};T)\to H^1(\Gamma_{{\mathbb R}};M)$$
 et
 $$H^1(\Gamma_{{\mathbb R}};T)\to H^1(\Gamma_{{\mathbb R}};\underline{M})$$
 sont \'egaux. Or ces applications s'inscrivent dans un diagramme commutatif
 $$\begin{array}{ccc}&&H^1(\Gamma_{{\mathbb R}};M)\\&\nearrow&\\H^1(\Gamma_{{\mathbb R}};T)&&\downarrow\\&\searrow&\\&&H^1(\Gamma_{{\mathbb R}};\underline{M})\\ \end{array}$$
 La fl\`eche de droite est injective car $M$ est un Levi de $\underline{M}$. L'\'egalit\'e cherch\'ee en r\'esulte. Cela prouve (1).

Pour tout $y\in \dot{{\cal Y}}^M(\eta)$, on pose $X[y]=ad_{y^{-1}}(X)$. Les d\'efinitions entra\^{\i}nent
 $$(2) \qquad S_{\tilde{M}}^{\tilde{G}}(exp(X)\eta,f)=\sum_{y\in \dot{{\cal Y}}^M(\eta)}I_{\tilde{M}}^{\tilde{G}}(exp(X[y])\eta[y],f)$$
 $$-\sum_{s\in Z(\hat{M})^{\Gamma_{{\mathbb R}}}/Z(\hat{G})^{\Gamma_{{\mathbb R}}}, s\not=1}i_{\tilde{M}}(\tilde{G},\tilde{G}'(s))S_{{\bf M}}^{{\bf G}'(s)}(exp(X)\eta,f^{{\bf G}'(s)}).$$
 Soit $s\in Z(\hat{M})^{\Gamma_{{\mathbb R}}}/Z(\hat{G})^{\Gamma_{{\mathbb R}}}$ tel que ${\bf G}'(s)$ soit elliptique. On fixe des donn\'ees auxiliaires $G'_{1}(s)$,...,$\Delta_{1}(s)$, on note $\tilde{M}_{1}(s)$, resp. $\tilde{T}_{1}(s)$, les images r\'eciproques de $\tilde{M}$, resp. $\tilde{T}$, dans $\tilde{G}'_{1}(s)$; on fixe $\eta_{1}(s)\in \tilde{T}_{1}(s;{\mathbb R})$ au-dessus de $\eta$. On doit aussi fixer une d\'ecomposition
 $$ \mathfrak{t}_{1}(s)=\mathfrak{c}_{1}(s)\oplus \mathfrak{t}$$
 selon la recette de 3.1, c'est-\`a-dire en fixant au pr\'ealable une d\'ecomposition d'alg\`ebres de Lie
 $$\mathfrak{g}'_{1}(s)=\mathfrak{c}_{1}(s)\oplus \mathfrak{g}'(s).$$
  Comme on l'a expliqu\'e en 3.1, il y a une fonction $X\mapsto \tau(s;X)$ sur $\mathfrak{t}(F)$ de sorte que l'on ait l'\'egalit\'e
  $$S_{{\bf M}}^{{\bf G}'(s)}(exp( X)\eta,f^{{\bf G}'(s)})=\tau(s;X)S_{\tilde{M}_{1}(s),\lambda_{1}(s)}^{\tilde{G}'_{1}(s)}(exp(X)\eta_{1}(s),f^{\tilde{G}'_{1}(s)}).$$
  Il r\'esulte des d\'efinitions que $\tau(s,X)=\Delta_{1}(s)(exp(X)\eta_{1}(s),exp(X)\eta)^{-1}$ si $exp(X)\eta$ est fortement $\tilde{G}$-r\'egulier.

 Supposons d'abord $\hat{\alpha}(s)\not=1$, autrement dit $s\not\in Z(\hat{M})^{\hat{\alpha}}$.. En utilisant encore une fois la description du syst\`eme de racines de $G'(s)_{\eta}$, cf. [W1] 3.3, on voit que $G'(s)_{\eta}=T$. Alors $\eta$ est $\tilde{G}'(s)$-\'equisingulier et $\eta_{1}(s)$ est $\tilde{G}'_{1}(s)$-\'equisingulier. D'apr\`es [V] 1.4(2), il existe une fonction $\varphi\in C_{c}^{\infty}(\tilde{M}_{1}(s);F)$ de sorte que l'on ait l'\'egalit\'e
 $$S_{\tilde{M}_{1}(s),\lambda_{1}(s)}^{\tilde{G}'_{1}(s)}(exp(X)\eta_{1}(s),f^{\tilde{G}'_{1}(s)})=S^{\tilde{M}'_{1}(s)}(exp(X)\eta_{1}(s),\varphi).$$
 Mais le point $\eta_{1}(s)$ est r\'egulier dans $\tilde{M}'_{1}(s)$. On sait qu'alors  l'int\'egrale ci-dessus est $C^{\infty}$ en $X$. Il en est de m\^eme de la fonction $X\mapsto \tau(s;X)$. Donc $X\mapsto S_{{\bf M}}^{{\bf G}'(s)}(exp(X)\eta,f^{{\bf G}'(s)})$ est $C^{\infty}$.  Les indices $s$ n'appartenant pas \`a $Z(\hat{M})^{\hat{\alpha}}$ contribuent donc au membre de droite de l'\'egalit\'e  (2) par une fonction $C^{\infty}$. Cela permet de r\'ecrire cette \'egalit\'e  sous la forme
  $$(3) \qquad S_{\tilde{M}}^{\tilde{G}}(exp(X)\eta,f)=\phi(X)+\sum_{y\in \dot{{\cal Y}}^M(\eta)}I_{\tilde{M}}^{\tilde{G}}(exp(X[y])\eta[y],f)$$
 $$-\sum_{s\in Z(\hat{M})^{ \hat{\alpha}}/Z(\hat{G})^{\Gamma_{{\mathbb R}}}, s\not=1}i_{\tilde{M}}(\tilde{G},\tilde{G}'(s))S_{{\bf M}}^{{\bf G}'(s)}(exp(X)\eta,f^{{\bf G}'(s)}),$$
 o\`u $\phi$ est une fonction $C^{\infty}$ de $X$.

 Supposons maintenant $\hat{\alpha}(s)=1$, autrement dit $s\in Z(\hat{M})^{\hat{\alpha}}$. Dans ce cas, l'ensemble de racines de $G'(s)_{\eta}$ relatif \`a $T$ est \'egal \`a $\{\pm \alpha\}$. Puisque l'action galoisienne sur $T$ n'a pas chang\'e, $\alpha$ est encore r\'eelle, c'est-\`a-dire fixe par $\Gamma_{{\mathbb R}}$. Il en r\'esulte que $G'(s)_{\eta,SC}$ est isomorphe \`a $SL(2)$. Autrement dit, $\eta$ et $\tilde{T}$ v\'erifient encore les hypoth\`eses de 4.1 o\`u l'on remplace $\tilde{G}$ par $\tilde{G}'(s)$. On dispose donc sur $\mathfrak{t}({\mathbb R})$ de l'int\'egrale orbitale stable modifi\'ee
 $$S_{\tilde{M}_{1}(s),\lambda_{1}(s)}^{\tilde{G}'_{1}(s),mod}(exp(X)\eta_{1}(s),f^{\tilde{G}'_{1}(s)})=S_{\tilde{M}_{1}(s),\lambda_{1}(s)}^{\tilde{G}'_{1}(s)}(exp(X)\eta_{1}(s),f^{\tilde{G}'_{1}(s)})$$
 $$+i^{\tilde{G}'_{1}(s)}_{\tilde{M}_{1}}(\eta_{1}(s))\vert \check{\alpha}\vert log(\vert \alpha(X)\vert )S_{\underline{\tilde{M}}_{1}(s),\lambda_{1}(s)}^{\tilde{G}'_{1}(s)}(exp(X)\eta_{1}(s),f^{\tilde{G}'_{1}(s)}).$$
 Expliquons davantage. Le fait que l'on travaille ici avec des fonctions qui se transforment selon le caract\`ere $\lambda_{1}(s)$ de $C_{1}(s;{\mathbb R})$ ne perturbe rien: puisqu'on travaille au voisinage de $\eta_{1}(s)$, on peut remplacer ces fonctions par des fonctions \`a support compact. 
La norme $\vert \check{\alpha}\vert $ ne d\'epend pas de $s$ car la norme provient d'une forme quadratique  qui, d'apr\`es nos d\'efinitions, ne change pas quand on remplace $\tilde{G}$ par $\tilde{G}'(s)$. L'espace $\underline{\tilde{M}}'_{1}(s)$ est l'espace de Levi de $\tilde{G}'_{1}(s)$ tel que ${\cal A}_{\underline{\tilde{M}}'_{1}(s)}$ soit le noyau de $\alpha$ dans ${\cal A}_{\tilde{M}'_{1}(s)}$.  Ou encore, notons $\underline{\tilde{M}}'(s)$ l'espace de Levi de $\tilde{G}'(s)$ tel que ${\cal A}_{\underline{\tilde{M}}'(s)}$ soit le noyau de $\alpha$ dans ${\cal A}_{\tilde{M}}$, autrement dit ${\cal A}_{\underline{\tilde{M}}'(s)}={\cal A}_{\underline{\tilde{M}}}$. Alors $\underline{\tilde{M}}'_{1}(s)$ est l'image r\'eciproque de $\underline{\tilde{M}}'(s)$ dans $\tilde{G}'_{1}(s)$. Montrons que

$$(4) \qquad i^{\tilde{G}'_{1}(s)}_{\tilde{M}_{1}(s)}(\eta_{1}(s))=i^{\tilde{G}'(s)}_{\tilde{M}} (\eta).$$

On a les suites exactes
$$1\to Z(\hat{M})\to Z(\hat{M}_{1}(s))\to \hat{C}_{1}(s)\to 1,$$
$$1\to Z(\underline{\hat{M}}'(s))\to Z(\underline{\hat{M}}'_{1}(s))\to \hat{C}_{1}(s)\to 1.$$
Parce que $C_{1}(s)$ est un tore induit, $\hat{C}_{1}(s)^{\Gamma_{{\mathbb R}}}$ est connexe. On en d\'eduit que les suites
$$1\to Z(\hat{M})^{\Gamma_{{\mathbb R}}}\to Z(\hat{M}_{1}(s))^{\Gamma_{{\mathbb R}}}\to \hat{C}_{1}(s)^{\Gamma_{{\mathbb R}}}\to 1,$$
$$(5) \qquad 1\to Z(\underline{\hat{M}}'(s))^{\Gamma_{{\mathbb R}}}\to Z(\underline{\hat{M}}'_{1}(s))^{\Gamma_{{\mathbb R}}}\to \hat{C}_{1}(s)^{\Gamma_{{\mathbb R}}}\to 1.$$
sont encore exactes. De la premi\`ere se d\'eduit une suite
$$(6) \qquad 1\to Z(\hat{M})^{\hat{\alpha}}\to Z(\hat{M}_{1}(s))^{\hat{\alpha}}\to \hat{C}_{1}(s)^{\Gamma_{{\mathbb R}}}\to 1.$$
Elle est encore exacte. En effet, le seul point non \'evident est la surjectivit\'e de l'homomorphisme de droite. Cette surjectivit\'e est assur\'ee par celle du dernier homomorphisme de (5)
 et par l'inclusion
$$ Z(\underline{\hat{M}}'_{1}(s))^{\Gamma_{{\mathbb R}}}\subset Z(\hat{M}_{1}(s))^{\hat{\alpha}}.$$
L'exactitude des deux suites (5) et (6) entra\^{\i}ne l'\'egalit\'e
$$[Z(\hat{M})^{\hat{\alpha}}:Z(\underline{\hat{M}}'(s))^{\Gamma_{{\mathbb R}}}]=[ Z(\hat{M}_{1}(s))^{\hat{\alpha}}:Z(\underline{\hat{M}}'_{1}(s))^{\Gamma_{{\mathbb R}}}],$$
d'o\`u (4).

On peut aussi exprimer les int\'egrales $I_{\tilde{M}}^{\tilde{G}}(exp(X[y])\eta[y],f)$ \`a l'aide des int\'egrales modifi\'ees $I_{\tilde{M}}^{\tilde{G},mod}(exp(X)\eta[y],f)$. On transforme ainsi l'expression (3) en
$$(7) \qquad S_{\tilde{M}}^{\tilde{G}}(exp(X)\eta,f) +\vert \check{\alpha}\vert log(\vert \alpha(X)\vert )B(X)=\phi(X)+\sum_{y\in \dot{{\cal Y}}^M(\eta)}I_{\tilde{M}}^{\tilde{G},mod}(exp(X[y])\eta[y],f)$$
 $$-\sum_{s\in Z(\hat{M})^{ \hat{\alpha}}/Z(\hat{G})^{\Gamma_{{\mathbb R}}}, s\not=1}i_{\tilde{M}}(\tilde{G},\tilde{G}'(s))\tau(s;X)S_{\tilde{M}_{1}(s),\lambda_{1}(s)}^{\tilde{G}'_{1}(s),mod}(exp(X)\eta_{1}(s),f^{\tilde{G}'_{1}(s)}),$$
 o\`u
 $$B(X)=\sum_{y\in \dot{{\cal Y}}^M(\eta)}I_{\underline{\tilde{M}}}^{\tilde{G}}(exp(X[y])\eta[y],f)$$
 $$-\sum_{s\in Z(\hat{M})^{ \hat{\alpha}}/Z(\hat{G})^{\Gamma_{{\mathbb R}}}, s\not=1}i_{\tilde{M}}(\tilde{G},\tilde{G}'(s))\tau(s;X)i^{\tilde{G}'(s)}_{\tilde{M}} (\eta)S_{\underline{\tilde{M}}'_{1}(s),\lambda_{1}(s)}^{\tilde{G}'_{1}(s)}(exp(X)\eta_{1}(s),f^{\tilde{G}'_{1}(s)}).$$
 Consid\'erons $B(X)$. On peut d\'ecomposer la somme en $s$ en une somme en $t\in Z(\hat{M})^{ \hat{\alpha}}/Z(\underline{\hat{M}})^{\Gamma_{{\mathbb R}}}$ suivie d'une somme en $s\in tZ(\underline{\hat{M}})^{\Gamma_{{\mathbb R}}}/Z(\hat{G})^{\Gamma_{{\mathbb R}}}$, avec $s\not=1$. C'est-\`a-dire
 $$B(X)=\sum_{y\in \dot{{\cal Y}}^M(\eta)}I_{\underline{\tilde{M}}}^{\tilde{G}}(exp(X[y])\eta[y],f)-\sum_{t\in Z(\hat{M})^{ \hat{\alpha}}/Z(\underline{\hat{M}})^{\Gamma_{{\mathbb R}}}}B(X;t),$$
 o\`u
 $$B(X;t)=\sum_{s\in tZ(\underline{\hat{M}})^{\Gamma_{{\mathbb R}}}/Z(\hat{G})^{\Gamma_{{\mathbb R}}}, s\not=1}i_{\tilde{M}}(\tilde{G},\tilde{G}'(s))\tau(s;X)i^{\tilde{G}'(s)}_{\tilde{M}} (\eta)S_{\underline{\tilde{M}}'_{1}(s),\lambda_{1}(s)}^{\tilde{G}'_{1}(s)}(exp(X)\eta_{1}(s),f^{\tilde{G}'_{1}(s)}).$$
 Fixons $t$. On voit que, pour $s\in tZ(\underline{\hat{M}})^{\Gamma_{{\mathbb R}}}/Z(\hat{G})^{\Gamma_{{\mathbb R}}}$, $\underline{\tilde{M}}'_{1}(s)$ est l'image r\'eciproque dans $\tilde{G}'_{1}(s)$ d'un espace de Levi 
$\underline{\tilde{M}}'(t)$ ind\'ependant de $s$. C'est l'espace de la donn\'ee endoscopique $\underline{{\bf M}}'(t)$ de $(\underline{M},\underline{\tilde{M}})$. Un calcul simple conduit \`a l'\'egalit\'e
$$(8) \qquad i_{\tilde{M}}(\tilde{G},\tilde{G}'(s))i^{\tilde{G}'(s)}_{\tilde{M}} (\eta)=i_{\underline{\tilde{M}}'(t)}(\tilde{G},\tilde{G}'(s))i^{\tilde{G}}_{\tilde{M}}(\eta).$$
Pour un moment, notons pour plus de clart\'e $\boldsymbol{\delta}$ l'int\'egrale orbitale stable sur $\tilde{M}({\mathbb R})$ associ\'ee \`a $exp(X)\eta$. Comme on l'a dit, pour tout $s$, elle s'identifie \`a l'int\'egrale orbitale stable sur $\tilde{M}_{1}(s;{\mathbb R})$ associ\'ee \`a $exp(X)\eta_{1}(s)$, multipli\'ee par $\tau(s,X)$. Notons $\boldsymbol{\delta}(s)$ cette distribution. La distribution qui intervient dans la d\'efinition de $B(X;t)$ est  la distribution induite $\boldsymbol{\delta}(s)^{\underline{\tilde{M}}'_{1}(s)}$. C'est la r\'ealisation de la distribution $\boldsymbol{\delta}^{\underline{{\bf M}}'(t)}\in D_{g\acute{e}om}^{st}(\underline{{\bf M}}'(t))$.  Si on oublie la restriction $s\not=1$ dans la d\'efinition de $B(X;t)$, on voit alors en utilisant (8) que
$$B(X;t)=i^{\tilde{G}}_{\tilde{M}}(\eta)I_{\underline{\tilde{M}}}^{\tilde{G},{\cal E}}(\underline{{\bf M}}'(t),\boldsymbol{\delta}^{\underline{{\bf M}}'(t)},f).$$
On r\'etablit cette restriction $s\not=1$ en retirant de la formule ci-dessus l'\'eventuel terme correspondant \`a $s=1$. On obtient  que $B(X,t)$ est donn\'e par la formule ci-dessus pour $t\not=1$, tandis que
$$B(X;1)=i^{\tilde{G}}_{\tilde{M}}(\eta)(I_{\underline{\tilde{M}}}^{\tilde{G},{\cal E}}(\underline{{\bf M}},\boldsymbol{\delta}^{\underline{{\bf M}}},f)-S_{\underline{\tilde{M}}}^{\tilde{G}}(\boldsymbol{\delta}^{\underline{{\bf M}}},f)).$$
La distribution $\boldsymbol{\delta}^{\underline{{\bf M}}'(t)}$ \'etant \`a support $\tilde{G}$-\'equisingulier, on peut utiliser la proposition [V] 1.13 qui nous dit que
$$I_{\underline{\tilde{M}}}^{\tilde{G},{\cal E}}(\underline{{\bf M}}'(t),\boldsymbol{\delta}^{\underline{{\bf M}}'(t)},f)=I_{\underline{\tilde{M}}}^{\tilde{G}}(transfert(\boldsymbol{\delta}^{\underline{{\bf M}}'(t)}),f).$$
Le transfert commutant \`a l'induction, le transfert de $\boldsymbol{\delta}^{\underline{{\bf M}}'(t)}$ est l'induite du transfert de $\boldsymbol{\delta}$. Ce dernier est l'int\'egrale orbitale stable  dans $\tilde{M}({\mathbb R})$ associ\'ee \`a $exp(X)\eta$ ou encore la somme sur $y\in \dot{{\cal Y}}^M(\eta)$ des int\'egrales orbitales associ\'ees \`a $exp(X[y])\eta[y]$.  Son induite est la m\^eme somme d'int\'egrales orbitales, cette fois sur $\underline{\tilde{M}}({\mathbb R})$. On obtient
$$(9)\qquad B(X;t)=i^{\tilde{G}}_{\tilde{M}}(\eta)\sum_{y\in \dot{{\cal Y}}^M(\eta)}I_{\underline{\tilde{M}}}^{\tilde{G}}(exp(X[y])\eta[y],f)$$
si $t\not=1$, tandis que $B(X;1)$ est le m\^eme terme moins 
$$i^{\tilde{G}}_{\tilde{M}}(\eta)S_{\underline{\tilde{M}}}^{\tilde{G}}(exp(X)\eta,f).$$
A ce dernier terme pr\`es, les $B(X;t)$ sont tous \'egaux. Leur somme  est donc le membre de gauche de (9) multipli\'e par le nombre d'\'el\'ements de l'ensemble de sommation, c'est-\`a-dire 
$$[Z(\hat{M})^{\hat{\alpha}}:Z(\underline{\hat{M}})^{\Gamma_{{\mathbb R}}}].$$
Ce facteur compense le terme $i^{\tilde{G}}_{\tilde{M}}(\eta)$ figurant dans (9). D'o\`u
$$\sum_{t\in Z(\hat{M})^{ \hat{\alpha}}/Z(\underline{\hat{M}})^{\Gamma_{{\mathbb R}}}}B(X;t)=(\sum_{y\in \dot{{\cal Y}}^M(\eta)}I_{\underline{\tilde{M}}}^{\tilde{G}}(exp(X[y])\eta[y],f))-i^{\tilde{G}}_{\tilde{M}}(\eta)S_{\underline{\tilde{M}}}^{\tilde{G}}(exp(X)\eta,f).$$
En se reportant \`a la d\'efinition de $B(X)$, on obtient
$$B(X)=i^{\tilde{G}}_{\tilde{M}}(\eta)S_{\underline{\tilde{M}}}^{\tilde{G}}(exp(X)\eta,f).$$
Le membre de gauche de (7) devient simplement $ S_{\tilde{M}}^{\tilde{G},mod}(exp(X)\eta,f)$ et on a obtenu l'\'egalit\'e
 $$  S_{\tilde{M}}^{\tilde{G},mod}(exp(X)\eta,f)  =\phi(X)+\sum_{y\in \dot{{\cal Y}}^M(\eta)}I_{\tilde{M}}^{\tilde{G},mod}(exp(X[y])\eta[y],f)$$
 $$-\sum_{s\in Z(\hat{M})^{ \hat{\alpha}}/Z(\hat{G})^{\Gamma_{{\mathbb R}}}, s\not=1}i_{\tilde{M}}(\tilde{G},\tilde{G}'(s))\tau(s;X)S_{\tilde{M}_{1}(s),\lambda_{1}(s)}^{\tilde{G}'_{1}(s),mod}(exp(X)\eta_{1}(s),f^{\tilde{G}'_{1}(s)}).$$
 Soit $U\in Sym(\mathfrak{t})$. Pour tout $y\in \dot{{\cal Y}}^M(\eta)$, $ad_{y^{-1}}$ envoie $T$ sur un tore $T[y]$ d\'efini sur ${\mathbb R}$ et envoie conform\'ement $U$ sur un \'el\'ement $U[y]\in Sym(\mathfrak{t}[y])$. Pour tout $s\in Z(\hat{M})^{ \hat{\alpha}}/Z(\hat{G})^{\Gamma_{{\mathbb R}}}$, on a introduit en 3.1 un isomorphisme  qui devient ici un automorphisme de $Sym(\mathfrak{t})$, que l'on note $\Xi(s)$. On pose $U(s)=\Xi(s)(U)$. De l'expression pr\'ec\'edente se d\'eduit l'\'egalit\'e
$$ \partial_{U} S_{\tilde{M}}^{\tilde{G},mod}(exp(X)\eta,f)  =\partial_{U}\phi(X)+\sum_{y\in \dot{{\cal Y}}^M(\eta)}\partial_{U[y]}I_{\tilde{M}}^{\tilde{G},mod}(exp(X[y])\eta[y],f)$$
 $$-\sum_{s\in Z(\hat{M})^{ \hat{\alpha}}/Z(\hat{G})^{\Gamma_{{\mathbb R}}}, s\not=1}i_{\tilde{M}}(\tilde{G},\tilde{G}'(s))\tau(s;X)\partial_{U(s)}S_{\tilde{M}_{1}(s),\lambda_{1}(s)}^{\tilde{G}'_{1}(s),mod}(exp(X)\eta_{1}(s),f^{\tilde{G}'_{1}(s)}).$$ 
 On applique maintenant cette \'egalit\'e au point $X=rH_{d}$, o\`u $r$ est r\'eel non nul et proche de $0$. Remarquons que, puisque le tore $T_{d}$ se projette dans le centre de $\tilde{M}$, l'\'el\'ement $H_{d}$ est fix\'e par $ad_{y^{-1}}$ pour tout $y\in \dot{{\cal Y}}^M(\eta)$. Remarquons aussi que, pour tout $s\in  Z(\hat{M})^{ \hat{\alpha}}$,  $H_{d}$ provient d'un \'el\'ement de $\underline{\mathfrak{m}}'(s)_{\eta,SC}({\mathbb R})$, a fortiori d'un \'el\'ement de $\mathfrak{g}'(s)_{SC}({\mathbb R})$. Comme on l'a dit en 3.1, il en r\'esulte que la fonction $r\mapsto \tau(s;rH_{d})$ est constante. On note $\tau(s;0)$ sa valeur constante. On obtient alors
$$(10) \qquad  \partial_{U} S_{\tilde{M}}^{\tilde{G},mod}(exp(rH_{d})\eta,f)  =\partial_{U}\phi(rH_{d})+\sum_{y\in \dot{{\cal Y}}^M(\eta)}\partial_{U[y]}I_{\tilde{M}}^{\tilde{G},mod}(exp(rH_{d})\eta[y],f)$$
 $$-\sum_{s\in Z(\hat{M})^{ \hat{\alpha}}/Z(\hat{G})^{\Gamma_{{\mathbb R}}}, s\not=1}i_{\tilde{M}}(\tilde{G},\tilde{G}'(s))\tau(s;0)\partial_{U(s)}S_{\tilde{M}_{1}(s),\lambda_{1}(s)}^{\tilde{G}'_{1}(s),mod}(exp(rH_{d})\eta_{1}(s),f^{\tilde{G}'_{1}(s)}).$$  
 
 Les deux premi\`eres assertions de l'\'enonc\'e en r\'esultent. En effet, ces assertions sont v\'erifi\'ees par la fonction $\phi$ qui est $C^{\infty}$. Elles le sont pour les int\'egrales orbitales de la premi\`ere somme d'apr\`es la proposition 4.1. En raisonnant par r\'ecurrence, elles le sont aussi pour les int\'egrales orbitales stables de la deuxi\`eme somme (puisqu'on somme sur $s\not=1$). 
 
 Prouvons l'assertion (iii). L'\'el\'ement $\eta$, vu comme un \'el\'ement de $\underline{\tilde{M}}({\mathbb R})$, est $\tilde{G}$-\'equisingulier. D'apr\`es [V] 1.4(2), il existe une fonction $\varphi\in C_{c}^{\infty}(\underline{\tilde{M}}({\mathbb R}))$ tel que, pour tout $X\in \underline{\mathfrak{m}}_{\eta}({\mathbb R})$ assez proche de $0$, on ait l'\'egalit\'e
 $$S_{\underline{\tilde{M}}}^{\tilde{G}}(exp(X)\eta,f)=S^{\underline{\tilde{M}}}(exp(X)\eta,\varphi).$$
 Alors l'assertion (iii) r\'esulte des propri\'et\'es bien connues des int\'egrales orbitales stables. Le m\^eme argument d\'emontre l'\'egalit\'e des deux derni\`eres limites de l'assertion (iv), l'\'element $\underline{U}$ intervenant  \'etant antisym\'etrique, c'est-\`a-dire tel que $w_{c}(\underline{U})=-\underline{U}$.
 
 Il reste \`a prouver la premi\`ere \'egalit\'e de cette assertion (iv). 
Pour $s\in Z(\hat{M})^{\hat{\alpha}}$, on a vu que $G'(s)_{\eta,SC}$ \'etait isomorphe \`a $SL(2)$. 
 L'ensemble de racines de $G'(s)_{\eta}$ relatif au sous-tore maximal $T$ est $\{\pm \alpha\}$. Les deux groupes $G_{\eta}$ et $G'(s)_{\eta}$ sont donc quasi-d\'eploy\'es. Ils ont un tore maximalement d\'eploy\'e commun (le tore $T$) et ont m\^eme ensemble de racines. Ils sont donc isomorphes. Plus pr\'ecis\'ement, on peut fixer un isomorphisme entre eux qui soit l'identit\'e sur $T$.  
Par  cet isomorphisme, $\underline{T}$ devient un sous-tore maximal de $G'(s)_{\eta}$. Rappelons que ce dernier groupe est \'egal \`a $ \underline{M}'(s)_{\eta}$. De m\^eme que l'on a introduit un isomorphisme $C:T\to \underline{T}$, on peut introduire un isomorphisme $C(s)$. Mais l'isomorphisme entre nos deux groupes permet de choisir $C(s)=C$.  De la d\'ecomposition d\'ej\`a fix\'ee de $\mathfrak{g}'_{1}(s)$ se d\'eduit une d\'ecomposition
 $$\mathfrak{\underline{t}}_{1}(s)=\mathfrak{c}_{1}(s)\oplus \mathfrak{\underline{t}}.$$
 Avec des notations compr\'ehensibles, on a un diagramme
 $$\begin{array}{ccc}Sym(\mathfrak{t})&\stackrel{\Xi(s)}{\to}&Sym(\mathfrak{t})\\ C\downarrow\,\,&&C(s)\downarrow\,\,\\ Sym(\mathfrak{\underline{t}})&\stackrel{\underline{\Xi}(s)}{\to}&Sym(\underline{\mathfrak{t}})\\ \end{array}.$$
 Montrons qu'il est commutatif. On a vu en 3.1 que $\Xi(s)$ envoyait un \'el\'ement $H\in \mathfrak{t}$ sur $H-<H,b(s)>$, o\`u $b(s)$ est un certain \'el\'ement de  $\mathfrak{z}(G'_{1}(s))^*$.  
 L'homomorphisme $\underline{\Xi}(s)$ envoie aussi un \'el\'ement $\underline{H}\in \underline{\mathfrak{t}}$ sur $\underline{H}-<\underline{H},b(s)>$. Le chemin nord-est du diagramme envoie $H$ sur $C(H)-<H,b(s)>$. Le chemin sud-ouest envoie $H$  sur $C(H)-<C(H),b(s)>$. Mais $C$ est une conjugaison par un \'el\'ement de $G'(s)_{SC}$ et fixe donc $b(s)$ qui est central. Donc $<H,b(s)>=<C(H),b(s)>$, d'o\`u la commutativit\'e cherch\'ee.  
 
 On est dans un cas particulier tr\`es simple de la th\'eorie locale rappel\'ee en [III] section 5 et [V] 4.1 . Celle-ci nous dit  que $G'(s)_{\eta,SC}$ est un "groupe endoscopique" de $G_{\eta,SC}$. Ces deux groupes sont ici \'egaux. En particulier, on peut prendre pour facteur de transfert pour cette donn\'ee endoscopique celui qui vaut $1$ sur la diagonale.    La relation [V] 4.1(1) implique alors l'existence d'une constante $d(s)\not=0$ telle que la propri\'et\'e suivante soit v\'erifi\'ee. Soit $X_{sc}\in \mathfrak{g}_{\eta,SC}({\mathbb R})$ un \'el\'ement r\'egulier et $Z\in \mathfrak{z}(G_{\eta};{\mathbb R})$ en position g\'en\'erale. Ces \'el\'ements sont tous deux proches de $0$. On consid\`ere $X_{sc}$ comme fix\'e et $Z$ comme variable. Alors
$$ lim_{Z\to 0}\Delta_{1}(s)(exp(Z+X_{sc})\eta_{1}(s),exp(Z+X_{sc})\eta)=d(s).$$
Le premier $Z+X_{sc}$ est identifi\'e \`a un \'el\'ement de $\mathfrak{g}'_{1}(s)$ par la section $\mathfrak{g}'(s)\to \mathfrak{g}'_{1}(s)$ que l'on a fix\'ee. 
 Appliquons cela \`a $X_{sc}\in \mathfrak{t}_{sc}$. Le facteur de transfert ci-dessus n'est autre que $\tau(s;Z+X_{sc})^{-1}$. Sa limite est $\tau(s;X_{sc})^{-1}$. On a d\'ej\`a dit que $\tau(s;X_{sc})$  \'etait constant, on a not\'e sa valeur $\tau(s;0)$. On obtient $d(s)=\tau(s;0)^{-1}$. Supposons $X_{sc}\in \underline{\mathfrak{t}}_{sc}({\mathbb R})$. Alors $Z+X_{sc}\in \underline{\mathfrak{t}}({\mathbb R})$. On a dit en 3.1 qu'alors $\Delta_{1}(s)(exp(Z+X_{sc})\eta_{1}(s),exp(Z+X_{sc})\eta)$ \'etait le produit de $e^{<b(s),Z>}$ et d'une fonction localement constante sur les \'el\'ements $Z+X_{sc}$ assez r\'eguliers. Ce qui pr\'ec\`ede montre que cette fonction est constante, de valeur $\tau(s;0)^{-1}$.  On peut donc \'etendre la d\'efinition de $\Delta_{1}(s)(exp(Z+X_{sc})\eta_{1}(s),exp(Z+X_{sc})\eta)$ ou, si l'on pr\'ef\`ere, de $\Delta_{1}(s)(exp(X)\eta_{1}(s),exp(X)\eta)$ pour $X\in \underline{\mathfrak{t}}({\mathbb R})$ (sans condition de r\'egularit\'e), par l'\'egalit\'e
 $$(11)\qquad \Delta_{1}(s)(exp(X)\eta_{1}(s),exp(X)\eta)=\tau(s;0)^{-1}e^{<b(s),X>}.$$
 
Fixons $U\in Sym(\mathfrak{t})$ tel que $w_{d}(U)=-U$, posons $\underline{U}=C(U)$. Pour $X\in \underline{\mathfrak{t}}({\mathbb R})\cap \mathfrak{g}_{\eta,reg}({\mathbb R})$, posons
$$(12)\qquad  \underline{S}(X)=\sum_{s\in Z(\hat{M})^{\hat{\alpha}}/Z(\hat{G})^{\Gamma_{{\mathbb R}}}}i_{\tilde{M}}(\tilde{G},\tilde{G}'(s))i_{\tilde{M}}^{\tilde{G}'(s)}(\eta) \Delta_{1}(s)(exp(X)\eta_{1}(s),exp(X)\eta)^{-1}$$
$$\partial_{\underline{U}(s)}S_{\underline{\tilde{M}}'_{1}(s),\lambda_{1}(s)}^{\tilde{G}'_{1}(s)}(exp(X)\eta_{1}(s),f^{\tilde{G}'_{1}(s)}),$$
o\`u $\underline{U}(s)=\underline{\Xi}(s)(\underline{U})$. 
Par la m\^eme d\'ecomposition d\'ej\`a utilis\'ee, on a
$$\underline{S}(X)=\sum_{t\in Z(\hat{M})^{\hat{\alpha}}/Z(\underline{\hat{M}})^{\Gamma_{{\mathbb R}}}}\underline{S}(X,t),$$
o\`u
$$\underline{S}(X,t)=\sum_{s\in tZ(\underline{\hat{M}})^{\Gamma_{{\mathbb R}}}/Z(\hat{G})^{\Gamma_{{\mathbb R}}}}i_{\tilde{M}}(\tilde{G},\tilde{G}'(s))i_{\tilde{M}}^{\tilde{G}'(s)}(\eta) \Delta_{1}(s)(exp(X)\eta_{1}(s),exp(X)\eta)^{-1}$$
$$\partial_{\underline{U}(s)}S_{\underline{\tilde{M}}'_{1}(s),\lambda_{1}(s)}^{\tilde{G}'_{1}(s)}(exp(X)\eta_{1}(s),f^{\tilde{G}'_{1}(s)}).$$
Fixons $t\in  Z(\hat{M})^{\hat{\alpha}}/Z(\underline{\hat{M}})^{\Gamma_{{\mathbb R}}}$. Pour simplifier la notation, on consid\`ere que $t$ est l'un des points de la sommation en $s$. Comme plus haut, pour $s\in tZ(\underline{\hat{M}})^{\Gamma_{{\mathbb R}}}/Z(\hat{G})^{\Gamma_{{\mathbb R}}}$, les donn\'ees $\underline{M}'_{1}(s)$,...,$\Delta_{1}(s)$ sont diff\'erentes donn\'ees auxiliaires pour une m\^eme donn\'ee endoscopique $\underline{{\bf M}}'(t)$ de $(\underline{M},\underline{\tilde{M}})$. Les distributions
$$\Delta_{1}(s)(exp(X)\eta_{1}(s),exp(X)\eta)^{-1}\partial_{\underline{U}(s)}S^{\underline{\tilde{M}}'_{1}(s)}_{\lambda_{1}(s)}(exp(X)\eta_{1}(s),.)$$
se recollent en une m\^eme distribution $\boldsymbol{\delta}_{\underline{U},X}\in D^{st}_{g\acute{e}om,\tilde{G}-\acute{e}qui}({\bf M}'(t))$. En utilisant (8), on obtient l'\'egalit\'e
$$\underline{S}(X,t)=i_{\tilde{M}}^{\tilde{G}}(\eta)I_{\underline{\tilde{M}}}^{\tilde{G},{\cal E}}({\bf M}'(t),\boldsymbol{\delta}_{\underline{U},X},f).$$
En utilisant la proposition [V] 1.13, on obtient
$$\underline{S}(X,t)=i_{\tilde{M}}^{\tilde{G}}(\eta)I_{\underline{\tilde{M}}}^{\tilde{G}}(transfert(\boldsymbol{\delta}_{\underline{U},X}),f).$$
Calculons le transfert de la distribution $\boldsymbol{\delta}_{\underline{U},X}$ en supposant d'abord $exp(X)\eta$ fortement r\'egulier. On r\'ealise maintenant cette distribution en utilisant les donn\'ees auxiliaires $\underline{M}'_{1}(t)$,...,$\Delta_{1}(t)$. Notons $\boldsymbol{\delta}'_{X}$ l'int\'egrale orbitale stable associ\'ee \`a $exp(X)\eta_{1}(t)$. Par d\'efinition, son transfert est la distribution
$$\varphi\mapsto \sum_{y\in \dot{{\cal Y}}^{\underline{\tilde{M}}}(exp(X)\eta)}\Delta_{1}(t)(exp(X)\eta_{1}(t),exp(X[y])\eta[y])I^{\underline{\tilde{M}}}(exp(X[y])\eta[y],\varphi)$$
sur $\underline{\tilde{M}}({\mathbb R})$, 
avec des notations similaires \`a celles utilis\'ees plus haut. En appliquant \`a $\boldsymbol{\delta}'_{X}$ l'op\'erateur diff\'erentiel $\partial_{\underline{U}(t)}$, on obtient une distribution que l'on note $\boldsymbol{\delta}'_{\underline{U},X}$. Pour tout $y\in   \dot{{\cal Y}}^{\underline{\tilde{M}}}(exp(X)\eta)$, notons $\underline{T}[y]=ad_{y^{-1}}(\underline{T})$. On a encore un isomorphisme $\underline{\Xi}(t,y):Sym(\underline{\mathfrak{t}}[y])\to Sym(\underline{\mathfrak{t}})$. En notant $\underline{U}[y]$ l'image r\'eciproque de $\underline{U}(t)$ par cet isomorphisme, le transfert de $\boldsymbol{\delta}'_{\underline{U},X}$ est
la distribution
$$\varphi\mapsto \sum_{y\in \dot{{\cal Y}}^{\underline{\tilde{M}}}(exp(X)\eta)}\Delta_{1}(t)(exp(X)\eta_{1}(t),exp(X[y])\eta[y])\partial_{\underline{U}[y]}I^{\underline{\tilde{M}}}(exp(X[y])\eta[y],\varphi).$$
Si $y=1$, on a $\underline{\Xi}(t,1)=\underline{\Xi}(t)$ donc $\underline{U}[1]=\underline{U}$. En fait, il r\'esulte des d\'efinitions que le diagramme
$$\begin{array}{ccc}Sym(\underline{\mathfrak{t}})&&\\ &\searrow \underline{\Xi}(t)&\\ ad_{y^{-1}}\downarrow&&Sym(\underline{\mathfrak{t}})\\ &\nearrow \underline{\Xi}(t,y)&\\ Sym(\underline{\mathfrak{t}}[y])&&\\ \end{array}$$
est commutatif. Donc $\underline{U}[y]$ est simplement l'image de $\underline{U}$ par l'isomorphisme $ad_{y^{-1}}$. Pour obtenir le transfert de $\boldsymbol{\delta}_{\underline{U},X}$, il reste \`a diviser par $\Delta_{1}(t)(exp(X)\eta_{1}(t),exp(X)\eta)$. On obtient que ce transfert est la distribution
$$\varphi\mapsto \sum_{y\in \dot{{\cal Y}}^{\underline{\tilde{M}}}(exp(X)\eta)}\frac{\Delta_{1}(t)(exp(X)\eta_{1}(t),exp(X[y])\eta[y])}{\Delta_{1}(t)(exp(X)\eta_{1}(t),exp(X)\eta)}\partial_{\underline{U}[y]}I^{\underline{\tilde{M}}}(exp(X[y])\eta[y],\varphi).$$ 
On en d\'eduit l'\'egalit\'e
$$(13) \qquad   \underline{S}(X,t)=i_{\tilde{M}}^{\tilde{G}}(\eta) \sum_{y\in \dot{{\cal Y}}^{\underline{\tilde{M}}}(exp(X)\eta)}\frac{\Delta_{1}(t)(exp(X)\eta_{1}(t),exp(X[y])\eta[y])}{\Delta_{1}(t)(exp(X)\eta_{1}(t),exp(X)\eta)}$$
$$\partial_{\underline{U}[y]}I_{\underline{\tilde{M}}}^{\tilde{G}}(exp(X[y])\eta[y],f).$$
Remarquons que l'on peut choisir l'ensemble $\dot{{\cal Y}}^{\underline{\tilde{M}}}(exp(X)\eta)$ ind\'ependant de $X$: il ne d\'epend que de $\underline{T}$. Notons-le plut\^ot $\dot{{\cal Y}}^{\underline{\tilde{M}}}(\underline{T})$. Les consid\'erations faites plus haut sur le facteur de transfert $\Delta_{1}(t)(exp(X)\eta_{1}(t),exp(X)\eta)$ s'appliquent aussi au facteur 
$$\Delta_{1}(t)(exp(X)\eta_{1}(t),exp(X[y])\eta[y])$$
 pour $y\in \dot{{\cal Y}}^{\underline{\tilde{M}}}(\underline{T})$. C'est-\`a-dire que $G'(t)_{\eta,SC}$ est un groupe endoscopique de $G_{\eta[y],SC}$. Ce dernier groupe n'est plus toujours $SL(2)$, ce peut \^etre une forme int\'erieure. Mais le premier groupe est encore $SL(2)$ donc est la forme quasi-d\'eploy\'ee du second. On peut encore choisir pour facteur de transfert  pour cette donn\'ee endoscopique celui qui vaut $1$ sur tout couple d'\'el\'ements qui se correspondent. Il existe alors une constante non nulle $d(t,y)$ v\'erifiant la propri\'et\'e suivante. Soient $X_{sc}[y]\in \mathfrak{g}_{\eta[y],SC}({\mathbb R})$ et $X'_{sc}\in \mathfrak{g}'(t)_{\eta,SC}({\mathbb R})$ deux \'el\'ements r\'eguliers dont les classes de conjugaison stables se correspondent et soit  $Z\in \mathfrak{z}(G_{\eta};{\mathbb R})$ en position g\'en\'erale. Ces \'el\'ements sont tous deux proches de $0$.  Alors
$$(14) \qquad  lim_{Z\to 0}\Delta_{1}(t)(exp(Z+X'_{sc})\eta_{1}(t),exp(Z[y]+X_{sc}[y])\eta)=d(t,y).$$
On peut prendre en particulier $X'_{sc}=rH_{c}$ , pour un r\'eel $r\not=0$ et proche de $0$, et $X_{sc}[y]=rH_{c}[y]$. 
Appliquons l'\'egalit\'e (13) \`a $X=Z+rH_{c}$  et pour un \'el\'ement $Z\in \mathfrak{z}(G_{\eta};{\mathbb R})$ en position g\'en\'erale. Faisons tendre $Z$ vers $0$. Tous les termes sont continus en $Z$ et on obtient l'\'egalit\'e
$$ \underline{S}(rH_{c},t)=i_{\tilde{M}}^{\tilde{G}}(\eta) \sum_{y\in \dot{{\cal Y}}^{\underline{\tilde{M}}}(\underline{T})}d(t,y)d(t,1)^{-1}\partial_{\underline{U}[y]}I_{\underline{\tilde{M}}}^{\tilde{G}}(exp(rH_{c}[y])\eta[y],f).$$
D'o\`u
$$(15)\qquad \underline{S}(rH_{c})=\sum_{y\in \dot{{\cal Y}}^{\underline{\tilde{M}}}(\underline{T})}d(y)\partial_{\underline{U}[y]}I_{\underline{\tilde{M}}}^{\tilde{G}}(exp(rH_{c}[y])\eta[y],f),$$
o\`u
$$d(y)=i_{\tilde{M}}^{\tilde{G}}(\eta)\sum_{t\in Z(\hat{M})^{\hat{\alpha}}/Z(\underline{\hat{M}})^{\Gamma_{{\mathbb R}}}}d(t,y)d(t,1)^{-1}.$$

 On a par d\'efinition l'inclusion ${\cal Y}^M(\eta)= {\cal Y}^{\underline{M}}(\eta)\cap M$. D'o\`u une application
 $$(16) \qquad T\backslash {\cal Y}^M(\eta)/M({\mathbb R})\to \underline{M}_{\eta}\backslash {\cal Y}^{\underline{M}}(\eta)/\underline{M}({\mathbb R}).$$
 Montrons que
 
 (17) cette application est injective; son image est l'image dans l'espace d'arriv\'ee de l'ensemble des $y\in{\cal Y}^{\underline{M}}(\eta)$ tels que le groupe $\underline{M}_{\eta[y],SC}$ soit isomorphe \`a $SL(2)$.
 
 Le membre de gauche de (16) s'identifie au noyau de l'application
 $$H^1(\Gamma_{{\mathbb R}};T)\to H^1(\Gamma_{{\mathbb R}};M)$$
 tandis que celui de droite s'identifie au noyau de l'application
 $$H^1(\Gamma_{{\mathbb R}};\underline{M}_{\eta})\to H^1(\Gamma_{{\mathbb R}};\underline{M}).$$
 L'injectivit\'e de l'application (16) r\'esulte de celle de l'application
 $$H^1(\Gamma_{{\mathbb R}};T)\to H^1(\Gamma_{{\mathbb R}};\underline{M}_{\eta}),$$
 laquelle provient du fait que $T$ est un Levi de $\underline{M}_{\eta}$. Si $y\in {\cal Y}^{\underline{M}}(\eta)$, $ad_{y^{-1}}$ se rel\`eve en un torseur int\'erieur de $\underline{M}_{\eta,SC}$ sur $\underline{M}_{\eta[y],SC}$.  Rappelons que le sous-tore $T_{d}$ de $\underline{M}_{\eta,SC}$ se projette sur un sous-tore d\'eploy\'e du centre de $M$. Si $y\in M$, le tore $ad_{y^{-1}}(T_{d})$ se projette sur le m\^eme tore. Il en r\'esulte que $ad_{y^{-1}}(T_{d})$ est d\'efini sur ${\mathbb R}$ et est d\'eploy\'e. Le groupe  $\underline{M}_{\eta[y],SC}$ est  donc une forme int\'erieure de $SL(2)$ qui contient un sous-tore d\'eploy\'e non trivial.  Cela implique qu'il est isomorphe \`a $SL(2)$. Inversement, supposons que $\underline{M}_{\eta[y],SC}$ soit isomorphe \`a $SL(2)$. Quitte \`a multiplier \`a gauche $y$ par un \'el\'ement de $\underline{M}_{\eta}$, on peut supposer que $ad_{y^{-1}}$ induit un isomorphisme d\'efini sur ${\mathbb R}$ de $\underline{M}_{\eta,SC}$ sur $\underline{M}_{\eta[y],SC}$. Alors $ad_{y^{-1}}$ induit aussi un isomorphisme d\'efini sur ${\mathbb R}$ de $\underline{M}_{\eta}$ sur $\underline{M}_{\eta[y]}$. Le groupe $A_{\tilde{M}}$ est inclus dans $\underline{M}_{\eta}$. Posons $A'=ad_{y^{-1}}(A_{\tilde{M}})$. Alors $A'$ est un tore d\'efini et d\'eploy\'e sur ${\mathbb R}$. Son commutant est le groupe $M'=ad_{y^{-1}}(M)$ et est d\'efini sur ${\mathbb R}$. Soit $P\in {\cal P}(M)$, posons $P'=ad_{y^{-1}}(P)$. Puisque $P$ est d\'etermin\'e par son ensemble associ\'e de racines positives, qui sont des caract\`eres de $A_{\tilde{M}}$,  le parabolique $P'$ est d\'etermin\'e par le m\^eme ensemble transport\'e \`a $A'$ par $ad_{y^{-1}}$. Donc $P'$ est d\'efini sur ${\mathbb R}$. Les paires $(P,M)$ et $(P',M')$ sont toutes deux d\'efinies sur ${\mathbb R}$ et sont conjugu\'ees par $y\in \underline{M}({\mathbb C})$. On sait qu'alors elles sont aussi conjugu\'ees par un \'el\'ement de  $\underline{M}({\mathbb R})$. Autrement dit, quitte \`a multiplier \`a droite $y$ par un \'el\'ement de $\underline{M}({\mathbb R})$, on peut supposer que les deux paires sont \'egales. On a alors $ad_{y^{-1}}(P,M)=(P,M)$, ce qui entra\^{\i}ne $y\in M$. Cela prouve (17).
 
 En vertu de (17), on peut supposer $\dot{{\cal Y}}^M(\eta)\subset \dot{{\cal Y}}^{\underline{\tilde{M}}}(\eta)$. La d\'emonstration montre que l'on peut supposer que, pour tout $y\in \dot{{\cal Y}}^M(\eta)$, $ad_{y^{-1}}$ se restreint en un isomorphisme d\'efini sur ${\mathbb R}$ de $\underline{M}_{\eta}$ sur $\underline{M}_{\eta[y]}$. 
 
  Puisque $\underline{T}\subset \underline{M}_{\eta}$, on a l'inclusion ${\cal Y}^{\underline{M}}(exp(X)\eta)\subset {\cal Y}^{\underline{M}}(\eta)$ pour tout $X\in \underline{\mathfrak{t}}({\mathbb R})\cap \mathfrak{g}_{\eta,reg}({\mathbb R})$. On en d\'eduit une application naturelle
 $$\underline{T}\backslash {\cal Y}^{\underline{M}}(exp(X)\eta)/\underline{M}({\mathbb R})\to \underline{M}_{\eta}\backslash {\cal Y}^{\underline{M}}(\eta)/\underline{M}({\mathbb R}).$$
 Cette application est surjective. Cela r\'esulte du fait que $\underline{T}$ est un sous-tore fondamental de $\underline{M}_{\eta}$ et se transf\`ere donc \`a toute forme int\'erieure de ce groupe. 
 Autrement dit, on peut choisir notre syst\`eme $\dot{{\cal Y}}^{\underline{\tilde{M}}}(\underline{T})$    v\'erifiant la propri\'et\'e suivante. Il y a une application surjective
 $$q:\dot{{\cal Y}}^{\underline{\tilde{M}}}(\underline{T})\to \dot{{\cal Y}}^{\underline{\tilde{M}}}(\eta)$$
 telle que, pour tout $y\in \dot{{\cal Y}}^{\underline{\tilde{M}}}(\underline{T})$, $q(y)$ est l'unique \'el\'ement de $\underline{M}_{\eta}y \cap  \dot{{\cal Y}}^{\underline{\tilde{M}}}(\eta)$. On d\'etermine les fibres de l'application $q$. Pour $y\in  \dot{{\cal Y}}^{\underline{\tilde{M}}}(\eta)$, la classe de conjugaison stable dans $\underline{\mathfrak{m}}_{\eta,SC}$  
d'un  $X$  comme ci-dessus se transf\`ere par $ad_{y^{-1}}$ en une classe de conjugaison stable dans $\underline{\mathfrak{m}}_{\eta[y],SC}$. L'ensemble $\{ad_{y'}^{-1}(X);y'\in q^{-1}(y)\}$ est un ensemble de repr\'esentants des classes de conjugaison par $\underline{M}_{\eta[y]}({\mathbb R})$ dans cette classe de conjugaison stable.

 Montrons que
 
 (18) on a $d(y)=1$ pour tout $y\in \dot{{\cal Y}}^{\underline{\tilde{M}}}(\underline{T})$ tel que $q(y)\in \dot{{\cal Y}}^M(\eta)$ tandis que $d(y)=0$ pour tout $y\in \dot{{\cal Y}}^{\underline{\tilde{M}}}(\underline{T})$ tel que $q(y)\not\in \dot{{\cal Y}}^M(\eta)$.
 
 C'est en fait un cas particulier de la proposition [III] 8.4. Reprenons partiellement le calcul. Si $q(y)\in  \dot{{\cal Y}}^M(\eta)$, on peut appliquer la d\'efinition (14) aux \'el\'ements $X_{sc}[y]=H_{d}[q(y)]$ et $X'_{sc}=H_{d}$. On obtient que, pour tout $t\in Z(\hat{M})^{\hat{\alpha}}/Z(\underline{\hat{M}})^{\Gamma_{{\mathbb R}}}$, on a l'\'egalit\'e
 $$d(t,y)d(t,1)^{-1}=lim_{Z\to 0}\frac{\Delta_{1}(t)(exp(Z+H_{d})\eta_{1}(t),exp(Z[y]+H_{d}[q(y)])\eta)}{\Delta_{1}(t)(exp(Z+H_{d})\eta_{1}(t),exp(Z+H_{d})\eta)}.$$
 Les \'el\'ements $exp(Z[y]+H_{d}[q(y)])\eta$ et $exp(Z+H_{d})\eta$ appartiennent \`a la m\^eme classe de conjugaison stable dans $\tilde{M}({\mathbb R})$. Le facteur $\Delta_{1}(t)$, restreint aux \'el\'ements de $\tilde{M}_{1}(t;{\mathbb R})\times \tilde{M}({\mathbb R})$, est un facteur pour la donn\'ee endoscopique triviale ${\bf M}$ de $(M,\tilde{M})$. Il est donc invariant par conjugaison stable en la deuxi\`eme variable. Le rapport ci-dessus vaut donc $1$. Puisque $i_{\tilde{M}}^{\tilde{G}}(\eta)$ est l'inverse du nombre d'\'el\'ements de l'ensemble de sommation $Z(\hat{M})^{\hat{\alpha}}/Z(\underline{\hat{M}})^{\Gamma_{{\mathbb R}}}$, on obtient la premi\`ere assertion. Supposons maintenant $q(y)\not\in \dot{{\cal Y}}^M(\eta)$. Le cocycle $\sigma\mapsto y\sigma(y)^{-1}$ \`a valeurs dans $G_{\eta}$ se pousse en un cocycle \`a valeurs dans $G_{\eta,ad}$, ce groupe \'etant l'image de $G_{\eta}$ dans $G_{AD}$. On note $\tau[y]$ sa classe. Introduisons le tore $\hat{T}$ dual de $T$ (ce n'est plus le tore d'une paire de Borel \'epingl\'ee $\hat{{\cal E}}$ comme au d\'ebut du paragraphe). Puisque $T$ est un sous-tore maximal de $M$ comme de $G_{\eta}$, on a des plongements $\hat{T}\subset \hat{M}$ et $\hat{T}\subset \hat{G}_{\eta}$. Le second r\'ealise $\hat{T}$ comme Levi de $\hat{G}_{\eta}$. Il est \'equivariant pour les actions galoisiennes. Du premier plongement se d\'eduit une inclusion $Z(\hat{M})\to \hat{T}$, qui est \'equivariante pour les actions galoisiennes. D'o\`u un plongement $Z(\hat{M})^{\Gamma_{{\mathbb R}}}\subset \hat{T}^{\Gamma_{{\mathbb R}}}\subset \hat{G}_{\eta}^{\Gamma_{{\mathbb R}}}$. Puisque $\hat{\alpha}$ est la seule racine (au signe pr\`es) de $\hat{G}_{\eta}$,  $Z(\hat{M})^{\hat{\alpha}}$ est exactement l'image r\'eciproque de $Z(\hat{G}_{\eta})^{\Gamma_{{\mathbb R}}}$ par le plongement pr\'ec\'edent.  Notons $\hat{M}_{sc}$ et $\hat{T}_{sc}$ les images r\'eciproques de $\hat{M}$ et $\hat{T}$ dans $\hat{G}$ et notons $\hat{G}_{\eta,sc}$ le groupe dual de $G_{\eta,ad}$. On a pour ces groupes des propri\'et\'es analogues aux pr\'ec\'edentes. Soit $t\in Z(\hat{M})^{\hat{\alpha}}$. On choisit un \'el\'ement $t_{sc}\in Z(\hat{M}_{sc})^{\Gamma_{{\mathbb R}},0}$ qui se projette dans $\hat{G}_{AD}$ sur le m\^eme \'el\'ement que $t$.  Alors $t_{sc}$ s'identifie \`a un \'el\'ement de $Z(\hat{G}_{\eta,sc})^{\Gamma_{{\mathbb R}}}$, que l'on projette dans le groupe des composantes connexes $\pi_{0}(Z(\hat{G}_{\eta,sc})^{\Gamma_{{\mathbb R}}})$. On note $u(t)$ son image. Remarquons qu'elle peut d\'ependre du choix de $t_{sc}$. On a un produit sur
 $$H^1(\Gamma_{{\mathbb R}};G_{\eta,ad})\times \pi_{0}(Z(\hat{G}_{\eta,sc})^{\Gamma_{{\mathbb R}}}).$$
 Le calcul de [III] 8.4 montre que
 $$d(t,y)d(t,1)^{-1}=<\tau[y],u(t)>.$$
 
 {\bf Remarque.} On pourrait aussi utiliser le th\'eor\`eme 5.1.D de [KS] pour obtenir cette formule.
 
Notons  $Ann[y]$ l'annulateur de $\tau[y]$ dans $\pi_{0}(Z(\hat{G}_{\eta,sc})^{\Gamma_{{\mathbb R}}})$. La formule montre que  l'image  $v(t)$ de $u(t)$ dans le quotient $\pi_{0}(Z(\hat{G}_{\eta,sc})^{\Gamma_{{\mathbb R}}})/Ann[y]$ est bien d\'etermin\'ee. On obtient ainsi une application
$$v: Z(\hat{M})^{\hat{\alpha}}/Z(\underline{\hat{M}})^{\Gamma_{{\mathbb R}}}\to \pi_{0}(Z(\hat{G}_{\eta,sc})^{\Gamma_{{\mathbb R}}})/Ann[y].$$
C'est un homomorphisme. Pour prouver la nullit\'e de $d(y)$, il suffit de prouver que l'image de $v$ n'est pas r\'eduite \`a $\{1\}$. Or cette image contient l'image naturelle du noyau de l'homomorphisme 
$$(19)\qquad  \pi_{0}(Z(\hat{G}_{\eta,sc})^{\Gamma_{{\mathbb R}}})\to \pi_{0}(\hat{T}_{sc}^{\Gamma_{{\mathbb R}}}).$$
En effet, un \'el\'ement $x$ de ce noyau se rel\`eve en un \'el\'ement  $t_{sc}\in Z(\hat{G}_{\eta,sc})^{\Gamma_{{\mathbb R}}}\cap \hat{T}_{sc}^{\Gamma_{{\mathbb R}},0}$. Parce que $T$ est elliptique dans $M$, on a l'\'egalit\'e $\hat{T}_{sc}^{\Gamma_{{\mathbb R}},0}=Z(\hat{M}_{sc})^{\Gamma_{{\mathbb R}},0}$. Donc $t_{sc}$ se projette en  un \'el\'ement $t\in Z(\hat{M})^{\Gamma_{{\mathbb R}}}$. Parce que $t_{sc}\in Z(\hat{G}_{\eta,sc})$, on a $t\in Z(\hat{M})^{\hat{\alpha}}$. On voit alors que $x=v(t)$, d'o\`u l'assertion. Il nous suffit de prouver que le noyau de (19) n'est pas contenu dans $Ann[y]$. On a un homomorphisme naturel
$$\pi_{0}(Z(\hat{G}_{\eta,SC})^{\Gamma_{{\mathbb R}}})\to \pi_{0}(Z(\hat{G}_{\eta,sc})^{\Gamma_{{\mathbb R}}}).$$
Son image est contenue dans le noyau de (19). En effet, puisque $\hat{G}_{\eta,SC}=SL(2,{\mathbb C})$, le groupe 
$\pi_{0}(Z(\hat{G}_{\eta,SC})^{\Gamma_{{\mathbb R}}})$ n'est autre que le centre $\{\pm 1\}$ de $ \hat{G}_{\eta,SC}$, qui est contenu dans $\hat{T}_{d}=\hat{T}_{d}^{\Gamma_{{\mathbb R}},0}$, donc se projette dans $\hat{T}_{sc}^{\Gamma_{{\mathbb R}},0}$. Il suffit donc de prouver que l'image de $ \pi_{0}(Z(\hat{G}_{\eta,SC})^{\Gamma_{{\mathbb R}}})$ n'est pas contenue dans $Ann[y]$. Cela \'equivaut \`a ce que l'image $\tau[y]_{AD}$ de $\tau[y]$ dans $H^1(\Gamma_{{\mathbb R}};G_{\eta,AD})$ ne soit pas dans le noyau de l'accouplement entre les deux ensembles 
$$H^1(\Gamma_{{\mathbb R}};G_{\eta,AD})\times \pi_{0}(Z(\hat{G}_{\eta,SC})^{\Gamma_{{\mathbb R}}}).$$
Il r\'esulte de [K] th\'eor\`eme 1.2 que ce dernier noyau est r\'eduit \`a l'\'el\'ement $1\in H^1(\Gamma_{{\mathbb R}};G_{\eta,AD})$. Or $\tau[y]_{AD}$ d\'etermine la forme int\'erieure $G_{\eta[y],SC}$ de $G_{\eta,SC}$. L'hypoth\`ese sur $y$ et l'assertion (17) entra\^{\i}nent que cette forme int\'erieure est non d\'eploy\'ee, donc que $\tau[y]_{AD}$ n'est pas \'egal \`a $1$. Cela ach\`eve la preuve de  (18).

En utilisant (18), on transforme l'\'egalit\'e (15) en
$$(20)\qquad \underline{S}(rH_{c})=\sum_{y\in \dot{{\cal Y}}^M(\eta)}\sum_{y'\in q^{-1}(y)}\partial_{\underline{U}[y']}I_{\underline{\tilde{M}}}^{\tilde{G}}(exp(rH_{c}[y'])\eta[y],f).$$

 Rempla\c{c}ons $X$ par $rH_{c}$ dans la formule (12) et utilisons (11). On isole le terme $s=1$ pour lequel on peut prendre pour ${\bf G}'(s)$ des donn\'ees auxiliaires triviales. On a pour celles-ci $\tau(1,0)=1$. On obtient
 $$\underline{S}(rH_{c})=i_{\tilde{M}}^{\tilde{G}}(\eta)\partial_{\underline{U}}S_{\underline{\tilde{M}}}^{\tilde{G}}(exp(rH_{c}\eta,f)
$$
$$+\sum_{s\in Z(\hat{M})^{\hat{\alpha}}/Z(\hat{G})^{\Gamma_{{\mathbb R}}}, s\not=1}i_{\tilde{M}}(\tilde{G},\tilde{G}'(s))i_{\tilde{M}}^{\tilde{G}'(s)}(\eta)\tau(s;0)\partial_{\underline{U}(s)}S_{\underline{\tilde{M}}'_{1}(s),\lambda_{1}(s)}^{\tilde{G}'_{1}(s)}(exp(rH_{c})\eta_{1}(s),f^{\tilde{G}'_{1}(s)}).$$
En utilisant cette \'egalit\'e et l'\'egalit\'e (10), on voit que
$$(21)\qquad \partial_{U}S_{\tilde{M}}^{\tilde{G},mod}(exp(rH_{d})\eta,f)-\partial_{U}S_{\tilde{M}}^{\tilde{G},mod}(exp(-rH_{d})\eta,f)$$
$$-2\pi i\vert  \check{\alpha}\vert i_{\tilde{M}}^{\tilde{G}}(\eta)\partial_{\underline{U}}S_{\underline{\tilde{M}}}^{\tilde{G}}(exp(rH_{c})\eta,f)$$
est la somme des expressions suivantes
 $$(22)\qquad \partial_{U}\phi(rH_{d})-\partial_{U}\phi(-rH_{d});$$
 $$(23)\qquad \sum_{y\in \dot{\cal Y}^M(\eta)}(\partial_{U[y]}I_{\tilde{M}}^{\tilde{G},mod}(exp(rH_{d})\eta[y],f)-\partial_{U[y]}I_{\tilde{M}}^{\tilde{G},mod}(exp(-rH_{d})\eta[y],f));$$
 $$(24) \qquad -2\pi i\vert  \check{\alpha}\vert \underline{S}(rH_{c});$$
 $$(25)\qquad \sum_{s\in Z(\hat{M})^{\hat{\alpha}}/Z(\hat{G})^{\Gamma_{{\mathbb R}}}, s\not=1}i_{\tilde{M}}(\tilde{G},\tilde{G}'(s)) \tau(s;0)C(r,s),$$
 o\`u
 $$C(r,s)=2\pi i\vert  \check{\alpha}\vert i_{\tilde{M}}^{\tilde{G}'(s)}(\eta)\partial_{\underline{U}(s)}S_{\underline{\tilde{M}}'_{1}(s),\lambda_{1}(s)}^{\tilde{G}'_{1}(s)}(exp(rH_{c})\eta_{1}(s),f^{\tilde{G}'_{1}(s)})$$
 $$-\partial_{U(s)}S_{\tilde{M}_{1}(s),\lambda_{1}(s)}^{\tilde{G}'_{1}(s)}(exp(rH_{d})\eta_{1}(s),f^{\tilde{G}'_{1}(s)})+\partial_{U(s)}S_{\tilde{M}_{1}(s),\lambda_{1}(s)}^{\tilde{G}'_{1}(s)}(exp(-rH_{d})\eta_{1}(s),f^{\tilde{G}'_{1}(s)}).$$
 Parce que $\phi$ est $C^{\infty}$, la limite de (22) quand $r$ tend vers $0$ est nulle. Pour $s\in Z(\hat{M})^{\hat{\alpha}}/Z(\hat{G})^{\Gamma_{{\mathbb R}}}$, $s\not=1$, on peut utiliser par r\'ecurrence l'assertion (iv) de l'\'enonc\'e. On a donc $lim_{r\to 0}C(r,s)=0$ et la limite de l'espression (25) quand $r$ tend vers $0$ est nulle. En utilisant (20), on voit que la somme de (23) et de (24) s'\'ecrit
 $$\sum_{y\in \dot{\cal Y}^M(\eta)}D(r,y),$$
 o\`u
 $$D(r,y)=\partial_{U[y]}I_{\tilde{M}}^{\tilde{G},mod}(exp(rH_{d})\eta[y],f)-\partial_{U[y]}I_{\tilde{M}}^{\tilde{G},mod}(exp(-rH_{d})\eta[y],f)$$
 $$-2\pi i\vert  \check{\alpha}\vert \sum_{y'\in q^{-1}(y)}\partial_{\underline{U}[y']}I_{\underline{\tilde{M}}}^{\tilde{G}}(exp(rH_{c}[y'])\eta[y],f).$$
 Fixons $y$. Comme on l'a dit, le nombre d'\'el\'ements de $q^{-1}(y)$ est le nombre de classes de conjugaison par $G_{\eta[y]}({\mathbb R})$ dans la classe de conjugaison stable de $exp(rH_{c}[y'])\eta[y])$, pour tout $y'\in q^{-1}(y)$. Autrement dit, il est \'egal \`a $2^{c(\eta[y])}$, avec la notation de la proposition 4.1. On peut donc r\'ecrire
 $$D(r,y)=2^{-c(\eta[y])}\sum_{y'\in q^{-1}(y)}D'(r,y'),$$
 o\`u
 $$D'(r,y')=\partial_{U[y]}I_{\tilde{M}}^{\tilde{G},mod}(exp(rH_{d})\eta[y],f)-\partial_{U[y]}I_{\tilde{M}}^{\tilde{G},mod}(exp(-rH_{d})\eta[y],f)$$
 $$-2^{1+c(\eta[y])}\pi i\vert  \check{\alpha}\vert \partial_{\underline{U}[y']}I_{\underline{\tilde{M}}}^{\tilde{G}}(exp(rH_{c}[y'])\eta[y],f).$$
  La proposition 4.1 nous dit que, pour tout $y'\in q^{-1}(y)$,  $lim_{r\to 0}D'(r,y')=0$. Il en r\'esulte que la limite quand $r$ tend vers $0$ de la somme des expressions (23) et (24)  est nulle. On conclut que la limite quand $r$ tend vers $0$ de l'expression (21) est nulle. C'est ce qu'affirme le (iv) de l'\'enonc\'e. Cela ach\`eve la d\'emonstration. $\square$

 \bigskip
 
 \subsection{Sauts des int\'egrales orbitales pond\'er\'ees endoscopiques}
 On conserve les hypoth\`eses de 4.1, le triplet $(G,\tilde{G},{\bf a})$ \'etant de nouveau g\'en\'eral. On suppose que $\tilde{G}$ est une composante connexe d'un $K$-espace $K\tilde{G}$ et que les espaces de Levi $\tilde{M}$ et $\underline{\tilde{M}}$ sont   des composantes connexes de $K$-espaces de Levi $K\tilde{M}$ et $K\tilde{G}$. Pour  $X\in \mathfrak{t}^{\theta}({\mathbb R})\cap \mathfrak{g}_{\eta,reg}({\mathbb R})$ assez proche de $0$, notons  $\boldsymbol{\gamma}(X)$ la distribution $\varphi\mapsto I^{K\tilde{M}}(exp(X)\eta,\omega,\varphi)$ pour $\varphi\in C_{c}^{\infty}(K\tilde{M}({\mathbb R}))$. Parce que
  l'\'el\'ement $exp(X)\eta$ de $\tilde{M}({\mathbb R})$ est $\tilde{G}$-\'equisingulier, $\boldsymbol{\gamma}(X)$ appartient \`a $D_{g\acute{e}om,\tilde{G}-\acute{e}qui}(K\tilde{M}({\mathbb R}),\omega)$. On a d\'efini en [V] 1.8 la distribution
  $$f\mapsto I_{K\tilde{M}}^{K\tilde{G},{\cal E}}(\boldsymbol{\gamma}(X),f)$$
  pour $f\in C_{c}^{\infty}(K\tilde{G}({\mathbb R}))$. On note maintenant ce terme $I_{K\tilde{M}}^{K\tilde{G},{\cal E}}(exp(X)\eta,\omega,f)$. De m\^eme, on d\'efinit $I_{K\underline{\tilde{M}}}^{K\tilde{G},{\cal E}}(exp(X)\eta,\omega,f)$. On d\'efinit alors comme en 4.1 l'int\'egrale modifi\'ee
  $$I_{K\tilde{M}}^{K\tilde{G},{\cal E},mod}(exp(X)\eta,\omega,f)=I_{K\tilde{M}}^{K\tilde{G},{\cal E}}(exp(X)\eta,\omega,f)+\vert \check{\alpha}\vert log(\vert \alpha(X)\vert )I_{K\underline{\tilde{M}}}^{K\tilde{G},{\cal E}}(exp(X)\eta,\omega,f).$$
 
  \ass{Proposition}{Soient $f\in C_{c}^{\infty}(K\tilde{G}({\mathbb R}))$.  
 
 (i) Pour tout  $U\in Sym(\mathfrak{t}^{\theta})$, les limites 
 
 $lim_{r\to 0+}\partial_{U}I_{K\tilde{M}}^{K\tilde{G},{\cal E},mod}(exp(rH_{d})\eta,\omega,f)$ et $lim_{r\to 0-}\partial_{U}I_{K\tilde{M}}^{K\tilde{G},{\cal E},mod}(exp(rH_{d})\eta,\omega,f)$
 
\noindent existent.
 
 (ii) Si $w_{d}(U)=U$, ces limites sont \'egales.
 
 (iii) Pour tout $\underline{U}\in Sym(\underline{\mathfrak{t}}^{\theta})$, les limites 
 
 $lim_{r\to 0+}\partial_{\underline{U}}I_{K\underline{\tilde{M}}}^{K\tilde{G},{\cal E}}(exp(rH_{c})\eta,\omega,f)$ et $lim_{r\to 0-}\partial_{\underline{U}}I_{K\underline{\tilde{M}}}^{K\tilde{G},{\cal E}}(exp(rH_{c})\eta,\omega,f)$ 
 
 \noindent existent.
 
 (iv) Soit $U\in Sym(\mathfrak{t}^{\theta})$, supposons $w_{d}(U)=-U$ et  posons $\underline{U}=C(U)$. Alors on a les \'egalit\'es
 $$lim_{r\to 0+}\partial_{U}I_{K\tilde{M}}^{K\tilde{G},{\cal E},mod}(exp(rH_{d})\eta,\omega,f)-lim_{r\to 0-}\partial_{U}I_{K\tilde{M}}^{K\tilde{G},{\cal E},mod}(exp(rH_{d})\eta,\omega,f)$$
 $$=2^{1+c(\eta)}\pi i\vert \check{\alpha}\vert \,lim_{r\to 0+}\partial_{\underline{U}}I_{K\underline{\tilde{M}}}^{K\tilde{G},{\cal E}}(exp(rH_{c})\eta,\omega,f)$$
 $$=-2^{1+c(\eta)}\pi i\vert \check{\alpha}\vert \,lim_{r\to 0-}\partial_{\underline{U}}I_{K\underline{\tilde{M}}}^{K\tilde{G},{\cal E}}(exp(rH_{c})\eta,\omega,f).$$}
 
 La preuve occupe les deux paragraphes suivants.
 
 \bigskip
 
 \subsection{Formules d'inversion}

 On note $(\tilde{G}_{p})_{p\in \Pi}$ le $K$-espace $K\tilde{G}$ et $p_{0}$ l'indice tel que $\tilde{G}=\tilde{G}_{p_{0}}$. On pose simplement $\phi_{p}=\phi_{p_{0},p}$, $\tilde{\phi}_{p}=\tilde{\phi}_{p_{0},p}$, $\nabla_{p}=\nabla_{p_{0},p}$, cf. [I] 1.11. Les $K$-espaces $K\tilde{M}$ et $K\underline{\tilde{M}}$ sont index\'es par des sous-ensembles $\Pi^M\subset \Pi^{\underline{M}}\subset \Pi$. On peut supposer que, pour $p\in \Pi^M$, resp. $p\in \Pi^{\underline{M}}$, $\nabla_{p}$ prend ses valeurs dans $M_{SC}$, resp. $\underline{M}_{SC}$, et $\tilde{\phi}_{p}(\tilde{M}_{p})=\tilde{M}$, resp. $\tilde{\phi}_{p}(\underline{\tilde{M}}_{p})=\underline{\tilde{M}}$. On fixe comme toujours une paire de Borel \'epingl\'ee $\hat{{\cal E}}$ de $\hat{G}$ conserv\'ee par l'action galoisienne, dont on note le tore $\hat{T}$ (pour ce qui est des actions galoisiennes, il ne s'agit pas du tore dual de $T$). On fixe un \'el\'ement $\tilde{P}\in {\cal P}(\tilde{M})$ et un \'el\'ement $\underline{\tilde{P}}\in {\cal P}(\underline{\tilde{M}})$ tels que $\tilde{P}\subset \underline{\tilde{P}}$. A l'aide de ces espaces paraboliques, on peut identifier les groupes duaux $\hat{M}$ et $\underline{\hat{M}}$ \`a des Levi de $\hat{G}$ qui sont standard pour $\hat{{\cal E}}$.
 
  Pour tout \'el\'ement semi-simple fortement r\'egulier $\gamma\in \tilde{M}({\mathbb R})$, notons ${\cal X}^M(\gamma)$ l'ensemble des couples   $({\bf  M}',\delta)$, o\`u ${\bf M}'$ est une donn\'ee endoscopique elliptique de $(M,\tilde{M},{\bf a})$ et $\delta\in \tilde{M}'({\mathbb R})$ est un \'el\'ement semi-simple (forc\'ement fortement r\'egulier) qui correspond \`a $\gamma$. Deux tels couples $({\bf M}'_{1},\delta_{1})$ et $({\bf M}'_{2},\delta_{2})$ sont dits \'equivalents s'il existe une \'equivalence entre ${\bf M}'_{1}$ et ${\bf M}'_{2}$, \`a laquelle est associ\'e un isomorphisme $\tilde{\iota}:\tilde{M}'_{1}\to \tilde{M}'_{2}$ d\'efini sur ${\mathbb R}$, de sorte que $\delta_{2}$ soit stablement conjugu\'e \`a $\tilde{\iota}(\delta_{1})$. On fixe un ensemble de repr\'esentants $\dot{{\cal X}}^M(\gamma)$ des classes d'\'equivalence dans ${\cal X}^M(\gamma)$.

 Consid\'erons un \'el\'ement $X_{0}\in \mathfrak{t}({\mathbb R})$ proche de $0$ tel que $exp(X_{0})\eta$ soit fortement r\'egulier dans $\tilde{ M}({\mathbb R})$.   On indexe $\dot{{\cal X}}^M(exp(X_{0})\eta)$ par un ensemble fini $J$: $\dot{{\cal X}}^M(exp(X_{0})\eta)=({\bf M}'_{j},\delta_{j})_{j\in J}$.  On note $T'_{j}$ le commutant de $\delta_{j}$ dans $M_{j}'$ et on pose $\tilde{T}'_{j}=T'_{j}\delta_{j}$. Le choix d'un diagramme  $(\delta_{j},B'_{j},T'_{j},B^M,T,exp(X_{0})\eta)$ d\'etermine un homomorphisme $\xi_{j}:T\to T_{j}'\simeq T/(1-\theta)(T)$. Il s'\'etend en une application compatible $\tilde{\xi}_{j}:\tilde{T}\to \tilde{T}'_{j}$ qui est d\'efinie sur ${\mathbb R}$ et v\'erifie $ \tilde{\xi}_{j}(exp(X_{0})\eta)=\delta_{j}$. On pose $\epsilon_{j}=\tilde{\xi}_{j}(\eta)$. Alors $(\epsilon_{j},B'_{j},T'_{j},B^M,T,\eta)$ est encore un diagramme. On v\'erifie que, pour tout $X\in \mathfrak{t}({\mathbb R})$  proche de $0$ et tel que $exp(X)\eta$ soit fortement r\'egulier dans $\tilde{M}({\mathbb R})$, on peut choisir pour ensemble $\dot{{\cal X}}^M(exp(X)\eta)$ l'ensemble $({\bf M}_{j}',exp(\xi_{j}(X))\epsilon_{j})_{j\in J}$. On fixe pour tout $j$ des donn\'ees auxiliaires $M'_{j,1}$,...,$\Delta_{j,1}$, un \'el\'ement $\epsilon_{j,1}\in \tilde{M}'_{j,1}({\mathbb R})$ au-dessus de $\epsilon_{j}$ et on note $\tilde{T}'_{j,1}$ l'image r\'eciproque de $\tilde{T}'_{j}$ dans $\tilde{M}'_{j,1}$. On fixe aussi une d\'ecomposition
 $$\mathfrak{t}'_{j,1}=\mathfrak{c}_{j,1}\oplus \mathfrak{t}'_{j}$$
 selon la recette de 3.1. 
 Rappelons que $\tilde{T}$ est un sous-tore tordu elliptique de $\tilde{M}({\mathbb R})$. La formule d'inversion [I] 4.9(5) se traduit ainsi. Pour $f\in C_{c}^{\infty}(K\tilde{M}({\mathbb R}))$, on a l'\'egalit\'e
 $$I^{K\tilde{M}}(exp(X)\eta,\omega,f)=[T^{\theta}({\mathbb R}):T^{\theta,0}({\mathbb R})]\vert J\vert ^{-1}d(\theta^*)^{-1/2}\sum_{j\in J}\Delta_{j,1}(exp(\xi_{j}(X))\epsilon_{1,j},exp(X)\eta)^{-1}$$
 $$S^{\tilde{M}'_{j,1}}_{\lambda_{j,1}}(exp(\xi_{j}(X))\epsilon_{1,j},f^{\tilde{M}'_{j,1}}).$$
 Tous les termes sont continus en un point $exp(X)\eta$ qui est seulement r\'egulier dans $\tilde{M}$ (et non fortement r\'egulier). Or, puisque $\eta$ lui-m\^eme est r\'egulier dans $\tilde{M}$, c'est le cas de tout $exp(X)\eta$ pour $X$ voisin de $0$. La formule ci-dessus est donc v\'erifi\'ee en tout $X$ proche de $0$.

On peut  traduire cette formule de la fa\c{c}on suivante. Notons $\boldsymbol{\gamma}(X)$ la distribution $I^{K\tilde{M}}(exp(X)\eta,\omega,.)$ et, pour  tout $j$, $\boldsymbol{\delta}_{j}(X)$ l'\'element de $D_{g\acute{e}om}^{st}({\bf M}'_{j})$ \`a laquelle s'identifie la distribution
$S^{\tilde{M}'_{j,1}}_{\lambda_{j,1}}(exp(\xi_{j}(X))\epsilon_{1,j},.)$, multipli\'ee par $\Delta_{j,1}(exp(\xi_{j}(X))\epsilon_{1,j},exp(X)\eta)^{-1}$. Remarquons que cette distribution est ind\'ependante des donn\'ees auxiliaires $M'_{j,1}$,...,$\Delta_{j,1}$. Alors
$$(1) \qquad \boldsymbol{\gamma}(X)=[T^{\theta}({\mathbb R}):T^{\theta,0}({\mathbb R})]\vert J\vert ^{-1}d(\theta^*)^{-1/2}\sum_{j\in J} transfert(\boldsymbol{\delta}_{j}(X)).$$

Pour tout \'el\'ement semi-simple $\gamma\in \tilde{G}({\mathbb R})$ et tout $p\in \Pi$,  on note ${\cal Y}_{p}(\gamma)$ l'ensemble des $y\in G$ tels que $y \nabla_{p}(\sigma)\sigma(y)^{-1}\in I_{\gamma}$ pour tout $\sigma\in \Gamma_{{\mathbb R}}$ (on rappelle que $I_{\gamma}=Z(G)^{\theta}G_{\gamma}$). Pour $y$ dans cet ensemble, on pose $\gamma[y]=\tilde{\phi}^{-1}_{p}(y^{-1}\gamma y)$. Alors $\gamma[y]\in \tilde{G}_{p}({\mathbb R})$. On fixe un ensemble de repr\'esentants $\dot{{\cal Y}}_{p}(\gamma)$ de l'ensemble des doubles classes $I_{\gamma}\backslash{\cal Y}_{p}(\gamma)/\phi_{p}(G_{p}({\mathbb R}))$. On note ${\cal Y}(\gamma)$, resp. 
$\dot{{\cal Y}}(\gamma)$, la r\'eunion disjointe des ${\cal Y}_{p}(\gamma)$, resp. $\dot{{\cal Y}}_{p}(\gamma)$, pour $p\in \Pi$. Si $\gamma$ est fortement r\'egulier, l'ensemble $\{\gamma[y]; y\in \dot{{\cal Y}}(\gamma)\}$ est un ensemble de repr\'esentants des classes de conjugaison dans la classe de conjugaison stable de $\gamma$. Si $\tilde{L}$ est un espace de Levi de $\tilde{G}$, composante connexe d'un $K$-espace de Levi $K\tilde{L}$ de $K\tilde{G}$, et si 
 $\gamma\in \tilde{L}({\mathbb R})$,  on peut effectuer les m\^emes constructions en rempla\c{c}ant $K\tilde{G}$ par $K\tilde{L}$. On affecte d'un exposant $L$ les objets obtenus, par exemple $\dot{{\cal Y}}^L(\gamma)$. 
 
 Soit comme ci-dessus $X\in \mathfrak{t}^{\theta}({\mathbb R})$ proche de $0$ et tel que $exp(X)\eta$ soit fortement r\'egulier dans $\tilde{M}(F)$. On a $I^M_{exp(X)\eta}=T^{\theta}$.  Les d\'efinitions des ensembles $\dot{{\cal Y}}^M(exp(X)\eta)$ et $\dot{{\cal Y}}_{p}^M(exp(X)\eta)$ ne d\'ependent que de ce groupe. On peut donc les choisir ind\'ependant de $X$ et on les note plut\^ot $\dot{{\cal Y}}^M(\tilde{T})$ et $\dot{{\cal Y}}_{p}^M(\tilde{T})$. Pour un \'el\'ement $y$ de cet ensemble, notons $p$ l'\'el\'ement de $\Pi^M$ tel que $y\in {\cal Y}_{p}^M( \tilde{T})$. On pose $\tilde{T}[y]=\tilde{\phi}_{p}^{-1}\circ ad_{y^{-1}}(\tilde{T})$, $\eta[y]=\tilde{\phi}_{p}^{-1}\circ ad_{y^{-1}}(\eta )$, $X[y]=\phi_{p}^{-1}\circ ad_{y^{-1}}(X)$. Pour $j\in J$, notons $\boldsymbol{\delta}'_{j}(X)$ l'\'element de $D_{g\acute{e}om}^{st}({\bf M}'_{j})$ \`a laquelle s'identifie la distribution
$S^{\tilde{M}'_{j,1}}_{\lambda_{j,1}}(exp(\xi_{j}(X))\epsilon_{1,j},.)$. La formule d'inversion [I] 4.9(4) se traduit par
$$transfert(\boldsymbol{\delta}'_{j}(X))=d(\theta^*)^{1/2}\sum_{y\in \dot{{\cal Y}}^M(\tilde{T})}[T[y]^{\theta}({\mathbb R}):T[y]^{\theta,0}({\mathbb R})]^{-1}$$
$$\Delta_{j,1}(exp(\xi_{j}(X))\epsilon_{j,1},exp(X[y])\eta[y])\boldsymbol{\gamma}(X,y),$$
o\`u $\boldsymbol{\gamma}(X,y)$ est la distribution $I^{K\tilde{M}}(exp(X[y])\eta[y],\omega,.)$. Pour obtenir le transfert de $\boldsymbol{\delta}_{j}(X)$, il suffit de diviser par $\Delta_{j,1}(exp(\xi_{j}(X))\epsilon_{1,j},exp(X)\eta)$. Consid\'erons le rapport
$$\frac{\Delta_{j,1}(exp(\xi_{j}(X))\epsilon_{j,1},exp(X[y])\eta[y])}{\Delta_{j,1}(exp(\xi_{j}(X))\epsilon_{1,j},exp(X)\eta)}.$$
Puisque $exp(\xi_{j}(X))\eta[y]$ et $exp(X)\eta$ sont stablement conjugu\'es, cette conjugaison stable \'etant r\'ealis\'ee par l'\'el\'ement $y$ qui est ind\'ependant de $X$, il r\'esulte de [KS] th\'eor\`eme 5.1.D(1) que ce rapport est d'une part ind\'ependant  de $X$, d'autre part ind\'ependant des donn\'ees auxiliaires $M'_{j,1}$,...,$\Delta_{j,1}$. On le note $d_{j}(y)$. On obtient
$$(2) \qquad transfert(\boldsymbol{\delta}_{j}(X))=d(\theta^*)^{1/2}\sum_{y\in \dot{{\cal Y}}^M(\tilde{T})}[T[y]^{\theta}({\mathbb R}):T[y]^{\theta,0}({\mathbb R})]^{-1} d_{j}(y)\boldsymbol{\gamma}(X,y).$$
On peut \'evidemment supposer que notre ensemble $\dot{{\cal Y}}^M(\tilde{T})$ contient l'\'el\'ement $y=1$. En comparant les formules (1) et (2), on obtient pour $y\in \dot{{\cal Y}}^M(\tilde{T})$ l'\'egalit\'e
$$(3) \qquad \vert J\vert ^{-1}\sum_{j\in J}d_{j}(y)=\left\lbrace\begin{array}{cc}1,&\text{ si }y=1,\\ 0,&\text{ si }y\not=1.\\  \end{array}\right.$$

 Fixons $j\in J$. On note ${\bf M}'_{j}=(M'_{j},{\cal M}'_{j},\tilde{\zeta}_{j})$. Pour $\tilde{s}\in\tilde{\zeta}_{j} Z(\hat{M})^{\Gamma_{{\mathbb R}},\hat{\theta}}/Z(\hat{G})^{\Gamma_{{\mathbb R}},\hat{\theta}}$, on introduit des donn\'ees suppl\'ementaires $G'_{1}(\tilde{s})$,...,$\Delta_{1}(\tilde{s})$ pour la donn\'ee endoscopique ${\bf G}'(\tilde{s})$. On fixe un point $\epsilon_{1}(\tilde{s})$ au-dessus de $\epsilon_{j}$ et on note  $\tilde{M}'_{1}(\tilde{s})$ et $\tilde{T}'_{1}(\tilde{s})$ les images r\'eciproques de $\tilde{M}'_{j}$ et $\tilde{T}'_{j}$ dans $\tilde{G}'_{1}(\tilde{s})$. On fixe aussi une d\'ecomposition
$$\mathfrak{t}'_{1}(\tilde{s})=\mathfrak{c}_{1}(\tilde{s})\oplus \mathfrak{t}'_{j}$$
selon la recette de 3.1.  On a fix\'e une paire de Borel \'epingl\'ee $\hat{{\cal E}}$ de $\hat{G}$, dont  on a  not\'e le tore $\hat{T}$. On peut supposer $\tilde{\zeta}_{j}=\zeta_{j}\hat{\theta}$ et on \'ecrit tout \'el\'ement $\tilde{s}$ sous la forme $s\hat{\theta}$.  On a  fix\'e un diagramme $(\epsilon_{j},B'_{j},T'_{j},B^M,T,\eta) $ reliant $\epsilon_{j}$ et $\eta$, d'espaces ambiants $\tilde{M}$ et $\tilde{M}'_{j}$. On a aussi fix\'e un \'el\'ement $\tilde{P}\in {\cal P}(\tilde{M})$. Notons $B$ le sous-groupe de   Borel   de $G$ tel que $B\subset P$ et $B\cap M=B^M$.  On v\'erifie que le diagramme s'\'etend pour tout $\tilde{s}$ en un diagramme $(\epsilon_{j},B'(\tilde{s}),T'_{j},B,T,\eta)$, d'espaces ambiants $\tilde{G}$ et $\tilde{G}'(\tilde{s})$. On compl\`ete $(B,T)$ en une paire de Borel \'epingl\'ee ${\cal E}$ de $G$. On fixe un \'el\'ement $e\in Z(\tilde{G},{\cal E})$ et on \'ecrit $\eta=\nu e$, avec $\nu\in T$. Le groupe $G_{\eta}$ a une unique racine (au signe pr\`es) que l'on a not\'ee $\alpha$. C'est un \'el\'ement de $X^*(T^{\theta,0})$.   Rappelons que l'on  note $\Sigma(T)$ l'ensemble des racines de $T$ dans $G$. Introduisons un \'el\'ement $\beta\in \Sigma(T)$ dont la restriction $\beta_{res}$ \`a $T^{\theta,0}$ soit $\alpha$. La coracine   $\check{\alpha}$ est alors $N(\check{\beta})$ si $\beta$ est de type 1 ou 3, $2N(\check{\beta})$ si $\beta$ est de type 2. D'apr\`es [W1] 3.3, l'ensemble des racines de $G_{\eta}$ est form\'e des $\beta'_{res}$, pour $\beta'\in \Sigma(T)$ tels que $N(\beta')(\nu)=1$ si $\beta'$ est de type 1 ou 2, $N(\beta')(\nu)=-1$ si $\beta'$ est de type 3. Notre hypoth\`ese est donc que cet ensemble de racines se r\'eduit aux racines de la forme $\pm \theta^k(\beta)$ pour $k\in {\mathbb N}$. Il correspond \`a $\beta$ une racine $\hat{\beta}$ de $\hat{T}$. De la description de l'ensemble de racines de $G'(\tilde{s})_{\epsilon_{j}}$ donn\'ee en [W1] 3.3 se d\'eduisent les r\'esultats suivants. Le groupe $G'(\tilde{s})_{\epsilon_{j}}$ est un tore ou de rang semi-simple $1$. Il est de rang semi-simple $1$ si et seulement si l'une des conditions suivantes est v\'erifi\'ee:

(a) $\beta$ de type 1 et $N(\hat{\beta})(s)=1$;

(b) $\beta$ de type 2 et $N(\hat{\beta})(s)=1$:

(c) $\beta$ de type 2 et $N(\hat{\beta})(s)=-1$;

(d) $\beta$ de type 3 et $N(\hat{\beta})(s)=1$.

L'unique racine $\alpha(\tilde{s})$ et l'unique coracine $\check{\alpha}(\tilde{s})$ de $G'(\tilde{s})_{\epsilon_{j}}$ (au signe pr\`es) se d\'ecrivent  dans chacun des cas de la fa\c{c}on suivante:

(a) $\alpha(\tilde{s})\circ \xi_{j}=N(\beta) $, $\check{\alpha}(\tilde{s})=\xi_{j}\circ \check{\beta}$;

(b) $\alpha(\tilde{s})\circ \xi_{j}=2N(\beta)$, $\check{\alpha}(\tilde{s})=\xi_{j}\circ \check{\beta}$;

(c) $\alpha(\tilde{s})\circ \xi_{j}=N(\beta)$, $\check{\alpha}(\tilde{s})=2\xi_{j}\circ \check{\beta}$;

(d) $\alpha(\tilde{s})\circ \xi_{j}=2N(\beta)$, $\check{\alpha}(\tilde{s})=\frac{1}{2}\xi_{j}\circ \check{\beta}$.

Puisque $\beta$ se restreint en un \'el\'ement non nul de $X_{*}(A_{\tilde{M}})$, $N\hat{\beta}$ se restreint en un caract\`ere non trivial de $Z(\hat{M})^{\Gamma_{{\mathbb R}},\hat{\theta},0}$. Quitte \`a multiplier $\tilde{\zeta}_{j}$ par un \'el\'ement  de ce groupe, on peut supposer $N(\hat{\beta})(\zeta_{j})=1$.  Posons $\hat{\alpha}^*=N(\hat{\beta})$ si $\beta$ est de type 1 ou 3, $\hat{\alpha}^*=2N(\hat{\beta})$ si $\beta$ est de type 2. 
Notons $Z(\hat{M})^{\hat{\alpha}^*}$ l'ensemble des $ t\in Z(\hat{M})^{\Gamma_{{\mathbb R}},\hat{\theta}}$ tels que $\hat{\alpha}^*(t)=1$. Pour $\tilde{s}\in\tilde{\zeta}_{j} Z(\hat{M})^{\Gamma_{{\mathbb R}},\hat{\theta}}/Z(\hat{G})^{\Gamma_{{\mathbb R}},\hat{\theta}}$, le groupe $G'(\tilde{s})_{\epsilon_{j}}$ est de rang semi-simple $1$ si et seulement si $\tilde{s}\in \tilde{\zeta}_{j}Z(\hat{M})^{\hat{\alpha}^*}/Z(\hat{G})^{\Gamma_{{\mathbb R}},\hat{\theta}}$. Soit $\tilde{s}\in \tilde{\zeta}_{j}Z(\hat{M})^{\hat{\alpha}^*}/Z(\hat{G})^{\Gamma_{{\mathbb R}},\hat{\theta}}$. Remarquons que, puisque $\beta_{res}$ est fixe par $\Gamma_{{\mathbb R}}$, il en est de m\^eme de la racine $ \alpha(\tilde{s})$   de $G'(\tilde{s})_{\epsilon_{j}}$. Donc $G'(\tilde{s})_{\epsilon_{j},SC}$ est isomorphe \`a $SL(2)$. On note $\underline{\tilde{M}}'(\tilde{s})$ l'espace de Levi de $\tilde{G}'(\tilde{s})$ tel que ${\cal A}_{\underline{\tilde{M}}'(\tilde{s})}$ soit le noyau de la racine $\alpha(\tilde{s})$ vue comme forme lin\'eaire sur ${\cal A}_{\tilde{M}'_{j}}$. On note $\underline{\tilde{M}}'_{1}(\tilde{s})$ son image r\'eciproque dans $\tilde{G}'_{1}(\tilde{s})$. On voit que $\underline{\tilde{M}}'(\tilde{s})$ n'est autre que l'espace endoscopique de la donn\'ee endoscopique $\underline{{\bf M}}'(\tilde{s})$ de $(\underline{M},\underline{\tilde{M}},{\bf a})$ associ\'ee \`a $\tilde{s}$. Cette donn\'ee est elliptique. En effet, on a les inclusions
$${\cal A}_{\underline{\tilde{M}}}\subset {\cal A}_{\underline{\tilde{M}}'(\tilde{s})}\subsetneq {\cal A}_{\tilde{M}'_{j}}={\cal A}_{\tilde{M}}.$$
Puisque la diff\'erence des dimensions des deux espaces extr\^emes est $1$, la premi\`ere inclusion est une \'egalit\'e, ce qui prouve l'ellipticit\'e.
D'autre part, d'apr\`es la d\'efinition de $\underline{\tilde{M}}'(\tilde{s})$, on a $\underline{M}'(\tilde{s})_{\epsilon_{j}}=G'(\tilde{s})_{\epsilon_{j}}$.

 Soit $\tilde{s}\in \tilde{\zeta}_{j}Z(\hat{M})^{\hat{\alpha}^*}$. Utilisons la th\'eorie locale rappel\'ee en [III] section 5.   On introduit la donn\'ee endoscopique $\bar{{\bf G}}'(\bar{s})$ du groupe $G_{\eta,SC}$. La paire $(\bar{G}'(\bar{s})_{SC},G'(\tilde{s})_{\epsilon_{j},SC} )$ se compl\`ete en un triplet endoscopique non standard.  Mais on sait que $G_{\eta,SC}$ et $G'(\tilde{s})_{\epsilon_{j},SC}$ sont tous deux isomorphes \`a $SL(2)$. Il en r\'esulte que $\bar{G}'(\bar{s})$ est lui-m\^eme simplement connexe et isomorphe \`a $SL(2)$. Le triplet endoscopique non standard est \'equivalent \`a un triplet "trivial". Notons toutefois que l'isomorphisme $j_{*}$ entre les alg\`ebres de Lie des tores d\'eploy\'es maximaux de ces groupes est en g\'en\'eral un multiple de l'isomorphisme naturel entre ces deux alg\`ebres (la notation $j_{*}$ se r\'ef\`ere \`a [III] 6.1, elle n'a rien \`a voir avec notre \'el\'ement $j\in J$). La donn\'ee $\bar{{\bf G}}'(\bar{s})$ est la donn\'ee endoscopique maximale  de $G_{\eta,SC}$. On peut prendre pour cette donn\'ee le facteur de transfert qui vaut $1$ sur des couples d'\'el\'ements qui se correspondent. De l'homomorphisme $\xi_{j}:T\to T'_{j}$ se d\'eduit un homomorphisme $\xi_{j}:\mathfrak{z}(G_{\eta})\to \mathfrak{z}(G'(\tilde{s})_{\epsilon_{j}})$. La relation [V] 4.1(1) se traduit  par l'existence d'une constante non nulle $d(\tilde{s})$ v\'erifiant la condition  suivante. Soit $X_{sc}\in \mathfrak{g}_{\eta,SC}({\mathbb R})$ un \'el\'ement r\'egulier. On le transf\`ere en un \'el\'ement $\bar{X}_{sc}\in \bar{\mathfrak{g}}'(\bar{s})({\mathbb R})$ et on transf\`ere celui-ci en un \'el\'ement $X'_{sc}\in \mathfrak{g}'(\tilde{s})_{\epsilon_{j},SC}({\mathbb R})$. On suppose ces \'el\'ements proches de $0$. Soit $Z\in \mathfrak{z}(G_{\eta},{\mathbb R})$ en position g\'en\'erale. Alors
 $$(4) \qquad     lim_{Z\to 0}\Delta_{1}(\tilde{s})(exp(\xi_{j}(Z)+X'_{sc})\epsilon_{1}(\tilde{s}),exp(Z+X_{sc})\eta)=d(\tilde{s}).$$
 Plus g\'en\'eralement, soit  $p\in \Pi$ et $y\in {\cal Y}_{p}(\eta)$. L'application $\phi_{p}^{-1}\circ ad_{y^{-1}}$ se restreint en un torseur int\'erieur de $G_{\eta}$ sur $G_{\eta[y]}$. Il se restreint en un isomorphisme d\'efini sur ${\mathbb R}$ de $\mathfrak{z}(G_{\eta})$ sur $\mathfrak{z}(G_{\eta[y]})$. On note $Z\mapsto Z[y]$ cet isomorphisme.  Le groupe $G_{\eta[y],SC}$ n'est plus en g\'en\'eral isomorphe \`a $SL(2)$, c'en est une forme int\'erieure. La donn\'ee $\bar{{\bf G}}'(\bar{s})$ est encore la donn\'ee endoscopique maximale de $G_{\eta[y],SC}$. On peut prendre pour cette donn\'ee le facteur de transfert qui vaut $1$ sur tout couple d'\'el\'ements qui se correspondent. La relation [V] 4.1(1) implique l'existence d'une constante non nulle $d(\tilde{s},y)$ v\'erifiant la condition  suivante. Soit $X_{sc}[y]\in \mathfrak{g}_{\eta[y],SC}({\mathbb R})$ un \'el\'ement r\'egulier. On le transf\`ere en un \'el\'ement $\bar{X}_{sc}\in \bar{\mathfrak{g}}'(\bar{s})({\mathbb R})$ et on transf\`ere celui-ci en un \'el\'ement $X'_{sc}\in \mathfrak{g}'(\tilde{s})_{\epsilon_{j},SC}({\mathbb R})$. On suppose ces \'el\'ements proches de $0$. Soit $Z\in \mathfrak{z}(G_{\eta},{\mathbb R})$ en position g\'en\'erale. Alors
 $$(5) \qquad   lim_{Z\to 0}\Delta_{1}(\tilde{s})(exp(\xi_{j}(Z)+X'_{sc})\epsilon_{1}(\tilde{s}),exp(Z[y]+X_{sc}[y])\eta)=d(\tilde{s},y).$$
 La constante $d(\tilde{s})$ introduite ci-dessus n'est autre que $d(\tilde{s},1)$. 
 
 On a introduit au d\'ebut du paragraphe des donn\'ees auxiliaires $M'_{j,1}$,...,$\Delta_{j,1}$. On peut prendre pour celles-ci les donn\'ees $M'_{1}(\tilde{s})$,...,$\Delta_{1}(\tilde{s})$ (quant \`a $\Delta_{1}(\tilde{s})$, il s'agit plut\^ot ici de la restriction de ce facteur de transfert aux couples d'\'el\'ements de $\tilde{M}'_{1}(\tilde{s},{\mathbb R})\times \tilde{M}({\mathbb R})$ qui se correspondent).  Dans le cas o\`u $y\in \dot{{\cal Y}}^M(\eta)$, la constante $d_{j}(y)$ d\'efinie plus haut n'est autre que $d(\tilde{s},y)d(\tilde{s},1)^{-1}$. Elle est ind\'ependante du choix de $\tilde{s}$.

\bigskip

\subsection{Preuve de la proposition 4.3}

Soit $X\in \mathfrak{t}^{\theta}({\mathbb R})\cap \mathfrak{g}_{\eta,reg}({\mathbb R})$.  Consid\'erons la formule 4.4(1). Toutes les distributions y intervenant sont $\tilde{G}$-\'equisinguli\`eres. On a alors par d\'efinition
$$(1) \qquad I_{K\tilde{M}}^{K\tilde{G},{\cal E}}(exp(X)\eta,\omega,f)=I_{K\tilde{M}}^{K\tilde{G},{\cal E}}(\boldsymbol{\gamma}(X),f)=[T^{\theta}({\mathbb R}):T^{\theta,0}({\mathbb R})]\vert J\vert ^{-1}d(\theta^*)^{-1/2}$$
$$\sum_{j\in J} I_{K\tilde{M}}^{K\tilde{G},{\cal E}}({\bf M}'_{j},\boldsymbol{\delta}_{j}(X),f)$$
pour tout $f\in C_{c}^{\infty}(K\tilde{G}({\mathbb R}))$. 

Fixons $j\in J$ et utilisons les notations introduites dans le paragraphe pr\'ec\'edent. Pour tout $\tilde{s}\in \tilde{\zeta}_{j}Z(\hat{M})^{\Gamma_{{\mathbb R}},\hat{\theta}}/Z(\hat{G})^{\Gamma_{{\mathbb R}},\hat{\theta}}$, la distribution $\boldsymbol{\delta}_{j}(X)$ s'identifie \`a
$$\Delta_{1}(\tilde{s})(exp(\xi_{j}(X))\epsilon_{1}(\tilde{s}),exp(X)\eta)^{-1}S^{\tilde{M}'_{1}(\tilde{s})}_{\lambda_{1}(\tilde{s})}(exp(\xi_{j}(X))\epsilon_{1}(\tilde{s}),.).$$
On a donc
$$(2)\qquad I_{K\tilde{M}}^{K\tilde{G},{\cal E}}({\bf M}'_{j},\boldsymbol{\delta}_{j}(X),f)=\sum_{\tilde{s}\in \tilde{\zeta}_{j}Z(\hat{M})^{\Gamma_{{\mathbb R}},\hat{\theta}}/Z(\hat{G})^{\Gamma_{{\mathbb R}},\hat{\theta}}}i_{\tilde{M}'_{j}}(\tilde{G},\tilde{G}'(\tilde{s}))$$
$$\Delta_{1}(\tilde{s})(exp(\xi_{j}(X))\epsilon_{1}(\tilde{s}),exp(X)\eta)^{-1}S_{\tilde{M}'_{1}(\tilde{s}),\lambda_{1}(\tilde{s})}^{\tilde{G}'_{1}(\tilde{s})}(exp(\xi_{j}(X))\epsilon_{1}(\tilde{s}),f^{\tilde{G}'_{1}(\tilde{s})}).$$
Soit $\tilde{s}\in \tilde{\zeta}_{j}Z(\hat{M})^{\Gamma_{{\mathbb R}},\hat{\theta}}/Z(\hat{G})^{\Gamma_{{\mathbb R}},\hat{\theta}}$, supposons d'abord $\tilde{s}\not\in \tilde{\zeta}_{j}Z(\hat{M})^{\hat{\alpha}^*}/Z(\hat{G})^{\Gamma_{{\mathbb R}},\hat{\theta}}$. Alors $G'(\tilde{s})_{\epsilon_{j}}$ est un tore et
$\epsilon_{j}$ est r\'egulier dans $\tilde{G}'(\tilde{s})$. Il en r\'esulte que la fonction
$$X\mapsto \Delta_{1}(\tilde{s})(exp(\xi_{j}(X))\epsilon_{1}(\tilde{s}),exp(X)\eta)^{-1}S_{\tilde{M}'_{1}(\tilde{s}),\lambda_{1}(\tilde{s})}^{\tilde{G}'_{1}(\tilde{s})}(exp(\xi_{j}(X))\epsilon_{1}(\tilde{s}),f^{\tilde{G}'_{1}(\tilde{s})})$$
est $C^{\infty}$ au voisinage de $0$. Supposons maintenant  $\tilde{s}\in \tilde{\zeta}_{j}Z(\hat{M})^{\hat{\alpha}^*}/Z(\hat{G})^{\Gamma_{{\mathbb R}},\hat{\theta}}$. Le groupe $G'(\tilde{s})_{\epsilon_{j},SC}=\underline{M}'(\tilde{s})_{\epsilon_{j},SC}$ est isomorphe \`a $SL(2)$. On peut introduire l'int\'egrale orbitale pond\'er\'ee stable modifi\'ee
$$(3) \qquad S_{\tilde{M}'_{1}(\tilde{s}),\lambda_{1}(\tilde{s})}^{\tilde{G}'_{1}(\tilde{s}),mod}(exp(\xi_{j}(X))\epsilon_{1}(\tilde{s}),f^{\tilde{G}'_{1}(\tilde{s})})=S_{\tilde{M}'_{1}(\tilde{s}),\lambda_{1}(\tilde{s})}^{\tilde{G}'_{1}(\tilde{s})}(exp(\xi_{j}(X))\epsilon_{1}(\tilde{s}),f^{\tilde{G}'_{1}(\tilde{s})})$$
$$+i_{\tilde{M}'_{1}(\tilde{s})}^{\tilde{G}'(\tilde{s})}(\epsilon_{1})\vert \check{\alpha}(\tilde{s})\vert log(\vert \alpha(\tilde{s})(\xi_{j}(X))\vert )S_{\underline{\tilde{M}}'_{1}(\tilde{s}),\lambda_{1}(\tilde{s})}^{\tilde{G}'_{1}(\tilde{s})}(exp(\xi_{j}(X))\epsilon_{1},f^{\tilde{G}'_{1}(\tilde{s})}).$$
D'apr\`es les formules de 4.4, on a $\alpha(\tilde{s})(\xi_{j}(X))=n(\tilde{s})\beta(X)=n(\tilde{s})\alpha(X)$, o\`u $n(\tilde{s})=n_{\beta}$ ou $2n_{\beta}$ selon le cas. Donc $ log(\vert \alpha(\tilde{s})(\xi_{j}(X))\vert )=log(n(\tilde{s}))+log(\vert \alpha(X)\vert )$. Le point $\epsilon_{j}$, consid\'er\'e comme un \'el\'ement de $\underline{\tilde{M}}'(\tilde{s};{\mathbb R})$, est $\tilde{G}'(\tilde{s})$-\'equisingulier. Il existe donc une fonction $\varphi\in C_{c}^{\infty}(\underline{\tilde{M}}'_{1}(\tilde{s},{\mathbb R}))$ telle que
$$S_{\underline{\tilde{M}}'_{1}(\tilde{s}),\lambda_{1}(\tilde{s})}^{\tilde{G}'_{1}(\tilde{s})}(exp(\xi_{j}(X))\epsilon_{1},f^{\tilde{G}'_{1}(\tilde{s})})=S^{\underline{\tilde{M}}'_{1}(\tilde{s})}_{\lambda_{1}(\tilde{s})}(exp(\xi_{j}(X))\epsilon_{1},\varphi)$$
pour tout $X$ proche de $0$. L'\'el\'ement $\epsilon_{j}$ n'est pas r\'egulier dans $\underline{\tilde{M}}(\tilde{s})$, mais la seule racine de $T'_{j}$ singuli\`ere en $\epsilon_{j}$ est "r\'eelle". Comme on le sait, cela entra\^{\i}ne que l'int\'egrale orbitale de droite ci-dessus est $C^{\infty}$ pour $X$ proche de $0$.  Donc remplacer $log(\vert \alpha(\tilde{s})(\xi_{j}(X))\vert )$ par $log(\vert \alpha(X)\vert )$ dans la formule (3) modifie le membre de droite par une fonction $C^{\infty}$ au voisinage de $X=0$.  
 Le facteur de transfert
$$X\mapsto  \Delta_{1}(\tilde{s})(exp(\xi_{j}(X))\epsilon_{1}(\tilde{s}),exp(X)\eta)^{-1}$$ 
est lui aussi $C^{\infty}$ au voisinage de $0$. En effet, $exp(X)\eta$ appartenant \`a $\tilde{M}({\mathbb R})$, il s'agit d'un facteur de transfert pour une donn\'ee auxiliaire de la donn\'ee endoscopique ${\bf M}'_{j}$. Il est $C^{\infty}$ car $\eta$ est r\'egulier dans $\tilde{M}({\mathbb R})$. Ces consid\'erations transforment l'expression (2) sous la forme suivante.
$$(4) \qquad  I_{K\tilde{M}}^{K\tilde{G},{\cal E}}({\bf M}'_{j},\boldsymbol{\delta}_{j}(X),f)=\phi_{j}(X)+B_{j}(X)-log(\vert \alpha(X)\vert )D_{j}(X)$$
o\`u $\phi_{j}$ est une fonction $C^{\infty}$ en $X=0$,
$$B_{j}(X)=\sum_{\tilde{s}\in\tilde{\zeta}_{j} Z(\hat{M})^{\hat{\alpha}^*}/Z(\hat{G})^{\Gamma_{{\mathbb R}},\hat{\theta}}}i_{\tilde{M}'_{j}}(\tilde{G},\tilde{G}'(\tilde{s}))\Delta_{1}(\tilde{s})(exp(\xi_{j}(X))\epsilon_{1}(\tilde{s}),exp(X)\eta)^{-1}$$
$$S_{\tilde{M}'_{1}(\tilde{s}),\lambda_{1}(\tilde{s})}^{\tilde{G}'_{1}(\tilde{s}),mod}(exp(\xi_{j}(X))\epsilon_{1}(\tilde{s}),f^{\tilde{G}'_{1}(\tilde{s})}),$$
$$D_{j}(X)=\sum_{\tilde{s}\in\tilde{\zeta}_{j} Z(\hat{M})^{\hat{\alpha}^*}/Z(\hat{G})^{\Gamma_{{\mathbb R}},\hat{\theta}}}i_{\tilde{M}'_{j}}(\tilde{G},\tilde{G}'(\tilde{s}))i_{\tilde{M}'_{1}(\tilde{s})}^{\tilde{G}'(\tilde{s})}(\epsilon_{1})\vert \check{\alpha}(\tilde{s})\vert\Delta_{1}(\tilde{s})(exp(\xi_{j}(X))\epsilon_{1}(\tilde{s}),exp(X)\eta)^{-1}$$
$$S_{\underline{\tilde{M}}'_{1}(\tilde{s}),\lambda_{1}(\tilde{s})}^{\tilde{G}'_{1}(\tilde{s})}(exp(\xi_{j}(X))\epsilon_{1}(\tilde{s}),f^{\tilde{G}'_{1}(\tilde{s})}).$$
On r\'ecrit
$$D_{j}(X)=\sum_{\tilde{t}\in\tilde{\zeta}_{j} Z(\hat{M})^{\hat{\alpha}^*}/Z(\underline{\hat{M}})^{\Gamma_{{\mathbb R}},\hat{\theta}}}D_{j}(X,\tilde{t}),$$
o\`u
$$D_{j}(X,\tilde{t})=\sum_{\tilde{s}\in \tilde{t}\in Z(\underline{\hat{M}})^{\Gamma_{{\mathbb R}},\hat{\theta}}/Z(\hat{G})^{\Gamma_{{\mathbb R}},\hat{\theta}}}i_{\tilde{M}'_{j}}(\tilde{G},\tilde{G}'(\tilde{s}))i_{\tilde{M}'_{1}(\tilde{s})}^{\tilde{G}'(\tilde{s})}(\epsilon_{1})\vert \check{\alpha}(\tilde{s})\vert\Delta_{1}(\tilde{s})(exp(\xi_{j}(X))\epsilon_{1}(\tilde{s}),exp(X)\eta)^{-1}$$
$$S_{\underline{\tilde{M}}'_{1}(\tilde{s}),\lambda_{1}(\tilde{s})}^{\tilde{G}'_{1}(\tilde{s})}(exp(\xi_{j}(X))\epsilon_{1}(\tilde{s}),f^{\tilde{G}'_{1}(\tilde{s})}).$$
Fixons $\tilde{t}\in\tilde{\zeta}_{j} Z(\hat{M})^{\hat{\alpha}^*}/Z(\underline{\hat{M}})^{\Gamma_{{\mathbb R}},\hat{\theta}}$. Cet \'el\'ement d\'etermine une donn\'ee endoscopique $\underline{{\bf M}}'(\tilde{t})$ de $(\underline{M},\underline{\tilde{M}},{\bf a})$, qui est elliptique comme on l'a dit en 4.4.   Soit $\tilde{s}\in\tilde{t} Z(\underline{\hat{M}})^{\Gamma_{{\mathbb R}},\hat{\theta}}/Z(\hat{G})^{\Gamma_{{\mathbb R}},\hat{\theta}}$. Alors l'espace $\underline{\tilde{M}}'_{1}(\tilde{s})$ fait partie de donn\'ees auxiliaires pour cette donn\'ee $\underline{{\bf M}}'(\tilde{t})$.  Posons
$$m=[Z(\hat{M})^{\hat{\alpha}^*}:Z(\underline{\hat{M}})^{\Gamma_{{\mathbb R}},\hat{\theta}}].$$
Montrons que
$$(5) \qquad i_{\tilde{M}'_{j}}(\tilde{G},\tilde{G}'(\tilde{s}))i_{\tilde{M}'_{1}(\tilde{s})}^{\tilde{G}_{1}'(\tilde{s})}(\epsilon_{1}(\tilde{s}))\vert \check{\alpha}(\tilde{s})\vert=m^{-1}i_{\underline{\tilde{M}}'(\tilde{t})}(\tilde{G},\tilde{G}'(\tilde{s}))\vert \check{\alpha}\vert .$$
Rappelons que $ i_{\tilde{M}'_{j}}(\tilde{G},\tilde{G}'(\tilde{s}))$ est l'inverse du nombre d'\'el\'ements du noyau $K_{1}$ de l'homomorphisme
$$ Z(\hat{M})^{\Gamma_{{\mathbb R}},\hat{\theta}}/Z(\hat{G})^{\Gamma_{{\mathbb R}},\hat{\theta}}\to Z(\hat{M}'_{j})^{\Gamma_{{\mathbb R}}}/Z(\hat{G}'(\tilde{s}))^{\Gamma_{{\mathbb R}}}.$$
Rappelons qu'un sous-tore maximal de $\hat{M}'_{j}$ s'identifie \`a $\hat{T}^{\hat{\theta},0}$. Notons $\hat{\alpha}'$ la restriction \`a ce tore de $N(\hat{\beta})$ si $\beta$ est de type 1 ou 3, de $2N(\hat{\beta})$ si $\beta$ est de type 2. Notons $Z(\hat{M}'_{j})^{\hat{\alpha}'}$ le groupe des $x\in  Z(\hat{M}'_{j})^{\Gamma_{{\mathbb R}}}$ tels que $\hat{\alpha}'(x)=1$. Alors
 $Z(\hat{M})^{\hat{\alpha}^*}/Z(\hat{G})^{\Gamma_{{\mathbb R}},\hat{\theta}}$ est l'image r\'eciproque par l'homomorphisme pr\'ec\'edent du groupe $Z(\hat{M}'_{j})^{\hat{\alpha}'}/Z(\hat{G}'(\tilde{s}))^{\Gamma_{{\mathbb R}}}$. On   a donc une suite exacte
 $$1\to K_{1}\to Z(\hat{M})^{\hat{\alpha}^*}/Z(\hat{G})^{\Gamma_{{\mathbb R}},\hat{\theta}}\to Z(\hat{M}'_{j})^{\hat{\alpha}'}/Z(\hat{G}'(\tilde{s}))^{\Gamma_{{\mathbb R}}}\to 1.$$
  Cette suite s'ins\`ere dans un diagramme
 $$\begin{array}{ccccccccc}&&1&&1&&1&&\\ &&\uparrow&&\uparrow&&\uparrow&&\\
 1&\to& K_{2}&\to &Z(\hat{M})^{\hat{\alpha}^*}/Z(\underline{\hat{M}})^{\Gamma_{{\mathbb R}},\hat{\theta}}&\to&Z(\hat{M}'_{j})^{\hat{\alpha}'}/Z(\underline{\hat{M}}'(\tilde{t}))^{\Gamma_{{\mathbb R}}}&\to 1\\&&\uparrow&&\uparrow&&\uparrow&&\\ 1&\to& K_{1}&\to& Z(\hat{M})^{\hat{\alpha}^*}/Z(\hat{G})^{\Gamma_{{\mathbb R}},\hat{\theta}}&\to& Z(\hat{M}'_{j})^{\hat{\alpha}'}/Z(\hat{G}'(\tilde{s}))^{\Gamma_{{\mathbb R}}}&\to& 1\\ &&\uparrow&&\uparrow&&\uparrow&&\\ 1&\to&K_{3}&\to&Z(\underline{\hat{M}})^{\Gamma_{{\mathbb R}},\hat{\theta}}/Z(\hat{G})^{\Gamma_{{\mathbb R}},\hat{\theta}}&\to&Z(\underline{\hat{M}}'(\tilde{t}))^{\Gamma_{{\mathbb R}}}/Z(\hat{G}'(\tilde{s}))^{\Gamma_{{\mathbb R}}}&\to&1\\ &&\uparrow&&\uparrow&&\uparrow&&\\ &&1&&1&&1&&\\ \end{array}$$
 Les groupes $K_{2}$ et $K_{3}$ sont d\'efinis de sorte que les suites horizontales soient exactes. Les deux suites verticales de droite sont exactes. Donc aussi celle de gauche. D'o\`u l'\'egalit\'e
 $$\vert K_{1}\vert =\vert K_{2}\vert \vert K_{3}\vert .$$
 Par d\'efinition, $\vert K_{3}\vert ^{-1}=i_{\underline{\tilde{M}}'(\tilde{t})}(\tilde{G},\tilde{G}'(\tilde{s}))$.  D'o\`u les \'egalit\'es
 $$ i_{\tilde{M}'_{j}}(\tilde{G},\tilde{G}'(\tilde{s}))=\vert K_{1}\vert ^{-1}=  \vert K_{2}\vert ^{-1}i_{\underline{\tilde{M}}'(\tilde{t})}(\tilde{G},\tilde{G}'(\tilde{s})).$$
 La premi\`ere suite horizontale ci-dessus est form\'ee de groupes finis. Donc $\vert K_{2}\vert $ est le quotient du nombre d'\'el\'em\'ent du groupe central de cette suite par le nombre d'\'el\'ements du groupe de droite. Le nombre d'\'el\'ements du groupe central est $m$.    D'o\`u l'\'egalit\'e
$$(6) \qquad  i_{\tilde{M}'_{j}}(\tilde{G},\tilde{G}'(\tilde{s}))= m ^{-1}[Z(\hat{M}'_{j})^{\hat{\alpha}'}:Z(\underline{\hat{M}}'(\tilde{t}))^{\Gamma_{{\mathbb R}}}]i_{\underline{\tilde{M}}'(\tilde{t})}(\tilde{G},\tilde{G}'(\tilde{s})).$$ 
  On a introduit deux sous-groupes de $Z(\hat{M}'_{j})^{\Gamma_{{\mathbb R}}}$: les sous-groupes $Z(\hat{M}'_{j})^{\hat{\alpha}(\tilde{s})}$ et $Z(\hat{M}'_{j})^{\hat{\alpha}'}$. On calcule la racine $\hat{\alpha}(\tilde{s})$: les formules pour les racines dans les groupes duaux sont identiques \`a celles pour les coracines dans les groupes sur ${\mathbb R}$. A l'aide des descriptions de 4.4, on obtient dans chaque cas
 
 (a) $\hat{\alpha}(\tilde{s})=\hat{\beta}_{res}$, $\hat{\alpha}'=n_{\alpha}\hat{\beta}_{res}$ (l'indice $res$ d\'esignant la restriction \`a $\hat{T}^{\hat{\theta},0}$);
 
 (b) $\hat{\alpha}(\tilde{s})=\hat{\beta}_{res}$, $\hat{\alpha}'=2n_{\alpha}\hat{\beta}_{res}$;
 
 (c) $\hat{\alpha}(\tilde{s})=2\hat{\beta}_{res}$, $\hat{\alpha}'=2n_{\alpha}\hat{\beta}_{res}$;
 
 (d) $\hat{\alpha}(\tilde{s})=\frac{1}{2}\hat{\beta}_{res}$, $\hat{\alpha}'=n_{\alpha}\hat{\beta}_{res}$.
 
 En comparant avec les formules donn\'ees  en 4.4 pour les coracines, on obtient l'\'egalit\'e
 $$\hat{\alpha}'=\vert \check{\alpha}\vert \vert \check{\alpha}(\tilde{s})\vert ^{-1} \hat{\alpha}(\tilde{s}).$$
 On a utilis\'e ici la compatibilit\'e de nos diff\'erentes normes: elles proviennent de formes quadratiques sur $X_{*}(T^{\theta,0})\otimes_{{\mathbb Z}}{\mathbb R}$, resp. $X_{*}(T_{j})\otimes_{{\mathbb Z}}{\mathbb R}$, qui se correspondent par l'isomorphisme naturel entre ces espaces.  On d\'eduit de ces calculs l'\'egalit\'e
 $$[Z(\hat{M}'_{j})^{\hat{\alpha}'}:Z(\underline{\hat{M}}'(\tilde{t}))^{\Gamma_{{\mathbb R}}}]=
  \vert \check{\alpha}\vert \vert \check{\alpha}(\tilde{s})\vert ^{-1} [Z(\hat{M}'_{j})^{\hat{\alpha}(\tilde{s})}:Z(\underline{\hat{M}}'(\tilde{t}))^{\Gamma_{{\mathbb R}}}].$$
  Par d\'efinition, le dernier terme ci-dessus est \'egal \`a $i_{\tilde{M}'_{j}}^{\tilde{G}'(\tilde{s})}(\epsilon_{j})^{-1}$.
  L'\'egalit\'e (6) se transforme en 
  $$i_{\tilde{M}'_{j}}(\tilde{G},\tilde{G}'(\tilde{s})) i_{\tilde{M}'_{j}}^{\tilde{G}'(\tilde{s})}(\epsilon_{j})\vert \check{\alpha}(\tilde{s})\vert=m ^{-1} \vert \check{\alpha}\vert i_{\underline{\tilde{M}}'(\tilde{t})}(\tilde{G},\tilde{G}'(\tilde{s})).$$
   Le m\^eme calcul qu'en 4.2(4) montre que  l'on a l'\'egalit\'e 
   $$i_{\tilde{M}'_{1}(\tilde{s})}^{\tilde{G}'_{1}(\tilde{s})}(\epsilon_{1}(\tilde{s}))= i_{\tilde{M}'_{j}}^{\tilde{G}'(\tilde{s})}(\epsilon_{j}).$$
   L'\'egalit\'e ci-dessus devient (5).
 
A l'aide de cette \'egalit\'e (5), on r\'ecrit
$$D_{j}(X,\tilde{t})=m^{-1}\vert \check{\alpha}\vert \sum_{\tilde{s}\in \tilde{t}\in Z(\underline{\hat{M}})^{\Gamma_{{\mathbb R}},\hat{\theta}}/Z(\hat{G})^{\Gamma_{{\mathbb R}},\hat{\theta}}}i_{\underline{\tilde{M}}'(\tilde{t})}(\tilde{G},\tilde{G}'(\tilde{s})) \Delta_{1}(\tilde{s})(exp(\xi_{j}(X))\epsilon_{1}(\tilde{s}),exp(X)\eta)^{-1}$$
$$S_{\underline{\tilde{M}}'_{1}(\tilde{s}),\lambda_{1}(\tilde{s})}^{\tilde{G}'_{1}(\tilde{s})}(exp(\xi_{j}(X))\epsilon_{1}(\tilde{s}),f^{\tilde{G}'_{1}(\tilde{s})}).$$
Il est clair que, pour tout $\tilde{s}$, l'int\'egrale orbitale stable associ\'ee \`a $exp(\xi_{j}(X))\epsilon_{1}(\tilde{s})$, multipli\'ee par $ \Delta_{1}(\tilde{s})(exp(\xi_{j}(X))\epsilon_{1}(\tilde{s}),exp(X)\eta)^{-1}$, s'identifie \`a une unique distribution appartenant \`a $D_{g\acute{e}om}^{st}(\underline{{\bf M}}'(\tilde{t}))$. C'est la distribution $\boldsymbol{\delta}_{j}(X)^{\underline{{\bf M}}'(\tilde{t})}$ induite de la distribution $\boldsymbol{\delta}_{j}(X)$ introduite en 4.4. L'\'egalit\'e ci-dessus se r\'ecrit
$$D_{j}(X,\tilde{t})=m^{-1}\vert \check{\alpha}\vert I_{K\underline{\tilde{M}}}^{K\tilde{G},{\cal E}}(\underline{{\bf M}}'(\tilde{t}),\boldsymbol{\delta}_{j}(X)^{\underline{{\bf M}}'(\tilde{t})},f).$$
Nos distributions sont \`a support $\tilde{G}$-\'equisingulier. On peut donc r\'ecrire
 $$D_{j}(X,\tilde{t})=m^{-1}\vert \check{\alpha}\vert I_{K\underline{\tilde{M}}}^{K\tilde{G},{\cal E}}(transfert(\boldsymbol{\delta}_{j}(X)^{\underline{{\bf M}}'(\tilde{t})}),f).$$
Ou encore, en utilisant la commutation du  transfert \`a l'induction,
$$D_{j}(X,\tilde{t})=m^{-1}\vert \check{\alpha}\vert I_{K\underline{\tilde{M}}}^{K\tilde{G},{\cal E}}((transfert(\boldsymbol{\delta}_{j}(X)))^{\underline{\tilde{M}}},f).$$
Ceci est ind\'ependant de $\tilde{t}$. Alors $D_{j}(X)$ est la m\^eme expression, multipli\'ee par le nombre d'\'el\'ements de la sommation en $\tilde{t}$, lequel n'est autre que $m$. D'o\`u
$$(7) \qquad D_{j}(X)=\vert \check{\alpha}\vert I_{K\underline{\tilde{M}}}^{K\tilde{G},{\cal E}}(transfert(\boldsymbol{\delta}_{j}(X))^{\underline{\tilde{M}}},f).$$

Reprenons les formules (1) et (4). On obtient
$$(8) \qquad I_{K\tilde{M}}^{K\tilde{G},{\cal E}}(exp(X)\eta,\omega,f)=\phi(X)+B(X)-\vert \check{\alpha}\vert log(\vert \alpha(X)\vert )D(X),$$
o\`u
$$B(X)=[T^{\theta}({\mathbb R}):T^{\theta,0}({\mathbb R})]\vert J\vert ^{-1}d(\theta^*)^{-1/2}\sum_{j\in J }B_{j}(X),$$
$$D(X)=\vert \check{\alpha}\vert ^{-1}[T^{\theta}({\mathbb R}):T^{\theta,0}({\mathbb R})]\vert J\vert ^{-1}d(\theta^*)^{-1/2}\sum_{j\in J }D_{j}(X),$$
et $\phi$ est une fonction $C^{\infty}$ au voisinage de $0$ dans $\mathfrak{t}^{\theta}({\mathbb R})$. Gr\^ace \`a (7), on obtient
$$D(X)= I_{K\underline{\tilde{M}}}^{K\tilde{G},{\cal E}}(\boldsymbol{\gamma}'(X)^{\underline{\tilde{M}}},f),$$
o\`u
$$\boldsymbol{\gamma}'(X)=[T^{\theta}({\mathbb R}):T^{\theta,0}({\mathbb R})]\vert J\vert ^{-1}d(\theta^*)^{-1/2}\sum_{j\in J }transfert(\boldsymbol{\delta}_{j}(X)).$$
D'apr\`es 4.4(1), on a $\boldsymbol{\gamma}'(X)=\boldsymbol{\gamma}(X)$. Rappelons que cette distribution est $I^{K\tilde{M}}(exp(X)\eta,\omega,.)$. Le point $exp(X)\eta$ \'etant r\'egulier dans $\tilde{G}$, l'induite de cette distribution  est $I^{K\underline{\tilde{M}}}(exp(X)\eta,\omega,.)$. Alors $D(X)=I_{K\underline{\tilde{M}}}^{K\tilde{G},{\cal E}}(exp(X)\eta,\omega,f)$. La formule (8) se r\'ecrit
$$I_{K\tilde{M}}^{K\tilde{G},{\cal E},mod}(exp(X)\eta,\omega,f)=\phi(X)+B(X).$$
Soit $U\in Sym(\mathfrak{t}^{\theta})$. Appliquons l'op\'erateur $\partial_{U}$. Pour tout $j\in J$, posons
$$B_{j}(X,U)=\sum_{\tilde{s}\in\tilde{\zeta}_{j} Z(\hat{M})^{\hat{\alpha}^*}/Z(\hat{G})^{\Gamma_{{\mathbb R}},\hat{\theta}}}i_{\tilde{M}'_{j}}(\tilde{G},\tilde{G}'(\tilde{s}))\Delta_{1}(\tilde{s})(exp(\xi_{j}(X))\epsilon_{1}(\tilde{s}),exp(X)\eta)^{-1}$$
$$\partial_{\Xi(\tilde{s},U)}S_{\tilde{M}'_{1}(\tilde{s}),\lambda_{1}(\tilde{s})}^{\tilde{G}'_{1}(\tilde{s}),mod}(exp(\xi_{j}(X))\epsilon_{1}(\tilde{s}),f^{\tilde{G}'_{1}(\tilde{s})}),$$
o\`u $\Xi(\tilde{s})$ est l'homomorphisme $Sym(\mathfrak{t}^{\theta})\to Sym(\mathfrak{t}'_{j})$ de 3.1. Posons
$$B(X,U)=[T^{\theta}({\mathbb R}):T^{\theta,0}({\mathbb R})]\vert J\vert ^{-1}d(\theta^*)^{-1/2}\sum_{j\in J }B_{j}(X,U).$$
On obtient alors
$$  \partial_{U}I_{K\tilde{M}}^{K\tilde{G},{\cal E},mod}(exp(X)\eta,\omega,f)=\partial_{U}\phi(X)+B(X,U).$$
Appliquons cela \`a $X=rH_{d}$ pour un r\'eel $r$ non nul et proche de $0$. On obtient
$$(9) \qquad \partial_{U}I_{K\tilde{M}}^{K\tilde{G},{\cal E},mod}(exp(rH_{d})\eta,\omega,f)=\partial_{U}\phi(rH_{d})+B(rH_{d},U).$$
Remarquons que, pour tout $j\in J$,  la d\'efinition de $B_{j}(X,U)$ et  la formule 4.4(4) conduisent \`a l'\'egalit\'e
 $$(10)\qquad B_{j}(rH_{d},U)=\sum_{\tilde{s}\in\tilde{\zeta}_{j} Z(\hat{M})^{\hat{\alpha}^*}/Z(\hat{G})^{\Gamma_{{\mathbb R}},\hat{\theta}}}i_{\tilde{M}'_{j}}(\tilde{G},\tilde{G}'(\tilde{s})) d(\tilde{s})^{-1}$$
$$\partial_{U_{\tilde{s}} }S_{\tilde{M}'_{1}(\tilde{s}),\lambda_{1}(\tilde{s})}^{\tilde{G}'_{1}(\tilde{s}),mod}(exp(\xi_{j}(rH_{d}))\epsilon_{1}(\tilde{s}),f^{\tilde{G}'_{1}(\tilde{s})}),$$
o\`u on a pos\'e $U_{\tilde{s}}=\Xi(\tilde{s},U)$.

Les deux premi\`eres assertions de l'\'enonc\'e r\'esultent de ces formules. En effet, elles sont v\'erifi\'ees pour la fonction $\phi$ qui est $C^{\infty}$. Elles le sont d'apr\`es la proposition 4.2 pour les int\'egrales orbitales stables modifi\'ees qui figurent dans la d\'efinition des fonctions $B_{j}(rH_{d},U)$.  Les assertions sont donc v\'erifi\'ees pour la fonction $B(rH_{d},U)$ et la conclusion s'ensuit.

La troisi\`eme assertion de l'\'enonc\'e se prouve comme en 4.2. Puisque $\eta$, vu comme un \'el\'ement de $\underline{\tilde{M}}({\mathbb R})$, est $\tilde{G}$-\'equisingulier, il existe une fonction $\varphi\in C_{c}^{\infty}(\underline{\tilde{M}}({\mathbb R}))$ telle que, pour tout $X\in \underline{\mathfrak{m}}_{\eta}({\mathbb R})$ assez proche de $0$, on ait l'\'egalit\'e
$$I_{K\underline{\tilde{M}}}^{K\tilde{G},{\cal E}}(exp(X)\eta,\omega,f)=I^{\underline{\tilde{M}}}(exp(X)\eta,\omega,\varphi).$$
L'assertion (iii) r\'esulte des propri\'et\'es bien connues des int\'egrales orbitales. Le m\^eme argument d\'emontre l'\'egalit\'e des deux derni\`eres limites de l'assertion (iv). 

Il reste \`a prouver la premi\`ere \'egalit\'e de cette assertion.  Soient  $j\in J$ et $\tilde{s}\in \tilde{\zeta}_{j}Z(\hat{M})^{\hat{\alpha}^*}$.  Le sous-tore $T_{c}$ de $G_{\eta,SC}$ se transf\`ere en un sous-tore de $\bar{G}'(\bar{s})$ (cf. 4.4), qui se transf\`ere lui-m\^eme en un sous-tore $T_{c}(\tilde{s})$ de $G'(\tilde{s})_{\epsilon_{j},SC}$. Notons $\underline{T}'(\tilde{s})$ le commutant dans $G'(\tilde{s})_{\epsilon_{j}}$ de l'image de $T_{c}(\tilde{s})$ dans ce groupe et posons $\underline{\tilde{T}}'(\tilde{s})=\underline{T}'(\tilde{s})\epsilon_{j}$.   Les paires $(\underline{T},\underline{\tilde{T}})$ et $(\underline{T}'(\tilde{s}),\underline{\tilde{T}}'(\tilde{s}))$ sont dans la m\^eme situation que les paires $(T,\tilde{T})$ et $(T'_{j},\tilde{T}'_{j})$. C'est-\`a-dire qu'il y a un homomorphisme $\underline{\xi}(\tilde{s}):\underline{T}\to \underline{T}'(\tilde{s})$ et une application compatible $\underline{\tilde{\xi}}(\tilde{s}):\underline{\tilde{T}}\to \underline{\tilde{T}}'(\tilde{s})$ qui v\'erifient les conditions suivantes. Ils sont d\'efinis sur ${\mathbb R}$. On a $\underline{\tilde{\xi}}(\tilde{s})(\eta)=\epsilon_{j}$. Pour tout $\gamma\in \underline{\tilde{T}}({\mathbb R})$, les \'el\'ements $\gamma$ et $\underline{\tilde{\xi}}(\tilde{s})(\gamma)$ se correspondent. De plus, $\underline{\xi}(\tilde{s})$ co\"{\i}ncide avec  $\xi_{j}$ sur $Z(G_{\eta})^0$.  De m\^eme que l'on a introduit un isomorphisme $C:T\to \underline{T}$ d\'efini sur ${\mathbb C}$, on introduit un isomorphisme $C(\tilde{s}):T'_{j}\to \underline{T}'(\tilde{s})$. Il est plus ou moins clair que l'on peut le choisir tel que le diagramme suivant soit commutatif
$$\begin{array}{ccc}T&\stackrel{\xi_{j}}{\to} &T'_{j}\\ C\downarrow&&C(\tilde{s})\downarrow\\ \underline{T}&\stackrel{\underline{\xi}(\tilde{s})}{\to}&\underline{T}'(\tilde{s})\\ \end{array}$$
On a d\'ej\`a fix\'e une d\'ecomposition
$$\mathfrak{t}'_{1}(\tilde{s})=\mathfrak{c}_{1}(\tilde{s})\oplus \mathfrak{t}'_{j}$$
provenant comme en 3.1 d'une d\'ecomposition d'alg\`ebres de Lie
$$\mathfrak{g}'_{1}(\tilde{s})=\mathfrak{c}_{1}(\tilde{s})\oplus \mathfrak{g}'(\tilde{s}).$$
On d\'eduit de cette derni\`ere une d\'ecomposition
$$\underline{\mathfrak{t}}'_{1}(\tilde{s})=\mathfrak{c}_{1}(\tilde{s})\oplus \underline{\mathfrak{t}}'(\tilde{s}).$$
Avec des notations compr\'ehensibles au moins par l'auteur, on a alors le diagramme
$$\begin{array}{ccc}Sym(\mathfrak{t}^{\theta})&\stackrel{\Xi(\tilde{s})}{\to}&Sym(\mathfrak{t}'_{j})\\ C\downarrow&&C(\tilde{s})\downarrow\\ Sym(\underline{\mathfrak{t}}^{\theta})&\stackrel{\underline{\Xi}(\tilde{s})}{\to}&Sym(\underline{\mathfrak{t}}'(\tilde{s}))\\ \end{array}$$
qui est commutatif pour la m\^eme raison qu'en 4.2

 Soit $U\in Sym(\mathfrak{t}^{\theta})$ tel que $w_{d}(U)=-U$, posons $\underline{U}=C(U)$. Pour $X\in \underline{\mathfrak{t}}({\mathbb R})\cap \mathfrak{g}_{\eta}({\mathbb R})$, posons
$$(11) \qquad \underline{B}_{j}(X,\underline{U})=\sum_{\tilde{s}\in\tilde{\zeta}_{j} Z(\hat{M})^{\hat{\alpha}^*}/Z(\hat{G})^{\Gamma_{{\mathbb R}},\hat{\theta}}}i_{\tilde{M}'_{j}}(\tilde{G},\tilde{G}'(\tilde{s}))i_{\tilde{M}'_{1}(\tilde{s})}^{\tilde{G}'_{1}(\tilde{s})}(\epsilon_{1}(\tilde{s}))\vert \check{\alpha}(\tilde{s})\vert \vert \check{\alpha}\vert ^{-1}$$
$$\Delta_{1}(\tilde{s})(exp(\underline{\xi}(\tilde{s},X))\epsilon_{1}(\tilde{s}),exp(X)\eta)^{-1}\partial_{\underline{U}_{\tilde{s}}}S_{\underline{\tilde{M}}'_{1}(\tilde{s}),\lambda_{1}(\tilde{s})}^{\tilde{G}'_{1}(\tilde{s}),mod}(exp(\underline{\xi}(\tilde{s},X))\epsilon_{1}(\tilde{s}),f^{\tilde{G}'_{1}(\tilde{s})}), $$ 
o\`u $\underline{U}_{\tilde{s}}=\underline{\Xi}(\tilde{s},\underline{U})=C(\tilde{s})\circ \Xi(\tilde{s})(U)$. Comme plus haut, on d\'ecompose
$$\underline{B}_{j}(X,\underline{U})=\sum_{\tilde{t}\in \tilde{\zeta}_{j}Z(\hat{M})^{\hat{\alpha}^*}/Z(\underline{\hat{M}})^{\Gamma_{{\mathbb R}},\hat{\theta}}}\underline{B}_{j}(X,\underline{U},\tilde{t}),$$
o\`u
$$\underline{B}_{j}(X,\underline{U},\tilde{t})=\sum_{\tilde{s}\in\tilde{t} Z(\underline{\hat{M}})^{\Gamma_{{\mathbb R}},\hat{\theta}}/Z(\hat{G})^{\Gamma_{{\mathbb R}},\hat{\theta}}}i_{\tilde{M}'_{j}}(\tilde{G},\tilde{G}'(\tilde{s}))i_{\tilde{M}'_{1}(\tilde{s})}^{\tilde{G}'_{1}(\tilde{s})}(\epsilon_{1}(\tilde{s}))\vert \check{\alpha}(\tilde{s})\vert \vert \check{\alpha}\vert ^{-1}$$
$$\Delta_{1}(\tilde{s})(exp(\underline{\xi}(\tilde{s},X))\epsilon_{1}(\tilde{s}),exp(X)\eta)^{-1}\partial_{\underline{U}_{\tilde{s}}}S_{\underline{\tilde{M}}'_{1}(\tilde{s}),\lambda_{1}(\tilde{s})}^{\tilde{G}'_{1}(\tilde{s}),mod}(exp(\underline{\xi}(\tilde{s},X))\epsilon_{1}(\tilde{s}),f^{\tilde{G}'_{1}(\tilde{s})}).$$
Fixons $\tilde{t}\in  \tilde{\zeta}_{j}Z(\hat{M})^{\hat{\alpha}^*}/Z(\underline{\hat{M}})^{\Gamma_{{\mathbb R}},\hat{\theta}}$. Pour simplifier la notation, on consid\`ere que $\tilde{t}$ est l'un des \'el\'ements de la sommation en $\tilde{s}$. Pour $\tilde{s}\in \tilde{t} Z(\underline{\hat{M}})^{\Gamma_{{\mathbb R}},\hat{\theta}}/Z(\hat{G})^{\Gamma_{{\mathbb R}},\hat{\theta}}$, les donn\'ees $\underline{M}'_{1}(\tilde{s})$,...,$\Delta_{1}(\tilde{s})$ sont diff\'erentes donn\'ees auxiliaires pour une m\^eme donn\'ee endospcopique elliptique $\underline{{\bf M}}'(\tilde{t})$ de $(\underline{M},\underline{\tilde{M}},{\bf a})$. Les distributions
$$\Delta_{1}(\tilde{s})(exp(\underline{\xi}(\tilde{s},X))\epsilon_{1}(\tilde{s}),exp(X)\eta)^{-1}\partial_{\underline{U}_{\tilde{s}}}S^{\underline{\tilde{M}}'_{1}(\tilde{s})}_{\lambda_{1}(\tilde{s})}(exp(\underline{\xi}(\tilde{s},X))\epsilon_{1}(\tilde{s}),.)$$
se recollent en une m\^eme distribution $\boldsymbol{\delta}_{X,\underline{U},\tilde{t}}\in D_{g\acute{e}om,\tilde{G}-\acute{e}qui}^{st}(\underline{{\bf M}}'(\tilde{t}))$. En utilisant (5), on obtient
$$(12)\qquad \underline{B}_{j}(X,U,\tilde{t})=m^{-1}I_{K\underline{\tilde{M}}}^{K\tilde{G},{\cal E}}(\underline{{\bf M}}'(\tilde{t}),\boldsymbol{\delta}_{X,\underline{U},\tilde{t}},f)=m^{-1}I_{K\underline{\tilde{M}}}^{K\tilde{G},{\cal E}}(transfert(\boldsymbol{\delta}_{X,\underline{U},\tilde{t}}),f).$$
On pose
$$\boldsymbol{\gamma}_{X,\underline{U},j}=m^{-1}\sum_{\tilde{t}\in  \tilde{\zeta}_{j}Z(\hat{M})^{\hat{\alpha}^*}/Z(\underline{\hat{M}})^{\Gamma_{{\mathbb R}},\hat{\theta}}}transfert(\boldsymbol{\delta}_{X,\underline{U},\tilde{t}}),$$
$$\boldsymbol{\gamma}_{X,\underline{U}}=[T^{\theta}({\mathbb R}):T^{\theta,0}({\mathbb R})]\vert J\vert ^{-1}d(\theta^*)^{-1/2}\sum_{j\in J}\boldsymbol{\gamma}_{X,\underline{U},j},$$
$$\underline{B}(X,\underline{U})=[T^{\theta}({\mathbb R}):T^{\theta,0}({\mathbb R})]\vert J\vert ^{-1}d(\theta^*)^{-1/2}\sum_{j\in J}\underline{B}_{j}(X,\underline{U}).$$
La formule (12) entra\^{\i}ne
$$(13)\qquad \underline{B}(X,\underline{U})=I_{K\underline{\tilde{M}}}^{K\tilde{G},{\cal E}}(\boldsymbol{\gamma}_{X,\underline{U}},f).$$

On va calculer la distribution $\boldsymbol{\gamma}_{X,U}$.   Supposons l'\'el\'ement $exp(X)\eta$ fortement r\'egulier. L'ensemble ${\cal Y}^{\underline{M}}(exp(X)\eta)$ introduit en 4.4  ne d\'epend par d\'efinition que de $\underline{\tilde{T}}$. D'autre part, l'inclusion $I^{\underline{M}}_{exp(X)\eta}\subset I^{\underline{M}}_{\eta}$ entra\^{\i}ne l'inclusion ${\cal Y}^{\underline{M}}(exp(X)\eta)\subset {\cal Y}^{\underline{M}}(\eta)$. On peut donc d'une part
  choisir l'ensemble $\dot{{\cal Y}}^{\underline{M}}(exp(X)\eta)$ ind\'ependant de $X$ et on le note plut\^ot $\dot{{\cal Y}}^{\underline{M}}(\underline{\tilde{T}})$; d'autre part supposer qu'il existe une application
  $$q:\dot{{\cal Y}}^{\underline{M}}(\underline{\tilde{T}})\to \dot{{\cal Y}}^{\underline{M}}(\eta)$$
  de sorte que, pour $y\in \dot{{\cal Y}}^{\underline{M}}(\underline{\tilde{T}})$, $q(y)$ soit l'unique \'el\'ement de $I^{\underline{M}}_{\eta}y\cap \dot{{\cal Y}}^{\underline{M}}(\eta)$. On utilise des notations similaires \`a celles de 4.1: pour $p\in \Pi^{\underline{M}}$ et $y\in  \dot{{\cal Y}}^{\underline{M}}(\underline{\tilde{T}})$, on pose $\underline{\tilde{T}}[y]=\tilde{\phi}_{p}^{-1}\circ ad_{y^{-1}}(\underline{\tilde{T}})$ et on note simplement $x\mapsto x[y]$ l'isomorphisme $\phi_{p}^{-1}\circ ad_{y^{-1}}$ de $\underline{T}$ sur $\underline{T}[y]$.

Soient $j\in J$ et $\tilde{t}\in \tilde{\zeta}_{j}Z(\hat{M})^{\hat{\alpha}^*}/Z(\underline{\hat{M}})^{\Gamma_{{\mathbb R}},\hat{\theta}}$.  Supposons comme ci-dessus l'\'el\'ement $exp(X)\eta$ fortement r\'egulier. Puisque $\underline{\tilde{T}}$ est un sous-tore tordu maximal et elliptique dans $\underline{\tilde{M}}$, l'\'el\'ement $exp(X)\eta$ est elliptique dans $\underline{\tilde{M}}({\mathbb R})$. Notons $\boldsymbol{\delta}'_{X,\tilde{t}}$ l'\'el\'ement de $ D_{g\acute{e}om,\tilde{G}-reg}^{st}(\underline{{\bf M}}'(\tilde{t}))$ auquel s'identifie la distribution $S^{\underline{\tilde{M}}'_{1}(\tilde{t})}_{\lambda_{1}(\tilde{t})}(exp(\underline{\xi}(\tilde{t},X))\epsilon_{1}(\tilde{t}),.)$. D'apr\`es [I] 4.9(4), on a l'\'egalit\'e
$$transfert(\boldsymbol{\delta}'_{X,\tilde{t}})=d(\theta^*)^{1/2}\sum_{y\in \dot{{\cal Y}}^{\underline{M}}(\underline{\tilde{T}})}\Delta_{1}(\tilde{t})(exp(\underline{\xi}(\tilde{t},X))\epsilon_{1}(\tilde{t}),exp(X[y])\eta[q(y)])$$
$$[\underline{T}[y]^{\theta}({\mathbb R}):\underline{T}[y]^{\theta,0}({\mathbb R})]^{-1} I^{K\underline{\tilde{M}}}(exp(X[y])\eta[q(y)]),\omega,.).$$
Notons  $\boldsymbol{\delta}'_{X,\underline{U},\tilde{t}}$  l'\'el\'ement de $ D_{g\acute{e}om,\tilde{G}-reg}^{st}(\underline{{\bf M}}'(\tilde{t}))$ auquel s'identifie la distribution 
$$\partial_{\underline{U}_{\tilde{t}}}S^{\underline{\tilde{M}}'_{1}(\tilde{t})}_{\lambda_{1}(\tilde{t})}(exp(\underline{\xi}(\tilde{t},X))\epsilon_{1}(\tilde{t}),.).$$
Alors
$$transfert(\boldsymbol{\delta}'_{X,\underline{U},\tilde{t}})=d(\theta^*)^{1/2}\sum_{y\in \dot{{\cal Y}}^{\underline{M}}(\underline{\tilde{T}})}\Delta_{1}(\tilde{t})(exp(\underline{\xi}(\tilde{t},X))\epsilon_{1}(\tilde{t}),exp(X[y])\eta[q(y)])$$
$$[\underline{T}[y]^{\theta}({\mathbb R}):\underline{T}[y]^{\theta,0}({\mathbb R})]^{-1} \partial_{\underline{U}[y]}I^{K\underline{\tilde{M}}}(exp(X[y])\eta[q(y)]),\omega,.).$$
Comme en 4.2, l'op\'erateur $\underline{U}[y]$ apparaissant ici n'est autre que le transport\'e de $\underline{U}$ par  l'isomorphisme $x\mapsto x[y]$ de $\underline{T}$ sur $\underline{T}[y]$. Pour obtenir le transfert de $\boldsymbol{\delta}_{X,\underline{U},\tilde{t}}$, il reste \`a diviser par $\Delta_{1}(\tilde{t})(exp(\underline{\xi}(\tilde{t},X))\epsilon_{1}(\tilde{t}),exp(X)\eta)$. D'o\`u
$$ transfert(\boldsymbol{\delta}_{X,\underline{U},\tilde{t}})=d(\theta^*)^{1/2}\sum_{y\in \dot{{\cal Y}}^{\underline{M}}(\underline{\tilde{T}})}\frac{\Delta_{1}(\tilde{t})(exp(\underline{\xi}(\tilde{t},X))\epsilon_{1}(\tilde{t}),exp(X[y])\eta[q(y)])}{\Delta_{1}(\tilde{t})(exp(\underline{\xi}(\tilde{t},X))\epsilon_{1}(\tilde{t}),exp(X)\eta)}$$
$$[\underline{T}[y]^{\theta}({\mathbb R}):\underline{T}[y]^{\theta,0}({\mathbb R})]^{-1} \partial_{\underline{U}[y]}I^{K\underline{\tilde{M}}}(exp(X[y])\eta[q(y)]),\omega,.).$$
 On peut maintenant  remplacer $X$ par $exp(rH_{c})$ pour $r$ r\'eel non nul et proche de $0$, tous les termes \'etant continus en un tel point. En utilisant les relations (4) et (5) de 4.4,  on obtient
 $$   transfert(\boldsymbol{\delta}_{rH_{c},\underline{U},\tilde{t}})=d(\theta^*)^{1/2}\sum_{y\in \dot{{\cal Y}}^{\underline{M}}(\underline{\tilde{T}})} d(\tilde{t},q(y))d(\tilde{t})^{-1}[\underline{T}[y]^{\theta}({\mathbb R}):\underline{T}[y]^{\theta,0}({\mathbb R})]^{-1} $$
 $$\partial_{\underline{U}[y]}I^{K\underline{\tilde{M}}}(exp(rH_{c}[y])\eta[q(y)]),\omega,.).$$
 Pour tout $y\in  \dot{{\cal Y}}^{\underline{M}}(\eta)$, posons
 $$d(y)=m^{-1}\vert J\vert ^{-1}\sum_{j\in J}\sum_{\tilde{t}\in \tilde{\zeta}_{j}Z(\hat{M})^{\hat{\alpha}^*}/Z(\underline{\hat{M}})^{\Gamma_{{\mathbb R}},\hat{\theta}}} d(\tilde{t},y)d(\tilde{t})^{-1}.$$
  Les d\'efinitions entra\^{\i}nent alors
 $$(14) \qquad \boldsymbol{\gamma}_{rH_{c},\underline{U}}=[T^{\theta}({\mathbb R}):T^{\theta,0}({\mathbb R})]\sum_{y\in \dot{{\cal Y}}^{\underline{M}}(\underline{\tilde{T}})}d(q(y))[\underline{T}[y]^{\theta}({\mathbb R}):\underline{T}[y]^{\theta,0}({\mathbb R})]^{-1} $$
 $$\partial_{\underline{U}[y]}I^{K\underline{\tilde{M}}}(exp(rH_{c}[y])\eta[q(y)]),\omega,.).$$

 Consid\'erons un \'el\'ement $X\in \mathfrak{t}({\mathbb R})\cap \mathfrak{g}_{\eta,reg}({\mathbb R})$ proche de $0$. On a $I^{M}_{\eta}=T^{\theta}=I^{M}_{exp(X)\eta}$. On peut donc supposer $\dot{{\cal Y}}^{M}(exp(X)\eta)=\dot{{\cal Y}}^{M}(\eta)$. On a l'\'egalit\'e $I^{M}_{\eta}=I^{\underline{M}}_{\eta}\cap M$. Il en r\'esulte que, pour tout $p\in \Pi^{M}$, on a une application naturelle
 $$(15)\qquad I^M_{\eta}\backslash {\cal Y}_{p}^M(\eta)/\phi_{p}(M_{p}({\mathbb R}))\to I^{\underline{M}}_{\eta}\backslash {\cal Y}_{p}^{\underline{M}}(\eta)/\phi_{p}(\underline{M}_{p}({\mathbb R})).$$
 Montrons que
 
 (16) cette application est injective; son image est l'image dans l'espace d'arriv\'ee de l'ensemble des $y\in {\cal Y}_{p}^{\underline{M}}(\eta)$ tels que $\underline{M}_{\eta[y],SC}$ soit isomorphe \`a $SL(2)$.
 
  Consid\'erons l'application qui, \`a $y\in {\cal Y}_{p}^M(\eta)$, associe la classe du cocycle $\sigma\mapsto y\nabla_{p}(\sigma)\sigma(y)^{-1}$. Elle se quotiente en une bijection $I^M_{\eta}\backslash {\cal Y}_{p}^M(\eta)/\phi_{p}(M_{p}({\mathbb R}))$ sur l'image r\'eciproque de la classe du cocycle $\nabla_{p}$ par l'application
  $$H^1(\Gamma_{{\mathbb R}};T^{\theta})\to H^1(\Gamma_{{\mathbb R}};M).$$
  De m\^eme, on a une bijection de $ I^{\underline{M}}_{\eta}\backslash {\cal Y}_{p}^{\underline{M}}(\eta)/\phi_{p}(\underline{M}_{p}({\mathbb R}))$ sur l'image r\'eciproque de la classe du cocycle $\nabla_{p}$ par l'application
  $$H^1(\Gamma_{{\mathbb R}};I^{\underline{M}}_{\eta})\to H^1(\Gamma_{{\mathbb R}};\underline{M}).$$
  Pour prouver l'injectivit\'e de (15), il suffit de prouver que l'application
  $$H^1(\Gamma_{{\mathbb R}};T^{\theta})\to H^1(\Gamma_{{\mathbb R}};I^{\underline{M}}_{\eta})$$
  est injective.  Notons $Z$ le centre de $I^{\underline{M}}_{\eta}$. Alors $I^{\underline{M}}_{\eta}/Z=\underline{M}_{\eta,AD}$ et on sait que ce groupe est $PGL(2)$. L'image $T_{ad}$ de $T^{\theta}$ dans ce groupe est un sous-tore maximal d\'eploy\'e. On a donc $H^1(\Gamma_{{\mathbb R}};T_{ad})=\{0\}$. Donc l'application naturelle $H^1(\Gamma_{{\mathbb R}};Z)\to H^1(\Gamma_{{\mathbb R}};T^{\theta})$ est surjective. Soient $u$ et $v$ deux cocycles \`a valeurs dans $Z$. Supposons que leurs images dans $H^1(\Gamma_{{\mathbb R}};I^{\underline{M}}_{\eta})$ soient \'egales. On doit prouver que leurs images dans $H^1(\Gamma_{{\mathbb R}};T^{\theta})$ le sont. L'hypoth\`ese signifie qu'il existe $x\in I^{\underline{M}}_{\eta}$ tel que $u(\sigma)=xv(\sigma)\sigma(x)^{-1}$ pour tout $\sigma\in \Gamma_{{\mathbb R}}$.  Puisque $u$ est \`a valeurs centrales, cela \'equivaut \`a $\sigma(x)u(\sigma)v(\sigma)^{-1}=x$. Introduisons un sous-groupe de Borel $B_{ad}=T_{ad}U$ de $\underline{M}_{\eta,AD}$ de tore $T_{ad}$ et relevons-le en un sous-groupe $B=T^{\theta}U$ de $I^{\underline{M}}_{\eta}$. Ces groupes sont d\'efinis sur ${\mathbb R}$.  Puisque $\underline{M}_{\eta,SC}\simeq SL(2)$, il existe un \'el\'ement $w\in \underline{M}_{\eta,SC}({\mathbb R})$ tel que l'\'el\'ement $x$ s'\'ecrive de fa\c{c}on unique $tn_{1}$ ou $n_{2}wtn_{1}$, avec $t\in T^{\theta}$ et $n_{1},n_{2}\in U$. Supposons par exemple $x=n_{2}wtn_{1}$. On a alors $\sigma(n_{2})w\sigma(t) u(\sigma)v(\sigma)^{-1}\sigma(n_{1})=x$. Cela fournit une autre d\'ecomposition de $x$. Par unicit\'e, on en d\'eduit $\sigma(t)u(\sigma)v(\sigma)^{-1}=t$, ou encore $u(\sigma)=tv(\sigma)\sigma(t)^{-1}$. Donc les images de $u$ et $v$ dans $H^1(\Gamma_{{\mathbb R}};T^{\theta})$ sont \'egales comme on le voulait. Cela d\'emontre l'injectivit\'e. La deuxi\`eme assertion de (16) se d\'emontre comme en 4.2(17). 
  
  Il r\'esulte de (16) que l'on peut supposer $\dot{{\cal Y}}^{M}(\eta)\subset \dot{{\cal Y}}^{\underline{M}}(\eta)$.   Soit $y\in \dot{{\cal Y}}^M(\eta)$. Pour $j\in J$ et $\tilde{t}\in \tilde{\zeta}_{j}Z(\hat{M})^{\hat{\alpha}^*}$, on a les \'egalit\'es
  $$d(\tilde{t},y)d(\tilde{t})^{-1}=d(\tilde{t},y)d(\tilde{t},1)^{-1}=d_{j}(y),$$
  avec les notations de 4.4. On en d\'eduit 
  $$d(y)=\vert J\vert ^{-1}\sum_{j\in J}d_{j}(y).$$
  L'assertion 4.4(3) nous dit que cette somme est nulle si $y\not=1$ et qu'elle vaut $1$ si $y=1$. 
     Montrons que
 
 (17) pour $y\in \dot{{\cal Y}}^{\underline{M}}(\eta)$, $y\not\in \dot{{\cal Y}}^M(\eta)$, et pour tout $j\in J$, on a l'\'egalit\'e 
 $$\sum_{\tilde{t}\in  \tilde{\zeta}_{j}Z(\hat{M})^{\hat{\alpha}^*}/Z(\underline{\hat{M}})^{\Gamma_{{\mathbb R}},\hat{\theta}}} d(\tilde{t},y)d(\tilde{t})^{-1}=0.$$
 
 Comme en 4.2(18), c'est un cas particulier de la proposition [III] 8.4. Soit $p\in \Pi^{\underline{M}}$ tel que  $y\in {\cal Y}_{p}^{\underline{M}}(\eta)$.  Le cocycle $\sigma\mapsto y\nabla_{p}(\sigma)\sigma(y)^{-1}$, \`a valeurs dans $I^{\underline{M}}_{\eta}$, se pousse en un cocycle \`a valeurs dans $\underline{M}_{\eta,ad}$, ce groupe \'etant l'image de $\underline{M}_{\eta}$ dans $G_{AD}$. On note $\tau[y]$ sa classe. Notons ici $\hat{T}$ le tore dual de $T$. On peut le consid\'erer comme celui de la paire de Borel \'epingl\'ee $\hat{{\cal E}}$, muni d'une action galoisienne tordue. Posons $\hat{\bar{T}}=\hat{T}/(1-\hat{\theta})(\hat{T})$. Puisque $T$ est un sous-tore maximal de $M$ et que $T^{\theta,0}$ est un sous-tore maximal de $\underline{M}_{\eta}$, on a des plongements $\hat{T}\to \hat{M}$ et $\hat{\bar{T}}\subset \underline{\hat{M}}_{\eta}$. En prenant les images r\'eciproques dans $\hat{G}_{SC}$, on a encore des plongements $\hat{T}_{sc}\to \hat{M}_{sc}$ et $\hat{\bar{T}}_{sc}\subset \underline{\hat{M}}_{\eta,sc}$. 
 Soit $t\in Z(\hat{M})^{\hat{\alpha}^*}$. On choisit $t_{sc}\in Z(\hat{M})_{sc}^{\Gamma_{{\mathbb R}},\hat{\theta},0}$ dont l'image dans $\hat{G}_{AD}$ soit la m\^eme que celle de $t$.   Parce que $\tilde{T}$ est elliptique dans $\tilde{M}$, on a l'\'egalit\'e $ Z(\hat{M})_{sc}^{\Gamma_{{\mathbb R}},\hat{\theta},0}=\hat{T}_{sc}^{\Gamma_{{\mathbb R}},\hat{\theta},0}$. L'\'el\'ement $t_{sc}$ s'envoie donc sur un \'el\'ement $\bar{t}_{sc}\in \hat{\bar{T}}_{sc}\subset \hat{\underline{M}}_{\eta,sc}$. On voit comme en 4.2 que l'hypoth\`ese $t\in Z(\hat{M})^{\hat{\alpha}^*}$ entra\^{\i}ne que $\bar{t}_{sc}\in Z(\hat{\underline{M}}_{sc})^{\Gamma_{{\mathbb R}}}$. On note $u(t)$ son image dans le groupe des composantes connexes $\pi_{0}(Z(\hat{\underline{M}}_{sc})^{\Gamma_{{\mathbb R}}})$. Le calcul de [III] 8.4 montre que
 $$d(\tilde{\zeta}_{j}t,y)d(\tilde{\zeta}_{j}t)^{-1}=d(\tilde{\zeta}_{j},y)d(\tilde{\zeta}_{j})^{-1}<\tau[y],u(t)>,$$
 le produit \'etant celui sur
 $$H^1(\Gamma_{{\mathbb R}};\underline{M}_{\eta,ad})\times \pi_{0}(Z(\hat{\underline{M}}_{sc})^{\Gamma_{{\mathbb R}}}).$$
 Il reste \`a prouver que l'application $t\mapsto <\tau[y],u(t)>$ est un caract\`ere non trivial de $Z(\hat{M})^{\hat{\alpha}^*}/Z(\underline{\hat{M}})^{\Gamma_{{\mathbb R}},\hat{\theta}}$. La preuve est similaire \`a celle de 4.2(18) et on la laisse au lecteur. 
 
 Il r\'esulte \'evidemment de (17) que $d(y)=0$ pour tout $y\not\in \dot{{\cal Y}}^M(\eta)$. Finalement $d(1)=1$ et $d(y)=0$ pour tout $y\not=1$. La formule (14) devient simplement
 $$  \boldsymbol{\gamma}_{rH_{c},\underline{U}}=[T^{\theta}({\mathbb R}):T^{\theta,0}({\mathbb R})]\sum_{y\in q^{-1}(1)}[\underline{T}[y]^{\theta}({\mathbb R}):\underline{T}[y]^{\theta,0}({\mathbb R})]^{-1}$$
 $$\partial_{\underline{U}[y]}I^{K\underline{\tilde{M}}}(exp(rH_{c}[y])\eta),\omega,.).$$
 La formule (13), \'evalu\'ee pour $X=exp(rH_{c})$, devient
 $$(18) \qquad \underline{B}(rH_{c},U)=[T^{\theta}({\mathbb R}):T^{\theta,0}({\mathbb R})]\sum_{y\in q^{-1}(1)}[\underline{T}[y]^{\theta}({\mathbb R}):\underline{T}[y]^{\theta,0}({\mathbb R})]^{-1}$$
 $$\partial_{\underline{U}[y]}I_{K\underline{\tilde{M}}}^{K\tilde{G},{\cal E}}(exp(rH_{c}[y])\eta,\omega,f).$$
 
 Montrons que
 
 (19)  on peut supposer que $\underline{T}[y]=\underline{T}$ pour tout $y\in q^{-1}(1)$;  le nombre d'\'el\'ements de $q^{-1}(1) $ est \'egal \`a
 $$2^{c(\eta)}[T^{\theta}({\mathbb R}):T^{\theta,0}({\mathbb R})]^{-1}[\underline{T}^{\theta}({\mathbb R}):\underline{T}^{\theta,0}({\mathbb R})];$$
 la limite
 $$lim_{r\to 0+}\partial_{\underline{U}[y]}I_{K\underline{\tilde{M}}}^{K\tilde{G},{\cal E}}(exp(rH_{c}[y])\eta,\omega,f)$$
 ne d\'epend pas de $y\in q^{-1}(1)$.

Notons ${\cal Y}_{0}$ l'ensemble des $y\in I^{\underline{M}}_{\eta}$ tels que $y\sigma(y)^{-1}\in \underline{T}^{\theta}$. D'apr\`es les d\'efinitions, $q^{-1}(y)$ est un ensemble de repr\'esentants de l'ensemble de doubles classes $\underline{T}^{\theta}\backslash{\cal Y}_{0}/I^{\underline{M}}_{\eta}({\mathbb R})$. 
 Pour $y\in q^{-1}(1)$, l'image r\'eciproque de $\underline{T}[y]$ dans $\underline{M}_{\eta,SC}$ est un sous-tore maximal de ce groupe, qui est d\'efini sur ${\mathbb R}$ et non d\'eploy\'e. Il est donc conjugu\'e \`a $T_{c}$ par un \'el\'ement de $\underline{M}_{\eta,SC}({\mathbb R})$. Quitte \`a multiplier $y$ \`a droite par un tel \'el\'ement, on peut donc supposer que cette image r\'eciproque est $T_{c}$. Mais cela entra\^{\i}ne que $\underline{T}[y]=\underline{T}$. D'o\`u la premi\`ere assertion. Notons ${\cal Y}_{1}={\cal Y}_{0}\cap \underline{M}_{\eta}$. Parce que $I^{\underline{M}}_{\eta}=\underline{T}^{\theta}\underline{M}_{\eta}$, on voit que l'injection ${\cal Y}_{1}\to {\cal Y}_{0}$ induit une surjection
 $${\cal Y}_{1}\to \underline{T}^{\theta}\backslash{\cal Y}_{0}/I^{\underline{M}}_{\eta}({\mathbb R}).$$
 Il s'en d\'eduit une surjection
 $$\underline{T}^{\theta,0}\backslash {\cal Y}_{1}/\underline{M}_{\eta}({\mathbb R})\to \underline{T}^{\theta}\backslash{\cal Y}_{0}/I^{\underline{M}}_{\eta}({\mathbb R}).$$
 Le premier ensemble classifie les classes de conjugaison par $\underline{M}_{\eta}({\mathbb R})$ dans la classe de conjugaison stable de $H_{c}$ dans $\underline{\mathfrak{m}}_{\eta}({\mathbb R})$. Son nombre d'\'el\'ements est $2^{c(\eta)}$. Ce nombre est au plus $2$ donc les fibres de la surjection ci-dessus ont forc\'ement le m\^eme nombre d'\'el\'ements. Le nombre d'\'el\'ements de $q^{-1}(1)$ est alors $2^{c(\eta)}$ divis\'e par le nombre d'\'el\'ements d'une fibre.
 La fibre au-dessus de la classe $\underline{T}^{\theta}I^{\underline{M}}_{\eta}({\mathbb R})$ est
 $${\mathbb  F}=\underline{T}^{\theta,0}\backslash(\underline{T}^{\theta}I^{\underline{M}}_{\eta}({\mathbb R})\cap \underline{M}_{\eta})/\underline{M}_{\eta}({\mathbb R}).$$
 On d\'efinit une application de $I^{\underline{M}}_{\eta}({\mathbb R})$ dans ce quotient de la fa\c{c}on suivante. Pour $x\in I^{\underline{M}}_{\eta}({\mathbb R})$, on choisit $t\in \underline{T}^{\theta}$ tel que $tx\in \underline{M}_{\eta}$ et on envoie $x$ sur l'image de $tx$ dans ${\mathbb F}$. Cette image ne d\'epend pas du choix de $t$. On v\'erifie que l'application ainsi d\'efinie se quotiente en une bijection
 $$\underline{T}^{\theta}({\mathbb R})\backslash I^{\underline{M}}_{\eta}({\mathbb R})/\underline{M}_{\eta}({\mathbb R})\to {\mathbb F}.$$
 Le nombre d'\'el\'ements du premier ensemble est
 $$[ I^{\underline{M}}_{\eta}({\mathbb R}):\underline{M}_{\eta}({\mathbb R})][\underline{T}^{\theta}({\mathbb R}):\underline{T}^{\theta,0}({\mathbb R})]^{-1}.$$
 On a une application naturelle
 $$T^{\theta}({\mathbb R})/T^{\theta,0}({\mathbb R})\to  I^{\underline{M}}_{\eta}({\mathbb R})/\underline{M}_{\eta}({\mathbb R}).$$
 Elle est clairement injective. Reprenons la preuve de (16). Tout \'el\'ement $x\in I^{\underline{M}}_{\eta}$ s'\'ecrit de fa\c{c}on unique $x=tn_{1}$ ou $x=n_{2}wtn_{1}$ avec les notations de cette preuve. Si $x\in I^{\underline{M}}_{\eta}({\mathbb R})$, l'unicit\'e entra\^{\i}ne que ces \'el\'ements $t$, $n_{1}$, $n_{2}$ sont tous d\'efinis sur ${\mathbb R}$. Les \'el\'ements $n_{1}$, $n_{2}$ et $w$ appartenant \`a $\underline{M}_{\eta}({\mathbb R})$ et ce groupe \'etant distingu\'e dans $I^{\underline{M}}_{\eta}({\mathbb R})$, on voit que $x$ et $t$ ont m\^eme image modulo $\underline{M}_{\eta}({\mathbb R})$. Donc l'application ci-dessus est surjective. D'o\`u l'\'egalit\'e
 $$[ I^{\underline{M}}_{\eta}({\mathbb R}):\underline{M}_{\eta}({\mathbb R})]=[T^{\theta}({\mathbb R}):T^{\theta,0}({\mathbb R})].$$
 En mettant ces calculs bout \`a bout, on obtient la deuxi\`eme assertion de (19).  Il r\'esulte de la description ci-dessus que, pour $y\in q^{-1}(1)$, on a soit $H_{c}[y]=H_{c}$ et alors $\underline{U}[y]=\underline{U}$, soit $H_{c}[y]=-H_{c}$ et alors $\underline{U}[y]$ est l'image de $\underline{U}$ par la sym\'etrie $w_{c}$ de $Sym(\underline{\mathfrak{t}})$.  On a suppos\'e $w_{d}(U)=-U$, donc $w_{c}(\underline{U})=-\underline{U}$. Dans ce cas, on a 
 $$lim_{r\to 0+}\partial_{\underline{U}[y]}I_{K\underline{\tilde{M}}}^{K\tilde{G},{\cal E}}(exp(rH_{c}[y])\eta,\omega,f)=-lim_{r\to 0-}\partial_{\underline{U}}I_{K\underline{\tilde{M}}}^{K\tilde{G},{\cal E}}(exp(rH_{c})\eta,\omega,f).$$
 Mais on a d\'ej\`a remarqu\'e que les deux derni\`eres limites de l'assertion (iv) de l'\'enonc\'e \'etaient \'egales. C'est-\`a-dire que la limite ci-dessus n'est autre que
 $$lim_{r\to 0+}\partial_{\underline{U}}I_{K\underline{\tilde{M}}}^{K\tilde{G},{\cal E}}(exp(rH_{c})\eta,\omega,f).$$
 Cela d\'emontre la troisi\`eme assertion de (19).

  Il r\'esulte de (18) et (19)  que
  $$(20) \qquad lim_{r\to 0+}\underline{B}(rH_{c},U)=2^{c(\eta)}lim_{r\to 0+}\partial_{\underline{U}}I_{K\underline{\tilde{M}}}^{K\tilde{G},{\cal E}}(exp(rH_{c})\eta,\omega,f).$$

 Soit $j\in J$. En rempla\c{c}ant $X$ par $rH_{c}$ dans la formule (11) et en utilisant 4.4(4), on obtient
 $$\underline{B}_{j}(rH_{c},\underline{U})=\sum_{\tilde{s}\in\tilde{\zeta}_{j} Z(\hat{M})^{\hat{\alpha}^*}/Z(\hat{G})^{\Gamma_{{\mathbb R}},\hat{\theta}}}i_{\tilde{M}'_{j}}(\tilde{G},\tilde{G}'(\tilde{s}))i_{\tilde{M}'_{1}(\tilde{s})}^{\tilde{G}'_{1}(\tilde{s})}(\epsilon_{1}(\tilde{s}))\vert \check{\alpha}(\tilde{s})\vert \vert \check{\alpha}\vert ^{-1}$$
$$ d(\tilde{s})^{-1}\partial_{\underline{U}_{\tilde{s}}}S_{\underline{\tilde{M}}'_{1}(\tilde{s}),\lambda_{1}(\tilde{s})}^{\tilde{G}'_{1}(\tilde{s}),mod}(exp(\underline{\xi}(\tilde{s},rH_{c}))\epsilon_{1}(\tilde{s}),f^{\tilde{G}'_{1}(\tilde{s})}).$$
  Jointe \`a la formule (10), cette \'egalit\'e entra\^{\i}ne
 $$B_{j}(rH_{d},U)-B_{j}(-rH_{d},U)-2\pi i \vert \alpha\vert \underline{B}_{j}(rH_{c},U)=\sum_{\tilde{s}\in\tilde{\zeta}_{j} Z(\hat{M})^{\hat{\alpha}^*}/Z(\hat{G})^{\Gamma_{{\mathbb R}},\hat{\theta}}}i_{\tilde{M}'_{j}}(\tilde{G},\tilde{G}'(\tilde{s}))d(\tilde{s})^{-1}X(\tilde{s},r),$$
 o\`u
$$X(\tilde{s},r)=\partial_{U_{\tilde{s}}} S_{\tilde{M}'_{1}(\tilde{s}),\lambda_{1}(\tilde{s})}^{\tilde{G}'_{1}(\tilde{s})}(exp(r\xi_{j}(H_{d}))\eta,f^{\tilde{G}'_{1}(\tilde{s})})-\partial_{U_{\tilde{s}}} S_{\tilde{M}'_{1}(\tilde{s}),\lambda_{1}(\tilde{s})}^{\tilde{G}'_{1}(\tilde{s})}(exp(-r\xi_{j}(H_{d}))\eta,f^{\tilde{G}'_{1}(\tilde{s})})$$
$$-2\pi i \vert \check{\alpha}(\tilde{s})\vert i_{\tilde{M}'_{1}(\tilde{s})}^{\tilde{G}'_{1}(\tilde{s})}(\epsilon_{1}(\tilde{s})  )\partial_{\underline{U}_{\tilde{s}}} S_{\underline{\tilde{M}}'_{1}(\tilde{s}),\lambda_{1}(\tilde{s})}^{\tilde{G}'_{1}(\tilde{s}),mod}(exp(r\underline{\xi}(\tilde{s},H_{c}))\epsilon_{1}(\tilde{s}),f^{\tilde{G}'_{1}(\tilde{s})}).$$
 On a l'\'egalit\'e $C(\tilde{s})(\xi_{j}(H_{d}))=i\underline{\xi}(\tilde{s},H_{c})$. On est donc dans la situation o\`u l'on peut appliquer la proposition 4.2. Elle  entra\^{\i}ne que 
$$lim_{r\to 0+}X(\tilde{s},r)=0$$
pour tout $\tilde{s}$ apparaissant ci-dessus. Donc
$$lim_{r\to 0+}\left(B_{j}(rH_{d},U)-B_{j}(-rH_{d},U)-2\pi i \vert \alpha\vert \underline{B}_{j}(rH_{c},U)\right)=0.$$
Cela \'etant vrai pour tout $j$, on a aussi
$$lim_{r\to 0+}\left(B(rH_{d},U)-B(-rH_{d},U)-2\pi i \vert \alpha\vert \underline{B}(rH_{c},U)\right)=0.$$
 Consid\'erons l'\'egalit\'e (9). Parce que $\phi$ est $C^{\infty}$, on a
 $$lim_{r\to 0+}\left(\partial_{U}\phi(rH_{d})-\partial_{U}\phi(-rH_{d})\right)=0.$$
 On en d\'eduit
 $$lim_{r\to 0+}(\partial_{U}I_{K\tilde{M}}^{K\tilde{G},mod}(exp(rH_{d})\eta,\omega,f)-\partial_{U}I_{K\tilde{M}}^{K\tilde{G},mod}(exp(-rH_{d})\eta,\omega,f)$$
 $$-2\pi i\vert \alpha\vert \underline{B}(rH_{c},U))=0.$$
 En utilisant (20), on obtient la premi\`ere \'egalit\'e de l'assertion (iv) de l'\'enonc\'e qu'il restait \`a prouver. $\square$
 
 \bigskip
 
 \section{Des variantes de l'application $\phi_{\tilde{M}}$}

\bigskip

\subsection{Normalisation partielle des op\'erateurs d'entrelacement}
Pour la suite de l'article, on suppose fix\'es une paire parabolique  minimale $(\tilde{P}_{0},\tilde{M}_{0})$
  et un sous-groupe compact maximal $K$ de $G({\mathbb R})$ en bonne position relativement \`a $M_{0}$.
  
Supposons fix\'ees des mesures de Haar sur tous les groupes intervenant. Soit $M\in {\cal L}(M_{0})$ et $\pi$ une repr\'esentation de $M({\mathbb R})$ irr\'eductible et temp\'er\'ee dans un espace $V_{\pi}$. Pour $\lambda\in {\cal A}_{M{\mathbb C}}^*$, on d\'efinit la repr\'esentation $\pi_{\lambda}$ par $\pi_{\lambda}(m)=e^{<\lambda,H_{M}(m)>}\pi(m)$. Soient $P,P'\in {\cal P}(M)$. On sait d\'efinir l'op\'erateur d'entrelacement $J_{P'\vert P}(\pi_{\lambda}):Ind_{P}^G(\pi_{\lambda})\to Ind_{P'}^G(\pi_{\lambda})$. Il est m\'eromorphe en $\lambda$. On sait normaliser cet op\'erateur, cf. [A4] th\'eor\`eme 2.1.   L'op\'erateur normalis\'e $R_{P'\vert P}(\pi_{\lambda})$ est \'egal \`a $r_{P'\vert P}(\pi_{\lambda})^{-1}J_{P'\vert P}(\pi_{\lambda})$, o\`u $\lambda\mapsto r_{P'\vert P}(\pi_{\lambda})$ est une fonction m\'eromorphe \`a valeurs dans ${\mathbb C}$. Nous allons d\'efinir une normalisation partielle que l'on pourra plus facilement stabiliser. 

Supposons que $\pi$ soit de la s\'erie discr\`ete. Cela implique que $M$ poss\`ede un sous-tore maximal elliptique et on fixe un tel sous-tore $T$. L'action de $\Gamma_{{\mathbb R}}\simeq\{\pm 1\}$ d\'ecompose $X_{*}(T)_{{\mathbb Z}}\otimes {\mathbb R}$ en
$$X_{*}(T)\otimes_{{\mathbb Z}}  {\mathbb R}=(X_{*}(T)\otimes_{{\mathbb Z}}{\mathbb R})^+\oplus (X_{*}(T)\otimes_{{\mathbb Z}} {\mathbb R})^-,$$
o\`u, pour $\epsilon=\pm 1$,  $X_{*}(T)\otimes_{{\mathbb Z}} {\mathbb R})^{\epsilon}$ est le sous-espace o\`u l'\'el\'ement non trivial $\boldsymbol{\sigma}$ de $\Gamma_{{\mathbb R}}$ agit par multiplication par $\epsilon$. L'hypoth\`ese d'ellipticit\'e signifie que 
$$(X_{*}(T)\otimes_{{\mathbb Z}} {\mathbb R})^+=\mathfrak{a}_{M}({\mathbb R}).$$
Posons $\mathfrak{h}^M_{{\mathbb R}}=(X_{*}(T)\otimes_{{\mathbb Z}} {\mathbb R})^- $. 
On a
$$\mathfrak{t}({\mathbb R})= \mathfrak{a}_{M}({\mathbb R})\oplus i\mathfrak{h}^M_{{\mathbb R}}$$
tandis que l'alg\`ebre $\mathfrak{h}_{{\mathbb R}}$ usuelle est isomorphe \`a
$$\mathfrak{h}_{{\mathbb R}}=\mathfrak{a}_{M}({\mathbb R})\oplus \mathfrak{h}^M_{{\mathbb R}}.$$
Le param\`etre $\mu(\pi)$ est la classe de conjugaison par le groupe de Weyl $W^M$ relatif \`a $T$ d'un \'el\'ement de $i\mathfrak{t}({\mathbb R})$. 
 Supposons que le caract\`ere central $\omega_{\pi}$ de $\pi$ soit trivial sur $\mathfrak{A}_{M}$ (on rappelle que ce groupe est la composante neutre pour la topologie r\'eelle de $A_{M}({\mathbb R})$).  Alors $\mu(\pi)$ est la classe de conjugaison par $W^M$ d'un \'el\'ement de $\mathfrak{h}^M_{{\mathbb R}}$. On fixe un \'el\'ement de cette classe. Notons $\Sigma(T)$ l'ensemble des racines de $T$ dans $G$. Pour $P,P'\in {\cal P}(M)$, notons  $\Sigma^{U_{P}}(T)$ le sous-ensemble des racines qui apparaissent dans l'alg\`ebre de Lie du radical unipotent de $P$ et posons
 $$\Sigma_{P'\vert P}(T)=\Sigma^{U_{P}}(T)-(\Sigma^{U_{P}}(T)\cap \Sigma^{U_{P'}}(T)).$$
 Le groupe $\Gamma_{{\mathbb R}}$ agit naturellement sur ces ensembles. Soit $\alpha\in \Sigma^{U_{P}}(T)$ dont l'orbite $(\alpha)$ pour cette action ait deux \'el\'ements, autrement dit $\boldsymbol{\sigma}(\alpha)\not=\alpha$. Quitte \`a \'echanger $\alpha$ et $\boldsymbol{\sigma}(\alpha)$, on peut supposer $<\mu(\pi),\boldsymbol{\sigma}(\check{\alpha})-\check{\alpha}>\leq 0$. Pour $\lambda\in {\cal A}_{M,{\mathbb C}}^*$, on pose
 $$r_{(\alpha)}(\pi_{\lambda})=2\pi <\mu(\pi)+\lambda,\check{\alpha}>^{-1}.$$
 Soit $\alpha\in \Sigma^{U_{P}}(T)$ dont l'orbite soit r\'eduite \`a un \'el\'ement, autrement dit $\boldsymbol{\sigma}(\alpha)=\alpha$ ($\alpha$ est "r\'eelle"). Remarquons qu'alors $<\mu(\pi),\check{\alpha}>=0$.  On pose
 $$r_{\alpha}(\pi_{\lambda})=\Gamma_{{\mathbb R}}(<\lambda,\check{\alpha}>+N_{\alpha})\Gamma_{{\mathbb R}}(<\lambda,\check{\alpha}>+N_{\alpha}+1)^{-1},$$
 o\`u $\Gamma_{{\mathbb R}}$ est la fonction usuelle
 $$\Gamma_{{\mathbb R}}(s)=\pi^{-s/2}\Gamma(s/2),$$
 et o\`u $N_{\alpha}$ est un \'el\'ement de $\{0,1\}$ qui d\'epend de $\pi$ et $\alpha$. D'apr\`es [A4] paragraphe 3, le facteur de normalisation $r_{P'\vert P}(\pi_{\lambda})$ est le produit des $r_{(\alpha)}(\pi_{\lambda})$ pour les orbites $(\alpha)$ \`a deux \'el\'ements  contenues dans $\Sigma_{P'\vert P}(T)$ et des $r_{\alpha}(\pi_{\lambda})$ pour les $\alpha$ r\'eelles contenues dans le m\^eme ensemble.

 D\'efinissons la fonction d'une variable complexe
 $$\Gamma_{0}(s)=\Gamma_{{\mathbb R}}(s)\Gamma_{{\mathbb R}}(s+1)^{-1}.$$
 D'apr\`es la formule usuelle $\Gamma(s+1)=s\Gamma(s)$, on a
$$\Gamma_{0}(s+1)=\frac{2\pi}{s}\Gamma_{0}(s)^{-1}.$$
 Consid\'erons une racine $\alpha$ r\'eelle. Posons $\epsilon(\alpha)=(-1)^{N_{\alpha}}$ et
 $$\rho_{\alpha}(\pi_{\lambda})=\Gamma_{0}(<\lambda,\check{\alpha}>)^{\epsilon(\alpha)} .$$
 D'apr\`es la formule  ci-dessus, on a
$$ r_{\alpha}(\pi_{\lambda})=\rho_{\alpha}(\pi_{\lambda})(\frac{2\pi}{<\lambda,\check{\alpha}>})^{N_{\alpha}}.$$
On d\'efinit $\rho_{P'\vert P}(\pi_{\lambda})$ comme les produit des $\rho_{\alpha}(\pi_{\lambda})$ sur les racines r\'eelles $\alpha$ contenues dans $\Sigma_{P'\vert P}(T)$. On pose
$$r^{rat}_{P'\vert P}(\pi_{\lambda})=\rho_{P'\vert P}(\pi_{\lambda})^{-1}r_{P'\vert P}(\pi_{\lambda}).$$
Il r\'esulte des formules ci-dessus que $r^{rat}_{P'\vert P}(\pi_{\lambda})$ est une fonction rationnelle en $\lambda$. Remarquons qu'au contraire, la fonction $\rho_{P'\vert P}(\pi_{\lambda})$ peut avoir une infinit\'e d'hyperplans polaires. On pose
$$R^{rat}_{P'\vert P}(\pi_{\lambda})=\rho_{P'\vert P}(\pi_{\lambda})^{-1}J_{P'\vert P}(\pi_{\lambda})=r^{rat}_{P'\vert P}(\pi_{\lambda})R_{P'\vert P}(\pi_{\lambda}).$$
Cet op\'erateur est rationnel en $\lambda$. On entend par l\`a la propri\'et\'e suivante.  On peut r\'ealiser les espaces des induites $Ind_{P}^G(\pi_{\lambda})$ et $Ind_{P'}^G(\pi_{\lambda})$ comme des espaces de fonctions sur $K$ ind\'ependants de $\lambda$. Ils sont munis de produits hermitiens d\'efinis positifs. Soient $e\in Ind_{P}^G(\pi_{\lambda})$ et $e'\in Ind_{P'}^G(\pi_{\lambda})$ deux \'el\'ements $K$-finis. Alors le produit hermitien
$$(e',R^{rat}_{P'\vert P}(\pi_{\lambda})(e))$$
est une fonction rationnelle en $\lambda$. On peut pr\'eciser que ses hyperplans polaires sont d'\'equations $<\lambda,\check{\alpha}>=c$ pour des racines $\alpha$ de $A_{M}$ dans $G$ (la d\'efinition exacte de $\check{\alpha}$ n'ayant pas beaucoup d'importance). On prendra garde toutefois que les p\^oles d\'ependent des $K$-types selon lesquels se transforment $e$ et $e'$. En g\'en\'eral, il n'y a pas un nombre fini d'hyperplans hors desquels la fonction ci-dessus n'a pas de p\^ole quels que soient $e$ et $e'$.
 
 Consid\'erons maintenant le cas g\'en\'eral o\`u $\pi$ est une repr\'esentation irr\'eductible temp\'er\'ee quelconque. On peut fixer 
 
 -un groupe de Levi $R$ de $M$ contenant $M_{0}$ et un sous-groupe parabolique $S\in {\cal P}^M(R)$;
 
 - une repr\'esentation $\sigma$ de $R({\mathbb R})$, irr\'eductible, de la s\'erie discr\`ete et telle que $\omega_{\sigma}$ soit trivial sur $\mathfrak{A}_{R}$;
 
 - un \'el\'ement $\lambda_{0}\in i{\cal A}_{R}^*$;
 
 \noindent de sorte que $\pi$ soit une sous-repr\'esentation de l'induite $Ind_{S}^M(\sigma_{\lambda_{0}})$. Soient $P,P'\in {\cal P}(M)$ et $\lambda\in {\cal A}_{M,{\mathbb C}}^*$. Notons
$Q$ et $Q'$ les \'el\'ements de ${\cal P}(R)$ tels que $Q\subset P$, $Q'\subset P'$, $Q\cap M=Q'\cap M=S$. On pose
$$\rho_{P'\vert P}(\pi_{\lambda})=\rho_{Q'\vert Q}(\sigma_{\lambda_{0}+\lambda}),$$
$$r^{rat}_{P'\vert P}(\pi_{\lambda})=r^{rat}_{Q'\vert Q}(\sigma_{\lambda_{0}+\lambda}),$$
et on d\'efinit comme ci-dessus
$$R^{rat}_{P'\vert P} (\pi_{\lambda})=\rho_{P'\vert P}(\pi_{\lambda})^{-1}J_{P'\vert P}(\pi_{\lambda})=r^{rat}_{P'\vert P}(\pi_{\lambda})R_{P'\vert P}(\pi_{\lambda}).$$
Ces termes ont les m\^emes propri\'et\'es que ci-dessus.

On a suppos\'e $\pi$ irr\'eductible. Remarquons que l'on peut poser les m\^emes d\'efinitions pour une repr\'esentation $\pi$ r\'eductible qui est une sous-repr\'esentation d'une induite  $Ind_{S}^M(\sigma_{\lambda_{0}})$ comme ci-dessus.

\bigskip

\subsection{Caract\`eres pond\'er\'es rationnels}

Soient $\tilde{M}\in {\cal L}(\tilde{M}_{0})$, $\tilde{R}$ un espace de Levi de $\tilde{M}$ contenant $\tilde{M}_{0}$ et $\tilde{\sigma}$ une $\omega$-repr\'esentation elliptique de $\tilde{R}({\mathbb R})$. Fixons un espace parabolique $\tilde{S}\in {\cal P}^{\tilde{M}}(\tilde{R})$ et posons $\tilde{\pi}=Ind_{\tilde{S}}^{\tilde{M}}(\tilde{\sigma})$. Fixons $\tilde{P}\in {\cal P}(\tilde{M})$ et soit $\lambda\in {\cal A}_{\tilde{M},{\mathbb C}}^*$. A la suite d'Arthur, on a d\'efini en [W2] 2.7 un op\'erateur ${\cal M}_{\tilde{M}}^{\tilde{G}}(\pi_{\lambda})$ de l'espace de la repr\'esentation $Ind_{\tilde{P}}^{\tilde{G}}(\pi_{\lambda})$, qui est m\'eromorphe en $\lambda$. Rappelons sa d\'efinition. La repr\'esentation sous-jacente $\pi$ de $M({\mathbb R})$ n'est pas irr\'eductible en g\'en\'eral, mais v\'erifie la derni\`ere propri\'et\'e du paragraphe pr\'ec\'edent. Pour $\tilde{P}'\in {\cal P}(\tilde{M})$, on peut donc d\'efinir les termes $r_{P'\vert P}(\pi_{\lambda})$ etc... On pose $\mu_{P'\vert P}(\pi_{\lambda})=r_{P\vert P'}(\pi_{\lambda})^{-1}r_{P'\vert P}(\pi_{\lambda})^{-1}$. Pour $\Lambda\in i{\cal A}_{\tilde{M}}^*$, on d\'efinit l'op\'erateur 
$${\cal M}(\pi_{\lambda};\Lambda,\tilde{P}')=\mu_{P'\vert P}(\pi_{\lambda})^{-1}\mu_{P'\vert P}(\pi_{\lambda+\Lambda/2})J_{P'\vert P}(\pi_{\lambda})^{-1}J_{P'\vert P}(\pi_{\lambda+\Lambda}).$$
La famille $({\cal M}(\pi_{\lambda};\Lambda,\tilde{P}'))_{\tilde{P}'\in {\cal P}(\tilde{M})}$ est une $(\tilde{G},\tilde{M})$-famille \`a valeurs op\'erateurs. L'op\'erateur ${\cal M}_{\tilde{M}}^{\tilde{G}}(\pi_{\lambda})$ est d\'eduit de cette famille de la fa\c{c}on habituelle.

Avec les m\^emes notations que ci-dessus, posons  $\mu^{rat}_{P'\vert P}(\pi_{\lambda})=r^{rat}_{P\vert P'}(\pi_{\lambda})^{-1}r^{rat}_{P'\vert P}(\pi_{\lambda})^{-1}$ et d\'efinissons l'op\'erateur
$${\cal M}^{rat}(\pi_{\lambda};\Lambda,\tilde{P}')=\mu^{rat}_{P'\vert P}(\pi_{\lambda})^{-1}\mu^{rat}_{P'\vert P}(\pi_{\lambda+\Lambda/2})R^{rat}_{P'\vert P}(\pi_{\lambda})^{-1}R^{rat}_{P'\vert P}(\pi_{\lambda+\Lambda}).$$
La  famille $({\cal M}^{rat}(\pi_{\lambda};\Lambda,\tilde{P}'))_{\tilde{P}'\in {\cal P}(\tilde{M})}$ est encore une $(\tilde{G},\tilde{M})$-famille \`a valeurs op\'erateurs. On en d\'eduit encore un op\'erateur que l'on note  ${\cal M}_{\tilde{M}}^{rat,\tilde{G}}(\pi_{\lambda})$. Il est clair qu'il est rationnel en $\lambda$. Ses hyperplans polaires sont d'\'equations $<\lambda,\check{\alpha}>=c$ pour les racines $\alpha$ de $A_{\tilde{M}}$ dans $G$. On a

(1) l'op\'erateur  ${\cal M}_{\tilde{M}}^{rat,\tilde{G}}(\pi_{\lambda})$ n'a pas de p\^ole pour $\lambda\in i{\cal A}_{\tilde{M}}^*$.

En effet, la d\'emonstration de [A5] proposition 2.3 s'applique. 

Rappelons que  l'on a fix\'e une fonction $H_{\tilde{M}}:\tilde{M}({\mathbb R})\to {\cal A}_{\tilde{M}}$ qui permet de d\'efinir la $\omega$-repr\'esentation $\tilde{\pi}_{\lambda}$, cf. [IV] 1.1.  Pour $f\in C_{c}^{\infty}(\tilde{G}({\mathbb R}))$, on a d\'efini le caract\`ere pond\'er\'e de $\tilde{\pi}_{\lambda}$ par
$$J_{\tilde{M}}^{\tilde{G}}(\tilde{\pi}_{\lambda},f)=trace( {\cal M}_{\tilde{M}}^{\tilde{G}}(\pi_{\lambda})Ind_{\tilde{P}}^{\tilde{G}}(\tilde{\pi}_{\lambda},f)).$$
On d\'efinit de m\^eme le caract\`ere pond\'er\'e "rationnel"
$$J_{\tilde{M}}^{rat,\tilde{G}}(\tilde{\pi}_{\lambda},f)=trace( {\cal M}_{\tilde{M}}^{rat, \tilde{G}}(\pi_{\lambda})Ind_{\tilde{P}}^{\tilde{G}}(\tilde{\pi}_{\lambda},f)).$$
Il v\'erifie les m\^emes propri\'et\'es que le pr\'ec\'edent et on utilisera ces propri\'et\'es sans plus de commentaires. Il n'est pas rationnel en $\lambda$ car l'op\'erateur $Ind_{\tilde{P}}^{\tilde{G}}(\tilde{\pi}_{\lambda},f)$ n'est pas rationnel. Mais ce dernier est holomorphe en $\lambda$ et \`a d\'ecroissance rapide dans les bandes verticales, au sens suivant.
 Pour un espace vectoriel r\'eel $V$ de dimension finie et pour une fonction $\psi$ m\'eromorphe  sur le complexifi\'e $V_{{\mathbb C}}$, on dit que $\psi$ est \`a d\'ecroissance rapide dans les bandes verticales si, 
 pour tout sous-ensemble compact $\Gamma\subset  V$ tel que cette fonction $\psi$ n'ait pas de p\^ole sur $\Gamma+i V$, celle-ci est \`a d\'ecroissance rapide sur cet ensemble.  Les propri\'et\'es de l'op\'erateur 
${\cal M}_{\tilde{M}}^{rat, \tilde{G}}(\pi_{\lambda})$  entra\^{\i}nent  les propri\'et\'es suivantes:

(2)  $J_{\tilde{M}}^{rat,\tilde{G}}(\tilde{\pi}_{\lambda},f)$ n'a qu'un nombre fini d'hyperplans polaires;

(3) chacun de ces hyperplans est d'\'equation $<\lambda,\check{\alpha}>=c$ o\`u  $\alpha$ est une racine de $A_{\tilde{M}}$ dans $G$ et $c\in {\mathbb C}$;

(4) $\lambda\mapsto J_{\tilde{M}}^{rat,\tilde{G}}(\tilde{\pi}_{\lambda},f)$ est \`a d\'ecroissance rapide dans les bandes verticales.

Les hyperplans polaires d\'ependent de $f$. Toutefois, on a le raffinement suivant. Rappelons que, pour un ensemble fini $\Omega$ de $K$-types, on note $C_{c}^{\infty}(\tilde{G}({\mathbb R}),\Omega)$ le sous-espace des $f\in C_{c}^{\infty}(\tilde{G}({\mathbb R}))$ qui se transforment \`a droite et \`a gauche selon des $K$-types appartenant \`a $\Omega$, cf. [IV] 3.4. Alors

(5) soit $\Omega$ un ensemble fini de $K$-types; alors il existe un ensemble  fini d'hyperplans de la forme  (3) tel que, pour $f\in C_{c}^{\infty}(\tilde{G}({\mathbb R}),\Omega)$, les hyperplans polaires de  $J_{\tilde{M}}^{rat,\tilde{G}}(\tilde{\pi}_{\lambda},f)$ appartiennent \`a cet ensemble.

 Le terme  $J_{\tilde{M}}^{rat,\tilde{G}}(\tilde{\pi}_{\lambda},f)$ n'a pas de p\^ole pour $\lambda\in i{\cal A}_{\tilde{M}}^*$.  Pour $X\in {\cal A}_{\tilde{M}}$, on peut d\'efinir
$$J_{\tilde{M}}^{rat,\tilde{G}}(\tilde{\pi},X,f)=\int_{i{\cal A}_{\tilde{M}}^*}J_{\tilde{M}}^{rat,\tilde{G}}(\tilde{\pi}_{\lambda},f)e^{<-\lambda,X>}\,d\lambda.$$
On obtient une fonction de Schwartz en $X$. On a d\'efini en [W2] 6.4 l'espace $C_{ac}^{\infty}(\tilde{G}({\mathbb R}))$. C'est celui des fonctions $f:\tilde{G}({\mathbb R})\to {\mathbb C}$ telles que $f(b\circ H_{\tilde{G}})\in C_{c}^{\infty}(\tilde{G}({\mathbb R}))$ pour toute fonction $b\in C_{c}^{\infty}({\cal A}_{\tilde{G}})$. Comme en [W2] 6.4, on peut \'etendre la d\'efinition de $J_{\tilde{M}}^{rat,\tilde{G}}(\tilde{\pi},X,f)$ par continuit\'e \`a $f\in C_{ac}^{\infty}(\tilde{G}({\mathbb R}))$.

On a impos\'e \`a $\tilde{\pi}$ d'\^etre une induite $\tilde{\pi}=Ind_{\tilde{S}}^{\tilde{M}}(\tilde{\sigma})$ o\`u $\tilde{\sigma}$ est une $\omega$-repr\'esentation elliptique de $\tilde{R}({\mathbb R})$. Par lin\'earit\'e, les termes $J_{\tilde{M}}^{rat,\tilde{G}}(\tilde{\pi}_{\lambda},f)$ et $J_{\tilde{M}}^{rat,\tilde{G}}(\tilde{\pi},X,f)$ s'\'etendent \`a $\tilde{\pi}\in Ind_{\tilde{R}}^{\tilde{M}}(D_{ell}(\tilde{R}({\mathbb R}),\omega))$, cf. [IV] 1.2.  Il est plus ou moins clair que ces termes sont invariants par l'action du groupe $W^{M}(\tilde{M}_{0})$ sur cet espace. En utilisant la version temp\'er\'ee de la d\'ecomposition [IV] 1.2(2), on peut alors \'etendre ces termes par lin\'earit\'e \`a tout $\tilde{\pi}\in D_{temp}(\tilde{M}({\mathbb R}),\omega)$.

\bigskip

\subsection{L'application $\phi_{\tilde{M}}^{rat,\tilde{G}}$}
On conserve la m\^eme situation. On vient de rappeler la d\'efinition de l'espace $C_{ac}^{\infty}(\tilde{G}({\mathbb R}))$. On note $C_{c}^{\infty}(\tilde{G}({\mathbb R}),K)$, resp. $C_{ac}^{\infty}(\tilde{G}({\mathbb R}),K)$, le sous-espace des \'el\'ements de $ C_{c}^{\infty}(\tilde{G}({\mathbb R}))$, resp. de $C_{ac}^{\infty}(\tilde{G}({\mathbb R}))$, qui sont $K$-finis \`a droite et \`a gauche. On note $I(\tilde{G}({\mathbb R}),\omega,K)$, resp. $I_{ac}(\tilde{G}({\mathbb R}),\omega)$,  $I_{ac}(\tilde{G}({\mathbb R}),\omega,K)$, le quotient de $C_{c}^{\infty}(\tilde{G}({\mathbb R}),K)$, resp. $C_{ac}^{\infty}(\tilde{G}({\mathbb R}))$, $C_{ac}^{\infty}(\tilde{G}({\mathbb R}),K)$, par le sous-espace des \'el\'ements $f$ tels que $I^{\tilde{G}}(\gamma,\omega,f)=0$ pour tout $\gamma\in \tilde{G}({\mathbb R})$ fortement r\'egulier. On pose $K^M=K\cap M({\mathbb R})$.

\ass{Proposition}{Il existe une unique application lin\'eaire
$$\phi_{\tilde{M}}^{rat,\tilde{G}}:C_{ac}^{\infty}(\tilde{G}({\mathbb R}))\to I_{ac}(\tilde{M}({\mathbb R}))$$
telle que, pour tout $\tilde{\pi}\in D_{temp}(\tilde{M}({\mathbb R}),\omega)$ et tout $X\in {\cal A}_{\tilde{M}}$, on ait l'\'egalit\'e
$$I^{\tilde{M}}(\tilde{\pi},X,\phi_{\tilde{M}}^{rat,\tilde{G}}(f))=J_{\tilde{M}}^{rat,\tilde{G}}(\tilde{\pi},X,f).$$
L'application $\phi_{\tilde{M}}^{rat,\tilde{G}}$ envoie $C_{ac}^{\infty}(\tilde{G}({\mathbb R}),K)$ dans $I_{ac}(\tilde{M}({\mathbb R}),K^M)$.}

Preuve. La preuve de la premi\`ere assertion est la m\^eme que celle de la proposition 6.4 de [W2]. Supposons $f$ $K$-finie \`a droite et \`a gauche. Elle se d\'ecompose donc selon un nombre fini de $K$-types, qui eux-m\^emes se d\'ecomposent en un nombre fini de $K^M$-types. Si  $\tilde{\pi}$ est l'induite d'une $\omega$-repr\'esentation elliptique d'un espace de Levi de $\tilde{M}$, il r\'esulte des constructions que $J_{\tilde{M}}^{rat,\tilde{G}}(\tilde{\pi},X,f)=0$ si aucun de ces $K^M$-types n'appara\^{\i}t dans $\pi$. Le th\'eor\`eme de Paley-Wiener de Delorme-Mezo (repris en [W2] 6.2) conduit alors \`a la seconde assertion. $\square$

De la d\'efinition r\'esulte que

(1) soient $f\in C_{c}^{\infty}(\tilde{G}({\mathbb R}))$ et $\tilde{\pi}\in D_{temp}(\tilde{M}({\mathbb R}),\omega)$; alors la fonction $X\mapsto I^{\tilde{M}}(\tilde{\pi},X,\phi_{\tilde{M}}^{rat,\tilde{G}}(f))$ est de Schwartz sur ${\cal A}_{\tilde{M}}$.

Comme toujours, on \'elimine les choix de mesures de Haar en d\'efinissant plus canoniquement une application
$$\phi_{\tilde{M}}^{rat,\tilde{G}}:C_{ac}^{\infty}(\tilde{G}({\mathbb R}))\otimes Mes(G({\mathbb R}))\to I_{ac}(\tilde{M}({\mathbb R}))\otimes Mes(M({\mathbb R})).$$

\bigskip

\subsection{Relation entre les applications $\phi_{\tilde{M}}^{\tilde{G}}$ et $\phi_{\tilde{M}}^{rat,\tilde{G}}$}
   Fixons un sous-tore maximal $T$ de $M$. Pour $\beta\in \Sigma(T)$, la coracine $\check{\beta}$ peut \^etre consid\'er\'ee comme une forme lin\'eaire sur $X^*(T)\otimes _{{\mathbb Z}}{\mathbb R}$ et se restreint en une forme lin\'eaire $\check{\beta}_{\tilde{M}}$ sur ${\cal A}_{\tilde{M}}^*$, autrement dit en un \'el\'ement $\check{\beta}_{\tilde{M}}\in {\cal A}_{\tilde{M}}$.  Notons $\check{\Sigma}_{\star }(A_{\tilde{M}})$ l'ensemble de ces restrictions $\check{\beta}_{\tilde{M}}$ qui sont non nulles. Cet ensemble ne d\'epend pas du choix de $T$. Soit $\tilde{\pi}$ une $\omega$-repr\'esentation de $\tilde{M}({\mathbb R})$, supposons comme en 5.1 que $\pi$ est une sous-repr\'esentation d'une induite $Ind_{S}^M(\sigma_{\lambda_{0}})$ o\`u $\sigma$ est de la s\'erie discr\`ete de $R({\mathbb R})$ et $\omega_{\sigma}=1$ sur $\mathfrak{A}_{R}$.  Supposons  que $T$ est un sous-tore maximal elliptique de $R$. Soit $\check{\alpha}\in\check{\Sigma}_{\star }(A_{\tilde{M}})$. Notons $E_{\check{\alpha}}(\pi)$ l'ensemble des $\beta\in \Sigma(T)$ qui sont r\'eelles et telles que $\check{\beta}_{\tilde{M}}=\check{\alpha}$. Cet ensemble d\'epend de $\pi$ car la notion de racine r\'eelle d\'epend de la structure sur ${\mathbb R}$ de $T$; ce tore lui-m\^eme d\'epend de $R$ qui depend de $\pi$. Pour tout $\beta\in E_{\check{\alpha}}(\pi)$, on a d\'efini la fonction $\rho_{\beta}(\pi_{\lambda})$ pour $\lambda\in {\cal A}_{M,{\mathbb C}}^*$, a fortiori pour $\lambda\in {\cal A}_{\tilde{M},{\mathbb C}}^*$. Pour un tel $\lambda$, elle est de la forme
   $$\rho_{\beta}(\pi_{\lambda})=\Gamma_{0}(<\lambda_{0},\check{\beta}>+<\lambda,\check{\alpha}>)^{\epsilon(\beta)}$$
   o\`u $\epsilon(\beta)\in \{\pm 1\}$. 
   Nous noterons maintenant $e$ au lieu de $\beta$ les \'el\'ements de $E_{\check{\alpha}}(\pi)$. Pour $e=\beta$, on pose $\epsilon(e)=\epsilon(\beta)$ et $\mu(e)=<\lambda_{0},\check{\beta}>$. Remarquons que $\mu(e)$ est imaginaire. La formule pr\'ec\'edente devient
   $$\rho_{e}(\pi_{\lambda})=\Gamma_{0}(\mu(e)+<\lambda,\check{\alpha}>)^{\epsilon(e)}.$$
   On pose
   $$\rho_{\check{\alpha}}(\pi_{\lambda})=\prod_{e\in E_{\check{\alpha}}(\pi)}\rho_{e}(\pi_{\lambda}).$$
   Remarquons qu'il y a une bijection naturelle de $E_{-\check{\alpha}}(\pi)$ sur $E_{\check{\alpha}}(\pi)$, que l'on peut noter $e\mapsto -e$, de sorte que $\mu(-e)=-\mu(e)$ et $\epsilon(-e)=\epsilon(e)$. On a ainsi la formule
   $$\rho_{-\check{\alpha}}(\pi_{\bar{\lambda}})=\overline{\rho_{\check{\alpha}}(\pi_{\lambda})}.$$
Pour $\tilde{P},\tilde{P}'\in {\cal P}(\tilde{M})$,  on d\'efinit les sous-ensembles $\check{\Sigma}_{\star }^{U_{P}}(A_{\tilde{M}})$ et $\check{\Sigma}_{\star ,P'\vert P}(A_{\tilde{M}})$ en imitant les d\'efinitions de 5.1. Il r\'esulte des constructions de ce paragraphe que 
$$(1) \qquad \rho_{P'\vert P}(\pi_{\lambda})=\prod_{\check{\alpha}\in \check{\Sigma}_{\star ,P'\vert P}(A_{\tilde{M}}) } \rho_{\check{\alpha}}(\pi_{\lambda}).$$

Consid\'erons $\tilde{P}$ comme fix\'e. Pour $\Lambda\in i{\cal A}_{\tilde{M}}^*$, posons
$$(2) \qquad \rho(\pi_{\lambda}:\Lambda,\tilde{P}')=\rho_{P\vert P'}(\pi_{\lambda}) \rho_{P\vert P'}(\pi_{\lambda+\Lambda/2})^{-1}\rho_{P'\vert P}(\pi_{\lambda+\Lambda/2})^{-1} \rho_{P'\vert P}(\pi_{\lambda+\Lambda}).$$
Il r\'esulte des d\'efinitions que l'on a l'\'egalit\'e
$${\cal M}(\pi_{\lambda};\Lambda,\tilde{P}')=\rho(\pi_{\lambda}:\Lambda,\tilde{P}'){\cal M}^{rat}(\pi_{\lambda};\Lambda,\tilde{P}').$$
La famille $(\rho(\pi_{\lambda}:\Lambda,\tilde{P}'))_{\tilde{P}'\in {\cal P}(\tilde{M})}$ est une $(\tilde{G},\tilde{M})$-famille. Les descriptions ci-dessus montrent qu'elle est d'une forme particuli\`ere pour laquelle le corollaire 6.5 de  [A6] s'applique. C'est-\`a-dire que l'on a l'\'egalit\'e
$${\cal M}_{\tilde{M}}^{\tilde{G}}(\pi_{\lambda})=\sum_{\tilde{L}\in {\cal L}(\tilde{M})}\rho_{\tilde{M}}^{\tilde{L}}(\pi_{\lambda}){\cal M}_{\tilde{L}}^{rat,\tilde{G}}(\pi_{\lambda}).$$
En utilisant les propri\'et\'es usuelles de commutation des op\'erateurs d'entrelacement \`a l'induction, on obtient la formule suivante, pour tout $f\in C_{c}^{\infty}(\tilde{G}({\mathbb R}))$:
$$(3) \qquad J_{\tilde{M}}^{\tilde{G}}(\tilde{\pi}_{\lambda},f)=\sum_{\tilde{L}\in {\cal L}(\tilde{M})}\rho_{\tilde{M}}^{\tilde{L}}(\pi_{\lambda})J_{\tilde{L}}^{rat,\tilde{G}}(Ind_{\tilde{Q}}^{\tilde{L}}(\tilde{\pi}_{\lambda}),f),$$
o\`u, pour tout $\tilde{L}$, on a fix\'e un espace parabolique $\tilde{Q}\in {\cal P}^{\tilde{L}}(\tilde{M})$. 

Calculons la fonction $\rho_{\tilde{M}}^{\tilde{G}}(\pi_{\lambda})$. Consid\'erons l'ensemble des ensembles $\underline{\check{\alpha}}=\{\check{\alpha}_{i};i=1,...,n\}$ o\`u les $\check{\alpha}_{i}$ sont des \'el\'ements lin\'eairements ind\'ependants de $\check{\Sigma}_{\star }(A_{\tilde{M}})$ et o\`u $n=a_{\tilde{M}}-a_{\tilde{G}}$. Disons que deux tels ensembles $\underline{\check{\alpha}}=\{\check{\alpha}_{i};i=1,...,n\}$ et $\underline{\check{\alpha}}'=\{\check{\alpha}'_{i};i=1,...,n\}$ sont \'equivalents si, quitte \`a changer leur num\'erotation, on a $\check{\alpha}'_{i}=\pm \check{\alpha}_{i}$ pour tout $i$. Notons ${\cal J}_{\tilde{M}}^{\tilde{G}}$ l'ensemble des classes d'\'equivalence. Puisqu'on a fix\'e un espace parabolique $\tilde{P}$, on peut identifier ${\cal J}_{\tilde{M}}^{\tilde{G}}$ \`a l'ensemble des $\underline{\check{\alpha}}=\{\check{\alpha}_{i};i=1,...,n\}$ comme ci-dessus  tels que $\check{\alpha}_{i}\in \check{\Sigma}_{\star }^{U_{P}}(A_{\tilde{M}})$ et c'est ce que nous faisons pour quelque temps. Pour un \'el\'ement $\underline{\check{\alpha}}=\{\check{\alpha}_{i};i=1,...,n\}$ de cet ensemble, on note $m(\underline{\check{\alpha}})$ le volume du quotient de ${\cal A}_{\tilde{M}}^{\tilde{G}}$ par le ${\mathbb Z}$-module engendr\'e par les $\check{\alpha}_{i}$. 
Pour $\check{\alpha}\in \check{\Sigma}_{\star }(A_{\tilde{M}})$, posons
$$c_{\check{\alpha}}(\pi_{\lambda};\Lambda)=\rho_{-\check{\alpha}}(\pi_{\lambda})\rho_{\check{\alpha}}(\pi_{\lambda+\Lambda/2})^{-1}\rho_{-\check{\alpha}}(\pi_{\lambda+\Lambda/2})^{-1}\rho_{\check{\alpha}}(\pi_{\lambda+\Lambda}).$$
Consid\'erons pour un instant $\lambda$ comme une constante. Pour $e\in E_{\check{\alpha}}(\pi)$, d\'efinissons la fonction $c_{e}$ d'une variable complexe $s$ par
$$c_{e}(s)=\Gamma_{0}(-\mu(e)-<\lambda,\check{\alpha}>)^{\epsilon(e)}\Gamma_{0}(\mu(e)+<\lambda,\check{\alpha}>+s/2)^{-\epsilon(e)}$$
$$\Gamma_{0}(-\mu(e)-<\lambda,\check{\alpha}>-s/2)^{-\epsilon(e)}\Gamma_{0}(\mu(e)+<\lambda,\check{\alpha}>+s)^{\epsilon(e)}.$$
Il r\'esulte des descriptions ci-dessus que
$$c_{\check{\alpha}}(\pi_{\lambda}:\Lambda)=\prod_{e\in E_{\check{\alpha}}(\pi)}c_{e}(<\Lambda,\check{\alpha}>).$$
Il r\'esulte de (1) et (2) que
$$ \rho(\pi_{\lambda}:\Lambda,\tilde{P}')=\prod_{\check{\alpha}\in \check{\Sigma}_{\star ,P'\vert P}(A_{\tilde{M}})}c_{\check{\alpha}}(\pi_{\lambda};\Lambda).$$
On peut alors utiliser le lemme 7.1 de [A7] qui calcule $\rho_{\tilde{M}}^{\tilde{G}}(\pi_{\lambda})$ sous la forme:
$$\rho_{\tilde{M}}^{\tilde{G}}(\pi_{\lambda})=\sum_{\underline{\check{\alpha}}=\{\check{\alpha}_{i};i=1,...,n\}\in {\cal J}_{\tilde{M}}^{\tilde{G}}}m(\underline{\check{\alpha}})\prod_{i=1,...,n}\sum_{e\in E_{\check{\alpha}_{i}}(\pi)}c'_{e}(0).$$
Le terme $c'_{e}(0)$ est la d\'eriv\'ee de $c_{e}$ \'evalu\'ee en $0$. On calcule
$$c'_{e}(0)=\frac{\epsilon(e)}{2}\left(\frac{\Gamma'_{0}(\mu(e)+<\lambda,\check{\alpha}_{i}>)}{\Gamma_{0}(\mu(e)+<\lambda,\check{\alpha}_{i}>)}+\frac{\Gamma'_{0}(-\mu(e)-<\lambda,\check{\alpha}_{i}>)}{\Gamma_{0}(-\mu(e)-<\lambda,\check{\alpha}_{i}>)}\right).$$
D\'efinissons la fonction $\Gamma_{1}$ d'une variable complexe $s$ par
$$\Gamma_{1}(s)=\frac{1}{4}\left(\frac{\Gamma'(s/2)}{\Gamma(s/2)}+\frac{\Gamma'(-s/2)}{\Gamma(-s/2)}-\frac{\Gamma'((1+s)/2)}{\Gamma((1+s)/2)}-\frac{\Gamma'((1-s)/2)}{\Gamma((1-s)/2)}\right).$$
On obtient
$$c'_{e}(0)=\epsilon(e)\Gamma_{1}(\mu(e)+<\lambda,\check{\alpha}_{i}>),$$
d'o\`u
$$(4) \qquad \rho_{\tilde{M}}^{\tilde{G}}(\pi_{\lambda})=\sum_{\underline{\check{\alpha}}=\{\check{\alpha}_{i};i=1,...,n\}\in {\cal J}_{\tilde{M}}^{\tilde{G}}}m(\underline{\check{\alpha}})\prod_{i=1,...,n}\sum_{e\in E_{\check{\alpha}_{i}}(\pi)} \epsilon(e)\Gamma_{1}(\mu(e)+<\lambda,\check{\alpha}_{i}>).$$
On a suppos\'e que les $\check{\alpha}_{i}$ \'etaient positifs relativement \`a $P$. En vertu des propri\'et\'es des ensembles $E_{\check{\alpha}_{i}}(\pi)$ et de la parit\'e de la fonction $\Gamma_{1}$, on voit que cette formule ne change pas si l'on remplace $\check{\alpha}_{i}$ par $-\check{\alpha}_{i}$. Elle reste donc correcte en consid\'erant ${\cal J}_{\tilde{M}}^{\tilde{G}}$ comme un ensemble de classes d'\'equivalence comme on l'a d\'efini plus haut. 

 Notons $U_{\tilde{M}}^{\tilde{G}}$ l'espace de fonctions sur ${\cal A}_{\tilde{M},{\mathbb C}}^*$ engendr\'e par les fonctions de la forme
 $$\lambda\mapsto \prod_{i=1,...,n}\Gamma_{1}(c_{i}+<\lambda,\check{\alpha}_{i}>)$$
 o\`u $\underline{\check{\alpha}}=\{\check{\alpha}_{i};i=1,...,n\}$ parcourt ${\cal J}_{\tilde{M}}^{\tilde{G}}$, ou plus exactement un ensemble de repr\'esentants de cet ensemble de classes, et o\`u, pour tout $i=1,...,n$, $c_{i}$ est un nombre imaginaire. La formule ci-dessus montre que la fonction $\lambda\mapsto  \rho_{\tilde{M}}^{\tilde{G}}(\pi_{\lambda})$ appartient \`a $U_{\tilde{M}}^{\tilde{G}}$. On d\'emontre classiquement la formule
 $$\Gamma_{1}(s)=\sum_{k\geq 1}(-1)^k\frac{k}{s^2-k^2}.$$
 Pour $u\in U_{\tilde{M}}^{\tilde{G}}$, on en d\'eduit  ais\'ement les propri\'et\'es suivantes:
 
 (5)  la fonction $u$ n'a pas de p\^ole sur 
$i{\cal A}_{\tilde{M}}^*$; pour tout  sous-ensemble compact $\Omega\subset {\cal A}_{\tilde{M}}^*$ tel que $u$ n'ait pas de p\^ole sur $ \Omega+i{\cal A}_{\tilde{M}}^*$, la fonction $u$ et ses d\'eriv\'ees sont born\'ees sur cet ensemble. 
 \bigskip
 
 \subsection{L'application $\theta_{\tilde{M}}^{rat,\tilde{G}}$}
 Nous allons d\'efinir une application lin\'eaire
 $$\theta_{\tilde{M}}^{rat,\tilde{G}}:C_{c}^{\infty}(\tilde{G}({\mathbb R}))\otimes Mes(G({\mathbb R}))\to I_{ac}(\tilde{M}({\mathbb R}),\omega)\otimes Mes(M({\mathbb R})).$$
 On doit admettre par r\'ecurrence certaines de ses propri\'et\'es. A savoir qu'elle   se prolonge \`a l'espace $C_{ac}^{\infty}(\tilde{G}({\mathbb R}))\otimes Mes(G({\mathbb R}))$, qu'elle est continue et se quotiente en une application lin\'eaire d\'efinie sur $I_{ac}(\tilde{G}({\mathbb R}),\omega)\otimes Mes(G({\mathbb R}))$. On pose alors la d\'efinition
 $$\theta_{\tilde{M}}^{rat,\tilde{G}}({\bf f})=\phi_{\tilde{M}}^{\tilde{G}}({\bf f})-\sum_{\tilde{L}\in {\cal L}(\tilde{M}),\tilde{L}\not=\tilde{G}}
\theta_{\tilde{M}}^{rat,\tilde{L}}(\phi_{\tilde{L}}^{rat,\tilde{G}}({\bf f})).$$
La preuve des propri\'etes \'evoqu\'ees ci-dessus  est formelle \`a partir des propri\'et\'es de l'application $\phi_{\tilde{M}}^{\tilde{G}}$, que l'on a prouv\'ees en [V] 1.2, et des propri\'et\'es analogues de l'application $\phi_{\tilde{M}}^{rat,\tilde{G}}$.  

\bigskip

\subsection{Un lemme auxiliaire}
Pour simplifier, on fixe des mesures de Haar sur tous les groupes intervenant.
\ass{Lemme}{Soient $\tilde{\pi}_{1},...,\tilde{\pi}_{n}$ un ensemble fini d'\'el\'ements de $D_{temp}(\tilde{M}({\mathbb R}),\omega)$. Il existe un sous-ensemble ${\cal H}$ de ${\cal A}_{\tilde{M},{\mathbb C}}^*$ qui est une r\'eunion finie  de sous-espaces affines propres invariants par translations par ${\cal A}_{\tilde{G},{\mathbb C}}^*$ de sorte que, pour tout $\lambda\in i{\cal A}_{\tilde{M}}^*$, $\lambda\not\in {\cal H}$, et pour tout $\phi\in I(\tilde{M}({\mathbb R}),\omega)$, il existe $f\in I(\tilde{G}({\mathbb R}),\omega)$ de sorte que l'on ait l'\'egalit\'e
$$I^{\tilde{M}}(\tilde{\pi}_{i,\lambda},f_{\tilde{M},\omega})=I^{\tilde{M}}(\tilde{\pi}_{i,\lambda},\phi)$$
pour tout $i=1,...,n$.}

Preuve.   Soit  $f\in I(\tilde{G}({\mathbb R}),\omega)$. Pour tout $\tilde{R}\in {\cal L}(\tilde{M}_{0})$ et tout $e\in D_{ell,0}(\tilde{R}({\mathbb R}),\omega)$, on d\'efinit la fonction $\varphi_{\tilde{R},e,f}$ sur ${\cal A}_{\tilde{R},{\mathbb C}}^*$ par $\varphi_{\tilde{R},e,f}(\mu)=I^{\tilde{R}}(e_{\mu},f_{\tilde{R},\omega})$. On rappelle que le th\'eor\`eme de Paley-Wiener affirme que l'application qui \`a $f$ associe la collection de fonctions $(\varphi_{\tilde{R},e,f})_{\tilde{R}\in {\cal L}(\tilde{M}_{0}),e\in D_{ell,0}(\tilde{R}({\mathbb R}),\omega)}$ est une bijection de $I(\tilde{G}({\mathbb R}),\omega)$ sur un espace $PW^{\infty}(\tilde{G},\omega)$ d\'ecrit en [IV] 1.4.  C'est l'espace des familles $(\varphi_{\tilde{R},e})_{\tilde{R}\in {\cal L}(\tilde{M}_{0}),e\in D_{ell,0}(\tilde{R}({\mathbb R}),\omega)}$ de fonctions qui v\'erifient certaines conditions d'analycit\'e, de croissance et  d'invariance relativement \`a l'action de $W(\tilde{M}_{0})$. Pr\'ecisons cette derni\`ere condition. Pour $w\in W(\tilde{M}_{0})$ et $\tilde{R}\in {\cal L}(\tilde{M}_{0})$, $w$ induit un isomorphisme de $D_{ell,0}(\tilde{R}({\mathbb R}),\omega)$ sur $D_{ell,0}(w(\tilde{R})({\mathbb R}),\omega)$, que l'on note encore $w$.  Consid\'erons la propri\'et\'e

$(1)_{\tilde{R},w}$ $\varphi_{w(\tilde{R}),w(e)}(w(\mu))=\varphi(\tilde{R},e,\mu)$
pour tout  $e\in D_{ell,0}(\tilde{R}({\mathbb R}),\omega)$ et tout $\mu\in {\cal A}_{\tilde{R},{\mathbb C}}^*$. 

La condition est que $(1)_{\tilde{R},w}$ doit \^etre v\'erifi\'ee pour tout $\tilde{R}\in {\cal L}(\tilde{M}_{0})$ et tout $w\in W(\tilde{M}_{0})$. L'espace $I(\tilde{M}({\mathbb R}),\omega)$ se d\'ecrit aussi comme un espace de familles $(\varphi_{\tilde{R},e})_{\tilde{R}\in {\cal L}(\tilde{M}_{0}),e\in D_{ell,0}(\tilde{R}({\mathbb R}),\omega)}$. Elles doivent v\'erifier les m\^emes conditions d'analycit\'e et de croissance que pr\'ec\'edemment. Ce qui change est la condition d'invariance. Cette condition est maintenant que $(1)_{\tilde{R},w}$ doit \^etre v\'erifi\'ee pour tout $\tilde{R}\in {\cal L}(\tilde{M}_{0})$ avec $\tilde{R}\subset \tilde{M}$ et pour tout $w\in W^M(\tilde{M}_{0})$. On doit de plus avoir

(2) $\varphi_{\tilde{R},e}=0$ si $\tilde{R}\not\subset \tilde{M}$. 

On ne perd rien \`a remplacer les $\tilde{\pi}_{i}$ par d'autres \'el\'ements de   $D_{temp}(\tilde{M}({\mathbb R}),\omega)$ qui engendrent le m\^eme sous-espace.   On peut donc supposer que, pour tout $i=1,...,n$, il existe

-  $\tilde{R}_{i}\in {\cal L}^{\tilde{M}}(\tilde{M}_{0})$;

- $e_{i}\in D_{ell,0}(\tilde{R}_{i}({\mathbb R}),\omega)$ et $\mu_{i}\in i{\cal A}_{\tilde{R}}^*$ (le dernier $i$ \'etant $\sqrt{-1}$);

\noindent de sorte que $\tilde{\pi}_{i}$ soit l'induite de $\tilde{R}_{i}$ \`a $\tilde{M}$ de $e_{i,\mu_{i}}$. L'\'egalit\'e de l'\'enonc\'e s'\'ecrit alors
$$(3) \qquad \varphi_{\tilde{R}_{i},e_{i},f}(\mu_{i}+\lambda)=\varphi_{\tilde{R}_{i},e_{i},\phi}(\mu_{i}+\lambda).$$

Consid\'erons l'ensemble des triplets $(i,j,w)$, o\`u $i,j\in \{1,...,n\}$ et $w\in W(\tilde{M}_{0})-W^{M}(\tilde{M}_{0})$, tels que  $w(\tilde{R}_{i})=\tilde{R}_{j}$. Puisque $w\not\in W^M(\tilde{M}_{0})$, l'ensemble des $\lambda\in {\cal A}_{\tilde{M},{\mathbb C}}^*$ v\'erifiant  l'\'equation $w(\mu_{i}+\lambda)=\mu_{j}+\lambda$ est soit vide, soit un  sous-espace affine propre invariant par translations par ${\cal A}_{\tilde{G},{\mathbb C}}^*$. On note ${\cal H}$ la r\'eunion de ces ensembles. Soit $\lambda\in i{\cal A}_{\tilde{M}}^*$ tel que $\lambda\not\in {\cal H}$. Soit $\tilde{R}\in {\cal L}(\tilde{M}_{0})$.   Consid\'erons les deux ensembles suivants

$E'(\tilde{R})$ est l'ensemble des $w(\mu_{i}+\lambda)$ pour $i=1,...,n$ et $w\in W^{M}(\tilde{M}_{0})$ tels que $w(\tilde{R}_{i})=\tilde{R}$;

 $E''(\tilde{R})$ est l'ensemble des $w(\mu_{i}+\lambda)$ pour $i=1,...,n$ et $w\in W(\tilde{M}_{0})-W^M(\tilde{M}_{0})$ tels que $w(\tilde{R}_{i})=\tilde{R}$.   

Puisque $\lambda\not\in {\cal H}$, ces ensembles sont disjoints. On peut donc fixer un polyn\^ome $p_{\tilde{R}}$ sur ${\cal A}_{\tilde{R},{\mathbb C}}^*$ de sorte que $p_{\tilde{R}}(\mu)=1$ pour tout $\mu\in E'(\tilde{R})$ et $p_{\tilde{R}}(\mu)=0$ pour tout $\mu\in E''(\tilde{R})$. Soit $\phi\in I(\tilde{M}({\mathbb R}),\omega)$, que l'on d\'ecrit par la collection $(\varphi_{\tilde{R},e,\phi})_{\tilde{R}\in {\cal L}(\tilde{M}_{0}),e\in D_{ell,0}(\tilde{R}({\mathbb R}),\omega)}$ comme ci-dessus. Pour $\tilde{R}\in {\cal L}(\tilde{M}_{0})$, $e\in D_{ell,0}(\tilde{R}({\mathbb R}),\omega)$ et $\mu\in {\cal A}_{\tilde{R},{\mathbb C}}^*$, posons
$$(4) \qquad \varphi_{\tilde{R},e}(\mu)=\vert W^M(\tilde{M}_{0})\vert ^{-1}\sum_{w\in W(\tilde{M}_{0})}\varphi_{w(\tilde{R}),w(e),\phi}(w(\mu))p_{w(\tilde{R})}(w(\mu)).$$
La famille $(\varphi_{\tilde{R},e})_{\tilde{R}\in {\cal L}(\tilde{M}_{0}),e\in D_{ell,0}(\tilde{R}({\mathbb R}),\omega)}$ v\'erifie les conditions d'analycit\'e et de croissance requises pour appartenir \`a l'espace $PW^{\infty}(\tilde{G},\omega)$. Elle v\'erifie aussi par d\'efinition la condition d'invariance $(1)_{\tilde{R},w}$ pour tout $\tilde{R}\in {\cal L}(\tilde{M}_{0})$ et tout $w\in W(\tilde{M}_{0})$. Elle appartient donc \`a cet espace de Paley-Wiener. Donc il existe $f\in I(\tilde{G}({\mathbb R}),\omega)$ telle que $\varphi_{\tilde{R},e,f}=\varphi_{\tilde{R},e}$ pour tous $\tilde{R}$, $e$. Il reste \`a prouver que la condition (3) est v\'erifi\'ee. Fixons $i\in \{1,...,n\}$. Soit $w\in W(\tilde{M}_{0})$. Si $w\not\in W^M(\tilde{M}_{0})$, l'\'el\'ement $w(\mu_{i}+\lambda)$ appartient \`a $E''(w(\tilde{R}_{i}))$. Donc $p_{w(\tilde{R})}(w(\mu_{i}+\lambda))=0$. Si $w\in W^M(\tilde{M}_{0})$, l'\'el\'ement $w(\mu_{i}+\lambda)$ appartient \`a $E'(w(\tilde{R}))$ donc $p_{w(\tilde{R})}(w(\mu_{i}+\lambda))=1$. La d\'efinition (4) donne alors 
$$\varphi_{\tilde{R}_{i},e_{i}}(\mu_{i}+\lambda)=\vert W^M(\tilde{M}_{0})\vert ^{-1}\sum_{w\in W^M(\tilde{M}_{0})}\varphi_{w(\tilde{R}_{i}),w(e_{i}),\phi}(w(\mu_{i}+\lambda)).$$
Mais tous les termes de cette somme sont \'egaux \`a $\varphi_{\tilde{R}_{i},e_{i},\phi}(\mu_{i}+\lambda)$ d'apr\`es la condition d'invariance v\'erifi\'ee par la famille $(\varphi_{\tilde{R},e,\phi})_{\tilde{R}\in {\cal L}(\tilde{M}_{0}),e\in D_{ell,0}(\tilde{R}({\mathbb R}),\omega)}$. L'\'egalit\'e ci-dessus devient (3), ce qui ach\`eve la d\'emonstration. $\square$

 On peut renforcer le lemme de la fa\c{c}on suivante.
\ass{Lemme bis}{Soient $\tilde{\pi}_{1},...,\tilde{\pi}_{n}$ un ensemble fini d'\'el\'ements de $D_{temp}(\tilde{M}({\mathbb R}),\omega)$. Il existe un sous-ensemble ${\cal H}$ de ${\cal A}_{\tilde{M},{\mathbb C}}^*$, qui est une r\'eunion finie  de sous-espaces affines propres invariants par translations par ${\cal A}_{\tilde{G},{\mathbb C}}^*$,  et il existe un ensemble fini $\Omega$ de $K$-types de sorte que, pour tout $\lambda\in i{\cal A}_{\tilde{M}}^*$, $\lambda\not\in {\cal H}$, et pour tout $\phi\in I(\tilde{M}({\mathbb R}),\omega)$, il existe $f\in I(\tilde{G}({\mathbb R}),\omega,\Omega)$ de sorte que l'on ait l'\'egalit\'e
$$I^{\tilde{M}}(\tilde{\pi}_{i,\lambda},f_{\tilde{M},\omega})=I^{\tilde{M}}(\tilde{\pi}_{i,\lambda},\phi)$$
pour tout $i=1,...,n$.}

Preuve. Notons $pw:I(\tilde{G}({\mathbb R}),\omega) \to PW^{\infty}(\tilde{G},\omega)$ l'isomorphisme du th\'eor\`eme de Renard, cf. [IV] 1.4. Avec la terminologie introduite en [IV] 3.5, il existe un ensemble fini $\Omega^{pw}$ de types spectraux de sorte que, pour tous $i$ et $\lambda$, $I^{\tilde{M}}(\tilde{\pi}_{i,\lambda},f_{\tilde{M},\omega})$ ne d\'epende que de la projection de $pw(f)$ dans le sous-espace $PW(\tilde{G},\omega,\Omega^{pw})$. D'apr\`es le th\'eor\`eme de Delorme et Mezo, on peut fixer un ensemble fini $\Omega$ de $K$-types de sorte que $PW(\tilde{G},\omega,\Omega^{pw})\subset pw(I(\tilde{G}({\mathbb R}),\omega,\Omega))$. Il suffit alors de remplacer la fonction $f$ fournie par le lemme pr\'ec\'edent par une fonction $f'\in I(\tilde{G}({\mathbb R}),\omega,\Omega)$ telle que la projection de $pw(f)$ dans $PW(\tilde{G},\omega,\Omega^{pw})$ soit \'egale \`a $pw(f')$. $\square$

Supposons $(G,\tilde{G},{\bf a})$ quasi-d\'eploy\'e et \`a torsion int\'erieure. On a d\'efini un sous-espace "instable" $D_{ell,0}^{inst}(\tilde{M})$, cf. [IV] 2.2.

\ass{Lemme ter}{Soient $\tilde{\pi}_{1},...,\tilde{\pi}_{n}$ un ensemble fini d'\'el\'ements de $D_{ell,0}^{inst}(\tilde{M}({\mathbb R}))$. Il existe un sous-ensemble ${\cal H}$ de ${\cal A}_{\tilde{M},{\mathbb C}}^*$, qui est une r\'eunion finie de sous-espaces affines propres invariants par translations par ${\cal A}_{\tilde{G},{\mathbb C}}^*$, et il existe un ensemble fini $\Omega$ de $K$-types de sorte que, pour tout $\lambda\in i{\cal A}_{\tilde{M}}^*$, $\lambda\not\in {\cal H}$, et pour tout $\phi\in I(\tilde{M}({\mathbb R}))$, il existe $f\in I(\tilde{G}({\mathbb R}),\Omega)$ de sorte que 

- l'image de $f$ dans $SI(\tilde{G}({\mathbb R}))$ soit nulle;

- on ait l'\'egalit\'e
$$I^{\tilde{M}}(\tilde{\pi}_{i,\lambda},f_{\tilde{M},\omega})=I^{\tilde{M}}(\tilde{\pi}_{i,\lambda},\phi)$$
pour tout $i=1,...,n$.}

Preuve. On a d\'ecompos\'e l'espace de Paley-Wiener $PW^{\infty}(\tilde{G})$ en somme de deux sous-espaces $PW^{\infty,st}(\tilde{G})$ et $PW^{\infty,inst}(\tilde{G})$, cf. [IV] 2.3. L'hypoth\`ese d'instabilit\'e des $\tilde{\pi}_{i}$ entra\^{\i}ne que la conclusion ne d\'epend que de la projection de $pw(f)$ dans $PW^{\infty,inst}(\tilde{G})$. On peut remplacer la fonction $f$ fournie par le lemme bis par une fonction $f'$ telle que $pw(f')$ a une projection nulle dans $PW^{\infty,st}(\tilde{G})$ et la m\^eme projection que $pw(f)$ dans $PW^{\infty,inst}(\tilde{G})$. La premi\`ere condition implique que l'image de $f'$ dans $SI(\tilde{G}({\mathbb R}))$ est nulle, d'apr\`es le th\'eor\`eme [IV] 2.3.  $\square$

\bigskip

\subsection{Propri\'et\'es de l'application $\theta_{\tilde{M}}^{rat,\tilde{G}}$}
Les propri\'et\'es d\'emontr\'ees en [VIII] 1.6 de l'application ${^c\theta}_{\tilde{M}}^{\tilde{G}}$ sur un corps de base non-archim\'edien valent aussi pour notre application $\theta_{\tilde{M}}^{rat,\tilde{G}}$. Les d\'emonstrations en sont les m\^emes. 

On a d\'efini l'espace $D_{temp}(\tilde{G}({\mathbb R}),\omega)$ engendr\'e par les caract\`eres de repr\'esentations temp\'er\'ees et son sous-espace $D_{temp,0}(\tilde{G}({\mathbb R}),\omega)$ engendr\'e par  les caract\`eres de repr\'esentations temp\'er\'ees dont le caract\`ere central $\omega_{\pi}$ est trivial sur $\mathfrak{A}_{\tilde{G}}$. Ainsi, l'application
$$\begin{array}{ccc}D_{temp,0}(\tilde{G}({\mathbb R}),\omega)\times i{\cal A}_{\tilde{G}}^*&\to&D_{temp}(\tilde{G}({\mathbb R}),\omega)\\ (\tilde{\pi},\lambda)&\mapsto&\tilde{\pi}_{\lambda}\\ \end{array}$$
est bijective.    
On a d\'efini en 5.4 l'espace $U_{\tilde{M}}^{\tilde{G}}$. Consid\'erons un \'el\'ement $v\in U_{\tilde{M}}^{\tilde{G}}\otimes D_{temp,0}(\tilde{M}({\mathbb R}),\omega)\otimes Mes(M({\mathbb R}))^*$, un \'el\'ement $\lambda\in {\cal A}_{\tilde{M},{\mathbb C}}^*$ et une fonction ${\bf f}\in I(\tilde{M}({\mathbb R}),\omega)\otimes Mes(M({\mathbb R}))$. On peut \'ecrire
$v=\sum_{j=1,...,k}u_{j}\tilde{\pi}_{j}$ o\`u les $u_{j}$ sont des \'el\'ements de $U_{\tilde{M}}^{\tilde{G}}$ et les $\tilde{\pi}_{j}$ appartiennent \`a $ D_{temp,0}(\tilde{M}({\mathbb R}),\omega)\otimes Mes(M({\mathbb R}))^*$. On peut alors d\'efinir
$$\sum_{j=1,...,k}u_{j}(\lambda)I^{\tilde{M}}(\tilde{\pi}_{\lambda},{\bf f}).$$
Cela ne d\'epend pas de la d\'ecomposition choisie de $v$. On note $I^{\tilde{M}}(v(\lambda),{\bf f})$ cette expression. Remarquons que, puisque les $u_{j}(\lambda)$ et ses d\'eriv\'ees sont born\'ees pour $\lambda\in i{\cal A}_{\tilde{M}}^*$ et que $I^{\tilde{M}}(\tilde{\pi}_{\lambda},{\bf f})$ est  de Schwartz sur cet ensemble, la fonction $I^{\tilde{M}}(v(\lambda),{\bf f})$ est elle aussi de Schwartz sur $i{\cal A}_{\tilde{M}}^*$.

\ass{Lemme}{  Il existe une unique application lin\'eaire
$$\rho_{\tilde{M}}^{\tilde{G}}:D_{temp,0}(\tilde{M}({\mathbb R}),\omega)\otimes Mes(M({\mathbb R}))^*\to U_{\tilde{M}}^{\tilde{G}}\otimes D_{temp,0}(\tilde{M}({\mathbb R}),\omega)\otimes Mes(M({\mathbb R}))^*$$
v\'erifiant la condition suivante. Soient ${\bf f}\in I(\tilde{G}({\mathbb R}),\omega)\otimes Mes(G({\mathbb R}))$, $\tilde{\pi}\in D_{temp,0}(\tilde{M}({\mathbb R}),\omega)\otimes Mes(M({\mathbb R}))^*$ et $X\in {\cal A}_{\tilde{M}}$. Alors on a l'\'egalit\'e
$$ I^{\tilde{M}}(\tilde{\pi},X,\theta_{\tilde{M}}^{rat,\tilde{G}}({\bf f}))=\int_{i{\cal A}_{\tilde{M}}^*}I^{\tilde{M}}(\rho_{\tilde{M}}^{\tilde{G}}(\tilde{\pi};\lambda),{\bf f}_{\tilde{M},\omega})e^{-<\lambda,X>}\, d\lambda.$$ }

{\bf Remarque.} Il r\'esulte de cette formule et de 5.4(5)  que la transform\'ee de Fourier $\lambda\mapsto I^{\tilde{M}}(\tilde{\pi},\lambda,\theta_{\tilde{M}}^{rat,\tilde{G}}({\bf f}))$ de la fonction $X\mapsto  I^{\tilde{M}}(\tilde{\pi},X,\theta_{\tilde{M}}^{rat,\tilde{G}}({\bf f}))$ s'\'etend en une fonction m\'eromorphe sur tout ${\cal A}_{\tilde{M},{\mathbb C}}^*$.  Elle est \`a d\'ecroissance rapide dans les bandes verticales. Ses p\^oles sont de la forme d\'ecrite en 5.2(3).  Par contre, le nombre de ces  p\^oles (ou plus exactement ces hyperplans polaires) n'est pas toujours fini. 

\bigskip

Preuve. Soit $\tilde{\pi} $ une $\omega$-repr\'esentation temp\'er\'ee  de $\tilde{M}({\mathbb R})$. Supposons comme en 5.1 que $\pi$ soit une sous-repr\'esentation d'une induite $Ind_{S}^M(\sigma_{\lambda_{0}})$, o\`u $\sigma$ est de la s\'erie discr\`ete de $R({\mathbb R})$ et $\omega_{\sigma}=1$ sur $\mathfrak{A}_{R}$. On suppose de plus que le caract\`ere central de $\pi$ est trivial sur ${\cal A}_{\tilde{M}}$, c'est-\`a-dire que $\lambda_{0}$ annule cet espace. On  d\'efinit $\rho_{\tilde{M}}^{\tilde{G}}(\tilde{\pi})$ comme le produit $u\otimes \tilde{\pi}$, o\`u $u$ est la fonction $u(\lambda)=\rho_{\tilde{M}}^{\tilde{G}}(\pi_{\lambda})$.  Pour simplifier, on fixe des mesures sur les groupes intervenant. Soit $f\in C_{c}^{\infty}(\tilde{G}({\mathbb R}))$.  L'\'egalit\'e  de l'\'enonc\'e s'\'ecrit plus explicitement
$$(1) \qquad  I^{\tilde{M}}(\tilde{\pi},X,\theta_{\tilde{M}}^{rat,\tilde{G}}( f))=\int_{i{\cal A}_{\tilde{M}}^*}\rho_{\tilde{M}}^{\tilde{G}}(\pi_{\lambda})I^{\tilde{M}}(\tilde{\pi}_{\lambda},f_{\tilde{M},\omega})e^{-<\lambda,X>}\,d\lambda.$$
On va v\'erifier cette \'egalit\'e. 
  On a par d\'efinition
$$(2) \qquad  I^{\tilde{M}}(\tilde{\pi},X,\theta_{\tilde{M}}^{rat,\tilde{G}}( f))=I^{\tilde{M}}(\tilde{\pi},X,\phi_{\tilde{M}}^{\tilde{G}}(f))-\sum_{\tilde{L}\in {\cal L}(\tilde{M}),\tilde{L}\not=\tilde{G}}I^{\tilde{M}}(\tilde{\pi},X,\theta_{\tilde{M}}^{rat,\tilde{L}}(\phi_{\tilde{L}}^{rat,\tilde{G}}(f))).$$
Le premier terme est par d\'efinition $J_{\tilde{M}}^{\tilde{G}}(\tilde{\pi},X,f)$, c'est-\`a-dire
$$(3) \qquad I^{\tilde{M}}(\tilde{\pi},X,\phi_{\tilde{M}}^{\tilde{G}}(f))=\int_{i{\cal A}_{\tilde{M}}^*}J_{\tilde{M}}^{\tilde{G}}(\tilde{\pi}_{\lambda},f)e^{-<\lambda,X>}\,d\lambda.$$
. Fixons $\tilde{L}\not=\tilde{G}$. On ne peut pas appliquer par r\'ecurrence la formule de l'\'enonc\'e pour calculer le terme index\'e par $\tilde{L}$ car la fonction $\phi_{\tilde{L}}^{rat,\tilde{G}}(f)$ n'est pas \`a support compact en g\'en\'eral. Mais, $X$ \'etant fix\'e, on choisit une fonction $b\in C_{c}^{\infty}({\cal A}_{\tilde{L}})$ qui vaut $1$ au voisinage de la projection $X_{\tilde{L}}$ de $X$ sur ${\cal A}_{\tilde{L}}$. Posons $\varphi_{1}=\theta_{\tilde{M}}^{rat,\tilde{L}}(\phi_{\tilde{L}}^{rat,\tilde{G}}(f) )(b\circ H_{\tilde{L}})$.  On a alors
$$ I^{\tilde{M}}(\tilde{\pi},X,\theta_{\tilde{M}}^{rat,\tilde{L}}(\phi_{\tilde{L}}^{rat,\tilde{G}}(f)))=I^{\tilde{M}}(\tilde{\pi},X,\varphi_{1}).$$
D'apr\`es [VIII] 1.6(3) qui est valable pour notre application $\theta_{\tilde{M}}^{rat,\tilde{L}}$, on a $\varphi_{1}=\theta_{\tilde{M}}^{rat,\tilde{L}}(\varphi_{2})$, o\`u $\varphi_{2}=\phi_{\tilde{L}}^{rat,\tilde{G}}(f)(b\circ H_{\tilde{L}})$. Maintenant, $\varphi_{2}$ est \`a support compact et on peut appliquer par r\'ecurrence la formule de l'\'enonc\'e. On obtient
 $$ I^{\tilde{M}}(\tilde{\pi},X,\theta_{\tilde{M}}^{rat,\tilde{L}}(\phi_{\tilde{L}}^{rat,\tilde{G}}(f)))=\int_{i{\cal A}_{\tilde{M}}^*}\rho_{\tilde{M}}^{\tilde{L}}(\pi_{\lambda})I^{\tilde{M}}(\tilde{\pi}_{\lambda},\varphi_{2,\tilde{M},\omega})e^{-<\lambda,X>}\,d\lambda.$$
 Il est clair que $\rho_{\tilde{M}}^{\tilde{L}}(\pi_{\lambda})$ ne d\'epend que de la projection de $\lambda$ sur $i{\cal A}_{\tilde{M}}^{\tilde{L},*}$. On peut r\'ecrire l'\'egalit\'e ci-dessus
  $$(4) \qquad  I^{\tilde{M}}(\tilde{\pi},X,\theta_{\tilde{M}}^{rat,\tilde{L}}(\phi_{\tilde{L}}^{rat,\tilde{G}}(f)))=\int_{i{\cal A}_{\tilde{M}}^{\tilde{L},*}}\rho_{\tilde{M}}^{\tilde{L}}(\pi_{\mu})B(\mu)e^{-<\mu,X>}\,d\mu,$$
  o\`u
  $$B(\mu)=\int_{i{\cal A}_{\tilde{L}}^*}I^{\tilde{M}}(\tilde{\pi}_{\lambda+\mu},\varphi_{2,\tilde{M},\omega})e^{-<\lambda,X>}\,d\lambda.$$
Fixons $\tilde{Q}\in {\cal P}^{\tilde{L}}(\tilde{M})$.  On a aussi
   $$B(\mu)=\int_{i{\cal A}_{\tilde{L}}^*}I^{\tilde{L}}(Ind_{\tilde{Q}}^{\tilde{L}}(\tilde{\pi}_{\mu+\lambda}),\varphi_{2})e^{-<\lambda,X>}\,d\lambda=I^{\tilde{L}}(\tilde{\pi}_{\mu},X_{\tilde{L}},\varphi_{2}).$$
   En se rappelant la d\'efinition $\varphi_{2}=\phi_{\tilde{L}}^{rat,\tilde{G}}(f)(b\circ H_{\tilde{L}})$, on voit que le terme ci-dessus est \'egal \`a $I^{\tilde{L}}(Ind_{\tilde{Q}}^{\tilde{L}}(\tilde{\pi}_{\mu}),X_{\tilde{L}},\phi_{\tilde{L}}^{rat,\tilde{G}}(f))$. Ou encore, par d\'efinition,
   $$B(\mu)=\int_{i{\cal A}_{\tilde{L}}^*}J_{\tilde{L}}^{rat,\tilde{G}}(Ind_{\tilde{Q}}^{\tilde{L}}(\tilde{\pi}_{\mu+\lambda}),f)e^{-<\lambda,X>}\,d\lambda.$$
   On peut alors reconstituer la formule (4) sous la forme
   $$(5) \qquad  I^{\tilde{M}}(\tilde{\pi},X,\theta_{\tilde{M}}^{rat,\tilde{L}}(\phi_{\tilde{L}}^{rat,\tilde{G}}(f)))=\int_{i{\cal A}_{\tilde{M}}^*}\rho_{\tilde{M}}^{\tilde{L}}(\pi_{\lambda})J_{\tilde{L}}^{rat,\tilde{G}}(Ind_{\tilde{Q}}^{\tilde{L}}(\tilde{\pi}_{\lambda}),f)e^{-<\lambda,X>}\,d\lambda,$$
   puisque la fonction que l'on int\`egre est \`a d\'ecroissance rapide. 
   
   Les formules (2), (3) et (5) conduisent \`a l'\'egalit\'e
 $$ I^{\tilde{M}}(\tilde{\pi},X,\theta_{\tilde{M}}^{rat,\tilde{G}}( f))=\int_{i{\cal A}_{\tilde{M}}^*}C(\lambda)e^{-<\lambda,X>}\, d\lambda,$$
 o\`u
 $$C(\lambda)=J_{\tilde{M}}^{\tilde{G}}(\tilde{\pi}_{\lambda},f)-\sum_{\tilde{L}\in {\cal L}(\tilde{M}),\tilde{L}\not=\tilde{G}}  \rho_{\tilde{M}}^{\tilde{L}}(\pi_{\lambda})J_{\tilde{L}}^{rat,\tilde{G}}(Ind_{\tilde{Q}}^{\tilde{L}}(\tilde{\pi}_{\lambda}),f).$$
 D'apr\`es 5.4(3), on a $C(\lambda)=\rho_{\tilde{M}}^{\tilde{G}}(\pi_{\lambda})I^{\tilde{G}}(Ind_{\tilde{Q}}^{\tilde{G}}(\tilde{\pi}_{\lambda}),f)$ o\`u $\tilde{Q}\in {\cal P}(\tilde{M})$, ou encore $C(\lambda)=\rho_{\tilde{M}}^{\tilde{G}}(\pi_{\lambda})I^{\tilde{M}}(\tilde{\pi}_{\lambda},f_{\tilde{M},\omega})$. La formule pr\'ec\'edente devient (1). Cela prouve cette relation. 
 
 On vient de prouver l'existence de $\rho_{\tilde{M}}^{\tilde{G}}(\tilde{\pi})$ pour un ensemble de $\tilde{\pi}$ qui engendre lin\'eairement l'espace $D_{temp,0}(\tilde{M}({\mathbb R}),\omega)$. On peut donc d\'efinir par lin\'earit\'e une application $\rho_{\tilde{M}}^{\tilde{G}}$ qui v\'erifie l'\'egalit\'e  de l'\'enonc\'e.    On doit prouver son unicit\'e. Celle-ci r\'esulte de l'assertion suivante
 
 (6) soit $v\in U_{\tilde{M}}^{\tilde{G}}\otimes D_{temp,0}(\tilde{M}({\mathbb R}),\omega)$; supposons que, pour tout $f\in I(\tilde{G}({\mathbb R}),\omega)$ et tout $X\in {\cal A}_{\tilde{M}}$, on ait l'\'egalit\'e
 $$\int_{i{\cal A}_{\tilde{M}}^*}I^{\tilde{M}}(v(\lambda),f_{\tilde{M},\omega})e^{-<\lambda,X>}\,d\lambda=0;$$
 alors $v=0$. 
 
 Par inversion de Fourier, l'hypoth\`ese \'equivaut \`a l'\'egalit\'e
   $$I^{\tilde{M}}(v(\lambda),f_{\tilde{M},\omega})=0$$
   pour tout $f\in I(\tilde{G}({\mathbb R}),\omega)$ et tout $\lambda\in i{\cal A}_{\tilde{M}}^*$. Ecrivons $v=\sum_{i=1,...,k}u_{i}\otimes \tilde{\pi}_{i}$, o\`u $u_{i}\in U_{\tilde{M}}^{\tilde{G}}$ et $\tilde{\pi}_{i}\in D_{temp,0}(\tilde{M}({\mathbb R}),\omega)$ pour tout $i$. On peut supposer  que $u_{i}\not=0$ pour tout $i$ et que la famille $(\tilde{\pi}_{i})_{i=1,...,k}$ est lin\'eairement ind\'ependante. Si $v\not=0$, on a $k\geq1$ et on peut fixer $\phi_{1}\in I(\tilde{M}({\mathbb R}),\omega)$ tel que $I^{\tilde{M}}(\tilde{\pi}_{1},\phi_{1})=1$ tandis que $I^{\tilde{M}}(\tilde{\pi}_{i},\phi_{1})=0$ pour tout 
$i\geq2$.  Introduisons l'ensemble ${\cal H}$ du lemme 5.6 pour notre famille   $(\tilde{\pi}_{i})_{i=1,...,k}$. Fixons $\lambda\in i{\cal A}_{\tilde{M}}^*$ tel que $\lambda\not\in {\cal H}$ et $u_{1}(\lambda)\not=0$.  D\'efinissons la fonction $\phi$ sur $\tilde{M}({\mathbb R})$ par $\phi(\gamma)=e^{-<\lambda,H_{\tilde{M}}(\gamma)>}\phi_{1}(\gamma)$. On a
$$I^{\tilde{M}}(\tilde{\pi}_{i,\lambda},\phi)=I^{\tilde{M}}(\tilde{\pi}_{i},\phi_{1})$$
pour tout $i$. Associons \`a $\phi$ une fonction $f\in I(\tilde{G}({\mathbb R}),\omega)$ satisfaisant la conclusion du lemme 5.6. On a alors $I^{\tilde{M}}(\tilde{\pi}_{1,\lambda},f_{\tilde{M},\omega})=1$ tandis que $I^{\tilde{M}}(\tilde{\pi}_{i,\lambda},f_{\tilde{M},\omega})=0$ pour tout $i\geq2$. Alors $I^{\tilde{M}}(v(\lambda),f_{\tilde{M},\omega})=u_{1}(\lambda)\not=0$, ce qui contredit l'hypoth\`ese. Cela prouve (6) et le lemme. $\square$

 Rappelons que, pour une $\omega$-repr\'esentation $\tilde{\pi}$ de $\tilde{G}({\mathbb R})$ telle que $\pi$ soit irr\'eductible, on d\'efinit son param\`etre infinit\'esimal $\mu(\tilde{\pi})$ qui est une orbite dans $\mathfrak{h}^*$ pour l'action du groupe de Weyl $W$. Pour une telle orbite, on note $D_{temp,\mu}(\tilde{G}({\mathbb R}),\omega)$ le sous-espace de $D_{temp}(\tilde{G}({\mathbb R}),\omega)$ engendr\'e par les $\tilde{\pi}$ de param\`etre infinit\'esimal $\mu$. On a les variantes $D_{temp,0,\mu}(\tilde{G}({\mathbb R}),\omega)$, $D_{ell,\mu}(\tilde{G}({\mathbb R}),\omega)$, $D_{ell,0,\mu}(\tilde{G}({\mathbb R}),\omega)$. Ce dernier espace est de dimension finie. Plus g\'en\'eralement, pour $\lambda\in i{\cal A}_{\tilde{G}}^*$, notons $D_{ell,0,\mu}(\tilde{G}({\mathbb R}),\omega)_{\lambda}$ l'ensemble des $\tilde{\pi}_{\lambda}$ pour $\tilde{\pi}\in D_{ell,0,\mu}(\tilde{G}({\mathbb R}),\omega)$. La construction explicite de l'appllication $\rho_{\tilde{M}}^{\tilde{G}}$ effectu\'ee dans la preuve ci-dessus entra\^{\i}ne la propri\'et\'e suivante. Soient $\tilde{R}\in {\cal L}^M(\tilde{M}_{0})$, $\mu$ une orbite dans $\mathfrak{h}^*$ pour l'action de $W^R$ et $\lambda\in i{\cal A}_{\tilde{R}}^{\tilde{M},*}$. Alors
 
 (7) $\rho_{\tilde{M}}^{\tilde{G}}$ envoie $Ind_{\tilde{R}}^{\tilde{M}}(D_{ell,0,\mu}(\tilde{R}({\mathbb R}),\omega)_{\lambda})\otimes Mes(R({\mathbb R}))^*$ dans $U_{\tilde{M}}^{\tilde{G}}\otimes Ind_{\tilde{R}}^{\tilde{M}}(D_{ell,0,\mu}(\tilde{R}({\mathbb R}),\omega)_{\lambda})\otimes Mes(R({\mathbb R}))^*$.
 
 En cons\'equence

 (8)   $\theta_{\tilde{M}}^{rat,\tilde{G}}$ envoie $C_{c}^{\infty}(\tilde{G}({\mathbb R}),K)\otimes Mes(G({\mathbb R}))$ dans $I_{ac}(\tilde{M}({\mathbb R}),\omega,K)\otimes Mes(M({\mathbb R}))$.
   
 Preuve. On oublie comme souvent les espaces de mesures. L'espace de Paley-Wiener $PW^{\infty}_{ell}(\tilde{G},\omega)$ est un espace de familles de fonctions holomorphes index\'ees par une base de $D_{ell,0}(\tilde{G}({\mathbb R}),\omega)$. On peut supposer que cette base est r\'eunion de bases des sous-espaces $D_{ell,0,\mu}(\tilde{G}({\mathbb R}),\omega)$ quand $\mu$ d\'ecrit tous les param\`etres possibles. On note $PW_{ell,\mu}(\tilde{G},\omega)$ le sous-espace des familles appartenant \`a $PW^{\infty}_{ell}(\tilde{G},\omega)$ dont les composantes sont nulles pour les indices n'appartenant pas \`a  $D_{ell,0,\mu}(\tilde{G}({\mathbb R}),\omega)$. Rappelons que l'on a des homomorphismes naturels
   $$\oplus_{\tilde{L}\in {\cal L}(\tilde{M}_{0})}PW^{\infty}_{ell}(\tilde{L},\omega)\stackrel{sym}{\to}PW^{\infty}(\tilde{G},\omega)\stackrel{pw}{\leftarrow}I(\tilde{G}({\mathbb R}),\omega).$$
    Un \'el\'ement $f\in I_{ac}(\tilde{G}({\mathbb R}),\omega)$ appartient \`a $I_{ac}(\tilde{G}({\mathbb R}),\omega,K)$ si et seulement s'il existe un nombre fini de couples $(\tilde{L}_{i},\mu_{i})_{i=1,...,k}$, o\`u $\tilde{L}_{i}\in {\cal L}(\tilde{M}_{0})$ et $\mu_{i}$ est une $W^L$-orbite dans $\mathfrak{h}^*$, de sorte que, pour toute fonction $b\in C_{c}^{\infty}({\cal A}_{\tilde{G}})$, on ait la relation
    $$pw(f(b\circ H_{\tilde{G}}))\in sym(\oplus_{i=1,...,k}PW_{ell,\mu_{i}}(\tilde{L}_{i},\omega)).$$
    Soit $f\in C_{c}^{\infty}(\tilde{G}({\mathbb R}),K)$. Alors $f_{\tilde{M},\omega}$ est $K^M$-finie \`a droite et \`a gauche et on peut fixer un ensemble fini de couples $(\tilde{L}_{i},\mu_{i})_{i=1,...,k}$ comme ci-dessus, avec $\tilde{L}_{i}\subset \tilde{M}$, de sorte que 
 $$pw^{M}(f_{\tilde{M},\omega})\in sym^{M}(\oplus_{i=1,...,k}PW_{ell,\mu_{i}}(\tilde{L}_{i},\omega)$$
(on a ajout\'e des exposants $M$ pour d\'esigner les objets relatifs \`a $\tilde{M}$ plut\^ot qu'\`a $\tilde{G}$). Quitte \`a accro\^{\i}tre la famille $(\tilde{L}_{i},\mu_{i})_{i=1,...,k}$, on peut la supposer invariante par l'action de $W^M$. Soit $b\in C_{c}^{\infty}({\cal A}_{\tilde{M}})$. Soient $\tilde{L}\in {\cal L}^M(\tilde{M}_{0})$ et $\mu$ une $W^L$-orbite dans $\mathfrak{h}^*$.
Pour tout $\lambda\in i{\cal A}_{\tilde{L}}^{\tilde{M},*}$,  tout $\tilde{\pi}\in D_{ell,\mu,0}(\tilde{L}({\mathbb R}),\omega)_{\lambda}$  et tout $X\in {\cal A}_{\tilde{M}}$, on a l'\'egalit\'e
$$I^{\tilde{M}}(Ind_{\tilde{L}}^{\tilde{M}}(\tilde{\pi}),X,\theta_{\tilde{M}}^{rat,\tilde{G}}(f)(b\circ H_{\tilde{M}}))=I^{\tilde{M}}(Ind_{\tilde{L}}^{\tilde{M}}(\tilde{\pi}),X,\theta_{\tilde{M}}^{rat,\tilde{G}}(f)) b(X).$$
L'assertion (7) et la formule du lemme montrent que ceci est nul si $(\tilde{L},\mu)$ n'est pas l'un des $(\tilde{L}_{i},\mu_{i})$. Il en r\'esulte que la composante de $pw^M(\theta_{\tilde{M}}^{rat,\tilde{G}}(f)(b\circ H_{\tilde{M}}))$ dans $PW_{ell,\mu}(\tilde{L},\omega)$ est nulle si  $(\tilde{L},\mu)$ n'est pas l'un des $(\tilde{L}_{i},\mu_{i})$. Cela prouve (8). $\square$

 \bigskip
 
 \subsection{L'application $^c\phi_{\tilde{M}}^{\tilde{G}}$}
 Pour tout $\tilde{M}\in {\cal L}(\tilde{M}_{0})$ et tout $\tilde{P}\in {\cal P}(\tilde{M})$, on fixe une fonction $\omega_{\tilde{P}}:{\cal A}_{\tilde{M}}\to [0,1]$. On suppose que ces fonctions v\'erifient les hypoth\`eses de [VIII] 1.1.
 
 {\bf Remarque.} Dans cette r\'ef\'erence, le corps de base \'etait non-archim\'edien et on avait d\'efini des fonctions $\omega_{\tilde{P}}$ sur les espaces $\tilde{{\cal A}}_{\tilde{M}}$.  Ici, le corps de base est archim\'edien et les fonctions $H_{\tilde{M}}:\tilde{M}({\mathbb R})\to {\cal A}_{\tilde{M}}$ que l'on a fix\'ees permettent d'identifier $\tilde{{\cal A}}_{\tilde{M}}$ \`a ${\cal A}_{\tilde{M}}$.
 \bigskip
 
 On suppose de plus que les fonctions $\omega_{\tilde{P}}$ sont $C^{\infty}$.

 Soient $\tilde{R}\in {\cal L}^{\tilde{M}}(\tilde{M}_{0})$, $\tilde{\pi}$ un \'el\'ement de $D_{ell}(\tilde{R}({\mathbb R}),\omega)\otimes Mes(R({\mathbb R}))^*$ et ${\bf f}\in C_{c}^{\infty}(\tilde{G}({\mathbb R}),K)\otimes Mes(G({\mathbb R}))$. Pour $X\in {\cal A}_{\tilde{R}}$, la fonction
 $$\lambda\mapsto J_{\tilde{M}}^{rat,\tilde{G}}(Ind_{\tilde{R}}^{\tilde{M}}(\tilde{\pi}_{\lambda}),{\bf f})e^{-<\lambda,X>}$$
 est rationnelle sur ${\cal A}_{\tilde{R},{\mathbb C}}^*$. Elle est \`a d\'ecroissance rapide dans les bandes verticales. Cela r\'esulte des propri\'et\'es de $J_{\tilde{M}}^{rat,\tilde{G}}$ et de l'hypoth\`ese que
 ${\bf f}$ est $K$-finie \`a droite et \`a gauche.    
 Pour $\nu\in {\cal A}_{\tilde{R}}^*$ tel que cette fonction n'ait pas de p\^ole sur $\nu+i{\cal A}_{\tilde{R}}^*$, on  peut former l'int\'egrale
 $$\int_{\nu+i{\cal A}_{\tilde{R}}^*}J_{\tilde{M}}^{rat,\tilde{G}}(Ind_{\tilde{R}}^{\tilde{M}}(\tilde{\pi}_{\lambda}),{\bf f})e^{-<\lambda,X>}\,d\lambda.$$
 Soit $\tilde{S}\in {\cal P}^{\tilde{G}}(\tilde{R})$. Fixons un point $\nu_{\tilde{S}}\in {\cal A}_{\tilde{R}}^*$ tel que $<\nu_{\tilde{S}},\check{\alpha}>$ soit assez grand pour tout $\alpha\in \Sigma^{\tilde{S}}(\tilde{R})$. L'int\'egrale ci-dessus pour $\nu=\nu_{\tilde{S}}$ ne d\'epend pas du choix de $\nu_{\tilde{S}}$. On pose
 $$^cJ_{\tilde{M}}^{rat,\tilde{G}}(Ind_{\tilde{R}}^{\tilde{M}}(\tilde{\pi}),{\bf f})=\int_{{\cal A}_{\tilde{R}}}\sum_{\tilde{S}\in {\cal P}^{\tilde{G}}(\tilde{R})}\omega_{\tilde{S}}(X)\int_{\nu_{\tilde{S}}+i{\cal A}_{\tilde{R}}^*}J_{\tilde{M}}^{rat,\tilde{G}}(Ind_{\tilde{R}}^{\tilde{M}}(\tilde{\pi}_{\lambda}){\bf f})e^{-<\lambda,X>}\,d\lambda\,dX.$$
 \ass{Proposition}{(i) La fonction en $X$ que l'on int\`egre ci-dessus est $C^{\infty}$ et \`a support compact. Le terme $^cJ_{\tilde{M}}^{rat,\tilde{G}}(Ind_{\tilde{R}}^{\tilde{M}}(\tilde{\pi}),{\bf f})$ ne d\'epend que de $Ind_{\tilde{R}}^{\tilde{M}}(\tilde{\pi})$ (et pas des choix de $\tilde{R}$ et $\tilde{\pi}$).
 
 (ii) Il existe une unique application lin\'eaire 
 $$^c\phi_{\tilde{M}}^{\tilde{G}} :C_{c}^{\infty}(\tilde{G}({\mathbb R}),K)\otimes Mes(G({\mathbb R}))\to I(\tilde{M}({\mathbb R}),\omega,K^M)\otimes Mes(M({\mathbb R}))$$
 de sorte que, pour tout $\tilde{R}\in {\cal L}^M(\tilde{M}_{0})$, tout $\tilde{\pi}\in D_{ell}(\tilde{R}({\mathbb R}),\omega)\otimes Mes(R({\mathbb R}))$ et tout ${\bf f}\in C_{c}^{\infty}(\tilde{G}({\mathbb R}),K)\otimes Mes(G({\mathbb R}))$, on ait l'\'egalit\'e
 $$I^{\tilde{M}}(Ind_{\tilde{R}}^{\tilde{M}}(\tilde{\pi}),{^c\phi}_{\tilde{M}}^{\tilde{G}}({\bf f}))={^cJ}_{\tilde{M}}^{rat,\tilde{G}}(Ind_{\tilde{R}}^{\tilde{M}}(\tilde{\pi}),{\bf f}).$$}
 
 La preuve est identique \`a celle  du cas non-archim\'edien, cf. [VIII] 1.3. 
 
 {\bf Remarque.} On peut pr\'eciser que, si $\Omega$ est un ensemble fini de $K$-types, il existe un ensemble fini $\Omega^M$ de $K^M$-types de sorte que
 $^c\phi_{\tilde{M}}^{\tilde{G}} $ envoie $C_{c}^{\infty}(\tilde{G}({\mathbb R}),\Omega)\otimes Mes(G({\mathbb R}))$ dans $ I(\tilde{M}({\mathbb R}),\omega,\Omega^M)\otimes Mes(M({\mathbb R}))$.
 
 \bigskip
 
 Les propri\'et\'es prouv\'ees  dans le cas non-archim\'edien en [VIII] 1.4 valent aussi  dans notre cas o\`u le corps de base est r\'eel. En particulier, la propri\'et\'e [VIII] 1.4(6) dit que $^c\phi_{\tilde{M}}^{\tilde{G}}$ s'\'etend en une application lin\'eaire
  $$C_{ac}^{\infty}(\tilde{G}({\mathbb R}),K)\otimes Mes(G({\mathbb R}))\to I_{ac}(\tilde{M}({\mathbb R}),\omega,K^M)\otimes Mes(M({\mathbb R})).$$

 \bigskip
 
 \subsection{L'application ${^c\theta}_{\tilde{M}}^{rat,\tilde{G}}$}
 On d\'efinit une application
 $${^c\theta}_{\tilde{M}}^{rat,\tilde{G}}:C_{c}^{\infty}(\tilde{G}({\mathbb R}),K)\otimes Mes(G({\mathbb R}))\to I_{ac}(\tilde{M}({\mathbb R}),\omega,K^M)\otimes Mes(M({\mathbb R}))$$
  par la formule de r\'ecurrence
 $${^c\theta}_{\tilde{M}}^{rat,\tilde{G}}({\bf f})=\phi_{\tilde{M}}^{rat,\tilde{G}}({\bf f})-\sum_{\tilde{L}\in {\cal L}(\tilde{M}),\tilde{L}\not=\tilde{G}}{^c\theta}_{\tilde{M}}^{rat,\tilde{L}}(^c\phi_{\tilde{L}}^{\tilde{G}}({\bf f})).$$
Pour que cette d\'efinition ait un sens, il faut admettre par r\'ecurrence la propri\'et\'e suivante.

\ass{Proposition }{L'application ${^c\theta}_{\tilde{M}}^{rat,\tilde{G}}$ se quotiente en une application lin\'eaire d\'efinie sur $I(\tilde{G}({\mathbb R}),\omega,K)\otimes Mes(G({\mathbb R}))$.}

 On d\'emontre facilement que l'application ${^c\theta}_{\tilde{M}}^{rat,\tilde{G}}$ est $\omega$-\'equivariante. C'est-\`a-dire que soit $h\in C_{c}^{\infty}(G({\mathbb R}))$ une fonction $K$-finie \`a droite et \`a gauche. Notons $h_{\omega}$ la fonction $h_{\omega}(g)=\omega(g)h(g)$. Fixons des mesures de Haar. Pour $f\in C_{c}^{\infty}(\tilde{G}({\mathbb R}),\omega,K)$, on d\'efinit les fonctions $f\star h$ et $h_{\omega}\star f$ sur $\tilde{G}({\mathbb R})$ par
 $$(f\star h)(\gamma)=\int_{G({\mathbb R})}f(\gamma g^{-1})h(g)\,dg,$$
 $$(h_{\omega}\star f)(\gamma)=\int_{G({\mathbb R})}h_{\omega}(g)f(g^{-1}\gamma)\, dg.$$
 Ces fonctions sont encore $K$-finies \`a droite et \`a gauche. 
 On montre que ${^c\theta}_{\tilde{M}}^{rat,\tilde{G}}(f\star h)={^c\theta}_{\tilde{M}}^{rat,\tilde{G}}(h_{\omega}\star f)$. Mais, dans notre cas o\`u le corps de base est r\'eel, cela ne suffit pas \`a prouver la proposition. Il faudrait de plus prouver une propri\'et\'e de continuit\'e de notre application. On pr\'ef\`ere revenir \`a la m\'ethode d'Arthur. La proposition sera prouv\'ee en 5.14. 
 
 Comme dans le paragraphe pr\'ec\'edent, on peut pr\'eciser que, si $\Omega$ est un ensemble fini de $K$-types, il existe un ensemble fini $\Omega^M$ de $K^M$-types de sorte que
 ${^c\theta}_{\tilde{M}}^{rat,\tilde{G}} $ envoie $C_{c}^{\infty}(\tilde{G}({\mathbb R}),\Omega)\otimes Mes(G({\mathbb R}))$ dans $ I_{ac}(\tilde{M}({\mathbb R}),\omega,\Omega^M)\otimes Mes(M({\mathbb R}))$.

 \bigskip
 
 \subsection{Propri\'et\'es de l'application ${^c\theta}_{\tilde{M}}^{rat,\tilde{G}}$}
 Les propri\'et\'es d\'emontr\'ees en [VIII] 1.6 et 1.8 dans le cas non-archim\'edien pour l'application not\'ee alors $^c\theta_{\tilde{M}}^{\tilde{G}}$ valent aussi sur le corps de base r\'eel pour l'application  ${^c\theta}_{\tilde{M}}^{rat,\tilde{G}}$. Rappelons la principale. Soit $\tilde{\pi}\in D_{temp}(\tilde{M}({\mathbb R}),\omega)\otimes Mes(M({\mathbb R}))^*$ et soit ${\bf f}\in C_{c}^{\infty}(\tilde{G}({\mathbb R}),K)\otimes Mes(G({\mathbb R}))$.  La fonction
 $$X\mapsto I^{\tilde{M}}(\tilde{\pi},X,\phi_{\tilde{M}}^{rat,\tilde{G}}({\bf f}))$$
 sur ${\cal A}_{\tilde{M}}$ est la transform\'ee de Fourier de la fonction 
 $\lambda\mapsto J_{\tilde{M}}^{rat,\tilde{G}}(\tilde{\pi}_{\lambda},{\bf f})$ sur $i{\cal A}_{\tilde{M}}^*$. Celle-ci  se prolonge en une fonction m\'eromorphe sur ${\cal A}_{\tilde{M},{\mathbb C}}^*$ qui v\'erifie les propri\'et\'es (2), (3), (4) et (5) de 5.2.

  On voit alors par r\'ecurrence que la fonction
 $$X\mapsto I^{\tilde{M}}(\tilde{\pi},X,{^c\theta}_{\tilde{M}}^{rat,\tilde{G}}({\bf f}))$$
 est la transform\'ee de Fourier d'une fonction sur $i{\cal A}_{\tilde{M}}^*$, laquelle se prolonge en une fonction m\'eromorphe sur ${\cal A}_{\tilde{M},{\mathbb C}}^*$  qui a les m\^emes propri\'et\'es que ci-dessus.  On note $\lambda\mapsto I^{\tilde{M}}(\tilde{\pi},\lambda,{^c\theta}_{\tilde{M}}^{rat,\tilde{G}}({\bf f}))$ cette fonction. Pour $\nu\in {\cal A}_{\tilde{M}}^*$ tel que cette fonction n'ait pas de p\^ole sur $\nu+i{\cal A}_{\tilde{M}}^*$ et pour $X\in {\cal A}_{\tilde{M}}$, on pose
 $$I^{\tilde{M}}(\tilde{\pi},\nu,X,{^c\theta}_{\tilde{M}}^{rat,\tilde{G}}({\bf f}))=\int_{\nu+i{\cal A}_{\tilde{M}}^*}I^{\tilde{M}}(\tilde{\pi},\lambda,{^c\theta}_{\tilde{M}}^{rat,\tilde{G}}({\bf f}))e^{-<\lambda,X>}\,d\lambda.$$
 
 Dans l'\'enonc\'e suivant, on fixe pour tout $\tilde{S}\in {\cal P}(\tilde{M})$ un point $\nu_{\tilde{S}}$ comme en 5.8.  
 
 \ass{Proposition}{Supposons $\tilde{M}\not=\tilde{G}$ et $\tilde{\pi}$ elliptique. Soit ${\bf f}\in C_{c}^{\infty}(\tilde{G}({\mathbb R}),K)\otimes Mes(G({\mathbb R}))$. Si chaque point $\nu_{\tilde{S}}$ est assez positif relativement \`a $\tilde{S}$, on a l'\'egalit\'e
 $$\sum_{\tilde{S}\in {\cal P}(\tilde{M})}\omega_{\tilde{S}}(X)I^{\tilde{M}}(\tilde{\pi},\nu_{\tilde{S}},X,{^c\theta}_{\tilde{M}}^{rat,\tilde{G}}({\bf f}))=0$$
 pour tout $X\in {\cal A}_{\tilde{M}}$.}
 
 {\bf Remarque.} La condition "assez positif" impos\'ee aux points $\nu_{\tilde{S}}$ d\'epend de $\tilde{\pi}$ et de ${\bf f}$. Toutefois, la propri\'et\'e 5.2(5) implique que, si l'on fixe un ensemble fini de $K$-types, on peut choisir des points $\nu_{\tilde{S}}$ assez positifs, cette notion ne d\'ependant plus que  de $\Omega$, de sorte que l'\'egalit\'e de l'\'enonc\'e soit v\'erifi\'ee pour tous $\tilde{\pi}$ et ${\bf f}$, pourvu que  ${\bf f}\in C_{c}^{\infty}(\tilde{G}({\mathbb R}),\Omega)\otimes Mes(G({\mathbb R}))$.
 \bigskip
 
 \subsection{L'application ${^c\theta}_{\tilde{M}}^{\tilde{G}}$}
  On d\'efinit une application
 $${^c\theta}_{\tilde{M}}^{\tilde{G}}:C_{c}^{\infty}(\tilde{G}({\mathbb R}),K)\otimes Mes(G({\mathbb R}))\to I_{ac}(\tilde{M}({\mathbb R}),\omega,K^M)\otimes Mes(M({\mathbb R}))$$
  par la formule de r\'ecurrence
 $${^c\theta}_{\tilde{M}}^{\tilde{G}}({\bf f})=\phi_{\tilde{M}}^{\tilde{G}}({\bf f})-\sum_{\tilde{L}\in {\cal L}(\tilde{M}),\tilde{L}\not=\tilde{G}}{^c\theta}_{\tilde{M}}^{\tilde{L}}(^c\phi_{\tilde{L}}^{\tilde{G}}({\bf f})).$$
 La diff\'erence avec la d\'efinition de  5.9 est que l'on a remplac\'e le premier terme $\phi_{\tilde{M}}^{rat,\tilde{G}}({\bf f})$ par $\phi_{\tilde{M}}^{\tilde{G}}({\bf f})$. Comme en 5.9, 
pour que cette d\'efinition ait un sens, il faut admettre par r\'ecurrence la propri\'et\'e suivante, qui sera prouv\'ee en 5.14.

\ass{Proposition }{L'application ${^c\theta}_{\tilde{M}}^{\tilde{G}}$ se quotiente en une application lin\'eaire d\'efinie sur $I(\tilde{G}({\mathbb R}),\omega,K)\otimes Mes(G({\mathbb R}))$.}

 Soit $\tilde{\pi}\in D_{temp}(\tilde{M}({\mathbb R}),\omega)\otimes Mes(M({\mathbb R}))^*$ et soit ${\bf f}\in C_{c}^{\infty}(\tilde{G}({\mathbb R}),K)\otimes Mes(G({\mathbb R}))$.  La fonction
 $$X\mapsto I^{\tilde{M}}(\tilde{\pi},X,{^c\theta}_{\tilde{M}}^{\tilde{G}}({\bf f}))$$
 sur ${\cal A}_{\tilde{M}}$ est la transform\'ee de Fourier d'une fonction $\lambda\mapsto I^{\tilde{M}}(\tilde{\pi},\lambda,{^c\theta}_{\tilde{M}}^{\tilde{G}}({\bf f}))$ sur $i{\cal A}_{\tilde{M}}^*$. Celle-ci s'\'etend en une fonction m\'eromorphe sur ${\cal A}_{\tilde{M},{\mathbb C}}^*$ qui est \`a d\'ecroissance rapide dans les bandes verticales. Ses p\^oles sont de la forme d\'ecrite en 5.2(3). Par contre, le nombre de ces p\^oles (ou plut\^ot de ces hyperplans polaires) n'est pas forc\'ement fini. Ces propri\'et\'es r\'esultent par r\'ecurrence de la d\'efinition ci-dessus et des propri\'et\'es analogues des fonctions $X\mapsto I^{\tilde{M}}(\tilde{\pi},X, \phi_{\tilde{M}}^{\tilde{G}}({\bf f}))$.

\bigskip

\subsection{Relation entre les applications $\theta_{\tilde{M}}^{rat,\tilde{G}}$, ${^c\theta}_{\tilde{M}}^{rat,\tilde{G}}$ et ${^c\theta}_{\tilde{M}}^{\tilde{G}}$}
\ass{Lemme}{Soit ${\bf f}\in C_{c}^{\infty}(\tilde{G}({\mathbb R},K)\otimes Mes(G({\mathbb R}))$. On a l'\'egalit\'e
$${^c\theta}_{\tilde{M}}^{\tilde{G}}({\bf f})=\sum_{\tilde{L}\in {\cal L}(\tilde{M})}\theta_{\tilde{M}}^{rat,\tilde{L}}\circ{^c\theta}_{\tilde{L}}^{rat,\tilde{G}}({\bf f}).$$}

Preuve. Si $\tilde{M}=\tilde{G}$, les trois applications $\theta_{\tilde{G}}^{rat,\tilde{G}}$, $^c\theta_{\tilde{G}}^{rat,\tilde{G}}$ et $^c\theta_{\tilde{G}}^{\tilde{G}}$ sont l'identit\'e et la relation de l'\'enonc\'e est claire. Supposons $\tilde{M}\not=\tilde{G}$.  En utilisant le fait que $\theta_{\tilde{M}}^{rat,\tilde{M}}$, $^c\theta_{\tilde{M}}^{rat,\tilde{M}}$ et $^c\theta_{\tilde{M}}^{\tilde{M}}$ sont l'identit\'e, les d\'efinitions de nos applications peuvent se reformuler de la fa\c{c}on suivante:
$$(1)\qquad \phi_{\tilde{M}}^{\tilde{G}}({\bf f})-{^c\phi}_{\tilde{M}}^{\tilde{G}}({\bf f})=\sum_{\tilde{L}\in {\cal L}(\tilde{M}),\tilde{L}\not=\tilde{M}}{^c\theta}_{\tilde{M}}^{\tilde{L}}({^c\phi}_{\tilde{L}}^{\tilde{G}}({\bf f})),$$
 $$(2) \qquad \phi_{\tilde{M}}^{\tilde{G}}({\bf f})-\phi_{\tilde{M}}^{rat, \tilde{G}}({\bf f})=\sum_{\tilde{L}\in {\cal L}(\tilde{M}),\tilde{L}\not=\tilde{M}}\theta_{\tilde{M}}^{rat,\tilde{L}}(\phi_{\tilde{L}}^{rat,\tilde{G}}({\bf f})),$$
$$(3) \qquad \phi_{\tilde{M}}^{rat,\tilde{G}}({\bf f})-{^c\phi}_{\tilde{M}}^{\tilde{G}}({\bf f})=\sum_{\tilde{L}\in {\cal L}(\tilde{M}),\tilde{L}\not=\tilde{M}}{^c\theta}_{\tilde{M}}^{rat,\tilde{L}}({^c\phi}_{\tilde{L}}^{\tilde{G}}({\bf f})).$$ 
Dans le terme du membre de droite de (2) index\'e par $\tilde{L}$, on utilise (3) avec $\tilde{M}$ remplac\'e par $\tilde{L}$. Cela transforme (2) en
$$\phi_{\tilde{M}}^{\tilde{G}}({\bf f})-\phi_{\tilde{M}}^{rat, \tilde{G}}({\bf f})=\sum_{\tilde{L}\in {\cal L}(\tilde{M}),\tilde{L}\not=\tilde{M}}\theta_{\tilde{M}}^{rat,\tilde{L}}({^c\phi}_{\tilde{L}}^{\tilde{G}}({\bf f}))$$
$$+\sum_{\tilde{L}\in {\cal L}(\tilde{M}),\tilde{L}\not=\tilde{M}}\sum_{\tilde{L}'\in {\cal L}(\tilde{L}),\tilde{L}'\not=\tilde{L}}
 \theta_{\tilde{M}}^{rat,\tilde{L}}\circ{^c\theta}_{\tilde{L}}^{rat,\tilde{L}'}({^c\phi}_{\tilde{L}'}^{\tilde{G}}({\bf f})).$$
 Dans la derni\`ere somme, on change la notation en rempla\c{c}ant $(\tilde{L},\tilde{L}')$ par $(\tilde{L}',\tilde{L})$.   Le membre de gauche de (1) est la somme de ceux de (2) et (3). Il en est donc de m\^eme des membres de droite. On obtient
 $$\sum_{\tilde{L}\in {\cal L}(\tilde{M}),\tilde{L}\not=\tilde{M}}{^c\theta}_{\tilde{M}}^{\tilde{L}}({^c\phi}_{\tilde{L}}^{\tilde{G}}({\bf f}))=\sum_{\tilde{L}\in {\cal L}(\tilde{M}),\tilde{L}\not=\tilde{M}}\theta_{\tilde{M}}^{rat,\tilde{L}}({^c\phi}_{\tilde{L}}^{\tilde{G}}({\bf f}))$$
 $$+\sum_{\tilde{L}'\in {\cal L}(\tilde{M}),\tilde{L}'\not=\tilde{M}}\sum_{\tilde{L}\in {\cal L}(\tilde{L}'),\tilde{L}\not=\tilde{L}'}
 \theta_{\tilde{M}}^{rat,\tilde{L}'}\circ{^c\theta}_{\tilde{L}'}^{rat,\tilde{L}}({^c\phi}_{\tilde{L}}^{\tilde{G}}({\bf f}))$$
 $$+\sum_{\tilde{L}\in {\cal L}(\tilde{M}),\tilde{L}\not=\tilde{M}}{^c\theta}_{\tilde{M}}^{rat,\tilde{L}}({^c\phi}_{\tilde{L}}^{\tilde{G}}({\bf f})).$$
 On peut inclure la premi\`ere somme du membre de droite dans la deuxi\`eme somme, en supprimant la condition $\tilde{L}'\not=\tilde{L}$. On peut de m\^eme inclure la troisi\`eme somme dans la deuxi\`eme en supprimant la condition $\tilde{L}'\not=\tilde{M}$. On obtient
 $$\sum_{\tilde{L}\in {\cal L}(\tilde{M}),\tilde{L}\not=\tilde{M}}{^c\theta}_{\tilde{M}}^{\tilde{L}}({^c\phi}_{\tilde{L}}^{\tilde{G}}({\bf f}))=\sum_{\tilde{L}\in {\cal L}(\tilde{M}), \tilde{L}\not=\tilde{M}}\sum_{\tilde{L}'\in {\cal L}^{\tilde{L}}(\tilde{M})}\theta_{\tilde{M}}^{rat,\tilde{L}'}\circ{^c\theta}_{\tilde{L}'}^{rat,\tilde{L}}({^c\phi}_{\tilde{L}}^{\tilde{G}}({\bf f})).$$
 En raisonnant par r\'ecurrence, on peut supposer la relation de l'\'enonc\'e prouv\'ee si l'on remplace $\tilde{G}$ par $\tilde{L}\not=\tilde{G}$. Alors les termes des deux membres ci-dessus index\'es par un tel $\tilde{L}$ sont \'egaux. Il ne reste que l'\'egalit\'e des termes index\'es par $\tilde{G}$. Cette \'egalit\'e est celle de l'\'enonce. $\square$

 \bigskip
 
 \subsection{Une variante des int\'egrales orbitales pond\'er\'ees $\omega$-\'equivariantes}
Soit $\boldsymbol{\gamma}\in D_{orb}(\tilde{M}({\mathbb R}),\omega)\otimes Mes(M({\mathbb R}))^*$, cf. [V] 1.3 pour cette notation. On d\'efinit une application lin\'eaire
$${\bf f}\mapsto {^cI}_{\tilde{M}}^{\tilde{G}}(\boldsymbol{\gamma},{\bf f})$$
sur $C_{c}^{\infty}(\tilde{G}({\mathbb R}),K)\otimes Mes(G({\mathbb R}))$ par la formule de r\'ecurrence
$$(1) \qquad {^cI}_{\tilde{M}}^{\tilde{G}}(\boldsymbol{\gamma},{\bf f})=J_{\tilde{M}}^{\tilde{G}}(\boldsymbol{\gamma},{\bf f})-\sum_{\tilde{L}\in {\cal L}(\tilde{M}),\tilde{L}\not=\tilde{G}}{^cI}_{\tilde{M}}^{\tilde{L}}(\boldsymbol{\gamma},{^c\phi}_{\tilde{L}}^{\tilde{G}}({\bf f})).$$
Elle v\'erifie les m\^emes propri\'et\'es que dans le cas non archim\'edien. Rappelons-les.

 (2) La distribution ${\bf f}\mapsto  {^cI}_{\tilde{M}}^{\tilde{G}}(\boldsymbol{\gamma},{\bf f})$ se quotiente en une forme lin\'eaire sur $I(\tilde{G}({\mathbb R}),\omega,K)\otimes Mes(G({\mathbb R}))$.
 
 Remarquons que cette propri\'et\'e est utilis\'ee par r\'ecurrence pour poser la d\'efinition (1).  Comme pour les propositions 5.9 et 5.11, la preuve est plus d\'elicate dans le cas archim\'edien. Elle sera faite en 5.14.  Les propri\'et\'es suivantes se prouvent, elles, comme dans le cas non-archim\'edien, cf. [VIII] 1.9.
 
 (3) Pour tout ${\bf f}\in C_{c}^{\infty}(\tilde{G}({\mathbb R}),K)\otimes Mes(G({\mathbb R}))$, il existe un sous-ensemble compact $\Gamma\subset \tilde{M}({\mathbb R})$ tel que ${^cI}_{\tilde{M}}^{\tilde{G}}(\boldsymbol{\gamma},{\bf f})=0$ si le support de $\boldsymbol{\gamma}$ ne coupe pas $\Gamma^M=\{m^{-1}\gamma m; m\in M({\mathbb R}), \gamma\in \Gamma\}$.
 
 (4) Soient $\tilde{R}\in {\cal L}^M(\tilde{M}_{0})$ et $\boldsymbol{\gamma}\in  D_{orb}(\tilde{R}({\mathbb R}),\omega)\otimes Mes(R({\mathbb R}))^*$. On a l'\'egalit\'e
 $$ {^cI}_{\tilde{M}}^{\tilde{G}}(\boldsymbol{\gamma}^{\tilde{M}},{\bf f})=\sum_{\tilde{L}\in {\cal L}(\tilde{R})}d_{\tilde{R}}^{\tilde{G}}(\tilde{M},\tilde{L}) {^cI}_{\tilde{R}}^{\tilde{L}}(\boldsymbol{\gamma},{\bf f}_{\tilde{L},\omega}).$$
 
 (5) On a l'\'egalit\'e
 $$ {^cI}_{\tilde{M}}^{\tilde{G}}(\boldsymbol{\gamma},{\bf f})=\sum_{\tilde{L}\in {\cal L}(\tilde{M})} I_{\tilde{M}}^{\tilde{L}}(\boldsymbol{\gamma},{^c\theta}_{\tilde{L}}^{\tilde{G}}({\bf f})).$$

 \bigskip
 
 \subsection{Preuve des propositions 5.9, 5.11 et de l'assertion 5.13(1)}
 Soit ${\bf f}\in C_{c}^{\infty}(\tilde{G}({\mathbb R}),K)\otimes Mes(G({\mathbb R}))$. Supposons que l'image de ${\bf f}$ dans $I(\tilde{G}({\mathbb R}),\omega,K)\otimes Mes(G({\mathbb R}))$ soit nulle. On veut prouver les assertions
 
 (1) ${^c\theta}_{\tilde{M}}^{rat,\tilde{G}}({\bf f})=0$;
 
 (2) ${^c\theta}_{\tilde{M}}^{\tilde{G}}({\bf f})=0$;
 
 (3) $^cI_{\tilde{M}}^{\tilde{G}}(\boldsymbol{\gamma},{\bf f})=0$  pour tout $\boldsymbol{\gamma}\in D_{orb}(\tilde{M}({\mathbb R}),\omega)\otimes Mes(M({\mathbb R}))^*$.
 
 Ces assertions sont tautologiques si $\tilde{M}=\tilde{G}$. Supposons $\tilde{M}\not=\tilde{G}$.
  Dans la formule 5.13(5), utilisons par r\'ecurrence l'assertion (2) o\`u l'on remplace $\tilde{M}$ par $\tilde{L}$ pour $\tilde{L}\in {\cal L}(\tilde{M})$, $\tilde{L}\not=\tilde{M}$. Cette formule 5.13(5) se simplifie en
 $$(4)\qquad ^cI_{\tilde{M}}^{\tilde{G}}(\boldsymbol{\gamma},{\bf f})=I^{\tilde{M}}(\boldsymbol{\gamma},{^c\theta}_{\tilde{M}}^{\tilde{G}}({\bf f})).$$
 Dans le lemme 5.12, utilisons de la m\^eme fa\c{c}on l'assertion (1) par r\'ecurrence (remarquons que, pour les applications $\theta_{\tilde{M}}^{rat,\tilde{L}}$, il n'y a pas de probl\`eme: on sait qu'elles se factorisent par $I_{ac}(\tilde{L}({\mathbb R}),\omega,K)\otimes Mes(L({\mathbb R}))$). On obtient
 $$(5) \qquad {^c\theta}_{\tilde{M}}^{\tilde{G}}({\bf f})={^c\theta}_{\tilde{M}}^{rat,\tilde{G}}({\bf f}).$$
 En cons\'equence de (4) et (5), on a
 $$^cI_{\tilde{M}}^{\tilde{G}}(\boldsymbol{\gamma},{\bf f})=I^{\tilde{M}}(\boldsymbol{\gamma},{^c\theta}_{\tilde{M}}^{rat,\tilde{G}}({\bf f})).$$
 La fonction ${^c\theta}_{\tilde{M}}^{rat,\tilde{G}}({\bf f})$ est un \'el\'ement de $I_{ac}(\tilde{M}({\mathbb R}),\omega,K)\otimes Mes(M({\mathbb R}))$. En apliquant la propri\'et\'e 5.13(3), l'\'egalit\'e ci-dessus entra\^{\i}ne que $I^{\tilde{M}}(\boldsymbol{\gamma},{^c\theta}_{\tilde{M}}^{rat,\tilde{G}}({\bf f}))=0$    si le support de $\boldsymbol{\gamma}$ ne coupe pas $\Omega^M$, o\`u $\Omega$ est un certain sous-ensemble compact de $\tilde{M}({\mathbb R})$. Il en r\'esulte que ${^c\theta}_{\tilde{M}}^{rat,\tilde{G}}({\bf f})$ est "\`a support compact", c'est-\`a-dire est un \'el\'ement de $I(\tilde{M}({\mathbb R}),\omega,K)\otimes Mes(M({\mathbb R}))$. Pour $\tilde{\pi}\in D_{temp}(\tilde{M}({\mathbb R}),\omega)\otimes Mes(M({\mathbb R}))^*$, la fonction $\lambda\mapsto I^{\tilde{M}}(\tilde{\pi},\lambda,{^c\theta}_{\tilde{M}}^{rat,\tilde{G}}({\bf f}))$ est donc holomorphe. Ses coefficients de Fourier $ I^{\tilde{M}}(\tilde{\pi},\nu,X,{^c\theta}_{\tilde{M}}^{rat,\tilde{G}}({\bf f}))$ ne d\'ependent pas du point $\nu$. Dans la proposition 5.10, on peut remplacer tous les $\nu_{\tilde{S}}$ par $0$. Puisque 
 $$\sum_{\tilde{S}\in {\cal P}(\tilde{M})}\omega_{\tilde{S}}(X)=1,$$
 on obtient que, pour $\tilde{\pi}$ elliptique, les coefficients de Fourier $ I^{\tilde{M}}(\tilde{\pi},0,X,{^c\theta}_{\tilde{M}}^{rat,\tilde{G}}({\bf f}))$ sont tous nuls. Cela entra\^{\i}ne que $I^{\tilde{M}}(\tilde{\pi},{^c\theta}_{\tilde{M}}^{rat,\tilde{G}}({\bf f}))=0$.  D'autre part, la fonction ${^c\theta}_{\tilde{M}}^{rat,\tilde{G}}({\bf f})$ est cuspidale. En effet, on a la formule de descente
 $$({^c\theta}_{\tilde{M}}^{rat,\tilde{G}}({\bf f}))_{\tilde{R},\omega}=\sum_{\tilde{L}\in {\cal L}(\tilde{R})}d_{\tilde{R}}^{\tilde{G}}(\tilde{M},\tilde{L}){^c\theta}_{\tilde{R}}^{rat,\tilde{L}}({\bf f}_{\tilde{Q},\omega})$$
 pour tout espace de Levi $\tilde{R}\in {\cal L}^M(\tilde{M}_{0})$. En utilisant l'assertion (1) par r\'ecurrence, on en d\'eduit que  $({^c\theta}_{\tilde{M}}^{rat,\tilde{G}}({\bf f}))_{\tilde{R},\omega}=0$ si $\tilde{R}$ est propre. Un \'el\'ement de $I(\tilde{M}({\mathbb R}),\omega,K)\otimes Mes(M({\mathbb R}))$ qui est cuspidal et annul\'e par toutes les repr\'esentations elliptiques est nul.  L'assertion (1) en r\'esulte.  L'\'egalit\'e (5) entra\^{\i}ne alors (2) et l'\'egalit\'e (4) entra\^{\i}ne (3). $\square$
 
 \bigskip
 
 \subsection{Une propri\'et\'e de l'espace $U_{\tilde{M}}^{\tilde{G}}$}
 On note ${\cal F}$ l'espace des fonctions holomorphes sur ${\cal A}_{\tilde{M},{\mathbb C}}^*$. Soit ${\cal H}$ un sous-ensemble de ${\cal A}_{\tilde{M},{\mathbb C}}^*$ qui est r\'eunion finie de sous-espaces affines propres  invariants par translations par ${\cal A}_{\tilde{G},{\mathbb C}}^*$. Soient $k\geq1$ un entier et $E$ un sous-ensemble de ${\cal F}^k$.  Un \'el\'ement de $E$ est une famille $e=(e_{j})_{j=1,...,k}$ o\`u $e_{j}\in {\cal F}$ pour tout $j$. On suppose
 
 (1) pour tout $\lambda\in {\cal A}_{\tilde{M},{\mathbb C}}^*-{\cal H}$, il existe $e=(e_{j})_{j=1,...,k}\in E$ tel que $e_{1}(\lambda)=1$ tandis que $e_{j}(\lambda)=0$ pour $j=2,...,k$.

 \ass{Proposition}{On suppose $\tilde{M}\not=\tilde{G}$. Soit $(u_{j})_{j=1,...,k}$ une famille d'\'el\'ements de $U_{\tilde{M}}^{\tilde{G}}$. On suppose que,
 pour tout $e=(e_{j})_{j=1,...,k}\in E$, la fonction $\lambda\mapsto \sum_{j=1,...,k}u_{j}(\lambda)e_{j}(\lambda)$ est holomorphe hors de ${\cal H}$. Alors $u_{1}=0$.}
 
 Preuve. On peut agrandir ${\cal H}$ et supposer que c'est une r\'eunion finie d'hyperplans affines.
  D\'efinissons une fonction $p:{\mathbb C}\to \{-1/2,0,1/2\}$ par
 $p(k)=\frac{(-1)^k sgn(k)}{2}$ pour $k\in {\mathbb Z}-\{0\}$ et $p(s)=0$ pour $s\not\in {\mathbb Z}$ ou $s=0$. La fonction $\Gamma_{1}(s)$ de 5.4 a pour p\^oles les points $s$ tels que $p(s)\not=0$, et son r\'esidu en $s$ est $p(s)$. Remarquons la propri\'et\'e suivante
 
 (2) soient $s\in {\mathbb C}$ et $k\in 2{\mathbb Z}$; si $Re(s)$ et $Re(s+k)$ sont tous deux  non nuls et s'ils sont de m\^eme signe, on a $p(s)=p(s+k)$; si $Re(s)$ et $Re(s+k)$ sont tous deux non nuls mais sont de signes oppos\'es, on a $p(s)=-p(s+k)$.

 Notons $n$ la dimension de ${\cal A}_{\tilde{M}}^{\tilde{G}}$. On a d\'efini l'ensemble ${\cal J}_{\tilde{M}}^{\tilde{G}}$ en 5.4. On fixe un \'el\'ement $\tilde{P}\in{\cal P}(\tilde{M})$, qui permet de supposer que les \'el\'ements de ${\cal J}_{\tilde{M}}^{\tilde{G}}$ sont des familles $\underline{\check{\alpha}}=(\check{\alpha}_{i})_{i=1,...,n}$ telles que, pour tout $i$,  $\check{\alpha}_{i}$ est positif pour $P$. Consid\'erons une telle famile, ainsi qu'une famille $\underline{y}=(y_{i})_{i=1,...,n}$ de nombres complexes  imaginaires. Soient $\lambda,\mu\in {\cal A}_{\tilde{M},{\mathbb C}}^*$, supposons $\mu$ en position g\'en\'erale (pr\'ecis\'ement $<\mu,\check{\alpha}_{i}>\not=0$ pour tout $i$). On voit que
 
 (3) la fonction 
 $$s\mapsto \prod_{i=1,...,n}\Gamma_{1}(y_{i}+<\lambda+s\mu,\check{\alpha}_{i}>)$$
 sur ${\mathbb C}$ a  un p\^ole d'ordre au plus  $n$ en $s=0$; le coefficient de $s^{-n}$ dans son d\'eveloppement en $0$ est \'egal \`a
 $$\prod_{i=1,..,n}p(y_{i}+<\lambda,\check{\alpha}_{i}>)<\mu,\check{\alpha}_{i}>^{-1}.$$

 Pour tout $j=1,...,k$, on peut fixer
 
 - un ensemble fini d'indices $B_{j}$;
 
 - pour $b\in B_{j}$, un nombre complexe $c_{b}$, une famille $\underline{y}_{b}=(y_{b,i})_{i=1,...,n}$ de nombres complexes imaginaires et une famille $\underline{\check{\alpha}}_{b}=(\check{\alpha}_{b,i})_{i=1,...,n}\in {\cal J}_{\tilde{M}}^{\tilde{G}}$;
 
 \noindent de sorte que, pour tout $\lambda\in {\cal A}_{\tilde{M},{\mathbb C}}^*$, on ait l'\'egalit\'e
 $$u_{j}(\lambda)=\sum_{b\in B_{j}}c_{b}\prod_{i=1,...,n}\Gamma_{1}(y_{b,i}+<\lambda,\check{\alpha}_{b,i}>).$$
 Cette formule ne  change pas si, pour un indice $j$ et pour \'el\'ement $b\in B_{j}$, on remplace les familles  $\underline{y}_{b}=(y_{b,i})_{i=1,...,n}$ et $\underline{\check{\alpha}}_{b}=(\check{\alpha}_{b,i})_{i=1,...,n}$ par les familles $\underline{y}'_{b}=(y_{b,\sigma(i)})_{i=1,...,n}$ et $\underline{\check{\alpha}}'_{b}=(\check{\alpha}_{b,\sigma(i)})_{i=1,...,n}$, o\`u $\sigma$ est une permutation de $\{1,...,n\}$. On s'autorisera \`a effectuer de telles permutations d'indices. 
 On peut supposer que $c_{b}\not=0$ pour tout $b\in B_{j}$ et que, si $b,b'$ sont deux \'el\'ements distincts de $B_{j}$,  les familles associ\'ees \`a $b$ et $b'$ sont distinctes ( c'est-\`a-dire qu'il n'y a pas de permutation $\sigma$ de $\{1,...,n\}$ telle que $y_{b,i}=y_{b',\sigma(i)}$ et $\check{\alpha}_{b,i}=\check{\alpha}_{b',\sigma(i)}$ pour tout $i$).
 
 Soient $\lambda\in {\cal A}_{\tilde{M},{\mathbb C}}^*-{\cal H}$, $\mu\in {\cal A}_{\tilde{M},{\mathbb C}}^*$ en position g\'en\'erale et $e=(e_{j})_{j=1,...,k}\in E$. D'apr\`es (3) le coefficient de $s^{-n}$ dans le d\'eveloppement en $0$ de la fonction
 $$s\mapsto \sum_{j=1,...,k}u_{j}(\lambda+s\mu)e_{j}(\lambda+s\mu)$$
 est
 $$\sum_{j=1,...,k}e_{j}(\lambda)\sum_{b\in B_{j}}c_{b}\prod_{i=1,...,n}p(y_{b,i}+<\lambda,\check{\alpha}_{b,i}>)<\mu,\check{\alpha}_{b,i}>^{-1}.$$
 L'hypoth\`ese d'holomorphie  et l'hypoth\`ese $n\geq1$ impliquent que cette expression est nulle. Par ailleurs, d'apr\`es la condition (1), on peut choisir $e$ de sorte que $e_{1}(\lambda)=1$ et $e_{j}(\lambda)=0$ pour $j=2,...,k$. En posant simplement $B=B_{1}$ (on peut \'evidemment supposer $k\geq1$), on obtient
 
 (4) pour tout $\lambda\in {\cal A}_{\tilde{M},{\mathbb C}}^*-{\cal H}$ et tout $\mu\in {\cal A}_{\tilde{M},{\mathbb C}}^*$ en position g\'en\'erale, on a
 $$\sum_{b\in B}c_{b}\prod_{i=1,...,n}p(y_{b,i}+<\lambda,\check{\alpha}_{b,i}>)<\mu,\check{\alpha}_{b,i}>^{-1}=0.$$
 
 Remarquons que l'on peut supprimer la condition que $\mu$ est en position g\'en\'erale: l'expression  ci-dessus est m\'eromorphe en $\mu$, la condition est que cette fonction m\'eromorphe est nulle. Remarquons aussi que nos fonctions ne d\'ependent que des projections de $\lambda$ et $\mu$ modulo ${\cal A}_{\tilde{G},{\mathbb C}}^*$. On ne perd rien \`a supposer que cet espace est nul.
 
 On veut prouver que $u_{1}=0$, c'est-\`a-dire que $B$ est vide. Raisonnons par l'absurde et supposons $B$ non vide. Par construction de l'ensemble ${\cal J}_{\tilde{M}}^{\tilde{G}}$, tous les \'el\'ements $\check{\alpha}_{b,i}$ appartiennent \`a un ${\mathbb Z}$-module de type fini contenu dans ${\cal A}_{\tilde{M}}^{\tilde{G}}$. Montrons que l'on peut choisir un \'el\'ement $b_{0}\in B$ tel que, quitte \`a permuter les indices, on ait
 
 (5) soient $m\in \{1,...,n\}$ et $b\in B$; supposons qu'il existe $c\in {\mathbb Q}$ tel que $\check{\alpha}_{b,i}=\check{\alpha}_{b_{0},i}$ pour $i=1,...,m-1$ et $\check{\alpha}_{b,m}=c\check{\alpha}_{b_{0},m}$; alors $0<c\leq 1$;
 
 En effet, soit $m\in \{1,...,n\}$,  supposons choisi $b_{0}$ tel que cette propri\'et\'e soit v\'erifi\'ee pour tout $m'<m$. Notons $B'$ le sous-ensemble des $b\in B$ tels que, quitte \`a permuter les indices, on ait $\check{\alpha}_{b,i}=\check{\alpha}_{b_{0},i}$ pour $i=1,...,m-1$. Consid\'erons l'ensemble des \'el\'ements $\check{\alpha}_{b,i}$ pour $b\in B'$ et $i\geq m$. On peut choisir un  \'el\'ement  $\check{\alpha}_{b_{1},i_{1}}$  de cet ensemble qui soit maximal, au sens que, si un autre \'el\'ement $\check{\alpha}_{b,i}$ v\'erifie $\check{\alpha}_{b,i}=c\check{\alpha}_{b_{1},i_{1}}$, avec $c\in {\mathbb Q}$, on ait $\vert c\vert \leq 1$. On a automatiquement $0< c\leq1$ par l'hypoth\`ese de positivit\'e impos\'ee \`a nos coracines. Quitte \`a permuter les indices, on peut supposer $i_{1}=m$. On remplace $b_{0}$ par $b_{1}$ et on voit que la propri\'et\'e (5) est alors v\'erifi\'ee pour tout $m'\leq m$. Par r\'ecurrence, on obtient (5).

 Fixons un \'el\'ement $b_{0}\in B$ v\'erifiant (5). Pour $i=1,...,n$, on pose simplement   $H_{i}=\check{\alpha}_{b_{0},i}$. La famille $(H_{i})_{i=1,...,n}$ est une base de ${\cal A}_{\tilde{M}}$ (puisqu'on a suppos\'e ${\cal A}_{\tilde{G}}=\{0\}$). On introduit les coordonn\'ees duales sur ${\cal A}_{\tilde{M},{\mathbb C} }^*$, c'est-\`a-dire que l'on note tout \'el\'ement $\lambda$ de cet espace sous la forme $\lambda=(\lambda_{i})_{i=1,...,n}$, o\`u $\lambda_{i}=<\lambda,H_{i}>$.  Pour tout $b\in B$ et tout $i\in \{1,...,n\}$, on peut \'ecrire
 $$(6) \qquad \check{\alpha}_{b,i}=\sum_{l=1,...,n}x_{b,i,l}H_{l},$$
 avec des coefficients $x_{b,i}\in {\mathbb Q}$. On fixe un d\'enominateur commun $D\in {\mathbb N}$, $D\not=0$. On note $C_{1}$ un majorant des valeurs absolues $\vert x_{b,i,l}\vert $ pour tous $b$, $i$ et $l$. Par hypoth\`ese, ${\cal H}$ est r\'eunion finie d'hyperplans, on note $h$ le nombre de ces hyperplans.
 Soit $C$ un r\'eel tel que
 $ C>5DnhC_{1}$. Pour un entier $m\in \{0,...,n\}$, notons $V[m]$ l'ensemble des $\lambda\in {\cal A}_{\tilde{M}}^*-{\cal H}$ tels que
 $C^{2i-1}\leq\vert Re(\lambda_{i})\vert \leq C^{2i}$ pour tout $i$ et $Re(\lambda_{i})>0$ pour $i\leq m$. Notons $B[m]$ le sous-ensemble des $b\in B$ tels qu'il existe des indices $i_{1},...,i_{m}\in \{1,...,n\}$ tous distincts de sorte que 
 
 - $\check{\alpha}_{b,i_{1}}$ appartient au sous-espace de ${\cal A}_{\tilde{M}}$ engendr\'e par $H_{1}$;
 
 - $\check{\alpha}_{b,i_{2}}$ appartient au sous-espace engendr\'e par $H_{1}$ et $H_{2}$;
 
 -...
 
 - $\check{\alpha}_{b,i_{m}}$ appartient au sous-espace engendr\'e par $H_{1}$,...,$H_{m}$.
 
 Ces conditions d\'eterminent les indices et, quitte \`a permuter l'ensemble d'indices, on peut supposer $i_{1}=1$,... $i_{m}=m$. On va prouver par r\'ecurrence sur $m$ que 
 
 (7) pour tout $\lambda\in V[m]$ et tout $\mu\in {\cal A}_{\tilde{M},{\mathbb C}}^*$ , on a l'\'egalit\'e 
 $$\sum_{b\in B[m]}c_{b}\prod_{i=1,...,n}p(y_{b,i}+<\lambda,\check{\alpha}_{b,i}>)<\mu,\check{\alpha}_{b,i}>^{-1}=0.$$
 
 Pour $m=0$, $B[m]=B$ et l'\'egalit\'e est un cas particulier de (4). Supposons $m>0$ et la relation prouv\'ee pour $m-1$. Soit $\lambda\in V[m]$. Les conditions impos\'ees \`a $C$ impliquent qu'il y a au moins $h+1$ points $z\in {\mathbb C}$ tels que
 $$\lambda_{m}-z\in 2D{\mathbb Z},\,\, -C^{2m}\leq Re(z)\leq -C^{2m-1}.$$
 Pour un tel point, notons $\lambda[z]$ l'\'el\'ement de ${\cal A}_{\tilde{M},{\mathbb C}}^*$ tel que $\lambda[z]_{i}=\lambda_{i}$ pour $i\not=m$ et $\lambda[z]_{m}=z$. Montrons que l'on  peut choisir $z$ tel que $\lambda[z]\not\in {\cal H}$. En effet, supposons que tous les $\lambda[z]$ appartiennent \`a ${\cal H}$. Puisqu'on a au moins $h+1$ points $z$, ,il y en a au moins deux, disons $z$ et $z'$, tels que $\lambda[z]$ et $\lambda[z']$ appartiennent \`a un m\^eme hyperplan contenu dans ${\cal H}$. Alors $\lambda=\frac{z'-\lambda_{m}}{z'-z}\lambda[z]+\frac{z-\lambda_{m}}{z-z'}\lambda[z']$ appartient aussi \`a cet hyperplan, contrairement \`a l'hypoth\`ese.  Cela prouve l'assertion. On choisit un  $z$ tel que $\lambda[z]\not\in {\cal H}$ et on pose simplement $\lambda'=\lambda[z]$. Les deux \'el\'ements $\lambda$ et $\lambda'$ appartiennent \`a $V[m-1]$. On applique l'hypoth\`ese de r\'ecurrence \`a chacun d'eux et on obtient par diff\'erence
 $$(8) \qquad \sum_{b\in B[m-1]}c_{b}X_{b}(\lambda,\lambda')\prod_{i=1,...,n}<\mu,\check{\alpha}_{b,i}>^{-1}=0,$$
 o\`u
 $$X_{b}(\lambda,\lambda')=\prod_{i=1,...,n}p(y_{b,i}+<\lambda,\check{\alpha}_{b,i}>)-\prod_{i=1,...,n}p(y_{b,i}+<\lambda',\check{\alpha}_{b,i}>).$$
 Montrons que
 
 (9 ) pour $b\in B[m-1]$, on a les \'egalit\'es 
 $$X_{b}(\lambda,\lambda')=\left\lbrace\begin{array}{cc}2\prod_{i=1,...,n}p(y_{b,i}+<\lambda,\check{\alpha}_{b,i}>),&\text{ si }b\in B[m],\\ 0,& \text{ si }b\not\in B[m].\\ \end{array}\right.$$
 
 Pour $i\leq m-1$, $<\lambda,\check{\alpha}_{b,i}>$ est combinaison lin\'eaire des $\lambda_{l}$ pour $l\leq m-1$ et remplacer $\lambda$ par $\lambda'$ ne change rien. Soit $i\geq m$. Notons $l_{i,max}$ le plus grand entier $l$ tel que $x_{b,i,l}\not=0$, avec la notation de (6). On a $l_{i,max}\geq m$. On a
 $$<\lambda,\check{\alpha}_{b,i}>=\sum_{l=1,...,l_{i,max}}x_{b,i,l}\lambda_{l}.$$
 La condition $\lambda\in V[m-1]$ et la condition impos\'ee \`a $C$ impliquent que
 $$\vert Re(x_{b,i,l_{i,max}}\lambda_{l_{i,max}})\vert >\sum_{l=1,...,l_{i,max}-1}\vert x_{b,i,l}\lambda_{l} \vert .$$
 Donc $Re(<\lambda,\check{\alpha}_{b,i}>)$ est du m\^eme signe que $ Re(x_{b,i,l_{i,max}}\lambda_{l_{i,max}})$. La m\^eme propri\'et\'e  vaut avec $\lambda$ remplac\'e par $\lambda'$. En cons\'equence, $Re(<\lambda,\check{\alpha}_{b,i}>)$ et $Re(<\lambda',\check{\alpha}_{b,i}>)$ sont de m\^eme signe si $l_{i,max}>m$, tandis qu'ils sont de signes oppos\'es si $l_{i,max}=m$. Par ailleurs, la condition $\lambda'_{m}-\lambda_{m}\in 2D{\mathbb Z}$ implique que 
 $$< \lambda,\check{\alpha}_{b,i}>-<\lambda',\check{\alpha}_{b,i}>\in 2{\mathbb Z}.$$
 D'apr\`es (2), on a donc 
 $$p(y_{b,i}+<\lambda,\check{\alpha}_{b,i}>)=p(y_{b,i}+<\lambda',\check{\alpha}_{b,i}>)$$
si $l_{i,max}>m$ tandis que
 $$p(y_{b,i}+<\lambda,\check{\alpha}_{b,i}>)=-p(y_{b,i}+<\lambda',\check{\alpha}_{b,i}>)$$
 si $l_{i,max}=m$. Les deux produits figurant dans la d\'efinition de $X_{b}(\lambda,\lambda')$ sont donc soit \'egaux, soit oppos\'es. Ils sont \'egaux exactement si le nombre d'indices $i\geq m$ tels que $l_{i,max}=m$ est impair. Or il r\'esulte des d\'efinitions que ce nombre d'indices est au plus $1$ et qu'il est \'egal \`a $1$ exactement si $b\in B[m]$. Cela prouve (9).
 
 Les relations (8) et (9) entra\^{\i}nent la conclusion de (7). Cela prouve cette assertion (7). 
 
 Pour $b\in B[n]$, l'ordre des indices des familles associ\'ees \`a $b$ est impos\'e ($\check{\alpha}_{b,i}$ appartient au sous-espace engendr\'e par $H_{1},...,H_{i}$). Pour un entier $q\in \{1,...,n+1\}$, notons $B[n,q]$ le sous-ensemble des $b\in B[n]$ tels que $x_{b,i,l}=0$ pour $i\geq q$ et $l\not=i$. Autrement dit, en posant simplement $x_{b,i}=x_{b,i,i}$, on a $\check{\alpha}_{b,i}=x_{b,i}H_{i}$ pour $i\geq q$. Montrons que
 
 (10)   pour tout $\lambda\in V[n]$ et tout $\mu\in {\cal A}_{\tilde{M},{\mathbb C}}^*$ , on a l'\'egalit\'e 
 $$\sum_{b\in B[n,q]}c_{b}\left(\prod_{i=1,...,n}p(y_{b,i}+<\lambda,\check{\alpha}_{b,i}>)\right)\left(\prod_{i=q,...,n}x_{b,i}^{-1}\right)\left(\prod_{i=1,...,q-1} <\mu,\check{\alpha}_{b,i}>^{-1}\right)=0.$$
 
 On raisonne par r\'ecurrence descendante sur $q$. Pour $q=n+1$, c'est l'assertion (7) pour $m=n$. Supposons $q\leq n$ et l'assertion prouv\'ee pour $q+1$. Remarquons que    l'assertion \`a prouver  ne d\'epend de $\mu$ que via les coordonn\'ees $\mu_{1}$,...,$\mu_{q-1}$. Puisque l'expression est m\'eromorphe en $\mu$, on peut se limiter au cas o\`u ces coordonn\'ees sont en position g\'en\'erale. Consid\'erons ces coordonn\'ees comme fix\'ees et consid\'erons la relation ci-dessus pour $q+1$. Elle d\'epend de la coordonn\'ee $\mu_{q}$ que l'on consid\`ere comme variable. On regarde son r\'esidu en $\mu_{q}=0$. Les seuls termes pouvant cr\'eer un r\'esidu sont les $<\mu,\check{\alpha}_{b,q}>^{-1}$ pour $b\in B[n,q+1]$. Plus explicitement, ce terme est
 $$(\sum_{i=1,...,q}\mu_{i}x_{b,q,i})^{-1}.$$
 S'il y a un $i<q$ tel que  $x_{b,q,i}\not=0$, l'hypoth\`ese que $\mu_{1}$,...,$\mu_{q-1}$ sont en position g\'en\'erale implique que le terme ci-dessus n'a pas de p\^ole en $\mu_{q}=0$. Si $x_{b,q,i}=0$ pour tout $i<q$, il y a un p\^ole et le r\'esidu est $x_{b,q,q}^{-1}$. Or la condition que $x_{b,q,i}=0$ pour tout $i<q$ \'equivaut \`a $b\in B[n,q]$. On voit alors que l'expression \`a prouver est \'egale au r\'esidu en $\mu_{q}=0$ de la m\^eme expression pour $q+1$. Cela prouve (10). 
 
 Posons simplement $B^*=B[n,1]$. C'est l'ensemble des $b\in B$ tels que, pour tout $i$, $\check{\alpha}_{b,i}=x_{b,i}H_{i}$. La relation (10) pour $q=1$ se r\'ecrit
 $$(11) \qquad \sum_{b\in B^*}c_{b}\prod_{i=1,...,n}x_{b,i}^{-1}p(y_{b,i}+\lambda_{i}x_{b,i})=0$$
 pour tout $\lambda\in V[n]$. Montrons que
 
 (12) il existe $\lambda\in V[n]$ tel que $y_{b_{0},i}+\lambda_{i}\in 1+2D{\mathbb Z}$ pour tout $i$. 
 
 Rappelons que $y_{b_{0},i}$ est imaginaire pour tout $i$. Pour tout $i$, notons $Z_{i} $ l'ensemble des  points $z\in {\mathbb C}$ tels que $y_{b_{0},i}+z\in 1+2D{\mathbb Z}$ et $C^{2i-1}\leq Re(z)\leq C^{2i}$. La condition impos\'ee \`a $C$ implique que $Z_{i}$ a au moins $h+1$ \'el\'ements. Notons $Z$ l'ensemble des $\lambda$ tels que $\lambda_{i}\in Z_{i}$ pour tout $i$. Un \'el\'ement $\lambda$ de cet ensemble convient, pourvu qu'il n'appartienne pas \`a ${\cal H}$. Il s'agit de prouver que $Z$ n'est pas inclus dans ${\cal H}$. Soit $H$ un hyperplan contenu dans ${\cal H}$. Il y a un indice $i$ tel que l'hyperplan vectoriel sous-jacent \`a $H$ ne contienne pas la droite port\'ee par $H_{i}$. Pour cet indice, la projection $\lambda\mapsto (\lambda_{l})_{l\not=i}$ est injective sur $H$. Elle envoie $Z\cap H$ dans $\prod_{l\not=i}Z_{l}$. Donc le nombre d\'el\'em\'ents de $Z\cap H$ est au plus $\prod_{l\not=i}\vert Z_{l}\vert$, ou encore $\vert Z\vert \vert Z_{i}\vert ^{-1}$, et ce nombre est major\'e par $(h+1)^{-1}\vert Z\vert $. Puisqu'il y a $h$ hyperplans, on en d\'eduit que $Z\cap {\cal H}$ a au plus $\frac{h}{h+1}\vert Z\vert $ \'el\'ements. Donc $Z$ n'est pas contenu dans ${\cal H}$, ce qui prouve (12).
 
 Appliquons (11) \`a un $\lambda$ v\'erifiant (12). Pour $b\in B^*$ et $i\in \{1,...,n\}$, on a $x_{b,i}>0$ d'apr\`es la condition de positivit\'e impos\'ee \`a $\check{\alpha}_{b,i}$. On a aussi $Re(\lambda_{i})>0$ par d\'efinition de $V[n]$. Donc les termes $y_{b,i}+\lambda_{i}x_{b,i}$  et $y_{b,i}+(1-y_{b_{0},i})x_{b,i}$ ont des parties r\'eelles de m\^eme signe. Leur diff\'erence appartient \`a $2{\mathbb Z}$. D'apr\`es (2),  la fonction $p$ prend la m\^eme valeur sur ces deux termes. Puisque $Re(y_{b,i}+(1-y_{b_{0},i})x_{b,i})=x_{b,i}$, la d\'efinition de cette fonction entra\^{\i}ne
 $p(y_{b,i}+\lambda_{i}x_{b,i})=0$ si $x_{b,i}<1$.  Supposons $b\not=b_{0}$. Prenons pour $i$ le plus petit indice tel que $x_{b,i}\not=1$. La condition (5) implique $x_{b,i}<1$. La contribution de $b$ \`a la formule (11) est donc nulle. Par contre, le m\^eme calcul montre que la contribution de $b_{0}$ est $p(1)^n$, c'est-\`a-dire $c_{b_{0}}(-2)^{-n}$. D'o\`u $(-2)^{-n}=0$, ce qui est  contradictoire puisque $c_{b_{0}}\not=0$. Cela ach\`eve la d\'emonstration. $\square$

  \bigskip
 
 \section{Endoscopie et applications $\theta_{\tilde{M}}^{rat,\tilde{G}}$, ${^c\theta}_{\tilde{M}}^{rat,\tilde{G}}$, ${^c\theta}_{\tilde{M}}^{\tilde{G}}$}
 
 \bigskip
 
 \subsection{Les applications stables}
 Toute les constructions de cette section sont  similaires \`a celles du cas non-archim\'edien. On se contentera la plupart du temps de rappeler les d\'efinitions et \'enonc\'es.
 
 Supposons $(G,\tilde{G},{\bf a})$ quasi-d\'eploy\'e et \`a torsion int\'erieure. Soit $\tilde{M}\in {\cal L}(\tilde{M}_{0})$. On note
 $$p^{st}_{\tilde{M}}:I_{ac}(\tilde{M}({\mathbb R}),K^M)\otimes Mes(M({\mathbb R}))\to S I_{ac}(\tilde{M}({\mathbb R}),K^M)\otimes Mes(M({\mathbb R}))$$
 la projection naturelle. Notons $\theta_{\tilde{M}}^{\tilde{G}}$ l'une de nos applications $\theta_{\tilde{M}}^{rat,\tilde{G}}$, ${^c\theta}_{\tilde{M}}^{rat,\tilde{G}}$, ${^c\theta}_{\tilde{M}}^{\tilde{G}}$. On d\'efinit  une application lin\'eaire
 $$S\theta_{\tilde{M}}^{\tilde{G}}:I(\tilde{G}({\mathbb R}),K)\otimes Mes(G({\mathbb R}))\to S I_{ac}(\tilde{M}({\mathbb R}),K^M)\otimes Mes(M({\mathbb R}))$$
 par la formule de r\'ecurrence
 $$S\theta_{\tilde{M}}^{\tilde{G}}({\bf f})=p_{\tilde{M}}^{st}\circ\theta_{\tilde{M}}^{\tilde{G}}({\bf f})-\sum_{s\in Z(\hat{M})^{\Gamma_{{\mathbb R}}}/Z(\hat{G})^{\Gamma_{{\mathbb R}}},s\not=1} i_{\tilde{M}}(\tilde{G},\tilde{G}'(s))S\theta_{{\bf M}}^{{\bf G}'(s)}({\bf f}^{{\bf G}'(s)}).$$
 Selon ce qu'est  l'application $\theta_{\tilde{M}}^{\tilde{G}}$, on notera  $S\theta_{\tilde{M}}^{rat,\tilde{G}}$, ${^cS}\theta_{\tilde{M}}^{rat,\tilde{G}}$, ${^cS}\theta_{\tilde{M}}^{\tilde{G}}$ l'application $S\theta_{\tilde{M}}^{\tilde{G}}$.
 
 \ass{Proposition  }{Les applications $S\theta_{\tilde{M}}^{rat,\tilde{G}}$, ${^cS}\theta_{\tilde{M}}^{rat,\tilde{G}}$ et ${^cS}\theta_{\tilde{M}}^{\tilde{G}}$ se quotientent en des applications lin\'eaires d\'efinies sur $SI(\tilde{G}({\mathbb R}),K)\otimes Mes(G({\mathbb R}))$.}
 
 Cela sera prouv\'e en 7.3.
 
 La proposition 2.3 de [VIII]  est valable pour chacune de nos applications.  
 
 {\bf Remarques.} (1) L'application $S\theta_{\tilde{M}}^{rat,\tilde{G}}$ peut \^etre d\'efinie comme une application lin\'eaire de $I(\tilde{G}({\mathbb R}))\otimes Mes(G({\mathbb R}))$ dans $ S I_{ac}(\tilde{M}({\mathbb R}))\otimes Mes(M({\mathbb R}))$. C'est-\`a-dire que l'on n'a pas besoin de supposer nos fonctions $K$-finies. L'application pr\'eserve n\'eanmoins la propri\'et\'e de $K$-finitude, puisqu'il en est de m\^eme de $\theta_{\tilde{M}}^{rat,\tilde{G}}$.
 
 (2) Il r\'esulte par r\'ecurrence des propri\'et\'es des applications $\theta_{\tilde{M}}^{\tilde{G}}$ que si l'on fixe un ensemble fini $\Omega$ de $K$-types, il existe un ensemble fini $\Omega^M$ de $K^M$-types tel que $S\theta_{\tilde{M}}^{\tilde{G}}$ envoie $I(\tilde{G}({\mathbb R}),\Omega)\otimes Mes(G({\mathbb R}))$ dans $ S I_{ac}(\tilde{M}({\mathbb R}),\Omega^M)\otimes Mes(M({\mathbb R}))$.
 
 \bigskip
 
 \subsection{Propri\'et\'es de l'application ${^cS}\theta_{\tilde{M}}^{rat,\tilde{G}}$}
 On suppose encore $(G,\tilde{G},{\bf a})$ quasi-d\'eploy\'e et \`a torsion int\'erieure.
 Soit $\tilde{\pi}\in D_{temp}^{st}(\tilde{M}({\mathbb R}))\otimes Mes(M({\mathbb R}))^*$ et ${\bf f}\in I(\tilde{G}({\mathbb R}),K)\otimes Mes(G({\mathbb R}))$.  On voit par r\'ecurrence que la fonction
 $$X\mapsto S^{\tilde{M}}(\tilde{\pi},X,{^cS}\theta_{\tilde{M}}^{rat,\tilde{G}}({\bf f}))$$
 est la transform\'ee de Fourier d'une fonction sur $i{\cal A}_{\tilde{M}}^*$, laquelle se prolonge en une fonction m\'eromorphe sur ${\cal A}_{\tilde{M},{\mathbb C}}^*$ qui v\'erifie les propri\'et\'es (2) \`a (5) de 5.2. On note cette fonction $\lambda\mapsto S^{\tilde{M}}(\tilde{\pi},\lambda,{^cS}\theta_{\tilde{M}}^{rat,\tilde{G}}({\bf f}))$. Pour $\nu\in {\cal A}_{\tilde{M}}^*$, on d\'efinit comme en 5.10  la transform\'ee de Fourier $X\mapsto S^{\tilde{M}}(\tilde{\pi},\nu,X,{^cS\theta}_{\tilde{M}}^{rat,\tilde{G}}({\bf f}))$. On voit  par r\'ecurrence que l'on a l'analogue suivant de la proposition 5.10.
 
 \ass{Proposition}{Supposons $\tilde{M}\not=\tilde{G}$ et $\tilde{\pi}$ elliptique. Soit ${\bf f}\in I(\tilde{G}({\mathbb R}),K)\otimes Mes(G({\mathbb R}))$. Si chaque point $\nu_{\tilde{S}}$ est assez positif relativement \`a $\tilde{S}$, on a l'\'egalit\'e
 $$\sum_{\tilde{S}\in {\cal P}(\tilde{M})}\omega_{\tilde{S}}(X)S^{\tilde{M}}(\tilde{\pi},\nu_{\tilde{S}},X,{^cS\theta}_{\tilde{M}}^{rat,\tilde{G}}({\bf f}))=0$$
 pour tout $X\in {\cal A}_{\tilde{M}}$.}
 
 Comme en 5.10, si l'on fixe un ensemble fini de $K$-types, on peut choisir les $\nu_{\tilde{S}}$ ind\'ependants de ${\bf f}$ et $\tilde{\pi}$, pourvu que ${\bf f}\in I(\tilde{G}({\mathbb R}),\Omega)\otimes Mes(G({\mathbb R}))$. 
 
 \bigskip
 
  \bigskip
 
 \subsection{Propri\'et\'es de l'application $S\theta_{\tilde{M}}^{rat,\tilde{G}}$}

 \ass{Lemme}{  Il existe une unique application lin\'eaire
$$\sigma_{\tilde{M}}^{\tilde{G}}:D_{temp,0}^{st}(\tilde{M}({\mathbb R}))\otimes Mes(M({\mathbb R}))^*\to U_{\tilde{M}}^{\tilde{G}}\otimes D_{temp,0}(\tilde{M}({\mathbb R}))\otimes Mes(M({\mathbb R}))^*$$
v\'erifiant la condition suivante. Soient ${\bf f}\in I(\tilde{G}({\mathbb R}),K)\otimes Mes(G({\mathbb R}))$, $\tilde{\pi}\in D_{temp,0}^{st}(\tilde{M}({\mathbb R}),\omega)\otimes Mes(M({\mathbb R}))^*$ et $X\in {\cal A}_{\tilde{M}}$. Alors on a l'\'egalit\'e
$$ S^{\tilde{M}}(\tilde{\pi},X,S\theta_{\tilde{M}}^{rat,\tilde{G}}({\bf f}))=\int_{i{\cal A}_{\tilde{M}}^*}I^{\tilde{M}}(\sigma_{\tilde{M}}^{\tilde{G}}(\tilde{\pi};\lambda),{\bf f}_{\tilde{M}})e^{-<\lambda,X>}\, d\lambda.$$}

Preuve. On utilise la d\'efinition de $S\theta_{\tilde{M}}^{rat,\tilde{G}}({\bf f})$ du paragraphe pr\'ec\'edent. Pour le premier terme, on a
$$S^{\tilde{M}}(\tilde{\pi},X,p^{st}_{\tilde{M}}\circ \theta_{\tilde{M}}^{rat,\tilde{G}}({\bf f}))=I^{\tilde{M}}(\tilde{\pi},X,\theta_{\tilde{M}}^{rat,\tilde{G}}({\bf f}))$$
puisque $\tilde{\pi}$ est stable, d'o\`u
$$S^{\tilde{M}}(\tilde{\pi},X,p^{st}_{\tilde{M}}\circ \theta_{\tilde{M}}^{rat,\tilde{G}}({\bf f}))=\int_{i{\cal A}_{\tilde{M}}^*}I^{\tilde{M}}(\rho_{\tilde{M}}^{\tilde{G}}(\tilde{\pi};\lambda),{\bf f}_{\tilde{M}})e^{-<\lambda,X>}\, d\lambda$$
d'apr\`es le lemme 5.7. Pour les autres termes, on utilise l'\'enonc\'e par r\'ecurrence:
$$S^{\tilde{M}}(\tilde{\pi},X,S\theta_{\tilde{M}}^{rat,{\bf G}'(s)}({\bf f}^{{\bf G}'(s)}))=\int_{i{\cal A}_{\tilde{M}}^*}I^{{\bf M}}(\sigma_{{\bf M}}^{{\bf G}'(s)}(\tilde{\pi};\lambda),({\bf f}^{{\bf G}'(\tilde{s})})_{{\bf M}})e^{-<\lambda,X>}\, d\lambda.$$
Evidemment, en identifiant l'espace de la donn\'ee endoscopique maximale ${\bf M}$ \`a $\tilde{M}$, on a simplement $({\bf f}^{{\bf G}'(\tilde{s})})_{{\bf M}})={\bf f}_{\tilde{M}}$. En d\'efinissant  $\sigma_{\tilde{M}}^{\tilde{G}}(\tilde{\pi})$ par
$$\sigma_{\tilde{M}}^{\tilde{G}}(\tilde{\pi})=\rho_{\tilde{M}}^{\tilde{G}}(\tilde{\pi})-\sum_{s\in Z(\hat{M})^{\Gamma_{{\mathbb R}}}/Z(\hat{G})^{\Gamma_{{\mathbb R}}},s\not=1} i_{\tilde{M}}(\tilde{G},\tilde{G}'(s))\sigma_{{\bf M}}^{{\bf G}'(s)}(\tilde{\pi}),$$
les formules pr\'ec\'edentes conduisent \`a l'\'egalit\'e de l'\'enonc\'e. L'unicit\'e de l'application $\sigma_{\tilde{M}}^{\tilde{G}}$ se prouve comme au lemme 5.7. $\square$

Avec les m\^emes notations qu'en 5.7(7), la d\'efinition par r\'ecurrence donn\'ee ci-dessus entra\^{\i}ne

(1) $\sigma_{\tilde{M}}^{\tilde{G}}$ envoie $Ind_{\tilde{R}}^{\tilde{M}}(D^{st}_{ell,0,\mu}(\tilde{R}({\mathbb R}),\omega)_{\lambda}\otimes Mes(R({\mathbb R}))^*)$ dans $U_{\tilde{M}}^{\tilde{G}}\otimes Ind_{\tilde{R}}^{\tilde{M}}(D_{ell,0,\mu}(\tilde{R}({\mathbb R}),\omega)_{\lambda}\otimes Mes(R({\mathbb R}))^*)$.

\bigskip

\subsection{Stabilit\'e de l'application $\sigma_{\tilde{M}}^{\tilde{G}}$}
On conserve les m\^emes hypoth\`eses.
\ass{Lemme }{L'application $\sigma_{\tilde{M}}^{\tilde{G}}$ prend ses valeurs dans $U_{\tilde{M}}^{\tilde{G}}\otimes D^{st}_{temp,0}(\tilde{M}({\mathbb R}))\otimes Mes(M({\mathbb R}))^*$.}
 
 Cela sera prouv\'e en 7.3.
 
 \bigskip
 
 \subsection{Une variante des int\'egrales orbitales pond\'er\'ees stables}
On suppose encore $(G,\tilde{G},{\bf a})$ quasi-d\'eploy\'e et \`a torsion int\'erieure. Soit $\boldsymbol{\delta}\in D_{orb}^{st}(\tilde{M}({\mathbb R}))\otimes Mes(M({\mathbb R}))^*$. On d\'efinit une forme lin\'eaire ${\bf f}\mapsto{^cS}_{\tilde{M}}^{\tilde{G}}(\boldsymbol{\delta},{\bf f})$ sur $I(\tilde{G}({\mathbb R}),K)\otimes Mes(G({\mathbb R}))$ par la formule habituelle
$$^cS_{\tilde{M}}^{\tilde{G}}(\boldsymbol{\delta},{\bf f})={^cI}_{\tilde{M}}^{\tilde{G}}(\boldsymbol{\delta},{\bf f})-\sum_{s\in Z(\hat{M})^{\Gamma_{{\mathbb R}}}/Z(\hat{G})^{\Gamma_{{\mathbb R}}},s\not=1} i_{\tilde{M}}(\tilde{G},\tilde{G}'(s)){^cS}_{{\bf M}}^{{\bf G}'(s)}(\boldsymbol{\delta},{\bf f}^{{\bf G}'(s)}).$$
La propri\'et\'e de compacit\'e 5.13(3) reste valable pour cette distribution.

\ass{Proposition }{Pour tout $\boldsymbol{\delta}\in D_{orb}^{st}(\tilde{M}({\mathbb R}))\otimes Mes(M({\mathbb R}))^*$, la forme lin\'eaire ${\bf f}\mapsto{^cS}_{\tilde{M}}^{\tilde{G}}(\boldsymbol{\delta},{\bf f})$ se quotiente en une forme lin\'eaire sur $SI(\tilde{G}({\mathbb R}),K)\otimes Mes(G({\mathbb R}))$.}

Cela sera prouv\'e en 7.3. 

\bigskip

\subsection{Les applications endoscopiques}
Le triplet $(G,\tilde{G},{\bf a})$ est quelconque. Soient $\tilde{M}\in {\cal L}(\tilde{M}_{0})$ et ${\bf M}'=(M',{\cal M}',\tilde{\zeta})$ une donn\'ee endoscopique elliptique et relevante de $(M,\tilde{M},{\bf a})$. Consid\'erons l'une de nos applications $\theta_{\tilde{M}}^{rat,\tilde{G}}$, ${^c\theta}_{\tilde{M}}^{rat,\tilde{G}}$, ${^c\theta}_{\tilde{M}}^{\tilde{G}}$ que l'on note $\theta_{\tilde{M}}^{\tilde{G}}$. On d\'efinit une application lin\'eaire
$$\theta_{\tilde{M}}^{\tilde{G},{\cal E}}({\bf M}'):I(\tilde{G}({\mathbb R}),\omega,K)\otimes Mes(G({\mathbb R}))\to SI_{ac}({\bf M}',K^{M'})\otimes Mes(M'({\mathbb R}))$$
par la formule
$$\theta_{\tilde{M}}^{\tilde{G},{\cal E}}({\bf M}',{\bf f})=\sum_{\tilde{s}\in \tilde{\zeta}Z(\hat{M})^{\Gamma_{{\mathbb R}},\hat{\theta}}/Z(\hat{G})^{\Gamma_{{\mathbb R}},\hat{\theta}}}i_{\tilde{M}'}(\tilde{G},\tilde{G}'(\tilde{s}))S\theta_{{\bf M}'}^{{\bf G}'(\tilde{s})}({\bf f}^{{\bf G}'(\tilde{s})}).$$
Dans le cas particulier o\`u $(G,\tilde{G},{\bf a})$ est quasi-d\'eploy\'e et \`a torsion int\'erieure et o\`u ${\bf M}'={\bf M}$, il faut remplacer le terme index\'e par $s=1$ par $S\theta_{\tilde{M}}^{\tilde{G}}({\bf f})$, cela parce que l'on n'a pas encore d\'emontr\'e que l'application $S\theta_{\tilde{M}}^{\tilde{G}}$ \'etait stable. 
On note suivant les cas cette application $\theta_{\tilde{M}}^{rat,\tilde{G},{\cal E}}({\bf M}')$, $^c\theta_{\tilde{M}}^{rat,\tilde{G},{\cal E}}({\bf M}')$, $^c\theta_{\tilde{M}}^{\tilde{G},{\cal E}}({\bf M}')$.

Consid\'erons maintenant un $K$-espace $(KG,K\tilde{G},{\bf a})$, un $K$-espace de Levi $K\tilde{M}\in {\cal L}(K\tilde{M}_{0})$ et une donn\'ee endoscopique ${\bf M}'$ de $(KM,K\tilde{M},{\bf a})$  qui est elliptique et relevante.    On d\'efinit par la m\^eme formule une application lin\'eaire 
$$\theta_{K\tilde{M}}^{K\tilde{G},{\cal E}}({\bf M}'):I(K\tilde{G}({\mathbb R}),\omega,K)\otimes Mes(G({\mathbb R}))\to SI_{ac}({\bf M}',K^{M'})\otimes Mes(M'({\mathbb R})).$$

\ass{Proposition}{Il existe une unique application lin\'eaire
$$\theta_{K\tilde{M}}^{K\tilde{G},{\cal E}}:I(K\tilde{G}({\mathbb R}),\omega,K)\otimes Mes(G({\mathbb R}))\to I_{ac}(K\tilde{M}({\mathbb R}),\omega,K^M)\otimes Mes(M({\mathbb R}))$$
telle que, pour toute donn\'ee endoscopique ${\bf M}'$ de $(KM,K\tilde{M},{\bf a})$ qui est elliptique et relevante et pour tout ${\bf f}\in I(K\tilde{G}({\mathbb R}),\omega,K)\otimes Mes(G({\mathbb R}))$, on ait l'\'egalit\'e
$$(\theta_{K\tilde{M}}^{K\tilde{G},{\cal E}}({\bf f}))^{{\bf M}'}=\theta_{K\tilde{M}}^{K\tilde{G},{\cal E}}({\bf M}',{\bf f}).$$}

L'application sera not\'ee selon le cas $\theta_{K\tilde{M}}^{rat,K\tilde{G},{\cal E}}$, $^c\theta_{K\tilde{M}}^{rat,K\tilde{G},{\cal E}}$, $^c\theta_{K\tilde{M}}^{K\tilde{G},{\cal E}}$. La preuve est la m\^eme que dans le cas non-archim\'edien. Il faut utiliser la version "$K$-finie" de la proposition 4.11 de [I], \`a savoir le corollaire 3.5 de [IV]. 

Ces applications v\'erifient les m\^emes propri\'et\'es que dans le cas non-archim\'edien, cf. en particulier [VIII] 3.6.  Soit $\tilde{M}$ une composante connexe de $K\tilde{M}$, soit $\tilde{\pi}\in D_{ell}(\tilde{M}({\mathbb R}),\omega)\otimes Mes(M({\mathbb R}))^*$ et ${\bf f}\in I(K\tilde{G}({\mathbb R}),\omega,K)\otimes Mes(G({\mathbb R}))$. L'application $X\mapsto I^{K\tilde{M}}(\tilde{\pi},X,\theta^{K\tilde{G},{\cal E}}_{K\tilde{M}}({\bf f}))$ est de Schwartz sur ${\cal A}_{\tilde{M}}$. Sa transform\'ee de Fourier  $\lambda\mapsto  I^{K\tilde{M}}(\tilde{\pi},\lambda,\theta^{K\tilde{G},{\cal E}}_{K\tilde{M}}({\bf f}))$ s'\'etend en une fonction m\'eromorphe sur ${\cal A}_{\tilde{M},{\mathbb C}}^*$. Elle est \`a d\'ecroissance rapide dans les bandes verticales. Ses p\^oles sont de la forme d\'ecrite en 5.2(3). Dans le cas de l'application $^c\theta_{K\tilde{M}}^{rat,K\tilde{G},{\cal E}}$ (et, en g\'en\'eral, seulement dans ce cas), les hyperplans polaires sont en nombre fini.

Rappelons une propri\'et\'e importante: 

(1) supposons $K\tilde{M}\not=K\tilde{G}$;  soit $\tilde{M}$ une composante connexe de $K\tilde{M}$, soit $\tilde{\pi}\in D_{ell}(\tilde{M}({\mathbb R}),\omega)\otimes Mes(M({\mathbb R}))^*$ et ${\bf f}\in I(K\tilde{G}({\mathbb R}),\omega,K)\otimes Mes(G({\mathbb R}))$; si chaque point $\nu_{\tilde{S}}$ est assez positif relativement \`a $\tilde{S}$, on a l'\'egalit\'e
$$\sum_{\tilde{S}\in {\cal P}(\tilde{M})}\omega_{\tilde{S}}(X)I^{K\tilde{M}}(\tilde{\pi},\nu_{\tilde{S}},X,{^c\theta}_{K\tilde{M}}^{rat,K\tilde{G},{\cal E}}({\bf f}))=0$$
pour tout $X\in {\cal A}_{\tilde{M}}$. 

Les propri\'et\'es de $K$-finitude du transfert entra\^{\i}nent que, si l'on fixe un ensemble fini $\Omega$ de $K$-types, il existe un ensemble fini $\Omega^M$ de $K^M$-types de sorte que $\theta_{K\tilde{M}}^{K\tilde{G},{\cal E}}$ envoie $I(K\tilde{G}({\mathbb R}),\omega,\Omega)\otimes Mes(G({\mathbb R}))$ dans $ I_{ac}(K\tilde{M}({\mathbb R}),\omega,\Omega^M)\otimes Mes(M({\mathbb R}))$. La remarque qui suit la proposition 5.10 vaut pour l'application ${^c\theta}_{K\tilde{M}}^{rat, K\tilde{G},{\cal E}}$. C'est-\`a-dire que  les p\^oles de la  fonction $\lambda\mapsto I^{K\tilde{M}}(\tilde{\pi},\lambda,{^c\theta}_{K\tilde{M}}^{rat,K\tilde{G},{\cal E}}({\bf f}))$ restent dans un nombre fini d'hyperplans ind\'ependants de ${\bf f}$, pourvu que ${\bf f}\in I(K\tilde{G}({\mathbb R}),\omega,\Omega)\otimes Mes(G({\mathbb R}))$. 

{\bf Remarque.}  Ici encore, l'application $\theta_{K\tilde{M}}^{rat,K\tilde{G},{\cal E}}$ peut \^etre d\'efinie sur des fonctions qui ne sont pas $K$-finies.

\bigskip

\subsection{Egalit\'e d'applications lin\'eaires}
Consid\'erons un $K$-espace $(KG,K\tilde{G},{\bf a})$ et un $K$-espace de Levi $K\tilde{M}\in {\cal L}(K\tilde{M}_{0})$. Notons $K\tilde{G}=(\tilde{G}_{p})_{p\in \Pi}$, $K\tilde{M}=(\tilde{M}_{p})_{p\in \Pi^M}$. Consid\'erons l'une de nos applications g\'en\'eriques $\theta_{\tilde{M}}^{rat,\tilde{G}}$, ${^c\theta}_{\tilde{M}}^{rat,\tilde{G}}$ ou ${^c\theta}_{\tilde{M}}^{\tilde{G}}$, que l'on note simplement $\theta_{\tilde{M}}^{\tilde{G}}$. On d\'efinit l'application lin\'eaire
$$\theta_{K\tilde{M}}^{K\tilde{G}} :I(K\tilde{G}({\mathbb R}),\omega,K)\otimes Mes(G({\mathbb R}))\to I_{ac}(K\tilde{M}({\mathbb R}),\omega,K^M)\otimes Mes(M({\mathbb R}))$$
comme la somme directe des applications $\theta_{\tilde{M}_{p}}^{\tilde{G}_{p}}$ pour $p\in \Pi^M$ et des applications nulles pour $p\in \Pi-\Pi^M$. 

\ass{Proposition (\`a prouver)}{On a l'\'egalit\'e $\theta_{K\tilde{M}}^{K\tilde{G},{\cal E}}=\theta_{K\tilde{M}}^{K\tilde{G}}$.}

On donnera en 7.4 une preuve soumise \`a une hypoth\`ese qui ne sera v\'erifi\'ee que plus tard. 

\bigskip

\subsection{Propri\'et\'es de l'application $\theta_{K\tilde{M}}^{rat,K\tilde{G},{\cal E}}$}
Consid\'erons un $K$-espace $(KG,K\tilde{G},{\bf a})$ et un $K$-espace de Levi $K\tilde{M}$. Ecrivons $K\tilde{G}=(\tilde{G}_{p})_{p\in \Pi}$, $K\tilde{M}=(\tilde{M}_{p})_{p\in \Pi^M}$. On d\'efinit $D_{temp}(K\tilde{M}({\mathbb R}),\omega)$ comme la somme directe des espaces $D_{temp}(\tilde{M}_{p}({\mathbb R}),\omega)$ pour $p\in \Pi^M$. Pour $\tilde{\pi}=\sum_{p\in \Pi^M}\tilde{\pi}_{p}$ dans cet espace et pour ${\bf f}=({\bf f}_{p})_{p\in \Pi^M}\in I(\tilde{M}({\mathbb R}),\omega)\otimes Mes(M({\mathbb R}))$, on pose
$$I^{K\tilde{M}}(\tilde{\pi},{\bf f})=\sum_{p\in \Pi^M}I^{\tilde{M}_{p}}(\tilde{\pi}_{p},{\bf f}_{p}).$$
On a diverses variantes: $D_{temp,0}(K\tilde{M}({\mathbb R}),\omega)$ etc...
\ass{Lemme}{  Il existe une unique application lin\'eaire
$$\rho_{K\tilde{M}}^{K\tilde{G},{\cal E}}:D_{temp,0}(K\tilde{M}({\mathbb R}),\omega)\otimes Mes(M({\mathbb R}))^*\to U_{\tilde{M}}^{\tilde{G}}\otimes D_{temp,0}(K\tilde{M}({\mathbb R}),\omega)\otimes Mes(M({\mathbb R}))^*$$
v\'erifiant la condition suivante. Soient ${\bf f}\in I(K\tilde{G}({\mathbb R}),K)\otimes Mes(G({\mathbb R}))$, $\tilde{\pi}\in D_{temp,0}(K\tilde{M}({\mathbb R}),\omega)\otimes Mes(M({\mathbb R}))^*$ et $X\in {\cal A}_{\tilde{M}}$. Alors on a l'\'egalit\'e
$$ I^{K\tilde{M}}(\tilde{\pi},X,\theta_{K\tilde{M}}^{rat,K\tilde{G},{\cal E}}({\bf f}))=\int_{i{\cal A}_{\tilde{M}}^*}I^{K\tilde{M}}(\rho_{K\tilde{M}}^{K\tilde{G}}(\tilde{\pi};\lambda),{\bf f}_{K\tilde{M}})e^{-<\lambda,X>}\, d\lambda.$$}

Preuve.  On introduit un ensemble ${\cal E}(K\tilde{M},{\bf a})$ de repr\'esentants des classes d'\'equivalence de donn\'ees endoscopiques de $(KM,K\tilde{M},{\bf a})$ qui sont elliptiques et relevantes. Pour tout ${\bf M}'\in {\cal E}(K\tilde{M},{\bf a})$, on peut introduire un \'el\'ement $\tilde{\pi}_{{\bf M}'}\in D^{st}_{temp,0}({\bf M}')\otimes Mes(M'({\mathbb R}))^*$ de sorte que
$$\tilde{\pi}=\sum_{{\bf M}'\in {\cal E}(K\tilde{M},{\bf a})}transfert(\tilde{\pi}_{{\bf M}'}).$$
On n'a aucun mal \`a prouver l'\'egalit\'e
$$ I^{K\tilde{M}}(\tilde{\pi},X,\theta_{K\tilde{M}}^{rat,K\tilde{G},{\cal E}}({\bf f}))=\sum_{{\bf M}'\in {\cal E}(K\tilde{M},{\bf a})}S^{{\bf M}'}(\tilde{\pi}_{{\bf M}'},X,(\theta_{K\tilde{M}}^{rat,K\tilde{G},{\cal E}}({\bf f}))^{{\bf M}'}).$$
A droite, on a identifi\'e  $X$ \`a un \'el\'ement de ${\cal A}_{\tilde{M}'}$ par l'isomorphisme entre cet espace et ${\cal A}_{\tilde{M}}$.
Fixons ${\bf M}'\in {\cal E}(K\tilde{M},{\bf a})$, que l'on \'ecrit  ${\bf M}'=(M',{\cal M}',\tilde{\zeta})$. Par d\'efinition de $\theta_{K\tilde{M}}^{rat,K\tilde{G},{\cal E}}({\bf f})$, le dernier terme ci-dessus  est aussi
$$S^{{\bf M}'}(\tilde{\pi}_{{\bf M}'},X,\theta_{K\tilde{M}}^{rat,K\tilde{G},{\cal E}}({\bf M}',{\bf f})).$$
Ou encore 
$$\sum_{\tilde{s}\in \tilde{\zeta}Z(\hat{M})^{\Gamma_{{\mathbb R}},\hat{\theta}}/Z(\hat{G})^{\Gamma_{{\mathbb R}},\hat{\theta}}}i_{\tilde{M}'}(\tilde{G},\tilde{G}'(\tilde{s}))S^{{\bf M}'}(\tilde{\pi}_{{\bf M}'},X,S\theta_{{\bf M}'}^{rat,{\bf G}'(\tilde{s})}({\bf f}^{{\bf G}'(\tilde{s})})).$$
Fixons $\tilde{s}$ et appliquons le lemme 6.2. Nos hypoth\`eses de r\'ecurrence nous autorisent \`a utiliser le lemme 6.3. On obtient
$$S^{{\bf M}'}(\tilde{\pi}_{{\bf M}'},X,S\theta_{{\bf M}'}^{rat,{\bf G}'(\tilde{s})}({\bf f}^{{\bf G}'(\tilde{s})}))=\int_{i{\cal A}_{\tilde{M}'}^*}S^{{\bf M}'}(\sigma_{{\bf M}'}^{{\bf G}'(\tilde{s})}(\tilde{\pi}_{{\bf M}'};\lambda),({\bf f}^{{\bf G}'(\tilde{s})})_{{\bf M}'}) e^{-<\lambda,X>}\,d\lambda.$$
    On a $({\bf f}^{{\bf G}'(\tilde{s})})_{{\bf M}'}=({\bf f}_{K\tilde{M},\omega})^{{\bf M}'}$. Alors
  $$S^{{\bf M}'}(\sigma_{{\bf M}'}^{{\bf G}'(\tilde{s})}(\tilde{\pi}_{{\bf M}'};\lambda),({\bf f}^{{\bf G}'(\tilde{s})})_{{\bf M}'})=I^{K\tilde{M}}(transfert(\sigma_{{\bf M}'}^{{\bf G}'(\tilde{s})}(\tilde{\pi}_{{\bf M}'};\lambda)),{\bf f}_{K\tilde{M},\omega}).$$
  On peut identifier ${\cal A}_{\tilde{M}'}^*$ \`a ${\cal A}_{\tilde{M}}^*$. L'application $\sigma_{{\bf M}'}^{{\bf G}'(\tilde{s})}$ prend ses valeurs dans $U_{\tilde{M}'}^{\tilde{G}'(\tilde{s})}\otimes D_{temp,0}^{st}({\bf M}')$. Ecrivons $\sigma_{{\bf M}'}^{{\bf G}'(\tilde{s})}(\tilde{\pi}_{{\bf M}'})=\sum_{j=1,...,k}u_{j}\otimes \tilde{\tau}_{j}$. Alors 
  $$transfert(\sigma_{{\bf M}'}^{{\bf G}'(\tilde{s})}(\tilde{\pi}_{{\bf M}'};\lambda))=\sum_{j=1,...,k}u_{j}(\lambda)transfert(\tilde{\tau}_{j,\lambda}).$$
  Cela s'interpr\`ete comme la valeur en $\lambda$ de l'\'el\'ement $\sum_{j=1,...,k}u_{j}\otimes transfert(\tilde{\tau}_{j})$. Pour interpr\'eter ce dernier terme comme un \'el\'ement de notre espace $U_{\tilde{M}}^{\tilde{G}}\otimes D_{temp,0}(K\tilde{M}({\mathbb R}),\omega)$, il faut remarquer que 
  
  (1) $U_{\tilde{M}'}^{\tilde{G}'(\tilde{s})}$ est inclus dans $U_{\tilde{M}}^{\tilde{G}}$. 
  
  Introduisons des paires de Borel $(B,T)$ de $K\tilde{G}$ et $(B',T')$ de $G'(\tilde{s})$ de sorte que $M$  soit standard pour $(B,T)$ et $M'$ soit standard pour $(B',T')$. L'ensemble $\check{\Sigma}_{\star }(A_{\tilde{M}})$ de 5.4 est celui des restrictions non nulles \`a ${\cal A}_{\tilde{M}}^*$ d'\'el\'ements de $\check{\Sigma}(T)$. L'ensemble similaire $\check{\Sigma}_{\star }(A_{\tilde{M}'})$ est celui des restrictions non nulles \`a ${\cal A}_{\tilde{M}'}^*$ d'\'el\'ements de $\check{\Sigma}(T')$.  On a un isomorphisme $T'\simeq T/(1-\theta)(T)$, d'o\`u un homomorphisme $X_{*}(T)\to X_{*}(T')$. Les descriptions usuelles des syst\`emes de racines montrent que $\check{\Sigma}(T')$ est contenu dans l'image de $\check{\Sigma}(T)$ par cet homomorphisme. Modulo l'identification de ${\cal A}_{\tilde{M}'}^*$ \`a ${\cal A}_{\tilde{M}}^*$, on en d\'eduit ais\'ement l'inclusion $\check{\Sigma}_{\star }(A_{\tilde{M}'})\subset \check{\Sigma}_{\star }(A_{\tilde{M}})$. L'assertion (1) r\'esulte alors de la d\'efinition donn\'ee en 5.4 des espaces $U_{\tilde{M}'}^{\tilde{G}'(\tilde{s})}$ et $U_{\tilde{M}}^{\tilde{G}}$.

  D\'efinissons $\rho_{K\tilde{M}}^{K\tilde{G},{\cal E}}({\bf M}',\tilde{\pi}_{{\bf M}'})$ par la formule
  $$(2) \qquad \rho_{K\tilde{M}}^{K\tilde{G},{\cal E}}({\bf M}',\tilde{\pi}_{{\bf M}'})=\sum_{\tilde{s}\in \tilde{\zeta}Z(\hat{M})^{\Gamma_{{\mathbb R}},\hat{\theta}}/Z(\hat{G})^{\Gamma_{{\mathbb R}},\hat{\theta}}}i_{\tilde{M}'}(\tilde{G},\tilde{G}'(\tilde{s}))transfert(\sigma_{{\bf M}'}^{{\bf G}'(\tilde{s})}(\tilde{\pi}_{{\bf M}'})).$$
  On obtient l'\'egalit\'e
  $$(3) \qquad S^{{\bf M}'}(\tilde{\pi}_{{\bf M}'},X,(\theta_{K\tilde{M}}^{rat,K\tilde{G},{\cal E}}({\bf f}))^{{\bf M}'})=\int_{i{\cal A}_{\tilde{M}}^*}I^{K\tilde{M}}(\rho_{K\tilde{M}}^{K\tilde{G},{\cal E}}({\bf M}',\tilde{\pi}_{{\bf M}'},\lambda),{\bf f}_{K\tilde{M},\omega})e^{-<\lambda,X>}\,d\lambda.$$
  D\'efinissons $\rho_{K\tilde{M}}^{K\tilde{G},{\cal E}}(\tilde{\pi})$ par
  $$\rho_{K\tilde{M}}^{K\tilde{G},{\cal E}}(\tilde{\pi})=\sum_{{\bf M}'\in {\cal E}(K\tilde{M},{\bf a})}\rho_{K\tilde{M}}^{K\tilde{G},{\cal E}}({\bf M}',\tilde{\pi}_{{\bf M}'}).$$
  On obtient alors la formule de l'\'enonc\'e. L'unicit\'e se prouve comme au lemme 5.7. $\square$
  
  \bigskip
  
  \subsection{Egalit\'e des fonctions $\rho_{K\tilde{M}}^{K\tilde{G}}$ et $\rho_{K\tilde{M}}^{K\tilde{G},{\cal E}}$}
  On conserve la m\^eme situation.
  On d\'efinit 
  $$\rho_{K\tilde{M}}^{K\tilde{G}}:D_{temp,0}(K\tilde{M}({\mathbb R}))\otimes Mes(M({\mathbb R}))^*\to U_{\tilde{M}}^{\tilde{G}}\otimes D_{temp,0}(K\tilde{M}({\mathbb R}))\otimes Mes(M({\mathbb R}))^*$$
   comme la somme directe des applications $\rho_{\tilde{M}_{p}}^{\tilde{G}_{p}}$ pour $p\in \Pi^M$. 
   \ass{Lemme (\`a prouver)}{On a l'\'egalit\'e $\rho_{K\tilde{M}}^{K\tilde{G},{\cal E}}=\rho_{K\tilde{M}}^{K\tilde{G}}$.}
   
   On donnera en 7.4 une preuve soumise \`a une hypoth\`ese qui ne sera v\'erifi\'ee que plus tard. 

\bigskip

\subsection{Variante des int\'egrales orbitales pond\'er\'ees elliptiques}
Soit $(G,\tilde{G},{\bf a})$ un triplet g\'en\'eral et $\tilde{M}$ un espace de Levi de $\tilde{G}$. Soit ${\bf M}'=(M',{\cal M}',\tilde{\zeta})$ une donn\'ee endoscopique elliptique et relevante de $(M,\tilde{M},{\bf a})$. Pour $\boldsymbol{\delta}\in D_{orb,\tilde{G}-reg}^{st}({\bf M}')\otimes Mes(M'({\mathbb R}))^*$ et ${\bf f}\in I(\tilde{G}({\mathbb R}),\omega,K)\otimes Mes(G({\mathbb R}))$, on pose
$$^cI_{\tilde{M}}^{\tilde{G},{\cal E}}({\bf M}',\boldsymbol{\delta},{\bf f})=\sum_{\tilde{s}\in \tilde{\zeta}Z(\hat{M})^{\Gamma_{{\mathbb R}},\hat{\theta}}/Z(\hat{G})^{\Gamma_{{\mathbb R}},\hat{\theta}}}i_{\tilde{M}'}(\tilde{G},\tilde{G}'(\tilde{s})) {^cS}_{{\bf M}'}^{{\bf G}'(\tilde{s})}(\boldsymbol{\delta},{\bf f}^{{\bf G}'(\tilde{s})}).$$
Ici encore, dans le cas particulier o\`u $(G,\tilde{G},{\bf a})$ est quasi-d\'eploy\'e et \`a torsion int\'erieure et o\`u ${\bf M}'={\bf M}$, on doit remplacer le terme index\'e par $s=1$ par $^cS_{\tilde{M}}^{\tilde{G}}(\boldsymbol{\delta},{\bf f})$.

Consid\'erons maintenant un $K$-espace $(KG,K\tilde{G},{\bf a})$ et un $K$-espace de Levi $K\tilde{M}$. Soit ${\bf M}'=(M',{\cal M}',\tilde{\zeta})$ une donn\'ee endoscopique elliptique et relevante de $(KM,K\tilde{M},{\bf a})$. Pour $\boldsymbol{\delta}\in D_{orb,\tilde{G}-reg}^{st}({\bf M}')\otimes Mes(M'({\mathbb R}))^*$ et ${\bf f}\in I(K\tilde{G}({\mathbb R}),\omega,K)\otimes Mes(G({\mathbb R}))$, on d\'efinit $^cI_{K\tilde{M}}^{K\tilde{G},{\cal E}}({\bf M}',\boldsymbol{\delta},{\bf f})$ par la m\^eme formule que ci-dessus. 

 Soit maintenant $\boldsymbol{\gamma}\in D_{orb,\tilde{G}-reg}(K\tilde{M}({\mathbb R}),\omega)\otimes Mes(M({\mathbb R}))^*$. On peut d\'ecomposer $\boldsymbol{\gamma}$ en $\sum_{p\in \Pi^M}\boldsymbol{\gamma}_{p}$, o\`u $\boldsymbol{\gamma}_{p}\in D_{orb,\tilde{G}-reg}(\tilde{M}_{p}({\mathbb R}),\omega)\otimes Mes(M({\mathbb R}))^*$. Pour ${\bf f}=({\bf f}_{p})_{p\in \Pi}\in 
 I(K\tilde{G}({\mathbb R}),\omega,K)\otimes Mes(G({\mathbb R}))$, on pose
 $${^cI}_{K\tilde{M}}^{K\tilde{G}}(\boldsymbol{\gamma},{\bf f})=\sum_{p\in \Pi^M}{^cI}_{\tilde{M}_{p}}^{\tilde{G}_{p}}(\boldsymbol{\gamma}_{p},{\bf f}_{p}).$$
 
 On introduit un ensemble ${\cal E}(K\tilde{M},{\bf a})$ de repr\'esentants des classes d'\'equivalence de donn\'ees endoscopiques de $(KM,K\tilde{M},{\bf a})$ qui sont elliptiques et relevantes. On \'ecrit
$$\boldsymbol{\gamma}=\sum_{{\bf M}'\in {\cal E}(K\tilde{M},{\bf a})}transfert(\boldsymbol{\delta}_{{\bf M}'}),$$
o\`u les $\boldsymbol{\delta}_{{\bf M}'}$ appartiennent \`a $D_{orb,\tilde{G}-reg}^{st}({\bf M}')\otimes Mes(M'({\mathbb R}))^*$. On pose
$$^cI_{K\tilde{M}}^{K\tilde{G},{\cal E}}(\boldsymbol{\gamma},{\bf f})=\sum_{{\bf M}'\in {\cal E}(K\tilde{M},{\bf a})}{^cI}_{K\tilde{M}}^{K\tilde{G},{\cal E}}({\bf M}',\boldsymbol{\delta}_{{\bf M}'},{\bf f}).$$
Cela ne d\'epend pas de la d\'ecomposition de $\boldsymbol{\gamma}$ choisie (la preuve est analogue \`a celle de [II] 1.15).

\ass{Proposition (\`a prouver)}{Pour tout $\boldsymbol{\gamma}\in D_{orb,\tilde{G}-reg}(K\tilde{M}({\mathbb R}),\omega)\otimes Mes(M({\mathbb R}))^*$ et tout ${\bf f}\in I(K\tilde{G}({\mathbb R}),\omega,K)\otimes Mes(G({\mathbb R}))$, on a l'\'egalit\'e
$$^cI_{K\tilde{M}}^{K\tilde{G},{\cal E}}(\boldsymbol{\gamma},{\bf f})={^cI}_{K\tilde{M}}^{K\tilde{G}}(\boldsymbol{\gamma},{\bf f}).$$}

On donnera en 7.4 une preuve soumise \`a une hypoth\`ese qui ne sera v\'erifi\'ee que plus tard. 

\bigskip

\subsection{Reformulation des \'enonc\'es dans le cas quasi-d\'eploy\'e et \`a torsion int\'erieure}
Soit $(G,\tilde{G},{\bf a})$ un triplet quasi-d\'eploy\'e et \`a torsion int\'erieure. On n'a pas envie de le plonger dans un $K$-espace car cela perturberait nos hypoth\`eses de r\'ecurrence. Mais on peut reformuler les  assertions des paragraphes 6.6 \`a 6.10 en se passant d'un tel plongement. 

Ainsi, soit $\tilde{M}$ un espace de Levi de $\tilde{G}$ et soit ${\bf M}'=(M',{\cal M}',\zeta)$ une donn\'ee endoscopique elliptique et relevante de $(M,\tilde{M})$.  Notons $\theta_{\tilde{M}}^{\tilde{G}}$ l'une de nos applications $\theta_{\tilde{M}}^{rat,\tilde{G}}$, ${^c\theta}_{\tilde{M}}^{rat,\tilde{G}}$, ${^c\theta}_{\tilde{M}}^{\tilde{G}}$. La proposition 6.7 se reformule par:

(A) pour tout ${\bf f}\in I(\tilde{G}({\mathbb R}),K)\otimes Mes(G({\mathbb R}))$, on a l'\'egalit\'e
$$\theta_{\tilde{M}}^{\tilde{G},{\cal E}}({\bf M}',{\bf f})=(\theta_{\tilde{M}}^{\tilde{G}}({\bf f}))^{{\bf M}'}.$$

Soit  $\boldsymbol{\delta}\in D_{orb,\tilde{G}-reg}^{st}({\bf M}')\otimes Mes(M'({\mathbb R}))^*$. La proposition 6.10 se reformule par

(B) pour tout ${\bf f}\in I(\tilde{G}({\mathbb R}),K)\otimes Mes(G({\mathbb R}))$, on a l'\'egalit\'e 
$$^cI_{\tilde{M}}^{\tilde{G},{\cal E}}({\bf M}',\boldsymbol{\delta},{\bf f})={^cI}_{\tilde{M}}^{\tilde{G}}(transfert(\boldsymbol{\delta}),{\bf f}).$$

Soit $\tilde{\pi}_{{\bf M}'}\in D_{temp,0}^{st}({\bf M}')\otimes Mes(M'({\mathbb R}))^*$. On peut poser comme en 6.8(2):
$$ \rho_{\tilde{M}}^{\tilde{G},{\cal E}}({\bf M}',\tilde{\pi}_{{\bf M}'})=\sum_{s\in \zeta Z(\hat{M})^{\Gamma_{{\mathbb R}}}/Z(\hat{G})^{\Gamma_{{\mathbb R}}}}i_{\tilde{M}'}(\tilde{G},\tilde{G}'(s))transfert(\sigma_{{\bf M}'}^{{\bf G}'(s)}(\tilde{\pi}_{{\bf M}'})).$$
Pour ${\bf f}\in I(\tilde{G}({\mathbb R}),K)\otimes Mes(G({\mathbb R}))$, l'\'egalit\'e 6.8(3) (qui est d\'ej\`a prouv\'ee) peut se r\'ecrire
  $$(1) \qquad S^{{\bf M}'}(\tilde{\pi}_{{\bf M}'},X,\theta_{\tilde{M}}^{rat,\tilde{G},{\cal E}}({\bf M}',{\bf f}))=\int_{i{\cal A}_{\tilde{M}}^*}I^{\tilde{M}}(\rho_{\tilde{M}}^{\tilde{G},{\cal E}}({\bf M}',\tilde{\pi}_{{\bf M}'},\lambda),{\bf f}_{\tilde{M},\omega})e^{-<\lambda,X>}\,d\lambda.$$
  Le lemme 6.9 peut se r\'ecrire
  
(C) pour tout $\tilde{\pi}_{{\bf M}'}\in D_{temp,0}^{st}({\bf M}')\otimes Mes(M'({\mathbb R}))^*$, on a l'\'egalit\'e
$$\rho_{\tilde{M}}^{\tilde{G}}(transfert(\tilde{\pi}_{{\bf M}'}))=  \rho_{\tilde{M}}^{\tilde{G},{\cal E}}({\bf M}',\tilde{\pi}_{{\bf M}'}).$$

Les assertions (A), (B) et (C) seront prouv\'ees en 7.5.

On aura besoin d'un substitut de la propri\'et\'e 6.6(1) (qui est d\'ej\`a prouv\'ee). Soit $\tilde{\pi}_{{\bf M}'}\in D_{ell}({\bf M}')\otimes Mes(M'({\mathbb R}))^*$ et ${\bf f}\in I(\tilde{G}({\mathbb R}),K)\otimes Mes(G({\mathbb R}))$. On n'a aucun mal \`a d\'efinir les termes $S^{{\bf M}'}(\tilde{\pi}_{{\bf M}'},\nu,X,{^c\theta}_{\tilde{M}}^{rat,\tilde{G},{\cal E}}({\bf M}',{\bf f}))$ pour $\nu\in {\cal A}_{\tilde{M}}^*$ et $X\in {\cal A}_{\tilde{M}}$. Alors

(2) supposons $\tilde{M}\not=\tilde{G}$; si chaque point $\nu_{\tilde{S}}$ est assez positif relativement \`a $\tilde{S}$, on a l'\'egalit\'e
$$\sum_{\tilde{S}\in {\cal P}(\tilde{M})}\omega_{\tilde{S}}(X)S^{{\bf M}'}(\tilde{\pi}_{{\bf M}'},\nu,X,{^c\theta}_{\tilde{M}}^{rat,\tilde{G},{\cal E}}({\bf M}',{\bf f}))=0$$
pour tout $X\in {\cal A}_{\tilde{M}}$.

\bigskip

\section{Les preuves des assertions de la section 6 }

\bigskip

\subsection{Lien entre les int\'egrales orbitales pond\'er\'ees endoscopiques et leurs variantes}
Dans ce paragraphe et le suivant, on consid\`ere l'une des trois situations

(A)   on se donne un triplet $(KG,K\tilde{G},{\bf a})$ et un $K$-espace de Levi $K\tilde{M}\in {\cal L}(K\tilde{M}_{0})$;

(B) on se donne un triplet  $(G,\tilde{G},{\bf a})$ quasi-d\'eploy\'e et \`a torsion int\'erieure et un espace de Levi $\tilde{M}$;

(C) on se donne un triplet  $(G,\tilde{G},{\bf a})$ quasi-d\'eploy\'e et \`a torsion int\'erieure, un espace de Levi $\tilde{M}$ et une donn\'ee endoscopique ${\bf M}'=(M',{\cal M}',\zeta)$ de $(M,\tilde{M})$ qui est elliptique et relevante. 

\ass{Lemme}{(i) Dans le cas (A), soient $\boldsymbol{\gamma }\in D_{orb,\tilde{G}-reg}(K\tilde{M}({\mathbb R}),\omega)\otimes Mes(M({\mathbb R}))^*$ et ${\bf f}\in I(K\tilde{G}({\mathbb R}),\omega,K)\otimes Mes(G({\mathbb R}))$. On a l'\'egalit\'e
$$^cI_{K\tilde{M}}^{K\tilde{G},{\cal E}}(\boldsymbol{\gamma},{\bf f})=\sum_{K\tilde{L}\in {\cal L}(K\tilde{M})}I_{K\tilde{M}}^{K\tilde{L},{\cal E}}(\boldsymbol{\gamma},{^c\theta}_{K\tilde{L}}^{K\tilde{G},{\cal E}}({\bf f})).$$

(ii) Dans le cas (B),  soient $\boldsymbol{\delta }\in D^{st}_{orb,\tilde{G}-reg}(\tilde{M}({\mathbb R}))\otimes Mes(M({\mathbb R}))^*$ et ${\bf f}\in I(\tilde{G}({\mathbb R}),K)\otimes Mes(G({\mathbb R}))$. On a l'\'egalit\'e
$$^cS_{\tilde{M}}^{\tilde{G}}(\boldsymbol{\delta},{\bf f})=\sum_{\tilde{L}\in {\cal L}(\tilde{M})}S_{\tilde{M}}^{\tilde{L}}(\boldsymbol{\delta},{^cS\theta}_{\tilde{L}}^{\tilde{G}}({\bf f})).$$

(iii) Dans le cas (C), soient $\boldsymbol{\delta}\in D^{st}_{orb,\tilde{G}-reg}({\bf M}')\otimes Mes(M'({\mathbb R}))^*$ et ${\bf f}\in I(\tilde{G}({\mathbb R}),K)\otimes Mes(G({\mathbb R}))$. On a l'\'egalit\'e
$$^cI_{\tilde{M}}^{\tilde{G},{\cal E}}({\bf M}',\boldsymbol{\delta},{\bf f})=S^{{\bf M}'}(\boldsymbol{\delta},{^c\theta}_{\tilde{M}}^{\tilde{G},{\cal E}}({\bf M}',{\bf f}))+\sum_{\tilde{L}\in {\cal L}(\tilde{M}),\tilde{L}\not=\tilde{M}}I_{\tilde{M}}^{\tilde{L},{\cal E}}({\bf M}',\boldsymbol{\delta}, {^c\theta}_{\tilde{L}}^{\tilde{G}}({\bf f})).$$}

Preuve.
La preuve du (i) est similaire \`a celle du cas non-archim\'edien ([VIII] proposition 4.1).  Traitons (iii). Supposons d'abord que ${\bf M}'$ est la donn\'ee triviale ${\bf M}$. La distribution ${\bf f}\mapsto {^cS}_{\tilde{M}}^{\tilde{G}}(\boldsymbol{\delta},{\bf f})$ est d\'efinie de sorte que l'on ait l'\'egalit\'e 
$$ ^cI_{\tilde{M}}^{\tilde{G},{\cal E}}({\bf M},\boldsymbol{\delta},{\bf f})={^cI}_{\tilde{M}}^{\tilde{G}}(\boldsymbol{\delta},{\bf f}).$$
De m\^eme, pour $\tilde{L}\in {\cal L}(\tilde{M})$ avec $\tilde{L}\not=\tilde{M}$, la distribution $\boldsymbol{\varphi}\mapsto S_{\tilde{M}}^{\tilde{L}}(\boldsymbol{\delta},\boldsymbol{\varphi})$ est d\'efinie de sorte que l'on ait l'\'egalit\'e
$$I_{\tilde{M}}^{\tilde{L},{\cal E}}({\bf M},\boldsymbol{\delta},\boldsymbol{\varphi})=I_{\tilde{M}}^{\tilde{L}}(\boldsymbol{\delta},\boldsymbol{\varphi}).$$
L'application ${^cS}\theta_{\tilde{M}}^{\tilde{G}}$ est d\'efinie de sorte que l'on ait l'\'egalit\'e
$$^c\theta_{\tilde{M}}^{\tilde{G},{\cal E}}({\bf M},{\bf f})=p_{\tilde{M}}^{st}\circ{^c\theta}_{\tilde{M}}^{\tilde{G}}({\bf f}).$$
Puisque  $\boldsymbol{\delta}$ est stable, on a $S^{\tilde{M}}(\boldsymbol{\delta},p_{\tilde{M}}^{st}(\boldsymbol{\varphi}))=I^{\tilde{M}}(\boldsymbol{\delta},\boldsymbol{\varphi})$ pour tout $\boldsymbol{\varphi}\in I(\tilde{M}({\mathbb R}),K^M)\otimes Mes(M({\mathbb R}))$. L'\'egalit\'e \`a prouver se transforme en
$$^cI_{\tilde{M}}^{\tilde{G}}(\boldsymbol{\delta},{\bf f})=\sum_{\tilde{L}\in {\cal L}(\tilde{M})}I_{\tilde{M}}^{\tilde{L}}(\boldsymbol{\delta},{^c\theta}_{\tilde{L}}^{\tilde{G}}({\bf f})),$$
ce qui est la relation 5.13(5). Cela conclut le cas ${\bf M}'={\bf M}$. Supposons maintenant ${\bf M}'\not={\bf M}$. 
Le m\^eme calcul qu'en [VIII] 4.1 conduit \`a l'\'egalit\'e
$$(1) \qquad ^cI_{\tilde{M}}^{\tilde{G},{\cal E}}({\bf M}',\boldsymbol{\delta},{\bf f})=\sum_{\tilde{L}\in {\cal L}(\tilde{M})}\sum_{s\in \zeta Z(\hat{M}^{\Gamma_{{\mathbb R}}}/Z(\hat{G})^{\Gamma_{{\mathbb R}}}}
i_{\tilde{M}'}(\tilde{L},\tilde{L}'(s))S_{{\bf M}'}^{{\bf L}'(s)}(\boldsymbol{\delta},{^c\theta}_{\tilde{L}}^{\tilde{G},{\cal E}}({\bf L}'(s),{\bf f})).$$
Pour $\tilde{L}=\tilde{M}$, la somme en $s$ est r\'eduite \`a l'\'el\'ement $s=\zeta$ et le terme dans la somme est le premier terme du membre de droite de l'\'egalit\'e de l'\'enonc\'e. 
Pour $\tilde{L}\not=\tilde{M}$, on peut utiliser par r\'ecurrence l'assertion 6.11(A): ${^c\theta}_{\tilde{L}}^{\tilde{G},{\cal E}}({\bf L}'(s),{\bf f})=(^c\theta_{\tilde{L}}^{\tilde{G}}({\bf f}))^{{\bf L}'(s)}$. Le terme index\'e par $\tilde{L}$ dans l'expression ci-dessus devient alors $I_{\tilde{M}}^{\tilde{L},{\cal E}}({\bf M}',\boldsymbol{\delta},^c\theta_{\tilde{L}}^{\tilde{G}}({\bf f}))$.  C'est le terme index\'e par $\tilde{L}$ dans le membre de droite de l'\'egalit\'e de l'\'enonc\'e. L'\'egalit\'e (1) devient donc celle de l'\'enonc\'e.

 Le (ii) se prouve comme dans le cas non-archim\'edien en appliquant au cas ${\bf M}'={\bf M}$  la preuve que l'on a donn\'ee du (iii) dans le  cas ${\bf M}'\not={\bf M}$. On renvoie le lecteur \`a la preuve  de [VIII] proposition 4.1(ii).  $\square$

\bigskip

\subsection{Relation entre les applications $\theta_{K\tilde{M}}^{rat,K\tilde{G},{\cal E}}$, $^c\theta_{K\tilde{M}}^{rat,K\tilde{G},{\cal E}}$, $^c\theta_{K\tilde{M}}^{K\tilde{G},{\cal E}}$}
On se r\'ef\`ere aux cas (A), (B), (C) du paragraphe pr\'ec\'edent. 
\ass{Lemme}{(i) Dans le cas (A), soit ${\bf f}\in I(K\tilde{G}({\mathbb R}),\omega,K)\otimes Mes(G({\mathbb R}))$. On a l'\'egalit\'e
$$^c\theta_{K\tilde{M}}^{K\tilde{G},{\cal E}}({\bf f})=\sum_{K\tilde{L}\in {\cal L}(K\tilde{M})}\theta_{K\tilde{M}}^{rat,K\tilde{L},{\cal E}}\circ {^c\theta}_{K\tilde{L}}^{rat,K\tilde{G},{\cal E}}({\bf f}).$$

(ii) Dans le cas (B), soit ${\bf f}\in I(\tilde{G}({\mathbb R}),K)\otimes Mes(G({\mathbb R}))$. On a l'\'egalit\'e
$${^cS}\theta_{\tilde{M}}^{\tilde{G}}({\bf f})=S\theta_{\tilde{M}}^{rat,\tilde{G}}({\bf f})+ \sum_{\tilde{L}\in {\cal L}(\tilde{M}),\tilde{L}\not=\tilde{G}}S\theta_{\tilde{M}}^{rat,\tilde{L}}\circ {^cS\theta}_{\tilde{L}}^{rat,\tilde{G}}({\bf f}).$$

(iii) Dans le cas (C), soit ${\bf f}\in I(\tilde{G}({\mathbb R}),K)\otimes Mes(G({\mathbb R}))$. On a l'\'egalit\'e
$$^c\theta_{\tilde{M}}^{\tilde{G},{\cal E}}({\bf M}',{\bf f})={^c\theta}_{\tilde{M}}^{rat,\tilde{G},{\cal E}}({\bf M}',{\bf f})+ \sum_{\tilde{L}\in {\cal L}(\tilde{M}),\tilde{L}\not=\tilde{M}}
\theta_{\tilde{M}}^{rat,\tilde{L},{\cal E}}({\bf M}',{^c\theta}_{\tilde{L}}^{rat,\tilde{G}}({\bf f})).$$}

{\bf Remarque.} Quand on aura prouv\'e la stabilit\'e de l'application $S\theta_{\tilde{M}}^{rat,\tilde{G}}$, le premier terme du membre de droite de (ii) rentrera dans la somme qui le suit, o\`u la restriction $\tilde{L}\not=\tilde{G}$ dispara\^{\i}tra.
\bigskip

 Preuve. Consid\'erons (i). Il suffit de prouver que, pour tout ${\bf M}'\in {\cal E}(K\tilde{M},{\bf a})$, les transferts \`a ${\bf M}'$  des deux membres de l'\'egalit\'e sont \'egaux. En appliquant les d\'efinitions, cela revient \`a prouver l'\'egalit\'e
 $$(1) \qquad ^c\theta_{K\tilde{M}}^{K\tilde{G},{\cal E}}({\bf M}',{\bf f})=\sum_{K\tilde{L}\in {\cal L}(K\tilde{M})}\theta_{K\tilde{M}}^{rat,K\tilde{L},{\cal E}}({\bf M}',{^c\theta}_{K\tilde{L}}^{rat,K\tilde{G},{\cal E}}({\bf f})).$$
 Notons ${\bf M}'=(M',{\cal M}',\tilde{\zeta})$. Utilisons la d\'efinition du membre de gauche
 $$^c\theta_{K\tilde{M}}^{K\tilde{G},{\cal E}}({\bf M}',{\bf f})=\sum_{\tilde{s}\in \tilde{\zeta}Z(\hat{M})^{\Gamma_{{\mathbb R}},\hat{\theta}}/Z(\hat{G})^{\Gamma_{{\mathbb R}},\hat{\theta}}}i_{\tilde{M}'}(\tilde{G},\tilde{G}'(\tilde{s})){^cS\theta}_{{\bf M}'}^{{\bf G}'(\tilde{s})}({\bf f}^{{\bf G}'(\tilde{s})}).$$
 Nos hypoth\`eses de r\'ecurrence et quelques formalit\'es nous autorisent \`a appliquer  aux termes du membre de droite le (ii) du pr\'esent \'enonc\'e, simplifi\'e par la remarque qui suit cet \'enonc\'e.  Le terme index\'e par $\tilde{s}$ se d\'eveloppe en une somme index\'ee par $\tilde{L}'_{\tilde{s}}\in {\cal L}^{\tilde{G}'(\tilde{s})}(\tilde{M}')$. Comme toujours, un tel $\tilde{L}'_{\tilde{s}}$ d\'etermine un $K$-espace de Levi $K\tilde{L}$. Il s'identifie alors \`a l'espace de la donn\'ee endoscopique ${\bf L}'(\tilde{s})$ de $(KL,K\tilde{L},{\bf a})$. On regroupe les couples $(\tilde{s},\tilde{L}_{\tilde{s}})$ selon le $K$-espace de Levi $K\tilde{L}$. On obtient
  $$(2)\qquad ^c\theta_{K\tilde{M}}^{K\tilde{G},{\cal E}}({\bf M}',{\bf f})=\sum_{K\tilde{L}\in {\cal L}(K\tilde{M})}\sum_{\tilde{s}\in \tilde{\zeta}Z(\hat{M})^{\Gamma_{{\mathbb R}},\hat{\theta}}/Z(\hat{L})^{\Gamma_{{\mathbb R}},\hat{\theta}},{\bf L}'(\tilde{s})\text{ elliptique }}$$
  $$\sum_{\tilde{t}\in \tilde{s}Z(\hat{L})^{\Gamma_{{\mathbb R}},\hat{\theta}}/Z(\hat{G})^{\Gamma_{{\mathbb R}},\hat{\theta}}}i_{\tilde{M}'}(\tilde{G},\tilde{G}'(\tilde{t}))
S\theta_{{\bf M}'}^{rat,{\bf L}'(\tilde{s})}\circ {^cS}\theta_{{\bf L}'(\tilde{s})}^{rat,{\bf G}'(\tilde{t})}({\bf f}^{{\bf G}'(\tilde{t})}).$$
On a comme toujours l'\'egalit\'e
$$i_{\tilde{M}'}(\tilde{G},\tilde{G}'(\tilde{t}))=i_{\tilde{M}'}(\tilde{L},\tilde{L}'(\tilde{s}))i_{\tilde{L}'(\tilde{s})}(\tilde{G},\tilde{G}'(\tilde{t})).$$
Sous cette forme, la restriction que ${\bf L}'(\tilde{s})$ doit \^etre elliptique devient superflue. Pour tout $\tilde{s}$, on a par d\'efinition
$$\sum_{\tilde{t}\in \tilde{s}Z(\hat{L})^{\Gamma_{{\mathbb R}},\hat{\theta}}/Z(\hat{G})^{\Gamma_{{\mathbb R}},\hat{\theta}}}i_{\tilde{L}'(\tilde{s})}(\tilde{G},\tilde{G}'(\tilde{t})){^cS}\theta_{{\bf L}'(\tilde{s})}^{rat,{\bf G}'(\tilde{t})}({\bf f}^{{\bf G}'(\tilde{t})})={^c\theta}_{K\tilde{L}}^{rat,K\tilde{G},{\cal E}}({\bf L}'(\tilde{s}),{\bf f})= {^c\theta}_{K\tilde{L}}^{rat,K\tilde{G},{\cal E}}({\bf f})^{{\bf L}'(\tilde{s})}.$$
Le terme index\'e par $K\tilde{L}$ de l'expression (1) devient
$$\sum_{\tilde{s}\in \tilde{\zeta}Z(\hat{M})^{\Gamma_{{\mathbb R}},\hat{\theta}}/Z(\hat{L})^{\Gamma_{{\mathbb R}},\hat{\theta}}}i_{\tilde{M}'}(\tilde{L},\tilde{L}'(\tilde{s}))S\theta_{{\bf M}'}^{rat,{\bf L}'(\tilde{s})}({^c\theta}_{K\tilde{L}}^{rat,K\tilde{G},{\cal E}}({\bf f})^{{\bf L}'(\tilde{s})}).$$
Par d\'efinition, c'est $\theta_{K\tilde{M}}^{rat,K\tilde{L},{\cal E}}({\bf M}',{^c\theta}_{K\tilde{L}}^{rat,K\tilde{G},{\cal E}}({\bf f}))$. Mais alors le membre de droite de (2) devient \'egal \`a celui de (1). Cela d\'emontre cette \'egalit\'e (1) et l'assertion (i) de l'\'enonc\'e. 

Traitons (iii). Supposons d'abord ${\bf M}'={\bf M}$. Comme on l'a dit, en notant $\theta_{\tilde{M}}^{\tilde{G}}$ l'une de nos application ${^c\theta}_{\tilde{M}}^{rat,\tilde{G}}$ etc..., l'application  $S\theta_{\tilde{M}}^{\tilde{G}}$ est d\'efinie de sorte que l'on ait l'\'egalit\'e
$$\theta_{\tilde{M}}^{\tilde{G},{\cal E}}({\bf M},{\bf f})=p_{\tilde{M}}^{st}\circ{\theta}_{\tilde{M}}^{\tilde{G}}({\bf f}).$$
Alors l'\'egalit\'e du (iii) n'est autre que l'\'egalit\'e du lemme 5.12 \`a laquelle on applique la projection $p_{\tilde{M}}^{st}$. Supposons maintenant ${\bf M}'\not={\bf M}$. On reprend la preuve du (i) en supprimant les $K$ et en rempla\c{c}ant $\tilde{s}$ par $s$. Le calcul marche jusqu'au point o\`u on a utilis\'e l'\'egalit\'e $^c\theta_{\tilde{L}}^{rat,\tilde{G},{\cal E}}({\bf L}'(s),{\bf f})=(^c\theta_{\tilde{L}}^{rat,\tilde{G},{\cal E}}({\bf f}))^{{\bf L}'(s)}$. Ici, on ne dispose plus de l'application $^c\theta_{\tilde{L}}^{rat,\tilde{G},{\cal E}}$. Mais, si $\tilde{L}\not=\tilde{M}$, on peut utiliser l'assertion 6.11(A): on a
 $^c\theta_{\tilde{L}}^{rat,\tilde{G}}({\bf L}'(s),{\bf f})=(^c\theta_{\tilde{L}}^{rat,\tilde{G}}({\bf f}))^{{\bf L}'(s)}$. Le calcul des termes index\'es par un $\tilde{L}\not=\tilde{M}$ se poursuit comme ci-dessus et donne le terme voulu du membre de droite de l'\'egalit\'e du (iii). Il reste le terme index\'e par $\tilde{M}$. C'est simplement $S\theta_{{\bf M}'}^{rat,{\bf M}'}(^c\theta_{\tilde{M}}^{rat,\tilde{G},{\cal E}}({\bf M}',{\bf f}))$. Mais $S\theta_{{\bf M}'}^{rat,{\bf M}'}$ est simplement l'identit\'e de l'espace $SI_{ac}({\bf M}')\otimes Mes(M'({\mathbb R}))$. Le terme pr\'ec\'edent est donc $^c\theta_{\tilde{M}}^{rat,\tilde{G},{\cal E}}({\bf M}',{\bf f})$, qui est le premier terme du membre de droite de l'\'egalit\'e du (iii). Cela prouve (iii).
 
 Le (ii) se prouve comme toujours en appliquant au cas ${\bf M}'={\bf M}$  la preuve que l'on a donn\'ee du (iii) dans le  cas ${\bf M}'\not={\bf M}$. On renvoie le lecteur \`a la preuve similaire de [VIII] proposition 4.1(ii).  $\square$

\subsection{Preuves des  propositions 6.1, 6.5 et du lemme 6.4}
On suppose $(G,\tilde{G},{\bf a})$ quasi-d\'eploy\'e et \`a torsion int\'erieure. Comme toujours, les propri\'et\'es \`a d\'emontrer sont tautologiques si $\tilde{M}=\tilde{G}$. On suppose $\tilde{M}\not=\tilde{G}$. 

Soit ${\bf f}\in I(\tilde{G}({\mathbb R}),K)\otimes Mes(G({\mathbb R}))$ dont l'image dans $SI(\tilde{G}({\mathbb R}),K)\otimes Mes(G({\mathbb R}))$ est nulle. Soit $\boldsymbol{\delta}\in D^{st}_{orb,\tilde{G}-reg}(\tilde{M}({\mathbb R}))\otimes Mes(M({\mathbb R}))^*$. Pour $\tilde{L}\in {\cal L}(\tilde{M})$, $\tilde{L}\not=\tilde{M}$, nos hypoth\`eses de r\'ecurrence assurent que $^cS\theta_{\tilde{L}}^{\tilde{G}}({\bf f})=0$. Le lemme 7.1(ii) se simplifie en
$$(1) \qquad {^cS}\theta_{\tilde{M}}^{rat,\tilde{G}}(\boldsymbol{\delta},{\bf f})=S^{\tilde{M}}(\boldsymbol{\delta},{^cS}\theta_{\tilde{M}}^{\tilde{G}}({\bf f})).$$
On a aussi $^cS\theta_{\tilde{L}}^{rat,\tilde{G}}({\bf f})=0$ pour $\tilde{L}\not=\tilde{G}$ et le lemme 7.2(ii) se simplifie en
$$(2) \qquad {^cS}\theta_{\tilde{M}}^{\tilde{G}}({\bf f})={^cS}\theta_{\tilde{M}}^{rat,\tilde{G}}({\bf f})+S\theta_{\tilde{M}}^{rat,\tilde{G}}({\bf f}).$$
Le membre de gauche de (1) poss\`ede la propri\'et\'e de compacit\'e 5.13(3). Il en r\'esulte que ${^cS}\theta_{\tilde{M}}^{\tilde{G}}({\bf f})$, qui est a priori un \'el\'ement de $SI_{ac}(\tilde{M}({\mathbb R}),K^M)\otimes Mes(M({\mathbb R}))$, est en fait "\`a support compact", c'est-\`a-dire appartient \`a $SI(\tilde{M}({\mathbb R}),K^M)\otimes Mes(M({\mathbb R}))$. Soit $\tilde{\pi}\in D_{ell,0}^{st}(\tilde{M}({\mathbb R}))\otimes Mes(M({\mathbb R}))^*$. Puisque  ${^cS}\theta_{\tilde{M}}^{\tilde{G}}({\bf f})$ est \`a support compact, la fonction $\lambda\mapsto S^{\tilde{M}}(\tilde{\pi}_{\lambda}, {^cS}\theta_{\tilde{M}}^{\tilde{G}}({\bf f}))$ est holomorphe pour $\lambda\in {\cal A}_{\tilde{M},{\mathbb C}}^*$.  Sa transform\'ee de Fourier sur $i{\cal A}_{\tilde{M}}^*$ est la fonction $X\mapsto S^{\tilde{M}}(\tilde{\pi},X, {^cS}\theta_{\tilde{M}}^{\tilde{G}}({\bf f}))$ sur ${\cal A}_{\tilde{M}}$. D'apr\`es (2), c'est la fonction
$$X\mapsto S^{\tilde{M}}(\tilde{\pi},X, {^cS}\theta_{\tilde{M}}^{rat,\tilde{G}}({\bf f}))+S^{\tilde{M}}(\tilde{\pi},X, S\theta_{\tilde{M}}^{rat,\tilde{G}}({\bf f})).$$
Or la premi\`ere fonction  est la transform\'ee de Fourier de la fonction $\lambda\mapsto S^{\tilde{M}}(\tilde{\pi},\lambda,{^cS}\theta_{\tilde{M}}^{rat, \tilde{G}}({\bf f}))$. La seconde  est la transform\'ee de Fourier de $\lambda\mapsto I^{\tilde{M}}( \sigma_{\tilde{M}}^{\tilde{G}}(\tilde{\pi},\lambda),{\bf f}_{\tilde{M}})$ d'apr\`es le lemme 6.2. Il en r\'esulte que
$$(3) \qquad S^{\tilde{M}}(\tilde{\pi}_{\lambda}, {^cS}\theta_{\tilde{M}}^{\tilde{G}}({\bf f}))=S^{\tilde{M}}(\tilde{\pi},\lambda,{^cS}\theta_{\tilde{M}}^{rat,\tilde{G}}({\bf f}))+I^{\tilde{M}}(\sigma_{\tilde{M}}^{\tilde{G}}(\tilde{\pi},\lambda),{\bf f}_{\tilde{M}}).$$
La fonction de gauche est holomorphe. La premi\`ere fonction de droite n'a qu'un nombre fini d'hyperplans polaires. Donc $\lambda\mapsto I^{\tilde{M}}(\sigma_{\tilde{M}}^{\tilde{G}}(\tilde{\pi},\lambda),{\bf f})$ n'a qu'un nombre fini d'hyperplans polaires.  On peut pr\'eciser que, si on fixe un ensemble fini $\Omega$ de $K$-types, ces hyperplans restent dans un ensemble fini ind\'ependant de ${\bf f}$ pourvu que ${\bf f}\in I(\tilde{G}({\mathbb R}),\Omega)\otimes Mes(G({\mathbb R}))$. On sait d\'ecomposer $D_{ell,0}(\tilde{M}({\mathbb R}))$ en somme de sa partie stable $D_{ell,0}^{st}(\tilde{M}({\mathbb R}))$ et de sa partie instable $D_{ell,0}^{inst}(\tilde{M}({\mathbb R}))$, cf. [IV] 2.2. On peut d\'ecomposer $\sigma_{\tilde{M}}^{\tilde{G}}(\tilde{\pi})$ sous la forme
$$\sigma_{\tilde{M}}^{\tilde{G}}(\tilde{\pi})=\sum_{j=1,...,k}u_{j}\otimes \tilde{\pi}_{j}$$
de sorte que les propri\'et\'es suivantes soient v\'erifi\'ees:

- $u_{j}\in U_{\tilde{M}}^{\tilde{G}}$ et $u_{j}\not=0$ pour tout $j$;

- $\tilde{\pi}_{j}\in D_{ell,0}(\tilde{M}({\mathbb R}))\otimes Mes(M({\mathbb R}))^*$ et la famille $(\tilde{\pi}_{j})_{j=1,...,k}$ est lin\'eairement ind\'ependante;

- il existe un entier $l\in \{0,...,k\}$ de sorte que $\tilde{\pi}_{j}\in D_{ell,0}^{inst}(\tilde{M}({\mathbb R}))\otimes Mes(M({\mathbb R}))^*$ si $j\leq l$, tandis que $\tilde{\pi}_{j}\in D_{ell,0}^{st}(\tilde{M}({\mathbb R}))\otimes Mes(M({\mathbb R}))^*$ si $j>l$.

L'hypoth\`ese sur ${\bf f}$ implique que l'image de ${\bf f}_{\tilde{M}}$ dans $SI(\tilde{M}({\mathbb R}))\otimes Mes(M({\mathbb R}))$ est nulle. Donc $I^{\tilde{M}}(\tilde{\pi}_{j,\lambda},{\bf f}_{\tilde{M}})=0$ pour tout $j>l$ et tout $\lambda$.  Supposons $l>0$. Appliquons le lemme 5.6 ter pour les $\tilde{\pi}_{1}$,...,$\tilde{\pi}_{l}$. On obtient un ensemble fini $\Omega$ de $K$-types et une r\'eunion finie de sous-espaces affines ${\cal H}$ v\'erifiant la conclusion de ce lemme. Comme on l'a dit, pour cet ensemble $\Omega$, les hyperplans polaires de la fonction $\lambda\mapsto  I^{\tilde{M}}(\sigma_{\tilde{M}}^{\tilde{G}}(\tilde{\pi},\lambda),{\bf f})$ restent dans un ensemble fini pourvu que 
${\bf f}\in I(\tilde{G}({\mathbb R}),\Omega)\otimes Mes(G({\mathbb R}))$. Quitte \`a accro\^{\i}tre ${\cal H}$, on peut supposer que ${\cal H}$ contient tous ces hyperplans. Notons ${\cal F}$ l'espace des fonctions holomorphes sur ${\cal A}_{\tilde{M},{\mathbb C}}^*$. Notons $I^{inst}(\tilde{G}({\mathbb R}),\Omega)$ le sous-espace des \'el\'ements ${\bf f}\in I(\tilde{G}({\mathbb R}),\Omega)$ d'image nulle dans $SI(\tilde{G}({\mathbb R}))$. A tout \'el\'ement  ${\bf f}\in I^{inst}(\tilde{G}({\mathbb R}),\Omega)\otimes Mes(G({\mathbb R}))$, associons la famille $\underline{e}_{{\bf f}}=(e_{{\bf f},j})_{j=1,...,l}\in {\cal F}^l$ d\'efinie par $e_{{\bf f},j}(\lambda)=I^{\tilde{M}}(\tilde{\pi}_{j,\lambda},{\bf f}_{\tilde{M}})$. Puisque
$$I^{\tilde{M}}(\sigma_{\tilde{M}}^{\tilde{G}}(\tilde{\pi},\lambda),{\bf f})=\sum_{j=1,...,l}u_{j}(\lambda)I^{\tilde{M}}(\tilde{\pi}_{j,\lambda},{\bf f}_{\tilde{M}}),$$
on obtient

(4) pour tout ${\bf f}\in  I^{inst}(\tilde{G}({\mathbb R}),\Omega)\otimes Mes(G({\mathbb R}))$, la fonction
$\lambda\mapsto \sum_{j=1,...,l}u_{j}(\lambda)e_{{\bf f},j}(\lambda)$ n'a pas de p\^ole hors de ${\cal H}$. 

On a aussi

(5) pour tout $\lambda\in {\cal A}_{\tilde{M},{\mathbb C}}^*-{\cal H}$, il existe ${\bf f}\in I^{inst}(\tilde{G}({\mathbb R}),\Omega)\otimes Mes(G({\mathbb R}))$ tel que $e_{{\bf f},1}(\lambda)=1$ et $e_{{\bf f},j}(\lambda)=0$ pour $j=2,...,l$.

En effet, on peut choisir $\boldsymbol{\phi}\in I(\tilde{M}({\mathbb R}))\otimes Mes(M({\mathbb R}))$ tel que $I^{\tilde{M}}(\tilde{\pi}_{1,\lambda},\boldsymbol{\phi})=1$ tandis que $I^{\tilde{M}}(\tilde{\pi}_{j,\lambda},\boldsymbol{\phi})=0$ pour $j=2,...,l$. Le lemme 5.6 ter associe \`a $\boldsymbol{\phi}$ une fonction ${\bf f}$ v\'erifiant (5).

Les propri\'et\'es (4) et (5) et la proposition 5.15 entra\^{\i}nent que $u_{1}=0$. Cela contredit l'hypoth\`ese.  Cela d\'emontre que $l=0$, autrement dit, $\sigma_{\tilde{M}}^{\tilde{G}}(\tilde{\pi})$ appartient \`a $U_{\tilde{M}}^{\tilde{G}}\otimes D_{ell,0}^{st}(\tilde{M}({\mathbb R}))\otimes Mes(M({\mathbb R}))^*$. Il en r\'esulte que $I^{\tilde{M}}(\sigma_{\tilde{M}}^{\tilde{G}}(\tilde{\pi},\lambda),{\bf f}_{\tilde{M}})=0$ pour tout ${\bf f}$ comme au d\'ebut de la preuve. En remontant le calcul, cela entra\^{\i}ne $S^{\tilde{M}}(\tilde{\pi},X,S\theta_{\tilde{M}}^{rat,\tilde{G}}({\bf f}))=0$ pour tout $X$. Mais la fonction  $S\theta_{\tilde{M}}^{rat,\tilde{G}}({\bf f})$ est cuspidale: il r\'esulte de la formule de descente [VIII] 2.3 et de nos hypoth\`eses de r\'ecurrence que $S\theta_{\tilde{M}}^{rat,\tilde{G}}({\bf f})_{\tilde{R}}=0$ pour tout espace de Levi $\tilde{R}\subsetneq \tilde{M}$. Une fonction cuspidale qui v\'erifie la nullit\'e ci-dessus pour tout $\tilde{\pi}$ elliptique est nulle. Cela prouve que $S\theta_{\tilde{M}}^{rat,\tilde{G}}({\bf f})=0$.

Revenons \`a la situation de la  formule (3). On vient de montrer que le deuxi\`eme terme du membre de droite \'etait nul.  Donc la fonction $\lambda\mapsto S^{\tilde{M}}(\tilde{\pi},\lambda,{^cS}\theta_{\tilde{M}}^{rat,\tilde{G}}({\bf f}))$ est holomorphe. La transform\'ee de Fourier $X\mapsto S^{\tilde{M}}(\tilde{\pi},\nu,X,{^cS}\theta_{\tilde{M}}^{rat,\tilde{G}}({\bf f}))$ ne d\'epend pas du point $\nu\in {\cal A}_{\tilde{M}}^*$. Dans la proposition 6.2, on peut remplacer tous les points $\nu_{\tilde{S}}$ par $0$. Puisque $\sum_{\tilde{S}\in {\cal P}(\tilde{M})}\omega_{\tilde{S}}(X)=1$ pour tout $X$, on obtient que $S^{\tilde{M}}(\tilde{\pi},X,{^cS}\theta_{\tilde{M}}^{rat,\tilde{G}}({\bf f}))=0$ pour tout $X$.  Pour la m\^eme raison que ci-dessus, la fonction  ${^cS}\theta_{\tilde{M}}^{rat,\tilde{G}}({\bf f})$ est cuspidale et la nullit\'e ci-dessus pour tout $\tilde{\pi}$ elliptique entra\^{\i}ne ${^cS}\theta_{\tilde{M}}^{rat,\tilde{G}}({\bf f})=0$. L'\'egalit\'e (2) entra\^{\i}ne alors ${^cS}\theta_{\tilde{M}}^{\tilde{G}}({\bf f})=0$. Cela d\'emontre la proposition 6.1. La proposition 6.5 en r\'esulte gr\^ace \`a l'\'egalit\'e (1). 

On vient de  prouver le lemme 6.4 restreint aux repr\'esentations elliptiques.  Pour $\tilde{\pi}\in D_{temp,0}^{st}(\tilde{M}({\mathbb R}))\otimes Mes(M({\mathbb R}))^*$, la formule du lemme 6.3 et la proposition 6.1 maintenant d\'emontr\'ee entra\^{\i}nent que $I^{\tilde{M}}(\sigma_{\tilde{M}}^{\tilde{G}}(\tilde{\pi};\lambda),{\bf f}_{\tilde{M}})=0$ pour tout ${\bf f}\in I(\tilde{G}({\mathbb R}),K)\otimes Mes(G({\mathbb R}))$ d'image nulle dans $SI(\tilde{G}({\mathbb R}))\otimes Mes(G({\mathbb R}))$. L'usage du lemme 5.7 ter et de la proposition 5.15 permet d'en d\'eduire comme ci-dessus que $\sigma_{\tilde{M}}^{\tilde{G}}(\tilde{\pi})$ appartient \`a $U_{\tilde{M}}^{\tilde{G}}\otimes D_{ell,0}^{st}(\tilde{M}({\mathbb R}))\otimes Mes(M({\mathbb R}))^*$. $\square$

\bigskip

\subsection{Preuve conditionnelle des propositions 6.7 et 6.10 et du lemme 6.9}
On consid\`ere un triplet $(KG,K\tilde{G},{\bf a})$ et un $K$-espace de Levi $K\tilde{M}\in {\cal L}(K\tilde{M}_{0})$. Les propri\'et\'es \`a prouver sont tautologiques si $K\tilde{M}=K\tilde{G}$. On suppose $K\tilde{M}\not=K\tilde{G}$. On impose la condition suivante:

 (Hyp) soient ${\bf f}\in I(K\tilde{G}({\mathbb R}),\omega,K)\otimes Mes(G({\mathbb R}))$ et $\boldsymbol{\gamma}\in D_{orb,\tilde{G}-reg}(K\tilde{M}({\mathbb R}),\omega)\otimes Mes(M({\mathbb R}))^*$; alors on a l'\'egalit\'e
 $$I_{K\tilde{M}}^{K\tilde{G},{\cal E}}(\boldsymbol{\gamma},{\bf f})=I_{K\tilde{M}}^{K\tilde{G}}(\boldsymbol{\gamma},{\bf f}).$$
 
 Sous cette hypoth\`ese, nous allons prouver les propositions 6.7, 6.10 et le lemme 6.9.   Soient ${\bf f}\in I(K\tilde{G}({\mathbb R}),\omega,K)\otimes Mes(G({\mathbb R}))$ et $\boldsymbol{\gamma}\in D_{orb,\tilde{G}-reg}(K\tilde{M}({\mathbb R}),\omega)\otimes Mes(M({\mathbb R}))^*$.
 Posons
 $${^c\varphi}={^c\theta}_{K\tilde{M}}^{K\tilde{G},{\cal E}}({\bf f})-{^c\theta}_{K\tilde{M}}^{K\tilde{G}}({\bf f}),$$
 $${^c\varphi}^{rat}={^c\theta}_{K\tilde{M}}^{rat,K\tilde{G},{\cal E}}({\bf f})-{^c\theta}_{K\tilde{M}}^{rat,K\tilde{G}}({\bf f}),$$
 $$\varphi^{rat}=\theta_{K\tilde{M}}^{rat,K\tilde{G},{\cal E}}({\bf f})-\theta_{K\tilde{M}}^{rat,K\tilde{G}}({\bf f}).$$
 
 Faisons la diff\'erence entre la formule du lemme 7.1(i) et la formule 5.13(5), \'etendue de fa\c{c}on \'evidente \`a notre $K$-espace. On obtient
 $$^cI_{K\tilde{M}}^{K\tilde{G},{\cal E}}(\boldsymbol{\gamma},{\bf f})-{^c}I_{K\tilde{M}}^{K\tilde{G},{\cal E}}(\boldsymbol{\gamma},{\bf f})=\sum_{K\tilde{L}\in {\cal L}(K\tilde{M})}I_{K\tilde{M}}^{K\tilde{L},{\cal E}}(\boldsymbol{\gamma},{^c\theta}_{K\tilde{L}}^{K\tilde{G},{\cal E}}({\bf f}))-I_{K\tilde{M}}^{K\tilde{L}}(\boldsymbol{\gamma},{^c\theta}_{K\tilde{L}}^{K\tilde{G}}({\bf f})).$$
 Pour $K\tilde{L}\not=K\tilde{G}$, les hypoth\`eses de r\'ecurrence assurent que $I_{K\tilde{M}}^{K\tilde{L},{\cal E}}(\boldsymbol{\gamma},.)=I_{K\tilde{M}}^{K\tilde{L}}(\boldsymbol{\gamma},.)$. Si $K\tilde{L}=K\tilde{G}$, cette \'egalit\'e est assur\'ee par notre hypoth\`ese (Hyp). Si $K\tilde{L}\not=K\tilde{M}$, nos hypoth\`eses de r\'ecurrence assurent que ${^c\theta}_{K\tilde{L}}^{K\tilde{G},{\cal E}}({\bf f})={^c\theta}_{K\tilde{L}}^{K\tilde{G}}({\bf f})$. Il reste
 $$(1) \qquad ^cI_{K\tilde{M}}^{K\tilde{G},{\cal E}}(\boldsymbol{\gamma},{\bf f})-{^c}I_{K\tilde{M}}^{K\tilde{G},{\cal E}}(\boldsymbol{\gamma},{\bf f})=I^{K\tilde{M}}(\boldsymbol{\gamma}, {^c\varphi}).$$
 Faisons la diff\'erence entre la formule du lemme 7.2(ii) et celle du lemme 5.12. De la m\^eme fa\c{c}on, les hypoth\`eses de r\'ecurrence conduisent \`a l'\'egalit\'e
  $$(2) \qquad {^c\varphi}={^c\varphi}^{rat}+\varphi^{rat}.$$
La suite de la preuve est similaire \`a celle du paragraphe pr\'ec\'edent. Le membre de gauche de (1) poss\`ede la propri\'et\'e de compacit\'e 5.13(3). Il en r\'esulte que ${^c\varphi}$, qui est a priori un \'el\'ement de $I_{ac}(K\tilde{M}({\mathbb R}),\omega,K^M)\otimes Mes(M({\mathbb R}))$, est en fait \`a support compact, c'est-\`a-dire appartient \`a $I(K\tilde{M}({\mathbb R}),\omega,K^M)\otimes Mes(M({\mathbb R}))$. Fixons une composante connexe $\tilde{M}$ de $K\tilde{M}$. Soit $\tilde{\pi}\in D_{ell,0}(\tilde{M}({\mathbb R}),\omega)\otimes Mes(M({\mathbb R}))^*$. De la compacit\'e du support de ${^c\varphi}$ r\'esulte que la fonction $\lambda\mapsto I^{K\tilde{M}}(\tilde{\pi}_{\lambda},{^c\varphi})$ est holomorphe sur ${\cal A}_{\tilde{M},{\mathbb C}}^*$. On calcule cette fonction gr\^ace \`a (2) comme en 7.3. On a
$$(3) \qquad I^{K\tilde{M}}(\tilde{\pi}_{\lambda},{^c\varphi})=I^{K\tilde{M}}(\tilde{\pi},\lambda,{^c\varphi}^{rat})+I^{K\tilde{M}}(\rho_{K\tilde{M}}^{K\tilde{G},{\cal E}}(\tilde{\pi},\lambda)-\rho_{K\tilde{M}}^{K\tilde{G}}(\tilde{\pi},\lambda),{\bf f}_{K\tilde{M},\omega}).$$
La premi\`ere fonction de droite n'a qu'un nombre fini d'hyperplans polaires. Si on fixe un ensemble fini $\Omega$ de $K$-types, ces hyperplans restent dans un ensemble fini ind\'ependant de ${\bf f}$ pourvu que ${\bf f}\in I(K\tilde{G}({\mathbb R}),\omega,\Omega)\otimes Mes(G({\mathbb R}))$. Il en est donc de m\^eme de la deuxi\`eme fonction. On d\'ecompose $\rho_{K\tilde{M}}^{K\tilde{G},{\cal E}}(\tilde{\pi})-\rho_{K\tilde{M}}^{K\tilde{G}}(\tilde{\pi})$ en
$$\rho_{K\tilde{M}}^{K\tilde{G},{\cal E}}(\tilde{\pi})-\rho_{K\tilde{M}}^{K\tilde{G}}(\tilde{\pi})=\sum_{j=1,...,k}u_{j}\otimes \tilde{\pi}_{j},$$
o\`u

- $u_{j}\in U_{\tilde{M}}^{\tilde{G}}$ et $u_{j}\not=0$ pour tout $j$;

- $\tilde{\pi}_{j}\in D_{ell,0}(K\tilde{M}({\mathbb R}),\omega)\otimes Mes(M({\mathbb R}))^*$ pour tout $j$ et la famille $(\tilde{\pi}_{j})_{j=1,...,k}$ est lin\'eairement ind\'ependante.

Le lemme 5.6 bis se g\'en\'eralise imm\'ediatement \`a notre $K$-espace $K\tilde{G}$.  Appliqu\'e \`a la famille $(\tilde{\pi}_{j})_{j=1,...,k}$, il nous fournit un ensemble fini $\Omega$ de $K$-types et une r\'eunion finie ${\cal H}$ de sous-espaces affines de ${\cal A}_{\tilde{M},{\mathbb C}}^*$.  Quitte \`a accro\^{\i}tre cet ensemble, on peut supposer que les fonctions du membre de droite de (3) n'ont pas de p\^ole hors de ${\cal H}$. A tout ${\bf f}\in I(K\tilde{G}({\mathbb R}),\omega,\Omega)\otimes Mes(G({\mathbb R}))$, on associe la famille $\underline{e}_{{\bf f}}=(e_{{\bf f},j})_{j=1,...,k}$ de fonctions holomorphes sur ${\cal A}_{\tilde{M},{\mathbb C}}^*$ d\'efinie par
$e_{{\bf f},j}(\lambda)=I^{K\tilde{M}}(\tilde{\pi}_{j},{\bf f}_{\tilde{M},\omega})$.  Le lemme 5.6 assure que, pour tout $\lambda\in {\cal A}_{\tilde{M},{\mathbb C}}^*-{\cal H}$, on peut trouver ${\bf f}\in I(K\tilde{G}({\mathbb R}),\omega,\Omega)\otimes Mes(G({\mathbb R}))$ de sorte que $e_{{\bf f},1}(\lambda)=1$ tandis que $e_{{\bf f},j}(\lambda)=0$ si $j=2,...,k$. Puisque
$$I^{K\tilde{M}}(\rho_{K\tilde{M}}^{K\tilde{G},{\cal E}}(\tilde{\pi},\lambda)-\rho_{K\tilde{M}}^{K\tilde{G}}(\tilde{\pi},\lambda),{\bf f}_{K\tilde{M},\omega})=\sum_{j=1,...,k}u_{j}(\lambda)e_{{\bf f},j}(\lambda),$$
cette derni\`ere somme est holomorphe hors de ${\cal H}$ pour tout ${\bf f}\in  I(K\tilde{G}({\mathbb R}),\omega,\Omega)\otimes Mes(G({\mathbb R}))$. La proposition 5.15 assure alors que $u_{1}=0$, ce qui est contradictoire avec nos hypoth\`eses sauf si $k=0$. On a donc $k=0$, ce qui \'equivaut \`a l'\'egalit\'e $\rho_{K\tilde{M}}^{K\tilde{G},{\cal E}}(\tilde{\pi})=\rho_{K\tilde{M}}^{K\tilde{G}}(\tilde{\pi})$. En remontant le calcul, cela assure que $I^{K\tilde{M}}(\tilde{\pi},X,\varphi^{rat})=0$. Les hypoth\`eses de r\'ecurrence et les formules habituelles de descente impliquent que $\varphi^{rat}$ est cuspidale. La nullit\'e pr\'ec\'edente pour toute composante connexe $\tilde{M}$ et tout $\tilde{\pi}$ elliptique impliquent alors $\varphi^{rat}=0$. 

Revenons \`a la situation de la formule (3). Le deuxi\`eme terme du membre de droite est maintenant nul. Donc le premier est holomorphe puisque le membre de gauche l'est. Il en r\'esulte que la transform\'ee de Fourier $X\mapsto I^{K\tilde{M}}(\tilde{\pi},\nu,X,{^c\varphi}^{rat})$ ne d\'epend pas du point $\nu$ choisi. Par diff\'erence entre la formule 6.6(1) et celle de la proposition 5.10, on a
$$\sum_{\tilde{S}\in {\cal P}(\tilde{M})}\omega_{\tilde{S}}(X)I^{K\tilde{M}}(\tilde{\pi},\nu_{\tilde{S}},X,{^c\varphi}^{rat})=0.$$
On peut remplacer chaque $\nu_{\tilde{S}}$ par $0$ et, par un raisonnement d\'ej\`a fait plusieurs fois, on en d\'eduit que $I^{K\tilde{M}}(\tilde{\pi},X,{^c\varphi}^{rat})=0$ pour tout $X$. Pour la m\^eme raison que ci-dessus, la fonction ${^c\varphi}^{rat}$ est cuspidale et la nullit\'e pr\'ec\'edente pour toute composante $\tilde{M}$ et tout $\tilde{\pi}$ elliptique implique que la fonction ${^c\varphi}^{rat}$ est nulle. Enfin l'\'egalit\'e (2) entra\^{\i}ne que ${^c\varphi}=0$. Cela prouve la proposition 6.7.
La proposition 6.10 r\'esulte maintenant de (1).

On a prouv\'e ci-dessus l'\'egalit\'e du lemme 6.9 pour les $\tilde{\pi}$ elliptiques. Pour 

\noindent $\tilde{\pi}\in D_{temp,0}(K\tilde{M}({\mathbb R}),\omega)\otimes Mes(M({\mathbb R}))^*$, les formules des lemmes 5.7 et 6.8 et l'\'egalit\'e maintenant prouv\'ee $\theta_{K\tilde{M}}^{rat,\tilde{G},{\cal E}}=\theta_{K\tilde{M}}^{rat,\tilde{G}}$ entra\^{\i}nent que
$$I^{K\tilde{M}}(\rho_{K\tilde{M}}^{K\tilde{G},{\cal E}}(\tilde{\pi},\lambda)-\rho_{K\tilde{M}}^{K\tilde{G}}(\tilde{\pi},\lambda),{\bf f})=0$$
pour tout $\lambda\in {\cal A}_{\tilde{M},{\mathbb C}}^*$ et tout ${\bf f}\in I(K\tilde{G}({\mathbb R}),\omega,K)\otimes Mes(G({\mathbb R}))$. Comme ci-dessus, l'usage du lemme 5.6 et de la proposition 5.15 d\'emontrent l'\'egalit\'e  $\rho_{K\tilde{M}}^{K\tilde{G},{\cal E}}(\tilde{\pi})=\rho_{K\tilde{M}}^{K\tilde{G}}(\tilde{\pi})$.  $\square$

\bigskip

\subsection{Variante dans le cas quasi-d\'eploy\'e et \`a torsion int\'erieure}
On consid\`ere un triplet $(G,\tilde{G},{\bf a})$ quasi-d\'eploy\'e et \`a torsion int\'erieure et un espace de Levi $\tilde{M}$. Soit ${\bf M}'=(M',{\cal M}',\zeta)$ une donn\'ee endoscopique elliptique et relevante de $(M,\tilde{M})$. On va prouver les assertions (A), (B) et (C) de 6.11. Soit ${\bf f}\in I(\tilde{G}({\mathbb R}),K)\otimes Mes(G({\mathbb R}))$.  Comme toujours, ces assertions sont tautologiques si $\tilde{M}=\tilde{G}$. On suppose $\tilde{M}\not=\tilde{G}$. On pose
$${^c\varphi}={^c\theta}_{\tilde{M}}^{\tilde{G},{\cal E}}({\bf M}',{\bf f})-({^c\theta}_{\tilde{M}}^{\tilde{G}}({\bf f}))^{{\bf M}'},$$
$${^c\varphi}^{rat}={^c\theta}_{\tilde{M}}^{rat,\tilde{G},{\cal E}}({\bf M}',{\bf f})-({^c\theta}_{\tilde{M}}^{rat,\tilde{G}}({\bf f}))^{{\bf M}'},$$
$$\varphi^{rat}=\theta_{\tilde{M}}^{rat,\tilde{G},{\cal E}}({\bf M}',{\bf f})-(\theta_{\tilde{M}}^{rat,\tilde{G}}({\bf f}))^{{\bf M}'}.$$
Soit $\boldsymbol{\delta}\in D_{orb,\tilde{G}-reg}^{st}({\bf M}')\otimes Mes(M'({\mathbb R}))^*$.
D'apr\`es 5.13(5) et le lemme 7.1(iii), on a
$${^cI}_{\tilde{M}}^{\tilde{G},{\cal E}}({\bf M}',\boldsymbol{\delta},{\bf f})=S^{{\bf M}'}(\boldsymbol{\delta},{^c\theta}_{\tilde{M}}^{\tilde{G},{\cal E}}({\bf M}',{\bf f}))+\sum_{\tilde{L}\in {\cal L}(\tilde{M}),\tilde{L}\not=\tilde{M}}I_{\tilde{M}}^{\tilde{L},{\cal E}}({\bf M}',\boldsymbol{\delta},{^c\theta}_{\tilde{L}}^{\tilde{G}}({\bf f}))$$
et
$${^cI}_{\tilde{M}}^{\tilde{G}}(transfert(\boldsymbol{\delta}),{\bf f})=\sum_{\tilde{L}\in {\cal L}(\tilde{M})}I_{\tilde{M}}^{\tilde{L}}(transfert(\boldsymbol{\delta}),{^c\theta}_{\tilde{L}}^{\tilde{G}}({\bf f})).$$
On a prouv\'e en [V] proposition 1.13 que l'on avait l'\'egalit\'e
$$I_{\tilde{M}}^{\tilde{L},{\cal E}}({\bf M}',\boldsymbol{\delta},\boldsymbol{\varphi})=I_{\tilde{M}}^{\tilde{L}}(transfert(\boldsymbol{\delta}),\boldsymbol{\varphi})$$
pour tout $\boldsymbol{\varphi}\in I(\tilde{L}({\mathbb R}))\otimes Mes(L({\mathbb R}))$. Cette \'egalit\'e s'\'etend imm\'ediatement \`a $\boldsymbol{\varphi}\in I_{ac}(\tilde{L}({\mathbb R}))\otimes Mes(L({\mathbb R}))$. Par diff\'erence, on obtient
$${^cI}_{\tilde{M}}^{\tilde{G},{\cal E}}({\bf M}',\boldsymbol{\delta},{\bf f})-{^cI}_{\tilde{M}}^{\tilde{G}}(transfert(\boldsymbol{\delta}),{\bf f})=S^{{\bf M}'}(\boldsymbol{\delta},{^c\theta}_{\tilde{M}}^{\tilde{G},{\cal E}}({\bf M}',{\bf f}))-I^{\tilde{M}}(transfert(\boldsymbol{\delta}),{^c\theta}_{\tilde{M}}^{\tilde{G}}({\bf f})).$$
Par d\'efinition du transfert des distributions, le dernier terme est aussi $S^{{\bf M}'}(\boldsymbol{\delta},{^c\theta}_{\tilde{M}}^{\tilde{G}}({\bf f})^{{\bf M}'})$. On obtient alors
$$(1)\qquad {^cI}_{\tilde{M}}^{\tilde{G},{\cal E}}({\bf M}',\boldsymbol{\delta},{\bf f})-{^cI}_{\tilde{M}}^{\tilde{G}}(transfert(\boldsymbol{\delta}),{\bf f})=S^{{\bf M}'}(\boldsymbol{\delta},{^c\varphi}).$$
Transf\'erons la formule du lemme 5.12. On obtient
$${^c\theta}_{\tilde{M}}^{\tilde{G}}({\bf f})^{{\bf M}'}=\sum_{\tilde{L}\in {\cal L}(\tilde{M})}(\theta_{\tilde{M}}^{rat,\tilde{L}}\circ{^c\theta}_{\tilde{L}}^{rat, \tilde{G}}({\bf f}))^{{\bf M}'}.$$
Pour $\tilde{L}=\tilde{M}$,  $\theta_{\tilde{M}}^{rat,\tilde{M}}$ est l'identit\'e. Pour $\tilde{L}=\tilde{G}$, $^c\theta_{\tilde{G}}^{rat,\tilde{G}}$ est l'identit\'e. 
Pour $\tilde{L}\not=\tilde{M}$, $\tilde{L}\not=\tilde{G}$, on peut appliquer la formule 6.11(A) par r\'ecurrence:
$$(\theta_{\tilde{M}}^{rat,\tilde{L}}\circ{^c\theta}_{\tilde{L}}^{rat, \tilde{G}}({\bf f}))^{{\bf M}'}=\theta_{\tilde{M}}^{rat,\tilde{L},{\cal E}}({\bf M}',{^c\theta}_{\tilde{L}}^{rat, \tilde{G}}({\bf f})^{{\bf M}'}).$$
 La formule ci-dessus se r\'ecrit
$${^c\theta}_{\tilde{M}}^{\tilde{G}}({\bf f})^{{\bf M}'}={^c\theta}_{\tilde{M}}^{rat,\tilde{G}}({\bf f})^{{\bf M}'}+\theta_{\tilde{M}}^{rat,\tilde{G}}({\bf f})^{{\bf M}'}+\sum_{\tilde{L}\in {\cal L}(\tilde{M}), \tilde{L}\not=\tilde{M},\tilde{L}\not=\tilde{G}}\theta_{\tilde{M}}^{rat,\tilde{L},{\cal E}}({\bf M}',{^c\theta}_{\tilde{L}}^{rat, \tilde{G}}({\bf f})^{{\bf M}'}).$$
Utilisons la formule du lemme 7.2(iii). Dans cette formule, on peut remplacer le terme index\'e par $\tilde{L}=\tilde{G}$ par $\theta_{\tilde{M}}^{rat,\tilde{G},{\cal E}}({\bf M}',{\bf f})$ puisque $^c\theta_{\tilde{G}}^{rat,\tilde{G}}$ est l'identit\'e. On obtient
$${^c\theta}_{\tilde{M}}^{\tilde{G}}({\bf f})^{{\bf M}'}={^c\theta}_{\tilde{M}}^{rat,\tilde{G},{\cal E}}({\bf M}',{\bf f})+\theta_{\tilde{M}}^{rat,\tilde{G},{\cal E}}({\bf M}',{\bf f})+\sum_{\tilde{L}\in {\cal L}(\tilde{M}), \tilde{L}\not=\tilde{M},\tilde{L}\not=\tilde{G}}\theta_{\tilde{M}}^{rat,\tilde{L},{\cal E}}({\bf M}',{^c\theta}_{\tilde{L}}^{rat, \tilde{G}}({\bf f})^{{\bf M}'}).$$
En comparant les deux formules ci-dessus, on obtient l'\'egalit\'e
$$(2) \qquad {^c\varphi}={^c\varphi}^{rat}+\varphi^{rat}.$$
A partir des \'egalit\'es (1) et (2), la preuve devient similaire \`a celles des deux paragraphes pr\'ec\'edents. On doit utiliser la variante 6.11(2) de la propri\'et\'e 6.6(1). On  laisse les d\'etails au lecteur. $\square$

\bigskip

\section{L'application $\epsilon_{\tilde{M}}$}

\bigskip

 \subsection{Un lemme \'el\'ementaire}
 On rappelle le lemme 6.2 de [A8] dont nous ferons usage plusieurs fois. Consid\'erons un espace $V={\mathbb R}^n$, notons $B\subset V$ la boule ferm\'ee  de centre $0$ et de rayon $1$. Notons $Diff^{cst}(V)$ l'espace des op\'erateurs diff\'erentiels \`a coefficients constants sur $V$. Consid\'erons une famille finie $(l_{i})_{i=1,...,k}$ de formes lin\'eaires non  nulles sur $V$.  Notons $V_{reg}=\{v\in V; \forall i=1,...,k, l_{i}(v)\not=0\}$. 
 
 \ass{Lemme}{Soit $q\in {\mathbb N}$. Il existe une famille $(\tau_{j})_{j\in {\mathbb N}}$ d'applications $\tau_{j}$ de $Diff^{cst}(V) $ dans lui-m\^eme v\'erifiant les propri\'et\'es suivantes:
 
 (i) pour tout $D$,  il n'y a qu'un nombre fini de $j$ tels que $\tau_{j}(D)\not=0$;
 
 (ii) soit $c:Diff^{cst}(V)\to {\mathbb R}_{\geq0}$ une fonction  telle que $c(0)=0$. D\'efinissons $c^*:Diff^{cst}(V)\to {\mathbb R}_{\geq0}$ par $c^*(D)=\sum_{j\in {\mathbb N}}c\circ \tau_{j}(D)$. Soit $\varphi\in C_{c}^{\infty}(V_{reg})$. Supposons que, pour tout $D\in Diff^{cst}(V)$ et tout $v\in B\cap V_{reg}$, on ait l'in\'egalit\'e
 $$\vert D\varphi(v)\vert\leq c(D)\prod_{i=1,...,k}\vert l_{i}(v)\vert^{-q}.$$
 Alors, pour tout $D\in Diff^{cst}(V)$ et tout $v\in B\cap V_{reg}$, on a l'in\'egalit\'e
 $$\vert D\varphi(v)\vert \leq c^*(D).$$} 

\bigskip

\subsection{D\'efinition locale}
On consid\`ere un triplet $(KG,K\tilde{G},{\bf a})$ et un $K$-espace de Levi $K\tilde{M}\in {\cal L}(K\tilde{M}_{0})$ de $K\tilde{G}$. Dans ce paragraphe, on fixe une composante connexe $\tilde{M}$ de $K\tilde{M}$, on note $\tilde{G}$ la composante correspondante de $K\tilde{G}$.

\ass{Proposition}{Soient ${\bf f}\in I(K\tilde{G}({\mathbb R}),\omega)\otimes Mes(G({\mathbb R}))$ et $\eta$ un \'el\'ement semi-simple de $\tilde{M}({\mathbb R})$. Alors il existe un voisinage $\tilde{U}$ de $\eta $ dans $\tilde{M}({\mathbb R})$ et une fonction $\boldsymbol{\varphi}\in I_{cusp}(\tilde{M}({\mathbb R}),\omega)\otimes Mes(M({\mathbb R}))$ telle que 
$$I^{\tilde{M}}(\boldsymbol{\gamma},\boldsymbol{\varphi})=I_{K\tilde{M}}^{K\tilde{G},{\cal E}}(\boldsymbol{\gamma},{\bf f})-I_{K\tilde{M}}^{K\tilde{G}}(\boldsymbol{\gamma},{\bf f})$$
pour tout $\boldsymbol{\gamma}\in D_{orb}(\tilde{M}({\mathbb R}),\omega)\otimes Mes(M({\mathbb R}))^*$ dont le support est form\'e d'\'el\'ements de $\tilde{U}$ qui sont fortement r\'eguliers dans $\tilde{G}$.}

Preuve. Pour simplifier, on fixe des mesures de Haar $dm$ sur $M({\mathbb R})$, $dg$ sur $G({\mathbb R})$ et on \'ecrit ${\bf f}=f\otimes dg$ avec $f\in  I(K\tilde{G}({\mathbb R}),\omega)$. Soit $\tilde{T}$ un sous-tore tordu maximal de $\tilde{M}$ tel que $\eta\in \tilde{T}({\mathbb R})$. On fixe une mesure de Haar $dt$ sur $T^{\theta,0}({\mathbb R})$. Pour $\gamma\in \tilde{T}({\mathbb R})\cap \tilde{G}_{reg}({\mathbb R})$, l'int\'egrale orbitale sur $\tilde{M}({\mathbb R})$  associ\'ee \`a  $\gamma$ et $dt$ est bien d\'efinie. Notons-la $\boldsymbol{\gamma}$ pour un instant.   On pose simplement
$$I_{K\tilde{M}}^{K\tilde{G}}(\gamma,\omega,f)=I_{K\tilde{M}}^{K\tilde{G}}(\boldsymbol{\gamma},f\otimes dg).$$
On d\'efinit de fa\c{c}on similaire $I_{K\tilde{M}}^{K\tilde{G},{\cal E}}(\gamma,\omega,f)$. On d\'efinit une fonction $\varphi_{f,\tilde{T}}$ presque partout sur $\tilde{T}({\mathbb R})$ par
$$\varphi_{f,\tilde{T}}(\gamma)=I_{K\tilde{M}}^{K\tilde{G},{\cal E}}(\gamma,\omega,f)-I_{K\tilde{M}}^{K\tilde{G}}(\gamma,\omega,f).$$

Soit $m\in  Z_{M}(\eta;{\mathbb R})$. Posons $\tilde{T}'=ad_{m}(\tilde{T})$ et transportons par $ad_{m}$ la mesure sur $T^{\theta,0}({\mathbb R})$ en une mesure sur $T^{_{'},\theta,0}({\mathbb R})$. Il est clair que, pour presque tout $\gamma\in \tilde{T}({\mathbb R})$,  on a l'\'egalit\'e
$$(1) \qquad \varphi_{f,\tilde{T}'}(m\gamma m^{-1})=\omega(m)\varphi_{f,\tilde{T}}(\gamma).$$

On a

(2) si $\tilde{T}$ n'est pas elliptique dans $\tilde{M}$, alors $\varphi_{f,\tilde{T}}=0$.

En effet, il existe un espace de Levi $\tilde{R}\subsetneq \tilde{M}$ tel que $\tilde{T}\subset \tilde{R}$. Quitte \`a conjuguer $\tilde{T}$, cet espace donne naissance \`a un $K$-espace de Levi $K\tilde{R}\in {\cal L}(K\tilde{M}_{0})$. D'apr\`es les formules de descente [II] lemme 1.7 et 1.15(1) , on a l'\'egalit\'e
$$\varphi_{f,\tilde{T}}(\gamma)=\sum_{K\tilde{L}\in {\cal L}(K\tilde{R})}d_{\tilde{R}}^{\tilde{G}}(\tilde{M},\tilde{L})\left(I_{K\tilde{R}}^{K\tilde{L},{\cal E}}(\gamma,\omega,f_{K\tilde{L},\omega})-I_{K\tilde{R}}^{K\tilde{L}}(\gamma,\omega,f_{K\tilde{L},\omega})\right).$$
Les espaces $K\tilde{L}$ intervenant sont des $K$-espaces propres de $K\tilde{G}$. Nos hypoth\`eses de r\'ecurrence assurent que tous les termes de ce d\'eveloppement sont nuls. D'o\`u (2).

Supposons maintenant $\tilde{T}$ elliptique. Notons $\Sigma^{G_{\eta}}(T^{\theta,0})$ et $\Sigma^{M_{\eta}}(T^{\theta,0})$ les ensembles de racines de $T^{\theta,0}$ dans $G_{\eta}$, resp. $M_{\eta}$. L'hypoth\`ese d'ellipticit\'e entra\^{\i}ne que tous les \'el\'ements de $\Sigma^{M_{\eta}}(T^{\theta,0})$ sont imaginaires. Par contre, aucun \'el\'ement de $\Sigma^{G_{\eta}}(T^{\theta,0})-\Sigma^{M_{\eta}}(T^{\theta,0})$ n'est imaginaire: un tel \'el\'ement se restreint non trivialement \`a $A_{M_{\eta}}$. On fixe un sous-ensemble positif dans $\Sigma^{M_{\eta}}(T^{\theta,0})$ et on d\'efinit une fonction $\Delta_{\eta}$ presque partout sur $\mathfrak{t}^{\theta}({\mathbb R})$ par
$$\Delta_{\eta}(X)=\prod_{\alpha\in \Sigma^{M_{\eta}}(T^{\theta,0}), \alpha>0}sgn(i\alpha(X)).$$
On va prouver

(3) si $\tilde{T}$ est elliptique dans $\tilde{M}$, la fonction $X\mapsto \Delta_{\eta}(X)\varphi_{f,\tilde{T}}(exp(X)\eta)$ se prolonge en une fonction $C^{\infty}$ dans un voisinage de $0$ dans $\mathfrak{t}^{\theta}({\mathbb R})$.

Admettons ce r\'esultat. La  th\'eorie de la descente et les r\'esultats de Bouaziz et Shelstad caract\'erisent les fonctions sur les tores $\tilde{T}$ comme ci-dessus qui sont au voisinage de $\eta$ les int\'egrales orbitales d'une fonction $C^{\infty}$ et cuspidale sur $\tilde{M}({\mathbb R})$. Ce sont pr\'ecis\'ement celles qui v\'erifient  les propri\'et\'es (1), (2) et (3). On obtient alors   l'assertion de la proposition.

Prouvons (3). On suppose donc $\tilde{T}$ elliptique dans $\tilde{M}$. Supposons d'abord que le point $\eta$ est $\tilde{G}$-\'equisingulier. D'apr\`es [V] 1.3 et lemme 1.9, il existe une fonction $\varphi\in I(\tilde{M}({\mathbb R}),\omega)$ telle que, pour tout  sous-tore tordu maximal $\tilde{T}'$ de $\tilde{M}$  tel que $\eta\in \tilde{T}'({\mathbb R})$ et pour tout $\gamma\in \tilde{T}'({\mathbb R})\cap \tilde{G}_{reg}({\mathbb R})$ assez voisin de $\eta$, on ait l'\'egalit\'e
$$\varphi_{f,\tilde{T}'}(\gamma)=I^{\tilde{M}}(\gamma,\omega,\varphi).$$
La fonction $\varphi$ est cuspidale d'apr\`es (2). Alors la propri\'et\'e (3) r\'esulte des propri\'et\'es habituelles des int\'egrales orbitales. 

 Notons $\mathfrak{t}_{1}$, resp. $\mathfrak{t}_{2}$,  le sous-ensemble des $X\in \mathfrak{t}^{\theta}$ tels que $\alpha(X)\not=0$ pour tout $\alpha\in \Sigma^{G_{\eta}}(T^{\theta,0})-\Sigma^{M_{\eta}}(T^{\theta,0})$, resp. $\alpha\in \Sigma^{G_{\eta}}(T^{\theta,0})$. Soit $X_{0}\in \mathfrak{t}_{1}({\mathbb R})$, supposons $X_{0}$ proche de $0$. Le point $\eta_{0}=exp(X_{0})\eta$ est $\tilde{G}$-\'equisingulier. On vient de prouver que la fonction $X'\mapsto \Delta_{\eta_{0}}(X')\varphi_{f,\tilde{T}}(exp(X')\eta_{0})$ se prolongeait en une fonction $C^{\infty}$ au voisinage de $0$. Mais on a l'\'egalit\'e $\Delta_{\eta_{0}}(X')=\epsilon\Delta_{\eta}(X_{0}+X')$ pour $X'$ proche de $0$, avec un $\epsilon\in \{\pm 1\}$ constant. Donc la fonction $X\mapsto \Delta_{\eta}(X)\varphi_{f,\tilde{T}}(exp(X)\eta)$ se prolonge en une fonction $C^{\infty}$ au voisinage de $X_{0}$. Soit $\Omega$ une composante connexe de $\mathfrak{t}_{1}({\mathbb R})$. On obtient

(4) il existe un voisinage $\mathfrak{u}$ de $0$ dans $\mathfrak{t}^{\theta}({\mathbb R})$ tel que la fonction 
$X\mapsto \Delta_{\eta}(X)\varphi_{f,\tilde{T}}(exp(X)\eta)$ se prolonge en une fonction $C^{\infty}$ sur $\Omega\cap \mathfrak{u}$. 

Consid\'erons le groupe  des $g\in G$ tels que $ad_{g}(\eta)=\eta$ et  $ad_{g}(\tilde{T})=\tilde{T}$. Il contient $T^{\theta}$ comme sous-groupe distingu\'e d'indice fini. Notons ${\cal W}$ le quotient. Le groupe ${\cal W}$ agit naturellement dans $\mathfrak{t}^{\theta}$. Pour $w\in {\cal W}-\{1\}$, l'ensemble des points fixes est un sous-espace propre. Notons  $\mathfrak{t}_{3}$  le compl\'ementaire dans $\mathfrak{t}_{2}$ de la r\'eunion de ces  sous-espaces. On voit que, pour $X\in  \mathfrak{t}_{3}({\mathbb R})$ assez proche de $0$, l'\'el\'ement $exp(X)\eta$ est fortement r\'egulier dans $\tilde{G}$.  Quitte \`a restreindre $\mathfrak{u}$, la propri\'et\'e 3.2(4) et le lemme 3.4 entra\^{\i}nent
que

(5) pour tout $U\in Sym(\mathfrak{t}^{\theta})$, il existe un entier $N$ et, pour tout $f\in I(\tilde{G}({\mathbb R}),\omega)$, il existe $c>0$ de sorte que l'on ait la majoration
$$\vert \partial_{U}\varphi_{f,\tilde{T}}(exp(X)\eta)\vert \leq  c D^{G_{\eta}}(X)^{-N}$$
pour tout $X\in \mathfrak{t}_{3}({\mathbb R})\cap \mathfrak{u}$. 

Rappelons l'homomorphisme d'Harish-Chandra, que l'on peut interpr\'eter ici comme un homomorphisme
$$\begin{array}{ccc} \mathfrak{Z}(G)&\to&Sym(\mathfrak{t}^{\theta})\\ z&\mapsto &z_{T^{\theta,0}}.\\ \end{array}$$
Montrons que

(6) on a  l'\'egalit\'e
$$\varphi_{zf,\tilde{T}}(exp(X)\eta)=\partial_{z_{T^{\theta,0}}}\varphi_{f,\tilde{T}}(exp(X)\eta)$$
pour tout $X\in \mathfrak{t}_{3}({\mathbb R})$.

En effet, la proposition 2.5 conduit \`a l'\'egalit\'e
$$\varphi_{zf,\tilde{T}}(exp(X)\eta)=\sum_{K\tilde{L}\in {\cal L}(K\tilde{M})}\delta_{\tilde{M}}^{\tilde{L}}(exp(X)\eta,z_{L})(I_{K\tilde{L}}^{K\tilde{G},{\cal E}}(exp(X)\eta,\omega,f)-I_{K\tilde{L}}^{K\tilde{G}}(exp(X)\eta,\omega,f)).$$
Les hypoth\`eses de r\'ecurrence assurent que tous les termes sont nuls, sauf celui pour $\tilde{L}=\tilde{M}$. Celui-ci est \'egal \`a $\partial_{z_{T^{\theta,0}}}\varphi_{f,\tilde{T}}(exp(X)\eta)$. D'o\`u (6).

L'alg\`ebre $Sym(\mathfrak{t}^{\theta})$ est de type fini sur l'image de l'homomorphisme d'Harish-Chandra. On peut donc fixer des \'el\'ements $U_{1},...,U_{k}\in Sym(\mathfrak{t}^{\theta})$ de sorte que tout $U\in Sym(\mathfrak{t}^{\theta})$ s'\'ecrive $U=\sum_{i=1,...,k}U_{i}z_{i,T^{\theta,0}}$, avec des $z_{i}\in  \mathfrak{Z}(G)$. On a alors
$$\partial_{U}\varphi_{f,\tilde{T}}(exp(X)\eta)=\sum_{i=1,...,k}\partial_{U_{i}}\varphi_{z_{i}f,\tilde{T}}(exp(X)\eta).$$
L'assertion (5) associe \`a chaque $U_{i}$ un entier $N_{i}$. En prenant pour $N$ un majorant de ces entiers, on obtient qu'il existe $N$ tel que, pour  tout $U\in Sym(\mathfrak{t}^{\theta})$ et pour tout $f\in I(\tilde{G}({\mathbb R}),\omega)$, il existe $c>0$ de sorte que l'on ait la majoration
$$\vert \partial_{U}\varphi_{f,\tilde{T}}(exp(X)\eta)\vert \leq c D^{G_{\eta}}(X)^{-N}$$
pour tout $X\in \mathfrak{t}_{3}({\mathbb R})\cap \mathfrak{u}$. On applique alors le lemme 8.1. On obtient
que,  pour tout $U\in Sym(\mathfrak{t}^{\theta})$ et pour tout $f\in I(\tilde{G}({\mathbb R}),\omega)$, la fonction $\partial_{U}\varphi_{f,\tilde{T}}(exp(X)\eta)$ est born\'ee sur $\mathfrak{t}_{3}({\mathbb R})\cap \mathfrak{u}$. Joint \`a (4), ce r\'esultat implique

(7) pour tout $U\in Sym(\mathfrak{t}^{\theta})$ et pour tout $f\in I(\tilde{G}({\mathbb R}),\omega)$, la fonction $\partial_{U}\varphi_{f,\tilde{T}}(X)$ est born\'ee sur $\Omega\cap \mathfrak{u}$. 

Il r\'esulte de (4) et (7) que cette fonction se prolonge continuement \`a l'adh\'erence de $\Omega\cap \mathfrak{u}$ et m\^eme \`a un voisinage de cette adh\'erence, cf. [B] remarque 3.2.

Soit $\alpha\in \Sigma^{G_{\eta}}(T^{\theta,0})$ une racine r\'eelle. Soit $X_{0}\in \mathfrak{u}$ tel que l'ensemble des  \'el\'ements de $ \Sigma^{G_{\eta}}(T^{\theta,0})$ annulant $X_{0}$ soit $\{\pm \alpha\}$. On va prouver

(8)   la fonction $X\mapsto \Delta_{\eta}(X)\varphi_{f,\tilde{T}}(exp(X)\eta)$ se  prolonge en une fonction $C^{\infty}$ au voisinage de $X_{0}$.

  Puisque $\alpha$ est r\'eelle, la fonction $\Delta_{\eta}(X)$ est constante au voisinage  de $X_{0}$ et on peut l'oublier.  Le point $X_{0}$ appartient aux adh\'erences de deux composantes connexes $\Omega_{1}$ et $\Omega_{2}$ s\'epar\'ees par l'hyperplan annul\'e par $\alpha$. Soit $U\in Sym(\mathfrak{t}^{\theta})$. Il suffit de montrer que les deux prolongements de 
$\partial_{U}\varphi_{f,\tilde{T}}(exp(X)\eta)$ dans les adh\'erences de ces deux composantes co\"{\i}ncident sur cet hyperplan. Soit $X_{1}$ un \'el\'ement 
 de cet hyperplan assez voisin de $X_{0}$. Posons $\eta_{1}=exp(X_{1})\eta$.  Le couple $(\eta_{1},\tilde{T})$ v\'erifie les hypoth\`eses de la section 4 et on utilise les notations de cette section.Les valeurs en $X_{1}$ des prolongements ci-dessus sont les deux limites 
$$(9) \qquad lim_{r\to 0+}\partial_{U}\varphi_{f,\tilde{T}}(exp(rH_{d})\eta_{1}) \,\,\text{ et }lim_{r\to 0-}\partial_{U}\varphi_{f,\tilde{T}}(exp(rH_{d})\eta_{1}).$$
Il s'agit de prouver qu'elles sont \'egales.  Pour $X$ voisin de $0$, on a 
$$\partial_{U}\varphi_{f,\tilde{T}}(exp(X)\eta_{1})=\partial_{U}I_{K\tilde{M}}^{K\tilde{G},{\cal E},mod}(exp(X)\eta_{1},\omega,f)-\partial_{U}I_{K\tilde{M}}^{K\tilde{G},mod}(exp(X)\eta_{1},\omega,f).$$
En effet, si on oublie les exposants $mod$, c'est la d\'efinition. Ajouter ces exposants ajoute au terme de droite l'image par $\partial_{U}$ de la fonction
$$\vert \check{\alpha}\vert log(\vert \alpha(X)\vert )(I_{K\underline{\tilde{M}}}^{K\tilde{G},{\cal E}}(exp(X)\eta_{1},\omega,f)-I_{K\underline{\tilde{M}}}^{K\tilde{G}}(exp(X)\eta_{1},\omega,f)).$$
Or cette fonction est nulle par hypoth\`ese de r\'ecurrence puisque $dim(A_{\underline{\tilde{M}}})< dim(A_{\tilde{M}})$.  La diff\'erence entre les deux limites (9) est calcul\'ee par les propositions 4.1 et 4.3. Le m\^eme argument montre qu'elle est nulle. Cela prouve (8).

Supposons maintenant que $\alpha$ soit imaginaire au lieu d'\^etre r\'eelle.  Le m\^eme r\'esultat vaut: la racine $\alpha$ appartient \`a $\Sigma^{M_{\eta}}(T^{\theta,0})$, le point $X_{0}$ appartient \`a $\mathfrak{t}_{1}({\mathbb R})$ et l'assertion r\'esulte de (4). Pour une racine $\alpha$ qui n'est ni r\'eelle, ni imaginaire, le sous-espace des \'el\'ements de $\mathfrak{t}^{\theta}({\mathbb R})$ annul\'es par $\alpha$ est de codimension deux: un tel \'el\'ement est annul\'e par $\alpha$ et par son conjugu\'e.  Mais alors, la fonction $X\mapsto \Delta_{\eta}(X)\varphi_{f,\tilde{T}}(exp(X)\eta)$ et l'ensemble de racines $\Sigma^{G_{\eta}}(T^{\theta,0})$ v\'erifient les hypoth\`eses du lemme 21 de [Va]. Ce lemme conclut que cette fonction se prolonge en une fonction $C^{\infty}$ au voisinage de $0$. Cela prouve (3) et la proposition. $\square$

\bigskip

\subsection{D\'efinition globale}
Pour ${\bf f}\in I(K\tilde{G}({\mathbb R}),\omega,K)\otimes Mes(G({\mathbb R}))$, posons
$$^c\boldsymbol{\varphi}_{{\bf f}}={^c\theta}_{K\tilde{M}}^{K\tilde{G},{\cal E}}({\bf f})-{^c\theta}_{K\tilde{M}}^{K\tilde{G}}({\bf f}).$$
C'est un \'el\'ement de $I_{ac}(K\tilde{M}({\mathbb R}),\omega,K)\otimes Mes(M({\mathbb R}))$. Un argument de descente d\'ej\`a utilis\'e plusieurs fois montre qu'il est cuspidal.

\ass{Proposition}{Soit ${\bf f}\in I(K\tilde{G}({\mathbb R}),\omega,K)\otimes Mes(G({\mathbb R}))$.  Il existe une fonction $\boldsymbol{\phi}\in I_{cusp}(K\tilde{M}({\mathbb R}),\omega)\otimes Mes(M({\mathbb R}))$ de sorte que l'on ait l'\'egalit\'e
$$I^{\tilde{M}}(\boldsymbol{\gamma},\boldsymbol{\phi}-{^c\boldsymbol{\varphi}}_{{\bf f}})=I_{K\tilde{M}}^{K\tilde{G},{\cal E}}(\boldsymbol{\gamma},{\bf f})-I_{K\tilde{M}}^{K\tilde{G}}(\boldsymbol{\gamma},{\bf f})$$
pour tout $\boldsymbol{\gamma}\in D_{orb}(K\tilde{M}({\mathbb R}),\omega)\otimes Mes(M({\mathbb R}))^*$ dont le support est form\'e d'\'el\'ements  fortement r\'eguliers dans $\tilde{G}$.}

Preuve. On reprend le d\'ebut de la preuve de 7.4. On n'impose plus l'hypoth\`ese (Hyp) de ce paragraphe. Le terme qui \'etait annul\'e par cette hypoth\`ese ne l'est plus et l'\'egalit\'e (1) de ce paragraphe se transforme en
$$^cI_{K\tilde{M}}^{K\tilde{G},{\cal E}}(\boldsymbol{\gamma},{\bf f})-{^cI}_{K\tilde{M}}^{K\tilde{G}}(\boldsymbol{\gamma},{\bf f})=I^{K\tilde{M}}(\boldsymbol{\gamma},{^c\boldsymbol{\varphi}}_{{\bf f}})+(I_{K\tilde{M}}^{K\tilde{G},{\cal E}}(\boldsymbol{\gamma},{\bf f})-I_{K\tilde{M}}^{K\tilde{G}}(\boldsymbol{\gamma},{\bf f})).$$
Le deuxi\`eme terme du membre de droite est localement une int\'egrale orbitale d'apr\`es la proposition 8.2. Il en est trivialement de m\^eme du premier terme. Donc le membre de gauche est localement une int\'egrale orbitale. Or il est \`a support dans un ensemble $\Gamma^M$ o\`u $\Gamma$ est compact. Par partition de l'unit\'e, on peut trouver $\boldsymbol{\phi}\in I(K\tilde{M}({\mathbb R}),\omega)\otimes Mes(M({\mathbb R}))$ de sorte que ce membre de gauche soit \'egal \`a $I^{\tilde{M}}(\boldsymbol{\gamma},\boldsymbol{\phi})$. On obtient la formule de l'\'enonc\'e. Les formules de descente habituelles et les hypoth\`eses de r\'ecurrence entra\^{\i}nent que le membre de droite de cette formule de l'\'enonc\'e est nul si le support de $\boldsymbol{\gamma}$ n'est pas form\'e d'\'el\'ements elliptiques. Cela entra\^{\i}ne que $\boldsymbol{\phi}-{^c\boldsymbol{\varphi}}_{{\bf f}}$ est cuspidale, donc aussi $\boldsymbol{\phi}$. $\square$

La fonction $\boldsymbol{\phi}$ de l'\'enonc\'e est uniquement d\'etermin\'ee. On pose
$$\epsilon_{K\tilde{M}}({\bf f})=\boldsymbol{\phi}-{^c\boldsymbol{\varphi}}_{{\bf f}}.$$
On a ainsi d\'efini une application lin\'eaire
$$\epsilon_{K\tilde{M}}: I(K\tilde{G}({\mathbb R}),\omega,K)\otimes Mes(G({\mathbb R}))\to I_{ac,cusp}(K\tilde{M}({\mathbb R}),\omega)\otimes Mes(M({\mathbb R})).$$
A ce point, on peut pr\'eciser que $\epsilon_{K\tilde{M}}$ prend ses valeurs dans 
$$(I_{cusp}(\tilde{M}({\mathbb R}),\omega)+I_{ac,cusp}(K\tilde{M}({\mathbb R}),\omega,K))\otimes Mes(M({\mathbb R})).$$
On prouvera dans les paragraphes suivants que l'on peut supprimer le premier terme  $I_{cusp}(\tilde{M}({\mathbb R}),\omega)$, autrement dit que $\epsilon_{K\tilde{M}}$ prend ses valeurs dans le sous-espace des fonctions $K$-finies.

Notons une propri\'et\'e importante:

(1) on a l'\'egalit\'e $\epsilon_{K\tilde{M}}(z{\bf f})=z_{\tilde{M}}\epsilon_{K\tilde{M}}({\bf f})$ pour tout ${\bf f}\in I(K\tilde{G}({\mathbb R}),\omega,K)\otimes Mes(G({\mathbb R}))$ et tout $z\in  \mathfrak{Z}(G)$.

Cela r\'esulte de 8.1(6). 

On a aussi

(2) pour toute $\omega$-repr\'esentation elliptique $\tilde{\pi}$ de $K\tilde{M}({\mathbb R})$, la fonction $X\mapsto I^{K\tilde{M}}(\tilde{\pi},X,\epsilon_{\tilde{M}}(f))$ est \`a d\'ecroissance rapide; sa  transform\'ee de Fourier   $\lambda\mapsto I^{K\tilde{M}}(\tilde{\pi},\lambda,\epsilon_{\tilde{M}}(f))$ sur $i{\cal A}_{\tilde{M}}^*$ se prolonge en une fonction   m\'eromorphe sur ${\cal A}_{\tilde{M},{\mathbb C}}^*$ qui est \`a d\'ecroissance rapide dans les bandes verticales; ses p\^oles sont de la forme d\'ecrite en 5.2(3). 

Preuve. On rappelle qu'un $\omega$-repr\'esentation elliptique de $K\tilde{M}({\mathbb R})$ est une telle $\omega$-repr\'esentation d'une composante connexe $\tilde{M}_{p}({\mathbb R})$. L'assertion r\'esulte du fait que  $\epsilon_{K\tilde{M}}(f)$ est la somme d'une fonction  \`a support compact et de ${^c\theta}_{K\tilde{M}}^{K\tilde{G}}(f)-{^c\theta}_{K\tilde{M}}^{K\tilde{G},{\cal E}}(f)$.  On a d\'ej\`a prouv\'e (cf. 5.11 et 6.6) que cette derni\`ere fonction v\'erifiait la propri\'et\'e (2). Cette propri\'et\'e est aussi v\'erifi\'ee pour une fonction \`a support compact. D'o\`u l'assertion. $\square$

Soulignons que l'on n'affirme pas que les hyperplans polaires sont en nombre fini.

\bigskip

\subsection{Retour sur la formule des traces locale sym\'etrique}
Pour simplifier, on fixe des mesures de Haar pour la fin de l'article. Dans ce paragraphe, on consid\`ere un triplet $(G,\tilde{G},{\bf a})$. Pour un espace de Levi $\tilde{L}\in {\cal L}(\tilde{M}_{0})$, on  introduit un ensemble de repr\'esentants $\underline{{\cal E}}_{ell,0}(\tilde{L},\omega)$ des classes d'\'equivalence des $\omega$-repr\'esentations elliptiques $\tilde{\pi}$ de $\tilde{L}({\mathbb R})$ dont le caract\`ere central $\omega_{\pi}$ est trivial sur $\mathfrak{A}_{\tilde{L}}$. Introduisons de m\^eme un ensemble de repr\'esentants $\underline{{\cal E}}_{disc,0}(\tilde{L},\omega)$ des classes d'\'equivalence des $\omega$-repr\'esentations discr\`etes $\tilde{\pi}$ de $\tilde{L}({\mathbb R})$ dont le caract\`ere central $\omega_{\pi}$ est trivial sur $\mathfrak{A}_{\tilde{L}}$. On renvoie \`a [W2] 2.11 pour la notion de $\omega$-repr\'esentation discr\`ete. Soit $\tilde{\pi}\in \underline{{\cal E}}_{disc,0}(\tilde{L},\omega)$. Il r\'esulte du lemme 2.11 de [W2] que c'est l'induite d'une repr\'esentation elliptique d'un espace de Levi de $\tilde{L}$. C'est-\`a-dire que l'on peut \'ecrire
$$(1) \qquad \tilde{\pi}= c \,Ind_{\tilde{R}}^{\tilde{L}}(\tilde{\sigma}_{\nu}),$$
o\`u $c\in {\mathbb C}$, $\tilde{R}\in {\cal L}(\tilde{M}_{0})$ est un espace de Levi contenu dans $\tilde{L}$, $\tilde{\sigma}$ appartient \`a $\underline{{\cal E}}_{ell,0}(\tilde{R},\omega)$ et $\nu\in i{\cal A}_{\tilde{R}}^*$.

Inversement, soient  $\tilde{R}\in {\cal L}(\tilde{M}_{0})$ et $\tilde{\sigma}\in \underline{{\cal E}}_{ell,0}(\tilde{R},\omega)$. Montrons que 

(2) il  n'existe qu'un nombre fini de couples $(\tilde{L},\tilde{\pi})$, o\`u $\tilde{L}\in {\cal L}(\tilde{R})$ et  $\tilde{\pi}\in \underline{{\cal E}}_{disc,0}(\tilde{L},\omega)$,  tels que le couple $(\tilde{R},\tilde{\sigma})$ soit  celui intervenant dans la d\'ecomposition (1). 

Preuve.   On peut \'evidemment fixer $\tilde{L}$ et prouver la finitude de l'ensemble des $\tilde{\pi}$. La repr\'esentation $\tilde{\sigma}$   est issue d'un triplet elliptique $(M,\tau,\tilde{{\bf r}})$, cf. [W2] 2.11. Le terme $M$ est un Levi de $R$ contenant $M_{0}$ et $\tau$ est une repr\'esentation de la s\'erie discr\`ete de $M({\mathbb R})$. Soit $\nu\in i{\cal A}_{\tilde{R}}^*$. On veut que le caract\`ere central de $Ind_{\tilde{R}}^{\tilde{L}}(\sigma_{\nu})$ soit trivial sur $\mathfrak{A}_{\tilde{L}}$. Cela impose que $\nu$ appartient \`a $ i{\cal A}_{\tilde{R}}^{\tilde{L},*}$. 
Pour que $Ind_{\tilde{R}}^{\tilde{L}}(\tilde{\sigma}_{\nu})$ soit une $\omega$-repr\'esentation discr\`ete, il est n\'ecessaire qu'il existe $\gamma\in \tilde{L}({\mathbb R})$ v\'erifiant les conditions suivantes:

-  $ad_{\gamma}$ conserve $M$;

- $\tau_{\nu}\circ ad_{\gamma}\simeq \tau_{\nu}\otimes \omega$;

- l'automorphisme $w_{\gamma}$ de ${\cal A}_{M}^{\tilde{L}}$ d\'eduit de $ad_{\gamma}$ n'a pas de points fixes non nuls. 

Il n'y a qu'un nombre fini d'automorphismes $w_{\gamma}$ possibles. Pour chacun d'eux, la deuxi\`eme relation d\'etermine $(w_{\gamma}^{-1}-1)(\nu)\in i{\cal A}_{M}^{\tilde{L}*}$. D'apr\`es la troisi\`eme relation, cela revient \`a d\'eterminer $\nu$. 
 Cela prouve (2). $\square$

  Soient $\tilde{R}\in {\cal L}(\tilde{M}_{0})$  et   $\tilde{\sigma}\in \underline{{\cal E}}_{ell,0}(\tilde{R},\omega)$.  A chaque couple $(\tilde{L},\tilde{\pi})$  v\'erifiant (2),   
    associons le sous-espace affine $H=\nu+i{\cal A}_{\tilde{L}}^*$ de $i{\cal A}_{\tilde{R}}^*$. On obtient un ensemble fini ${\cal H}_{\tilde{R},1}^{\tilde{G}}(\tilde{\sigma})$ de sous-espaces affines de $i{\cal A}_{\tilde{R}}^*$.  Cet ensemble contient $i{\cal A}_{\tilde{R}}^*$: cet espace est associ\'e au couple $(\tilde{L},\tilde{\pi})=(\tilde{R},\tilde{\sigma})$ lui-m\^eme. On le sym\'etrise de la fa\c{c}on suivante.  Un \'el\'ement  $w\in W(\tilde{M}_{0})$ agit en envoyant $\tilde{R}$ sur $w(\tilde{R})$ et $\tilde{\sigma}$, sinon sur un \'el\'ement de $\underline{{\cal E}}(w(\tilde{R}),\omega)$, du moins sur un multiple d'un tel \'el\'ement, c'est-\`a-dire  $w(\tilde{\sigma})=z_{w}\tilde{\sigma}_{w}$, o\`u $z_{w}\in {\mathbb C}^{\times}$ et $\tilde{\sigma}_{w}\in \underline{{\cal E}}(w(\tilde{R}),\omega)$. On note ${\cal H}_{\tilde{R}}^{\tilde{G}}(\tilde{\sigma})$ la r\'eunion des $w^{-1}({\cal H}_{w(\tilde{R}),1}^{\tilde{G}}(\tilde{\sigma}_{w}))$ sur $w\in W(\tilde{M}_{0})$.

Pour tout $\tilde{L}\in {\cal L}(\tilde{M}_{0})$ et pour tout $\tilde{\pi}\in \underline{{\cal E}}_{disc,0}(\tilde{L},\omega)$, soit $\phi_{\tilde{L},\tilde{\pi}}$ une fonction $C^{\infty}$ et \`a croissance polynomiale sur $i{\cal A}_{\tilde{L}}^*$. Pour $f\in I(\tilde{G}({\mathbb R}),\omega,K)$, l'expression
$$J(f)=\sum_{\tilde{L}\in {\cal L}(\tilde{M}_{0})}\sum_{\tilde{\pi}\in \underline{{\cal E}}_{disc,0}(\tilde{L},\omega)} \int_{i{\cal A}_{\tilde{L}}^*}\phi_{\tilde{L},\tilde{\pi}}(\lambda)I^{\tilde{L}}(\tilde{\pi}_{\lambda},f_{\tilde{L},\omega})\,d\lambda.$$
est convergente. Le r\'esultat ci-dessus montre qu'\`a la famille $(\phi_{\tilde{L},\tilde{\pi}})_{\tilde{L}\in {\cal L}(\tilde{M}_{0}),\tilde{\pi}\in \underline{{\cal E}}_{disc,0}(\tilde{L},\omega)}$, on peut associer une famille $(\phi_{\tilde{R},\tilde{\sigma},H})_{\tilde{R}\in {\cal L}(\tilde{M}_{0}),\tilde{\sigma}\in \underline{{\cal E}}_{ell,0}(\tilde{R},\omega),H\in {\cal H}_{\tilde{R}}^{\tilde{G}}(\tilde{\sigma})}$, o\`u $\phi_{\tilde{R},\tilde{\sigma},H}$ est une fonction $C^{\infty}$ \`a croissance polynomiale sur $H$, de sorte que
$$J(f)=\sum_{\tilde{R}\in {\cal L}(\tilde{M}_{0})}\sum_{\tilde{\sigma}\in \underline{{\cal E}}_{ell,0}(\tilde{R},\omega)}\sum_{H\in {\cal H}_{\tilde{R}}^{\tilde{G}}(\tilde{\sigma})}\int_{H}\phi_{\tilde{R},\tilde{\sigma},H}(\lambda)I^{\tilde{R}}(\tilde{\sigma}_{\lambda},f_{\tilde{R},\omega})\,d\lambda.$$
En vertu des propri\'et\'es de sym\'etries de la famille $(f_{\tilde{R},\omega})_{\tilde{R}\in {\cal L}(\tilde{M}_{0})}$, on peut sym\'etriser la famille $(\phi_{\tilde{R},\tilde{\sigma},H})_{\tilde{R}\in {\cal L}(\tilde{M}_{0}),\tilde{\sigma}\in \underline{{\cal E}}_{ell,0}(\tilde{R},\omega),H\in {\cal H}_{\tilde{R}}^{\tilde{G}}(\tilde{\sigma})}$ et supposer

(3) pour $\tilde{R}\in {\cal L}(\tilde{M}_{0})$, $\tilde{\sigma}\in \underline{{\cal E}}_{ell,0}(\tilde{R},\omega)$, $H\in {\cal H}_{\tilde{R}}^{\tilde{G}}(\tilde{\sigma})$, $\lambda\in H$ et $w\in W(\tilde{M}_{0})$, on a l'\'egalit\'e $\phi_{\tilde{R},\tilde{\sigma},H}(\lambda)=z_{w}\phi_{w(\tilde{R}),\tilde{\sigma}_{w},w(H)}(w\lambda)$, o\`u on a \'ecrit $w(\tilde{\sigma})=z_{w}\tilde{\sigma}_{w}$ comme ci-dessus. 

En imposant cette condition, la famille $(\phi_{\tilde{R},\tilde{\sigma},H})_{\tilde{R}\in {\cal L}(\tilde{M}_{0}),\tilde{\sigma}\in \underline{{\cal E}}_{ell,0}(\tilde{R},\omega),H\in {\cal H}_{\tilde{R}}^{\tilde{G}}(\tilde{\sigma})}$ est enti\`erement d\'etermin\'ee par la forme lin\'eaire $f\mapsto J(f)$.

Soit $\tilde{M}\in {\cal L}(\tilde{M}_{0})$ et soit $\tilde{T}$ un sous-tore tordu maximal elliptique de $\tilde{M}$. Fixons un sous-ensemble ouvert $\tilde{U}\subset \tilde{T}({\mathbb R})$ d'adh\'erence compacte et  invariant par multiplication par $(1-\theta)(T({\mathbb R}))$. Posons $\tilde{U}_{reg}=\tilde{U}\cap \tilde{G}_{reg}({\mathbb R})$. Notons $C_{c}^{\infty}(\tilde{U}_{reg})^{\omega^{-1}-inv}$ l'espace des fonctions  $\varphi$ sur  $\tilde{T}({\mathbb R})$ qui  v\'erifient les conditions suivantes:

- $\varphi$ est $C^{\infty}$;

- on a
$$\varphi(t^{-1}\gamma t)=\omega(t)\varphi(\gamma)$$
pour tous $\gamma\in \tilde{T}({\mathbb R})$ et $t\in T({\mathbb R})$;

- le support de $\varphi$ est contenu dans $\tilde{U}_{reg}$ et est d'image compacte dans $\tilde{T}({\mathbb R})/(1-\theta)(T({\mathbb R}))$.  

Soulignons que la condition de transformation n'est pas la m\^eme qu'en 1.1: on a invers\'e $\omega$. Rappelons que l'on suppose $\omega$ unitaire (cf. introduction) donc $\omega^{-1}=\bar{\omega}$.  L'alg\`ebre d'op\'erateurs diff\'erentiels $Diff^{cst}(\tilde{T}({\mathbb R}))^{\omega^{-1}-inv}$ de 1.1 (avec $\omega$ invers\'e) agit sur $C_{c}^{\infty}(\tilde{U}_{reg})^{\omega^{-1}-inv}$.    Pour $\varphi\in C_{c}^{\infty}(\tilde{U}_{reg})^{\omega^{-1}-inv}$, posons
  $$\delta(\varphi)=inf_{\gamma\in Supp(\varphi)} D^{\tilde{G}}(\gamma) .$$

\ass{Proposition}{ Pour tout $\varphi\in  C_{c}^{\infty}(\tilde{U}_{reg})^{\omega^{-1}-inv}$, il existe une unique famille 
$$(\phi_{\tilde{R},\tilde{\sigma},H}(\varphi))_{\tilde{R}\in {\cal L}(\tilde{M}_{0}),\tilde{\sigma}\in \underline{{\cal E}}_{ell,0}(\tilde{R},\omega),H\in {\cal H}_{\tilde{R}}^{\tilde{G}}(\tilde{\sigma})}$$
 v\'erifiant les conditions suivantes.  
 
(i)  Soient  $\tilde{R}\in {\cal L}(\tilde{M}_{0})$, $\tilde{\sigma}\in \underline{{\cal E}}_{ell,0}(\tilde{R},\omega)$ et $H\in {\cal H}_{\tilde{R}}^{\tilde{G}}(\tilde{\sigma})$.
 Le terme $\phi_{\tilde{R},\tilde{\sigma},H}(\varphi)$ est une fonction $C^{\infty}$
  sur $H$. On note  $ \phi_{\tilde{R},\tilde{\sigma},H}(\varphi,\lambda)$  sa valeur en un point  $\lambda\in H$.

(ii)  Soient  $\tilde{R}\in {\cal L}(\tilde{M}_{0})$, $\tilde{\sigma}\in \underline{{\cal E}}_{ell,0}(\tilde{R},\omega)$ et $H\in {\cal H}_{\tilde{R}}^{\tilde{G}}(\tilde{\sigma})$.
  Pour $N\in {\mathbb N}$ et pour un op\'erateur diff\'erentiel \`a coefficients constants $\Delta$ sur $H$, il existe des op\'erateurs $D_{i}\in Diff^{cst}(\tilde{T}({\mathbb R}))^{\omega^{-1}-inv}$ pour $i=1,...,n$ ind\'ependants de $\varphi$  et un entier $d\in {\mathbb N}$ ind\'ependant de $\varphi$ de sorte que
$$\vert \Delta\phi_{\tilde{R},\tilde{\sigma},H}(\varphi,\lambda)\vert \leq (1+\vert \lambda\vert )^{-N}\delta(\varphi)^{-d}\sum_{i=1,...,n}sup_{\gamma\in \tilde{U}}\vert D_{i}\varphi(\gamma)\vert $$
pour tout $\lambda\in H$.

(iii) La famille v\'erifie la condition de sym\'etrie (3).

(iv) Pour tout $f\in I(\tilde{G}({\mathbb R}),\omega,K)$, on a l'\'egalit\'e
$$\int_{\tilde{T}({\mathbb R})/(1-\theta)(T({\mathbb R}))}\varphi(\gamma)I^{\tilde{G}}_{\tilde{M}}(\gamma,\omega,f)\,d\gamma=\sum_{\tilde{R}\in {\cal L}(\tilde{M}_{0})}\sum_{\tilde{\sigma}\in \underline{{\cal E}}_{ell,0}(\tilde{R},\omega)}\sum_{H\in {\cal H}_{\tilde{R}}^{\tilde{G}}(\tilde{\sigma})}\int_{H}\phi_{\tilde{R},\tilde{\sigma},H}(\varphi,\lambda)I^{\tilde{R}}(\tilde{\sigma}_{\lambda},f_{\tilde{R},\omega})\,d\lambda.$$}

Preuve. On peut \'evidemment supposer que $\omega$ est trivial sur $T^{\theta}({\mathbb R})$, sinon $ C_{c}^{\infty}(\tilde{U}_{reg})^{\omega^{-1}-inv}=\{0\}$.  On peut effectuer une partition de l'unit\'e et supposer que $\tilde{U}$ est de la forme suivante. Fixons $\eta\in \tilde{T}({\mathbb R})$.  Notons $Norm_{G}(\tilde{T})$ le normalisateur de $\tilde{T}$ dans $G$. Il contient $T$. Posons $N =Norm_{G}(\tilde{T};{\mathbb R})\cap Z_{G}(\eta;{\mathbb R})$  et $\Xi=N/T^{\theta}({\mathbb R})$. Ce groupe est fini et agit naturellement sur $\mathfrak{t}^{\theta}({\mathbb R})$. Fixons un voisinage $\mathfrak{u}$ de $0$ dans $\mathfrak{t}^{\theta}({\mathbb R})$, que l'on suppose assez petit et invariant par $\Xi$. On suppose que 
$\tilde{U}= \{ x^{-1} exp(X)\eta x; X\in \mathfrak{u}, x\in Norm_{G}(\tilde{T};{\mathbb  R})\}$. Notons $I(\tilde{U}_{reg},\omega^{-1})$ le sous-espace des \'el\'ements $\varphi\in C_{c}^{\infty}(\tilde{U}_{reg})^{\omega^{-1}-inv}$ qui v\'erifient la condition
$$\varphi(x^{-1} exp(X)\eta x)=\omega(x)\varphi(exp(X)\eta)$$
pour tout $X\in \mathfrak{u}$ et $x\in Norm_{G}(\tilde{T};{\mathbb R})$. Posons $W_{T}=Norm_{G}(\tilde{T};{\mathbb R})/T({\mathbb R})$. On a une application de sym\'etrisation 
$$\begin{array}{ccc}C_{c}^{\infty}(\tilde{U}_{reg})^{\omega^{-1}-inv}&\to&I(\tilde{U}_{reg},\omega^{-1})\\ \varphi&\mapsto& \varphi^{sym}\\ \end{array}$$ 
d\'efinie par la formule 
$$\varphi^{sym}(\gamma)=\vert W_{T}\vert ^{-1}\sum_{  x\in Norm_{G}(\tilde{T};{\mathbb R})/T({\mathbb R}) }\omega(x)^{-1}\varphi(x^{-1}\gamma x).$$
Remarquons que l'action de $Norm_{G}(\tilde{T};{\mathbb R})$ conserve  $A_{\tilde{T}}$, lequel  est \'egal \`a $A_{\tilde{M}}$ puisque  $\tilde{T}$ est un sous-tore elliptique de $\tilde{M}$. Donc cette action conserve $\tilde{M}$. Il en r\'esulte que la fonction $\gamma\mapsto I^{\tilde{G}}_{\tilde{M}}(\gamma,\omega,f)$ intervenant dans (iv) v\'erifie elle-m\^eme une condition de sym\'etrie relativement \`a cette action. Cela implique que le membre de gauche de (iv) ne change pas quand on remplace $\varphi$ par $\varphi^{sym}$. On peut donc se limiter \`a d\'emontrer la proposition pour des $\varphi\in I(\tilde{U}_{reg},\omega^{-1})$. On d\'efinit l'espace $I(\tilde{U}_{reg},\omega)$ comme on a d\'efini $I(\tilde{U}_{reg},\omega^{-1})$. Parce que $\omega$ est unitaire, l'application $\varphi\mapsto \bar{\varphi}$ est un automorphisme antilin\'eaire de $I(\tilde{U}_{reg},\omega^{-1})$ sur $I(\tilde{U}_{reg},\omega)$.   Notons   $I(\tilde{G}_{reg}({\mathbb R}),\omega)$ l'image dans $I(\tilde{G}({\mathbb R}),\omega)$ de l'espace des \'el\'ements de $C_{c}^{\infty}(\tilde{G}({\mathbb R}))$ \`a support fortement r\'egulier.  L'espace $I(\tilde{U}_{reg},\omega)$   s'identifie \`a un sous-espace de $I(\tilde{G}_{reg}({\mathbb R}),\omega)$ de la fa\c{c}on suivante. Une fonction $\varphi\in I(\tilde{U}_{reg},\omega)$ s'identifie \`a la fonction 
  $f\in I(\tilde{G}_{reg}({\mathbb R}),\omega)$ telle que $I^{\tilde{G}}(\gamma,\omega,f)=0$ pour tout $\gamma\in \tilde{G}({\mathbb R})$ qui n'est pas conjugu\'e \`a un \'el\'ement de $\tilde{T}({\mathbb R})$ et que $I^{\tilde{G}}(\gamma,\omega,f)=\varphi(\gamma)$ pour tout $\gamma\in \tilde{T}({\mathbb R})$.

 On utilise [Moe] proposition 2, dont on reprend les notations. Pour $\varphi\in I(\tilde{U}_{reg},\omega)$, on a l'\'egalit\'e
$$\int_{\tilde{T}({\mathbb R})/(1-\theta)(T({\mathbb R})}\overline{\varphi(\gamma)}I^{\tilde{G}}_{\tilde{M}}(\gamma,\omega,f)\,d\gamma=\sum_{\tilde{L}\in {\cal L}(\tilde{M}_{0})}\sum_{\tilde{\pi}\in \underline{{\cal E}}_{disc,0}(\tilde{L},\omega)}c_{\tilde{L},\tilde{\pi}} \int_{i{\cal A}_{\tilde{L}}^*}\overline{ I_{\tilde{L}}^{_{'}\tilde{G}}(\tilde{\pi}_{\lambda},\varphi)}I^{\tilde{L}}(\tilde{\pi}_{\lambda},f_{\tilde{L},\omega})\,d\lambda.$$ 
Pour $\varphi\in I(\tilde{U}_{reg},\omega^{-1})$, on en d\'eduit
$$\int_{\tilde{T}({\mathbb R})/(1-\theta)(T({\mathbb R})}\varphi(\gamma)I^{\tilde{G}}_{\tilde{M}}(\gamma,\omega,f)\,d\gamma=\sum_{\tilde{L}\in {\cal L}(\tilde{M}_{0})}\sum_{\tilde{\pi}\in \underline{{\cal E}}_{disc,0}(\tilde{L},\omega)}c_{\tilde{L},\tilde{\pi}} \int_{i{\cal A}_{\tilde{L}}^*}\overline{ I_{\tilde{L}}^{_{'}\tilde{G}}(\tilde{\pi}_{\lambda},\bar{\varphi})}I^{\tilde{L}}(\tilde{\pi}_{\lambda},f_{\tilde{L},\omega})\,d\lambda.$$

Les explications donn\'ees avant l'\'enonc\'e permettent de convertir cette formule en celle du (iv) de l'\'enonc\'e. On peut imposer (iii) et la famille obtenue est alors unique. Les conditions (i) et (ii) se d\'eduisent d'assertions analogues concernant les fonctions $\lambda\mapsto I_{\tilde{L}}^{_{'}\tilde{G}}(\tilde{\pi}_{\lambda},\varphi) $ intervenant ci-dessus, pour $\varphi\in I(\tilde{U}_{reg},\omega)$. Il reste \`a prouver ces assertions. On fixe $\tilde{L}\in {\cal L}(\tilde{M}_{0})$ et $\tilde{\pi}\in \underline{{\cal E}}_{disc,0}(\tilde{L},\omega)$.  Soit $f\in C_{c}^{\infty}(\tilde{G}({\mathbb R}))$, supposons $f$ \`a support fortement r\'egulier. On d\'efinit par r\'ecurrence un terme
$$(4) \qquad I_{\tilde{L}}^{_{'}\tilde{G}}(\tilde{\pi}_{\lambda},f)=J_{\tilde{L}}^{\tilde{G}}(\tilde{\pi}_{\lambda},f)-\sum_{\tilde{R}\in {\cal L}(\tilde{L}),\tilde{R}\not=\tilde{G}}
I_{\tilde{L}}^{_{'}\tilde{R}}(\tilde{\pi}_{\lambda},\phi_{\tilde{R}}^{_{'}\tilde{G}}(f)).$$
On rappellera plus loin la d\'efinition de $\phi_{\tilde{R}}^{_{'}\tilde{G}}(f)$. Moeglin montre que l'application lin\'eaire $f\mapsto  I_{\tilde{L}}^{_{'}\tilde{G}}(\tilde{\pi}_{\lambda},f)$ se quotiente en une application d\'efinie sur  $I(\tilde{G}_{reg}({\mathbb R}),\omega)$. Le terme $I_{\tilde{L}}^{_{'}\tilde{G}}(\tilde{\pi}_{\lambda},\varphi)$ intervenant ci-dessus est la valeur de cette application en $\varphi\in I(\tilde{G}_{reg}({\mathbb R}),\omega)$. Que ce terme soit $C^{\infty}$ en $\lambda$ est imm\'ediat par r\'ecurrence sur la formule (4) car le terme $J_{\tilde{L}}^{\tilde{G}}(\tilde{\pi}_{\lambda},f)$ est $C^{\infty}$ en $\lambda$. 
Nous voulons montrer que $I_{\tilde{L}}^{_{'}\tilde{G}}(\tilde{\pi}_{\lambda},\varphi)$ v\'erifie une majoration comme dans le (ii) de l'\'enonc\'e.  

L'ensemble $\tilde{U}^G$ des \'el\'ements de $\tilde{G}({\mathbb R})$ qui sont conjugu\'es \`a un \'el\'ement de $\tilde{U}$ est \'egal \`a $\{x^{-1} exp(X)\eta x; X\in \mathfrak{u}, x\in G({\mathbb R})\}$. Notons $\mathfrak{u}_{reg}$ le sous-ensemble des \'el\'ements de $\mathfrak{u}$ dont le stabilisateur dans $\Xi$ est r\'eduit \`a $\{1\}$. Consid\'erons l'application
$$\begin{array}{ccc}\mathfrak{u}_{reg}\times T^{\theta}({\mathbb R})\backslash G({\mathbb R})&\to& \tilde{G}({\mathbb R})\\ (X,x)&\mapsto& x^{-1} exp(X)\eta x.\\ \end{array}$$
Le groupe $N$ agit sur l'espace de d\'epart par $(y,X,x)\mapsto (ad_{y}(X),yx)$ pour $y\in N$. Cette action se quotiente en une action du groupe fini $\Xi$. L'application ci-dessus a pour image $\tilde{U}^{\tilde{G}}\cap \tilde{G}_{reg}({\mathbb R})$ et, pourvu que $\mathfrak{u}$ soit assez petit,  est un rev\^etement fini de groupe $\Xi$ au-dessus de cette image. 
  Fixons une fonction $\beta\in C_{c}^{\infty}(T^{\theta}({\mathbb R})\backslash G({\mathbb R}))$ dont l'int\'egrale vaut $1$ et qui est invariante par l'action de $\Xi$ par translation \`a gauche. Pour $\varphi\in I(\tilde{U}_{reg},\omega)$, notons $f'_{\varphi}$ la fonction sur $\mathfrak{u}_{reg}\times T^{\theta}({\mathbb R})\backslash G({\mathbb R})$ d\'efinie par
 $$f'_{\varphi}(X,x)=\omega(x)^{-1}\beta(x) D^{\tilde{G}}(exp(X)\eta) ^{-1/2}\varphi(exp(X)\eta).$$
 Elle est invariante par l'action de $\Xi$. Elle se descend donc en une fonction $f_{\varphi}$ sur $\tilde{G}({\mathbb R})$, \`a support dans $\tilde{U}^{\tilde{G}}\cap \tilde{G}_{reg}({\mathbb R})$.  
On peut fixer un sous-ensemble ouvert $\tilde{V}$ de $\tilde{G}({\mathbb R})$, d'adh\'erence compacte, qui contient 
 $x^{-1}exp(X)\eta x$ pour  tout $X\in \mathfrak{u}$ et  tout $x$ dans  le  support de $\beta$.   La fonction $f_{\varphi}$ est \`a support dans $\tilde{V}$. L'alg\`ebre enveloppante $\mathfrak{U}(G)$ agit \`a droite et \`a gauche sur $C_{c}^{\infty}(\tilde{V})$. On note ces actions $(X,f,Y)\mapsto XfY$. On munit  $C_{c}^{\infty}(\tilde{V})$ des semi-normes $f\mapsto sup_{\gamma\in \tilde{V}}\vert (XfY)(\gamma)\vert$ pour $X,Y\in \mathfrak{U}(G)$.  Montrons que
 
 (5) pour toute semi-norme ${\cal N}$ sur $C_{c}^{\infty}(\tilde{V})$, il existe un entier $d$ et un ensemble fini d'\'el\'ements $D_{i}\in Diff^{cst}(\tilde{T}({\mathbb R}))^{\omega-inv}$ pour $i=1,...,n$  de sorte que, pour tout $\varphi\in I(\tilde{U}_{reg},\omega)$, on ait 
 $${\cal N}(f_{\varphi})\leq \delta(\varphi)^{-d}\sum_{i=1,...,n}sup_{X\in \mathfrak{u}}\vert D_{i}\varphi(exp(X)\eta)\vert .$$
 
 Fixons une base $(U_{j})_{j\in {\mathbb N}}$ de $\mathfrak{U}(G)$ et une base $(H_{l})_{l\in {\mathbb N}}$ de $Sym(\mathfrak{t}^{\theta})$. L'alg\`ebre $  Sym(\mathfrak{t}^{\theta})\otimes\mathfrak{U}(G)$ agit naturellement sur l'espace des fonctions sur $\mathfrak{u}_{reg}\times T^{\theta}({\mathbb R})\backslash G({\mathbb R})$. Fixons $j_{1},j_{2}\in {\mathbb N}$. 
 En utilisant l'application d'Harish-Chandra de 1.5, on montre que l'on a une \'egalit\'e
 $$(U_{j_{1}}f_{\varphi}U_{j_{2}})(x^{-1}exp(X)\eta x)=\sum_{j,l\in {\mathbb N}}c_{j,l}(x,exp(X)\eta)(H_{l}U_{j}f'_{\varphi})(X,x)$$
 o\`u les fonctions $c_{j,l}(x,exp(X)\eta)$ sont presque toutes nulles, sont $C^{\infty}$ en $x$ et rationnelles en $exp(X)\eta$. Le (i) du lemme 1.5 montre que, multipli\'ees par une puissance convenable de $D_{\star}^{\tilde{G}}(exp(X)\eta)$, ces fonctions deviennent
polynomiales en $exp(X)\eta$. Pour prouver (5), il suffit alors d'appliquer la d\'efinition de $f'_{\varphi}$ et le fait que cette fonction est \`a support compact en $x$, ce support ne d\'ependant pas de $\varphi$. Cela prouve (5).

 L'image de $f_{\varphi}$ dans $I(\tilde{G}({\mathbb R}),\omega)$ est \'egale \`a $\varphi$.  Cela implique $I^{_{'}\tilde{G}}_{\tilde{L}}(\tilde{\pi}_{\lambda},\varphi)=I^{_{'}\tilde{G}}_{\tilde{L}}(\tilde{\pi}_{\lambda},f_{\varphi})$.  L'\'egalit\'e (4) devient
 $$(6) \qquad I_{\tilde{L}}^{_{'}\tilde{G}}(\tilde{\pi}_{\lambda},\varphi)=J_{\tilde{L}}^{\tilde{G}}(\tilde{\pi}_{\lambda},f_{\varphi})-\sum_{\tilde{R}\in {\cal L}(\tilde{L}),\tilde{R}\not=\tilde{G}}
I_{\tilde{L}}^{_{'}\tilde{R}}(\tilde{\pi}_{\lambda},\phi_{\tilde{R}}^{_{'}\tilde{G}}(f_{\varphi})).$$
 
 Rappelons  l'estimation du premier terme du membre de droite, que l'on reprend d'Arthur.     
 Pour tout $N\in {\mathbb N}$ et tout op\'erateur diff\'erentiel $\Delta$ sur $i{\cal A}_{\tilde{L}}^*$ \`a coefficients constants, il existe  un ensemble fini de semi-normes ${\cal N}_{j}$ pour  $j=1,...,m$ sur $C_{c}^{\infty}(\tilde{V})$ de sorte que
$$  \vert \Delta J_{\tilde{L}}^{\tilde{G}}(\tilde{\pi}_{\lambda},f)\vert \leq (1+\vert \lambda\vert )^{-N}\sum_{j=1,...,m} {\cal N}_{j}(f) $$
pour tout $\lambda\in i{\cal A}_{\tilde{L}}^*$ et tout $f\in C_{c}^{\infty}(\tilde{V})$. 

{\bf Remarque.} On a rappel\'e en [W2] 5.2 la preuve de cette majoration pour $\Delta=1$. Le m\^eme argument vaut pour tout $\Delta$. 

On applique cela  \`a $f_{\varphi}$ et on utilise (5). On obtient que le terme $J_{\tilde{L}}^{\tilde{G}}(\tilde{\pi}_{\lambda},f_{\varphi})$ v\'erifie une majoration comme dans le (ii) de l'\'enonc\'e (l'alg\`ebre $Diff^{cst}(\tilde{T}({\mathbb R}))^{\omega^{-1}-inv}$ \'etant remplac\'ee par $Diff^{cst}(\tilde{T}({\mathbb R}))^{\omega-inv}$ conform\'ement \`a l'isomorphisme antilin\'eaire que l'on a utilis\'e plus haut).

Soit $\tilde{R}\in {\cal L}(\tilde{L})$, $\tilde{R}\not=\tilde{G}$. Consid\'erons $\phi_{\tilde{R}}^{_{'}\tilde{G}}(f_{\varphi}) $ comme un \'el\'ement de $I(\tilde{R}({\mathbb R}),\omega)$. C'est un \'el\'ement de $I^{\tilde{R}}(\tilde{U}_{reg},\omega)$, l'analogue de $I(\tilde{U}_{reg},\omega)$ pour l'espace ambiant $\tilde{R}$. En raisonnant par r\'ecurrence, on peut  supposer que  l'application
$$\varphi'\mapsto I_{\tilde{L}}^{_{'}\tilde{R}}(\tilde{\pi}_{\lambda}, \varphi')$$
d\'efinie sur $I^{\tilde{R}}(\tilde{U}_{reg},\omega)$ v\'erifie les majorations souhait\'ees. Il intervient dans ces majorations un terme $\delta^{\tilde{R}}(\varphi')$ analogue de $\delta(\varphi')$ quand on remplace $\tilde{G}$ par $\tilde{R}$. Mais, d'apr\`es sa d\'efinition, on a $\delta^{\tilde{R}}(\varphi')\geq c \delta(\varphi')$ pour tout $\varphi'$, avec une constante $c>0$ ind\'ependante de $\varphi'$. 
Pour $X\in \mathfrak{u}$ et $n\in Norm_{G}(\tilde{T};{\mathbb  R})$, on a par d\'efinition
$$\phi_{\tilde{R}}^{_{'}\tilde{G}}(f_{\varphi}) (n^{-1}exp(X)\eta n)=J_{\tilde{R}}^{\tilde{G}}(n^{-1}exp(X)\eta n,\omega,f_{\varphi})$$
$$= \varphi(exp(X)\eta)\omega(n)^{-1}\int_{T^{\theta}({\mathbb R})\backslash G({\mathbb R})}\beta(nx)v_{\tilde{R}}^{\tilde{G}}(x)\,dx.$$
Il en r\'esulte que, pour tout $D\in Diff^{cst}(\tilde{T}({\mathbb R}))^{\omega-inv}$, il existe $c>0$ tel que
$$sup_{\gamma\in \tilde{U}}\vert D\phi_{\tilde{R}}^{_{'}\tilde{G}}(f_{\varphi};\gamma)\vert \leq c\,sup_{\gamma\in \tilde{U}}\vert D\varphi(\gamma)\vert$$
pour tout $\varphi\in I(\tilde{U}_{reg},\omega)$.  On en d\'eduit que le terme $I_{\tilde{L}}^{_{'}\tilde{R}}(\tilde{\pi}_{\lambda},\phi_{\tilde{R}}^{_{'}\tilde{G}}(f_{\varphi}))$ v\'erifie une majoration comme dans le (ii) de l'\'enonc\'e. 
   Alors  l'\'egalit\'e (6) entra\^{\i}ne  une telle majoration pour $ I_{\tilde{L}}^{_{'}\tilde{G}}(\tilde{\pi}_{\lambda},\varphi)$. Cela ach\`eve la d\'emonstration. $\square$

\bigskip

\subsection{Stabilisation de la formule pr\'ec\'edente}
On suppose $(G,\tilde{G},{\bf a})$ quasi-d\'eploy\'e et \`a torsion int\'erieure. 
On consid\`ere un espace de Levi $\tilde{M}\in {\cal L}(\tilde{M}_{0})$, un sous-tore tordu $\tilde{T}$ de $\tilde{M}$ maximal et elliptique et un ensemble $\tilde{U}\subset \tilde{T}({\mathbb R}) $. On suppose $\tilde{U}$ ouvert et d'adh\'erence compacte.    Pour tout espace de Levi $\tilde{L}\in {\cal L}(\tilde{M}_{0})$, on fixe une base $B(\tilde{L})$ de $D_{ell,0}(\tilde{L})$, poss\'edant les propri\'et\'es de [IV] 2.2. En particulier, elle est r\'eunion de bases $B(\tilde{L},\mu)$ des sous-espaces $D_{ell,0,\mu}(\tilde{L})$, pour tout param\`etre infinit\'esimal $\mu$. On suppose de plus que, pour tout $\mu$, $B(\tilde{L},\mu)$ est r\'eunion de bases $B^{st}(\tilde{L},\mu)$ de $D^{st}_{ell,0,\mu}(\tilde{L})$ et $B^{inst}(\tilde{L},\mu)$ de $D^{inst}_{ell,0,\mu}(\tilde{L})$. La base $B(\tilde{L})$ elle-m\^eme est donc r\'eunion de $B^{st}(\tilde{L})$ et $B^{inst}(\tilde{L})$. 

\ass{Proposition}{Pour tout espace de Levi $\tilde{R}\in {\cal L}(\tilde{M}_{0})$ et tout $\tilde{\sigma}\in B^{st}(\tilde{R})$, il existe un ensemble fini $S{\cal H}_{\tilde{R}}^{\tilde{G}}(\tilde{\sigma})$ de sous-espaces affines de $i{\cal A}_{\tilde{R}}^*$ de sorte que les propri\'et\'es suivantes soient v\'erifi\'ees pour tous $\tilde{M}$, $\tilde{T}$, $U$ comme ci-dessus.   Pour tout $\varphi\in  C_{c}^{\infty}(\tilde{U}_{reg})$, il existe une  famille 
$$(S\phi_{\tilde{R},\tilde{\sigma},H}(\varphi))_{\tilde{R}\in {\cal L}(\tilde{M}_{0}),\tilde{\sigma}\in  B^{st}(\tilde{R}),H\in S{\cal H}_{\tilde{R}}^{\tilde{G}}(\tilde{\sigma})}$$
 v\'erifiant les conditions (i) et (ii) de la proposition pr\'ec\'edente ainsi que
 
(iii) pour tout $f\in SI(\tilde{G}({\mathbb R}),K)$, on a l'\'egalit\'e
$$\int_{\tilde{T}({\mathbb R}) }\varphi(\delta)S^{\tilde{G}}_{\tilde{M}}(\delta,f)\,d\delta=\sum_{\tilde{R}\in {\cal L}(\tilde{M}_{0})}\sum_{\tilde{\sigma}\in  B^{st}(\tilde{R})}\sum_{H\in S{\cal H}_{\tilde{R}}^{\tilde{G}}(\tilde{\sigma})}\int_{H}S\phi_{\tilde{R},\tilde{\sigma},H}(\varphi,\lambda)S^{\tilde{R}}(\tilde{\sigma}_{\lambda},f_{\tilde{R}})\,d\lambda.$$}

Preuve.   Fixons   un ensemble de repr\'esentants ${\cal X}$ de l'ensemble de doubles classes 
$$T\backslash\{x\in M; \forall \sigma\in \Gamma_{{\mathbb R}}, x\sigma(x)^{-1}\in T\}/M({\mathbb R}).$$
Pour $\psi\in C_{c}^{\infty}(\tilde{M}({\mathbb R}))$, on a l'\'egalit\'e
$$S^{\tilde{M}}(\delta,\psi)=\sum_{x\in {\cal X}}I^{\tilde{M}}(x^{-1} \delta x,\psi)$$
pour tout $\delta\in \tilde{T}_{\tilde{G}-reg}({\mathbb R})$. Soit $f\in I(\tilde{G}({\mathbb R}),K)$.
  Le membre de gauche de l'\'egalit\'e (iii) est la somme de
$$(1) \qquad \int_{\tilde{T}({\mathbb R}) }\varphi(\delta)\sum_{x\in {\cal X}}I^{\tilde{G}}_{\tilde{M}}(x^{-1}\delta x,f)\,d\delta$$
et de  la somme sur $s\in Z(\hat{M})^{\Gamma_{{\mathbb R}}}/Z(\hat{G})^{\Gamma_{{\mathbb R}}}$, $s\not=1$, de
$$(2) \qquad -i_{\tilde{M}}(\tilde{G},\tilde{G}'(s)) \int_{\tilde{T}({\mathbb R}) }\varphi(\delta)S^{{\bf G}'(s)}_{{\bf M}}(\delta,f^{{\bf G}'(s)})\,d\delta.$$
 Un changement de variables transforme (1) en 
$$ \sum_{x\in {\cal X}}\int_{\tilde{T}({\mathbb R}) }\varphi(x\delta x^{-1}) I^{\tilde{G}}_{\tilde{M}}(\delta,f)\, d\delta. $$
Gr\^ace \`a la proposition pr\'ec\'edente, ce terme  s'\'ecrit sous la forme
$$(3) \qquad   \sum_{\tilde{R}\in {\cal L}(\tilde{M}_{0})}\sum_{\tilde{\sigma}'\in \underline{{\cal E}}_{ell,0}(\tilde{R})}\sum_{H\in {\cal H}_{\tilde{R}}^{\tilde{G}}(\tilde{\sigma}')}\sum_{x\in {\cal X}}\int_{H}\phi_{\tilde{R},\tilde{\sigma}',H}(ad_{x^{-1}}(\varphi),\lambda)I^{\tilde{R}}(\tilde{\sigma}'_{\lambda},f_{\tilde{R}})\,d\lambda.$$ 
Fixons $\tilde{R}$ et introduisons la matrice de changement de base reliant les deux bases $\underline{{\cal E}}_{ell,0}(\tilde{R})$ et $B(\tilde{R})$. C'est-\`a-dire que, pour $\tilde{\sigma}'\in \underline{{\cal E}}_{ell,0}(\tilde{R})$, on \'ecrit $\tilde{\sigma}'=\sum_{\tilde{\sigma}\in B(\tilde{R})}c_{\tilde{\sigma}',\tilde{\sigma}}\tilde{\sigma}$. Pour $\tilde{\sigma}'$ fix\'e, il n'y a qu'un nombre fini de $\tilde{\sigma}$ tels que $c_{\tilde{\sigma}',\tilde{\sigma}}\not=0$. Les propri\'et\'es des bases entra\^{\i}nent que, dans l'autre sens,  pour   $\tilde{\sigma}$ fix\'e, il n'y a qu'un nombre fini de $\tilde{\sigma}'$ tels que $c_{\tilde{\sigma}',\tilde{\sigma}}\not=0$. Pour $\tilde{\sigma}\in B(\tilde{R})$, posons ${\cal H}_{\tilde{R}}^{\tilde{G},0}(\tilde{\sigma})=\cup_{\tilde{\sigma}'\in \underline{{\cal E}}_{ell,0}(\tilde{R}), c_{\tilde{\sigma}',\tilde{\sigma}}\not=0}{\cal H}_{\tilde{R}}^{\tilde{G}}(\tilde{\sigma}')$. Pour $H\in  {\cal H}_{\tilde{R}}^{\tilde{G},0}(\tilde{\sigma})$, posons 
$$\phi_{\tilde{R},\tilde{\sigma},H}^0(\varphi)=\sum_{x\in {\cal X}}\sum_{\tilde{\sigma}'\in \underline{{\cal E}}_{ell,0}(\tilde{R}), H\in {\cal H}_{\tilde{R}}^{\tilde{G}}(\tilde{\sigma}')}c_{\tilde{\sigma}',\tilde{\sigma}}\phi_{\tilde{R},\tilde{\sigma}',H}(ad_{x^{-1}}(\varphi)).$$
On transforme facilement l'expression (3) en
$$(4) \qquad   \sum_{\tilde{R}\in {\cal L}(\tilde{M}_{0})}\sum_{\tilde{\sigma}\in  B(\tilde{R})}\sum_{H\in {\cal H}_{\tilde{R}}^{\tilde{G},0}(\tilde{\sigma})}\int_{H}\phi_{\tilde{R},\tilde{\sigma},H}^0(\varphi,\lambda)I^{\tilde{R}}(\tilde{\sigma}_{\lambda},f_{\tilde{R}})\,d\lambda.$$

  Fixons $s\in Z(\hat{M})^{\Gamma_{{\mathbb R}}}/Z(\hat{G})^{\Gamma_{{\mathbb R}}}$ tel que $s\not=1$ et $i_{\tilde{M}}(\tilde{G},\tilde{G}'(s))\not=0$. Par r\'ecurrence, le terme (2) s'\'ecrit comme une  multiple somme. L'un des indices est un espace de Levi $\tilde{R}'_{s}\in {\cal L}^{\tilde{G}'(s)}(\tilde{M})$. Comme on le sait, un tel espace est  associ\'e \`a une donn\'ee endoscopique elliptique ${\bf R}'(s)$ d'un espace de Levi $\tilde{R}\in {\cal L}(\tilde{M})$. L'espace $\tilde{R}'_{s}$ \'etant fix\'e, la somme int\'erieure est
  $$-i_{\tilde{M}}(\tilde{G},\tilde{G}'(s))\sum_{\tilde{\sigma}'\in B^{st}({\bf R}'(s))}\sum_{H'\in S{\cal H}_{{\bf R}'(s)}^{{\bf G}'(s)}(\tilde{\sigma}')}$$
  $$
 \int_{H'\in S{\cal H}_{{\bf R}'(s)}^{{\bf G}'(s)}(\tilde{\sigma}')}S\phi_{{\bf R}'(s),\tilde{\sigma}',H'}(\varphi,\lambda)S^{{\bf R}'(s)}(\tilde{\sigma}'_{\lambda},(f^{{\bf G}'(s)})_{{\bf R}'(s)})\,d\lambda.$$
Comme toujours, $f^{{\bf G}'(s)}_{{\bf R}'(s)}=(f_{\tilde{R}})^{{\bf R}'(s)}$, d'o\`u
$$S^{{\bf R}'(s)}(\tilde{\sigma}'_{\lambda},(f^{{\bf G}'(s)})_{{\bf R}'(s)})=I^{\tilde{R}}(transfert(\tilde{\sigma}'_{\lambda}),f_{\tilde{R}}).$$
On peut identifier $H'$ \`a un sous-espace affine de $i{\cal A}_{\tilde{R}}^*$. On peut aussi \'ecrire
$$transfert(\tilde{\sigma}'_{\lambda})=\sum_{\tilde{\sigma}\in B(\tilde{R})}c_{\tilde{\sigma}',\tilde{\sigma}}\tilde{\sigma}_{\lambda},$$
avec des coefficients $c_{\tilde{\sigma}',\tilde{\sigma}}$ nuls pour presque tout $\tilde{\sigma}$. L'int\'egrale intervenant ci-dessus est alors
$$\sum_{\tilde{\sigma}\in B(\tilde{R})}c_{\tilde{\sigma}',\tilde{\sigma}}\int_{H'}S\phi_{{\bf R}'(s),\tilde{\sigma}',H'}(\varphi,\lambda)I^{\tilde{R}}(\tilde{\sigma}_{\lambda},f_{\tilde{R}})\,d\lambda.$$
Les propri\'et\'es de finitude du transfert assurent que, pour tout $\tilde{\sigma}$, il n'existe qu'un nombre fini de $\tilde{\sigma}'$ tels que $c_{\tilde{\sigma}',\tilde{\sigma}}\not=0$. Pour $\tilde{\sigma}\in B(\tilde{R})$, posons ${\cal H}_{\tilde{R}}^{\tilde{G},s}(\tilde{\sigma})=\cup_{\tilde{\sigma}'\in  B^{st}({\bf R}'(s)), c_{\tilde{\sigma}',\tilde{\sigma}}\not=0}S{\cal H}_{{\bf R}'(s)}^{{\bf G}'(s)}(\tilde{\sigma}')$. Pour $H\in  {\cal H}_{\tilde{R}}^{\tilde{G},s}(\tilde{\sigma})$, posons 
$$\phi_{\tilde{R},\tilde{\sigma},H}^s(\varphi)=\sum_{\tilde{\sigma}'\in  B^{st}({\bf R}'(s)), H\in S{\cal H}_{{\bf R}'(s)}^{{\bf G}'(s)}(\tilde{\sigma}')}c_{\tilde{\sigma}',\tilde{\sigma}}S\phi_{{\bf R}'(s),\tilde{\sigma}',H}(\varphi).$$
Ces d\'efinitions valent si ${\bf R}'(s)$ est elliptique. Pour $\tilde{R}\in {\cal L}(\tilde{M})$ tel que ${\bf R}'(s)$ n'est pas elliptique, on pose ${\cal H}_{\tilde{R}}^{\tilde{G},s}(\tilde{\sigma})=\emptyset$ pour tout $\tilde{\sigma}\in B(\tilde{R})$. 
Alors l'expression (2) s'\'ecrit
$$(5) \qquad  -i_{\tilde{M}}(\tilde{G},\tilde{G}'(s))\vert \sum_{\tilde{R}\in {\cal L}(\tilde{M}_{0})}\sum_{\tilde{\sigma}\in  B(\tilde{R})}\sum_{H\in {\cal H}_{\tilde{R}}^{\tilde{G},s}(\tilde{\sigma})}\int_{H}\phi_{\tilde{R},\tilde{\sigma},H}^s(\varphi,\lambda)I^{\tilde{R}}(\tilde{\sigma}_{\lambda},f_{\tilde{R}})\,d\lambda.$$
On peut accro\^{\i}tre les ensembles ${\cal H}_{\tilde{R}}^{\tilde{G},0}(\tilde{\sigma})$ et ${\cal H}_{\tilde{R}}^{\tilde{G},s}(\tilde{\sigma})$ intervenant dans (4) et (5), quitte \`a introduire des fonctions nulles $\phi_{\tilde{R},\tilde{\sigma},H}^0(\varphi)$ ou $\phi_{\tilde{R},\tilde{\sigma},H}^s(\varphi)$.  Dans un premier temps, on peut donc remplacer chacun de ces ensembles par leur r\'eunion. On obtient ainsi un ensemble qui, a priori, d\'epend de $\tilde{M}$. On peut de nouveau accro\^{\i}tre cet ensemble en le rempla\c{c}ant par la r\'eunion sur les $\tilde{M}\in {\cal L}^{\tilde{R}}(\tilde{M}_{0})$ des ensembles relatifs \`a $\tilde{M}$. On obtient un ensemble $S{\cal H}_{\tilde{R}}^{\tilde{G}}(\tilde{\sigma})$ qui ne d\'epend plus de $\tilde{M}$. On peut  supposer que, 
 pour $\tilde{R}\in {\cal L}(\tilde{M})$ et $\tilde{\sigma}\in B(\tilde{R})$, les ensembles    ${\cal H}_{\tilde{R}}^{\tilde{G},0}(\tilde{\sigma})$ et  ${\cal H}_{\tilde{R}}^{\tilde{G},s}(\tilde{\sigma})$  intervenant sont tous \'egaux \`a $S{\cal H}_{\tilde{R}}^{\tilde{G}}(\tilde{\sigma})$.   Pour $H\in S{\cal H}_{\tilde{R}}^{\tilde{G}}(\tilde{\sigma})$, posons 
$$S\phi_{\tilde{R},\tilde{\sigma},H}(\varphi)= \phi_{\tilde{R},\tilde{\sigma},H}^0(\varphi)-\sum_{s\in Z(\hat{M})^{\Gamma_{{\mathbb R}}}/Z(\hat{G})^{\Gamma_{{\mathbb R}}},s\not=1}i_{\tilde{M}}(\tilde{G},\tilde{G}'(s))\phi_{\tilde{R},\tilde{\sigma},H}^s(\varphi).$$
 On obtient que le membre de gauche de l'\'egalit\'e (iii) de l'\'enonc\'e est \'egal \`a
$$(6) \qquad \sum_{\tilde{R}\in {\cal L}(\tilde{M}_{0})}\sum_{\tilde{\sigma}\in  B(\tilde{R})}\sum_{H\in S{\cal H}_{\tilde{R}}^{\tilde{G}}(\tilde{\sigma})}\int_{H}S\phi_{\tilde{R},\tilde{\sigma},H}(\varphi,\lambda)I^{\tilde{R}}(\tilde{\sigma}_{\lambda},f_{\tilde{R}})\,d\lambda.$$
D'apr\`es leur construction, les fonctions $\phi_{\tilde{R},\tilde{\sigma},H}(\varphi)$ v\'erifient les propri\'et\'es (i) et (ii) de la proposition 8.3. On doit remarquer que, dans les majorations de $\phi_{\tilde{R},\tilde{\sigma},H}^s(\varphi)$, il intervient un terme $\delta^{\tilde{G}'(s)}(\varphi)$ analogue \`a $\delta(\varphi)$ quand on remplace $\tilde{G}$ par $\tilde{G}'(s)$. Mais ce terme est minor\'e par le produit de $\delta(\varphi)$ et d'une constante ind\'ependante de $\varphi$.

   On d\'ecompose l'expression  (6) en une somme de deux expressions $J^{st}(\varphi,f)$ et $J^{inst}(\varphi,f)$. Dans la premi\`ere,  resp. la seconde, on regroupe les contributions des $\tilde{\sigma}\in B^{st}(\tilde{R})$, resp. $\tilde{\sigma}\in B^{inst}(\tilde{R})$. Pour $\tilde{\sigma}\in B^{st}(\tilde{R})$, on a par d\'efinition $I^{\tilde{R}}(\tilde{\sigma}_{\lambda},f_{\tilde{R}})=S^{\tilde{R}}(\tilde{\sigma}_{\lambda},f_{\tilde{R}})$. Donc $J^{st}(\varphi,f)$ est de la forme du membre de droite du (iii) de l'\'enonc\'e. Pour d\'emontrer la proposition, il reste \`a prouver que $J^{inst}(\varphi,f)=0$. 
   On a d\'ecompos\'e l'espace de Paley-Wiener en $PW^{st,\infty}(\tilde{G})\oplus PW^{inst,\infty}(\tilde{G})$, cf. [IV] 2.4.  Il r\'esulte des d\'efinitions que $J^{st}(\varphi,f)$ ne d\'epend que de la projection de $pw(f)$ sur $PW^{st,\infty}(\tilde{G})$, tandis que $J^{inst}(\varphi,f)$ ne d\'epend que de la projection de $pw(f)$ sur $PW^{inst,\infty}(\tilde{G})$. Mais on sait que les int\'egrales orbitales pond\'er\'ees $f\mapsto S_{\tilde{M}}^{\tilde{G}}(\delta,f)$ sont stables ([V] th\'eor\`eme 1.4). Donc le membre de gauche de (iii) est stable. Soit $f'$ l'\'el\'ement de $I(\tilde{G}({\mathbb R}))$ tel que $pw(f')$ soit \'egal \`a la projection de $pw(f)$ sur $PW^{st,\infty}(\tilde{G})$.  La fonction $f'$ est encore $K$-finie. Le membre de gauche pour $f$ est alors \'egal au m\^eme membre pour $f'$. Puisque ces termes sont calcul\'es par la formule (6), on obtient 
   $$J^{st}(\varphi,f)+J^{inst}(\varphi,f)=J^{st}(\varphi,f')+J^{inst}(\varphi,f').$$
   Par construction de $f'$, on a $J^{st}(\varphi,f)=J^{st}(\varphi,f')$ tandis que $J^{inst}(\varphi,f')=0$. Cela entra\^{\i}ne $J^{inst}(\varphi,f)=0$. Comme on l'a dit, cela prouve la proposition. 
   $\square$
 
 \bigskip
 
 \subsection{Version endoscopique de la proposition 8.4}
 On consid\`ere ici un triplet $(KG,K\tilde{G},{\bf a})$   et un $K$-espace $K\tilde{M}\in {\cal L}(K\tilde{M}_{0})$. On fixe une composante connexe $\tilde{M}$ de $K\tilde{M}$, contenue dans une composante $\tilde{G}$ de $K\tilde{G}$,  un sous-tore tordu maximal elliptique $\tilde{T}$ de $\tilde{M}$ et un ensemble $\tilde{U}\subset \tilde{T}({\mathbb R}) $ v\'erifiant les m\^emes propri\'et\'es qu'en 8.4. 
\ass{Proposition}{Pour tout $K\tilde{R}\in {\cal L}(K\tilde{M}_{0})$ et tout $\tilde{\sigma}\in \underline{{\cal E}}_{ell,0}(K\tilde{R},\omega)$, il existe un ensemble fini ${\cal H}_{K\tilde{R}}^{K\tilde{G},{\cal E}}(\tilde{\sigma})$ de sous-espaces affines de $i{\cal A}_{\tilde{R}}^*$ v\'erifiant les conditions suivantes.

(1) Pour $w\in W(K\tilde{M}_{0})$, $w({\cal H}_{K\tilde{R}}^{K\tilde{G},{\cal E}}(\tilde{\sigma}))={\cal H}_{w(K\tilde{R})}^{K\tilde{G},{\cal E}}(\tilde{\sigma}_{w})$,  o\`u on a pos\'e  $w(\tilde{\sigma})=z_{w}\tilde{\sigma}_{w}$ comme en  8.4.

(2) Pour tous $\tilde{M}$, $\tilde{T}$, $U$ comme ci-dessus et pour tout $\varphi\in C_{c}^{\infty}(\tilde{U}_{reg})^{\omega^{-1}-inv}$, il existe une unique famille 
$$(\phi_{K\tilde{R},\tilde{\sigma},H}^{{\cal E}}(\varphi))_{K\tilde{R}\in {\cal L}(K\tilde{M}_{0}),\tilde{\sigma}\in \underline{{\cal E}}_{ell,0}(K\tilde{R},\omega),H\in {\cal H}_{K\tilde{R}}^{K\tilde{G},{\cal E}}(\tilde{\sigma})}$$
 v\'erifiant les conditions  (i), (ii) et (iii) de la proposition 8.4 ainsi que

(iv) pour tout $f\in I(K\tilde{G}({\mathbb R}),\omega,K)$, on a l'\'egalit\'e
$$\int_{\tilde{T}({\mathbb R})/(1-\theta)(T({\mathbb R}))}\varphi(\gamma)I^{K\tilde{G},{\cal E}}_{K\tilde{M}}(\gamma,\omega,f)\,d\gamma=\sum_{K\tilde{R}\in {\cal L}(K\tilde{M}_{0})}\sum_{\tilde{\sigma}\in \underline{{\cal E}}_{ell,0}(K\tilde{R},\omega)}\sum_{H\in {\cal H}_{K\tilde{R}}^{K\tilde{G},{\cal E}}(\tilde{\sigma})}$$
$$\int_{H}\phi_{K\tilde{R},\tilde{\sigma},H}^{{\cal E}}(\varphi,\lambda)I^{K\tilde{R}}(\tilde{\sigma}_{\lambda},f_{K\tilde{R},\omega})\,d\lambda.$$}

Preuve. On affirme l'existence de familles d'espaces affines ind\'ependantes de $\tilde{M}$ et $\tilde{T}$. Mais, si on d\'emontre pour chaque couple $(\tilde{M},\tilde{T})$ l'existence d'une telle famille v\'erifiant les conditions requises pour ce couple, il suffit de prendre la r\'eunion de ces familles sur l'ensemble des couples $(\tilde{M},\tilde{T})$ pris \`a conjugaison pr\`es pour r\'esoudre le probl\`eme.   De m\^eme, on peut n\'egliger les conditions d'invariance par $W(\tilde{M}_{0})$: si on r\'esout le probl\`eme sans ces conditions, on peut ensuite accro\^{\i}tre les familles d'espaces affines de sorte \`a les rendre sym\'etriques et sym\'etriser les familles 
$$(\phi_{K\tilde{R},\tilde{\sigma},H}^{{\cal E}}(\varphi))_{K\tilde{R}\in {\cal L}(K\tilde{M}_{0}),\tilde{\sigma}\in \underline{{\cal E}}_{ell,0}(K\tilde{R},\omega),H\in {\cal H}_{K\tilde{R}}^{K\tilde{G},{\cal E}}(\tilde{\sigma})}.$$

Fixons donc $\tilde{M}$, $\tilde{T}$ et $\tilde{U}$. On peut supposer  comme dans la preuve de 8.3 que $\tilde{U}= \{x^{-1} exp(X)\eta x; X\in \mathfrak{u},x\in Norm_{G}(\tilde{T};{\mathbb R})\}$, o\`u $\eta$ est un \'el\'ement fix\'e de $\tilde{T}({\mathbb R}) $ et $\mathfrak{u}$ est un voisinage de $0$ dans $\mathfrak{t}^{\theta}({\mathbb R})$, assez petit.  Soit $f\in I(K\tilde{G}({\mathbb R}),\omega,K)$ et $\varphi\in C_{c}^{\infty}(\tilde{U}_{reg})^{\omega^{-1}-inv}$. Le membre de gauche de (iv) est \'egal \`a
$$\int_{\mathfrak{u}}\varphi(exp(X)\eta)I^{K\tilde{G},{\cal E}}_{K\tilde{M}}(exp(X)\eta,\omega,f)\,dX.$$
Ici comme dans la suite, les mesures doivent \^etre convenablement normalis\'ees pour que la formule soit exacte, mais ces normalisations importent peu puisqu'on ne se propose pas de calculer explicitement les fonctions $\phi_{K\tilde{R},\tilde{\sigma},H}^{{\cal E}}(\varphi)$. On reprend les constructions et notations de 4.4. La formule (1) de ce paragraphe entra\^{\i}ne l'\'egalit\'e
$$I^{K\tilde{G},{\cal E}}_{K\tilde{M}}(exp(X)\eta,\omega,f)=c\sum_{j\in J}\Delta_{j,1}(exp(\xi_{j}(X))\epsilon_{j,1},exp(X)\eta)^{-1}$$
$$S_{\tilde{M}'_{j,1},\lambda_{j,1}}^{\tilde{G}'_{j,1}}(exp(\xi_{j}(X))\epsilon_{j,1},f^{\tilde{G}'_{j,1}}),$$
o\`u $c$ est une constante non nulle. Pour chaque $j\in J$, posons $\mathfrak{u}_{j}=\xi_{j}(\mathfrak{u})$. Le membre de gauche de (iv) est donc
$$\sum_{j\in J}\int_{\mathfrak{u}_{j}}\varphi_{j}(exp(Y)\epsilon_{j,1})S_{\tilde{M}'_{j,1},\lambda_{j,1}}^{\tilde{G}'_{j,1}}(exp(Y))\epsilon_{j,1},f^{\tilde{G}'_{j,1}})\,dY,$$
o\`u
$$\varphi_{j}(exp(Y)\epsilon_{j,1})=c\Delta_{j,1}(exp(Y)\epsilon_{j,1},exp(\xi_{j}^{-1}(Y))\eta)^{-1}\varphi(exp(\xi_{j}^{-1}(Y))\eta).$$
Pour chaque $j$, on peut appliquer la proposition 8.5. Certes, il faut l'adapter \`a la situation des donn\'ees auxiliaires o\`u l'on consid\`ere des fonctions se transformant selon le caract\`ere $\lambda_{j,1}$ de $C_{1}({\mathbb R})$.  On laisse ces d\'etails techniques. Remarquons que la fonction $\varphi_{j}$ est $C^{\infty}$ car le facteur de transfert l'est sur l'ensemble des \'el\'ements fortement r\'eguliers.  Le r\'esultat est que le terme index\'e par $j$ dans la somme ci-dessus s'exprime sous la forme
$$\sum_{\tilde{R}'_{j}\in {\cal L}^{\tilde{G}'_{j}}(\tilde{M}'_{j})}\sum_{\tilde{\sigma}'\in B^{st}({\bf R}'_{j})}\sum_{H\in S{\cal H}_{\tilde{R}'_{j}}^{\tilde{G}'_{j}}(\tilde{\sigma}')}\int_{H}S\phi_{{\bf R}',\tilde{\sigma}',H}(\varphi_{j},\lambda)S^{{\bf R}'_{j}}(\tilde{\sigma}'_{\lambda},(f^{{\bf G}'_{j}})_{{\bf R}'_{j}})\,d\lambda.$$
On a utilis\'e la notation ${\bf R}'_{j}$: comme toujours, l'espace $\tilde{R}'_{j}$ appara\^{\i}t comme l'espace d'une donn\'ee endoscopique ${\bf R}'_{j}$ d'un $K$-espace de Levi $K\tilde{R}\in {\cal L}(K\tilde{M}_{0})$. La d\'emonstration se poursuit alors comme celle de 8.5. Fixons  $\tilde{R}'_{j}$ intervenant ci-dessus. On identifie ${\cal A}_{\tilde{R}'_{j}}^*$ \`a ${\cal A}_{\tilde{R}}^*$. Pour chaque $\tilde{\sigma}'$ intervenant, l'ensemble $S{\cal H}_{\tilde{R}'_{j}}^{\tilde{G}'_{j}}(\tilde{\sigma}')$ s'identifie \`a un ensemble de sous-espaces affines de $i{\cal A}_{\tilde{R}}$. On \'ecrit
$$transfert(\tilde{\sigma}'_{\lambda})=\sum_{\tilde{\sigma}\in \underline{{\cal E}}(K\tilde{R},\omega)}c_{\tilde{\sigma}',\tilde{\sigma}}\tilde{\sigma}_{\lambda} .$$
Alors
$$\int_{H}S\phi_{{\bf R}',\tilde{\sigma}',H}(\varphi_{j},\lambda)S^{{\bf R}'_{j}}(\tilde{\sigma}'_{\lambda},(f^{{\bf G}'_{j}})_{{\bf R}'_{j}})\,d\lambda=\sum_{\tilde{\sigma}\in \underline{{\cal E}}(K\tilde{R},\omega)}c_{\tilde{\sigma}',\tilde{\sigma}}$$
$$\int_{H}S\phi_{{\bf R}',\tilde{\sigma}',H}(\varphi_{j},\lambda) I^{K\tilde{R}}(\tilde{\sigma}_{\lambda}, f_{K\tilde{R},\omega})\,d\lambda.$$
Pour chaque $\tilde{\sigma}$, il n'y a qu'un nombre fini de $\tilde{\sigma}'$ pour lesquels $c_{\tilde{\sigma}',\tilde{\sigma}}\not=0$. En sommant les expressions obtenues, 
on obtient finalement une expression du membre de gauche de (iv) de la forme voulue. Les familles de sous-espaces affines d\'ependent des donn\'ees ${\bf M}'_{j}$ choisies, mais on peut supposer que celles-ci appartiennent \`a un ensemble fini de repr\'esentants des classes d'\'equivalence de donn\'ees endoscopiques elliptiques pour $(KM,K\tilde{M},{\bf a})$. Donc ces sous-espaces affines restent dans un ensemble fini. 

Parce que le facteur de transfert $\Delta_{j,1}(exp(Y)\epsilon_{j,1},exp(\xi_{j}^{-1}(Y))\eta)$ est $C^{\infty}$ en $Y$, on voit facilement que les fonctions $S\phi_{{\bf R}',\tilde{\sigma}',H}(\varphi_{j},\lambda)$ v\'erifient les majorations souhait\'ees. 
Cela ach\`eve la d\'emonstration. $\square$

\bigskip

\subsection{Expression de $\epsilon_{K\tilde{M}}(f)$}
La proposition 8.4 s'adapte \'evidemment aux $K$-espaces. Puisque l'on peut toujours accro\^{\i}tre nos familles de sous-espaces affines, on peut supposer que, pour tout $K\tilde{R}\in {\cal L}(K\tilde{M}_{0})$, toute composante connexe $\tilde{R}$ de $K\tilde{R}$,  tout $\tilde{\sigma}\in \underline{{\cal E}}_{ell,0}(\tilde{R},\omega)$, l'ensemble ${\cal H}_{\tilde{R}}^{\tilde{G}}(\tilde{\sigma})$ de 8.4 est contenu dans l'ensemble ${\cal H}_{K\tilde{R}}^{K\tilde{G},{\cal E}}(\tilde{\sigma})$ de 8.6. Pour $\tilde{M}$, $\tilde{T}$ et $\tilde{U}$ comme en 8.6, la conjonction des propositions 8.4 et 8.6 entra\^{\i}ne le r\'esultat suivant.  Pour tout $\varphi\in C_{c}^{\infty}(\tilde{U}_{reg})^{\omega^{-1}-inv}$, il existe une unique famille 
$$(\phi_{K\tilde{R},\tilde{\sigma},H}(\varphi))_{K\tilde{R}\in {\cal L}(K\tilde{M}_{0}),\tilde{\sigma}\in \underline{{\cal E}}_{ell,0}(K\tilde{R},\omega),H\in {\cal H}_{\tilde{R}}^{\tilde{G}}(\tilde{\sigma})}$$
 v\'erifiant les conditions (i), (ii) et (iii) de 8.4 ainsi que:

(1) pour tout $f\in I(K\tilde{G}({\mathbb R}),\omega,K)$, on a l'\'egalit\'e
$$\int_{\tilde{T}({\mathbb R})/(1-\theta)(T({\mathbb R}))}\varphi(\gamma)I^{K\tilde{M}}(\gamma,\omega,\epsilon_{K\tilde{M}}(f))\,d\gamma=\sum_{K\tilde{R}\in {\cal L}(\tilde{M}_{0})}\sum_{\tilde{\sigma}\in \underline{{\cal E}}_{ell,0}(K\tilde{R},\omega)}\sum_{H\in {\cal H}_{K\tilde{R}}^{K\tilde{G},{\cal E}}(\tilde{\sigma})}$$
$$\int_{H}\phi_{K\tilde{R},\tilde{\sigma},H}(\varphi,\lambda)I^{K\tilde{R}}(\tilde{\sigma}_{\lambda},f_{K\tilde{R},\omega})\,d\lambda.$$

\ass{Proposition}{Pour tout $K\tilde{R}\in {\cal L}(K\tilde{M}_{0})$, tout $\tilde{\sigma}\in \underline{{\cal E}}_{ell,0}(K\tilde{R},\omega)$, tout $H\in {\cal H}_{K\tilde{R}}^{K\tilde{G},{\cal E}}(\tilde{\sigma})$, il existe une unique fonction $\xi_{K\tilde{R},\tilde{\sigma},H}$ sur $(\tilde{T}({\mathbb R})\cap \tilde{G}_{reg}({\mathbb R}))\times H$ de sorte que les conditions suivantes soient v\'erifi\'ees. Dans les deux premi\`eres, on fixe $K\tilde{R}$, $\tilde{\sigma}$ et $H$.

(i) La fonction $\xi_{K\tilde{R},\tilde{\sigma},H}$ est $C^{\infty}$ sur $(\tilde{T}({\mathbb R})\cap \tilde{G}_{reg}({\mathbb R}))\times H$;

(ii) Pour tout op\'erateur diff\'erentiel  $D \in Diff^{cst}(\tilde{T}({\mathbb R}))^{\omega-inv}$, tout op\'erateur diff\'erentiel $\Delta$ \`a coefficients constants sur $H$ et tout sous-ensemble compact $\Gamma\subset \tilde{T}({\mathbb R})$, il existe $d, N\in {\mathbb N}$ et une constante $c>0$  de sorte que
$$\vert D\Delta \xi_{K\tilde{R},\tilde{\sigma},H}(\gamma,\lambda)\vert \leq c\vert  D^{\tilde{G}}(\gamma)\vert ^{-d}(1+\vert \lambda\vert )^N$$
pour tout $\gamma\in  \Gamma\cap \tilde{G}_{reg}({\mathbb R})$ et tout $\lambda\in H$.

(iii) La famille v\'erifie la condition de sym\'etrie (2) de 8.4.

(iv) Pour  tout $f\in I(K\tilde{G}({\mathbb R}),\omega,K)$ et tout $\gamma\in \tilde{T}({\mathbb R})\cap \tilde{G}_{reg}({\mathbb R})$ , on a l'\'egalit\'e
$$I^{K\tilde{M}}(\gamma,\omega,\epsilon_{K\tilde{M}}(f))=\sum_{K\tilde{R}\in {\cal L}(K\tilde{M}_{0})}\sum_{\tilde{\sigma}\in \underline{{\cal E}}_{ell,0}(\tilde{R},\omega)}\sum_{H\in {\cal H}_{K\tilde{R}}^{K\tilde{G},{\cal E}}(\tilde{\sigma})}$$
$$\int_{H}\xi_{K\tilde{R},\tilde{\sigma},H}(\gamma,\lambda)I^{K\tilde{R}}(\tilde{\sigma}_{\lambda},f_{K\tilde{R},\omega})\,d\lambda.$$}

Preuve. La preuve d'Arthur [A9] p. 198, 199, 200  s'applique. On va la reprendre car Arthur \'enonce une majoration plus faible que celle du (ii) de l'\'enonc\'e.   On a muni l'espace $\mathfrak{t}^{\theta}({\mathbb R})$ est d'une norme  euclidienne. Notons $\mathfrak{u}'$ une boule  ouverte centr\'ee en $0$ de rayon $c_{0}$ assez petit et $\mathfrak{u}$ la boule de rayon $2c_{0}$.    Fixons  $\eta\in \tilde{T}({\mathbb R})$ et  posons $\tilde{U}=\{(1-\theta)(t)exp(X)\eta; t\in T({\mathbb R}), X\in \mathfrak{u}\}$. Notons $\mathfrak{u}_{reg}$, resp. $\mathfrak{u}'_{reg}$, le sous-ensemble des $X\in \mathfrak{u}$, resp. $X\in \mathfrak{u}'$,  tels que $exp(X)\eta$ soit fortement r\'egulier dans $\tilde{G}$. Avec les notations de 8.3, $\mathfrak{u}_{reg}$ est  l'ensemble des $X\in \mathfrak{u}$ dont le fixateur dans $\Xi$ est trivial, pourvu que $c_{0}$ soit assez petit.  La propri\'et\'e suivante est claire:

(2) il existe $c_{1}>0$ tel que, pour tout $X\in \mathfrak{u}'_{reg}$ et pour tout $Y\in \mathfrak{t}({\mathbb R})$, la condition $\vert X-Y\vert \leq c_{1}D^{\tilde{G}}(exp(X)\eta)$ implique que $Y\in \mathfrak{u}_{reg}$ et  que $D^{\tilde{G}}(exp(Y)\eta)\geq D^{\tilde{G}}(exp(X)\eta)/2$.

  L'espace $C_{c}^{\infty}(\tilde{U}_{reg})^{\omega^{-1}-inv}$ s'identifie \`a $C_{c}^{\infty}(\mathfrak{u}_{reg})$. Les fonctions $\phi_{K\tilde{R},\tilde{\sigma},H}(\varphi)$ de (1) sont donc d\'efinies pour $\varphi\in C_{c}^{\infty}(\mathfrak{u}_{reg})$. Pour tout $m\in {\mathbb N}$, notons $C_{c}^m(\mathfrak{u}_{reg})$ l'espace des fonctions de classe $C^m$ sur $\mathfrak{t}^{\theta}({\mathbb R})$ dont le support est compact et inclus dans $\mathfrak{u}_{reg}$. Consid\'erons une suite $(\varphi_{i})_{i\in {\mathbb N}}$ d'\'el\'ements de $C_{c}^m(\mathfrak{u}_{reg})$. On dit qu'elle tend vers $0$ s'il existe un sous-ensemble compact de $\mathfrak{u}_{reg}$ tel que le support de $\varphi_{i}$ soit inclus dans ce compact pour tout $i$ et si, pour tout $H\in Sym(\mathfrak{t}^{\theta})$ de degr\'e inf\'erieur ou \'egal \`a $m$, la suite $(\partial_{H}\varphi_{i})_{i\in {\mathbb N}}$ tend uniform\'ement vers $0$. Cela munit $C_{c}^m(\mathfrak{u}_{reg})$ d'une topologie. L'application $\varphi\mapsto \delta(\varphi)$ de 8.3 se d\'efinit aussi bien pour $\varphi\in C_{c}^m(\mathfrak{u}_{reg})$. 
     Pour $H$ intervenant dans (1), notons $C^m(H)$ l'espace des fonctions de classe $C^m$ sur $H$.  On munit cet espace de la topologie usuelle. Montrons que
  
  (3) pour tout $m_{0}\in {\mathbb N}$, il existe un entier $m_{1}\in {\mathbb N}$ tel que, pour tous $K\tilde{R}$, $\tilde{\sigma}$, $H$ intervenant dans (1), l'application $\varphi\mapsto \phi_{K\tilde{R},\tilde{\sigma},H}(\varphi)$ s'\'etend par continuit\'e en une application de $C^{m_{1}}(\mathfrak{u}_{reg})$ dans $C^{m_{0}}(H)$;   pour tout op\'erateur diff\'erentiel \`a coefficients constants $\Delta$ sur $H$ de degr\'e $\leq m_{0}$, il existe  une famille finie $(H_{i})_{i=1,...,n}$  d'\'el\'ements  de $  Sym(\mathfrak{t}^{\theta})$ de degr\'e $\leq m_{1}$ et un entier $d\in {\mathbb N}$    de sorte que
  $$\vert \Delta\phi_{K\tilde{R},\tilde{\sigma},H}(\varphi,\lambda)\vert \leq c\delta(\varphi)^{-d}\sum_{i=1,...,n}sup_{X\in \mathfrak{u}}\vert \partial_{H_{i}}\varphi(X)\vert $$
  pour tout $\lambda\in H$ et tout $\varphi\in C^{m_{1}}(\mathfrak{u}_{reg})$; ces applications v\'erifient la propri\'et\'e de sym\'etrie 8.4(3) et l'\'egalit\'e (1) ci-dessus;
    
Pour tout $H$, on  fixe une base $(\Delta_{H,j})_{j=1,...,j_{H}}$ de l'espace des op\'erateurs diff\'erentiels \`a coefficients constants sur $H$ de degr\'e $\leq m_{0}$. En appliquant les majorations d\'ej\`a prouv\'ees, on voit que l'on peut trouver une
famille finie $(H_{i})_{i=1,...,n}$  d'\'el\'ements  de $  Sym(\mathfrak{t}^{\theta})$, un entier $d\in {\mathbb N}$ et une constante $c>0$ de sorte que
  $$(4) \qquad \vert \Delta_{H,j}\phi_{K\tilde{R},\tilde{\sigma},H}(\varphi,\lambda)\vert \leq c\delta(\varphi)^{-d}\sum_{i=1,...,n}sup_{X\in \mathfrak{u}}\vert \partial_{H_{i}}\varphi(X)\vert $$
  pour tous $K\tilde{R}$, $\tilde{\sigma}$, $H$, tout $j=1,...,j_{H}$, tout $\lambda\in H$ et tout $\varphi\in C_{c}^{\infty}(\mathfrak{u}_{reg})$.  Notons $m_{1}$ le plus grand des degr\'es des $H_{i}$. Quand une suite $(\varphi_{k})_{k\in {\mathbb N}}$ d'\'el\'ements de $C_{c}^{\infty}(\mathfrak{u}_{reg})$ tend vers une fonction $\varphi\in C_{c}^{m_{1}}(\mathfrak{u}_{reg})$, les termes $\delta(\varphi_{k})$ sont minor\'es puisque les supports des $\varphi_{k}$ restent dans un compact de $\mathfrak{u}_{reg}$. Il r\'esulte alors de la majoration (4) que la suite $(\phi_{K\tilde{R},\tilde{\sigma},H}(\varphi_{k}))_{k\in {\mathbb N}}$ est de Cauchy  dans $C^{m_{0}}(H)$.  Elle converge vers un \'el\'ement de cet espace que l'on note $\phi_{K\tilde{R},\tilde{\sigma},H}(\varphi)$.  Les propri\'et\'es de cette application r\'esultent par continuit\'e des propri\'et\'es de l'application initiale d\'efinie pour $\varphi\in C_{c}^{\infty}(\mathfrak{u}_{reg})$. On doit juste pr\'eciser un point. Pour une suite $(\varphi_{k})_{k\in {\mathbb N}}$ tendant vers  une fonction $\varphi$ comme ci-dessus, il n'y a pas de raison que $\delta(\varphi_{k})$ tende vers $\delta(\varphi)$. La majoration voulue de $\Delta_{H,j}\phi_{K\tilde{R},\tilde{\sigma},H}(\varphi,\lambda)$ ne r\'esulte pas des majorations (4) des $\Delta_{H,j}\phi_{K\tilde{R},\tilde{\sigma},H}(\varphi_{k},\lambda)$. Mais, pour tout $\varphi\in C^{m_{1}}_{c}(\mathfrak{u}_{reg})$, on peut choisir une suite $(\varphi_{k})_{k\in {\mathbb N}}$ d'\'el\'ements de $C_{c}^{\infty}(\mathfrak{u}_{reg})$ tendant vers $\varphi$ et v\'erifiant de plus la propri\'et\'e $\delta(\varphi_{k})\geq \delta(\varphi)/2$ pour tout $k$. En utilisant une telle suite, on obtient la majoration voulue. Cela prouve (3). 
  
  Notons $\varpi$ l'\'el\'ement de $Sym(\mathfrak{t}^{\theta})$ tel que $\partial_{\varpi}$ soit le laplacien relatif \`a notre m\'etrique. Puisque $Sym(\mathfrak{t}^{\theta})$ est un module de type fini fini sur l'image de $\mathfrak{Z}(G)$ par l'application d'Harish-Chandra $z\mapsto z_{T^{\theta}}$, on a
  
  (5) il existe un entier $r\geq1$ tel que, pour tout entier $m\geq1$, il existe des \'el\'ements $z_{j}$ de $\mathfrak{Z}(G)$ pour $j=1,...,r$ de sorte que
  $$\varpi^{mr}=\sum_{j=1,...,r}z_{j,T^{\theta}}\varpi^{m(r-j)}.$$
  
  Notons $\delta_{0}$ la mesure de Dirac en $0$ sur $\mathfrak{t}^{\theta}({\mathbb R})$. La th\'eorie des op\'erateurs elliptiques implique que
  
  (6) il existe $m_{2}\in {\mathbb N}$ tel que, pour tout $m\geq m_{2}/2$, il existe une fonction $\psi^m$ de classe $2m-m_{2}$ sur $\mathfrak{t}^{\theta}({\mathbb R})$ de sorte que
  $$\delta_{0}=\partial_{\varpi^m}\psi^m.$$
  
  Fixons une fonction $\alpha\in C_{c}^{\infty}(\mathfrak{u})$ qui soit constante de valeur $1$ au voisinage de $0$. Pour $\epsilon>0$  et $m\in {\mathbb N}$, posons $\psi^m_{\epsilon}(X)=\psi^m(X)\alpha(\epsilon^{-1}X)$. La fonction $\psi^m_{\epsilon}$ est de classe $2m-m_{2}$ et est \`a support dans la boule centr\'ee en $0$ et de rayon $2\epsilon c_{0}$. On a  une \'egalit\'e
 $$\delta_{0}=\partial_{\varpi^m}\psi_{\epsilon}^m+\chi^m_{\epsilon}$$
 o\`u $\chi^m_{\epsilon}$ est une fonction $C^{\infty}$ sur $\mathfrak{t}^{\theta}({\mathbb R})$ \`a support dans la boule centr\'ee en $0$ et de rayon $2\epsilon c_{0}$.
 
 Fixons maintenant un \'el\'ement $X\in \mathfrak{u}'_{reg}$ et un entier $m_{0}\in {\mathbb N}$.  On pose $\gamma=exp(X)\eta$. Fixons un $\epsilon>0$ tel que 
 $$\epsilon<\frac{c_{1}}{2c_{0}}D^{\tilde{G}}(\gamma)$$
  o\`u $c_{1}$ est comme en (2). Fixons un entier $m\in {\mathbb N}$ tel que $2m>m_{0}+m_{1}+m_{2}$ o\`u $m_{1}$ et $m_{2}$ sont d\'etermin\'es par (3) et (6). Posons $\psi=\psi_{\epsilon}^{mr}$ et $\chi=\chi_{\epsilon}^{mr}$. Pour toute fonction $\beta$ sur $\mathfrak{t}^{\theta}({\mathbb R})$ et tout $X'\in \mathfrak{t}^{\theta}({\mathbb R})$, d\'efinissons $\beta_{X'}$ par $\beta_{X'}(Y)=\beta(Y-X')$ pour tout $Y\in \mathfrak{t}^{\theta}({\mathbb R})$. 
   La propri\'et\'e (2) et le choix de $\epsilon$ entra\^{\i}nent que $\psi_{X}$ et $\chi_{X}$ sont \`a support dans l'ensemble des \'el\'ements  $Y\in \mathfrak{u}_{reg}$ qui v\'erifient l'in\'egalit\'e $D^{\tilde{G}}(exp(Y)\eta)\geq D^{\tilde{G}}(\gamma)/2$. La fonction $\chi_{X}$ est $C^{\infty}$ tandis que $\psi_{X}$ est de classe $2mr-m_{2}$. En notant $\delta_{X}$ la mesure de Dirac en $X$, on a 
 $$\delta_{X}=\partial_{\varpi^{mr}}\psi_{X}+\chi_{X}.$$
 Ou encore, en utilisant (5),
 $$\delta_{X}=\chi_{X}+\sum_{j=1,...,r}\partial_{z_{j,T^{\theta}}\varpi^{m(r-j)}}\psi_{X}.$$ 
 Soit $f\in I(K\tilde{G}({\mathbb R}),\omega)$. Appliquons $\delta_{X}$ \`a la fonction $Y\mapsto I^{K\tilde{M}}(exp(Y)\eta,\omega,\epsilon_{K\tilde{M}}(f))$. On obtient $I^{K\tilde{M}}(\gamma,\omega,\epsilon_{K\tilde{M}}(f))$. En utilisant la formule pr\'ec\'edente, on en d\'eduit
 $$(7) \qquad I^{K\tilde{M}}(\gamma,\omega,\epsilon_{K\tilde{M}}(f))=\sum_{j=0,...,r}I_{j},$$
 o\`u
 $$I_{0}=\int_{\mathfrak{t}^{\theta}({\mathbb R})}\chi_{X}(Y)I^{K\tilde{M}}(exp(Y)\eta,\omega,\epsilon_{K\tilde{M}}(f))\,dY$$
 et, pour $j=1,...,r$,
 $$I_{j}=\int_{\mathfrak{t}^{\theta}({\mathbb R})}\partial_{z_{j,T^{\theta}}\varpi^{m(r-j})}\psi_{X}(Y)I^{K\tilde{M}}(exp(Y)\eta,\omega,\epsilon_{K\tilde{M}}(f))\,dY.$$
 Pour tout $j=1,...,r$, la fonction $\partial_{\varpi^{m(r-j)}}\psi_{X}$ est de classe $2m-m_{2}$ (qui est $>0$ par le choix de $m$) et est \`a support compact.  Notons que l'application $Y\mapsto -Y$ fixe $\varpi$ donc aussi les $z_{j,T^{\theta}}$. Par int\'egration par parties,  $I_{j}$ se r\'ecrit
 $$I_{j}=\int_{\mathfrak{t}^{\theta}({\mathbb R})}\partial_{ \varpi^{m(r-j})}\psi_{X}(Y)\partial_{z_{j,T^{\theta}}}I^{K\tilde{M}}(exp(Y)\eta,\omega,\epsilon_{K\tilde{M}}(f))\,dY.$$
Ou encore, gr\^ace \`a 8.2(1)
$$I_{j}=\int_{\mathfrak{t}^{\theta}({\mathbb R})}\partial_{ \varpi^{m(r-j})}\psi_{X}(Y) I^{K\tilde{M}}(exp(Y)\eta,\omega,\epsilon_{K\tilde{M}}(z_{j}f))\,dY.$$
Les fonctions $\xi_{X}$ comme $\partial_{ \varpi^{m(r-j})}\psi_{X}$ appartiennent \`a $C_{c}^{m_{1}}(\mathfrak{u}_{reg})$. Gr\^ace \`a (3), on peut exprimer chaque int\'egrale $I_{j}$ sous la forme (1). La formule (7) devient alors
$$(8) \qquad  I^{K\tilde{M}}(\gamma,\omega,\epsilon_{K\tilde{M}}(f))=\sum_{K\tilde{R}\in {\cal L}(K\tilde{M}_{0})}\sum_{\tilde{\sigma}\in \underline{{\cal E}}_{ell,0}(\tilde{R},\omega)}\sum_{H\in {\cal H}_{K\tilde{R}}^{K\tilde{G},{\cal E}}(\tilde{\sigma})}\int_{H}\sum_{j=0,...,r}\Psi_{K\tilde{R},\tilde{\sigma},H,j}(\lambda)\,d\lambda,$$
 o\`u
 $$\Psi_{K\tilde{R},\tilde{\sigma},H,0}(\lambda)=\phi_{K\tilde{R},\tilde{\sigma},H}(\chi_{X},\lambda)I^{K\tilde{R}}(\tilde{\sigma}_{\lambda},f_{K\tilde{R},\omega})$$
 et, pour $j=1,...,r$,
 $$\Psi_{K\tilde{R},\tilde{\sigma},H,j}(\lambda)=\phi_{K\tilde{R},\tilde{\sigma},H}(\partial_{\varpi^{m(r-j)}}\psi_{X},\lambda)I^{K\tilde{R}}(\tilde{\sigma}_{\lambda},z_{j}f_{K\tilde{R},\omega}).$$
 On sait qu'\`a $\tilde{\sigma}$ est associ\'e un param\`etre infinit\'esimal $\mu(\tilde{\sigma})$ qui est une orbite dans $\mathfrak{h}^*$ pour l'action de $W^{R}$. Identifions-le \`a un point de cette orbite.  En identifiant $\mathfrak{Z}(G)$ \`a l'alg\`ebre des polyn\^omes sur $\mathfrak{h}^*\simeq \mathfrak{t}^*$, on  a l'\'egalit\'e 
 $$I^{K\tilde{R}}(\tilde{\sigma}_{\lambda},z_{j}f_{K\tilde{R},\omega})=z_{j}(\mu(\tilde{\sigma})+\lambda)I^{K\tilde{R}}(\tilde{\sigma}_{\lambda},f_{K\tilde{R},\omega}).$$
 Posons
 $$(9) \qquad \xi_{K\tilde{R},\tilde{\sigma},H}(\gamma,\lambda)=\phi_{K\tilde{R},\tilde{\sigma},H}(\chi_{X},\lambda)+\sum_{j=1,...,r}z_{j}(\mu(\tilde{\sigma}+\lambda)\phi_{K\tilde{R},\tilde{\sigma},H}(\partial_{\varpi^{m(r-j)}}\psi_{X},\lambda).$$
 Alors  
 $$\sum_{j=0,...,r}\Psi_{K\tilde{R},\tilde{\sigma},H,j}(\lambda)=\xi_{K\tilde{R},\tilde{\sigma},H}(\gamma,\lambda)I^{K\tilde{R}}(\tilde{\sigma}_{\lambda},f_{K\tilde{R},\omega})$$
 et la formule (8) devient celle du (iv) de l'\'enonc\'e. 
 
 D'apr\`es (3) et la d\'efinition (9), les fonctions $\lambda\mapsto \xi_{K\tilde{R},\tilde{\sigma},H}(\gamma,\lambda)$ sont de classe $C^{m_{0}}$.  Reprenons la construction en rempla\c{c}ant l'\'el\'ement $X$ par un \'el\'ement $X'$ voisin de $X$. Si $X'$ est assez proche de $X$, on peut utiliser le m\^eme $\epsilon$. Les fonctions $\chi_{X'}$ et $\psi_{X'}$ d\'ependent de $X'$ seulement par translation. En se rappelant la condition impos\'ee \`a $m$, on voit alors que les applications $X'\mapsto \chi_{X'}$ et $X'\mapsto \psi_{X'}$ sont des applications $m_{0}$ fois d\'erivables \`a valeurs dans $C^{m_{1}}( \mathfrak{u}_{reg})$. Pour $U\in Sym(\mathfrak{t}^{\theta})$ de degr\'e au plus $m_{0}$, les images de ces applications par $\partial_{U}$ sont \'egales \`a $X'\mapsto \partial_{U}\chi_{X'}$ et $X'\mapsto \partial_{U}\psi_{X'}$. D'autre part, consid\'erons un ouvert $E$ d'un espace ${\mathbb R}^{a}$, $b\in {\mathbb N}$ un entier et   $e\mapsto \varphi[e]$ une application $b$ fois d\'erivable de $E$ dans $C^{m_{1}}(\mathfrak{u}_{reg})$.  La majoration (3) entra\^{\i}ne que, pour tous $K\tilde{R}$, $\tilde{\sigma}$, $H$ et tout $\lambda\in H$, l'application $e\mapsto \phi_{K\tilde{R},\tilde{\sigma},H}(\varphi[e],\lambda)$ est $b$ fois d\'erivable. Pour un op\'erateur diff\'erentiel $D$ sur ${\mathbb R}^{a}$, \`a coefficients constants et de degr\'e au plus $b$, l'image par $D$ de cette application   est l'application $e\mapsto \phi_{K\tilde{R},\tilde{\sigma},H}(D\varphi[e],\lambda)$. De cela et de la formule (9) r\'esulte que $X'\mapsto \xi_{K\tilde{R},\tilde{\sigma},H}(exp(X')\eta,\lambda)$ est $m_{0}$ fois d\'erivable au point $X$ et que, pour $D\in Sym(\mathfrak{t}^{\theta})$ de degr\'e au plus $m_{0}$, on a l'\'egalit\'e
 $$(10) \qquad D\xi_{K\tilde{R},\tilde{\sigma},H}(exp(X)\eta,\lambda)=\phi_{K\tilde{R},\tilde{\sigma},H}(D\chi_{X},\lambda)$$
 $$+\sum_{j=1,...,r}z_{j}(\mu(\tilde{\sigma})+\lambda)\phi_{K\tilde{R},\tilde{\sigma},H}(\partial_{\varpi^{m(r-j)}}D\psi_{X},\lambda).$$
 On remarque que les valeurs de $\delta(D\chi_{X})$ et $\delta(D\psi_{X})$ sont minor\'ees par $D^{\tilde{G}}(exp(X)\eta)/2$ d'apr\`es le choix de $\epsilon$. 
 En utilisant (3), on en d\'eduit que l'application $(X,\lambda)\mapsto \xi_{K\tilde{R},\tilde{\sigma},H}(exp(X)\eta,\lambda) $ est $m_{0}$ fois d\'erivable en les deux variables. Pour tout op\'erateur diff\'erentiel \`a coefficients constants $\Delta$ sur $H$ de degr\'e au plus $m_{0}$  et tout \'el\'ement $D\in Sym(\mathfrak{t}^{\theta})$ de degr\'e au plus $m_{0}$, il existe  un entier $N\in {\mathbb N}$, une famille finie $(H_{i})_{i=1,...,n}$ d'\'el\'ements de $Sym(\mathfrak{t}^{\theta})$ de degr\'e au plus $m_{3}=2m(r-1)+m_{0}+m_{1}$ et un entier $d\in {\mathbb N}$ de sorte que
 $$(11) \qquad \vert D\Delta\xi_{K\tilde{R},\tilde{\sigma},H}(exp(X)\eta,\lambda)\vert \leq D^{\tilde{G}}(exp(X)\eta)^{-d}(1+\vert \lambda\vert )^N$$
 $$\sum_{i=1,...,n}(sup_{Y\in \mathfrak{t}^{\theta}({\mathbb R})}\vert \partial_{H_{i}}\chi(Y)\vert) +(sup_{Y\in \mathfrak{t}^{\theta}({\mathbb R})}\vert \partial_{H_{i}}\psi(Y)\vert).$$
Notons que $m_{3}< 2mr-m_{2}$ par d\'efinition de $m$. Le terme $(1+\vert \lambda\vert )^N$ s'introduit \`a cause des polyn\^omes $z_{j}(\mu(\tilde{\sigma})+\lambda)$ de la formule (10).  Reprenons la d\'efinition des fonctions $\chi$ et $\psi$. On voit que les termes $\vert \partial_{H_{i}}\chi(Y)\vert $ et $\vert \partial_{H_{i}}\psi(Y)\vert$ sont born\'es par une constante et une puissance n\'egative de $\epsilon$, l'exposant \'etant  au plus $m_{3}$. On a suppos\'e $\epsilon< \frac{c_{1}}{2c_{0}}D^{\tilde{G}}(exp(X)\eta)$ mais on peut aussi bien choisir $\epsilon=\frac{c_{1}}{4c_{0}}D^{\tilde{G}}(exp(X)\eta)$. Les termes ci-dessus sont alors born\'es par une constante et une puissance n\'egative de $D^{\tilde{G}}(exp(X)\eta)$, l'exposant \'etant au plus $m_{3}$. La majoration (11) prend alors la forme du (ii) de l'\'enonc\'e. 
 
 On obtient ainsi une forme affaiblie des assertions (i) et (ii) de l'\'enonc\'e:  pour (i), on prouve seulement que la fonction est $C^{m_{0}}$ en chaque variable; pour (ii), on impose que les degr\'es des op\'erateurs diff\'erentiels sont au plus $m_{0}$.  Mais le principe d'unicit\'e \'enonc\'e en 8.3 entra\^{\i}ne que les fonctions $\xi_{K\tilde{R},\tilde{\sigma},H}$ que l'on a construites sont ind\'ependantes de l'entier $m_{0}$ utilis\'e pour les construire. En faisant varier cet entier, les formes faibles de (i) et (ii) entra\^{\i}nent ces assertions telles qu'\'enonc\'ees dans la proposition. $\square$

$\square$

\bigskip

\subsection{Description des fonctions $\xi_{K\tilde{R},\tilde{\sigma},H}$}
On conserve la situation du paragraphe pr\'ec\'edent. Fixons $\eta\in \tilde{T}({\mathbb R})\cap \tilde{G}_{reg}({\mathbb R})$. On dispose de l'exponentielle
$$exp:\mathfrak{t}^{\theta}({\mathbb R})\to T^{\theta,0}({\mathbb R}).$$
 Son  noyau est un sous-${\mathbb Z}$-module de type fini de $\mathfrak{t}^{\theta}({\mathbb R})$ que l'on note $Ker_{T}$.  Il engendre le sous-espace $X_{*}(T^{\theta,0})^-\otimes i{\mathbb R}$ de $\mathfrak{t}^{\theta}({\mathbb R})$, o\`u l'exposant $^-$ d\'esigne le sous-espace sur lequel  $\Gamma_{{\mathbb R}}$ agit par son unique caract\`ere non trivial.
Notons $\mathfrak{t}'$ le sous-ensemble des \'el\'ements $X\in \mathfrak{t}^{\theta}({\mathbb R})$ tels que $exp(X)\eta$ soit r\'egulier dans $\tilde{M}({\mathbb R})$. C'est le compl\'ementaire d'un ensemble localement fini d'hyperplans affines. Chacun de ces hyperplans est invariant par translations par $\mathfrak{a}_{\tilde{M}}({\mathbb R})$ et ne contient pas $0$ puisque $\eta$ est fortement r\'egulier dans $\tilde{G}$. Cet ensemble d'hyperplans est conserv\'e par translations par $Ker_{T}$.  Pour $X\in \mathfrak{t}'$, notons $n(X)$ le nombres de ces hyperplans qui s\'eparent $X$ de $0$ et posons $\epsilon(X)=(-1)^{n(X)}$. Cela d\'efinit une fonction $\epsilon:\mathfrak{t}'\to \{\pm 1\}$.

Soient $K\tilde{R}\in {\cal L}(K\tilde{M}_{0})$, $\tilde{\sigma}\in \underline{{\cal E}}_{ell,0}(K\tilde{R},\omega)$ et $H\in {\cal H}_{K\tilde{R}}^{K\tilde{G},{\cal E}}(\tilde{\sigma})$. On sait d\'efinir le param\`etre infinit\'esimal $\mu(\tilde{\sigma})$ qui est une orbite dans $\mathfrak{h}^*$ pour l'action de $W$. On a vu en [IV] 1.2 que son intersection avec l'espace affine $\tilde{\mu}(\omega)+
\mathfrak{h}^{\theta,*}$ \'etait une unique orbite pour l'action de $W^{\theta}$. Fixons un point de cette intersection que, pour simplifier, on note encore $\mu(\tilde{\sigma})$. On identifie $\mathfrak{t}^*$ \`a $\mathfrak{h}^*$ et donc $H$ \`a un sous-ensemble de $\mathfrak{h}^*$, en fait de $i\mathfrak{h}_{{\mathbb R}}^*$.  Notons $W_{\tilde{\sigma},H}$ le sous-groupe des \'el\'ements de $W^{\theta}$ qui fixent tout \'el\'ement de $\mu(\tilde{\sigma})+H$. Notons $H'$ le sous-ensemble des $\lambda\in H$ tels que le fixateur de $\mu(\tilde{\sigma})+\lambda$ dans $W^{\theta}$ soit \'egal \`a $W_{\tilde{\sigma},H}$. C'est le compl\'ementaire dans $H$ d'un ensemble fini de sous-espaces affines propres. 

\ass{Proposition}{ (i) Pour tout $w\in W^{\theta}/W_{\tilde{\sigma},H}$, il existe une unique fonction
 $$\begin{array}{ccc}\mathfrak{t}^{\theta}({\mathbb R})\times H'&\to&{\mathbb C}\\ (X,\lambda)&\mapsto& P_{w}(X,\lambda)\\ \end{array}$$
  v\'erifiant les conditions suivantes:

- $P_{w}(X,\lambda)$ est $C^{\infty}$ en $\lambda$ et polynomiale en $X$ de degr\'e inf\'erieur ou \'egal \`a $\vert W_{\tilde{\sigma},H}\vert $;

- pour $X\in \mathfrak{t}^{\theta}({\mathbb R})$ tel que $exp(X)\eta$ soit fortement r\'egulier dans $\tilde{G}$ et pour $\lambda\in H'$, on a l'\'egalit\'e
$$\xi_{K\tilde{R},\tilde{\sigma},H}(exp(X)\eta,\lambda)=\epsilon(X)\sum_{w\in W^{\theta}/W_{\tilde{\sigma},H}}e^{<X,w(\mu(\tilde{\sigma})+\lambda)>}P_{w}(X,\lambda).$$

(ii) On a $\xi_{K\tilde{R},\tilde{\sigma},H}=0$ si $dim(H)>dim( {\cal A}_{\tilde{M}})$.}

Preuve.  On oublie nos donn\'ees fix\'ees $K\tilde{R}$, $\tilde{\sigma}$ et $H$ qui ne joueront pas de r\^ole particulier, ce qui lib\`ere ces symboles. Notons $\mathfrak{t}''$ le sous-ensemble des \'el\'ements $X\in \mathfrak{t}^{\theta}({\mathbb R})$ tels que $exp(X)\eta$ soit fortement r\'egulier dans $\tilde{G}$. C'est un sous-ensemble de $\mathfrak{t}'$ qui est le compl\'ementaire dans $\mathfrak{t}^{\theta}({\mathbb R})$ d'un ensemble localement fini de sous-espaces propres. Rappelons que l'on dispose de 
   l'homomorphisme d'Harish-Chandra $z\mapsto z_{T^{\theta,0}}$ de $\mathfrak{Z}(G)$ dans 
   $Sym(\mathfrak{t})_{\theta,\omega}^{W^{\theta}}\simeq Sym(\mathfrak{t}^{\theta})^{W^{\theta}}$.  On a l'\'egalit\'e
$$I^{K\tilde{M}}(exp(X)\eta,\omega,z_{M}\varphi)=\partial_{z_{T^{\theta,0}}}I^{K\tilde{M}}(exp(X)\eta,\omega,\varphi)$$
pour tout $z\in  \mathfrak{Z}(G)$, tout  $X\in \mathfrak{t}''$ et tout $\varphi\in I(K\tilde{M}({\mathbb R}),\omega)$. Soient $z\in  \mathfrak{Z}(G)$, $X\in \mathfrak{t}''$ et $f\in I(\tilde{G}({\mathbb R}),\omega)$. On a l'\'egalit\'e $\epsilon_{K\tilde{M}}(zf)=z_{M}\epsilon_{K\tilde{M}}(f)$, cf. 8.2(1). On a donc
$$I^{K\tilde{M}}(exp(X)\eta,\omega,\epsilon_{K\tilde{M}}(zf))=\partial_{z_{T^{\theta,0}}}I^{K\tilde{M}}(exp(X)\eta,\omega,\epsilon_{K\tilde{M}}(f)).$$
Exprimons les deux termes \`a l'aide du (iv) de la proposition 8.7. On obtient pour chacun d'eux une somme en $K\tilde{R}$, $\tilde{\sigma}$, $H$  d'int\'egrales
$$\int_{H}\partial_{z_{T^{\theta,0}}}\xi_{K\tilde{R},\tilde{\sigma},H}(exp(X)\eta,\lambda)I^{K\tilde{R}}(\tilde{\sigma}_{\lambda},f_{K\tilde{R},\omega})\,d\lambda$$
pour le membre de droite et
$$\int_{H}\xi_{K\tilde{R},\tilde{\sigma},H}(exp(X)\eta,\lambda)I^{K\tilde{R}}(\tilde{\sigma}_{\lambda},(zf)_{K\tilde{R},\omega})\,d\lambda$$
pour celui de gauche. On a  
$$I^{K\tilde{R}}(\tilde{\sigma}_{\lambda},(zf)_{K\tilde{R},\omega})=I^{K\tilde{R}}(\tilde{\sigma}_{\lambda},z_{R}f_{K\tilde{R},\omega})=z(\mu(\tilde{\sigma})+\lambda)I^{K\tilde{R}}(\tilde{\sigma}_{\lambda},f_{K\tilde{R},\omega}).$$
On obtient ainsi des expressions similaires, avec une fonction $\partial_{z_{T^{\theta,0}}}\xi_{K\tilde{R},\tilde{\sigma},H}(exp(X)\eta,\lambda)$ dans le membre de droite et une fonction $z(\mu(\tilde{\sigma})+\lambda)\xi_{K\tilde{R},\tilde{\sigma},H}(exp(X)\eta,\lambda)$. Ces fonctions forment deux familles index\'ees par $\tilde{R}$, $\tilde{\sigma}$, $H$. Ces deux familles v\'erifient la condition de sym\'etrie (2) de 8.4. Mais elles calculent le m\^eme terme. Comme on l'a dit en 8.4, elles sont donc \'egales. On obtient l'\'egalit\'e
$$(1) \qquad \partial_{z_{T^{\theta,0}}}\xi_{K\tilde{R},\tilde{\sigma},H}(exp(X)\eta,\lambda)=z(\mu(\tilde{\sigma})+\lambda)\xi_{K\tilde{R},\tilde{\sigma},H}(exp(X)\eta,\lambda).$$
L'alg\`ebre $Sym(\mathfrak{t}^{\theta})$ est un module de type fini sur l'image de $\mathfrak{Z}(G)$. On peut en fixer une base $(U_{i})_{i=1,...,n}$. Pour $U\in Sym(\mathfrak{t}^{\theta})$, \'ecrivons $U=\sum_{i=1,...,n}z_{i,T^{\theta,0}}U_{i}$ avec des $z_{i}\in \mathfrak{Z}(G)$. Alors
$$\partial_{U}\xi_{K\tilde{R},\tilde{\sigma},H}(exp(X)\eta,\lambda)=\sum_{i=1,...,n}\partial_{U_{i}}\partial_{z_{i,T^{\theta,0}}}\xi_{K\tilde{R},\tilde{\sigma},H}(exp(X)\eta,\lambda)$$
$$=\sum_{i=1,...,n}z_{i}(\mu(\tilde{\sigma})+\lambda)\partial_{U_{i}}\xi_{K\tilde{R},\tilde{\sigma},H}(exp(X)\eta,\lambda).$$
On en d\'eduit une premi\`ere am\'elioration du (ii) de la proposition 8.7:

(2) pour tout sous-ensemble compact $\Gamma\subset \mathfrak{t}^{\theta}({\mathbb R})$, il existe $d\in {\mathbb N}$  tel que, pour tout $U\in Sym(\mathfrak{t}^{\theta})$, il existe un entier $N\in {\mathbb N}$ et une constante $c>0$ de sorte que
$$\vert \partial_{U}\xi_{K\tilde{R},\tilde{\sigma},H}(exp(X)\eta,\lambda)\vert \leq c\vert D^{\tilde{G}}(exp(X)\eta)\vert ^{-d}(1+\vert \lambda\vert )^{N}$$
pour tout $X\in \Gamma\cap \mathfrak{t}''$ et tout $\lambda\in H$.  

L'entier $d$ est devenu ind\'ependant de $U$. En appliquant le lemme 8.1, on obtient la seconde am\'elioration suivante:

(3)  pour tout sous-ensemble compact $\Gamma\subset \mathfrak{t}^{\theta}({\mathbb R})$ et tout $U\in Sym(\mathfrak{t}^{\theta})$, il existe un entier $N\in {\mathbb N}$ et une constante $c>0$ de sorte que
$$\vert \partial_{U}\xi_{K\tilde{R},\tilde{\sigma},H}(exp(X)\eta,\lambda)\vert \leq c (1+\vert \lambda\vert )^{N}$$
pour tout $X\in \Gamma\cap \mathfrak{t}''$ et tout $\lambda\in H$.

Consid\'erons pour quelques instants que $\lambda$ est fix\'e  et qu'il appartient \`a $H'$. Notons simplement $\xi(X)=\xi_{K\tilde{R},\tilde{\sigma},H}(exp(X)\eta,\lambda)$.    L'\'egalit\'e (1) est une \'equation diff\'erentielle portant sur la fonction $\xi$. En faisant varier $z$, on obtient un syst\`eme d'\'equations diff\'erentielles. On conna\^{\i}t gr\^ace \`a [Va] th\'eor\`eme 11 la forme des solutions.   Dans tout ouvert connexe o\`u $\xi$ est $C^{\infty}$,  $\xi$ s'\'ecrit
$$\xi(X)=\sum_{w\in W^{\theta}/W_{\tilde{\sigma},H}}e^{<X,w(\mu(\tilde{\sigma})+\lambda)>}P_{w}(X),$$
o\`u $P_{w}$ est un polyn\^ome en $X$ de degr\'e au plus $\vert W_{\tilde{\sigma},H}\vert $. Ces polyn\^omes sont uniquement d\'etermin\'es. Evidemment, quand on consid\`ere de nouveau $\lambda$ comme variable, les polyn\^omes d\'ependent de $\lambda$. On v\'erifie ais\'ement  par interpolation que, puisque $\xi_{K\tilde{R},\tilde{\sigma},H}(exp(X)\eta,\lambda)$ est $C^{\infty}$ en $\lambda$, les polyn\^omes obtenus sont $C^{\infty}$ en tout point $\lambda$ en position g\'en\'erale. Par ailleurs, on peut aussi bien remplacer $\xi_{K\tilde{R},\tilde{\sigma},H}(exp(X)\eta,\lambda)$ par $\epsilon(X)\xi_{K\tilde{R},\tilde{\sigma},H}(exp(X)\eta,\lambda)$ puisque $\epsilon$ est localement constante sur $\mathfrak{t}''$.   En r\'esum\'e, pour toute composante connexe $\Omega$ de $\mathfrak{t}''$, il existe d'uniques fonctions $P^{\Omega}_{w}(X,\lambda)$, pour $w\in W^{\theta}/W_{\tilde{\sigma},H}$, qui sont polynomiales en $X$ de degr\'e au plus $\vert W_{\tilde{\sigma},H}\vert $ et qui sont $C^{\infty}$ en tout point $\lambda\in H'$, de sorte que
$$(4) \qquad \epsilon(X)\xi_{K\tilde{R},\tilde{\sigma},H}(exp(X)\eta,\lambda)=\sum_{w\in W^{\theta}/W_{\tilde{\sigma},H}}e^{<X,w(\mu(\tilde{\sigma})+\lambda)>}P^{\Omega}_{w}(X,\lambda)$$
pour tout $X\in \Omega$ et tout $\lambda\in H'$. Notons $\xi_{K\tilde{R},\tilde{\sigma},H}^{\Omega}(exp(X)\eta,\lambda)$ le membre de droite ci-dessus. On remarque que la fonction $X\mapsto \xi_{K\tilde{R},\tilde{\sigma},H}^{\Omega}(exp(X)\eta,\lambda)$ s'\'etend naturellement en une fonction $C^{\infty}$ sur $\mathfrak{t}^{\theta}({\mathbb R})$ d\'efinie par la m\^eme formule.

On va prouver que cette fonction ne d\'epend pas de $\Omega$. Le compl\'ementaire de $\mathfrak{t}''$ est une r\'eunion de sous-espaces propres. Mais les sous-espaces de codimension au moins $2$ ne cr\'eent pas de disconnexit\'e. Il nous suffit de prouver l'\'egalit\'e $\xi_{K\tilde{R},\tilde{\sigma},H}^{\Omega}(exp(X)\eta,\lambda)=\xi_{K\tilde{R},\tilde{\sigma},H}^{\Omega'}(exp(X)\eta,\lambda)$ quand $\Omega$ et $\Omega'$ sont deux composantes connexes s\'epar\'ees 
  par un unique hyperplan singulier.  Consid\'erons de telles $\Omega$ et $\Omega'$. Consid\'erons  un \'el\'ement $X_{0}$ de  l'hyperplan singulier qui s\'epare ces deux composantes, qui est en position g\'en\'erale dans cet hyperplan et qui appartient aux adh\'erences de $\Omega$ et $\Omega'$. On a d\'efini en 8.2 un signe $  \Delta_{exp(X_{0})\eta}(Y)$ pour un \'el\'ement $Y$ g\'en\'eral et proche de $0$.   On v\'erifie sur la d\'efinition que le rapport $\Delta_{exp(X_{0})\eta}(Y)\epsilon(X_{0}+Y)^{-1}$ est constant pour $Y$ g\'en\'eral et proche de $0$. 
Soit  $f\in I(K\tilde{G}({\mathbb R}),\omega,K)$. Parce que $\epsilon_{K\tilde{M}}(f)$ est une fonction cuspidale, la fonction $Y\mapsto \Delta_{exp(X_{0})\eta}(Y)I^{K\tilde{M}}(exp(Y+X_{0})\eta,\omega,\epsilon_{K\tilde{M}}(f))$ est $C^{\infty}$ au voisinage de $0$.   Ou encore la fonction $X\mapsto \epsilon(X)I^{K\tilde{M}}(exp(X)\eta,\omega,\epsilon_{K\tilde{M}}(f))$ est $C^{\infty}$ au voisinage de $X_{0}$. Soit $U\in Sym(\mathfrak{t}^{\theta})$. Pour $X\in \mathfrak{t}''$ proche de $X_{0}$, $\epsilon(X)\partial_{U}I^{K\tilde{M}}(exp(X)\eta,\omega,\epsilon_{K\tilde{M}}(f))$ se calcule par la formule 8.7(iv). On obtient une somme en $K\tilde{R}$, $\tilde{\sigma}$, $H$ d'int\'egrales
$$\epsilon(X)\int_{H}\partial_{U}\xi_{K\tilde{R},\tilde{\sigma},H}(exp(X)\eta,\lambda)I^{K\tilde{R}}(\tilde{\sigma}_{\lambda},f_{K\tilde{R},\omega})\, d\lambda.$$
Supposons $X\in \Omega$. Pour $\lambda\in H'$, on peut remplacer $\epsilon(X)\partial_{U}\xi_{K\tilde{R},\tilde{\sigma},H}(exp(X)\eta,\lambda)$ par $\partial_{U}\xi_{K\tilde{R},\tilde{\sigma},H}^{\Omega}(exp(X)\eta,\lambda)$.  Faisons tendre $X$ vers $X_{0}$. Pour tout $\lambda\in H'$,  $\partial_{U}\xi_{K\tilde{R},\tilde{\sigma},H}^{\Omega}(exp(X)\eta,\lambda)$ tend vers $\partial_{U}\xi_{K\tilde{R},\tilde{\sigma},H}^{\Omega}(exp(X_{0})\eta,\lambda)$. L'assertion (3) montre que $\partial_{U}\xi_{K\tilde{R},\tilde{\sigma},H}^{\Omega}(exp(X)\eta,\lambda)$ reste uniform\'ement born\'ee par le produit d'une constante et de $(1+\vert \lambda\vert )^{N}$ pour un entier $N$ convenable. On peut donc appliquer le th\'eor\`eme de convergence domin\'ee. La fonction  $\partial_{U}\xi_{K\tilde{R},\tilde{\sigma},H}^{\Omega}(exp(X_{0})\eta,\lambda)$ est encore born\'ee par une puissance de $1+\vert \lambda\vert $ et la limite quand $X$ tend vers $0$ de  $\epsilon(X)\partial_{U}I^{K\tilde{M}}(exp(X)\eta,\omega,\epsilon_{K\tilde{M}}(f))$  est la somme en $K\tilde{R}$, $\tilde{\sigma}$, $H$ des  int\'egrales
$$ \int_{H}\partial_{U}\xi_{K\tilde{R},\tilde{\sigma},H}^{\Omega}(exp(X_{0})\eta,\lambda)I^{K\tilde{R}}(\tilde{\sigma}_{\lambda},f_{K\tilde{R},\omega})\, d\lambda.$$
 On peut \'evidemment refaire le calcul en rempla\c{c}ant $\Omega$ par $\Omega'$. On obtient une expression similaire. Les deux expressions obtenues v\'erifient la condition de sym\'etrie 8.4(3) (parce que ce sont des limites d'expressions qui la v\'erifient). Il en r\'esulte que ces deux expressions co\"{\i}ncident. C'est-\`a-dire que, pour tous $K\tilde{R}$, $\tilde{\sigma}$, $H$ et pour tout $U$, on a l'\'egalit\'e
$$\partial_{U}\xi_{K\tilde{R},\tilde{\sigma},H}^{\Omega}(exp(X_{0})\eta,\lambda)=\partial_{U}\xi_{K\tilde{R},\tilde{\sigma},H}^{\Omega'}(exp(X_{0})\eta,\lambda).$$
Puisque les deux fonctions $\xi_{K\tilde{R},\tilde{\sigma},H}^{\Omega}(exp(X)\eta,\lambda)$ et $\xi_{K\tilde{R},\tilde{\sigma},H}^{\Omega'}(exp(X)\eta,\lambda)$ sont des exponentielles-polyn\^omes, l'\'egalit\'e de toutes leurs d\'eriv\'ees en un point entra\^{\i}ne leur \'egalit\'e. Donc, pour $\lambda\in H'$, $\xi_{K\tilde{R},\tilde{\sigma},H}^{\Omega}(exp(X)\eta,\lambda)$ ne d\'epend pas de $\Omega$.  Les termes $P_{w}^{\Omega}(X,\lambda)$ non plus. Mais alors, l'assertion (4) d\'emontre le (i) de l'\'enonc\'e. 

 Consid\'erons de nouveau des donn\'ees $K\tilde{R}$, $\tilde{\sigma}$ et $H$ fix\'ees. On voit facilement que la  fonction $\epsilon$ v\'erifie la condition de p\'eriodicit\'e: pour $Y\in Ker_{T}$, $\epsilon(X)/\epsilon(X+Y)$ ne d\'epend pas de $X\in \mathfrak{t}'$. Il en r\'esulte que $\epsilon(X+2Y)=\epsilon(X)$ pour $X\in \mathfrak{t}'$ et $Y\in Ker_{T}$.  Puisqu'on a aussi $exp(X+2Y)=exp(X)$, la formule du (i) de l'\'enonc\'e entra\^{\i}ne
 $$\sum_{w\in W^{\theta}/W_{\tilde{\sigma},H}}e^{<X+2Y,w(\mu(\tilde{\sigma})+\lambda)>}P_{w}(X+2Y,\lambda)=\sum_{w\in W^{\theta}/W_{\tilde{\sigma},H}}e^{<X,w(\mu(\tilde{\sigma})+\lambda)>}P_{w}(X,\lambda).$$
 Cette \'egalit\'e pour tout $X\in \mathfrak{t}'$ et tout $\lambda\in H'$ entra\^{\i}ne que
 $$(5) \qquad e^{<2Y,w(\mu(\tilde{\sigma})+\lambda)>}P_{w}(X+2Y,\lambda)=P_{w}(X,\lambda)$$
 pour tous $w\in W^{\theta}$, $X\in \mathfrak{t}'$, $\lambda\in H'$ et $Y\in Ker_{T}$. Si $P_{w}(X,\lambda)$ est identiquement nul pour tout $w$, la conclusion du (ii) de l'\'enonc\'e est claire. Sinon, fixons $w$, $X$ et $\lambda_{0}$ tels que $P_{w}(X,\lambda)\not=0$. Pour ces valeurs de $w$ et  $X$ et pour  $\lambda$, l'\'egalit\'e (5) dit qu'une fonction exponentielle co\"{\i}ncide avec une fonction rationnelle sur le ${\mathbb Z}$-module $Ker_{T}$. Il en r\'esulte ais\'ement que cette exponentielle est constante sur ce r\'eseau, ce qui implique que $<2Y,w(\mu(\tilde{\sigma})+\lambda)>\in 2\pi i{\mathbb Z}$ pour tout $Y\in Ker_{T}$. On se rappelle que $H$ est un sous-espace  affine de $i{\cal A}_{\tilde{R}}^*$. Notons $H^0$ le sous-espace lin\'eaire associ\'e. Les seules conditions impos\'ees \`a $\lambda$ sont $\lambda\in H'$ et $P_{w}(X,\lambda)\not=0$. Elles d\'efinissent un ouvert non vide et la relation ci-dessus est v\'erifi\'ee pour tout $\lambda$ dans cet ouvert.   On peut donc  remplacer $\lambda$ par  $\lambda+\nu$, o\`u $\nu$ est un \'el\'ement de $H^0$ assez voisin de $0$. On en d\'eduit que
 $<2Y,w(\nu)>\in 2\pi i{\mathbb Z}$ pour tout $Y\in Ker_{T}$ et tout $\nu\in H^0$ assez proche de $0$. Cela entra\^{\i}ne facilement que $<Y,w(\nu)>=0$ pour tout $\nu\in H^0\otimes_{{\mathbb R}}{\mathbb C}$ et tout $Y$ dans le ${\mathbb C}$-sous-espace de $\mathfrak{t}^{\theta}({\mathbb R})$  engendr\'e par $Ker_{T}$. Comme on l'a dit, ce sous-espace est $X_{*}(T^{\theta,0})^-\otimes_{{\mathbb Z}}{\mathbb C}$. Parce que $\tilde{T}$ est un sous-tore tordu elliptique de $\tilde{M}$, l'orthogonal de ce sous-espace est 
 $\mathfrak{a}_{\tilde{M}}^*$. La relation pr\'ec\'edent signifie que  $w(H^0\otimes_{{\mathbb R}}{\mathbb C})$ est inclus    dans  $\mathfrak{a}_{\tilde{M}}^*$. Cela implique que $dim(H)\leq dim(\mathfrak{a}_{\tilde{M}}({\mathbb R}))=dim({\cal A}_{\tilde{M}})$. Cela prouve le (ii) de l'\'enonc\'e. $\square$

\bigskip

\subsection{$K$-finitude,}
 
On r\'ealise l'espace de Paley-Wiener $PW^{\infty}(K\tilde{G},\omega)$ et son sous-espace $PW(K\tilde{G},\omega)$ en utilisant les bases  $\underline{{\cal E}}_{ell,0}(K\tilde{R},\omega)$ des espaces $D_{ell,0}(K\tilde{R},\omega)$ pour tout $K\tilde{R}\in {\cal L}(K\tilde{M}_{0})$. Fixons un tel $K\tilde{R}$ et un \'el\'ement $\tilde{\sigma}\in \underline{{\cal E}}_{ell,0}(K\tilde{R},\omega)$. Notons ${\cal F} $ l'espace des fonctions de Paley-Wiener sur ${\cal A}_{\tilde{R},{\mathbb C}}^*$.  On l'identifie au sous-espace des \'el\'ements  $(\varphi_{\tilde{\sigma}'})_{\tilde{\sigma}'\in \underline{{\cal E}}_{ell,0}(K\tilde{R},\omega)}\in PW_{ell}(K\tilde{R},\omega)$ dont toutes les composantes sont nulles sauf celle pour $\tilde{\sigma}'=\tilde{\sigma}$. On dispose alors de  l'homomorphisme de sym\'etrisation  $sym:{\cal F}\to {\cal F}^{W(\tilde{R})}\to PW(\tilde{G},\omega)$ et de  l'isomorphisme de Paley-Wiener $pw:  I(\tilde{G}({\mathbb R}),\omega,K)\to PW(\tilde{G},\omega)$ (cf.  [IV] 2.4 et [W2] 6.1). Pour $\varphi\in {\cal F}$, on pose $f_{\varphi}=pw^{-1}\circ sym(\varphi)$.   

\ass{Proposition}{(i) Supposons $dim({\cal A}_{\tilde{R}})< dim({\cal A}_{\tilde{M}})$. Alors $\epsilon_{K\tilde{M}}(f_{\varphi})=0$ pour tout $\varphi\in {\cal F}$.

(ii) Supposons $dim({\cal A}_{\tilde{R}})\geq dim({\cal A}_{\tilde{M}})$. Alors
$\epsilon_{K\tilde{M}}(f_{\varphi})$ est $K$-finie, c'est-\`a-dire appartient \`a $I_{ac,cusp}(K\tilde{M}({\mathbb R}),\omega,K)$,  pour tout $\varphi\in {\cal F}$.}

Preuve.   Fixons $\tilde{\pi}\in \underline{{\cal E}}_{ell,0}(K\tilde{M},\omega)$. Soit $\varphi\in {\cal F}$. Comme on l'a dit en 8.3(2), on dispose d'une fonction m\'eromorphe $\lambda\mapsto   I^{K\tilde{M}}(\tilde{\pi},\lambda,\epsilon_{\tilde{M}}(f_{\varphi}))$ sur ${\cal A}_{\tilde{M},{\mathbb C}}^*$. Montrons que

(1) pour $\lambda\in {\cal A}_{\tilde{M},{\mathbb C}}^*$ et $z\in \mathfrak{Z}(M)$, on a l'\'egalit\'e
$$ I^{K\tilde{M}}(\tilde{\pi},\lambda,z\epsilon_{\tilde{M}}(f_{\varphi}))=z(\mu(\tilde{\pi})+\lambda) I^{K\tilde{M}}(\tilde{\pi},\lambda,\epsilon_{\tilde{M}}(f_{\varphi})).$$

 Consid\'erons les espaces ${\cal S}(i{\cal A}_{\tilde{M}}^*)$ et ${\cal S}({\cal A}_{\tilde{M}})$ des fonctions de Schwartz sur $i{\cal A}_{\tilde{M}}^*$ et ${\cal A}_{\tilde{M}}$. On a la transformation de Fourier $\psi\mapsto \hat{\psi}$ de ${\cal S}(i{\cal A}_{\tilde{M}}^*)$ sur  ${\cal S}({\cal A}_{\tilde{M}})$. L'alg\`ebre $\mathfrak{Z}(M)$ est isomorphe \`a l'alg\`ebre des polyn\^omes sur $\mathfrak{h}^*$ invariants par $W^M$. On la fait agir sur ${\cal S}(i{\cal A}_{\tilde{M}}^*)$ par $(z\psi)(\lambda)=z(\mu(\tilde{\pi})+\lambda)\psi(\lambda)$. Par transformation de Fourier, on obtient une action sur ${\cal S}({\cal A}_{\tilde{M}})$ que l'on note $\rho$. C'est-\`a-dire que l'on a $(z\psi)^{\hat{}}=\rho(z)\hat{\psi}$. Il est facile d'expliciter cette action.   On d\'ecompose $\mathfrak{h}^*=\mathfrak{h}^{M,*}\oplus \mathfrak{a}_{\tilde{M}}$. Cela d\'ecompose $\mathfrak{Z}(M)$ en produit tensoriel de l'alg\`ebre des polyn\^omes sur $\mathfrak{h}^{M,*}$ invariants par $W^M$ et une alg\`ebre isomorphe \`a $Sym(\mathfrak{a}_{\tilde{M}})$. Pour $z$ dans la premi\`ere alg\`ebre, $\rho(z)$ est la multiplication par $z(\mu(\tilde{\pi}))$. Pour $X\in \mathfrak{a}_{\tilde{M}}$, l'action $\rho(X)$ est la d\'erivation $\partial_{X}$.  Soit $h\in I(K\tilde{M}({\mathbb R}),\omega)$. La fonction $X\mapsto I^{K\tilde{M}}(\tilde{\pi},X,h)$ est la transform\'ee de Fourier de $\lambda\mapsto I^{K\tilde{M}}(\tilde{\pi}_{\lambda},h)$. On sait que, pour $z\in \mathfrak{Z}(M)$, on a l'\'egalit\'e
$$ I^{K\tilde{M}}(\tilde{\pi}_{\lambda},zh)=z(\mu(\tilde{\pi})+\lambda) I^{K\tilde{M}}(\tilde{\pi}_{\lambda},h).$$
On en d\'eduit que
$$(2) \qquad I^{K\tilde{M}}(\tilde{\pi},X,zh)=\rho(z)I^{K\tilde{M}}(\tilde{\pi},X,h).$$
Mais cette formule se g\'en\'eralise \`a $h\in I_{ac}(K\tilde{M}({\mathbb R}),\omega)$. En effet, fixons $X$, puis une fonction $b\in C_{c}^{\infty}({\cal A}_{\tilde{M}})$ qui vaut $1$ dans un voisinage de $X$. On a
$$(3) \qquad I^{K\tilde{M}}(\tilde{\pi},X,z(h(b\circ H_{\tilde{M}})))=\rho(z)I^{K\tilde{M}}(\tilde{\pi},X,h(b\circ H_{\tilde{M}})).$$
Il r\'esulte de la description de l'action $\rho$ que $\rho(z)I^{K\tilde{M}}(\tilde{\pi},X,h(b\circ H_{\tilde{M}}))$ ne d\'epend que des valeurs de $I^{K\tilde{M}}(\tilde{\pi},X',h(b\circ H_{\tilde{M}}))$ pour $X'$ proche de $X$. Or, pour de tels $X'$, on a  $I^{K\tilde{M}}(\tilde{\pi},X',h(b\circ H_{\tilde{M}}))= I^{K\tilde{M}}(\tilde{\pi},X',h)$. D'autre part, la diff\'erence $ z(h(b\circ H_{\tilde{M}}))-(zh)(b\circ H_{\tilde{M}})$ est nulle en un point $\gamma$ tel que $H_{\tilde{M}}(\gamma)$ est proche de $X$. On en d\'eduit 
$$I^{K\tilde{M}}(\tilde{\pi},X,z(h(b\circ H_{\tilde{M}})))=I^{K\tilde{M}}(\tilde{\pi},X,(zh)(b\circ H_{\tilde{M}}))=I^{K\tilde{M}}(\tilde{\pi},X,zh).$$
L'\'egalit\'e (3) est donc identique \`a (2). En particulier, l'\'egalit\'e (2) est v\'erifi\'ee pour $h=\epsilon_{\tilde{M}}(f_{\varphi})$. Par inversion de Fourier, on en d\'eduit l'\'egalit\'e de l'assertion (1).

Soit $\nu\in {\cal A}_{\tilde{M},{\mathbb C}}^*$. D\'efinissons une forme lin\'eaire $l_{\nu}$ sur ${\cal F}$ par $l_{\nu}(\varphi)=I^{K\tilde{M}}(\tilde{\pi},\nu,\epsilon_{\tilde{M}}(f_{\varphi}))$.
 Les assertions (1) ci-dessus et 8.3(1) entra\^{\i}nent que $l_{\nu}(z\varphi)=z(\mu(\tilde{\pi})+\nu)l_{\nu}(\varphi)$. En notant $J_{\nu}$ l'id\'eal des \'el\'ements $z\in \mathfrak{Z}(G)$ tels que $z(\mu(\tilde{\pi})+\nu)=0$, on obtient que $l_{\nu}$ annule $J_{\nu}{\cal F}$. Rappelons que $ \mathfrak{Z}(G)$ agit sur ${\cal F}$ par $(z\varphi)(\lambda)=z(\mu(\tilde{\sigma})+\lambda)\varphi(\lambda)$. Supposons que $W(\mu(\tilde{\pi})+\nu)\cap( \mu(\tilde{\sigma})+{\cal A}_{\tilde{R},{\mathbb C}}^*)=\emptyset$. Le lemme [IV] 2.5 entra\^{\i}ne alors que $J_{\nu}{\cal F}={\cal F}$, donc $l_{\nu}$ est nulle.  Supposons  que $W(\mu(\tilde{\pi})+\nu)\cap (\mu(\tilde{\sigma})+{\cal A}_{\tilde{R},{\mathbb C}}^*)$ soit non vide, notons $\lambda_{1}$,...,$\lambda_{m}$ les \'el\'ements de ${\cal A}_{\tilde{R},{\mathbb C}}^*$ tels que cette intersection soit $\{\mu(\tilde{\sigma})+\lambda_{j};j=1,...,m\}$. Le m\^eme lemme entra\^{\i}ne l'existence d'un entier $N\geq0$ tel que $J_{\nu}{\cal F}$ contienne tous les \'el\'ements de ${\cal F}$ qui s'annulent \`a l'ordre $N$ en chaque $\lambda_{j}$.  Pour chaque espace $H\in {\cal H}_{K\tilde{R}}^{K\tilde{G},{\cal E}}(\tilde{\sigma})$, \'ecrivons $H=i\mu_{H}+iV_{H}$, o\`u $V_{H}$ est un sous-espace de ${\cal A}_{\tilde{R}}^*$ et $\mu_{H}\in {\cal A}_{\tilde{R}}^*$ est orthogonal \`a $V_{H}$ et  posons $H_{{\mathbb C}}=i\mu_{H}+V_{H,{\mathbb C}}$. Notons ${\cal H}_{K\tilde{R}}^{K\tilde{G},{\cal E}}(\tilde{\sigma},\leq a_{\tilde{M}})$ le sous-ensemble des $H\in {\cal H}_{K\tilde{R}}^{K\tilde{G},{\cal E}}(\tilde{\sigma})$ tels que $dim(H)\leq a_{\tilde{M}}=dim({\cal A}_{\tilde{M}})$. Montrons que
 
 (4) supposons que,  pour tout $j=1,...,m$ et tout $H\in {\cal H}_{K\tilde{R}}^{K\tilde{G},{\cal E}}(\tilde{\sigma},\leq a_{\tilde{M}})$,  $\lambda_{j}$ n'apparient pas \`a $ H_{{\mathbb C}}$; alors $l_{\nu}=0$. 
 
 Sous l'hypoth\`ese de (4), on peut trouver un polyn\^ome $Q$ sur ${\cal A}_{\tilde{R},{\mathbb C}}^*$ tel que $Q-1$ s'annule \`a l'ordre $N$ en tout $\lambda_{j}$ et $Q$ s'annule sur tout \'el\'ement de ${\cal H}_{K\tilde{R}}^{K\tilde{G},{\cal E}}(\tilde{\sigma},\leq a_{\tilde{M}})$. Ces conditions \'etant invariantes par $W(\tilde{R})$, on peut supposer $Q$ invariant par ce groupe.
 Pour $\varphi\in {\cal F}$, on a $\varphi=(1-Q)\varphi+Q\varphi$. Le premier terme appartient \`a $J_{\nu}{\cal F}$ donc est annul\'e par $l_{\nu}$.  On va montrer que
 
 (5) $\epsilon_{K\tilde{M}}(f_{Q\varphi})=0$.
 
 Cela entra\^{\i}ne a fortiori $l_{\nu}(Q\varphi)=0$, ce qui prouve (4). Prouvons (5). Si $K\tilde{M}$
 ne poss\`ede pas de sous-tore tordu maximal elliptique, toute fonction elliptique sur $K\tilde{M}({\mathbb R})$ est nulle, d'o\`u (5). Sinon, fixons un tel 
   sous-tore tordu maximal elliptique $\tilde{T}$ de $K\tilde{M}$. Pour $\gamma\in \tilde{T}({\mathbb R})\cap \tilde{G}_{reg}({\mathbb R})$, on va prouver que   $I^{K\tilde{M}}(\gamma,\omega,\epsilon_{K\tilde{M}}(f_{Q\varphi}))=0$, ce qui d\'emontrera (5). Cette int\'egrale orbitale est calcul\'ee par la proposition 8.7. Compte tenu de la d\'efinition de ${\cal F}$,  seules les paires conjugu\'ees \`a la paire fix\'ee $(K\tilde{R},\tilde{\sigma})$ interviennent dans les deux premi\`eres sommes.   Compte tenu des invariances par conjugaison de nos diff\'erents objets, on voit que, pour prouver la nullit\'e souhait\'ee, il suffit de prouver que
   $$\xi_{K\tilde{R},\tilde{\sigma},H}(\gamma,\lambda)I^{K\tilde{R}}(\tilde{\sigma}_{\lambda},(f_{Q\varphi})_{K\tilde{R},\omega})=0$$
   pour tout $H\in {\cal H}_{K\tilde{R}}^{K\tilde{G},{\cal E}}(\tilde{\sigma})$ et tout $\lambda\in H$ en position g\'en\'erale.
   On a
   $$I^{K\tilde{R}}(\tilde{\sigma}_{\lambda},(f_{Q\varphi})_{K\tilde{R},\omega})=\vert W(\tilde{R})\vert ^{-1}\sum_{w\in W(\tilde{R})}Q(w\lambda)\varphi(w\lambda).$$
   D'apr\`es l'invariance de $Q$, il suffit de prouver que
 $$\xi_{K\tilde{R},\tilde{\sigma},H}(\gamma,\lambda)Q(\lambda) =0$$
 pour tout $H\in {\cal H}_{K\tilde{R}}^{K\tilde{G},{\cal E}}(\tilde{\sigma})$ et tout $\lambda\in H$ en position g\'en\'erale. Or, si $dim(H)\leq a_{\tilde{M}}$, $Q(\lambda)=0$. Si $dim(H)> a_{\tilde{M}}$, c'est la premi\`ere fonction qui est nulle d'apr\`es la proposition 8.8(ii).  Cela prouve (5) et (4).
  
 Il r\'esulte de (4) que, pour que $l_{\nu}$ soit non nulle, il est n\'ecessaire qu'il existe $H\in {\cal H}_{K\tilde{R}}^{K\tilde{G},{\cal E}}(\tilde{\sigma},\leq a_{\tilde{M}})$ de sorte que 
 $$W(\mu(\tilde{\pi})+\nu)\cap (\mu(\tilde{\sigma})+i\mu_{H}+V_{H,{\mathbb C}})\not=\emptyset.$$
Notons $E$ l'ensemble des $\nu\in {\cal A}_{\tilde{M},{\mathbb C}}^*$ v\'erifiant cette condition. Supposons dor\'enavant que

(6) la fonction $(X,\varphi)\mapsto I^{K\tilde{M}}(\tilde{\pi},X,\epsilon_{K\tilde{M}}(f_{\varphi}))$ sur ${\cal A}_{\tilde{M}}\times {\cal F}$ est non nulle. 

Il en est de m\^eme de la fonction $(\nu,\varphi)\mapsto I^{K\tilde{M}}(\tilde{\pi},\nu,\varphi)$ sur $i{\cal A}_{\tilde{M}}^*\times {\cal F}$. Puisque cette fonction est $C^{\infty}$ en $\nu$, cela entra\^{\i}ne que $l_{\nu}$ est non nulle pour $\nu$ dans un ouvert non vide de $i{\cal A}_{\tilde{M}}^*$. Un tel ouvert est donc contenu dans $E$. A fortiori
 $E$ n'est pas inclus dans une r\'eunion finie d'hyperplans  affines de ${\cal A}_{\tilde{M},{\mathbb C}}^*$. 
Rempla\c{c}ons les orbites $\mu(\tilde{\pi})$ et $\mu(\tilde{\sigma})$ par des points dans ces orbites. L'ensemble $E$ est la r\'eunion sur $H\in {\cal H}_{K\tilde{R}}^{K\tilde{G},{\cal E}}(\tilde{\sigma},\leq a_{\tilde{M}})$, $w\in W$ et $w'\in W^M$ des ensembles
 $$E_{H,w,w'}= (w(\mu(\tilde{\sigma})+i\mu_{H}+ V_{H,{\mathbb C}})-w'(\mu(\tilde{\pi})))\cap {\cal A}_{\tilde{M},{\mathbb C}}^*.$$
 On peut donc fixer $H$, $w$ et $w'$ tels que $E_{H,w,w'}$ ne soit pas inclus dans un hyperplan affine de ${\cal A}_{\tilde{M},{\mathbb C}}^*$. Mais $E_{H,w,w'}$ est clairement contenu dans un tel hyperplan, sauf si c'est ${\cal A}_{\tilde{M},{\mathbb C}}^*$ tout entier. Donc
 $$(7) \qquad {\cal A}_{\tilde{M},{\mathbb C}}^*\subset (w(\mu(\tilde{\sigma})+i\mu_{H}+ V_{H,{\mathbb C}})-w'(\mu(\tilde{\pi}))).$$
 Cette condition implique  ${\cal A}_{\tilde{M},{\mathbb C}}^*\subset w( V_{H,{\mathbb C}})$. En vertu de l'hypoth\`ese $dim(H)\leq a_{\tilde{M}}$, cette inclusion est une \'egalit\'e. L'\'egalit\'e (7) entra\^{\i}ne alors 
 $wÕ(\mu(\tilde{\pi}))=w(\mu(\tilde{\sigma})+i\mu_{H})$. En consid\'erant de nouveau $\mu(\tilde{\pi})$ comme une orbite sous l'action de $W^M$, cette \'egalit\'e impose \`a $\mu(\tilde{\pi})$ d'appartenir \`a un ensemble fini d'orbites d\'etermin\'e par $\tilde{\sigma}$. La condition (6)  implique donc
 
(8) il existe $H\in {\cal H}_{K\tilde{R}}^{K\tilde{G},{\cal E}}(\tilde{\sigma})$ tel que $dim(H)=a_{\tilde{M}}$;

(9) l'orbite $\mu(\tilde{\pi})$ appartient \`a un ensemble fini dÕorbites d\'etermin\'e par $\tilde{\sigma}$. 

Si $dim({\cal A}_{\tilde{R}})<dim({\cal A}_{\tilde{M}})$, la condition (8) n'est jamais v\'erifi\'ee. En g\'en\'eral, la condition (9) n'est v\'erifi\'ee que pour un ensemble fini de $\tilde{\pi}\in \underline{{\cal E}}_{ell,0}(\tilde{M},\omega)$. Les m\^emes assertions valent donc pour la condition (6). 
   Cela prouve la proposition. $\square$
 
 La proposition entra\^{\i}ne imm\'ediatement le corollaire suivant. 

 \ass{Corollaire}{ On a $\epsilon_{K\tilde{M}}(f)\in I_{ac,cusp}(K\tilde{M}({\mathbb R}),\omega,K)$ pour tout $f\in C_{c}^{\infty}(\tilde{G}({\mathbb R}),K)$.}

\bigskip

 {\bf Bibliographie}
 
 [A1] J. Arthur: {\it The local behaviour of weighted orbital integrals}, Duke Math. J. 56 (1988), p. 223-293
 
 [A2] -----------: {\it Stabilization of a family of differential equations}, Proc. of Symp. in pure Math. 68 (2000)
 
 [A3] -----------: {\it The characters of discrete series as orbital integrals}, Inventiones Math. 32 (1976), p.205-261
 
 [A4] -----------: {\it Intertwining operators and residues I. Weighted characters}, J. of Funct. Analysis 84 (1989), p. 19-84
 
 [A5] -----------: {\it Canonical normalization of weighted characters and a transfer conjecture}, C. R. Math. Acad. Sci. Soc. R. Can. 20 (1998), p. 33-52
 
 [A6] -----------: {\it The trace formula in invariant form}, Annals of Math. 114 (1981), p. 1-74
 
 [A7] -----------: {\it On a family of distributions obtained from Eisenstein series II: explicit formulas}, Amer. J. of Math. 104 (1981), p. 1289- 1336
 
 [A8] -----------: {\it Parabolic transfer for real groups}, J. AMS 21 (2008), p. 171-234
 
 [A9] -----------: {\it On the Fourier transforms of weighted orbital integrals}, J. reine angew. Math. 452 (1994), p. 163-217
 
 [B] A. Bouaziz: {\it Int\'egrales orbitales sur les groupes de Lie r\'eductifs}, Ann. Sc. ENS 27 (1994), p. 573-609
 
 [HC] Harish-Chandra: {\it The characters of semi-simple Lie groups}, Trans. AMS 83 (1956), p. 98-163
 
 [K] R. Kottwitz: { \it Stable trace formula: elliptic singular terms}, Math. Ann. 275 (1986), p. 365-399
 
 [KS] -------------, D. Shelstad: {Foundations of twisted endoscopy}, Ast\'erisque 255 (1999)
 
 [L] J.-P. Labesse: {\it Stable twisted trace formula: elliptic terms}, Journal of the Inst. of Math. Jussieu 3 (2004), p. 473-530
 
 [LW] ---------------, J.-L. Waldspurger: {\it La formule des traces tordue d'apr\`es le Friday Morning Seminar}, CRM Monograph series 31 (2013)
 
 [Me] P. Mezo: {\it Spectral transfer in the twisted endoscopy for real groups}, pr\'epublication 2013
 
 [Moe] C. Moeglin: {\it Appendice \`a la formule des traces locale}, pr\'epublication 2013
 
 [S] D. Shelstad: {\it On geometric transfer in real twisted endoscopy}, Annals of Math. 176 (2012), p. 1919-1985 
 
 [V] V. Varadarajan: {\it Harmonic analysis on real reductive groups}, Lecture Notes in Math 576, Springer 1977
 
 [W1] J.-L. Waldspurger: {\it L'endoscopie tordue n'est pas si tordue}, Memoirs AMS 908 (2008)
 
 [W2] ------------------------: {\it La formule des traces locale tordue}, pr\'epublication 2012
 
 [I], [II], [III], [IV], [V], [VIII] ------------------------: {\it  Stabilisation de la formule des traces tordue I: endoscopie tordue sur un corps local}, {\it II: int\'egrales orbitales et endoscopie sur un corps local non-archim\'edien;  d\'efinitions et \'enonc\'es des r\'esultats}, {\it III: int\'egrales orbitales et endoscopie sur un corps local non-archim\'edien; r\'eductions et preuves}, {\it IV: transfert spectral archim\'edien}, {\it V: int\'egrales orbitales et endoscopie sur ${\mathbb R}$}, {\it VIII: l'application $\epsilon_{\tilde{M}}$ sur un corps local non-archim\'edien}, pr\'epublications 2014
 
      \bigskip
     
     Institut de Math\'ematiques de Jussieu
     
     2 place Jussieu, 75005 Paris
     
     e-mail: jean-loup.waldspurger@imj-prg.fr
 
 \end{document}